\documentclass[11pt]{article}
\usepackage[letterpaper,margin=1.05in,bottom=1.65in]{geometry}

\usepackage[utf8]{inputenc} 
\usepackage[T1]{fontenc}    
\usepackage[hidelinks]{hyperref}       
\usepackage{url}            
\usepackage{booktabs}       
\usepackage{amsfonts}       
\usepackage{nicefrac}       
\usepackage{microtype}      
\usepackage{xcolor}         
\usepackage[affil-it]{authblk} 
\usepackage{amsmath}
\usepackage{amssymb}
\usepackage{amsthm}
\usepackage{mathtools}
\mathtoolsset{showonlyrefs,showmanualtags}
\usepackage{mathabx} 
\allowdisplaybreaks
\usepackage{bbm}
\usepackage{enumitem}
\usepackage{graphicx}
\usepackage[ruled]{algorithm2e}
\usepackage{caption}
\usepackage{subcaption}
\usepackage{float}

\newcommand{\figpath}{./figs}

\DeclarePairedDelimiter\paren\lparen\rparen
\DeclarePairedDelimiter\bracket\lbrack\rbrack
\DeclarePairedDelimiter\braces\lbrace\rbrace

\DeclarePairedDelimiter\abs\lvert\rvert

\providecommand{\bbone}{\mathbf{1}}
\DeclarePairedDelimiterXPP\indicator[1]{\bbone}{\lbrack}{\rbrack}{}{#1}

\DeclarePairedDelimiterXPP\expf[1]{\exp}{\lparen}{\rparen}{}{#1}
\DeclarePairedDelimiterXPP\logf[1]{\log}{\lparen}{\rparen}{}{#1}
\DeclarePairedDelimiterXPP\maxf[1]{\max}{\lparen}{\rparen}{}{#1}
\DeclarePairedDelimiterXPP\minf[1]{\min}{\lparen}{\rparen}{}{#1}

\DeclareMathOperator*{\sgn}{sgn}
\DeclarePairedDelimiterXPP\sgnf[1]{\sgn}{\lparen}{\rparen}{}{#1}

\DeclarePairedDelimiterXPP\func[2]{#1}{\lparen}{\rparen}{}{#2}

\DeclareMathOperator*{\atan}{atan}
\DeclarePairedDelimiterXPP\atanf[1]{\atan}{\lparen}{\rparen}{}{#1}
\DeclarePairedDelimiterXPP\tanf[1]{\tan}{\lparen}{\rparen}{}{#1}

\DeclarePairedDelimiter\setb\lbrace\rbrace

\newcommand{\Reals}{\mathbb{R}}

\newcommand{\Naturals}{\mathbb{N}}

\DeclareMathOperator*{\argmin}{arg\,min}

\newcommand{\mat}[1]{\boldsymbol{#1}}
\renewcommand{\vec}[1]{\boldsymbol{#1}}

\newcommand{\unitM}{\mat{I}}

\makeatletter
\newcommand*{\tran}{{\mathpalette\@tran{}}}
\newcommand*{\@tran}[2]{\raisebox{\depth}{$\m@th#1\intercal$}}
\makeatother

\DeclarePairedDelimiter\iprod\langle\rangle

\DeclarePairedDelimiter\norm\lVert\rVert
\DeclarePairedDelimiterXPP\tnorm[1]{}{\lVert}{\rVert_{1}}{}{#1}
\DeclarePairedDelimiterXPP\enorm[1]{}{\lVert}{\rVert_{2}}{}{#1}
\DeclarePairedDelimiterXPP\inorm[1]{}{\lVert}{\rVert_{\infty}}{}{#1}
\DeclarePairedDelimiterXPP\pnorm[2]{}{\lVert}{\rVert_{#1}}{}{#2}
\DeclarePairedDelimiterXPP\opnorm[1]{}{\lVert}{\rVert_{op}}{}{#1}

\DeclarePairedDelimiterXPP\detf[1]{\det}{\lparen}{\rparen}{}{#1}

\DeclarePairedDelimiterXPP\pospart[1]{}{\lparen}{\rparen_+}{}{#1}

\DeclarePairedDelimiterXPP\kerf[1]{\ker}{\lparen}{\rparen}{}{#1}

\DeclareMathOperator{\trsym}{tr}
\DeclarePairedDelimiterXPP\tr[1]{\trsym}{\lparen}{\rparen}{}{#1}

\DeclareMathOperator{\vspansym}{span}
\DeclarePairedDelimiterXPP\vspan[1]{\vspansym}{\lparen}{\rparen}{}{#1}

\DeclareMathOperator{\diagsym}{diag}
\DeclarePairedDelimiterXPP\diag[1]{\diagsym}{\lparen}{\rparen}{}{#1}

\DeclareMathOperator{\ranksym}{rank}
\DeclarePairedDelimiterXPP\rank[1]{\ranksym}{\lparen}{\rparen}{}{#1}

\DeclareMathOperator{\vectorizesym}{vec}
\DeclarePairedDelimiterXPP\vectorize[1]{\vectorizesym}{\lparen}{\rparen}{}{#1}

\DeclareMathOperator*{\esssup}{ess\,sup}
\DeclarePairedDelimiterXPP\esssupf[1]{\esssup}{\lparen}{\rparen}{}{#1}

\let\Prsym\Pr
\let\Pr\relax
\DeclarePairedDelimiterXPP\Pr[1]{\Prsym}{\lparen}{\rparen}{}{%
	#1}

\DeclarePairedDelimiterXPP\Prsub[2]{\Prsym_{#1}}{\lparen}{\rparen}{}{%
	#2}

\DeclareMathOperator{\Esym}{E}
\DeclarePairedDelimiterXPP\E[1]{\Esym}{\lbrack}{\rbrack}{}{%
	#1}

\DeclarePairedDelimiterXPP\Esub[2]{\Esym_{#1}}{\lbrack}{\rbrack}{}{%
	#2}

\DeclareMathOperator{\Varsym}{Var}
\DeclarePairedDelimiterXPP\Var[1]{\Varsym}{\lparen}{\rparen}{}{%
	#1}

\DeclarePairedDelimiterXPP\Varsub[2]{\Varsym_{#1}}{\lparen}{\rparen}{}{%
	#2}

\DeclarePairedDelimiterXPP\EstVar[1]{\widehat{\Varsym}}{\lparen}{\rparen}{}{%
	#1}

\DeclareMathOperator{\Covsym}{Cov}
\DeclarePairedDelimiterXPP\Cov[1]{\Covsym}{\lparen}{\rparen}{}{%
	#1}

\DeclarePairedDelimiterXPP\Covsub[2]{\Covsym_{#1}}{\lparen}{\rparen}{}{%
	#2}

\DeclareMathOperator{\Corrsym}{Corr}
\DeclarePairedDelimiterXPP\Corr[1]{\Corrsym}{\lparen}{\rparen}{}{%
	#1}

\makeatletter
\newcommand{\indep}{\protect\mathpalette{\protect\@indep}{\perp}}
\newcommand*{\@indep}[2]{\mathrel{\rlap{$#1#2$}\mkern3mu{#1#2}}}
\makeatother

\newcommand{\bigOsym}{\mathcal{O}}
\DeclarePairedDelimiterXPP\bigO[1]{\bigOsym}{\lparen}{\rparen}{}{#1}
\DeclarePairedDelimiterXPP\bigOt[1]{\widetilde{\bigOsym}}{\lparen}{\rparen}{}{#1}

\newcommand{\littleOsym}{o}
\DeclarePairedDelimiterXPP\littleO[1]{\littleOsym}{\lparen}{\rparen}{}{#1}

\newcommand{\bigOpsym}{\bigOsym_p}
\DeclarePairedDelimiterXPP\bigOp[1]{\bigOpsym}{\lparen}{\rparen}{}{#1}

\newcommand{\littleOpsym}{\littleOsym_p}
\DeclarePairedDelimiterXPP\littleOp[1]{\littleOpsym}{\lparen}{\rparen}{}{#1}

\newcommand{\bigOmegasym}{\Omega}
\DeclarePairedDelimiterXPP\bigOmega[1]{\bigOmegasym}{\lparen}{\rparen}{}{#1}

\newcommand{\littleOmegasym}{\omega}
\DeclarePairedDelimiterXPP\littleOmega[1]{\littleOmegasym}{\lparen}{\rparen}{}{#1}

\newcommand{\bigThetasym}{\Theta}
\DeclarePairedDelimiterXPP\bigTheta[1]{\bigThetasym}{\lparen}{\rparen}{}{#1}

\DeclarePairedDelimiterXPP\cosf[1]{\cos}{\lparen}{\rparen}{}{#1}
\DeclarePairedDelimiterXPP\sinf[1]{\sin}{\lparen}{\rparen}{}{#1}

\newcommand{\quadtext}[1]{\quad\text{#1}\quad}

\newcommand{\quadand}{\quadtext{and}}

\newcommand{\newvar}[2]{
	\expandafter\newcommand\csname #1\endcsname{#2}
}

\newcommand{\newvars}[2]{
	\expandafter\newcommand\csname #1\endcsname[1]{#2_{##1}}
}

\newcommand{\newvarss}[2]{
	\expandafter\newcommand\csname #1\endcsname[2]{#2_{##1,##2}}
}

\newcommand{\newvarsss}[2]{
	\expandafter\newcommand\csname #1\endcsname[3]{#2_{##1,##2,##3}}
}

\newcommand{\newvarssss}[2]{
	\expandafter\newcommand\csname #1\endcsname[4]{#2_{##1,##2,##3,##4}}
}

\newcommand{\newfunconly}[2]{
	\expandafter\DeclarePairedDelimiterXPP\csname #1\endcsname[1]{#2}{\lparen}{\rparen}{}{##1}
}

\newcommand{\newfuncs}[2]{
	\expandafter\DeclarePairedDelimiterXPP\csname #1\endcsname[2]{#2_{##1}}{\lparen}{\rparen}{}{##2}
}

\newcommand{\newfuncss}[2]{
	\expandafter\DeclarePairedDelimiterXPP\csname #1\endcsname[3]{#2_{##1,##2}}{\lparen}{\rparen}{}{##3}
}

\newcommand{\newfuncsss}[2]{
	\expandafter\DeclarePairedDelimiterXPP\csname #1\endcsname[4]{#2_{##1,##2,##3}}{\lparen}{\rparen}{}{##4}
}

\newcommand{\newfuncssss}[2]{
	\expandafter\DeclarePairedDelimiterXPP\csname #1\endcsname[5]{#2_{##1,##2,##3,##4}}{\lparen}{\rparen}{}{##5}
}

\newcommand{\varpartialeval}[3]{
	\expandafter\newcommand\csname #2\endcsname{\csname #1\endcsname{#3}}
}

\newcommand{\funcpartialeval}[3]{
	\expandafter\newcommand\csname #2\endcsname[1][]{\csname #1\endcsname[##1]{#3}}
}

\newcommand{\varevali}[2][]{
	\varpartialeval{#2#1}{#2i}{i}
	\varpartialeval{#2#1}{#2j}{j}
}

\newcommand{\varevalk}[2][]{
	\varpartialeval{#2#1}{#2k}{k}
	\varpartialeval{#2#1}{#2l}{\ell}
}

\newcommand{\varevals}[2][]{
	\varpartialeval{#2#1}{#2s}{s}
	\varpartialeval{#2#1}{#2r}{r}
}

\newcommand{\funcevali}[2][]{
	\funcpartialeval{#2#1}{#2i}{i}
	\funcpartialeval{#2#1}{#2j}{j}
}

\newcommand{\funcevalk}[2][]{
	\funcpartialeval{#2#1}{#2k}{k}
	\funcpartialeval{#2#1}{#2l}{\ell}
}

\newcommand{\varevalii}[2][]{
	\varevali[#1]{#2}
	\varevali{#2i}
	\varevali{#2j}
}

\newcommand{\varevalik}[2][]{
	\varevali[#1]{#2}
	\varevalk{#2i}
	\varevalk{#2j}
}

\newcommand{\varevaliik}[2][]{
	\varevalii[#1]{#2}
	\varevalk{#2ii}
	\varevalk{#2ij}
	\varevalk{#2ji}
	\varevalk{#2jj}
}

\newcommand{\varevaliis}[2][]{
	\varevalii[#1]{#2}
	\varevals{#2ii}
	\varevals{#2ij}
	\varevals{#2ji}
	\varevals{#2jj}
}

\newcommand{\varevaliikk}[2][]{
	\varevaliik[#1]{#2}
	\varevalk{#2iik}
	\varevalk{#2ijk}
	\varevalk{#2jik}
	\varevalk{#2jjk}
	\varevalk{#2iil}
	\varevalk{#2ijl}
	\varevalk{#2jil}
	\varevalk{#2jjl}
}

\newcommand{\funcevalii}[2][]{
	\funcevali[#1]{#2}
	\funcevali{#2i}
	\funcevali{#2j}
}

\newcommand{\funcevalik}[2][]{
	\funcevali[#1]{#2}
	\funcevalk{#2i}
	\funcevalk{#2j}
}

\newcommand{\funcevaliik}[2][]{
	\funcevalii[#1]{#2}
	\funcevalk{#2ii}
	\funcevalk{#2ij}
	\funcevalk{#2ji}
	\funcevalk{#2jj}
}

\newcommand{\funcevaliikk}[2][]{
	\funcevaliik[#1]{#2}
	\funcevalk{#2iik}
	\funcevalk{#2ijk}
	\funcevalk{#2jik}
	\funcevalk{#2jjk}
	\funcevalk{#2iil}
	\funcevalk{#2ijl}
	\funcevalk{#2jil}
	\funcevalk{#2jjl}
}

\newcommand{\funcevalX}[3]{
	\expandafter\newcommand\csname #1v#2\endcsname{\csname #1v\endcsname{#3}}
}

\newcommand{\funcevaliX}[3]{
	\expandafter\newcommand\csname #1v#2\endcsname{\csname #1v\endcsname{#3}}
	\expandafter\newcommand\csname #1v#2e\endcsname[1]{\csname #1ve\endcsname{##1}{#3}}
	\expandafter\newcommand\csname #1v#2i\endcsname{\csname #1vi\endcsname{#3}}
	\expandafter\newcommand\csname #1v#2j\endcsname{\csname #1vj\endcsname{#3}}
}

\newcommand{\funcevalikX}[3]{
	\expandafter\newcommand\csname #1v#2\endcsname{\csname #1v\endcsname{#3}}
	\expandafter\newcommand\csname #1v#2e\endcsname[2]{\csname #1ve\endcsname{##1}{##2}{#3}}
	\expandafter\newcommand\csname #1v#2i\endcsname[1]{\csname #1vi\endcsname{##1}{#3}}
	\expandafter\newcommand\csname #1v#2j\endcsname[1]{\csname #1vj\endcsname{##1}{#3}}
	\expandafter\newcommand\csname #1v#2ik\endcsname{\csname #1vik\endcsname{#3}}
	\expandafter\newcommand\csname #1v#2il\endcsname{\csname #1vil\endcsname{#3}}
	\expandafter\newcommand\csname #1v#2jk\endcsname{\csname #1vik\endcsname{#3}}
	\expandafter\newcommand\csname #1v#2jl\endcsname{\csname #1vjl\endcsname{#3}}
}

\newcommand{\newvari}[2]{
	\newvar{#1}{#2}
	\newvars{#1e}{\csname #1\endcsname}
	\varevali[e]{#1}
}

\newcommand{\newvark}[2]{
	\newvar{#1}{#2}
	\newvars{#1e}{\csname #1\endcsname}
	\varevalk[e]{#1}
}

\newcommand{\newvarik}[2]{
	\newvar{#1}{#2}
	\newvarss{#1e}{\csname #1\endcsname}
	\varevalik[e]{#1}
}

\newcommand{\newvarii}[2]{
	\newvar{#1}{#2}
	\newvarss{#1e}{\csname #1\endcsname}
	\varevalii[e]{#1}
}

\newcommand{\newvariik}[2]{
	\newvar{#1}{#2}
	\newvarsss{#1e}{\csname #1\endcsname}
	\varevaliik[e]{#1}
}

\newcommand{\newvariis}[2]{
	\newvar{#1}{#2}
	\newvarsss{#1e}{\csname #1\endcsname}
	\varevaliis[e]{#1}
}

\newcommand{\newvariikk}[2]{
	\newvar{#1}{#2}
	\newvarssss{#1e}{\csname #1\endcsname}
	\varevaliikk[e]{#1}
}

\newcommand{\newfunc}[2]{
	\newvar{#1}{#2}
	\newfunconly{#1v}{\csname #1\endcsname}
}

\newcommand{\newfunci}[2]{
	\newvari{#1}{#2}
	\newfunconly{#1v}{\csname #1\endcsname}
	\newfuncs{#1ve}{\csname #1\endcsname}
	\funcevali[e]{#1v}
}

\newcommand{\newfunck}[2]{
	\newvark{#1}{#2}
	\newfunconly{#1v}{\csname #1\endcsname}
	\newfuncs{#1ve}{\csname #1\endcsname}
	\funcevalk[e]{#1v}
}

\newcommand{\newfuncik}[2]{
	\newvarik{#1}{#2}
	\newfunconly{#1v}{\csname #1\endcsname}
	\newfuncss{#1ve}{\csname #1\endcsname}
	\funcevalik[e]{#1v}
}

\newcommand{\newfuncii}[2]{
	\newvarii{#1}{#2}
	\newfunconly{#1v}{\csname #1\endcsname}
	\newfuncss{#1ve}{\csname #1\endcsname}
	\funcevalii[e]{#1v}
}

\newcommand{\newfunciik}[2]{
	\newvariik{#1}{#2}
	\newfunconly{#1v}{\csname #1\endcsname}
	\newfuncsss{#1ve}{\csname #1\endcsname}
	\funcevaliik[e]{#1v}
}

\newcommand{\newfunciikk}[2]{
	\newvariikk{#1}{#2}
	\newfunconly{#1v}{\csname #1\endcsname}
	\newfuncssss{#1ve}{\csname #1\endcsname}
	\funcevaliikk[e]{#1v}
}

\theoremstyle{plain}
\newtheorem{theorem}{Theorem}[section]

\newtheorem{corollary}[theorem]{Corollary}
\newtheorem{lemma}[theorem]{Lemma}
\newtheorem{proposition}[theorem]{Proposition}

\newtheorem{innerrefcorollary}{Corollary}
\newenvironment{refcorollary}[1]
{\renewcommand\theinnerrefcorollary{#1}\innerrefcorollary}
{\endinnerrefcorollary}

\newtheorem{innerreflemma}{Lemma}
\newenvironment{reflemma}[1]
{\renewcommand\theinnerreflemma{#1}\innerreflemma}
{\endinnerreflemma}

\newtheorem{innerrefproposition}{Proposition}
\newenvironment{refproposition}[1]
{\renewcommand\theinnerrefproposition{#1}\innerrefproposition}
{\endinnerrefproposition}

\newtheorem{innerreftheorem}{Theorem}
\newenvironment{reftheorem}[1]
{\renewcommand\theinnerreftheorem{#1}\innerreftheorem}
{\endinnerreftheorem}

\newtheorem{innerrefcondition}{Condition}
\newenvironment{refcondition}[1]
{\renewcommand\theinnerrefcondition{#1}\innerrefcondition}
{\endinnerrefcondition}

\theoremstyle{definition}
\newtheorem{assumption}{Assumption}
\newtheorem{condition}{Condition}
\newtheorem{definition}{Definition}

\theoremstyle{remark}

\newcommand{\xv}{\vec{x}}
\newcommand{\xM}{\mat{X}}

\newcommand{\designsym}{D}
\newcommand{\tarbdesign}{\designsym_{t}}
\newcommand{\fullarbdesign}{\mathcal{\designsym}}

\newcommand{\bv}{\vec{\beta}}
\newcommand{\optbv}{\bv^*}
\newcommand{\woptbv}{\widetilde{\bv}}

\newcommand{\olsres}[1]{\mathcal{E}(#1)}
\newcommand{\estolsres}[1]{\widehat{\mathcal{E}}(#1)}
\newcommand{\sqolsres}[1]{A(#1)}
\newcommand{\estsqolsres}[1]{\widehat{A}(#1)}

\newcommand{\ate}{\tau}

\newcommand{\filt}{\mathcal{F}}

\newcommand{\eate}{\hat{\ate}}

\newcommand{\regret}{\mathcal{R}}
\newcommand{\neymanregret}{\regret^{\textrm{Neyman}}}

\newcommand{\probloss}[1]{f_{#1}}
\newcommand{\regretprob}{\regret^{\mathrm{prob}}}
\newcommand{\newpredloss}[1]{ \ell_{#1} }
\newcommand{\newpredregret}{\regret^{\mathrm{pred}}}
\newcommand{\sigloss}[1]{h_{#1}}
\newcommand{\estsigloss}[1]{\widehat{h}_{#1}}

\newcommand{\newestpredloss}[1]{\widehat{\ell}_{#1}}
\newcommand{\newtpredloss}[1]{L_{#1} }
\newcommand{\newptpredloss}[1]{\widetilde{L}_{#1} }

\newcommand{\tprobloss}[1]{H_{#1}}
\newcommand{\ptprobloss}[1]{\widetilde{H}_{#1}}

\newcommand{\evb}{\widehat{\vb}}
\newcommand{\vb}{\textrm{VB}}
\newcommand{\optvar}{\textrm{V}^*}

\newcommand{\ourdesign}{\textsc{Sigmoid-FTRL}{}}

\newcommand{\ores}[2]{ A_{#1}(#2) }
\newcommand{\estores}[2]{ \widehat{A}_{#1}(#2) }

\newcommand{\breg}[3]{\mathcal{B}_{#1}( #2 \vert #3 )}
\newcommand{\bregf}[1]{\mathcal{B}_{#1}}

\newcommand{\grad}{\nabla}

\newcommand{\ci}{\widehat{\mathrm{CI}}_{\alpha}}

\newcommand{\bb}{b}
\newcommand{\e}{\mathrm{E}}

\newcounter{spacesave}
\usepackage[style=authoryear-comp,sorting=nyt,sortcites=false, natbib=true, backend=bibtex]{biblatex}

\usepackage[toc,page,header]{appendix}
\usepackage{minitoc}

\usepackage{expcmd}
\newcommand\mainref\ref
\newcommand\suppref\ref


\title{Sigmoid-FTRL: Design-Based \\ Adaptive Neyman Allocation for AIPW Estimators}
\author[1]{Fangyi Chen}
\author[2]{Shu Ge}
\author[3]{Jian Qian}
\author[4]{Christopher Harshaw} 
\affil[1]{Harvard University}
\affil[2]{Massachusetts Institute of Technology}
\affil[3]{The University of Hong Kong}
\affil[4]{Columbia University}
\date{\today}

\begin{document}
	
	\makeatletter%
	
	\begin{NoHyper}\gdef\@thefnmark{}\@footnotetext{\hspace{-1em}We thank 
P.M. Aronow,
Alexander Rakhlin,
Fredrik S{\"a}vje,
and
Stefan Wager
for insightful discussions which helped to shape this work.
Christopher Harshaw gratefully acknowledges support from Foundations of Data Science Institute (FODSI) NSF grant DMS2023505 and NSF grant MMS2316335.}\end{NoHyper}%
	\makeatother%
	
	\maketitle
	\thispagestyle{empty}
	
	\abstract{
		We consider the problem of Adaptive Neyman Allocation for the class of AIPW estimators for estimating the average treatment effect in a design-based setting, where potential outcomes and covariates are deterministic.
As each subject arrives, an adaptive procedure must select both a treatment assignment probability and a pair of linear predictors to be used in the AIPW estimator.
Our goal is to construct an adaptive procedure that minimizes the Neyman Regret, which is the difference between the variance of the adaptive procedure and an oracle variance which uses the optimal non-adaptive choice of assignment probabilities and linear predictors.
While previous work has drawn fruitful connections between Neyman Regret and online convex optimization for the Horvitz--Thompson estimator, one of the central challenges for the AIPW estimator is that the underlying optimization is non-convex.
In this paper, we propose Sigmoid-FTRL, an adaptive experimental design which addresses the non-convexity via simultaneous minimization of two convex regrets.
We prove that under standard regularity conditions, the Neyman Regret of Sigmoid-FTRL converges at a $\sqrt{d/T}$ rate, where $d$ is the dimension of the covariates and $T$ is the number of subjects in the experiment.
We establish that this $T$ dependence is optimal by proving a $T^{-1/2}$ lower bound on the Neyman Regret for any adaptive design under our regularity conditions.
Finally, we establish a central limit theorem and a consistently conservative variance estimator which facilitate the construction of asymptotically valid Wald-type confidence intervals.
	}

	\newpage
	
	\pagenumbering{roman}
	
	\doparttoc 
	\faketableofcontents 
	\part{} 
	\parttoc 
	
	\newpage
	
	\pagenumbering{arabic}
	
	\section{Introduction} \label{sec:intro}

Randomized experiments are used to investigate causal effects in virtually all of the social sciences, from economics and political science to sociology and public health.
In a classical randomized experiment, the experimental design does not depend on the observed outcomes: the subjects enter the study, treatment is assigned, and only afterwards are the outcomes observed and the treatment effects consequently estimated.
In recent years, there has been a growing interest in adaptive randomized experiments, where subjects arrive sequentially and the experimenter can incorporate previously observed outcomes in the experimental design.

In this paper, we study the problem of Adaptive Neyman Allocation for AIPW estimators in the design-based framework.
In this context, an adaptive experimental design must select both the treatment assignment probability and the linear predictors used in the AIPW estimator.
Roughly speaking, the goal of Adaptive Neyman Allocation is to construct an adaptive experiment design under which the variance of an effect estimator is nearly equal to its optimal variance under the best non-adaptive design that has oracle access to all potential outcomes.
We focus specifically on constructing an adaptive experiment design that minimizes this difference, which is known as the Neyman Regret \citep{Dai2023Clip, Kato2025Efficient}.
A formal description of the problem is deferred to Section~\ref{sec:adaptive-neyman-allocation}.

We focus on the design-based framework, where potential outcomes and covariates of each subject are considered to be deterministic and treatment assignment is the sole source of randomness.
The design-based framework stands in contrast to a super-population framework where subjects are assumed to be independent and identically distributed draws from an unknown distribution.
The super-population assumption may be difficult to interpret and justify in settings where subjects were not literally randomly selected into the study; moreover, it precludes the possibility of drift or otherwise systematic change in the potential outcomes of subjects over time.
For these reasons, the design-based framework is sometimes seen as more robust and assumption-lean \citep{Harshaw2025Why}.

The paper which is most closely related to ours is \citet{Dai2023Clip}, who consider Adaptive Neyman Allocation for the unadjusted Horvitz--Thompson estimator.
They introduce Clip-OGD, an experimental design based on online gradient descent with probability clipping, which guarantees that the Neyman Regret converges at a rate of $T^{-1/2} \exp(\sqrt{\log(T)})$.
Their results are based upon a connection between Adaptive Neyman Allocation and online convex optimization, which we further explore here.

In the offline setting, AIPW estimators leverage available covariate information to gain precision improvements over the unadjusted Horvitz--Thompson estimator \citep{Lin2013Agnostic, Lei2020Regression}.
In the online setting, an adaptive design which approaches the optimal AIPW performance will also have significant precision improvements compared to one that targets the unadjusted Horvitz--Thompson estimator.
In this sense, the shift from Horvitz--Thompson to AIPW is not merely of theoretical interest, but also offers practical improvements in statistical precision.
Nevertheless, this shift presents several theoretical challenges: for example, the optimization problem underlying Neyman Regret for AIPW estimation is non-convex, which precludes the possibility of directly using techniques from online convex optimization.

The main contributions of this paper are as follows:
\begin{itemize}
	\item \textbf{Adaptive Design and Analysis}: We present \ourdesign{}, a new adaptive experimental design under which the AIPW Neyman Regret converges at a rate of $\sqrt{d / T}$, where $d$ is the covariate dimension and $T$ is the number of subjects.
	To overcome the issue of non-convexity, the design simultaneously minimizes two convex objectives corresponding to the selection of treatment assignment probability and linear predictors, respectively.
	We derive a $T^{-1/2}$ lower bound which demonstrates that the design has optimal $T$ dependence under our regularity assumptions.
	In order to obtain these rates, \ourdesign{} employs a sigmoidal transformation of the domain, which may be of independent interest to the online optimization community.
	
	\item \textbf{Inferential Methods}: We provide a central limit theorem for the AIPW estimator under \ourdesign{}, which requires further technical developments. We also construct a consistent estimator for Neyman's variance bound. Together, these enable the development of Wald-type intervals that asymptotically cover at least at the nominal level.
\end{itemize}

The \ourdesign{} design and its analysis highlight an underlying tension when deciding which covariates to use in the AIPW estimator: as more covariates are added, there is an improvement in the oracle variance that \ourdesign{} approaches, but at the cost of a slower rate of convergence.
A similar trade-off for constructing AIPW estimators also occurs in the offline setting.
We explore this trade-off in greater detail thru numerical simulations in Section~\ref{sec:simulations}, where we use real data collected from a microfinance experiment.

Another interesting conclusion of this work is the distinction between Adaptive Neyman Allocation in design-based and super-population frameworks.
We show that $T^{-1/2}$ is the optimal rate of Neyman Regret in a design-based setting, whereas prior work has shown that $T^{-1} \log(T)$ is the optimal rate in a super-population setting \citep{Neopane2025Logarithmic,Neopane2025Optimistic}.
This difference mirrors results in the bandit literature, where adversarial and stochastic settings have minimax regret $T^{1/2}$ and $\log(T)$, respectively \citep{Lattimore2020Bandit}.
In both literatures, treating data as being deterministic offers more robustness, but at the cost of slower convergence.

\subsection{Related Work} \label{sec:related-work}

The foundations of Adaptive Neyman Allocation go back nearly a century.
\citet{Neyman1934Two} was the first to consider optimal allocation strategies, demonstrating that sampling from treatments proportional to the within-treatment standard deviation will minimize the variance of standard estimators.
To the best of our knowledge, \citet{Robbins1958Some} was the first to propose the sequential problem of constructing adaptive procedures which attain nearly the same variance as their optimal non-adaptive counterparts.

Adaptive experiments have seen a resurgence of interest from the causal inference community in the last twenty years.
We focus our attention on the potential outcomes framework \citep{Neyman1923, Rubin1980Randomization}.
From the causal inference perspective, early work focused on estimation and inference under adaptive treatment assignment \citep{vanderlann2008Construction, Hahn2011Adaptive}.
The focus in these papers is on estimating parameters (e.g. treatment probabilities) of the optimal design, but this does not directly guarantee that the resulting effect estimator obtains the correspondingly optimal variance.
More recent work has focused on adaptive experiment design for obtaining efficiency bounds \citep{vanderlann2014Online, Cook2024Semiparametric, Kato2024Active}.
These results are for the asymptotic variance and do not provide an analysis of how close the finite sample variance is to the efficiency bound.

The Neyman Regret is a non-asymptotic quantification of the gap between the actual variance and the optimal variance, which has only been investigated more recently.
To the best of our knowledge, the concept of Neyman Regret was proposed independently by \citet{Dai2023Clip} in the design-based setting and \citet{Kato2025Efficient} in the super-population setting.
In the super-population setting, \citet{Kato2025Efficient} develops adaptive AIPW estimators aimed at targeting an oracle obtained from minimizing the semi-parametric efficiency bound.
Using the assumption of consistent nuisance estimators, they establish that the Neyman Regret converges to zero, albeit at an unspecified rate.
In a series of work, \citet{Neopane2025Logarithmic,Neopane2025Optimistic} show that $T^{-1} \log(T)$ Neyman regret can be attained when focusing on Horvitz--Thompson and AIPW estimators in a super-population setting.
\citet{Li2024Optimal} have shown that vanishing Neyman regret is achievable in these settings even when the design parameters can only be adapted a few times, i.e. low-switching designs.

In the context of design-based inference, \citet{Dai2023Clip} use techniques from online optimization to construct an experimental design which attains $T^{-1/2} \exp( \sqrt{\log T} )$ Neyman Regret for the Horvitz--Thompson estimator.
\citet{Noarov2025Stronger} obtain $T^{-1} \log(T)$ Neyman Regret for the Horvitz--Thompson estimator, but under non-standard assumptions where the experimenter knows a constant lower bound on the absolute value of each individual potential outcome\footnote{In the context of AIPW estimators, this assumption would translate into a constant and known lower bound on each subject's absolute residual under the optimal regression.}.
\citet{Blackwell2022Batch} propose a two-stage approach to reduce the variance of the difference-in-means estimator, but the Neyman Regret is not explicitly analyzed.

This paper also contributes to the online optimization literature and particularly to the study of online methods for variance reduction.
There is a series of work, motivated by the training of machine learning models, which has studied regret minimization of importance sampling objectives $f_t(p) = \sum_{k=1}^K a_{t,k} / p_k$ over probability distributions $p \in \Delta_K$.
When $K=2$, the importance sampling objective reduces to the conditional variance objectives considered in this paper, albeit only when the predictors $\bv_t(1)$ and $\bv_t(0)$ are held fixed.
One idea in this literature has been to place a lower bound on the individual probabilities and use standard first-order methods to minimize regret \citep{Namkoong2017Adaptive,Zhao2025Adaptive}.
In an effort to avoid placing an artificial lower bound constraint, several works have sought alternative procedures for regret minimization.
\citet{Boros2018Online} propose an FTRL method with an inverse probability regularizer which obtains $\bigO{T^{-1/3}}$ regret and later extend this approach to a mixture setting \citep{Boros2019Online}.
\citet{ElHanchi2020Adaptive} propose a clipping technique for Stochastic Gradient Descent and show that it also achieves $\bigO{T^{-1/3}}$ regret.
The sigmoidal regularizer we present in this paper allows for an improved $\bigO{T^{-1/2}}$ regret bound in these settings; however, we wish to emphasize that the objective in this paper has an additional predictor $\bv_t(k)$ component which introduces additional complexities such as non-convexity.
Additionally, online techniques have been used in constructing more efficient Monte Carlo integrations \citep{Carpentier2011Finite, Carpentier2012Adaptive, Carpentier2012Minimax}.

Outside of Adaptive Neyman Allocation, other aspects of adaptive experiments have been studied from a number of perspectives.
A recent line of work has provided methods for estimation and inference of causal effects when treatment is assigned via a bandit algorithm \citep{Hadad2021Confidence, Zhang2020Inference, Zhang2021Statistical}.
\citet{OfferWestort2021Adaptive,Chen2023Optimal} focus on selective inference under adaptive treatment assignment.
Finally, anytime inference and data-dependent stopping are central areas of study in sequential analysis \citep{Wald1945Sequential, Howard2021Time, WaudbySmith2024Time}.

	\section{Design-based Adaptive Experiments} \label{sec:preliminaries}

\subsection{Design-Based Potential Outcomes Framework}

We consider a sequential experiment with $T$ experimental subjects denoted by integers $t \in [T]$.
The experimenter must assign each unit to exactly one of two treatment conditions.
For each subject $t \in [T]$, we denote their treatment assignment as $Z_t \in \setb{0,1}$.
Each experimental subject $t \in [T]$ is presumed to have two potential outcomes, $y_t(1)$ and $y_t(0)$, which correspond to the outcomes that would be measured under treatment ($Z_t = 1$) and control ($Z_t = 0$), respectively.
The causal estimand of interest is the \emph{average treatment effect} (ATE) which is defined as
\[
\ate = \frac{1}{T} \sum_{t=1}^T \paren[\big]{ y_t(1) - y_t(0) } \enspace.
\]
The experimenter also measures a vector of covariates $\xv_t \in \Reals^d$ for each subject $t \in [T]$, which are not affected by treatment assignment and may be used to improve estimates of the average treatment effect.
Implicit in the above is the standard Stable Unit Treatment Value Assumption (SUTVA) which posits that there are no hidden versions of treatment and that subjects do not interfere with each other \citep{Holland1986Statistics, ImbensRubin2015, Hernan2020What}.

In this paper, we work in a design-based framework where the subjects, their potential outcomes, and their covariates are all considered to be deterministic and treatment is the sole source of randomness.
In such a setting, randomization of treatment serves as the sole basis for statistical inference, e.g. no i.i.d. assumptions are placed on the subjects.

\subsection{Adaptive Experiment Designs}

The sequential experimental procedure proceeds in $T$ rounds.
At each round $t \in [T]$, the experimenter observes the covariates $\xv_t \in \Reals^d$, then assigns treatment $Z_t \in \setb{0,1}$ and consequently observes the outcome
\[
Y_t = \indicator{Z_t = 1} y_t(1) + \indicator{Z_t = 0} y_t(0) \enspace.
\]
The experiment may be adaptive in the sense that the randomization of treatment assignment $Z_t$ and the choice of linear predictors $\bv_t(1)$ and $\bv_t(0)$ (to be defined in the next section) can depend on previously observed outcomes $Y_1, \dots, Y_{t-1}$, treatment assignments $Z_1, \dots, Z_{t-1}$, and covariates $\xv_1 , \dots, \xv_{t}$.
Formally, an \emph{adaptive experiment design} is a sequence of mappings $\fullarbdesign = \setb{ \tarbdesign }_{t=1}^T$ with signature $\tarbdesign : (\setb{0,1}\times \Reals^d \times \Reals)^{t-1} \times \Reals^d \to [0,1] \times \Reals^d \times \Reals^d$ that encode the conditional treatment assignment probability and linear predictors, i.e.
\[
\Pr{ Z_t = 1 \mid \filt_{t-1} }, \bv_t(1), \bv_t(0) = \tarbdesign( Z_1, \xv_1, Y_1, \dots , Z_{t-1}, \xv_{t-1}, Y_{t-1}, \xv_t)
\enspace,
\]
where $\filt_{t-1}$ conditions on past observations, i.e. formally, $\filt_{t-1}$ is the $\sigma$-algebra generated by $Z_1, \dots, Z_{t-1}$.
Recall that in the design-based context, the potential outcomes are deterministic so that all randomness (including observed outcomes $Y_1 ,\dots, Y_{t-1}$) are functions of treatment assignments $Z_1, \dots, Z_{t-1}$.
The linear predictors $\bv_t(1)$ and $\bv_t(0)$ are used in the estimator, as described in Section~\ref{sec:adaptive-AIPW-estimators}.
We focus on the setting where the order of subjects in the sequence is fixed and arbitrary (i.e. cannot be chosen by the experimenter), which reflects the conditions often arising in practice.
For this reason, we refer to the \emph{sequence} of potential outcomes $\setb{ (y_t(1), y_t(0) ) }_{t=1}^T$.
We will not assume that the sequence of potential outcomes satisfies any type of stationary condition, i.e. the outcomes are allowed to arrive in an arbitrary order.

\subsection{Technical Assumptions} \label{sec:technical-assumptions}

Our main technical assumptions on the potential outcomes and covariates are given below.
We emphasize that all constants in the assumptions are presumed to exist but are not known to the experimenter.
Indeed, the fact that these constants are not known before the experiment is part of what makes the problem more challenging.

\begin{assumption}[Bounded Moments] \label{assumption:moments}
	There exist constants $0 < c_0 \leq c_1$ such that for both treatments $k \in \setb{0,1}$,
	\[
	c_0 
	\leq
	\min_{\bv \in \Reals^d} \paren[\Big]{ \frac{1}{T} \sum_{t=1}^T \braces[\Big]{y_t(k) - \iprod{ \xv_t , \bv } }^2 }^{1/2}
	\leq 
	\paren[\Big]{ \frac{1}{T} \sum_{t=1}^T y_t(k)^4 }^{1/4}
	\leq 
	c_1
	\enspace.
	\]
\end{assumption}

The moment conditions in Assumption~\ref{assumption:moments} ensure two things.
First, the second moments of the OLS residuals are assumed to be bounded from below.
This assumption is generally plausible, unless the outcomes are suspected to be exactly a linear function of the covariates, which is rarely (if ever) the case in practice.
The assumption also places a bound on the fourth moments of the potential outcomes.
Because the problem of Adaptive Neyman Allocation essentially requires the estimation of squared residuals, it is unlikely that either of these moment assumptions can be weakened.

\begin{assumption}[Covariate Regularity] \label{assumption:covariate-regularity}
	There exist constants $c_2 > 0$ and $\gamma_0 \in (0, \sqrt{T/d})$ such that for all $t \geq \sqrt{T d} \cdot \gamma_0$, the covariate matrix is well-invertible:
	\[
	\frac{1}{c_2} \leq \sigma_{\min}\paren[\Big]{ \frac{1}{t} \sum_{s \leq t} \xv_s \xv_s^\tran }
	\enspace.
	\]
\end{assumption}

Assumption~\ref{assumption:covariate-regularity} ensures that after the ``early iterations'' (i.e. iterations $t = \bigO{\sqrt{Td} }$), the empirical covariance matrix $\frac{1}{t} \sum_{s \leq t} \xv_s \xv_s^\tran$ is well-invertible.
This invertibility condition guarantees that adaptively estimated regression coefficients are well-behaved.
The condition that $\gamma_0 \leq \sqrt{T/d}$ ensures that the covariate regularity in Assumption~\ref{assumption:covariate-regularity} holds for the entire set of covariates, i.e. $t=T$. 
We will not require that the covariate matrices in the sequence are well-conditioned, in the sense that the largest singular value may at times be much larger than its smallest one.
Moreover, our asymptotic analyses will not assume that this covariate matrix converges to any limiting quantity.

The next assumption bounds the maximum radius, defined as $R = \max_{t \in [T]} \norm{\xv_t}$.
We do not presume that the maximum radius $R$ is known a priori to the experimenter.

\begin{assumption}[Maximum Radius]
	\label{assumption:maximum-radius}
	There exists a constant $c_3 > 0$ such that $R \leq c_3 T^{1/2}$.
\end{assumption}

The most salient aspect of Assumptions~\ref{assumption:moments}-\ref{assumption:maximum-radius} is that they allow for non-stationarity in the sequence of potential outcomes and covariates.
Nearly all quantities are allowed to drift arbitrarily throughout the experiment, including individual treatment effects and the residuals of best linear predictors.
The only substantive restriction on the order of the subjects is the well-invertible condition of Assumption~\ref{assumption:covariate-regularity}, which is fairly mild.

The covariate dimension $d$ and the maximum radius $R$ are implicitly related through these assumptions in various ways.
First, we should note that while Assumptions~\ref{assumption:moments} and \ref{assumption:covariate-regularity} require that $d < T$, they permit $d$ to grow linearly with $T$.
Next, we note that these assumptions imply a non-degeneracy condition where the maximum radius $R$ is at least order $\sqrt{d}$.
In particular, Assumption~\ref{assumption:covariate-regularity} implies that 
\begin{equation} \label{eq:R_lowerbound}
c_{2}^{-1} 
\leq \lambda_{\min}\paren[\Big]{ \frac{1}{T} \sum_{t \leq T} \xv_t \xv_t^\tran }
\leq \frac{1}{d} \tr[\Big]{ \frac{1}{T} \sum_{t \leq T} \xv_t \xv_t^\tran }
= \frac{1}{d T} \sum_{t \leq T} \norm{\xv_t}^2 
\leq \frac{R^2}{d}
\enspace,
\end{equation}
which yields that $R \geq c_2^{-1/2} \sqrt{d}$.
This does rule out settings with sparse, high-dimensional covariate vectors; however, OLS-based AIPW estimators considered in this paper would be an inappropriate choice in such settings.

Under the natural regime where each entry of every covariate vector is of constant order (i.e. $\max_{t \in [T], k \in [d]} \abs{ \xv_t(k) } = \bigO{1} $ ), our assumptions imply the relation $R^2 \asymp d$.
This is arguably the most relevant regime to consider because covariates are often deemed to be bounded in this way in many practical situations, e.g. age, blood pressure, etc.
For this reason, the main results stated in the abstract and introduction have been presented specifically in this regime.

However, the analysis in the rest of our paper requires only Assumptions~\ref{assumption:moments}-\ref{assumption:maximum-radius}, which allow for greater flexibility.
In particular, our assumptions allow for the existence of outlier covariates, which are forbidden in the natural regime considered above.
For example, Assumption~\ref{assumption:maximum-radius} allows that $R^2 \gg d$ in low-dimensional settings where $d = \littleO{T}$.
In a similar spirit, Assumption~\ref{assumption:covariate-regularity} allows that the empirical covariance matrix might have a large eigenvalue due to an outlier covariate.
The most attractive aspect of these assumptions is that neither the magnitude of the outlier nor its position in the sequence is required to be known by the experimenter in advance.

To reason about the practicality of the assumptions above, it may help to view them through the lens of a super-population.
If the outcomes and covariates were sampled i.i.d., then Assumptions~\ref{assumption:moments} and \ref{assumption:covariate-regularity} would hold with probability tending to 1 under the conditions that (i) the fourth moments of the outcomes existed, (ii) the covariates were sampled according to a subgaussian distribution with $d = \littleO{ T^{ 1/2 } }$, and (iii) the conditional variance of the outcomes is positive almost surely.

\paragraph{Triangular Array Asymptotics}
Although the majority of our results are finite-sample in nature, we will also consider asymptotic analyses.
We follow the convention in the design-based literature of using triangular array asymptotics.
In the triangular array asymptotics, we consider a sequence of experiments indexed by $T \in \Naturals$.
For each $T$, there is a sequence of potential outcomes and covariates $\setb{ y_t^{(T)}(1) , y_t^{(T)}(0), \xv_t^{(T)} }_{t=1}^T$ and an adaptive experimental design $\mathcal{D}^{(T)}$.
This yields a sequence of (deterministic) average treatment effects $\ate^{(T)}$ and estimators $\eate^{(T)}$.
All limiting statements, e.g. $\ate^{(T)} - \eate^{(T)} \xrightarrow{p} 0$, are made with respect to this sequence.
For notational clarity, we drop the superscript $T$ in the asymptotic statements.

%

\subsection{Adaptive AIPW Estimators} \label{sec:adaptive-AIPW-estimators}

We focus on the class of adaptive Augmented Inverse Propensity Weighted (AIPW) estimators for the average treatment effect.
The AIPW estimator is widespread in causal inference, where a regression model is used to improve the efficiency of the standard IPW estimator.
In this work, we focus on adaptive AIPW estimators with linear regression models.

At each round, the experimenter updates a linear regression model of the outcomes under treatment and control using the observed history.
Let $\bv_t(1)$, $\bv_t(0) \in \Reals^d$ denote the coefficients of the linear regression models for the outcomes under treatment and control at time $t$, respectively, i.e. $\bv_t(1)$ and $\bv_t(0)$ are determined completely by the history $\filt_{t-1}$.
For notational convenience, we use $p_t$ to denote the conditional probability of treatment assignment at round $t$ given history, i.e. $p_t = \Pr{ Z_t = 1 \mid \filt_{t-1} }$.
At each round $t \in [T]$, the adaptive AIPW estimator proceeds by estimating the individual treatment effect using the regression models, then correcting via an IPW estimate of the residuals.
Formally, the adaptive AIPW estimator is given as
\begin{align*}
\eate &= \frac{1}{T} \sum_{t=1}^T \braces[\Big]{ 
	\iprod{\xv_t , \bv_t(1)} - \iprod{ \xv_t , \bv_t(0) } } \\
	&\quad 
	+ \frac{1}{T} \sum_{t=1}^T \braces[\Big]{
	\frac{\indicator{Z_t = 1}}{p_t} \paren[\Big]{ Y_t - \iprod{\xv_t , \bv_t(1)} }
	- \frac{\indicator{Z_t = 0}}{1-p_t} \paren[\Big]{ Y_t - \iprod{\xv_t , \bv_t(0)} }
}
\enspace.
\end{align*}

\expcommand{\aipwbias}{%
	If $p_t \in (0,1)$ for all $t \in [T]$ a.s. then the adaptive AIPW estimator is unbiased: $\E{ \eate } = \ate$.
}
\begin{proposition}[AIPW Bias] \label{prop:aipw-unbiased}
	\aipwbias
\end{proposition}

Proposition~\ref{prop:aipw-unbiased} shows that the adaptive AIPW estimator is unbiased, regardless of how well the online regression model fits the data.
The next proposition shows that better fitted models typically correspond to smaller variance, so long as the conditional treatment probabilities are not too extreme.

\expcommand{\aipwvar}{
The normalized variance $T \cdot \Var{\eate}$ of the AIPW estimator is given as
\begin{align*}
	 \E[\Bigg]{ \frac{1}{T} \sum_{t=1}^T \paren[\Bigg]{ 
			\braces[\big]{ y_t(1) - \iprod{ \xv_t , \bv_t(1) } } \cdot \sqrt{ \frac{1-p_t}{p_t} }
			+  	\braces[\big]{ y_t(0) - \iprod{ \xv_t , \bv_t(0) } } \cdot \sqrt{ \frac{p_t}{1-p_t} }
		}^2 }
	\enspace.
\end{align*}
}
\begin{proposition}[AIPW Variance] \label{prop:aipw-variance}
	\aipwvar
\end{proposition}
	
	\section{Adaptive Neyman Allocation} \label{sec:adaptive-neyman-allocation}

\subsection{Formulation of the Neyman Regret}

Introduced by \citet{Robbins1958Some}, the problem of Adaptive Neyman Allocation is to design an adaptive protocol which has nearly the same performance as the optimal non-adaptive protocol which has access to all of the data.
In the recent literature on Adaptive Neyman Allocation, the performance of an experimental design is measured by the Neyman Regret, which is the gap between the adaptive variance and the non-adaptive oracle variance.

We begin by deriving the oracle variance which serves as the relevant comparator when constructing an adaptive design.
In this context, the \emph{oracle variance} $\optvar$ is defined as the minimal variance of the AIPW estimator when selecting the best fixed linear predictors and assignment probability:
\[
\optvar = \min_{ p, \bv(1), \bv(0) } \Var[\big]{ \eate ; p, \bv(1), \bv(0) }
\enspace.
\]
The \emph{Neyman Allocation} refers to the optimal choice of linear predictors $\optbv(1)$, $\optbv(0)$ and assignment probability $p^*$ which attain this oracle variance.
Note that the oracle variance $\optvar$ and the Neyman Allocation $\optbv(1)$, $\optbv(0)$, and $p^*$ depend on the entire sequence of potential outcomes and covariates.
For each treatment $k \in \setb{0,1}$, define the \emph{optimal residuals} as
\[
\olsres{k} = \min_{\bv \in \Reals^d} \paren[\Bigg]{ \frac{1}{T} \sum_{t=1}^T \paren{ y_t(k) - \iprod{ \xv_t, \bv } }^2 }^{1/2}
\enspace,
\]
and denote the corresponding OLS predictor as $\bv_{\textrm{OLS}}(k)$.
Define the \emph{residual correlation} $\rho$ as
\[
\rho = \frac{ \frac{1}{T} \sum_{t=1}^T \paren[\big]{ y_t(1) - \iprod{ \xv_t , \bv_{\textrm{OLS}}(1) } } \paren[\big]{ y_t(0) - \iprod{ \xv_t , \bv_{\textrm{OLS}}(0) } } }{ \olsres{0} \olsres{1} }
\enspace.
\]
The following proposition derives the oracle variance and the Neyman Allocation in terms of these parameters.

\expcommand{\oracle}{
	The oracle variance is given by $T \cdot \optvar = 2 (1+\rho) \olsres{1} \olsres{0}$
	and the Neyman allocation is given by the least squares predictors $\optbv(k) = \bv_{\textrm{OLS}}(k)$ and assignment probability $p^*=(1 + \olsres{0} / \olsres{1})^{-1}$.
}

\begin{proposition} \label{prop:oracle}
	\oracle
\end{proposition}

At first glance, it may appear that the OLS predictors cannot be the minimizers of the variance due to the presence of the cross term, i.e. the correlation between the residuals.
However, this intuition is false{\textemdash}it turns out that the minimizers are indeed given by the OLS predictors.
Some care is needed to establish this fact because the function $(p, \bv(1), \bv(0)) \mapsto \Var{\eate ; p, \bv(1), \bv(0)}$ is non-convex.
In Section~\suppref{section:C1} of the appendix, we give the full proof which uses orthogonality properties of the OLS predictors as well as several judicious applications of AM-GM and Cauchy-Schwarz.

\begin{definition}
	Given an adaptive design, the \emph{Neyman Regret} $\neymanregret_T$ is defined as the difference between the (normalized) adaptive and oracle variances:
	\[
	\neymanregret_T = T \cdot \Var{\eate} - T \cdot \optvar \enspace.
	\]
\end{definition}

We use a subscript $T$ in the Neyman Regret $\neymanregret_T$ to reflect the dependence on the sample size $T$.
If the Neyman Regret is decreasing to $0$ with the sample size, then the true adaptive variance is essentially upper bounded by the oracle variance.
Roughly speaking, this means that the experimenter will be able to use estimates of the oracle variance $T \cdot \optvar$ in their inferential procedures.

Our goal will be to construct an adaptive experimental design which approaches oracle Neyman Allocation performance in the sense that $\neymanregret_T \to 0$.
But why might the oracle described above{\textemdash}in particular, one that considers fixed treatment probability $p^*$ and predictors $\optbv(1)$ and $\optbv(0)${\textemdash}be an appropriate target?
The strongest argument in favor of targeting this oracle is practical: by doing so, we obtain an adaptive design which provides a substantial improvement in variance over existing techniques within the design-based framework.
While considerably stronger oracles will provide smaller variance, they will not serve as practical targets for adaptive designs.
For example, consider the strongest oracle which is allowed to select $p^*_t$, $\optbv_t(1)$, and $\optbv_t(0)$ uniquely for each subject in an unrestricted manner.
This oracle can achieve 0 variance by exactly predicting both potential outcomes.
In fact, the weaker oracle which fixes $\optbv(1) = \optbv(0) = 0$ and selects only individual $p_t^*$ uniquely can still achieve 0 variance when $y_t(1)$ and $y_t(0)$ have different signs.
Clearly, no adaptive design can meaningfully compete with these stronger oracles.
A complete and thorough investigation into which stronger oracles can serve as meaningful targets lies outside the scope of this paper.

Before we begin constructing an adaptive experimental design under which $\neymanregret_T \to 0$, we should first ask ourselves how quickly we can expect this convergence to happen.
The following theorem shows that under our assumptions, the Neyman regret cannot converge faster than $T^{-1/2}$.

\expcommand{\lowerbound}{
	There exists a constant $c > 0$ such that
	for all integers $T$ and any adaptive experimental design $\fullarbdesign$,
	there exists a sequence $\{ y_t(1), y_t(0), \xv_t \}_{t=1}^T$ satisfying Assumptions~\mainref{assumption:moments}-\mainref{assumption:maximum-radius} under which the corresponding Neyman regret is at least $\neymanregret_T \geq c \cdot T^{-1/2} $.
}

\begin{theorem}[Lower Bound] \label{theorem:neyman-regret-lower-bound}
	\lowerbound
\end{theorem}

Theorem~\ref{theorem:neyman-regret-lower-bound} demonstrates a $T^{-1/2}$ lower bound under the conditions of bounded moments (Assumption~\ref{assumption:moments}), covariate regularity (Assumption~\ref{assumption:covariate-regularity}), and bounded maximum radius (Assumption~\ref{assumption:maximum-radius}).
The constant $c > 0$ in Theorem~\ref{theorem:neyman-regret-lower-bound} depends on the constants $c_0$, $c_1$, $c_2$, $c_3$, and $\gamma_0$ appearing in these assumptions.
In Section~\suppref{sec:supp-lower-bound} of the appendix, we show that weakening any of these assumptions even moderately will imply that the Neyman regret is at least constant (i.e. not converging to zero) for any adaptive experimental design.

\subsection{Online Optimization and Technical Challenges}

Following \citet{Dai2023Clip}, we view the problem of minimizing the Neyman Regret through the lens of online optimization.
Define the objective functions $g_t : (0,1) \times \Reals^d \times \Reals^d \to \Reals_+$ as
\[
g_t\paren[\big]{ p, \bv(1), \bv(0) }
= \paren[\Bigg]{ 
	\braces[\big]{ y_t(1) - \iprod{\xv_t, \bv(1)} } \cdot \sqrt{ \frac{1-p}{p} }
	+ \braces[\big]{ y_t(0) - \iprod{\xv_t, \bv(0)} } \cdot \sqrt{ \frac{p}{1-p} }
  }^2
\enspace.
\]
Then, using Propositions~\ref{prop:aipw-variance} and \ref{prop:oracle}, we can write the Neyman regret in terms of online optimization, i.e. as the difference between the expected objective under the adaptive design and the objective for the optimal set of parameters:
\[
\neymanregret_T 
= \E[\Big]{ \frac{1}{T} \sum_{t=1}^T g_t\paren[\big]{ p_t, \bv_t(1), \bv_t(0) } }
- \min_{p, \bv(1), \bv(0)} \frac{1}{T} \sum_{t=1}^T g_t \paren[\big]{ p, \bv(1), \bv(0) } 
\enspace.
\]
The underlying online optimization problem will guide our construction of an adaptive experimental design.
However, the online optimization itself presents two challenges which require new technical developments:

\begin{enumerate}
	\item \textbf{Nonconvexity}: The objective functions $g_t$ are non-convex.
	Broadly speaking, convexity is considered to be the watershed for whether a problem can be efficiently solved.
	In order to guarantee convergence of the Neyman Regret, we need to overcome the issues presented by non-convexity in the objective.
	
	\item \textbf{Ill-conditioned}: The gradients of $g_t$ become arbitrarily large as $p$ approaches the boundary of $(0,1)$.
	Analysis of conventional optimization methods (e.g. gradient descent) often requires that the gradients are bounded over the entire domain.
	The tension is that the optimal $p^*$ may be close to the boundary $(0,1)$ and so we should allow the adaptive $p_t$ to get close to the boundary, albeit in a controlled manner.
	In order to get clean $T^{-1/2}$ rates of Neyman Regret, we must also address the fundamental ill-conditioning in the objective functions.
\end{enumerate}

\subsection{Decomposition of the Neyman Regret}

In this section, we show how to overcome the non-convexity underlying the Neyman Regret.
We show that the non-convex Neyman Regret can be decomposed into the sum of two convex regrets: probability regret and prediction regret.
Roughly speaking, the probability regret measures the extent to which the adaptively chosen probability balances the online residuals while the prediction regret measures the performance of the adaptively chosen linear predictors.

We begin by showing how to evaluate the performance of the adaptively chosen treatment probability.
For each iteration $t$, we define the \emph{probability loss} $\probloss{t} : (0,1) \to \Reals_+$ as the convex function
\[
\probloss{t}(p) = 
\paren[\Big]{ y_t(1) - \iprod{ \xv_t , \bv_t(1) } }^2 \cdot \frac{1}{p}
+ \paren[\Big]{ y_t(0) - \iprod{ \xv_t , \bv_t(0) } }^2 \cdot \frac{1}{1-p}
\enspace.
\]
The probability loss measures how well a treatment assignment probability $p$ balances the squared error of the predictions made at iteration $t$.
We define the \emph{probability regret} as
\[
\regretprob_T = \sum_{t=1}^T \probloss{t}(p_t) - \sum_{t=1}^T \probloss{t}(p^*)
\enspace.
\]
The probability regret measures how well the adaptively chosen assignment probabilities $p_1, \dots, p_T$ balance the online residuals via comparison to the Neyman Allocation assignment probability $p^*$.
Because the adaptively chosen probabilities $p_t$ and linear predictors $\bv_t(1)$ and $\bv_t(0)$ are random, the probability regret $\regretprob_T$ is also a random variable.
Note that the scaling of the probability regret is such that $\littleO{T}$ is considered to be asymptotically vanishing.

We now focus on how to evaluate the performance of the adaptively chosen linear predictors.
For each iteration $t$, we define the \emph{prediction loss} function $\ell_t : \Reals^d \times \Reals^d \to \Reals_+$ as the convex function
\[
\ell_t \paren[\big]{ \bv(1), \bv(0) }
= \paren[\Bigg]{
	\braces[\big]{ y_t(1) - \iprod{ \xv_t , \bv(1) } } \cdot \sqrt{ \frac{\olsres{0}}{\olsres{1}} }
	+
	\braces[\big]{ y_t(0) - \iprod{ \xv_t , \bv(0) } } \cdot \sqrt{ \frac{\olsres{1}}{\olsres{0}} }
}^2
\enspace.
\]
The prediction loss measures the error of linear predictions of the potential outcomes at a given iteration, weighted by the ratio of the optimal residuals.
The ratio of optimal residuals can also be understood in terms of the Neyman Allocation probability $p^*$ because $\olsres{0} / \olsres{1} = 1/p^* - 1$ and $\olsres{1} / \olsres{0} = 1/(1-p^*) - 1$.
From this perspective, the prediction loss measures the error of the linear predictors $\bv_t(1), \bv_t(0)$, weighted by the Neyman Allocation $p^*$.
We define the \emph{prediction regret} as
\[
\newpredregret_T 
= \sum_{t=1}^T \newpredloss{t} \paren[\big]{ \bv_t(1), \bv_t(0) } 
- \sum_{t=1}^T \newpredloss{t} \paren[\big]{ \optbv(1),  \optbv(0) }
\enspace.
\]
The prediction regret measures the overall prediction error of the adaptively chosen predictors by comparing to the prediction error under the optimal least squares predictors $\bv^*(1), \bv^*(0)$ -- again, both errors are weighted accordingly by $p^*$.
Because the linear predictors are random, the prediction regret $\newpredregret_T $ is also a random variable.
The prediction regret is on the same scale as the probability regret, i.e. a sublinear $\littleO{T}$ regret is considered to be asymptotically vanishing.

\expcommand{\regretdecomposition}{%
	The Neyman Regret can be decomposed as the $T$-normalized sum of the probability regret and the prediction regret:
	\[
	\neymanregret_T = \frac{1}{T} \E{ \regretprob_T } + \frac{1}{T} \E{ \newpredregret_T }
	\enspace.
	\]
}
\begin{lemma}\label{lemma:regret-decomposition}
	\regretdecomposition
\end{lemma}

Lemma~\ref{lemma:regret-decomposition} demonstrates that the Neyman Regret can be decomposed into the sum of the expected probability and prediction regrets.
The proof, which appears in Section~\suppref{section:C1} of the supplement, is simple: roughly speaking, simply add and subtract the AIPW variance that would be obtained when selecting the \emph{actual} choice of linear predictors $\bv_t(1)$ and $\bv_t(0)$ in the experiment together with the \emph{hypothetical} oracle choice of assignment probability $p^*$.
What appears is the sum of the probability and prediction regrets.
This explains the different roles that the two regrets have: the probability regret measures the balancing of actual residuals by selection of $p_t$ whereas the prediction regret measures the accuracy of the predictors weighted by the optimal allocation $p^*$.

The decomposition of the Neyman regret in Lemma~\ref{lemma:regret-decomposition} effectively replaces the issue of non-convexity with the issue of multiple objectives; that is, we still need the expected probability regret and the expected prediction regret to be small \emph{simultaneously}.
The main benefit of Lemma~\ref{lemma:regret-decomposition} is that it identifies individual convex objectives on which techniques from convex analysis may be used.
	
	\section{The Sigmoid-FTRL Design} \label{sec:our-design}

\subsection{Formal Description of the Sigmoid-FTRL Design} \label{sec:presenting-our-design}

\begin{algorithm}[t] \caption{\ourdesign}\label{alg:sigmoid-ftrl}
	\KwIn{Sigmoid $\phi : \Reals \to (0,1)$ satisfying Condition~\mainref{condition:sigmoid}. }
	Define the sigmoid regularizer $\Psi  = \psi \circ \phi^{-1}$ with $\psi(u) = \frac{1}{2}u^2 + |u|^3$. \\
	Initialize max radius $R_0 = \sqrt{d}$. \\
	\For{ $t=1, 2, \dots T$}{
		Observe covariate $\xv_t \in \Reals^d$ and update $R_t = \maxf{ R_{t-1}, \norm{\xv_t} }$ and $\eta_t = T^{-1/2} R_t^{-1}$. \\
		Construct the regression coefficients $\bv_t(1)$ and $\bv_t(0)$ as
		\begin{align*}
			\bv_t(1) &= \argmin_{ \bv \in \Reals^d } \sum_{s=1}^{t-1} \paren[\Big]{ Y_s \cdot \frac{\indicator{Z_s=1}}{p_s} - \iprod{ \xv_s , \bv } }^2 
			+ \eta_t^{-1}  \norm{\bv}^2 \\
			\bv_t(0) &= \argmin_{ \bv \in \Reals^d } \sum_{s=1}^{t-1} \paren[\Big]{ Y_s \cdot \frac{\indicator{Z_s=0}}{1-p_s} - \iprod{ \xv_s , \bv } }^2 + \eta_t^{-1}   \norm{\bv}^2
		\end{align*} \\
		Construct estimates of the online squared residuals $\estores{t-1}{1}$ and $\estores{t-1}{0}$ as
		\begin{align*}
			\estores{t-1}{1} &= \sum_{s=1}^{t-1} \frac{\indicator{Z_s = 1}}{p_s} \paren[\Big]{ Y_s - \iprod{ \xv_s , \bv_s(1) } }^2 \\
			\estores{t-1}{0} &= \sum_{s=1}^{t-1} \frac{\indicator{Z_s = 0}}{1-p_s} \paren[\Big]{ Y_s - \iprod{ \xv_s , \bv_s(0) } }^2
		\end{align*} \\
		Select assignment probability $p_t$ as
		\[
		p_t = \argmin_{p \in (0,1)} 
		\frac{ \estores{t-1}{1} }{p} 
		+ \frac{ \estores{t-1}{0} }{1-p}
		+ \eta_t^{-1} \Psi(p)
		\enspace.
		\]\\
		Sample treatment assignment $Z_t = 1$ with probability $p_t$, and $Z_t=0$ otherwise. \\
		Observe outcome $Y_t = \indicator{Z_t=1} \cdot y_t(1) + \indicator{Z_t = 0} \cdot y_t(0)$.
	}
\end{algorithm}

The \ourdesign{} design is formally described in Algorithm~\ref{alg:sigmoid-ftrl}.
The design can be understood as applying the \emph{follow-the-regularized leader} (FTRL) principle separately to the probability regret and the prediction regret.
FTRL is a commonly used technique in online convex optimization which is an alternative to gradient-descent-based methods \citep{Hazan2016Introduction}.
As each subject arrives in the experiment, the \ourdesign{} design proceeds in two main steps: first, the linear predictors $\bv_t(1)$, $\bv_t(0)$ are computed and then the treatment assignment probability $p_t$ is computed.

The linear predictors  $\bv_t(1)$ and $\bv_t(0)$ are chosen to minimize the estimated squared prediction errors on the previously observed units.
The squared prediction errors on previously observed units are estimated using adaptive IPW weighting, which is necessary because we observe only one outcome for each unit.
To ensure that predictors are sufficiently regularized, we add a ridge regularizer with penalty term $\eta_t^{-1}$, which depends on the number of units $T$ and the largest norm of any covariate vector seen so far.
As we demonstrate in Section~\ref{sec:prediction-regret}, this step aims to minimize the prediction regret.

The treatment probability $p_t$ is chosen so that the subject is more likely to be assigned to the treatment where the online predictions have so far resulted in larger errors.
For each treatment $k\in \setb{0,1}$, define the \emph{online residuals up to time $t-1$} as $\ores{t-1}{k} = \sum_{s \leq t-1} (y_s(k) - \iprod{ \xv_s , \bv_s(k) } )^2$.
The online residuals are not directly observable because in all previous iterations $s < t$, we have seen either $y_s(1)$ or $y_s(0)$, but not both.
For this reason, we use adaptive IPW weighting to obtain estimates $\estores{t-1}{k}$.
The treatment probability $p_t$ is then chosen to minimize the weighted sum $\estores{t-1}{1} / p + \estores{t-1}{0} /(1-p)$ with an additional regularization term $\eta_t^{-1} \cdot \Psi(p)$.
As we demonstrate in Section~\ref{sec:probability-regret}, this step aims to minimize the probability regret.

We use a \emph{sigmoid regularizer} $\Psi : (0,1) \to \Reals_+$ to regularize the selected treatment probability $p_t$.
The purpose of the sigmoidal regularizer is to ensure that the treatment probability does not get too close to the boundary of the interval $[0,1]$, which would increase the variance of the AIPW estimator.
We refer to this penalty function as sigmoidal because it takes the form $\Psi  = \psi \circ \phi^{-1}$, where $\phi : \Reals \to (0,1)$ is a sigmoid function and $\psi(u) = \frac{1}{2} u^2 + |u|^3$.
In other words, the penalty $\Psi$ can be understood as using a quadratic + cubic penalty on the transformed variable $u_t$ which yields the treatment probability, i.e. $p_t = \phi(u_t)$.

The sigmoidal penalty is one of the novel techniques introduced in our work and is the namesake of the \ourdesign{} design.
The sigmoid penalty is a key ingredient in achieving the clean $T^{-1/2}$ term in our regret analysis.
This improves upon the probability clipping design of \citet{Dai2023Clip}, which features an additional sub-polynomial $\expf{ \sqrt{\log(T)} }$ factor as a result of the probability clipping.
One can interpret the sigmoidal regularizer as a gentler alternative to the harsh regularization implicitly introduced by probability clipping.

The choice of an appropriate sigmoid function $\phi$ is crucial for obtaining the performance guarantees of the \ourdesign{} design.
Many choices of sigmoid will not work, including the logistic function, the hyperbolic tangent, and the error function.
Fortunately, there are several good choices, such as the arctangent function $\phi(u) = \frac{1}{\pi} \paren{ \arctan(u) + \pi / 2}$ and the algebraic sigmoid $\phi(u) = \frac{1}{2} \paren{ \frac{u}{1+\abs{u}} + 1 }$.
Our theoretical results will hold for any sigmoid satisfying the following condition:

\expcommand{\sigmoidcondition}{%
	The sigmoid $\phi : \Reals \to (0,1)$ satisfies the following:
	\begin{enumerate}
		\item $\phi(u)$ is strictly monotone increasing with $\phi(u)+\phi(-u)=1$ and $\phi(+\infty)=1$.
		\item $u \mapsto 1/\phi(u)$ and $u \mapsto 1/(1-\phi(u))$ are convex functions.
		\item There exist constants $b_1, b_2, b_3>0$ such that
		\begin{enumerate}
			\item $-\left(\frac{1}{\phi(u)}\right)^{\prime} \leq b_1$ for all $u \in \Reals$.
			\item $ \left(\frac{1}{\phi(u)}\right)^{\prime\prime}\leq b_2\cdot \frac{1}{(1+|u|)^3}$ for all $u \in \Reals$ except for a finite number of points.
			\item $ \left(\frac{1}{\phi(u)}\right)^{\prime\prime}\geq b_3\cdot \frac{1}{(1+u)^3}$ for all $u \geq 0$ except for a finite number of points.
		\end{enumerate}
	\end{enumerate}
}

\begin{condition}[Sigmoid Condition] \label{condition:sigmoid}
	\sigmoidcondition
\end{condition}

Another key aspect of the \ourdesign{} design is the introduction of an adaptive step size $\eta_t$.
In the context of FTRL, the step size may be equivalently viewed as a penalty term.
The penalty term is set adaptively as $\eta_t = (T^{1/2} R_t)^{-1}$ where $R_t$ is the largest covariate norm seen so far.
The variable is initialized as $R_0 = \sqrt{d}$ and updated iteratively as $R_t = \maxf{ R_{t-1}, \norm{\xv_t} }$.
From a theoretical point of view, the adaptively chosen penalty term $\eta_t$ ensures that the regularization appropriately scales with the magnitude of the covariates.
From a practical point of view, the primary benefit of an adaptive penalty term $\eta_t$ is that the experimenter is not required to correctly specify the magnitude of outlier covariates a priori.
In keeping with the conventions of online optimization, we define $\eta_{T+1} = \eta_T$.

\ourdesign{} is computationally practical to implement.
In particular, the number of arithmetic operations per iteration scales like $\bigO{d^3}$ and the total storage required for the algorithm is $\bigO{d^2}$.
The dominant computational cost is solving the ridge regression.
By maintaining and updating a few intermediate variables, the solution to the ridge regression may be obtained by solving a $d$-by-$d$ linear system, requiring $\bigO{d^3}$ arithmetic operations and $\bigO{d^2}$ storage.
The estimated online squared residuals $\estores{t}{1}$ and $\estores{t}{0}$ can be updated using $\bigO{d}$ operations.
The minimization required to obtain $p_t$ is convex and one-dimensional, so simple root finding algorithms (e.g. bisection) may be used for this purpose.

\subsection{Neyman Regret Guarantee}

The first main result of this paper is the following theorem, which shows that the Neyman regret converges to zero at a $T^{-1/2} R$ rate under \ourdesign{}.

\expcommand{\Neymanregret}{%
	Under Assumptions~\mainref{assumption:moments}-\mainref{assumption:maximum-radius} and Condition \mainref{condition:sigmoid}, there exists a constant $C>0$ such that the Neyman Regret under \ourdesign{} is bounded as
	\[
	\neymanregret_T \leq C\cdot T^{-1/2} R \enspace.
	\]
}
\begin{theorem} \label{thm:neyman-regret}
	\Neymanregret
\end{theorem}

In the theorem above, $C$ is a constant depending only on the constants $c_0$, $c_1$, $c_2$, $c_3$, $\gamma_0$ and $b_1$, $b_2$, $b_3$ appearing in Assumptions~\mainref{assumption:moments}-\mainref{assumption:maximum-radius} and Condition \mainref{condition:sigmoid}, respectively.
We remark that $C$ is a small polynomial function of these constants, unlike the exponential dependence in the analysis of \citet{Dai2023Clip}.
In a slight abuse of notation, we use $C$ throughout the paper to express any constant which depends only on these constants in the assumptions; however, all constants are made explicit in the appendix.

In the regime where the entries of each covariate vector are considered to be uniformly bounded, i.e. $R^2 \asymp d$, then Theorem~\ref{thm:neyman-regret} demonstrates a regret of $\sqrt{d/T}$, as stated in the abstract and introduction.
In the context of our lower bound (Theorem~\ref{theorem:neyman-regret-lower-bound}),
Theorem~\ref{thm:neyman-regret} establishes that $T^{-1/2}$ is optimal $T$ dependence for Neyman regret under our regularity conditions.
In particular, this removes the sub-polynomial factors in the probability clipping approach of \citet{Dai2023Clip}.

Although the proof of Theorem~\ref{thm:neyman-regret} is lengthy and technical, the central ideas of the proof can be conveyed by several key lemmas.
In the remainder of Section~\ref{sec:our-design}, we will present certain lemmas with the aim of providing the intuition underlying our approach.
The full proof of Theorem~\ref{thm:neyman-regret} appears in Section~\suppref{sec:supp-neyman-regret-analysis} of the appendix.

\subsection{Probability Regret} \label{sec:probability-regret}

%
%

\expcommand{\probregretinitialbound}{
	The expected probability regret can be bounded as
	\[
	\E{ \regretprob_T }
	\leq \frac{1}{\eta_{T+1}} \psi(u^*)
	+ \sum_{t=1}^T \eta_t \E[\Bigg]{ \frac{ \paren[\big]{\grad \estsigloss{t}(u_t)}^2 }{2(1+|u_t|)} }
	\enspace.
	\] 
}

\begin{figure}[t]
	\centering  
	\begin{subfigure}{0.31\textwidth}
		\centering
		\includegraphics[width=\linewidth]{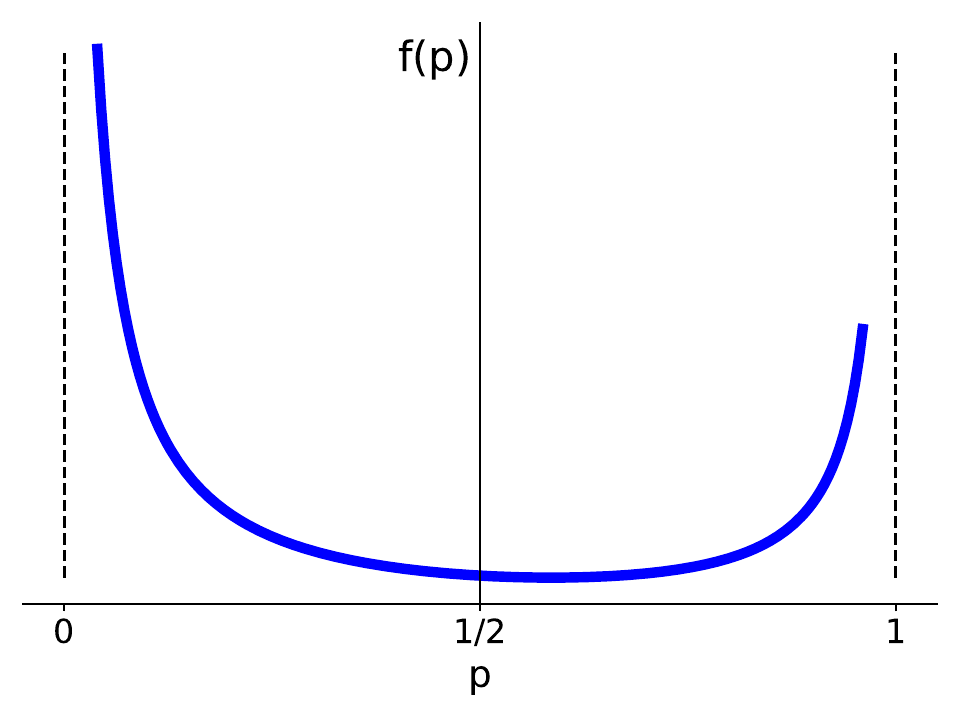}
		\caption{Probability Loss}
		\label{fig:prob-loss}
	\end{subfigure}
	\hfill
	\begin{subfigure}{0.31\textwidth}
		\centering
		\includegraphics[width=\linewidth]{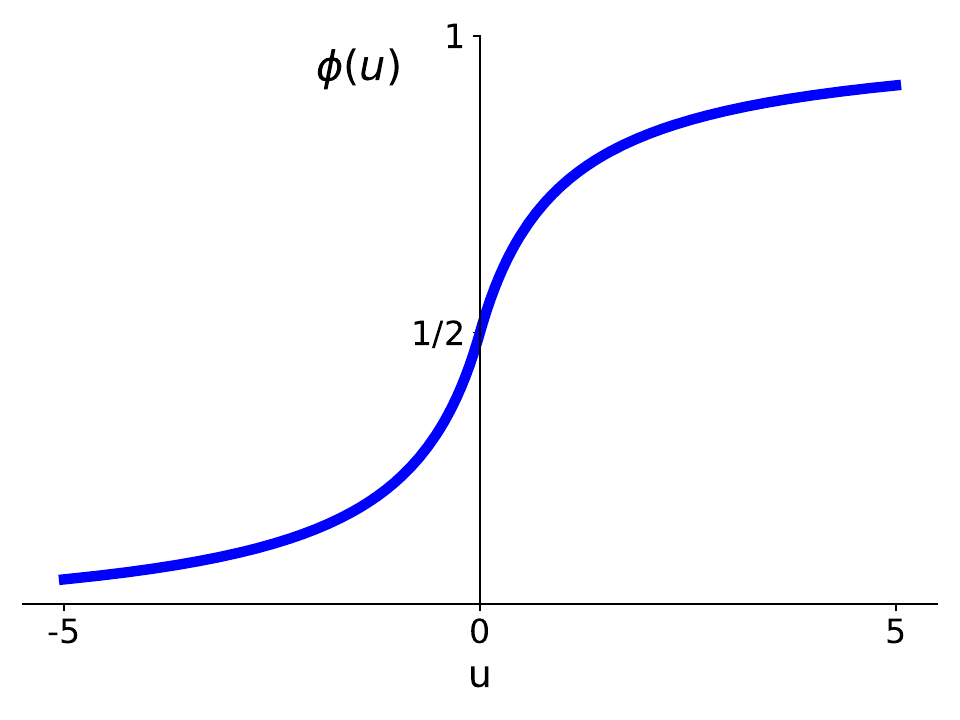}
		\caption{Sigmoid Transformation}
		\label{fig:sigmoid}
	\end{subfigure}
	\hfill
	\begin{subfigure}{0.31\textwidth}
		\centering
		\includegraphics[width=\linewidth]{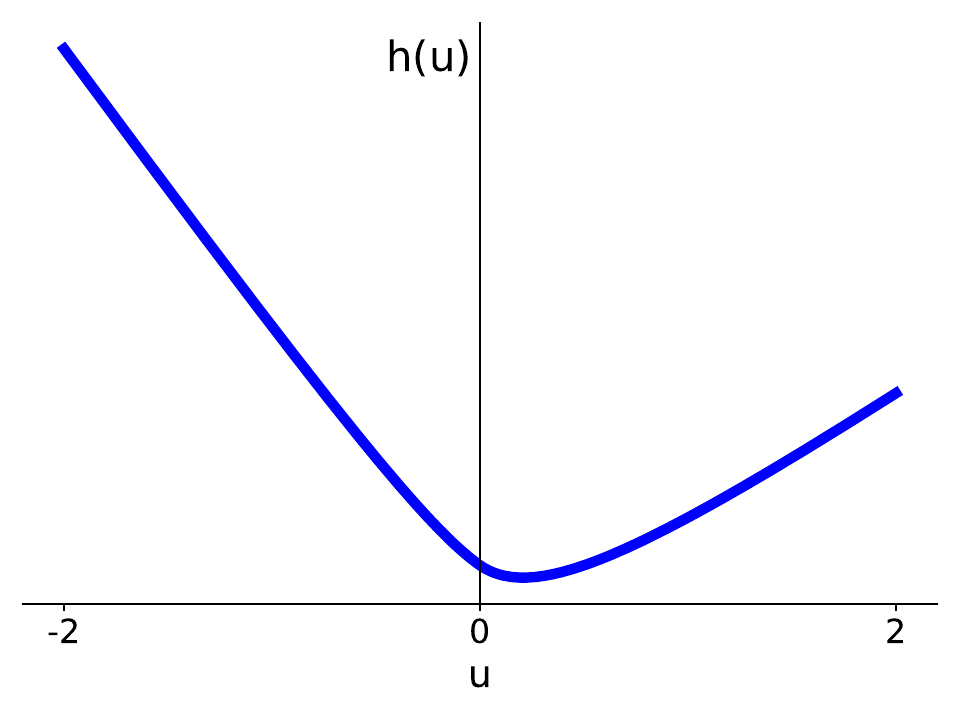}
		\caption{Transformed Loss}
		\label{figtransformed-loss}
	\end{subfigure}
	\caption{
		\textbf{Effect of Sigmoidal Transformation}:
		Fig~\ref{fig:prob-loss} shows the probability loss $f_t(p)$, which has gradients that explode at the boundary.
		Fig~\ref{fig:sigmoid} shows the algebraic sigmoid transformation $\phi(u)$.
		Fig~\ref{figtransformed-loss} shows the transformed loss function $h_t(u)$, which has bounded gradients.
	}
	\label{fig:loss-funcs}
\end{figure}

We begin by showing how to bound the expected probability regret $\E{ \regretprob_T }$ under the \ourdesign{} design.
The key idea is the use of a \emph{sigmoidal transformation} of the regret problem.
Rather than selecting a probability $p_t \in (0,1)$, we consider the equivalent selection of a decision variable $u_t \in \Reals$ through the sigmoid function $\phi : \Reals \to (0,1)$, i.e. $p_t = \phi(u_t)$.
For each iteration $t \in [T]$, we define the \emph{sigmoid probability loss} function $\sigloss{t} : \Reals \to \Reals_+$ as the composition $\sigloss{t} = \probloss{t} \circ \phi$ which may be expressed as:
\[
\sigloss{t}(u) = 
\paren[\Big]{ y_t(1) - \iprod{ \xv_t , \bv_t(1) } }^2 \cdot \frac{1}{\phi(u)} 
+
\paren[\Big]{ y_t(0) - \iprod{ \xv_t , \bv_t(0) } }^2 \cdot \frac{1}{1-\phi(u)} 
\enspace.
\]
Using the equivalence $p_t = \phi(u_t)$, we can re-express the probability regret as
\[
\regretprob_T = \sum_{t=1}^T \sigloss{t}(u_t) - \sum_{t=1}^T \sigloss{t}(u^*)
\enspace,
\]
where $u^* = \phi^{-1}(p^*)$ is the transformed Neyman Allocation.
Parts 1-3 of Condition~\mainref{condition:sigmoid} ensure that the sigmoid probability loss functions $\sigloss{t}$ are convex and have uniformly bounded gradients.
In contrast, the original probability loss functions $f_t$ have gradients which blow up at the boundary of the interval, i.e. $\abs{ \grad f_t(p) } \to \infty$ as $p \to 1$ or $p \to 0$.
In this way, the benefit of the sigmoidal transformation is to transform an ill-conditioned constrained problem into a well-conditioned unconstrained problem.
See Figure~\ref{fig:loss-funcs} for an illustration of the well-conditioning properties of the sigmoidal transformation.

The sigmoid probability loss functions $\sigloss{t}$ depend on both of the potential outcomes, so they cannot be directly observed.
In light of this fact, the \ourdesign{} design uses the adaptive IPW weighting technique to construct \emph{estimated sigmoid loss functions} $\estsigloss{t}$ as
\[
\estsigloss{t}(u)
= \paren[\Big]{ Y_t - \iprod{ \xv_t , \bv_t(1) } }^2 \cdot \frac{\indicator{Z_t = 1}}{p_t} \cdot \frac{1}{\phi(u)} 
+
\paren[\Big]{ Y_t - \iprod{ \xv_t , \bv_t(0) } }^2 \cdot \frac{\indicator{Z_t = 0}}{1-p_t} \cdot \frac{1}{1-\phi(u)}
\]
and then selects the treatment assignment probability $p_t = \phi(u_t)$ where $u_t$ is the minimizer of the estimated sigmoid loss functions
\[
u_t 
= \argmin_{u \in \Reals} \sum_{s < t} \estsigloss{s}(u) + \eta_t^{-1} \psi(u) 
:= \argmin_{u \in \Reals} \tprobloss{t}(u)
\enspace,
\]
with regularization $\psi(u) = \frac{1}{2} u^2 + |u|^3$.
In this way, we can understand the \ourdesign{} design as working directly in the sigmoidal transformation.
The following lemma shows that these estimated sigmoid loss functions $\estsigloss{t}$ are conditionally unbiased for the sigmoid loss functions $\sigloss{t}$.

\expcommand{\unbiasedsigmoidloss}{
	The estimated sigmoid loss functions are conditionally unbiased: for all $\filt_{t-1}$-measurable random variables $u$, $\E{ \estsigloss{t}(u) \mid \filt_{t-1} } = \sigloss{t}(u)$ almost surely.
}

\begin{lemma} \label{lemma:unbiased-sigmoid-loss}
	\unbiasedsigmoidloss
\end{lemma}

An immediate corollary is the following: while the probability regret is defined in terms of the true sigmoid loss functions $\sigloss{t}$, its expectation can be expressed via the expected regret under the (random) estimated sigmoid loss functions $\estsigloss{t}$. 

\expcommand{\expprobloss}{
	The expected probability regret under \ourdesign{} is
	\[
	\E{ \regretprob_T }
	= \E[\Big]{ \sum_{t=1}^T \estsigloss{t}(u_t) - \sum_{t=1}^T \estsigloss{t}(u^*) }
	\enspace.
	\]
}
\begin{corollary} \label{corollary:exp-prob-loss}
	\expprobloss
\end{corollary}

With Corollary~\ref{corollary:exp-prob-loss} in hand, we may use standard techniques in the online optimization literature to bound the regret inside the expectation, e.g. see Chapter~7 of \citep{Orabona2023Modern}.
Given that we are working with time-varying step sizes $\eta_t$, it will be convenient to define
\[
\ptprobloss{t}(u) = \sum_{s=1}^{t-1} \estsigloss{s}(u) + \eta_{t-1}^{-1} \psi(u)
\quadand
\widetilde{u}_{t} = \argmin_{u \in \Reals} \ptprobloss{t}(u) 
\enspace.
\]
The function $\ptprobloss{t}$ is similar to the objective $\tprobloss{t}$ which is minimized at the $t$th iteration in \ourdesign{}, except that the previous step size $\eta_{t-1}$ is used in the regularization term.
Note that $\ptprobloss{t}$ and $\widetilde{u}_{t}$ exist only for the purpose of analysis, i.e. they are not computed in \ourdesign{}.
The following lemma bounds the expected probability regret in terms of these quantities.

\expcommand{\expproblossub}{
	The expected probability regret is bounded as
	\[
		\E{ \regretprob_T }
		\leq \frac{1}{\eta_{T+1}} \psi(u^*)
		+ \sum_{t=1}^T \E[\Big]{ \ptprobloss{t+1}(u_t) - \ptprobloss{t+1}( \widetilde{u}_{t+1} ) }
		\enspace.
	\]
}

\begin{lemma}\label{lemma:exp-prob-regret-initial-oco-bound}
	\expproblossub
\end{lemma}

The first term in the bound is the regularization term, which features both the final step size $\eta_{T+1}^{-1} = T^{1/2} R_T$ and the regularizer $\psi$ evaluated at the optimal $u^*$.
By Assumptions~\mainref{assumption:moments} and \ref{assumption:covariate-regularity}, this first term is of the desired order $T^{1/2} R$.
Thus, we can focus our attention on the second term.

The second term contains differences $\ptprobloss{t+1}(u_t) - \ptprobloss{t+1}( \widetilde{u}_{t+1} )$ which measure the suboptimality gap of $u_t$ chosen by \ourdesign{} on the modified loss $\ptprobloss{t+1}$.
We present the usual analysis from online convex optimization for these suboptimality terms \citep[see, e.g.][]{Hazan2016Introduction, Orabona2023Modern}.
Using first order convexity conditions, rearranging terms, and the fact that $u_t$ minimizes $\tprobloss{t}$ so that $\grad \tprobloss{t}(u_t) = 0$, we have that
\begin{align*}
	& \ptprobloss{t+1}(u_t) - \ptprobloss{t+1}( \widetilde{u}_{t+1} ) \\
	&\quad = \sum_{s=1}^{t} \paren[\Big]{\estsigloss{s}( u_t ) - \estsigloss{s}( \widetilde{u}_{t+1} )}
	+ \eta_{t}^{-1} \paren[\Big]{ \psi(u_t) - \psi(\widetilde{u}_{t+1}) } \\
	&\quad \leq \iprod[\Big]{ \sum_{s=1}^t \grad \estsigloss{s}(u_t) , u_t - \widetilde{u}_{t+1} }
	+ \eta_{t}^{-1} \paren[\Big]{ \psi(u_t) - \psi(\widetilde{u}_{t+1}) } \\
	&\quad = \iprod[\Big]{ \grad \tprobloss{t}(u_t) + \grad \estsigloss{t}(u_t) ,  u_t - \widetilde{u}_{t+1} }
	+ \eta_{t}^{-1} \paren[\Big]{ \psi(u_t) - \psi(\widetilde{u}_{t+1}) - \iprod{ \grad \psi(u_t), u_t-\widetilde{u}_{t+1} }} \\
	&\quad = \iprod{ \grad \estsigloss{t}(u_t) ,  u_t - \widetilde{u}_{t+1} }
	- \eta_{t}^{-1} \paren[\Big]{ \psi(\widetilde{u}_{t+1})- \psi(u_t)  - \iprod{ \grad \psi(u_t), \widetilde{u}_{t+1}-u_t }} \\
	&\quad = \iprod{ \grad \estsigloss{t}(u_t) ,  u_t - \widetilde{u}_{t+1} } - \frac{\eta_{t}^{-1}}{2} (u_t - \widetilde{u}_{t+1})^2 \psi''(\widebar{u}_t) \\
	&\quad \leq \frac{\eta_t}{2} \cdot \frac{ \paren{ \grad \estsigloss{t}(u_t) }^2} { \psi''(
		\widebar{u}_t) }
	\enspace,
\end{align*}
where the second to last line follows from Taylor's theorem, i.e. there exists a $\widebar{u}_t$ between $u_t$ and $\widetilde{u}_{t+1}$ such that 
$
\psi(\widetilde{u}_{t+1})- \psi(u_t)  - \iprod{ \grad \psi(u_t), \widetilde{u}_{t+1}-u_t } 
= \frac{1}{2}(u_t - \widetilde{u}_{t+1})^2 \psi''(\widebar{u}_t) 
$.
The last line follows by taking the supremum over all possible values of $u_t - \widetilde{u}_{t+1}$.
If we further make the approximation that $\psi''(\widebar{u}_t) \approx \psi''(u_t)$, then we can approximate the conditional expectation of the suboptimality terms by
\begin{equation} \label{eq:regularizer-intuition}
\E[\Big]{ \ptprobloss{t+1}(u_t) - \ptprobloss{t+1}( \widetilde{u}_{t+1} ) \mid \filt_{t-1} }
\approx 
\frac{\eta_t}{2} \cdot  \frac{ \E[\Big]{ \paren{ \grad \estsigloss{t}(u_t) }^2 \mid \filt_{t-1}} } { \psi''(u_t) }
\enspace.
\end{equation}
Formally speaking, the approximation \eqref{eq:regularizer-intuition} is not rigorous because there is no justification for the approximation $\psi''(\widebar{u}_t) \approx \psi''(u_t)$.
However, let us now consider such an approximation for the sake of developing the intuition of our analysis.
Once the intuition is made clear, we will formally justify why we can use the approximation \eqref{eq:regularizer-intuition} as a valid upper bound, up to constants.

The approximation \eqref{eq:regularizer-intuition} provides guidance on the choice of regularizer: we will choose $\psi$ so that $\psi''$ offsets any $u$ dependence in $\grad \estsigloss{t}$.
Although the gradients of the true sigmoid loss $\sigloss{t}$ are uniformly bounded, the estimated sigmoid loss $\estsigloss{t}$ uses an inverse probability weighting, so its gradients will depend on $p_t$ and thus $u_t$.
The expected squared gradients of the estimated sigmoid loss can be bounded
\[
\E[\Big]{  \paren{ \grad \estsigloss{t}(u_t) }^2 \mid \filt_{t-1} }
\leq \paren{ C_1 + C_2 \abs{ u_t } } \cdot \paren[\Big]{
	\braces[\big]{ y_t(1) - \iprod{ \xv_t, \bv_t(1) } }^4
	+\braces[\big]{ y_t(0) - \iprod{ \xv_t, \bv_t(0) } }^4
}
,
\]
where the constants $C_1$ and $C_2$ depend on the choice of sigmoid, i.e. constants in Condition~\mainref{condition:sigmoid}.
The $C_1 + C_2 \abs{u_t}$ term arises because it is proportional to $\max \setb{ 1/p_t, 1/(1-p_t) }$ which appears from the inverse probability weighting.
If we were to use the usual squared regularizer $\psi(u) = \frac{1}{2} u^2$, then $\psi''(u) = 1$ so that the right hand side of approximation \eqref{eq:regularizer-intuition} has a numerator which grows with $u_t$ and a denominator which is constant.
In this way, the expected suboptimality gap term $\E{\ptprobloss{t+1}(u_t) - \ptprobloss{t+1}( \widetilde{u}_{t+1} )}$ could depend on the magnitude of $u_t$.
On the other hand, our choice of squared + cubic regularizer $\psi(u) = \frac{1}{2} u^2 + \abs{u}^3$ results in $\psi''(u) = 1 + 6 \abs{u}$ which cancels out (up to constants) the $u$ dependence from the expected squared gradient.
Using the stronger regularizer, we have that the expected suboptimality gap $\E{\ptprobloss{t+1}(u_t) - \ptprobloss{t+1}( \widetilde{u}_{t+1} )}$ will not directly depend on $u_t$, and instead depend only on the fourth moments of the online residuals.

In order to make these arguments rigorous, we need to formalize the approximation that $\psi''(\widebar{u}_t) \approx \psi''(u_t)$.
To this end, we use the geometry of Bregman divergences.
Bregman divergences play a crucial role in the analysis of optimization algorithms because they express the relevant curvature of the objective functions \citep{Hazan2016Introduction,Orabona2023Modern}.
Given a differentiable convex function $f : \mathcal{X} \to \Reals$ on a convex set $\mathcal{X}$, the \emph{Bregman divergence} $\bregf{f} : \mathcal{X} \times \mathcal{X} \to \Reals$ is defined as $\breg{f}{u}{v} = f(u) - f(v) - \iprod{ \grad f(v) , u-v}$.
The Bregman divergence measures local geometry of the function $f$ at different pairs of points through the gap between $f(u)$ and the convex lower bound $f(v) + \iprod{ \grad f(v), u-v }$.

The analysis of the suboptimality gap above actually used Taylor's theorem to characterize the Bregman divergence, i.e. $\breg{\psi}{\widetilde{u}_{t+1}}{u_t} = \frac{1}{2}(u_t - \widetilde{u}_{t+1})^2 \psi''(\widebar{u}_t)$.
This is known as the local norm interpretation of the Bregman divergence.
The downside of local norm interpretation is that nothing is known about the midpoint $\widebar{u}_t$, and so $\psi''(\widebar{u}_t)$ is not clearly related to $\psi''(u_t)$.
The following lemma provides an alternative approach: the Bregman divergence of our regularizer $\psi$ is lower bounded in a way that depends only on the two points in question, without any reference to a midpoint.

\expcommand{\breglb}{
	The Bregman divergence is lower bounded:
	$\breg{\psi}{v}{u} \geq \frac{1}{2} (v-u)^2 (1+|u|)$.
}
\begin{lemma} \label{lemma:breg-lb}
	\breglb
\end{lemma}

Lemma~\ref{lemma:breg-lb} shows that the Bregman divergence is lower bounded by a quantity which, up to constants, actually mimics the desired behavior, i.e. $(u_t - \widetilde{u}_{t+1})^2 \psi''(u_t)$.
This justifies the non-rigorous approximation $\psi''(\widebar{u}_t) \approx \psi''(u_t)$ made earlier in the section.
In fact, using Lemma~\ref{lemma:breg-lb} in the analysis of the suboptimality terms above provides a fully rigorous argument.
For these reasons, Lemma~\ref{lemma:breg-lb} is a key technical result which facilitates the formalization of approximation \eqref{eq:regularizer-intuition} and is ultimately responsible for the improved $T^{-1/2} R$ rates of Neyman regret.

At first sight, \ourdesign{} may appear similar to imposing a heavy regularizer, e.g.
\[
\Psi(p) = p^{-2} + (1-p)^{-2} + p^{-3} + (1-p)^{-3}.
\]
However, the key advantage of the sigmoid transformation is different: Lemma~\ref{lemma:breg-lb} shows that the Bregman divergence $\breg{\psi}{u_{t+1}}{u_t}$ in the transformed space can be controlled \emph{globally}, even when the corresponding probabilities $p_t = \phi(u_t)$ and $p_{t+1} = \phi(u_{t+1})$ are far apart.
In contrast, working with such a regularizer $\Psi$ directly in probability space would force us to rely on $\breg{\Psi}{p_{t+1}}{p_t}$, which could only be controlled \emph{locally}, due to the classical localization constraint of FTRL: the divergence can only be well-approximated by its second-order expansion when $p_{t+1}$ remains close to $p_t$.
The sigmoid geometry circumvents this limitation by expressing potentially large moves in probability space as well-behaved movements in the unconstrained $u$-space.
We believe this transformation technique will be of more general interest beyond the problem of Adaptive Neyman Allocation considered here.

The following lemma provides the formal statement that follows directly from the arguments above.
We see that the second term in the regret analysis is bounded in terms of the fourth moments of the online residuals.

\expcommand{\probabilitygradientfourthmomentbound}{
	There exists a constant $C>0$ depending on the sigmoid Condition~\mainref{condition:sigmoid} such that the suboptimality term is bounded by the expected fourth moment of online residuals:
	\[
	\sum_{t=1}^T \E[\Big]{ \ptprobloss{t+1}(u_t) - \ptprobloss{t+1}( \widetilde{u}_{t+1} )}
	\leq 
	C \cdot \sum_{k \in \setb{0,1}} \E[\Bigg]{ \sum_{t=1}^T \eta_t \paren[\big]{ y_t(k) - \iprod{ \xv_t, \bv_t(k) } }^4 }
	\enspace.
	\]
}

\begin{lemma} \label{lemma:probability-gradient-fourth-moment-bound}
	\probabilitygradientfourthmomentbound
\end{lemma}

In order to bound the probability regret, the last step is to ensure that the fourth moment of the online residuals is bounded.
To bound the fourth moments of the online residuals, we develop a specialized ``prediction tracking'' technique.
To the best of our knowledge, this technique is new to this work and does not exist in the online optimization literature.
The proof sketch of the prediction tracking technique is delayed until Section~\ref{sec:prediction-tracking}, but the following lemma contains the main result.

\expcommand{\boundedfourthmoment}{
	Under Assumptions~\mainref{assumption:moments}-\mainref{assumption:maximum-radius} and Condition~\mainref{condition:sigmoid}, there exists a constant $C > 0$ such that
	for each treatment $k \in \setb{0,1}$, the fourth moment of the online residuals can be bounded as
	\[
	\E[\Bigg]{ \sum_{t=1}^T \eta_t \braces[\big]{ y_t(k) - \iprod{ \xv_t, \bv_t(k) } }^4 }
	\leq C \cdot T^{1/2} R\enspace.
	\]
}
\begin{lemma} \label{lemma:bounded-fourth-moments}
	\boundedfourthmoment
\end{lemma}

We are now ready to provide a bound on the probability regret.
Combining Lemmata~\ref{lemma:exp-prob-regret-initial-oco-bound}, \ref{lemma:probability-gradient-fourth-moment-bound}, and \ref{lemma:bounded-fourth-moments}, we arrive at the following bound on the expected probability regret:

\expcommand{\probregretbound}{
	Under Assumptions \mainref{assumption:moments}-\mainref{assumption:maximum-radius} and Condition \mainref{condition:sigmoid}, there exists a constant $C>0$ such that the expected probability regret is bounded as $\E{ \regretprob_T} \leq C\cdot T^{1/2} R$.
}
\begin{proposition} \label{prop:prob-regret-bound}
	\probregretbound
\end{proposition}

\subsection{Prediction Regret} \label{sec:prediction-regret}

In this section, we show how to bound the prediction regret.
The prediction loss $\newpredloss{t}$ depends on both of the potential outcomes, which are not simultaneously observed.
Using the adaptive inverse probability weighting technique, we define the \emph{estimated prediction loss} $\newestpredloss{t} : \Reals^d \times \Reals^d \to \Reals_+$ as
\begin{align*}
&\newestpredloss{t}\paren[\big]{ \bv(1), \bv(0) } \\
& = \paren[\Bigg]{
	\braces[\Big]{ Y_t \cdot \frac{\indicator{Z_t=1}}{p_t} - \iprod{ \xv_t , \bv(1) } } \cdot \sqrt{ \frac{\olsres{0}}{\olsres{1}} }
	+
	\braces[\Big]{ Y_t \cdot \frac{\indicator{Z_t=0}}{1-p_t} - \iprod{ \xv_t , \bv(0) } } \cdot \sqrt{ \frac{\olsres{1}}{\olsres{0}} }
}^2
.
\end{align*}
Define also the regularizer function $m : \Reals^d \times \Reals^d \to \Reals^+$ as
\[
m \paren[\big]{ \bv(1), \bv(0) }
= \norm[\Bigg]{ 
	\sqrt{\frac{\olsres{0}}{\olsres{1}}}  \cdot \bv(1)
	+ \sqrt{\frac{\olsres{1}}{\olsres{0}}}  \cdot \bv(0)
}^2
\enspace.
\]
Recall that \ourdesign{} selects the two linear predictors $\bv_t(1)$ and $\bv_t(0)$ separately via a ridge regression with IPW-weighted outcomes, i.e. see Algorithm~\ref{alg:sigmoid-ftrl}.
What is perhaps surprising is that the separate selection of linear predictors $\bv_t(1)$ and $\bv_t(0)$ actually corresponds to joint selection via the FTRL principle applied to the estimated prediction loss functions using the regularizer defined above:
\begin{equation} \label{eq:joint-minimizers-are-separate-minimizers}
\bv_t(1), \bv_t(0) =
\argmin_{\bv(1), \bv(0)} \sum_{s=1}^{t-1} \newestpredloss{s}\paren[\big]{ \bv(1), \bv(0) }
+ \eta_t^{-1}  m \paren[\big]{ \bv(1), \bv(0) }
\enspace.
\end{equation}
Note that the program above includes both a cross term and the optimal residuals $\olsres{1}$ and $\olsres{0}$ in its definition, neither of which appear in the usual ridge loss.
In the same way that the individual OLS predictors jointly minimize the AIPW variance  (Proposition~\ref{prop:oracle}), it is also true that the individual ridge predictors jointly minimize the program above.
Establishing \eqref{eq:joint-minimizers-are-separate-minimizers} is the main argument in the proof of Lemma~\ref{lemma:standard-oco-for-pred-regret} which appears in Section~\suppref{section:C4} of the supplement.

The next steps of the analysis of the prediction regret mirror the standard FTRL arguments also used in the analysis of the probability regret.
The following lemma demonstrates that the estimated prediction loss functions $\newestpredloss{t}$ are equivalent, in expectation, to working with the actual prediction loss functions $\newpredloss{t}$:

\expcommand{\predregretunbiased}{%
	For each iteration $t$, the following conditional unbiasedness holds:
	\[
	\E[\big]{ \newestpredloss{t}\paren[\big]{ \bv_t(1), \bv_t(0) } - \newestpredloss{t}\paren[\big]{ \optbv(1), \optbv(0) }\mid \filt_{t-1} }
	= \newpredloss{t}\paren[\big]{ \bv_t(1), \bv_t(0) } - \newpredloss{t}\paren[\big]{ \optbv(1), \optbv(0) } \enspace.
	\]
}
\begin{lemma} \label{lemma:pred-regret-unbiased}
	\predregretunbiased
\end{lemma}

\noindent
As in Section~\ref{sec:probability-regret}, we introduce auxiliary functions and variables
\begin{align*}
\newptpredloss{t}\paren[\big]{ \bv(1), \bv(0) } 
&= \sum_{s=1}^{t-1} \newestpredloss{s} \paren{ \bv(1), \bv(0) } + \eta_{t-1}^{-1} m \paren[\big]{ \bv(1), \bv(0) } \quadand \\
\woptbv_t(1), \woptbv_t(0) &= \argmin_{\bv(1), \bv(0)} \newptpredloss{t} \paren[\big]{ \bv(1), \bv(0) } 
\enspace.
\end{align*}
The auxiliary objective $\newptpredloss{t}$ is similar to the regularized objective function $\newtpredloss{t}$ except that it uses the previous iteration's step size $\eta_{t-1}^{-1}$.
The pair of linear predictors $\woptbv_t(1)$, $\woptbv_t(0)$ is the minimizer of $\newptpredloss{t}$.
We emphasize that the auxiliary loss $\newptpredloss{t}$ and predictors $\woptbv_t(1)$, $\woptbv_t(0)$ are defined for the purposes of analysis and not constructed in \ourdesign{}.
As before, standard FTRL arguments yield the following regret bound:

\expcommand{\standardoco}{
	The expected prediction regret is bounded as
	\[
		\E{ \newpredregret_T }
		\leq 
		\frac{ m \paren[\big]{ \optbv(1), \optbv(0) } }{\eta_{T+1}}
		+ \E[\Bigg]{ \sum_{t=1}^T 
			\newptpredloss{t+1}\paren[\big]{ \bv_t(1), \bv_t(0) } 
			- \newptpredloss{t+1}\paren[\big]{ \woptbv_{t+1}(1), \woptbv_{t+1}(0) }  
		 } 
	 \enspace.
	\]
}
\begin{lemma}\label{lemma:standard-oco-for-pred-regret}
	\standardoco	
\end{lemma}

The upper bound in Lemma~\ref{lemma:standard-oco-for-pred-regret} has two parts.
The first term represents the cost of regularization.
By Assumptions~\mainref{assumption:moments} and \ref{assumption:covariate-regularity}, the norms of the optimal linear predictors $\optbv(1)$ and $\optbv(0)$ will be of constant order.
Moreover, the final step size of the algorithm is chosen to be $\eta_{T+1}^{-1} = T^{1/2} R_T$.
Thus, the regularization term will be order $T^{1/2} R$, as desired.

The second term again contains the differences which measure suboptimality of the predictors $\bv_t(1), \bv_t(0)$ on the next objective function $\newptpredloss{t+1}$.
Such analyses{\textemdash}including the analysis of probability regret in Section~\ref{sec:probability-regret}{\textemdash}usually proceed via an application of first order convexity of the objective, i.e., using the shorthand $\bv_t = \paren{ \bv_t(1), \bv_t(0) }$  for concatenation:
\[
\newptpredloss{t+1}( \bv_{t} )
- \newptpredloss{t+1}( \woptbv_{t+1} )
\leq \iprod{ \grad \newptpredloss{t+1}( \bv_{t} ) , \bv_{t} - \woptbv_{t+1} }
\enspace.
\]
However, this will not suffice for our analysis because it will not allow us to leverage the regularity in the covariates.
Indeed, using the first order convexity bound will result in each covariate vector $\xv_t$ being considered separately, whereas the covariates are only assumed to be regular in aggregate as measured by the invertibility of the covariate matrix (Assumption~\ref{assumption:covariate-regularity}).
For this reason, we must analyze this difference of terms appearing in Lemma~\ref{lemma:standard-oco-for-pred-regret} directly.
To this end, define the \emph{per iteration leverage scores} as
\[
\Pi_{t,s} = \xv_t^\tran \paren[\Big]{ \xM_{t-1}^\tran \xM_{t-1} + \eta_t^{-1} \unitM }^{-1} \xv_s
\enspace,
\]
where $\xM_{t-1} = \bracket{ \xv_1, \dots, \xv_{t-1} }^\tran $, so that $\xM_{t-1}^\tran \xM_{t-1} = \sum_{s \leq t-1} \xv_s \xv_s^\tran$.
Unlike the usual leverage score, the per iteration leverage score $\Pi_{t,s}$ is not symmetric.
The per-iteration leverage scores also differ from the usual leverage score because $\xv_t$ is not part of the matrix $\xM_{t-1}$ which defines the quadratic form.
Nevertheless, the per-iteration leverage score $\Pi_{t,s}$ arises naturally when expressing the online ridge predictors as linear estimators:

\expcommand{\ridgepredexactform}{
	At each iteration $t$, the online ridge predictions admit the decompositions
	\[
	\iprod{ \xv_t, \bv_t(1) } = \sum_{s=1}^{t-1} \Pi_{t,s} y_s(1) \frac{\indicator{Z_s=1}}{p_s}
	\quadand
	\iprod{ \xv_t, \bv_t(0) } = \sum_{s=1}^{t-1} \Pi_{t,s} y_s(0) \frac{\indicator{Z_s=0}}{1-p_s}
	\enspace.
	\]
}
\begin{lemma} \label{lemma:ridge-pred-exact-form}
	\ridgepredexactform
\end{lemma}

Using matrix algebra, we can write the suboptimality gap in terms of the per-iteration leverage score and the estimated prediction loss as follows:

\expcommand{\predlossdiff}{
	The suboptimality gap can be written as:
	\[
	\newptpredloss{t+1}\paren[\big]{ \bv_t(1), \bv_t(0) } 
	- \newptpredloss{t+1}\paren[\big]{ \woptbv_{t+1}(1), \woptbv_{t+1}(0) } 
	= \frac{ \Pi_{t,t} }{1 + \Pi_{t,t}} \cdot \newestpredloss{t} \paren[\big]{ \bv_t(1), \bv_t(0) }
	\enspace.
	\]
}
\begin{lemma}\label{lemma:ptpredloss-diff}
	\predlossdiff
\end{lemma}

Similar derivations of this type have appeared in the early literature on online optimization \citep[see e.g.][]{Azoury2001Relative, Zhdanov2010Competing}.
The benefit of Lemma~\ref{lemma:ptpredloss-diff} is that it allows us to make use of the regularity of the entire sequence of covariates (Assumption~\ref{assumption:covariate-regularity}) through the per-iteration leverage scores.

The remainder of the proof for bounding the prediction regret can be carried out using various properties of the per-iteration leverage scores (e.g. bounds on the partial sums of higher moments) and control of the moments of inverse probabilities, which are  presented in the next section.
We defer to Section \suppref{section:C4} of the appendix for details and summarize the main result here:

\expcommand{\predregret}{
	Under Assumptions \mainref{assumption:moments}-\mainref{assumption:maximum-radius} and Condition~\mainref{condition:sigmoid}, there exists a constant $C>0$ such that the expected prediction regret is bounded as $\E{ \newpredregret_T} \leq C \cdot  T^{1/2} R$.
}
\begin{proposition} \label{prop:pred-regret}
	\predregret
\end{proposition}

The proof of Theorem~\ref{thm:neyman-regret} is completed by applying the Neyman Regret decomposition (Lemma~\ref{lemma:regret-decomposition}) with the bounds on the probability regret (Proposition~\ref{prop:prob-regret-bound}) and prediction regret (Proposition~\ref{prop:pred-regret}).

\subsection{Fourth Moments: Prediction Tracking} \label{sec:prediction-tracking}

In this section, we provide a proof sketch for Lemma~\ref{lemma:bounded-fourth-moments}, which bounds the fourth moments of the online residuals.
Our high level approach will be as follows: we introduce a sequence of deterministic linear predictors and argue that (i) the online residuals of this deterministic predictor sequence are bounded, (ii) the difference between the deterministic and the adaptive residuals is a lower-order term.
We refer to this proof technique as \emph{prediction tracking}.

For each treatment $k \in \setb{0,1}$, we define the sequence of \emph{full information predictors} $\setb{ \optbv_t(k) }_{t=1}^{T}$ as the solutions to the ridge regression:
\[
\optbv_t(k) = \argmin_{ \bv \in \Reals^d } \sum_{s < t} \paren[\Big]{ y_s(k) - \iprod{\xv_s, \bv} }^2 + \eta_t^{-1} \norm{\bv}^2 \enspace,
\]
where the penalty term $\eta_t$ is as chosen in \ourdesign{}.
We remark that the sequence of full-information predictors $\setb{ \optbv_t(k) }_{t=1}^{T}$ is deterministic because it depends on the potential outcomes, not the observed outcomes, hence the name ``full-information''.
As such, we cannot hope to compute the full-information predictors directly{\textemdash}they are defined merely for the purposes of analysis.
The purpose of introducing these full-information predictors is that they are a deterministic sequence around which the adaptively chosen predictors will concentrate.
The following unbiasedness property is a direct consequence of the adaptive IPW weighting used in the ridge regression step of \ourdesign{}:

\expcommand{\predictorexpectation}{
	For each iteration $t \in [T]$ and treatment $k \in \setb{0,1}$,
	the adaptive linear predictors are unbiased for the full-information predictors: $\E{\bv_t(k)} = \optbv_t(k)$.
}
\begin{lemma} \label{lemma:predictor-expectation}
	\predictorexpectation
\end{lemma}

In Section~\suppref{section:C5} of the appendix, we show that the fourth moments of the full-information residuals are bounded as $\sum_{t=1}^T \eta_t \braces{ y_t(k) - \iprod{ \xv_t , \optbv_t(k) } }^4 = \bigO{T^{1/2}R}$ using the moment conditions and covariate regularity, i.e. Assumptions~\mainref{assumption:moments} and \ref{assumption:covariate-regularity}.
What remains to be shown is that the adaptively selected predictors' residuals ``track'' those of the deterministic full-information predictors, in the sense that they have roughly the same fourth moment.

We formalize this idea in the following way:
as a corollary of Lemma~\ref{lemma:ridge-pred-exact-form}, we can obtain a closed form expression of the \emph{prediction tracking error}:
\[
\E[\Big]{ \sum_{t=1}^T \eta_t \iprod{ \xv_t, \bv_t(1) - \optbv_t(1)}^4}
=\E[\Big]{ \sum_{t=1}^T \eta_t \paren[\Big]{ \sum_{s=1}^{t-1} \Pi_{t,s} y_s(1) \braces[\Big]{ \frac{\indicator{Z_s=1}}{p_s} - 1 } }^4 }
\enspace.
\]
The expression above demonstrates that the expected prediction tracking error depends on three things: the per-iteration leverage scores, the potential outcomes, and the fourth moments of the inverse probabilities.
Throughout the remainder of the section, we focus on the role of the moments of the inverse probabilities; in the appendix, all three of these quantities are carefully dealt with.

The inverse probabilities are controlled, in part, by our choice of sigmoidal regularization.
The following lemma illustrates the effect of the sigmoidal regularization for any sigmoid satisfying Condition~\mainref{condition:sigmoid} on the resulting inverse probabilities.

\expcommand{\pregularization}{
	Consider $A, B \geq 0$ and define $p^*$ as the minimizer of the following program:
	\[
	p^* = \argmin_{p \in (0,1)} \frac{A}{p} + \frac{B}{1-p} + \eta^{-1} \Psi(p)
	\enspace.
	\]
	Then, the minimizer $p^*$ is bounded away from $0$ and $1$ in the following sense: there exists a universal constant $C>0$ depending only on $b_1, b_2, b_3$ in Condition~\mainref{condition:sigmoid} such that
	\[
	\frac{1}{p^*} \leq 2 + C\cdot \minf[\Big]{ \paren{ \eta B }^{1/4}, \paren[\Big]{\frac{B}{A}}^{1/2} } 
	\quadand
	\frac{1}{1-p^*} \leq 2 + C\cdot \minf[\Big]{ \paren{  \eta A }^{1/4}, \paren[\Big]{\frac{A}{B}}^{1/2} } 
	\enspace.
	\]
}
\begin{lemma} \label{lemma:effect-of-p-regularization}
	\pregularization
\end{lemma}

Lemma~\ref{lemma:effect-of-p-regularization} provides two types of upper bounds on the inverse probabilities based on the quantities $A$ and $B$ in the objective and the regularization parameter $\eta^{-1}$.
If the regularization parameter $\eta^{-1}$ is large relative to $A$ and $B$, then the inverse probabilities are bounded.
Alternatively, if $A$ and $B$ are of the same magnitude so that their ratio is not too large, then the inverse probabilities are also bounded.
Lemma~\ref{lemma:effect-of-p-regularization} relies on the curvature conditions of the sigmoidal regularizer (Condition~\mainref{condition:sigmoid}).
In Lemma~\ref{lemma:effect-of-p-regularization}, we adopt the convention that $x/0 = \infty$ for $x \in \Reals$, i.e. the bound remains non-vacuous when $A$ or $B$ is $0$.

In each iteration of the \ourdesign{} design, the treatment assignment probability $p_t$ is chosen according to the minimization in Lemma~\ref{lemma:effect-of-p-regularization} where $\eta$ is the step size $\eta_t$ and $A$ and $B$ are the estimated online residuals $\estores{t-1}{1}$ and $\estores{t-1}{0}$ under treatment and control, respectively.
The first type of bounds in Lemma~\ref{lemma:effect-of-p-regularization} imply that if we want $\E{ 1/p_t }$ and $\E{ 1 /(1-p_t) }$ to be small, we need that the estimated residuals are small in expectation:

\expcommand{\pmoment}{
	There exists a universal constant $C>0$ depending on $b_1$ and $b_2$ in Condition~\mainref{condition:sigmoid} such that
	for each iteration $t \in [T]$ and any $0 \leq k \leq 4$, the $k$th moments of the inverse probabilities are bounded as
	\begin{align*}
		\E[\Big]{ \paren[\Big]{ \frac{1}{p_t} }^k } &\leq \braces[\Big]{ 2 + C\cdot\paren[\Big]{  \eta_t \E{ \estores{t-1}{0}  }  }^{1/4} }^k 
		\\
		\E[\Big]{ \paren[\Big]{ \frac{1}{1-p_t} }^k } &\leq \braces[\Big]{ 2 + C\cdot\paren[\Big]{ \eta_t \E{ \estores{t-1}{1}  }  }^{1/4} }^k 
		\enspace.
	\end{align*}
}
\begin{corollary}\label{corollary:p-moment}
	\pmoment
\end{corollary}

Corollary~\ref{corollary:p-moment} demonstrates that in order to bound the inverse probabilities, we must bound the expectations of the estimated online residuals $\E{ \estores{t-1}{1} }$ and $\E{ \estores{t-1}{0} }$.
The following lemma computes these expectations directly.
For each treatment $k \in \setb{0,1}$, we define the full-information online residuals as $A^*_t(k) = \sum_{s \leq t} \paren{ y_s(k) - \iprod{ \xv_s, \optbv_s(k) } }^2$. 

\expcommand{\expectationestore}{
	For any $t \in [T]$, the expectation of the estimated squared residuals is equal to 
	\begin{align*}
		\E{ \estores{t}{1} } &= A^*_t(1) + \sum_{s=1}^t \sum_{r=1}^{s-1} \Pi_{s,r}^2 y_r(1)^2 \E[\Big]{ \frac{1}{p_r} - 1 } \\
		\E{ \estores{t}{0} } &= A^*_t(0) + \sum_{s=1}^t \sum_{r=1}^{s-1} \Pi_{s,r}^2 y_r(0)^2 \E[\Big]{ \frac{1}{1-p_r} - 1 }
		\enspace.
	\end{align*}
}
\begin{lemma} \label{lemma:expectation-of-estores}
	\expectationestore
\end{lemma}

The magnitude of the estimated online residuals depends on all of the inverse probabilities in previous iterations.
The reason for this is that $\estores{t}{1}$ and $\estores{t}{0}$ are estimated using the adaptive IPW technique.
At first glance, it appears that we have run into a circular problem: in order to bound the inverse probabilities, we must bound the estimated online residuals, but doing so requires bounding the inverse probabilities.

Fortunately, both the inverse probabilities and the estimated online residuals can be bounded through an inductive method.
At the first iteration, the estimated online residuals are $\estores{0}{1} = \estores{0}{0} = 0$ so the treatment probability is $p_1 = 1/2$, and thus both of these quantities are well-controlled.
At later iterations, we use the inductive hypothesis that $\E{\estores{t-1}{1}}$ and $\E{\estores{t-1}{0}}$ cannot be too large together with Corollary~\ref{corollary:p-moment} to establish that the inverse probabilities $\E{1/p_t}$ and $\E{1/(1-p_t)}$ are bounded.
Using this together with the inductive hypothesis that previous inverse probabilities $\E{1/p_s}$ and $\E{1/(1-p_s)}$ are bounded for $s < t$, we are able to bound $\E{\estores{t}{1}}$ and $\E{\estores{t}{0}}$, thereby completing the inductive proof.
This method of argument, described in full detail in Sections~\suppref{section:C3} and \suppref{section:C5} of the appendix, establishes the bound on the prediction tracking error:

\expcommand{\expectedtracking}{
	Under Assumptions \mainref{assumption:moments}-\mainref{assumption:maximum-radius} and Condition \mainref{condition:sigmoid}, there exists a constant $C>0$ such that the prediction tracking error for both treatments $k \in \setb{0,1}$ can be bounded as
	 $\E{ \sum_{t=1}^T \eta_t \iprod{ \xv_t , \bv_t(k) - \optbv_t(k) }^4 } \leq C\cdot T^{1/2} R$.
}
\begin{proposition}\label{prop:expected-tracking}
	\expectedtracking
\end{proposition}
	
	\section{Asymptotically Valid Inference} \label{sec:inference}

In this section, we propose methods for asymptotically valid confidence intervals using \ourdesign{}.
While the results on Neyman regret in Section~\ref{sec:our-design} were stated in terms of finite-sample analysis, the results in this section will be solely focused on asymptotic analysis using the triangular array asymptotics discussed in Section~\ref{sec:technical-assumptions}.

\subsection{Non-Superefficiency} \label{sec:non-superefficiency}

In design-based inference, it is necessary to place assumptions on the potential outcomes which ensure that the variance of an effect estimator is not superefficient, i.e. the variance does not decay faster than the usual $T^{-1}$ rate \citep{Aronow2017Estimating, Leung2022Causal, Harshaw2022Design}.
For example, if the potential outcomes are perfectly correlated so that the AIPW estimator is constant, then the variance can be equal to zero.
When potential outcomes are presumed to be sampled i.i.d. from a superpopulation, such an event would occur with probability zero; however, it has so far been permissible under our assumptions.
In particular, note that the convergence of the Neyman regret ensures only an upper bound on the variance, not a lower bound.

Throughout most of the design-based literature which focuses on complex designs, non-superefficiency is guaranteed by direct assumption{\textemdash}that is, by assuming that $\liminf_{T \to \infty} T \cdot \Var{\eate} > 0$.
In this work, we are able to identify precise conditions when the adaptive AIPW estimator will be non-superefficient under \ourdesign{}.
Recall that the oracle variance is given by $T \cdot \optvar = 2(1+\rho) \olsres{1} \olsres{0}$, where $\olsres{1}$ and $\olsres{0}$ are the minimal residuals in each treatment condition and $\rho$ is the correlation of residuals between the conditions.
The following theorem shows that the asymptotic variance of the adaptive estimator is exactly the oracle variance, provided that we strengthen Assumption~\ref{assumption:maximum-radius} from $R \lesssim T^{1/2}$ to $R = \littleO{T^{1/2}}$.

\begin{assumption}[Maximum Radius Strengthened]
	\label{assumption:maximum-radius-strengthened}
	The maximum radius is bounded as $R = o(T^{1/2})$.
\end{assumption}

\expcommand{\asymptoticvariance}{
	Under Assumptions \mainref{assumption:moments}-\mainref{assumption:maximum-radius-strengthened} and Condition \mainref{condition:sigmoid}, the asymptotic variance of the adaptive AIPW estimator under \ourdesign{} is the oracle variance:
	\[
	T \cdot \Var{\eate} = 2 (1 + \rho) \mathcal{E}(1) \mathcal{E}(0) + \littleO{1}
	\enspace.
	\]
}

\begin{theorem} \label{thm:exact-asymptotic-variance}
	\asymptoticvariance
\end{theorem}

While Theorem~\ref{thm:exact-asymptotic-variance} derives the exact asymptotic variance of the adaptive AIPW estimator under \ourdesign{}, it does not imply the Neyman regret guarantee of Theorem~\ref{thm:neyman-regret}.
Theorem~\ref{thm:neyman-regret} shows that the difference between the (normalized) adaptive and oracle variances goes to zero at the rate $T^{-1/2} R$.
On the other hand, Theorem~\ref{thm:exact-asymptotic-variance} establishes the exact asymptotic variance, albeit at slower rates of convergence hidden in the $\littleO{1}$ term.
In this sense, Theorems~\ref{thm:neyman-regret} and \ref{thm:exact-asymptotic-variance} are complementary but ultimately incomparable.

Theorem~\ref{thm:exact-asymptotic-variance} shows that the exact asymptotic variance depends only on the optimal residuals $\olsres{1}$ and $\olsres{0}$ as well as the residual correlation $\rho$.
Assumption~\mainref{assumption:moments} ensures that these $\olsres{1}$ and $\olsres{0}$ are bounded from below.
The following assumption ensures that the correlation is bounded away from $-1$, which implies that the estimator is not superefficient.

\begin{assumption}[Bounded Correlation] \label{assumption:bounded-correlation}
	$\liminf_{T \to \infty} \rho > -1$.
\end{assumption}

\expcommand{\nonsuperefficiency}{
	Under Assumptions~\mainref{assumption:moments}-\mainref{assumption:bounded-correlation} and Condition \mainref{condition:sigmoid}, $\liminf_{T \to \infty} T \cdot \Var{\eate} > 0$.
}

\begin{corollary} \label{corollary:non-superefficiency}
	\nonsuperefficiency
\end{corollary}

\subsection{Central Limit Theorem} \label{sec:clt}

To facilitate the development of asymptotically valid inference, we derive a central limit theorem for the adaptive AIPW estimator under \ourdesign{}.
The Central Limit theorem holds under the same conditions required for non-superefficiency:

\expcommand{\ourclt}{
Under Assumptions~\mainref{assumption:moments}-\mainref{assumption:bounded-correlation} and Condition \mainref{condition:sigmoid}, the estimator is asymptotically normal under \ourdesign{}: $(\eate - \ate) / \sqrt{\Var{\eate}} \xrightarrow{d} \mathcal{N}(0,1)$.
}

\begin{theorem} \label{thm:clt}
	\ourclt
\end{theorem}

Because the estimator has a martingale structure with respect to iterations in the adaptive design, we will use a standard martingale CLT.
While the martingale CLT itself is standard, verifying the necessary conditions will require new technical developments.

\expcommand{\martingaleclt}{
	If $\{X_{t,T}\}$ is a triangular array of martingale difference sequences with respect to filtrations $\filt_{t,T}$, i.e. $\E{ X_{t,T} \mid \filt_{t-1,T} } = 0$, then if
	\begin{enumerate}
		\item Conditional Lyapunov Condition: $\exists$ $\delta > 0$ such that $\sum_{t=1}^T \E{ |X_{t,T}|^{2+\delta} \mid \filt_{t-1,T} } \xrightarrow{p} 0$
		\item Conditional Variance Convergence: $V_T^2 := \sum_{t=1}^T \E{ X_{t,T}^2 \mid \filt_{t-1,T} } \xrightarrow{p} 1$
	\end{enumerate}
	then the sum $S_T = \sum_{t=1}^T X_{t,T}$ converges to a standard normal in distribution, $S_T \xrightarrow{d} \mathcal{N}(0,1)$.
}

\begin{lemma}[\cite{helland1982central}] \label{lemma:martingale-clt}
	\martingaleclt
\end{lemma}

In order to use the martingale CLT, we will consider the martingale difference sequence given by
\[
X_{t,T} = \frac{ \eate_t - \ate_t }{T \sqrt{ \Var{\eate} } }
\enspace,
\]
where $\ate_t = y_t(1) - y_t(0)$ is the individual treatment effect and 
\[
\eate_t = \iprod{\xv_t,\bv_t(1)} - \iprod{\xv_t,\bv_t(0)} +  \paren[\big]{ Y_t - \iprod{ \xv_t, \bv_t(1) } }  \frac{ \indicator{Z_t=1} }{p_t}
- \paren[\big]{ Y_t - \iprod{ \xv_t, \bv_t(0) } }  \frac{ \indicator{Z_t=0} }{1-p_t}
\]
is the individual effect estimator.
One can readily verify that $\sum_{t=1}^T X_{t,T} = (\eate - \ate) / \sqrt{\Var{\eate}}$, so that it remains to argue for the two conditions in Lemma~\ref{lemma:martingale-clt}: the conditional Lyapunov condition and the conditional variance convergence.
Due to space considerations, we defer the full proof to Section \suppref{section:D3} in the appendix.
Instead, we highlight two technical developments in the proof which involve controlling the inverse probabilities.

\paragraph{Almost Sure Bounds on Inverse Probabilities}
In the proof of the central limit theorem, we bound the inverse probabilities almost surely.
Because Assumption~\mainref{assumption:moments} bounds only the fourth moments of outcomes, an almost sure bound on the inverse probabilities is necessary to verify the conditional Lyapunov condition.
The following proposition contains our almost sure bound:

\expcommand{\asinvprobbound}{
	Under Assumptions \mainref{assumption:moments}-\mainref{assumption:maximum-radius-strengthened} and Condition \mainref{condition:sigmoid}, there exists a constant $K>0$ so that
	\[
	\Pr[\Big]{ \max \braces[\Big]{ \frac{1}{p_t}, \frac{1}{1-p_t} } \leq K \cdot T^{7/26} R_t^{-2/11} \ \text{ for all } t \in [T]} = 1 \enspace.
	\]
}

\begin{proposition} \label{prop:as-inv-prob-bound}
	\asinvprobbound
\end{proposition}

While most adaptive designs obtain bounded inverse probabilities through explicit truncation, Proposition~\ref{prop:as-inv-prob-bound} demonstrates that such bounds are actually a \emph{consequence of the sigmoidal regularizer} used in \ourdesign{}.
To the best of our knowledge, both this type of result and its proof are unique in this literature and may be of independent interest.
For this reason, we sketch the key arguments and refer readers to supplement Section~\suppref{section:D2} for the full proof.

The proof is based on applying successive refinements of an inductive argument.
Roughly speaking, each inductive argument proceeds in the following manner: if $\max_{s \leq t} \maxf{1/p_s, 1/(1-p_s)}$ is bounded by $\mu$, then the estimated online residuals $\estores{t}{1}$ and $\estores{t}{0}$ are also suitably bounded, which further implies that $\maxf{1/p_{t+1}, 1/(1-p_{t+1})} \leq \mu$ is bounded by the probability regularization lemma (Lemma~\ref{lemma:effect-of-p-regularization}).
We remark that the curvature conditions on the sigmoidal regularizer are crucial for this argument.
Applying this inductive argument once would yield an almost sure bound of order $\mu = \bigO{T^{1/2}}$, which is too weak to be useful in the CLT.
So, we refine this uniform bound by recursively applying the inductive argument above to obtain a sequence of refined uniform bounds $\setb{\mu_r}_{r=1}^\infty$ for which $\mu_r \to \mu^* = \bigO{T^{7/26}}$.
Further time-dependent improvements can now be made by using this $\mu^*$ bound in the inductive argument and taking into account the time-dependent step size $\eta_t = T^{-1/2} R_t^{-1}$.
This final refinement yields the bound appearing in Proposition~\ref{prop:as-inv-prob-bound}.

\paragraph{Stability of Inverse Probabilities}
In proving the conditional variance convergence in Lemma~\ref{lemma:martingale-clt}, we need to go beyond almost sure bounds and establish the stability of the inverse probabilities; stability of the linear predictors is also required, but this is established through the prediction tracking technique described in Section~\ref{sec:prediction-tracking}.
For this purpose, we introduce the \emph{stabilized probability} $\widebar{p}_t$, which we define as
\begin{equation} \label{eq:stabilized-prob}
\widebar{p}_t = \argmin_{p \in (0,1)} 
\frac{ \E{\estores{t-1}{1}} }{p} 
+ \frac{ \E{\estores{t-1}{0}} }{1-p}
+ \eta_t^{-1} \Psi(p)
\enspace.
\end{equation}
The stabilized probability $\widebar{p}_t$ is a deterministic quantity which is defined similarly to the random adaptively chosen probability $p_t$, except that the stabilized probability $\widebar{p}_t$ replaces the estimated squared residuals with their expected values.

In the proof of the central limit theorem, we investigate the expected absolute differences between the inverse probabilities and the inverse stabilized probabilities: $\E{\abs{ 1/p_t - 1/\widebar{p}_t }}$ and $\E{\abs{ 1/(1-p_t) - 1/(1-\widebar{p}_t) }}$.
Through the first order conditions of \eqref{eq:stabilized-prob}, these expected differences of inverse stabilized probabilities can be expressed as differences in the estimated online residuals, i.e. $\estores{t-1}{1} - \E{ \estores{t-1}{1} }$ and $\estores{t-1}{0} - \E{ \estores{t-1}{0} }$.
The subtle aspect is that the variance $\Var{ \estores{t-1}{1} }$ depends on the expectation $\E{ \estores{t-1}{0}}$ in the following way: when $\E{ \estores{t-1}{0}}$ is large, then $1/p_t$ is large and thus $\Var{ \estores{t-1}{1} }$ is large due to the inverse probability weighting.
Likewise, the variance $\Var{ \estores{t-1}{0} }$ depends on the expectation $\E{ \estores{t-1}{1}}$.
In essence, the terms $\estores{t-1}{1}$ and $\estores{t-1}{0}$ have a \emph{mutually normalizing} property.
This mutually normalizing property is used to prove several delicate aspects of the central limit theorem, including the following proposition:

\begin{proposition} \label{prop:interesting-variance-stabilizing-term}
Under Assumptions~\mainref{assumption:moments}-\mainref{assumption:maximum-radius-strengthened} and Condition \mainref{condition:sigmoid},
\[
\frac{1}{T} \sum_{t=1}^T \paren[\big]{ y_t(1) - \iprod{ \xv_t, \optbv_t(1) } }^2 \E[\Big]{ \abs[\Big]{ \frac{1}{p_t} - \frac{1}{\widebar{p}_t} } } 
\to 0 
\enspace.
\]
\end{proposition}

Proving Proposition~\ref{prop:interesting-variance-stabilizing-term} is challenging because the residuals and the difference of the inverse probabilities cannot be analyzed separately, e.g. using H{\"o}lder's inequality.
The issue is that under \ourdesign{}, a small proportion of the inverse probabilities $1/p_t$ could be large; in fact, some inverse probability terms $1/p_t$ will not concentrate around the inverse stabilized probability $1/\widebar{p}_t$.
This means that some of the terms $\E{ \abs{ 1/p_t - 1/\widebar{p}_t } }$ can be large.
The saving grace is that when the difference of these probabilities is large, the corresponding residuals must be small, so that the product in the sum is controlled.
The formal proof of Proposition~\ref{prop:interesting-variance-stabilizing-term} relies heavily on the mutually normalizing property described above.
Given its central role in the proof, we believe that the mutually normalizing property might be of independent interest.

\subsection{Variance Estimation and Confidence Intervals} \label{sec:var-estimation}

In order to construct confidence intervals, we first need to construct an estimator for the oracle variance $T \cdot \optvar = 2 (1+\rho) \olsres{1} \olsres{0}$.
However, this oracle variance cannot be consistently estimated because the residual correlation $\rho$ does not admit a consistent estimator.
The impossibility of consistent variance estimation is a well-known problem in the design-based framework \citep[see e.g.][]{ImbensRubin2015, Harshaw2021Optimized}.
In light of this issue, experimenters tend to opt for conservative (i.e. upwardly biased) estimates of the variance, which in turn yield conservative inferential procedures.
Conservative variance estimates are typically constructed by first identifying an estimable upper bound on the variance, or \emph{variance bound} for short.

The most common variance bound is due to \citet{Neyman1923} and was originally derived in the setting of a (non-adaptive) completely randomized experiment.
We proceed by adapting his variance bound to the current setting.
Because the correlation is at most $\rho \leq 1$, we have that
\[
T \cdot \optvar = 2(1+\rho) \olsres{1} \olsres{0} \leq 4 \olsres{1} \olsres{0} := T \cdot \vb
\enspace,
\]
which is referred to as the Neyman variance bound.
Unlike the variance itself, the Neyman variance bound is estimable because it depends only on the optimal residuals $\olsres{1}$ and $\olsres{0}$, each of which may be estimated to high precision.
Interestingly, selecting the Neyman allocation $p^*$, $\bv^*(1)$, and $\bv^*(0)$ which minimizes the variance and then applying Neyman's bound is equivalent to first applying Neyman's bound and then selecting the allocation to minimize the bound.
In other words, by constructing an adaptive design which minimizes the variance of the point estimator, we are also minimizing the expected width of the resulting confidence interval.

In order to construct a consistent estimator for the variance bound $\vb$, we focus first on constructing consistent estimates of the squared residuals.
For each treatment $k \in \setb{0,1}$, the squared residual is given by
\[
\olsres{k}^2 = \min_{\bv \in \Reals^d} \frac{1}{T} \sum_{t=1}^T \paren[\big]{ y_t(k) - \iprod{ \xv_t , \bv } }^2
= \frac{1}{T} \vec{y}(k)^\tran \mat{Q} \vec{y}(k) \enspace,
\]
where $\mat{Q}$ is the orthogonal projection matrix $I - \xM (\xM^\tran \xM)^{-1} \xM^\tran$, $\xM$ is the $T$-by-$d$ matrix whose rows are covariate vectors, and $\vec{y}(k)$ is the vector of potential outcomes $\vec{y}(k) = \paren{ y_1(k) \dots y_T(k) }^{\tran}$.
To simplify the notation, we introduce $\sqolsres{k} = \olsres{k}^2$.
In order to estimate the optimal squared residuals, we use an adaptive IPW estimator applied to the quadratic form,
\begin{align*}
\estsqolsres{1}
&= \frac{1}{T} \sum_{t=1}^T Q_{t,t} Y_t^2 \frac{ \indicator{Z_t =1} }{p_t} 
+ \frac{1}{T} \sum_{t=1}^T \sum_{s \neq t} Q_{t,s} Y_t Y_s \frac{ \indicator{Z_t = 1, Z_s = 1} }{p_tp_s} 
\quadand \\
\estsqolsres{0}
&=  \frac{1}{T} \sum_{t=1}^T Q_{t,t} Y_t^2 \frac{ \indicator{Z_t = 0} }{1-p_t} 
+ \frac{1}{T} \sum_{t=1}^T \sum_{s \neq t} Q_{t,s} Y_t Y_s \frac{ \indicator{Z_t = 0, Z_s = 0} }{(1-p_t)(1-p_s)}
\enspace.
\end{align*}
The adaptive inverse probability weighting is chosen so that $\estsqolsres{k}$ is an unbiased estimate of $\sqolsres{k}$.
Even though this estimate $\estsqolsres{k}$ is unbiased for the positive value $\sqolsres{k}$, there is a small chance that it will take negative values.
To account for this possibility, we define the estimated OLS residuals as $\estolsres{k} = \maxf{  \estsqolsres{k}, 0}^{1/2}$.
However, the following theorem shows that $\estsqolsres{k}$ is a consistent estimator, and so it takes negative values with vanishingly small probability.

\expcommand{\sqolsresestimates}{
	Under Assumptions~\mainref{assumption:moments}-\mainref{assumption:maximum-radius-strengthened} and Condition \mainref{condition:sigmoid}, for each treatment $k \in \setb{0,1}$, we have that the estimated squared OLS residuals satisfy
	\[
	\E{ \estsqolsres{k} } = \sqolsres{k}
	\quadand 
	\Var{ \estsqolsres{k} } = \bigO{ \braces{ T^{-5/12} R^{5/6} }^2 }     
	\enspace.
	\]
}

\begin{theorem} \label{thm:sq-ols-res-estimates}
	\sqolsresestimates
\end{theorem}

Theorem~\ref{thm:sq-ols-res-estimates} shows that $\estolsres{k}$ is consistent at a rate $T^{-5/12} R^{5/6}$, which depends on both the number of samples $T$ and the maximum radius of the covariates $R$.
Under Assumption~\mainref{assumption:maximum-radius-strengthened}, we have that $R = \littleO{T^{1/2}}$ so that consistency is ensured, i.e. $\estolsres{k} - \olsres{k} = \littleOp{1}$.
In the regime of uniformly bounded covariates where $R \asymp d^{1/2}$, Theorem~\ref{thm:sq-ols-res-estimates} demonstrates that the estimator $\estolsres{k}$ converges at a rate of $\bigOp{(d/T)^{5/12}}$.

It's interesting to note that \ourdesign{} uses adaptive experimentation to ensure maximal efficiency of the point estimator, but this seems to come at the cost of slightly worsened rates of consistency of the variance bound estimator.
We find this trade-off to be unproblematic because efficiency of the point estimator is the first-order concern for confidence intervals.
In contrast, the estimated variance bound is required only to be consistent, and the precise rate is inconsequential for asymptotic coverage.
It is an interesting question whether this trade-off is necessary, or merely an artifact of our approach.

The proof of Theorem~\ref{thm:sq-ols-res-estimates} proceeds by analyzing the covariance of the individual terms in the quadratic estimator.
This requires specialized case analysis of subject pairs $(s,t)$ and $(k,\ell)$.
For example, if the subject pairs are distinct, i.e. $\{s,t\} \cap \{k,\ell\} = \emptyset$, then the corresponding covariance term is zero.
We refer the reader to Section \suppref{section:D5} of the appendix for the full case analysis.

Given the estimates of the OLS residuals, we construct an estimator of the Neyman variance bound as $T \cdot \evb = 4 \estolsres{1} \estolsres{0}$.
The estimated variance bound inherits the same rates of consistency:

\expcommand{\vbconsistent}{
	Under Assumptions~\mainref{assumption:moments}-\mainref{assumption:maximum-radius-strengthened} and Condition \mainref{condition:sigmoid},
	$T \cdot \evb - T \cdot \vb = \bigOp{T^{-5/12} R^{5/6}}$.
}

\begin{corollary} \label{corollary:vb-is-consistent}
	\vbconsistent
\end{corollary}

The estimator of the variance bound, together with the central limit theorem, yields the following Wald-type confidence interval: for a given level $\alpha \in (0,1)$, the confidence interval $\ci$ is
\[
\ci = \eate \pm \Phi^{-1}(1 - \alpha/2) \cdot \sqrt{ \evb } \enspace,
\]
where $\Phi^{-1} : (0,1) \to \Reals$ is the quantile function of a standard normal.

\expcommand{\wald}{
	Under Assumptions~\mainref{assumption:moments}-\mainref{assumption:bounded-correlation} and Condition \mainref{condition:sigmoid}, the Wald-type intervals cover at least at the nominal level: $\liminf_{T \to \infty} \Pr{ \ate \in \ci } \geq 1 - \alpha$.
}
\begin{corollary} \label{corollary:wald-coverage}
	\wald
\end{corollary}
	
	\section{Numerical Simulations}\label{sec:simulations}

We assess the empirical performance of \ourdesign{} using data collected in the (non-adaptive) field experiment of \citet{groh2016macroinsurance}.
The original study examined whether providing protection against macroeconomic risks could benefit microenterprise owners in Egypt.
It included 2961 clients of the country's largest microfinance institution and evaluated a newly introduced insurance product on a collection of outcomes.
Our objective is not to revisit the findings of the original study, but rather to use this experimentally collected data to investigate the behavior of our adaptive experimental design.

Due to the fundamental problem of causal inference, the dataset does not contain both potential outcomes for any of the measured outcomes.
One way to overcome this issue is to impute missing outcomes via some generative working model.
Rather than taking this approach, we use two different outcomes which were actually observed in the experiment to act as the two potential outcomes in our simulation.
Specifically, average monthly business profit is taken as the potential outcome under treatment for every subject, irrespective of the subject's actual treatment assignment in the original study.
Similarly, average monthly business revenue is taken as the potential outcome under control.
This approach to selecting the two potential outcomes allows us to preserve the irregularities and idosyncrasies of real data as much as possible.
However, the causal effect induced by this construction does not necessarily yield a substantive interpretation.

We include all available baseline variables in the original dataset as covariates in our simulations, including an intercept term.
The covariates are standardized to have mean zero and unit standard deviation.
We drop ``profits2'' and ``missprofits2'', since these two variables are highly collinear with ``profits1'' and ``missprofits1'', respectively.
After removing observations with missing data, the resulting sample contains 2,906 subjects with 12 covariates.
We further rescale the two potential-outcome vectors jointly to normalize the two residual sums of squares so that $\min\{\olsres{1},\olsres{0}\}=1$.

When forming the sequential experiment, we retain the ordering of the subjects appearing in the original dataset.
To increase the sample size, we select larger subsequences of the entire data sequence.
We set the \ourdesign{} step size as $\eta_t=c\cdot T^{-1/2}R_t^{-1}$ with $c=10$.
In Section~\suppref{sec:supp-simulations} of the supplement, we run additional simulations which demonstrate that performance is not sensitive to the choice of $c$.
All statistics presented here are based on 100,000 Monte Carlo simulations.

\begin{figure}[t]
	\centering  
	\begin{subfigure}{0.45\textwidth}
		\centering
		\includegraphics[width=\linewidth]{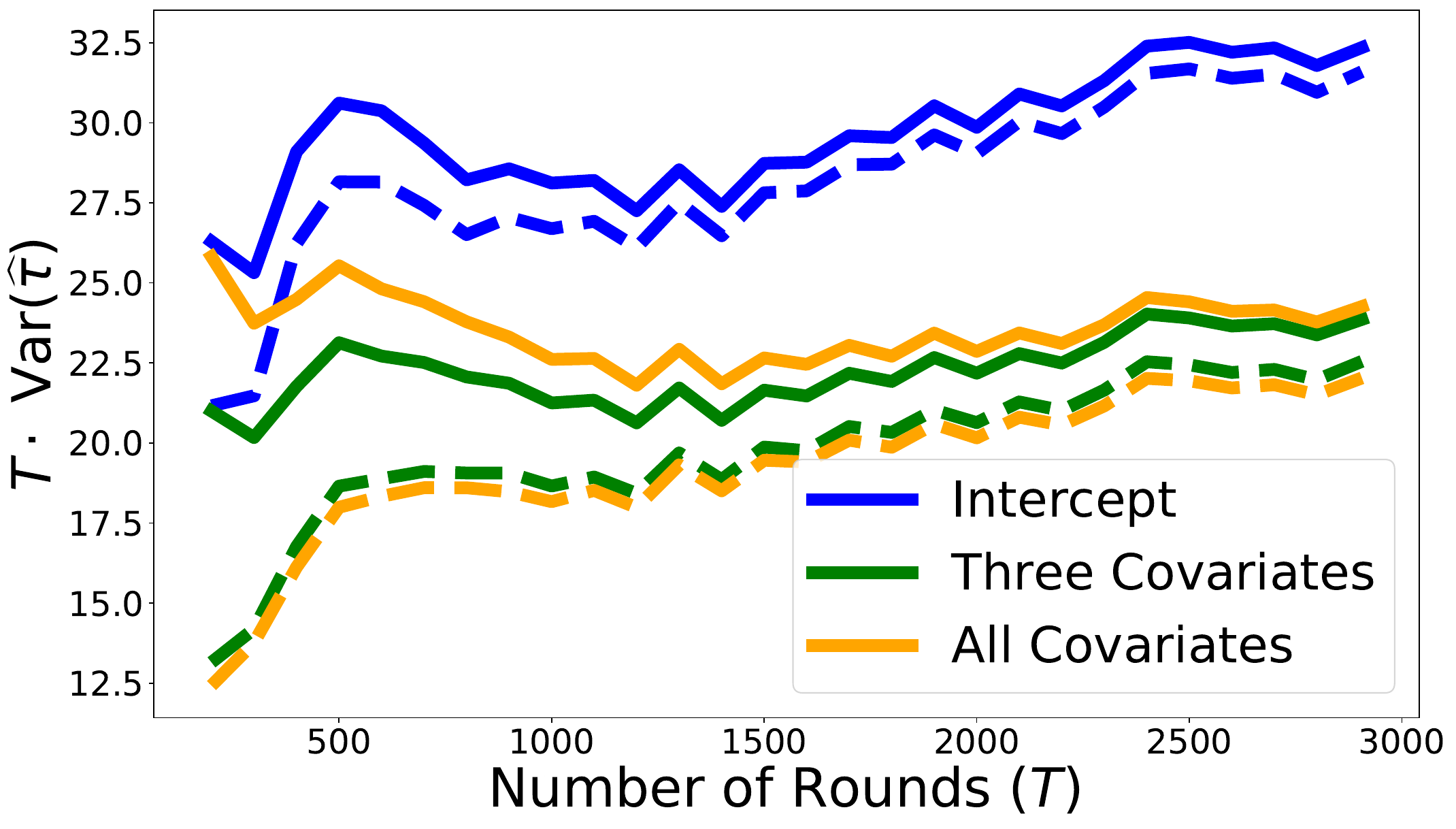}
		\caption{Normalized Variance}
		\label{fig:comparison-our-design-var}
	\end{subfigure}
	\hfill
	\begin{subfigure}{0.45\textwidth}
		\centering
		\includegraphics[width=\linewidth]{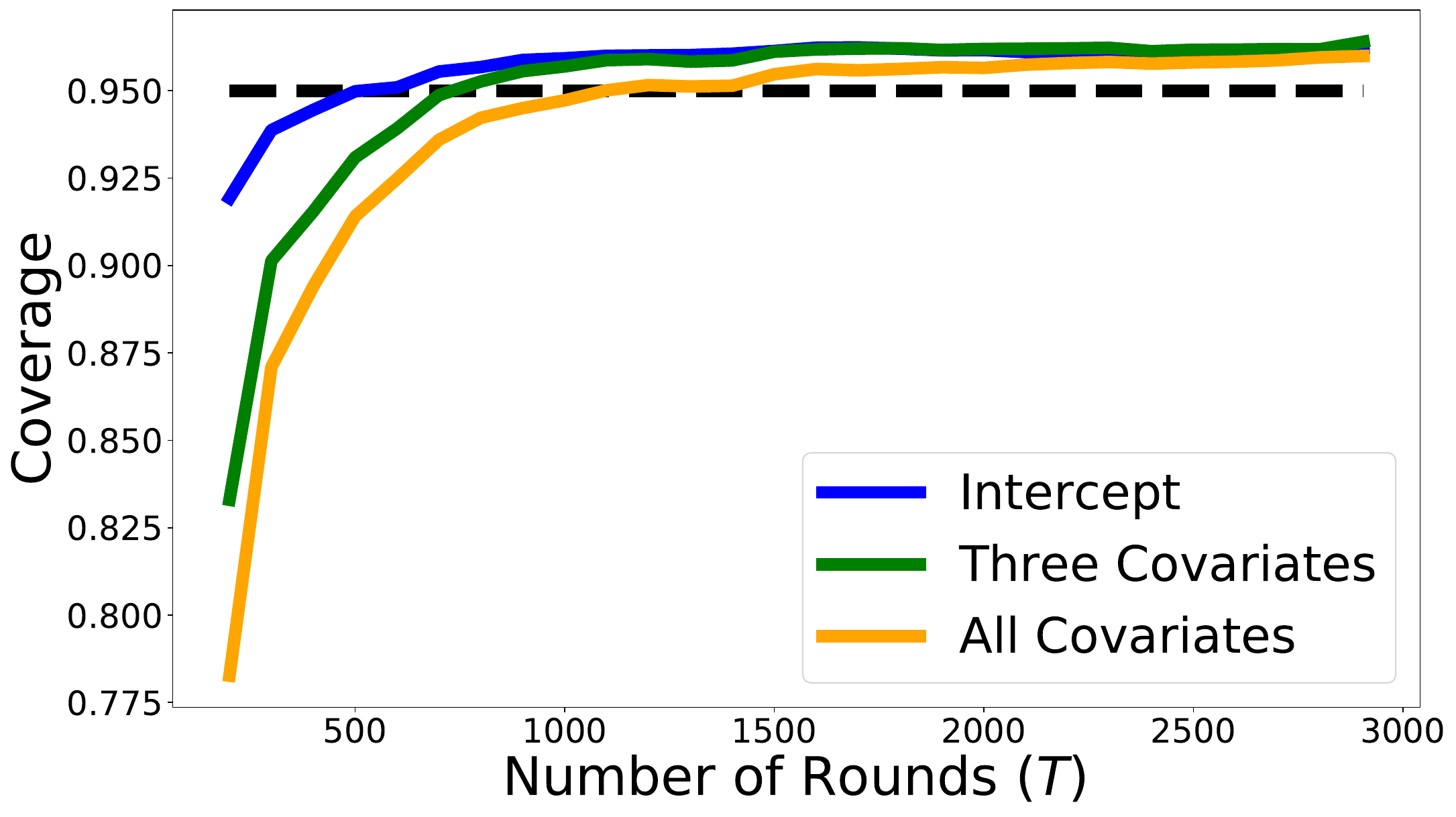}
		\caption{Coverage of Intervals}
		\label{fig:comparison-our-design-coverage}
	\end{subfigure}
	\caption{
	\ourdesign{} was run with three different choices of covariates: (1) intercept only (2) three most informative covariates (3) all covariates.
	The effect of these choices on the (normalized) variance of the estimator and coverage of intervals are shown in Figures~\ref{fig:comparison-our-design-var} and \ref{fig:comparison-our-design-coverage}, respectively. 
	}
	\label{fig:covarite-comparison}
\end{figure}

We begin by empirically investigating the trade-off inherent in adding more covariates: as additional covariates are added, the targeted oracle variance grows smaller, but the convergence of the Neyman Regret is slower.
Figure~\ref{fig:covarite-comparison} demonstrates the effect of this trade-off by comparing the performance of \ourdesign{} across three choices of covariates: (1) the intercept alone, (2) the three most informative covariates (``expectloslessthan3000'', ``gender'', ``profits1''), and (3)  all 12 covariates.
Figure~\ref{fig:comparison-our-design-var} plots the adaptive and oracle variances corresponding to these three choices of covariates.
The solid lines correspond to the actual adaptive variance, the dotted lines correspond to the oracle Neyman variance, and thus the difference between the two is the Neyman Regret.
For each choice of covariates, we see that the Neyman regret is indeed shrinking as the sample size grows, i.e. the solid and dotted lines are getting closer.
In addition, we see that adding more covariates decreases the oracle Neyman variance, although little oracle precision is gained by adding the entire set of covariates in addition to the three most informative ones.
The most interesting aspect of Figure~\ref{fig:comparison-our-design-var} is the illustration of the aforementioned trade-off: \ourdesign{} with only the three most informative covariates achieves a lower variance than \ourdesign{} with all covariates, at least for the sample sizes we consider here.
It is worth noting, however, that this difference is less pronounced as the sample size grows: indeed, the performance is almost identical at sample sizes above $T \geq 2,500$.

Figure~\ref{fig:comparison-our-design-coverage} contains the finite-sample coverage of our proposed confidence intervals under these three choices of covariates.
While all intervals eventually cover at the nominal rate, the sample size at which this happens grows as more covariates are added.
Regardless of how many covariates are included, the proposed intervals have valid coverage for sample sizes larger than $T = 1,000$.
At first glance, this might seem to contradict the results in Figure~\ref{fig:comparison-our-design-var}: how can two different covariate choices both cover at a sample size of $T = 1,000$ when their variances are not comparable until larger sample sizes of $T \geq 2,500$?
The answer is that the confidence intervals are not constructed using a consistent variance estimator, but instead one which is conservatively consistent, i.e. converges to an upper bound.
Roughly speaking, coverage occurs once the estimated oracle variance bound is larger than the actual variance.
In particular, this can occur at sample sizes which are smaller than what is required for the actual variances to become comparable.
This numerical simulation highlights that the considerations for point estimation and inference might differ in this way.

The practical recommendation from these simulations is that researchers ought to carefully choose covariates which they believe will be most informative, rather than simply using all available covariates.
However, this consideration is most important for smaller sample sizes.

\begin{figure}[t]
	\centering  
	\begin{subfigure}{0.32\textwidth}
		\centering
		\includegraphics[width=\linewidth]{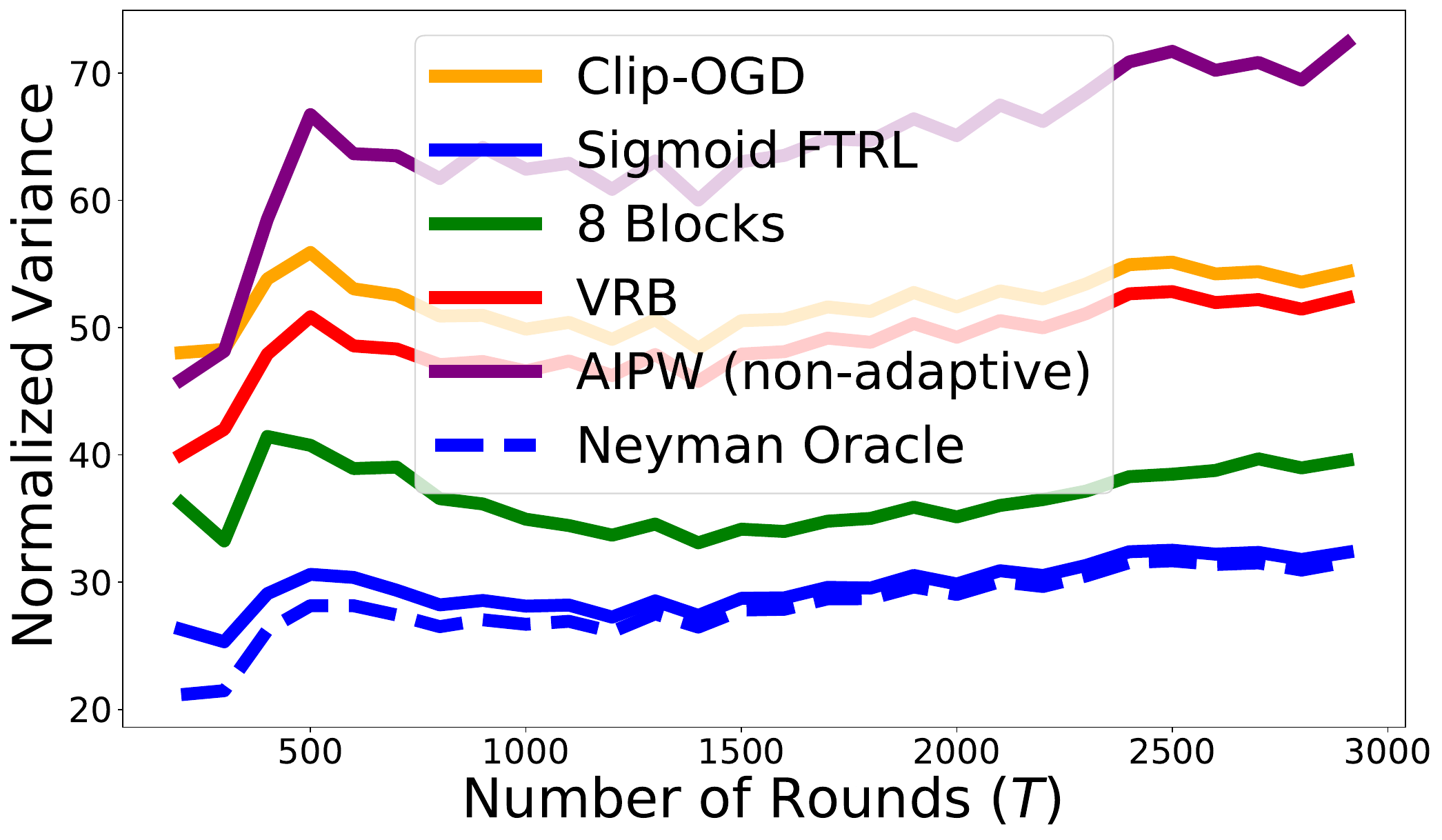}
		\caption{Intercept only}\label{fig:intercept_variance}
	\end{subfigure}
	\hfill
	\begin{subfigure}{0.32\textwidth}
		\centering
		\includegraphics[width=\linewidth]{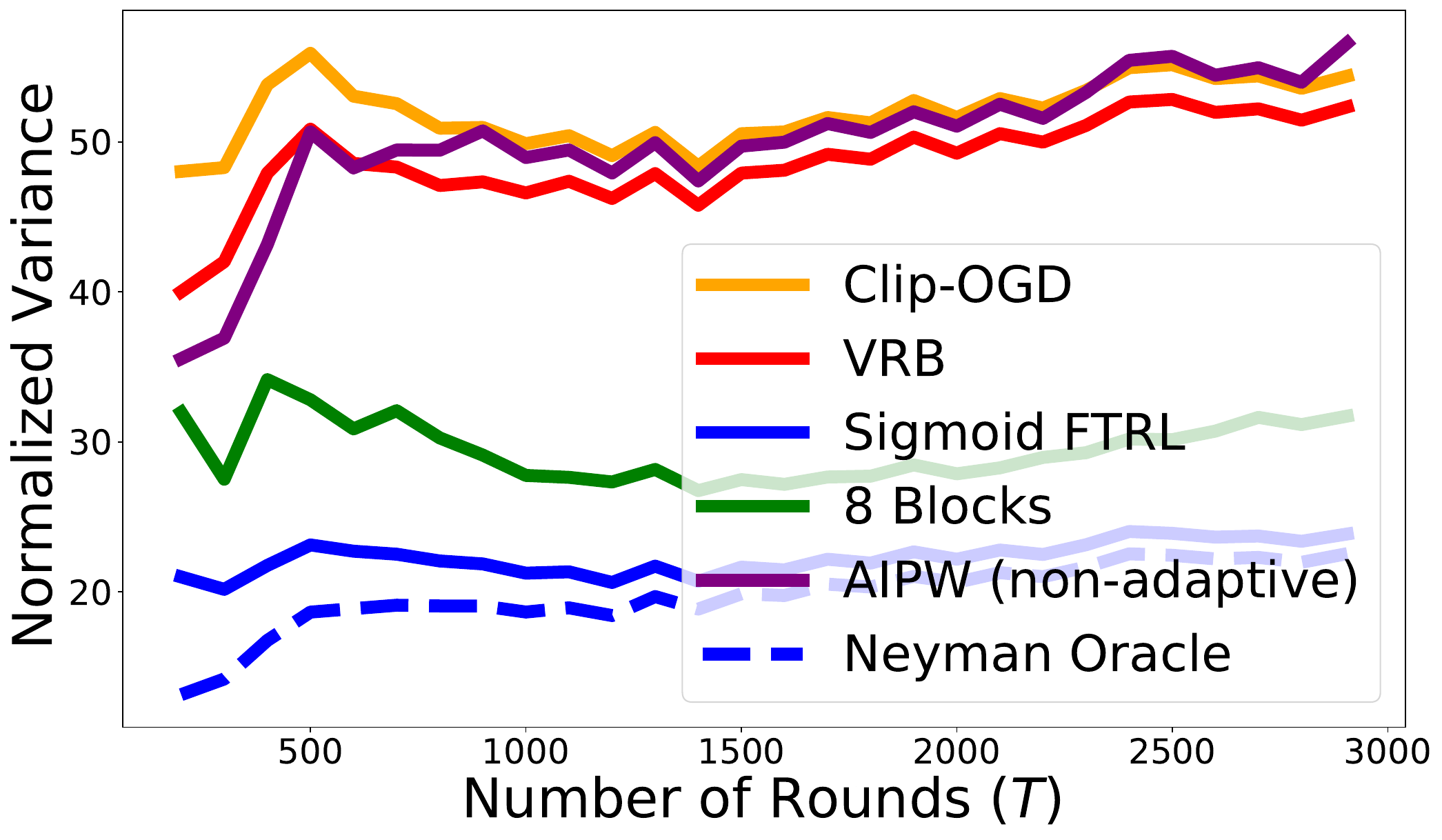}
		\caption{Three covariates}\label{fig:three_covariates_variance}
	\end{subfigure}
	\hfill
	\begin{subfigure}{0.32\textwidth}
		\centering
		\includegraphics[width=\linewidth]{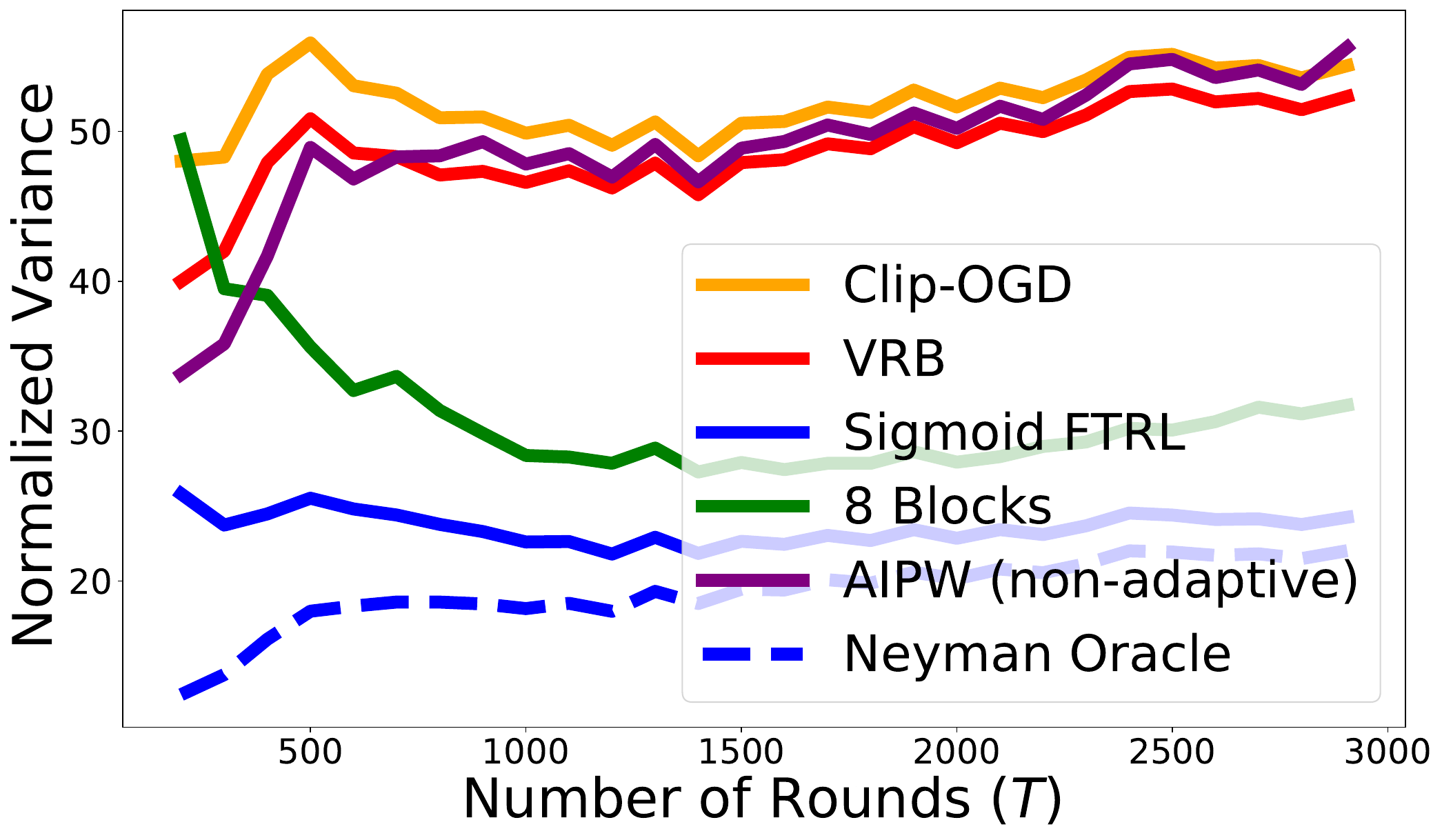}
		\caption{All covariates}\label{fig:all_covariates_variance}
	\end{subfigure}
	\caption{
		We compare the performance of \ourdesign{} against four benchmark designs across the three choices of covariates.
		Each figure shows the normalized variance of each method plotted as sample size increases.
		Figure~\ref{fig:intercept_variance} shows the intercept alone, Figure~\ref{fig:three_covariates_variance} shows three most informative covariates, and Figure~\ref{fig:all_covariates_variance} shows all the covariates.}
	\label{fig:variance-comparison}
\end{figure}

Next, we compare the performance of \ourdesign{} relative to different benchmark designs across the three different choices of covariates described above.
We include the non-adaptive AIPW estimator with $p=1/2$ \citep{Lin2013Agnostic}.
This serves as a baseline against which we can compare the gains of adaptivity of \ourdesign{}.
Next, we include adaptive designs that target the IPW estimator, including \textsc{Clip-OGD} of \citet{Dai2023Clip} and \textsc{Variance-Reduced-Bandit} of \citet{Boros2018Online}.
These serve as baselines against which we can compare the gains of incorporating covariate information, i.e. targeting AIPW rather than IPW.
We also include a heuristic ``block'' design that selects the probabilities $p_t$ and predictors $\bv_t(1)$, $\bv_t(0)$ in 8 block updates via un-regularized minimization of the loss functions we propose.

Figure~\ref{fig:variance-comparison} contains the results of this comparison where the $T$-normalized variance is plotted for each of these methods as the sample size increases.
We also plot the oracle Neyman Variance as a relevant benchmark for the performance of \ourdesign{}.
We see that \ourdesign{} is converging to the oracle Neyman Variance in these settings, although the convergence is slower when more covariates are included.
We also see that \ourdesign{} design attains the smallest variance at all sample sizes in each setting.
This is to be expected because unlike the methods that either target the unadjusted IPW estimator or implement covariate adjustment, \ourdesign{} is doing both simultaneously.
While the block design achieves the second-smallest variance, it does not converge to the oracle Neyman Variance because only a few block updates are used.
These comparisons highlight that in order to improve precision, experimenters ought to use adaptive designs which leverage adaptivity in the treatment assignment and covariate adjustment, when possible.
	
	\section{Conclusion} \label{sec:conclusion}

In this paper, we have investigated the problem of Adaptive Neyman Allocation for AIPW estimators in the design-based framework.
We have presented \ourdesign{}, an adaptive experimental design which selects treatment probabilities and linear predictors together.
Assuming mild regularity conditions, the Neyman Regret under \ourdesign{} converges at a rate of $\sqrt{d/T}$, which constitutes an optimal $T$ dependence under the regularity conditions we consider.
A central limit theorem and a variance estimator facilitate the construction of asymptotically valid Wald-type confidence intervals.

There are several avenues for future work.
One direction is to study adaptive designs where the probability of treatment can depend on the observed covariates, which would likely yield further precision improvements.
It is also an interesting direction to explore anytime valid confidence sequences in conjunction with Adaptive Neyman Allocation in the design-based setting.
For example, it would be interesting to understand the extent to which minimizing the width of the Wald-type interval (i.e. the central goal of Adaptive Neyman Allocation) and minimizing the width of the anytime valid intervals may be similar or conflicting goals.
	
	
	\printbibliography
	
	\newpage
	
	\newgeometry{letterpaper,margin=1in,bottom=1in}
	
	\appendix
	
	\addcontentsline{toc}{section}{Appendix} 
	\part{Appendix} 
	\parttoc 
	\newpage
	
	\section{Preliminary Results}\label{sec:supp-preliminary}

In this section, we prove preliminary results regarding the adaptive AIPW estimator and the sigmoidal transformation.
In Section \ref{section:A1}, we verify the unbiasedness of the adaptive AIPW estimator and provide an explicit form for its variance under an arbitrary adaptive design.
In Section \ref{section:A2}, we provide two examples of commonly used sigmoid functions that satisfy Condition \mainref{condition:sigmoid}.
In Section \ref{section:A3}, we provide a standard decomposition lemma for the regret in an FTRL algorithm and establish an upper bound on it.

\subsection{Properties of AIPW Estimator}\label{section:A1}

Throughout our analysis, it will be convenient to decompose both the average treatment effect and the adaptive AIPW estimator into individual-level quantities.
For each subject $t \in [T]$, we define the \emph{individual treatment effect} $\ate_t = y_t(1) - y_t(0)$ and the \emph{individual treatment effect estimator} $\eate_t$ as
\[
\eate_t = \iprod{\xv_t , \bv_t(1)} - \iprod{ \xv_t , \bv_t(0) }  
+ \frac{\indicator{Z_t = 1}}{p_t} \cdot(Y_t - \iprod{\xv_t , \bv_t(1)} )
- \frac{\indicator{Z_t = 0}}{1-p_t}\cdot (Y_t - \iprod{\xv_t , \bv_t(0)} )
\enspace.
\]
The average treatment effect can be decomposed as $\ate = (1/T) \sum_{t=1}^T \ate_t$ and, likewise, the adaptive AIPW estimator can be decomposed as $\eate = (1/T) \sum_{t=1}^T \eate_t$.
We will see shortly that the individual estimator $\eate_t$ can be understood as an unbiased (but noisy) estimate of the individual effect $\ate_t$.

For each $t$, let $\filt_0 \subseteq \filt_1 \subseteq \cdots \subseteq \filt_T$ be the filtration where $\filt_t$ conditions on all information collected in the first $t$ rounds, i.e., $\filt_t=\sigma\{Z_1,\ldots,Z_t\}$.
The following lemma shows that the sequence $\setb{ \eate_t - \ate_t}_{t=1}^T$ forms a martingale difference sequence with respect to this filtration.

\begin{lemma}\label{lemma:MDS}
	If $p_t \in (0,1)$ almost surely, then for each $t \in [T]$, the individual estimator is conditionally unbiased: $\E{ \eate_t \mid \filt_{t-1}} = \ate_t$.
\end{lemma}

\begin{proof}
	Recall that $p_t = \Pr{ Z_t=1 \mid \filt_{t-1} }$ is the treatment assignment probability conditioned on the observed history $\filt_{t-1}$.
	The central insight in this proof is to observe that because $p_t \in (0,1)$ is $\filt_{t-1}$-measurable, by the law of iterated expectations, we have that
    \[\E[\Big]{ \frac{\indicator{Z_t = 1}}{p_t}  }=\E[\Big]{\E[\Big]{ \frac{\indicator{Z_t = 1}}{p_t}  |\filt_{t-1}}}=\E[\Big]{\frac{p_t}{p_t}}=1
    \]
    and likewise
	\[\E[\Big]{ \frac{\indicator{Z_t = 0}}{1-p_t}  }=\E[\Big]{\E[\Big]{ \frac{\indicator{Z_t = 0}}{1-p_t}  |\filt_{t-1}}}=\E[\Big]{\frac{1-p_t}{1-p_t}}=1\enspace.
    \]
	By definition, the assignment probability $p_t$ and the linear predictors $\bv_{t}(1)$, $\bv_{t}(0)$ are measurable with respect to the observed history $\filt_{t-1}$.
	Using this together with linearity of expectation and the above calculation, we have that
	\begin{align*}
		&\E{ \eate_t \mid \filt_{t-1} } \\
		&= 
			\iprod{\xv_t , \bv_t(1)} - \iprod{ \xv_t , \bv_t(0) }  \\
			&\quad + \E[\Big]{ \frac{\indicator{Z_t = 1}}{p_t} \mid \filt_{t-1}} \paren[\big]{ y_t(1) - \iprod{\xv_t , \bv_t(1)} }- \E[\Big]{ \frac{\indicator{Z_t = 0}}{1-p_t} \mid \filt_{t-1}}  \paren[\big]{ y_t(0) - \iprod{\xv_t , \bv_t(0)} }
		 \\
		 &= \iprod{\xv_t , \bv_t(1)} - \iprod{ \xv_t , \bv_t(0) }  
		 + \paren[\big]{ y_t(1) - \iprod{\xv_t , \bv_t(1)} }
		 - \paren[\big]{ y_t(0) - \iprod{\xv_t , \bv_t(0)} } \\
		 &= y_t(1) - y_t(0) = \ate_t
		\enspace.
		\qedhere
	\end{align*}
\end{proof}

Given the property of a martingale difference sequence, we can immediately verify the unbiasedness of the adaptive AIPW estimator for the average treatment effect.

\begin{refproposition}{\mainref{prop:aipw-unbiased}}
	\aipwbias
\end{refproposition}

\begin{proof}
    Using the law of iterated expectations and Lemma \ref{lemma:MDS}, we have
    \[
    \E{\eate} 
    = \E[\Big]{ \frac{1}{T} \sum_{t=1}^T \eate_t }
    = \E[\Big]{ \frac{1}{T} \sum_{t=1}^T \E{ \eate_t \mid \filt_{t-1} } }
    = \frac{1}{T} \sum_{t=1}^T \ate_t
    = \ate 
    \enspace.
    \qedhere
    \]
\end{proof}

We can also use the martingale property of the individual estimators to derive a simple expression for the variance of the AIPW estimator under an arbitrary design.

\begin{refproposition}{\mainref{prop:aipw-variance}}
    \aipwvar
\end{refproposition}

\begin{proof}
    By decomposing the estimator and the estimand and applying the martingale property of the individual estimators (Lemma~\ref{lemma:MDS}), we have that
	\[
	T \cdot \Var{ \eate }
	= T \cdot \E{ ( \eate - \ate )^2 }
	= T \cdot \E[\Big]{ \paren[\Big]{ \frac{1}{T} \sum_{t=1}^T (\eate_t - \ate_t) }^2 }
	= \E[\Big]{ \frac{1}{T} \sum_{t=1}^T \E[\big]{ \paren{  \eate_t - \ate_t }^2 \mid \filt_{t-1} } }
	\]
	Thus, it remains to calculate $\E{ \paren{  \eate_t - \ate_t }^2 \mid \filt_{t-1} }$.
	To this end, observe that the difference $\eate_t - \ate_t$ can be written as
	\[
	\eate_t - \ate_t = 
	\paren[\Big]{ \frac{\indicator{Z_t = 1}}{p_t} - 1 } 
	\cdot \underbrace{\paren[\Big]{ y_t(1) - \iprod{\xv_t, \bv_t(1)} }}_{:= \Delta_t(1)} 
	-
	\paren[\Big]{ \frac{\indicator{Z_t = 0}}{1-p_t} - 1 } 
	\cdot \underbrace{\paren[\Big]{ y_t(0) - \iprod{\xv_t, \bv_t(0)} }}_{:= \Delta_t(0)}
	\enspace,
	\]
	where we have introduced the shorthand $\Delta_t(k) = y_t(k) - \iprod{  \xv_t, \bv_t(k) }$.
	Because the linear predictors $\bv_t(1)$ and $\bv_t(0)$ are measurable with respect to the observed history $\filt_{t-1}$, so too are $\Delta_t(1)$ and $\Delta_t(0)$.
	Before continuing, let us present a few conditional expectations which will be useful:
	\begin{align*}
		\E[\Big]{ \paren[\Big]{ \frac{\indicator{Z_t = 1}}{p_t} - 1 }^2 \mid \filt_{t-1} } &= \frac{1}{p_t} - 1 \\
		\E[\Big]{ \paren[\Big]{ \frac{\indicator{Z_t = 0}}{1-p_t} - 1 }^2 \mid \filt_{t-1} } &= \frac{1}{1-p_t} - 1 \\
		\E[\Big]{ \paren[\Big]{ \frac{\indicator{Z_t = 1}}{p_t} - 1 } \cdot \paren[\Big]{ \frac{\indicator{Z_t = 0}}{1-p_t} - 1 }  \mid \filt_{t-1} } &= -1
	\end{align*}
	The above can be verified via direct calculation, keeping in mind that the assignment probability $p_t$ is measurable with respect to the observed history $\filt_{t-1}$.
	Putting these together, we can now calculate  $\E{ \paren{  \eate_t - \ate_t }^2 \mid \filt_{t-1} }$ as
	\begin{align*}
		 \E{ \paren{  \eate_t - \ate_t }^2 \mid \filt_{t-1} }
		 &= \E[\Bigg]{ \braces[\Bigg]{ \paren[\Big]{ \frac{\indicator{Z_t = 1}}{p_t} - 1 } \Delta_t(1)  - \paren[\Big]{ \frac{\indicator{Z_t = 0}}{1-p_t} - 1 } \Delta_t(0) }^2 \mid \filt_{t-1} } \\
		 &= \E[\Big]{ \paren[\Big]{ \frac{\indicator{Z_t = 1}}{p_t} - 1 }^2 \mid \filt_{t-1} } \cdot \Delta_t(1)^2
		 + \E[\Big]{ \paren[\Big]{ \frac{\indicator{Z_t = 0}}{1-p_t} - 1 }^2 \mid \filt_{t-1} } \cdot \Delta_t(0)^2 \\
		 &\quad - 2 \E[\Big]{ \paren[\Big]{ \frac{\indicator{Z_t = 1}}{p_t} - 1 } \cdot \paren[\Big]{ \frac{\indicator{Z_t = 0}}{1-p_t} - 1 }  \mid \filt_{t-1} } \Delta_t(1) \Delta_t(0)  \\
		 &= \paren[\Big]{ \frac{1}{p_t} - 1  } \cdot \Delta_t(1)^2
		 + \paren[\Big]{ \frac{1}{1-p_t} - 1  }  \cdot \Delta_t(0)^2
		 + 2 \Delta_t(1) \Delta_t(0)  \\
		 &= \paren[\Bigg]{ \sqrt{\frac{1-p_t}{p_t}} \Delta_t(1) + \sqrt{ \frac{p_t}{1-p_t} } \Delta_t(0)  }^2 \enspace.
	\end{align*}
	The proof is completed by substituting the calculation into the expression for the variance and recalling the definitions of $\Delta_t(1)$ and $\Delta_t(0)$.
\end{proof}

\subsection{Verification of the Sigmoidal Condition}\label{section:A2}

In this section, we verify that two commonly used sigmoid functions satisfy Condition~\mainref{condition:sigmoid}, which is required for our theoretical analysis.
We further provide some basic inequalities on the transformed space.
Namely, we focus on the \emph{trigonometric sigmoid} $\phi(u)=\frac{1}{\pi}(\operatorname{arctan}(u)+\pi/2)$ and the \emph{algebraic sigmoid} $\phi(u)=\frac{1}{2} \paren{ \frac{u}{1+\abs{u}} + 1 }$.
For completeness, we recall Condition~\mainref{condition:sigmoid} below:

\begin{refcondition}{\mainref{condition:sigmoid}}
	\sigmoidcondition
\end{refcondition}

First, we show that the trigonometric sigmoid $\phi(u)=\frac{1}{\pi}(\operatorname{arctan}(u)+\pi/2)$ satisfies these conditions.

\begin{lemma}\label{lemma:arctan}
    The sigmoid function $\phi(u)=\frac{1}{\pi}(\operatorname{arctan}(u)+\pi/2)$ satisfies Condition \mainref{condition:sigmoid} with constants $\bb_1=\pi$, $\bb_2=\frac{2^{5/2}\pi}{3}$ and $\bb_3=\frac{2}{\pi}$.
\end{lemma}

\begin{proof}
    We verify the requirements in Condition \mainref{condition:sigmoid} one by one.
    \begin{enumerate}
    \item[(1)] Due to the fact that $\operatorname{arctan}(u)$ is a strictly increasing odd function and $\lim_{u\rightarrow +\infty}\operatorname{arctan}(u)=\pi/2$, condition (1) is easily verified.
    \item[(2)] We first verify the convexity of the function $1/\phi(u)$ by calculating its second derivative.
    The derivative of $1/\phi(u)$ is calculated as:
    \begin{align*}
        &\left(\frac{1}{\phi(u)}\right)^{\prime}=-\frac{\pi}{(1+u^2)\left(\operatorname{arctan}(u)+\frac{1}{2}\pi\right)^2}\enspace.
    \end{align*}
    The second derivative of $1/\phi(u)$ is hence calculated as:
    \begin{align*}
        \left(\frac{1}{\phi(u)}\right)^{\prime\prime}=\frac{2}{\pi^2(1+u^2)^2\left(\frac{\operatorname{arctan}(u)}{\pi}+\frac{1}{2}\right)^3}\left[1+u\left(\operatorname{arctan}(u)+\frac{1}{2}\pi\right)\right]\enspace.
    \end{align*}
    Since $\frac{\operatorname{arctan}(u)}{\pi}+\frac{1}{2}>0$, it suffices to verify that $u\left(\operatorname{arctan}(u)+\frac{1}{2}\pi\right)>-1$ to prove the convexity. 
    When $u\geq 0$, $u\left(\operatorname{arctan}(u)+\frac{1}{2}\pi\right)\geq 0$.
    When $u<0$, by basic properties of the arctangent function we have
    \begin{align*}
        u\left(\operatorname{arctan}(u)+\frac{1}{2}\pi\right)=&-(-u)\left(-\operatorname{arctan}(-u)+\frac{1}{2}\pi\right)\\
        =&-(-u)\operatorname{arctan}(-1/u)\\
        >&-(-u)\cdot (-1/u)\\
        =&-1\enspace.
    \end{align*}
    Hence $1/\phi(u)$ is convex. By condition (1), we have $1/(1-\phi(u))=1/\phi(-u)$.
    The convexity of $1/(1-\phi(u))$ is then implied by the convexity of $1/\phi(u)$.
    \item[(3a)] By the convexity of $1/\phi(u)$, we only need to lower bound its gradient when $u$ tends to $-\infty$ to verify condition 3(a).
    By L'H\^opital's rule, we have
    \begin{align*}
        \lim_{u\rightarrow -\infty}\sqrt{1+u^2}\left(\operatorname{arctan}(u)+\frac{1}{2}\pi\right)=&\lim_{u\rightarrow -\infty}\frac{\left(\operatorname{arctan}(u)+\frac{1}{2}\pi\right)}{1/\sqrt{1+u^2}}\\
        =&\lim_{u\rightarrow -\infty}\frac{1/(1+u^2)}{-u(1+u^2)^{-3/2}}\quad(\text{L'H\^opital's rule})\\
        =&\lim_{u\rightarrow -\infty}\frac{(1+u^2)^{1/2}}{-u}\\
        =&1\enspace.
    \end{align*}
    By the explicit form of $\left(1/\phi(u)\right)^{\prime}$ that we have calculated in Part 2, it holds that $-\left(1/\phi(u)\right)^{\prime}\leq \pi$ since $\left(1/\phi(u)\right)^{\prime}$ is monotone increasing by its convexity.
    \item[(3b)] We first prove a stronger result than condition 3(b): For any $u\in\mathbb{R}$,
    \begin{align}\label{lemma:arctan_eq1}
        \left(\frac{1}{\phi(u)}\right)^{\prime\prime}\leq \frac{2\pi}{3}(1+u^2)^{-3/2}\enspace.
    \end{align}
	To see that \eqref{lemma:arctan_eq1} implies part 3b of Condition~\mainref{condition:sigmoid}, observe that it implies that for any $u \in \Reals$,
	\[
	\left(\frac{1}{\phi(u)}\right)^{\prime\prime}
	\leq \frac{2\pi}{3}(1+u^2)^{-3/2}
	\leq \frac{2\pi}{3}\cdot 2^{3/2}(1+|u|)^{-3}
	=\frac{2^{5/2}\pi}{3}(1+|u|)^{-3}
	\enspace,
	\]
	where the second inequality is implied by inequality $1+u^2\geq \frac{1}{2}(1+|u|)^2$.
	
	Now we will prove the stronger condition \eqref{lemma:arctan_eq1}.
	First, consider a change of variables $t=\operatorname{arctan}(u)+\frac{1}{2}\pi$ and note that $t\in (0,\pi)$.
	Then we have
    \begin{align*}
        1+u^2=&1+\tan^2(t-\pi/2)=1+\frac{\cos^2 t}{\sin^2 t}=\frac{\cos^2 t+\sin^2 t}{\sin^2 t}=\frac{1}{\sin^2 t}\enspace,\\
        u=&\tan(t-\pi/2)=-\frac{\cos t}{\sin t}\enspace.
    \end{align*}
    By the explicit form of $\left(1/\phi(u)\right)^{\prime\prime}$ calculated in Part 2, we can transform the function as
    \begin{align*}
        \frac{1}{2\pi}(1+u^2)^{3/2}\left(\frac{1}{\phi(u)}\right)^{\prime\prime}
        =&\frac{1}{\pi^3(1+u^2)^{1/2}\left(\frac{\operatorname{arctan}(u)}{\pi}+\frac{1}{2}\right)^3}\left[1+u\left(\operatorname{arctan}(u)+\frac{1}{2}\pi\right)\right]\\
        =&\frac{\sin t}{t^3}\left(1-\frac{t\cos t}{\sin t}\right)\\
        =&\frac{\sin t-t\cos t}{t^3}\\
        :=& f(t)\enspace.
    \end{align*}
    In order to verify \eqref{lemma:arctan_eq1}, it suffices to prove that for any $t\in(0,\pi)$, $f(t)\leq \frac{1}{3}$. 
    In order to verify this, we prove that $\lim_{t\downarrow 0}f(t)=\frac{1}{3}$ and $f(t)$ is monotone decreasing.
    By L'H\^opital's rule, we have
    \begin{align*}
        \lim_{t\downarrow 0}f(t)=&\lim_{t\downarrow 0}\frac{\sin t-t\cos t}{t^3}=\lim_{t\downarrow 0}\frac{\cos t-\cos t+t\sin t}{3t^2}=\lim_{t\downarrow 0}\frac{\sin t}{3t}=\frac{1}{3}\enspace.
    \end{align*}
    To verify the monotonicity of $f(t)$, we prove the negativity of
    \begin{align*}
        f^{\prime}(t)=\frac{\sin t(t^2-3)+3t\cos t}{t^4}:= \frac{k(t)}{t^4}\enspace,
    \end{align*}
    which is equivalent to proving the negativity of $k(t)$ on $(0,\pi)$.
    We then prove that $\lim_{t\downarrow 0}k(t)=0$ and $k(t)$ is monotone decreasing.
    Proving $\lim_{t\downarrow 0}k(t)=0$ is straightforward.
    To verify the monotonicity of $k(t)$, we prove the negativity of $k^{\prime}(t)=t(t\cos t-\sin t)$, which can be easily seen since $t\cos t-\sin t$ equals $0$ at $t=0$ and $(t\cos t-\sin t)^{\prime}=-t\sin t<0$ on $(0,\pi)$.
    This completes the proof of inequality \eqref{lemma:arctan_eq1}.
    \item[(3c)] We use the function $f(t)$ defined in Part 3b.
    Since $f(t)$ is proved to be monotone decreasing and $\lim_{t\uparrow \pi}f(t)=1/\pi^2$, we have $f(t)\geq 1/\pi^2$ for any $t\in(0,\pi)$.
    This implies that
    \begin{align*}
        \frac{1}{2\pi}(1+u^2)^{3/2}\left(\frac{1}{\phi(u)}\right)^{\prime\prime}=f(t)\geq 1/\pi^2\enspace.
    \end{align*}
    Hence for any $u\geq 0$, by inequality $1+u^2\leq (1+u)^2$, we have
    \begin{align*}
        \left(\frac{1}{\phi(u)}\right)^{\prime\prime}\geq \frac{2}{\pi}(1+u^2)^{-3/2}\geq \frac{2}{\pi}(1+u)^{-3}\enspace.
    \end{align*}
    \end{enumerate}
    Hence the sigmoid function $\phi(u)$ satisfies Condition \mainref{condition:sigmoid} with constants $\bb_1=\pi$, $\bb_2=\frac{2^{5/2}\pi}{3}$ and $\bb_3=\frac{2}{\pi}$.
\end{proof}

Next, we show that the algebraic sigmoid $\phi(u)=\frac{1}{2} \paren{ \frac{u}{1+\abs{u}} + 1 }$ also satisfies Condition~\mainref{condition:sigmoid}.

\begin{lemma}\label{lemma:algebraic}
    The sigmoid function $\phi(u)=\frac{1}{2} \paren{ \frac{u}{1+\abs{u}} + 1 }$ satisfies Condition \mainref{condition:sigmoid} with constants $\bb_1=2$, $\bb_2=8$ and $\bb_3=1$.
\end{lemma}

\begin{proof}
    We verify the requirements in Condition \mainref{condition:sigmoid} one by one.
    \begin{enumerate}
    \item[(1)] The strict monotonicity of $\phi(u)$ can be seen from its derivative $\phi^{\prime}(u)=\frac{1}{2(1+|u|)^2}$.
    Since $\frac{u}{1+\abs{u}}$ is an odd function and $\lim_{u\rightarrow\infty}\phi(u)=1$, we have verified the other two conditions.
    \item[(2)] We first verify the convexity of the function $1/\phi(u)$ by showing that its derivative is monotone increasing.
    The derivative of $1/\phi(u)$ is calculated as:
    \begin{align*}
        \left(\frac{1}{\phi(u)}\right)^{\prime}=&\begin{cases}
            -\frac{2}{(2u+1)^2}&\text{if $u\geq 0$}\\
            -2&\text{if $u\leq 0$}
        \end{cases}\enspace.
    \end{align*}
    Because $\left(1/\phi(u)\right)^{\prime}$ is monotone increasing, we have that $1/\phi(u)$ is convex. 
    By condition (1), we have $1/(1-\phi(u))=1/\phi(-u)$.
    The convexity of $1/(1-\phi(u))$ is then implied by the convexity of $1/\phi(u)$.
    \item[(3a)] By the explicit form of $\left(1/\phi(u)\right)^{\prime}$ calculated in Part 2, we have that for any $u\in\mathbb{R}$,
    \begin{align*}
        -\left(\frac{1}{\phi(u)}\right)^{\prime}\leq 2\enspace.
    \end{align*}
    \item[(3b)] The second derivative of $1/\phi(u)$ is calculated as:
    \begin{align*}
        \left(\frac{1}{\phi(u)}\right)^{\prime\prime}=&\begin{cases}
            \frac{8}{(2u+1)^3}&\text{if $u> 0$}\\
            0&\text{if $u< 0$}
        \end{cases}\enspace.
    \end{align*}
    We have that for any $u\in\mathbb{R}\setminus\{0\}$,
    \begin{align*}
        \left(\frac{1}{\phi(u)}\right)^{\prime\prime}\leq \frac{8}{(1+2|u|)^3}\leq \frac{8}{(1+|u|)^3}\enspace.
    \end{align*}
    \item[(3c)] By the explicit form of $\left(1/\phi(u)\right)^{\prime\prime}$ calculated in Part 3b, we have that for any $u> 0$,
    \begin{align*}
        \left(\frac{1}{\phi(u)}\right)^{\prime\prime}=\frac{8}{(2u+1)^3}\geq \frac{1}{(1+u)^3}\enspace.
    \end{align*}
\end{enumerate}
Hence the sigmoid function $\phi(u)$ satisfies Condition \mainref{condition:sigmoid} with constants $\bb_1=2$, $\bb_2=8$ and $\bb_3=1$.
\end{proof}

Recall that the sigmoid $\phi : \Reals \to (0,1)$ transforms the ill-conditioned problem of choosing $p \in (0,1)$ into the well-conditioned problem of choosing $u \in \Reals$.
Throughout our analysis, we will need to move back and forth between what we call the \emph{$p$-space}, i.e., the open interval $(0,1)$, and the \emph{$u$-space}, i.e., the real line $\Reals$, via this sigmoidal transformation.
For example, if $u$ and $u'$ are close, then how does this ensure closeness of $p = \phi(u)$ and $p' = \phi(u')$?
The following lemma provides a variety of quantitative answers to these questions which we use throughout our analysis.

\begin{lemma}\label{lemma:u}
    Fix $u,\widetilde u\ge 0$ and let $p=\phi(u)$, $\widetilde p=\phi(\widetilde u)$.
    Under Condition \mainref{condition:sigmoid}, we have the following results:
    \begin{enumerate}
        \item[(1)] $\left(\frac{1}{1-\phi(u)}\right)^{\prime}\big|_{u=0}\geq \frac{\bb_3}{2}$.
        \item[(2)] $|\left(\frac{1}{\phi(u)}\right)^{\prime}|\leq \frac{\bb_2}{2}\cdot \frac{1}{(1+u)^2}$.
        \item[(3)] $\frac{\bb_3}{2}\cdot|u-\widetilde{u}|\leq|\frac{1}{1-p}-\frac{1}{1-\widetilde{p}}|\leq \bb_1\cdot|u-\widetilde{u}|$. The second inequality holds for any $u,\widetilde{u}\in\mathbb{R}$.
        \item[(4)] $|\frac{1}{p}-\frac{1}{\widetilde{p}}|\leq \frac{\bb_2}{2}\cdot \frac{|u-\widetilde{u}|}{(1+u)(1+\widetilde{u})}$.
        \item[(5)] $|\left(\frac{1}{\phi(u)}\right)^{\prime}-\left(\frac{1}{\phi(\widetilde{u})}\right)^{\prime}|\geq \frac{\bb_3}{2}\cdot\frac{|u-\widetilde{u}|}{(1+u)(1+\widetilde{u})(1+\min\{u,\widetilde{u}\})}$.
        \item[(6)] $\max\left\{\frac{1}{\phi(u)}-\frac{1}{\phi(\widetilde{u})}-\left(\frac{1}{\phi(\widetilde{u})}\right)^{\prime}(u-\widetilde{u}),\frac{1}{1-\phi(u)}-\frac{1}{1-\phi(\widetilde{u})}-\left(\frac{1}{1-\phi(\widetilde{u})}\right)^{\prime}(u-\widetilde{u})\right\}\leq \frac{\bb_2}{2}\cdot \frac{\widetilde{u}-u}{(1+u)(1+\widetilde{u})}$ for $\widetilde{u}\geq u$.
    \end{enumerate}
\end{lemma}

\begin{proof}
	
	\noindent
	\textbf{Part 1}:
	By the symmetry property in Condition \mainref{condition:sigmoid}-1, we have
    \begin{align*}
        \frac{\mathrm{d}^2}{\mathrm{d}u^2} \left(\frac{1}{1-\phi(u)}\right)=\frac{\mathrm{d}^2}{\mathrm{d}u^2} \left(\frac{1}{\phi(-u)}\right)=\left(\frac{1}{\phi}\right)^{\prime\prime}(-u)\enspace.
    \end{align*}    
    Recall that the monotonicity property of Condition \mainref{condition:sigmoid}-1 requires that $\left(1/(1-\phi(u))\right)^{\prime}\geq 0$ for any $u$ and Condition \mainref{condition:sigmoid}-3c requires that $(1/\phi(u))'' \geq \bb_3/(1+u)^3$ for all $u \geq 0$ except for a finite number of points.
	Using these parts of Condition \mainref{condition:sigmoid}, we obtain that
	\begin{align*}
		 \left(\frac{1}{1-\phi(u)}\right)^{\prime}\Big|_{u=0}
		 &= \int_{-\infty}^{0}\left(\frac{1}{1-\phi(u)}\right)^{\prime\prime}\mathrm{d}u+\lim_{u\rightarrow-\infty }\left(\frac{1}{1-\phi(u)}\right)^{\prime} \\
		 &\geq 
		 \int_{-\infty}^{0}\left(\frac{1}{1-\phi(u)}\right)^{\prime\prime}\mathrm{d}u
		 	&\text{(Condition \mainref{condition:sigmoid}-1)} \\
            &=\int_{-\infty}^{0}\left(\frac{1}{\phi}\right)^{\prime\prime}(-u)\mathrm{d}u\\
		 &\geq \bb_3\int_{-\infty}^{0}\frac{1}{(1-u)^3}\mathrm{d}u
		 	&\text{(Condition \mainref{condition:sigmoid}-3c)} \\
		 &= \frac{\bb_3}{2}\enspace.
	\end{align*}

	\noindent
    \textbf{Part 2}: 
    By the convexity property in Condition \mainref{condition:sigmoid}-2, $\left(1/\phi(u)\right)^{\prime}$ is nondecreasing.
    On the other hand, $\lim_{u\rightarrow\infty }1/\phi(u)=1$ by Condition \mainref{condition:sigmoid}-1, which implies that $\lim_{u\rightarrow\infty }\left(1/\phi(u)\right)^{\prime}=0$.
    Using this together with Condition \mainref{condition:sigmoid}-3b, we have that
    \begin{align*}
        \left|\left(\frac{1}{\phi(u)}\right)^{\prime}\right|
        &=-\left(\frac{1}{\phi(u)}\right)^{\prime} \\
        &=\int_{u}^{\infty}\left(\frac{1}{\phi(t)}\right)^{\prime\prime}\mathrm{d}t-\lim_{u\rightarrow\infty }\left(\frac{1}{\phi(u)}\right)^{\prime} \\
        &=\int_{u}^{\infty}\left(\frac{1}{\phi(t)}\right)^{\prime\prime}\mathrm{d}t \\
        &\leq \bb_2\int_{u}^{\infty}\frac{1}{(1+t)^3}\mathrm{d}t
        	&\text{(Condition \mainref{condition:sigmoid}-3b)} \\
        &= \frac{\bb_2}{2(1+u)^2}\enspace.
    \end{align*}

	\noindent
   \textbf{Part 3}: 
 	First, let us establish this upper bound.
 	We will use the symmetry property in Condition \mainref{condition:sigmoid}-1 and the bound in Condition \mainref{condition:sigmoid}-3a.
 	By the mean value theorem, there exists $\widebar{u}$ between $u$ and $\widetilde{u}$ such that
    \begin{align*}
        \left|\frac{1}{1-p}-\frac{1}{1-\widetilde{p}}\right|
        =\left|\frac{1}{1-\phi(u)}-\frac{1}{1-\phi(\widetilde{u})}\right| 
        &=\left|\left(\frac{1}{1-\phi(\widebar{u})}\right)^{\prime}\right||u-\widetilde{u}| 
        	&\text{(mean value theorem)}\\
        &=\left|\left(\frac{1}{\phi(-\widebar{u})}\right)^{\prime}\right||u-\widetilde{u}|
        	&\text{(Condition \mainref{condition:sigmoid}-1)} \\
        & \leq \bb_1|u-\widetilde{u}|
        \enspace.
        &\text{(Condition \mainref{condition:sigmoid}-3a)}
    \end{align*}
    Next, we prove the lower bound.
    By the result in Part 1 of this proof and the convexity property in Condition \mainref{condition:sigmoid}-2, we have
    \begin{align*}
        \left|\frac{1}{1-p}-\frac{1}{1-\widetilde{p}}\right|=&\left|\frac{1}{1-\phi(u)}-\frac{1}{1-\phi(\widetilde{u})}\right|\\
        =&\left|\int_{\widetilde{u}}^{u}\left(\frac{1}{1-\phi(t)}\right)^{\prime}\mathrm{d}t\right|\\
        \geq&\left|\int_{\widetilde{u}}^{u}\left(\frac{1}{1-\phi(0)}\right)^{\prime}\mathrm{d}t\right|
        	&\text{(convexity)} \\
        \geq&\frac{\bb_3}{2}|u-\widetilde{u}|
        &\text{(Part 1)} 
        \enspace.
    \end{align*}

	\noindent
    \textbf{Part 4}:
    By the result in Part 2 of the proof, we have
    \begin{align*}
        \left|\frac{1}{p}-\frac{1}{\widetilde{p}}\right|=\left|\frac{1}{\phi(u)}-\frac{1}{\phi(\widetilde{u})}\right|=\left|\int_{u}^{\widetilde{u}}\left(\frac{1}{\phi(t)}\right)^{\prime}\mathrm{d}t\right|\leq\frac{\bb_2}{2}\left|\int_{u}^{\widetilde{u}}\frac{1}{(1+t)^2}\mathrm{d}t\right|=\frac{\bb_2|u-\widetilde{u}|}{2(1+u)(1+\widetilde{u})}\enspace.
    \end{align*}

	\noindent
   \textbf{Part 5}:
    By Condition \mainref{condition:sigmoid}-3c, we have
    \begin{align*}
        \left|\left(\frac{1}{\phi(u)}\right)^{\prime}-\left(\frac{1}{\phi(\widetilde{u})}\right)^{\prime}\right|
        =&\left|\int_{\widetilde{u}}^{u}\left(\frac{1}{\phi(t)}\right)^{\prime\prime}\mathrm{d}t\right| \\
        \geq& \bb_3\left|\int_{\widetilde{u}}^{u}\frac{1}{(1+t)^3}\mathrm{d}t\right|
        \quad(\text{Condition \mainref{condition:sigmoid}-3c})\notag\\
        =&\frac{\bb_3}{2}\left|\frac{1}{(1+\widetilde{u})^2}-\frac{1}{(1+u)^2}\right|\notag\\
        =&\frac{\bb_3}{2}\frac{(2+\widetilde{u}+u)}{(1+\widetilde{u})^2(1+u)^2}|u-\widetilde{u}|\notag\\
        \geq&\frac{\bb_3}{2}\frac{1}{(1+\widetilde{u})(1+u)(1+\min\{u,\widetilde{u}\})}|u-\widetilde{u}| 
        \enspace,
    \end{align*}
	where the last inequality follows because $2+\widetilde{u}+u\geq \max\{1+\widetilde{u},1+u\}$.\\

	\noindent
    \textbf{Part 6}: 
    We upper bound each term separately.
    By the monotonicity property in Condition \mainref{condition:sigmoid}-1 and the convexity property in Condition \mainref{condition:sigmoid}-2, for $\widetilde{u}\geq u$ we have
    \begin{align*}
        \frac{1}{\phi(u)}-\frac{1}{\phi(\widetilde{u})}-\left(\frac{1}{\phi(\widetilde{u})}\right)^{\prime}(u-\widetilde{u})\leq\frac{1}{\phi(u)}-\frac{1}{\phi(\widetilde{u})}\enspace.
    \end{align*}
    Then by the result in Part 4, we have
    \begin{align*}
        \frac{1}{\phi(u)}-\frac{1}{\phi(\widetilde{u})}-\left(\frac{1}{\phi(\widetilde{u})}\right)^{\prime}(u-\widetilde{u})\leq\frac{1}{\phi(u)}-\frac{1}{\phi(\widetilde{u})}\leq\frac{\bb_2(\widetilde{u}-u)}{2(1+u)(1+\widetilde{u})}\enspace.
    \end{align*}
    On the other hand, by the symmetry property in Condition \mainref{condition:sigmoid}-1, and Condition \mainref{condition:sigmoid}-3b, we have
    \begin{align*}
        \frac{1}{1-\phi(u)}-\frac{1}{1-\phi(\widetilde{u})}- &\left(\frac{1}{1-\phi(\widetilde{u})}\right)^{\prime}(u-\widetilde{u})\\
        =&\int_{u}^{\widetilde{u}}\int_{t}^{\widetilde{u}}\left(\frac{1}{1-\phi(s)}\right)^{\prime\prime}\mathrm{d}s\mathrm{d}t\\
        =&\int_{u}^{\widetilde{u}}\int_{t}^{\widetilde{u}}\left(\frac{1}{\phi(-s)}\right)^{\prime\prime}\mathrm{d}s\mathrm{d}t
        	&\text{(symmetry)}\\
        \leq&\bb_2\int_{u}^{\widetilde{u}}\int_{t}^{\widetilde{u}}\frac{1}{(1+s)^3}\mathrm{d}s\mathrm{d}t
        	&(\text{Condition \mainref{condition:sigmoid}-3b})\\
        =&\frac{\bb_2}{2}\int_{u}^{\widetilde{u}}\left[\frac{1}{(1+t)^2}-\frac{1}{(1+\widetilde{u})^2}\right]\mathrm{d}t\\
        \leq&\frac{\bb_2}{2}\left[\frac{1}{1+u}-\frac{1}{1+\widetilde{u}}-\frac{1}{(1+\widetilde{u})^2}(\widetilde{u}-u)\right]\\
        =&\frac{\bb_2}{2}\cdot\frac{(1+\widetilde{u})(\widetilde{u}-u)-(1+u)(\widetilde{u}-u)}{(1+u)(1+\widetilde{u})^2}\\
        =&\frac{\bb_2}{2}\cdot\frac{(\widetilde{u}-u)^2}{(1+u)(1+\widetilde{u})^2}\\
        \leq&\frac{\bb_2}{2}\cdot \frac{\widetilde{u}-u}{(1+u)(1+\widetilde{u})}
        	&\text{($\widetilde{u}\geq u\geq 0$)}
        \enspace,
    \end{align*}
    which completes the proof.
\end{proof}

\subsection{General FTRL analysis}\label{section:A3}
In this section, we provide a standard analysis of regret for an FTRL algorithm in a general OCO problem.
Here is the general online optimization framework. At each iteration, the player selects an action $x_t$ from a convex set $\mathcal{X}$.
Then, the convex loss function $f_t: \mathcal{X} \rightarrow \mathbb{R}$ is revealed and the player incurs a loss $f_t\left(x_t\right)$.
We allow for the loss function $f_t$ to depend on previously chosen actions, i.e., $x_1, \ldots, x_{t-1}$.

The Follow-the-Regularized-Leader (FTRL) algorithm requires the specification of a regularizer function $R: \mathcal{X} \rightarrow \mathbb{R}_+$, which we assume is differentiable.
Assume that $\min_{x\in\mathcal{X}}R(x)=0$, i.e., the regularizer can attain a zero value for certain actions. 
Assume further that we use a sequence of decreasing step sizes $\{\eta_t\}_{t\geq 1}$. For notational simplicity, let $\eta_0=\eta_1$.
The iterates of the FTRL algorithm are defined as
\begin{align*}
    x_{t}\in\arg\min_{x\in\mathcal{X}}\sum_{s<t}f_s(x)+\eta_{t}^{-1}R(x)\enspace.
\end{align*}
For any $t\geq 1$, we denote
\begin{align*}
    F_t(x):=&\sum_{s<t}f_s(x)+\eta_{t}^{-1}R(x)\enspace,\\
    \widetilde{F}_t(x):=&\sum_{s<t}f_s(x)+\eta_{t-1}^{-1}R(x)\enspace.
\end{align*}
By definition we have $x_t\in\arg\min_{x\in\mathcal{X}} F_t(x)$.
We further denote $\widetilde{x}_t\in \arg\min_{x\in\mathcal{X}}\widetilde{F}_t(x)$.
The following lemma decomposes the regret in an FTRL algorithm into three parts: the regularizer, the sum of telescoping differences, and the difference in the last step.

\begin{lemma}\label{lemma:standard-FTRL-decomposition}
    For any $x^*\in\mathcal{X}$, the iterates of FTRL satisfy the following equality:
    \begin{align*}
        \sum_{t=1}^Tf_t(x_t)-\sum_{t=1}^Tf_t(x^*)=\eta_{T+1}^{-1}R(x^*)+\sum_{t=1}^T\left[\widetilde{F}_{t+1}(x_t)-F_{t+1}(x_{t+1})\right]+F_{T+1}(x_{T+1})-F_{T+1}(x^*)\enspace.
    \end{align*}
\end{lemma}

\begin{proof}
    By definition, we have $F_{T+1}(x^*)=\sum_{t=1}^Tf_t(x^*)+\eta_{T+1}^{-1}R(x^*)$ and $F_1(x_1)=\eta_1^{-1}R(x_1)=\eta_1^{-1}\min_{x\in\mathcal{X}}R(x)=0$.
    Hence we obtain
    \begin{align*}
        &\sum_{t=1}^Tf_t(x_t)-\sum_{t=1}^Tf_t(x^*)\\
        =&\eta_{T+1}^{-1}R(x^*)+\sum_{t=1}^Tf_t(x_t)+F_1(x_1)-F_{T+1}(x^*)\\
        =&\eta_{T+1}^{-1}R(x^*)+\sum_{t=1}^Tf_t(x_t)+F_1(x_1)-F_{T+1}(x_{T+1})+F_{T+1}(x_{T+1})-F_{T+1}(x^*)\quad(\text{add and subtract})\\
        =&\eta_{T+1}^{-1}R(x^*)+\sum_{t=1}^T\left[f_t(x_t)+F_t(x_t)-F_{t+1}(x_{t+1})\right]+F_{T+1}(x_{T+1})-F_{T+1}(x^*)\\
        \intertext{By definition $f_t(x_t)+F_t(x_t)=\widetilde{F}_{t+1}(x_t)$. Hence}
        =&\eta_{T+1}^{-1}R(x^*)+\sum_{t=1}^T\left[\widetilde{F}_{t+1}(x_t)-F_{t+1}(x_{t+1})\right]+F_{T+1}(x_{T+1})-F_{T+1}(x^*)\enspace,
    \end{align*}
    which completes the proof.
\end{proof}

Using the optimality of $\widetilde{x}_t$ and $x_t$ in their definitions, the following corollary directly follows from Lemma \ref{lemma:standard-FTRL-decomposition}.

\begin{corollary}\label{corollary:standard-FTRL-upper-bound}
    For any $x^*\in\mathcal{X}$, we have the upper bound:
    \begin{align*}
        \sum_{t=1}^Tf_t(x_t)-\sum_{t=1}^Tf_t(x^*)\leq \eta_{T+1}^{-1}R(x^*)+\sum_{t=1}^T\left[\widetilde{F}_{t+1}(x_t)-\widetilde{F}_{t+1}(\widetilde{x}_{t+1})\right]\enspace.
    \end{align*}
\end{corollary}

\begin{proof}
    By Lemma \ref{lemma:standard-FTRL-decomposition}, we have the following equality:
    \begin{align*}
        \sum_{t=1}^Tf_t(x_t)-\sum_{t=1}^Tf_t(x^*)=\eta_{T+1}^{-1}R(x^*)+\sum_{t=1}^T[\widetilde{F}_{t+1}(x_t)-F_{t+1}(x_{t+1})]+F_{T+1}(x_{T+1})-F_{T+1}(x^*)\enspace.
    \end{align*}
    Since $x_{T+1}$ minimizes $F_{T+1}$, we have
    \begin{align*}
        F_{T+1}(x_{T+1})-F_{T+1}(x^*)\leq 0\enspace.
    \end{align*}
    Since $\eta_{t+1}\leq \eta_{t}$, for any $t\in[T]$ we have
    \begin{align*}
        \widetilde{F}_{t+1}(x_t)-F_{t+1}(x_{t+1})\leq \widetilde{F}_{t+1}(x_t)-\widetilde{F}_{t+1}(x_{t+1})\enspace.
    \end{align*}
    Since $\widetilde{x}_{t+1}$ minimizes $\widetilde{F}_{t+1}$, for any $t\in[T]$ we have
    \begin{align*}
        \widetilde{F}_{t+1}(x_t)-\widetilde{F}_{t+1}(x_{t+1})\leq \widetilde{F}_{t+1}(x_t)-\widetilde{F}_{t+1}(\widetilde{x}_{t+1})\enspace.
    \end{align*}
    Combining the results, we have
    \begin{align*}
        \sum_{t=1}^Tf_t(x_t)-\sum_{t=1}^Tf_t(x^*)\leq \eta_{T+1}^{-1}R(x^*)+\sum_{t=1}^T\left[\widetilde{F}_{t+1}(x_t)-\widetilde{F}_{t+1}(\widetilde{x}_{t+1})\right]\enspace,
    \end{align*}
    which completes the proof.
\end{proof}
	
	\section{Neyman Regret Lower Bound}\label{sec:supp-lower-bound}

In this section, we focus on constructing lower bounds for the Neyman Regret under the assumptions considered in the main body.
In Section \ref{section:B1}, we prove the lower bound on the Neyman regret stated in Theorem \mainref{theorem:neyman-regret-lower-bound}.
In Section \ref{section:B2}, we justify the necessity of the imposed assumptions in order to obtain a regret lower bound of order $\Omega(T^{-1/2})$.

\subsection{Lower Bound (Theorem \mainref{theorem:neyman-regret-lower-bound})}\label{section:B1}
We first provide a simple lemma regarding the sample variance for independent but nonidentical random variables, which simplifies the calculation in the construction of the regret lower bound.

\begin{lemma}\label{lemma:sample-variance}
    Suppose $X_1,\ldots,X_T$ are independent random variables and let $\widebar{X} = \frac{1}{T} \sum_{t=1}^T X_t$ denote their sample mean.
    Then we have the following equation:
    \begin{align*}
        \e\left[\sum_{t=1}^T(X_t-\widebar{X})^2\right]=\frac{T-1}{T}\sum_{t=1}^T\operatorname{Var}(X_t)+\sum_{t=1}^{T}\e[X_t]^2-T\,\e[\widebar{X}]^2\enspace.
    \end{align*}
\end{lemma}
\begin{proof}
	Observe that by definition of $\widebar{X}$, we have the identity
	\[
	\sum_{t=1}^T (X_t - \widebar{X})^2 
	= \sum_{t=1}^T X_t^2 - T \cdot \widebar{X}^2
	\enspace.
	\]
	Using this, we can calculate the expectation:
    \begin{align*}
        \e\left[\sum_{t=1}^T(X_t-\widebar{X})^2\right]=&\sum_{t=1}^T\e[X_t^2]-T\cdot\e[\widebar{X}^2]\\
        =&\sum_{t=1}^T\operatorname{Var}(X_t)+\sum_{t=1}^T\e[X_t]^2-T\cdot \operatorname{Var}(\widebar{X})-T\cdot \e[\widebar{X}]^2\\
        =&\sum_{t=1}^T\operatorname{Var}(X_t)+\sum_{t=1}^T\e[X_t]^2-T\cdot \frac{1}{T^2}\sum_{t=1}^T\operatorname{Var}({X}_t)-T\cdot \e[\widebar{X}]^2 \\
        =&\frac{T-1}{T}\sum_{t=1}^T\operatorname{Var}(X_t)+\sum_{t=1}^{T}\e[X_t]^2-T\,\e[\widebar{X}]^2\enspace,
    \end{align*}
	where the first equality follows from the above identity and linearity of expectation, the second equality comes from using $\E{ \widebar{X}^2 }  = \Var{\widebar{X}} + \E{ \widebar{X} }^2$, the third equality follows from independence of the $X_1, \dots, X_T$, and the fourth equality comes from collecting terms.
\end{proof}

We now provide the proof of Theorem \mainref{theorem:neyman-regret-lower-bound}.
The proof follows a standard lower-bound construction technique in the online convex optimization literature (e.g., see \cite{Hazan2016Introduction} Chapter 3.2). 
The main idea is to construct a distribution over random problem instances, i.e. sequences $\setb{ y_t(1), y_t(0), \xv_t}_{t=1}^T$ satisfying Assumptions~\mainref{assumption:moments}-\mainref{assumption:maximum-radius}, and show that, for an arbitrary adaptive design, the expected Neyman Regret over this distribution of problem instances is at least $\bigOmega{T^{-1/2}}$.
To do this, we must show that there is a \emph{separation} that occurs: an arbitrary adaptive algorithm must incur a large expected AIPW variance over this distribution of problem instances, whereas the expected oracle Neyman variance should be small.

\begin{reftheorem}{\mainref{theorem:neyman-regret-lower-bound}}[Lower Bound]
	\lowerbound
\end{reftheorem}

\begin{proof}
For expositional clarity, we will consider a simple setting where $d=1$ and $R=1$; however, the same argument can be extended to larger $d$ and $R$.
For notational convenience, we define $\Delta_{t,k}(\bv)=y_t(k)-\iprod{\xv_t,\bv}$ for $k\in\{0,1\}$.
For simplicity, we suppose that $T$ is an even number.

Our construction of the hard instance proceeds as follows:
let $\xv_{t}=1$ for all $t\in[T]$.
We construct the random sequence of potential outcomes $\{y_t(1),y_t(0):t \in[T] \}$ according to the following procedure.

\begin{enumerate}
    \item Generate noise variables $\epsilon_1,\ldots,\epsilon_T$ where the even noise variables $\epsilon_{2k}$ are i.i.d. $\pm 1$ Rademacher random variables and the odd noise variables $\epsilon_{2k+1}$ are $0$.
    \item Generate a random vector $\vec{U}=(U_1,U_0)$ that is independent of $\epsilon_1,\ldots,\epsilon_T$ which takes value $(2,4)$ with probability $1/2$ and takes value $(-4,-2)$ with probability $1/2$.
    \item Let $y_1(1)=U_1(T^{1/4}+\epsilon_1)$ and $y_1(0)=U_0(T^{1/4}+\epsilon_1)$.
    \item For any $t=2,\ldots,T$, let $y_t(1)=U_1(1+\epsilon_t)$ and $y_t(0)=U_0(1+\epsilon_t)$.
\end{enumerate}

The high-level idea behind the proof is to compare (i) the expected conditional variance under an arbitrary adaptive design $\mathcal{D}$ and (ii) the expected conditional variance under the oracle design. The construction is a two-point mixture in $\vec U$ with a single ``large'' observation at $t=1$ (of order $T^{1/4}$) and Rademacher perturbations on even times.
This creates a nontrivial gap.
Any adaptive design $\mathcal{D}$ must pay for the uncertainty in $\vec U$ at least in the first round, while the oracle regression fit partially cancels it.
The inclusion of $\epsilon_1,\ldots,\epsilon_T$ is not essential for the regret bound, but is needed to ensure that the sums of squared OLS residuals are of order $\Omega(T)$ (Assumption \mainref{assumption:moments}).

We define the filtration as: $\filt_t=\sigma\{ Z_s, Y_s:1\leq s\leq t\}$.
Note that this filtration is different from the definition in the main paper since we are generating random potential outcomes in the lower bound construction.
By definition, we have that $p_t$ is $\filt_{t-1}$-measurable in the design $\fullarbdesign$ and does not depend on any randomness independent of $\filt_{t-1}$.
We write $\e_y[\cdot]$ for expectation over the random potential outcomes.
$\operatorname{Var}_{-y}(\cdot)$ denotes the conditional variance given a fixed realization of the potential outcomes.\\[3mm]
\textbf{Step 1:} Lower-bound the expected (conditional) variance under design $\fullarbdesign$.\\[3mm]
For any $k\in\{0,1\}$, denote $e_{1,k}=y_1(k)-\e[y_1(k)|\filt_0]$, $r_{1,k}=\e[y_1(k)|\filt_0]-\iprod{ \xv_1 , \bv_1(k) }$. 
Then $\e[e_{1,k}|\filt_0]=0$ and $r_{1,k}$ is $\filt_0$-measurable.
For any $k\in\{0,1\}$ and $2\leq t\leq T$, denote $e_{t,k}=y_t(k)-\e[y_t(k)|\filt_{t-1},\vec{U}]$ and $r_{t,k}=\e[y_t(k)|\filt_{t-1},\vec{U}]-\iprod{ \xv_t , \bv_t(k) }$.
Then $\e[e_{t,k}|\filt_{t-1},\vec{U}]=0$, and $r_{t,k}$ is $\sigma\{\filt_{t-1},\vec{U}\}$-measurable.

We separate the case $t=1$ because $\vec{U}$ is completely unknown at the beginning and is independent of $\filt_0$. 
Starting from $t\geq 2$, $\vec{U}$ may be dependent on $\filt_{t-1}$.
In the worst case the learner may have already inferred $\vec{U}$ exactly from past observations. Therefore, to obtain bounds that remain valid in all scenarios, we condition on $\vec{U}$ for $t\geq 2$. 

By Proposition \mainref{prop:aipw-variance}, the conditional variance of $\hat\tau$ given the realized potential outcomes can be written as a sum of per-time contributions.
Taking the expectation with respect to the potential outcomes yields the expression below:    
\begin{align}\label{theorem:neyman-regret-lower-bound_eq1}
    &\e_{y}\left[T\cdot \operatorname{Var}_{-y}(\hat{\tau})\right]\notag\\
    =&\frac{1}{T} \sum_{t=1}^T\e\left[\left(\Delta_{t,1}(\bv_t(1)) \cdot  \sqrt{\frac{1-p_t}{p_t}}+\Delta_{t,0}(\bv_t(0)) \cdot \sqrt{\frac{p_t}{1-p_t}} \right)^2\right]\notag\\
    =&\frac{1}{T}\sum_{t=1}^T\e\left[\left(e_{t,1}+r_{t,1}\right)^2 \cdot  \frac{1-p_t}{p_t}+\left(e_{t,0}+r_{t,0}\right)^2 \cdot \frac{p_t}{1-p_t} +2\left(e_{t,1}+r_{t,1}\right)\cdot\left(e_{t,0}+r_{t,0}\right) \right]\notag\\
    \intertext{Reorganizing the terms yields the decomposition}
    =&
    \frac{1}{T} 
    \underbrace{\sum_{t=1}^T \e\left[\left(r_{t,1}\cdot \sqrt{\frac{1-p_t}{p_t}}+r_{t,0}\cdot \sqrt{\frac{p_t}{1-p_t}}\right)^2\right]}_{:= S_1}
        \notag\\
        &+\frac{1}{T} \underbrace{\sum_{t=1}^T\e\left[e_{t,1}^2\cdot  \frac{1-p_t}{p_t}+e_{t,0}^2\cdot \frac{p_t}{1-p_t}+2e_{t,0}e_{t,1}\right]}_{:= S_2}
        \notag \\
    &+\frac{2}{T}\underbrace{\sum_{t=1}^T\e\left[e_{t,1}r_{t,1}\cdot\frac{1-p_t}{p_t}+e_{t,0}r_{t,0}\cdot\frac{p_t}{1-p_t}+e_{t,1}r_{t,0}+e_{t,0}r_{t,1}\right]}_{:= S_3}
    \enspace.
\end{align}
We then lower bound $S_1$, $S_2$ and $S_3$ separately.
Because $S_{1}$ is a nonnegative term, we can lower bound it by zero, i.e. $S_1 \geq 0$.

Next, we turn to analyzing $S_2$.
When $t=1$, since $\epsilon_1=0$ and $\vec{U}$ is independent of $\filt_0$, we can calculate the following:
\begin{align}\label{theorem:neyman-regret-lower-bound_eq2}
    \e[e_{1,k}^2|\filt_{0}]=&\operatorname{Var}(e_{1,k}|\filt_{0})=\operatorname{Var}(U_kT^{1/4}|\filt_{0})=T^{1/2}\operatorname{Var}(U_k)=9T^{1/2}\enspace,\notag\\
    \e[e_{1,0}e_{1,1}|\filt_{0}]=&\operatorname{Cov}(U_1T^{1/4},U_0T^{1/4}|\filt_{0})=T^{1/2}\operatorname{Cov}\left(U_1,U_0\right)=9T^{1/2}\enspace.
\end{align}
By \eqref{theorem:neyman-regret-lower-bound_eq2} and the law of iterated expectations, we obtain
\begin{align}\label{theorem:neyman-regret-lower-bound_eq4}
    &\e\left[e_{1,1}^2\cdot \frac{1-p_1}{p_1}+e_{1,0}^2\cdot \frac{p_1}{1-p_1}+2e_{1,0}e_{1,1}\right]\notag\\
    =&\e\left[\frac{1-p_1}{p_1}\cdot \e[e_{1,1}^2|\filt_{0}]+\frac{p_1}{1-p_1}\cdot \e[e_{1,0}^2|\filt_{0}]+2\e[e_{1,0}e_{1,1}|\filt_{0}]\right]\notag\\
    =&\e\left[9T^{1/2}\left(\frac{p_1}{1-p_1}+\frac{1-p_1}{p_1}\right)+2\cdot 9T^{1/2}\right]\notag\\
    \geq&2\cdot 9T^{1/2}+2\cdot 9T^{1/2}\quad(\text{AM-GM inequality})\notag\\
    =&36T^{1/2}\enspace.
\end{align}
On the other hand, for $2\leq t\leq T$, note that $e_{t,k}=y_t(k)-\e[y_t(k)|\filt_{t-1},\vec{U}]=U_k(1+\epsilon_t)-U_k=U_k\epsilon_t$.
Hence we can calculate the following conditional expectations:
\begin{align*}
    \e[e_{t,k}^2|\filt_{t-1},\vec{U}]=&\e[U_k^2\epsilon_t^2|\filt_{t-1},\vec{U}]=U_k^2\e[\epsilon_t^2]=U_k^2\cdot \indicator{\text{$t$ even}}\enspace,\\
    \e[e_{t,1}e_{t,0}|\filt_{t-1},\vec{U}]=&U_1U_0\cdot \indicator{\text{$t$ even}}\enspace,
\end{align*}
which, together with the law of iterated expectations, implies that
\begin{align}\label{theorem:neyman-regret-lower-bound_eq5}
    &\e\left[e_{t,1}^2\cdot  \frac{1-p_t}{p_t}+e_{t,0}^2\cdot \frac{p_t}{1-p_t}+2e_{t,0}e_{t,1}\right]\notag\\
    =&\e\left[\e\left[e_{t,1}^2\cdot  \frac{1-p_t}{p_t}+e_{t,0}^2\cdot \frac{p_t}{1-p_t}+2e_{t,0}e_{t,1}\Big|\filt_{t-1},\vec{U}\right]\right]\notag\\
    =&\e\left[ \frac{1-p_t}{p_t}\cdot \e[e_{t,1}^2|\filt_{t-1},\vec{U}]+\frac{p_t}{1-p_t}\cdot \e[e_{t,0}^2|\filt_{t-1},\vec{U}]+2\e[e_{t,0}e_{t,1}|\filt_{t-1},\vec{U}]\right]\notag\\
    =&\e\left[\indicator{\text{$t$ even}}\cdot\left(U_1^2\cdot \frac{1-p_t}{p_t}+U_0^2\cdot \frac{p_t}{1-p_t}+2U_1U_0\right)\right]\notag\\
    \geq&\e\left[\indicator{\text{$t$ even}}\cdot 4U_1U_0\right]\notag\\
    =&32\cdot \indicator{\text{$t$ even}}\quad(\text{$U_1U_0=8$ by construction})
\end{align}
Combining the results in \eqref{theorem:neyman-regret-lower-bound_eq4} and \eqref{theorem:neyman-regret-lower-bound_eq5}, we can obtain the following lower bound on $S_2$:
\begin{align}\label{theorem:neyman-regret-lower-bound_eq6}
    S_2=&\e\left[e_{1,1}^2\cdot \frac{1-p_1}{p_1}+e_{1,0}^2\cdot \frac{p_1}{1-p_1}+2e_{1,0}e_{1,1}\right]\notag\\
    &+\sum_{t=2}^T\e\left[e_{t,1}^2\cdot  \frac{1-p_t}{p_t}+e_{t,0}^2\cdot \frac{p_t}{1-p_t}+2e_{t,0}e_{t,1}\right]\notag\\
    \geq&36T^{1/2}+32\sum_{t=2}^T\indicator{\text{$t$ even}}\notag\\
    =&36T^{1/2}+16T\enspace.
\end{align}
Finally, we verify that $S_3$ equals 0.
For $t=1$, by the law of iterated expectations, we have
\begin{align}\label{theorem:neyman-regret-lower-bound_eq7}
    &\e\left[e_{1,1}r_{1,1}\cdot\frac{1-p_1}{p_1}+e_{1,0}r_{1,0}\cdot\frac{p_1}{1-p_1}+e_{1,1}r_{1,0}+e_{1,0}r_{1,1}\right]\notag\\
    =&\e\left[\e\left[e_{1,1}r_{1,1}\cdot\frac{1-p_1}{p_1}+e_{1,0}r_{1,0}\cdot\frac{p_1}{1-p_1}+e_{1,1}r_{1,0}+e_{1,0}r_{1,1}\Big|\mathcal{F}_0\right]\right]\notag\\
    =&\e\left[\left(r_{1,1}\cdot\frac{1-p_1}{p_1}+r_{1,0}\right)\e\left[e_{1,1}|\mathcal{F}_0\right]+\left(r_{1,0}\cdot\frac{p_1}{1-p_1}+r_{1,1}\right)\e\left[e_{1,0}|\mathcal{F}_0\right]\right]\notag\\
    =&0\enspace.
\end{align}
Similarly, for $2\leq t\leq T$, by the law of iterated expectations, we have
\begin{align}\label{theorem:neyman-regret-lower-bound_eq8}
    \quad&\e\left[e_{t,1}r_{t,1}\cdot\frac{1-p_t}{p_t}+e_{t,0}r_{t,0}\cdot\frac{p_t}{1-p_t}+e_{t,1}r_{t,0}+e_{t,0}r_{t,1}\right]\notag\\
    \quad=&\e\left[\e\left[e_{t,1}r_{t,1}\cdot\frac{1-p_t}{p_t}+e_{t,0}r_{t,0}\cdot\frac{p_t}{1-p_t}+e_{t,1}r_{t,0}+e_{t,0}r_{t,1}\Big|\mathcal{F}_{t-1},\vec{U}\right]\right]\notag\\
    \quad=&\e\left[\left(r_{t,1}\cdot\frac{1-p_t}{p_t}+r_{t,0}\right)\e\left[e_{t,1}|\mathcal{F}_{t-1},\vec{U}\right]+\left(r_{t,0}\cdot\frac{p_t}{1-p_t}+r_{t,1}\right)\e\left[e_{t,0}|\mathcal{F}_{t-1},\vec{U}\right]\right]\notag\\
    \quad=&0\enspace.
\end{align}
Combining the results in \eqref{theorem:neyman-regret-lower-bound_eq7} and \eqref{theorem:neyman-regret-lower-bound_eq8}, we obtain $S_3=0$.
Hence by \eqref{theorem:neyman-regret-lower-bound_eq1} and \eqref{theorem:neyman-regret-lower-bound_eq6}, we have verified that
\begin{align}\label{theorem:neyman-regret-lower-bound_eq12}
    \e_{y}\left[T\cdot \operatorname{Var}_{-y}(\hat{\tau})\right]\geq 16+36T^{-1/2}\enspace.
\end{align}
\textbf{Step 2:} Upper-bound the expected (conditional) variance under the oracle design.\\[3mm]
Using the explicit form of the oracle variance established in Proposition \mainref{prop:oracle} as well as $\rho \leq 1$ implied by the Cauchy-Schwarz inequality, we have that 
\begin{align} \label{theorem:neyman-regret-lower-bound_eq13}
    \Esub[\Big]{y}{ T \cdot \optvar_{-y} }
    &= \Esub[\Big]{y}{ 2 (1 + \rho) \olsres{1} \olsres{0} } 
    \notag \\
    &\leq \Esub[\Big]{y}{ 4 \olsres{1} \olsres{0} } 
    \notag
    \intertext{Using the law of total expectation and conditioning on the values of the vector $\vec{U}$, we decompose the expectation as}
    &= \Esub[\Big]{y}{ 4 \olsres{1} \olsres{0} \mid \vec{U} = (-4, -2) } \Pr{ \vec{U} = (-4, -2) } 
    \notag \\
    &\quad + \Esub[\Big]{y}{ 4 \olsres{1} \olsres{0} \mid \vec{U} = (2, 4) } \Pr{ \vec{U} = (2, 4) } 
    \notag
    \intertext{
        Note that $\setb{y_t(1),y_t(0)}_{t=1}^T$ has the same distribution as $\setb{-y_t(0),-y_t(1)}_{t=1}^T$ by construction, so that the two conditional expectations above are equal.
        Moreover, $\vec{U}$ takes each value with equal probability.
        Thus, the expectation is 
    }
    &= \Esub[\Big]{y}{ 4 \olsres{1} \olsres{0} \mid \vec{U} = (-4, -2) } 
    \notag \\
    \intertext{Applying the AM-GM inequality $4ab \leq a^2 + 4b^2$ yields}
    &\leq \Esub[\Big]{y}{ \olsres{1}^2 \mid \vec{U} = (-4, -2) } 
    + 4 \Esub[\Big]{y}{ \olsres{0}^2 \mid \vec{U} = (-4, -2) }
    \enspace.
\end{align}

Next, we turn to calculating the expected squared OLS residuals above under the mixture.
This is a simple task because we have only included intercept covariates, i.e. $\xv_t=1$ for all $t\in[T]$.
This allows us to compute the OLS coefficient $\bv^*(k)=\widebar{y}(k):= \frac{1}{T}\sum_{t=1}^Ty_t(k)$ for each $k\in\{0,1\}$.
Note that $y_1(k),\ldots,y_T(k)$ are conditionally independent given $\vec{U}$.
Hence, by Lemma \ref{lemma:sample-variance} we have
\begin{align}\label{theorem:neyman-regret-lower-bound_eq14}
    &\e_y[\olsres{1}^2|\vec{U}=(-4,-2)]\notag\\
    =&\e_{y}\left[\frac{1}{T}\sum_{t=1}^T\left(y_t(1)-\widebar{y}(1)\right)^2\Bigg|\vec{U}=(-4,-2)\right]\notag\\
    =&\frac{T-1}{T^2}\sum_{t=1}^T\operatorname{Var}(y_t(1)|\vec{U}=(-4,-2))+\frac{1}{T}\sum_{t=1}^T(\e[y_t(1)|\vec{U}=(-4,-2)])^2\notag\\
    &-(\e[\widebar{y}(1)|\vec{U}=(-4,-2)])^2\notag\\
    =&\frac{T-1}{T^2}\cdot \frac{4^2T}{2}+\frac{(4T^{1/4})^2}{T}+\frac{T-1}{T}\cdot 4^2-\left(\frac{1}{T}\cdot 4T^{1/4}+\frac{T-1}{T}\cdot 4\right)^2\notag\\
    =&8+16T^{-1/2}+8T^{-1}-16T^{-3/2}-16T^{-2}-32T^{-3/4}+32T^{-7/4}
\end{align}
and
\begin{align}\label{theorem:neyman-regret-lower-bound_eq15}
    &\e_y[\olsres{0}^2|\vec{U}=(-4,-2)]\notag\\
    =&2+4T^{-1/2}+2T^{-1}-4T^{-3/2}-4T^{-2}-8T^{-3/4}+8T^{-7/4}\enspace.
\end{align}
Hence by \eqref{theorem:neyman-regret-lower-bound_eq13}, \eqref{theorem:neyman-regret-lower-bound_eq14} and \eqref{theorem:neyman-regret-lower-bound_eq15}, we obtain the following upper bound:
\begin{align}\label{theorem:neyman-regret-lower-bound_eq16}
    \e_{y}[T\cdot \mathrm{V}^{*}_{-y}]\leq& 16+32T^{-1/2}+16T^{-1}-32T^{-3/2}-32T^{-2}-64T^{-3/4}+64T^{-7/4}\notag\\
    \leq&16+32T^{-1/2}+16T^{-1}-32T^{-3/2}-32T^{-2}\enspace.
\end{align}
By combining the results in \eqref{theorem:neyman-regret-lower-bound_eq12} and \eqref{theorem:neyman-regret-lower-bound_eq16}, we have
\begin{align*}
    \e_{y}\left[T\cdot \operatorname{Var}_{-y}(\hat{\tau})-T\cdot \mathrm{V}^{*}_{-y}\right]\geq& 4T^{-1/2}-16T^{-1}+32T^{-3/2}+32T^{-2}\enspace.
\end{align*}
All that remains to be shown is that this lower bound is at least $c \cdot T^{-1/2}$ for some constant $c > 0$.
To this end, observe that when $T\geq 25$, we have $\frac{16}{5}T^{-1/2}\geq 16T^{-1}$, which implies that
\begin{align*}
    4T^{-1/2}-16T^{-1}+32T^{-3/2}+32T^{-2}\geq \frac{4}{5}T^{-1/2}+32T^{-3/2}+32T^{-2}\geq \frac{4}{5}T^{-1/2}\enspace.
\end{align*}
When $2\leq T\leq 25$, direct calculation shows that this term is still strictly positive:
\begin{align*}
    4T^{-1/2}-16T^{-1}+32T^{-3/2}+32T^{-2}>0\enspace.
\end{align*}
Therefore, there exists a constant $c>0$, such that for any $T\geq 2$,
\begin{align}\label{theorem:neyman-regret-lower-bound_eq17}
    \e_{y}\left[T\cdot \operatorname{Var}_{-y}(\hat{\tau})-T\cdot \mathrm{V}^{*}_{-y}\right]\geq&c\cdot T^{-1/2}\enspace.
\end{align}
\textbf{Step 3:} Verify Assumptions \mainref{assumption:moments}-\mainref{assumption:maximum-radius}.\\[3mm]
We then prove that all configurations of $\{ y_t(1), y_t(0), \xv_t \}_{t=1}^T$ generated by this mixture satisfy Assumptions \mainref{assumption:moments}-\mainref{assumption:maximum-radius}.
Assumptions \mainref{assumption:covariate-regularity} and \mainref{assumption:maximum-radius} hold automatically under the choice $\xv_t\equiv 1$.
For any generated configuration of $\{ y_t(1), y_t(0)\}_{t=1}^T$, it can be readily shown that $\sum_{t=1}^Ty_t(1)^4$ and $\sum_{t=1}^Ty_t(0)^4$ are upper bounded by
\begin{align*}
    (4T^{1/4}+4)^4+(T-1)\cdot (4+4)^4\leq 8^{4}T+8^4T=8192T\enspace.
\end{align*}
We only verify the case where $\vec{U}=(2,4)$. 
Denote $\widebar{\epsilon}=\frac{1}{T-1}\sum_{t=2}^T\epsilon_t$, then
\begin{align*}
    \sum_{t=1}^T(y_t(1)-\iprod{\xv_t,\bv^*(1)})^2=&\sum_{t=1}^T(y_t(1)-\bv^*(1))^2\\
    \geq &\min_{\beta\in\mathbb{R}}\sum_{t=2}^T(y_t(1)-\beta)^2\\
    =&\min_{\beta\in\mathbb{R}}\sum_{t=2}^T(2+2\epsilon_t-\beta)^2\\
    =&\min_{\beta\in\mathbb{R}}\sum_{t=2}^T(2\epsilon_t-\beta)^2\quad(\text{shift by constant})\\
    =&4\left(\sum_{t=2}^T\epsilon_t^2-(T-1)\widebar{\epsilon}^2\right)\\
    =&4\left(\frac{T}{2}-(T-1)\widebar{\epsilon}^2\right)\quad(\text{since $\epsilon_t^2=1$ for $t$ even and $\epsilon_t^2=0$ for $t$ odd})\enspace.
\end{align*}
Similarly we can prove that $\sum_{t=1}^T(y_t(0)-\iprod{\xv_t,\bv^*(0)})^2\geq 16\left(T/2-(T-1)\widebar{\epsilon}^2\right)\geq 4\left(T/2-(T-1)\widebar{\epsilon}^2\right)$.
Note that by construction, it naturally holds that $|\widebar{\epsilon}|\leq T/(2(T-1))$ since $\epsilon_3=\epsilon_5=\ldots=\epsilon_{T-1}=0$.
Hence when $T\geq 3$, we have
\begin{align*}
    &\min\left\{\sum_{t=1}^T(y_t(1)-\iprod{\xv_t,\bv^*(1)})^2,\sum_{t=1}^T(y_t(0)-\iprod{\xv_t,\bv^*(0)})^2\right\}\\
    \geq&4\left(\frac{T}{2}-(T-1)\widebar{\epsilon}^2\right)\\
    \geq&4\left(\frac{T}{2}-(T-1)\cdot\frac{T^2}{4(T-1)^2}\right)\\
    \geq&4\left(\frac{T}{2}-\frac{3T}{8}\right)\\
    \geq&\frac{1}{2}T\enspace.
\end{align*}
When $T=2$, we can also derive a strictly positive lower bound for these two residuals since $y_1(1)\neq y_2(1)$ and $y_1(0)\neq y_2(0)$ by construction.
Hence Assumption \mainref{assumption:moments} is satisfied for any constructed configurations.
This, together with \eqref{theorem:neyman-regret-lower-bound_eq17}, indicates that for any adaptive design $\mathcal{D}$, there exists a certain configuration of $\{ y_t(1), y_t(0), \xv_t \}_{t=1}^T$ where the Neyman regret is at least $c\cdot T^{-1/2}$.
\end{proof}

\subsection{Necessity of Regularity Assumptions \mainref{assumption:moments}-\mainref{assumption:maximum-radius}}\label{section:B2}

In this section, we construct counterexamples to illustrate the necessity of the imposed Assumptions \mainref{assumption:moments}-\mainref{assumption:maximum-radius}.
Lemma \ref{lemma:lower-bound-1} and Lemma \ref{lemma:lower-bound-2} demonstrate the necessity of the fourth-moment condition in Assumption \mainref{assumption:moments} and the covariate regularity condition in Assumption \mainref{assumption:covariate-regularity} in order to obtain a regret lower bound of order $\Omega(T^{-1/2}R)$, respectively.

By a construction similar to that in Theorem \mainref{theorem:neyman-regret-lower-bound}, we can show the necessity of the bounded fourth-moment condition in Assumption \mainref{assumption:moments}.
This is formally stated in the following lemma.

\begin{lemma}\label{lemma:lower-bound-1}
    There exists a constant $c > 0$ such that
	for all integers $T$ and any adaptive experimental design $\fullarbdesign$,
	there exists a sequence $\{ y_t(1), y_t(0), \xv_t \}_{t=1}^T$ satisfying Assumptions \mainref{assumption:moments}-\mainref{assumption:maximum-radius} except for the bounded fourth-moment condition in Assumption \mainref{assumption:moments}, under which the corresponding Neyman regret is at least $\neymanregret_T \geq c$.
\end{lemma}

\begin{proof}
    For simplicity, we only prove the case where $d=1$ and $R=1$.
    The proof is similar to Theorem \mainref{theorem:neyman-regret-lower-bound}. 
    Let $\xv_{t}=1$ for all $t\in[T]$. 
    We construct the random sequence $\{y_t(1),y_t(0):t\in[T]\}$ as follows:
    \begin{enumerate}
        \item[1.] Generate $\epsilon_1,\ldots,\epsilon_T$ where $\epsilon_{2},\epsilon_{4},\ldots,\epsilon_{T}$ are i.i.d. Rademacher random variables and $\epsilon_1=\epsilon_3=\ldots=\epsilon_{T-1}=0$.
        \item[2.] Generate a random vector $\vec{U}=(U_1,U_0)$ that is independent of $\epsilon_1,\ldots,\epsilon_T$.
        It takes value $(2,4)$ with probability $1/2$ and takes value $(-4,-2)$ with probability $1/2$.
        \item[3.] Let $y_1(1)=U_1(T^{1/2}+\epsilon_1)$ and $y_1(0)=U_0(T^{1/2}+\epsilon_1)$.
        \item[4.] For any $t=2,\ldots,T$, let $y_t(1)=U_1(1+\epsilon_t)$ and $y_t(0)=U_0(1+\epsilon_t)$.
    \end{enumerate} 
    This construction differs from that in Theorem \mainref{theorem:neyman-regret-lower-bound} only in how $y_1(1)$ and $y_1(0)$ are specified, where their orders $\Theta(T^{1/4})$ are replaced with $\Theta(T^{1/2})$.
    The construction satisfies the bounded second moment condition, but its empirical fourth moment is of order $\Omega(T)$, which fails to satisfy the bounded fourth-moment condition in Assumption \mainref{assumption:moments}.
    By the same proof as in Theorem \mainref{theorem:neyman-regret-lower-bound}, we can prove that there exists a constant $c>0$ such that $\e_y\left[T\cdot \operatorname{Var}_{-y}(\hat{\tau})-T\cdot \mathrm{V}_{-y}^*\right]\geq c$ for any $T$.
    Moreover, it is straightforward to verify Assumptions \mainref{assumption:covariate-regularity} and \mainref{assumption:maximum-radius}. 
    This completes the proof.
\end{proof}

The following lemma states that the lower bound on the Neyman regret can also be of order $\Omega(1)$ if the covariate regularity condition (Assumption \mainref{assumption:covariate-regularity}) is not imposed.

\begin{lemma}\label{lemma:lower-bound-2}
    There exists a constant $c > 0$ such that
	for all integers $T$ and any adaptive experimental design $\fullarbdesign$,
	there exists a sequence $\{ y_t(1), y_t(0), \xv_t \}_{t=1}^T$ satisfying Assumptions \mainref{assumption:moments} and \mainref{assumption:maximum-radius}, under which the corresponding Neyman regret is at least $\neymanregret_T \geq c$.
\end{lemma}

\begin{proof}
    For simplicity, we assume that $T$ is even.
    Suppose $\xv_{t}\in\mathbb{R}^d$, where $d=T/2$. 
    We construct the random sequence $\{y_t(1),y_t(0),\xv_t\}$ as follows: $\xv_t=\vec{e}_t$ for $t=1,\ldots,T/2$ and $\xv_{T/2+1}=\ldots=\xv_{T}=\vec{0}$.
    Let $y_1(1),\ldots,y_T(1)$ be independent Rademacher random variables and let $y_t(0)=y_t(1)$ for $t=1,\ldots,T$. 
    It can be readily shown that Assumption \mainref{assumption:maximum-radius} is satisfied with $R=1$.
    Assumption \mainref{assumption:moments} also holds since $\sum_{t=1}^{T}y_t(k)^4=T$ and $\sum_{t=1}^T(y_t(k)-\iprod{\xv_t,\bv^*(k)})^2\geq \sum_{t=T/2+1}^T(y_t(k)-\iprod{\xv_t,\bv^*(k)})^2=\sum_{t=T/2+1}^Ty_t(k)^2=T/2$ for $k\in\{0,1\}$.
    Moreover, the $y$'s have marginal expectation 0.
    Now we prove that
    \begin{align}\label{lemma:lower-bound-2_eq1}
        \e_{y}\left[T\cdot \operatorname{Var}_{-y}(\hat{\tau})-T\cdot \mathrm{V}_{-y}^*\right]=\Omega(1).
    \end{align}
    By a similar argument to that in Theorem \mainref{theorem:neyman-regret-lower-bound}, we can decompose the expected conditional variance under any given adaptive design and derive the following lower bound:
    \begin{align}\label{lemma:lower-bound-2_eq2}
        &\e_{y}\left[T\cdot \operatorname{Var}_{-y}(\hat{\tau})\right]\notag\\
        =&\frac{1}{T}\sum_{t=1}^T\e\left[\braces[\big]{ y_t(1) - \iprod{ \xv_t , \bv_t(1) } }^2 \cdot  \frac{1-p_t}{p_t}+\braces[\big]{ y_t(0) - \iprod{ \xv_t , \bv_t(0) } }^2 \cdot \frac{p_t}{1-p_t} \right]\notag\\
        &+\frac{2}{T}\sum_{t=1}^T \e\left[ \braces[\big]{ y_t(1) - \iprod{ \xv_t , \bv_t(1) } }\cdot\braces[\big]{ y_t(0) - \iprod{ \xv_t , \bv_t(0) } } \right]\quad(\text{Proposition \mainref{prop:aipw-variance}})\notag\\
        =&\frac{1}{T}\sum_{t=1}^T\e\left[\left(\iprod{\xv_t,\bv_t(1)}\cdot \sqrt{\frac{1-p_t}{p_t}}+\iprod{\xv_t,\bv_t(0)}\cdot \sqrt{\frac{p_t}{1-p_t}}\right)^2\right]\notag\\
        &+\frac{1}{T}\sum_{t=1}^T\e\left[y_t(1)^2\cdot  \frac{1-p_t}{p_t}+y_t(0)^2\cdot \frac{p_t}{1-p_t}+2y_t(1)y_t(0)\right]\notag\\
        &-\frac{2}{T}\sum_{t=1}^T\e\Bigg[y_t(1)\left(\iprod{\xv_t,\bv_t(1)}\cdot\frac{1-p_t}{p_t}+\iprod{\xv_t,\bv_t(0)}\right)\notag\\
        &\quad\quad\quad\quad+y_t(0)\left(\iprod{\xv_t,\bv_t(0)}\cdot\frac{p_t}{1-p_t}+\iprod{\xv_t,\bv_t(1)}\right)\Bigg]\notag\\
        \intertext{The last term equals 0 since $y_t(1),y_t(0)$ have mean zero and are independent of $p_t,\bv_t(1),\bv_t(0)$.}
        \geq&\frac{1}{T}\sum_{t=1}^T\e\left[\frac{1-p_t}{p_t}+\frac{p_t}{1-p_t}+2\right]\quad(\text{since $y_t(1)^2=y_t(0)^2=1$ for all $t\in[T]$})\notag\\
        \geq&4\quad(\text{AM-GM inequality})\enspace.
    \end{align}
    Denote $\vec{Y}=(y_1(1),\ldots,y_T(1))^{\tran}$ and $\vec{X}=(\xv_1,\ldots,\xv_T)^{\tran}$.
    By a similar argument as in Theorem \mainref{theorem:neyman-regret-lower-bound}, we can derive the following upper bound:
    \begin{align}\label{lemma:lower-bound-2_eq3}
        \e_{y}[T\cdot \mathrm{V}^{*}_{-y}]\leq&4\e_{y}\left[\olsres{1}\olsres{0}\right]\quad(\text{similar proof as in Theorem \mainref{theorem:neyman-regret-lower-bound}})\notag\\
        =&4\e_{y}\left[\olsres{1}^2\right]\quad(\text{since $y_t(1)=y_t(0)$ for all $t\in[T]$})\notag\\
        =&\frac{4}{T}\e_y\left[\vec{Y}^{\tran}\left(\vec{I}_T-\vec{X}(\vec{X}^{\tran}\vec{X})^{-1}\vec{X}^{\tran}\right)\vec{Y}\right]\notag\\
        =&\frac{4}{T}\e_y\left[\sum_{t=1}^Ty_t(1)^2\right]-\frac{4}{T}\e_y\left[\operatorname{tr}(\vec{Y}^{\tran}\vec{X}(\vec{X}^{\tran}\vec{X})^{-1}\vec{X}^{\tran}\vec{Y})\right]\notag\\
        =&4-\frac{4}{T}\operatorname{tr}(\vec{X}(\vec{X}^{\tran}\vec{X})^{-1}\vec{X}^{\tran}\e_y[\vec{Y}\vec{Y}^{\tran}])\notag\\
        =&4-\frac{4}{T}\operatorname{tr}(\vec{X}(\vec{X}^{\tran}\vec{X})^{-1}\vec{X}^{\tran}\vec{I}_T)\notag\\
        =&4-\frac{4}{T}\operatorname{tr}((\vec{X}^{\tran}\vec{X})^{-1}\vec{X}^{\tran}\vec{X})\notag\\
        =&4-\frac{4}{T}\operatorname{tr}(\vec{I}_d)\notag\\
        =&2\quad(\text{since $d=T/2$})\enspace.
    \end{align}
    Hence \eqref{lemma:lower-bound-2_eq1} follows from \eqref{lemma:lower-bound-2_eq2} and \eqref{lemma:lower-bound-2_eq3}.
\end{proof}

Finally, we remark on the necessity of Assumption \mainref{assumption:maximum-radius}.
Consider the case where every covariate is of constant order, i.e. $\xv_{t,j} = \bigO{1}$, so that $R = \bigO{d^{1/2}}$.
In this case, the invertibility stipulated by Assumption \mainref{assumption:covariate-regularity}, i.e. that $\frac{1}{t} \xM_t^\tran \xM_t \succ 0$ for $t \geq \bigOmega{T^{1/2}d^{1/2}}$, requires that $d \leq \bigO{T}$, which is the content of Assumption \mainref{assumption:maximum-radius}.
In other words, when $R = \bigO{d^{1/2}}$ then Assumption \mainref{assumption:maximum-radius} is implied by Assumption \mainref{assumption:covariate-regularity}.
In this sense, Assumption \mainref{assumption:maximum-radius} is necessary for Assumption \mainref{assumption:covariate-regularity} to hold in the case where $\xv_{t,j} = \bigO{1}$.

\section{Neyman Regret Upper Bound} \label{sec:supp-neyman-regret-analysis}
In this section, we aim to perform a comprehensive analysis of the upper bound on the Neyman regret.
In Section \ref{section:C1}, we provide a proof of Lemma \mainref{lemma:regret-decomposition}, which decomposes the Neyman regret into two components: the expected probability regret and the expected prediction regret.
In Section \ref{section:C2}, we provide an upper bound on the expected probability regret.
Section \ref{section:C3} establishes an upper bound on the moments of inverse probabilities.
In Section \ref{section:C4}, we provide an upper bound on the expected prediction regret.
Section \ref{section:C5} establishes the upper bound on the fourth moment of the online residuals, a crucial component in upper-bounding the expected probability regret.
By combining the results from these sections, we establish the main theorem on Neyman regret (Theorem \mainref{thm:neyman-regret}) in Section \ref{section:C6}.

\subsection{Regret Decomposition}\label{section:C1}
In this section, we derive a decomposition of the Neyman regret into two components: the expected probability regret and the expected prediction regret. 
We start by deriving the explicit form of the oracle variance in the following proposition.

\begin{refproposition}{\mainref{prop:oracle}}
    \oracle
\end{refproposition}

\begin{proof}
    For $k\in\{0,1\}$, denote $\Delta_{t,k}(\bv)=y_t(k)-\iprod{\xv_t,\bv}$.
    By direct calculation, the explicit form of $T^2\cdot\operatorname{Var}(\hat{\tau};\bv(1),\bv(0),p)$ is given by:
    \begin{align*}
	    T^2\cdot\operatorname{Var}(\hat{\tau};\bv(1),\bv(0),p)&= \left(\sum_{t=1}^T\Delta_{t,1}(\bv(1))^2\right) \cdot \paren[\Big]{ \frac{1}{p} - 1}  
	    + \left(\sum_{t=1}^T \Delta_{t,0}(\bv(0))^2\right) \cdot \paren[\Big]{ \frac{1}{1-p} - 1}\\
	    &+ 2\sum_{t=1}^T \Delta_{t,1}(\bv(1)) \cdot \Delta_{t,0}(\bv(0))  
	    \enspace.
    \end{align*}
    The subsequent proof consists of two steps.
    We first fix $\bv(1)$ and $\bv(0)$ and derive the optimal choice of $p$.
    Then we prove that $\bv^*(1)$ and $\bv^*(0)$ minimize the resulting variance.\\[3mm]
    \textbf{Step 1: }For fixed $\bv(1)$ and $\bv(0)$, obtain the optimal $p$ and its corresponding variance.\\[3mm]
    By first-order optimality conditions, it can be readily shown that the optimal choice of $p$ (as a function of $\bv(1)$, $\bv(0)$ and the data) is
    \begin{align*}
        p^*(\bv(1),\bv(0))=\left(1+\sqrt{\frac{\sum_{t=1}^T\Delta_{t,0}(\bv(0))^2}{\sum_{t=1}^T\Delta_{t,1}(\bv(1))^2}}\right)^{-1}\enspace.
    \end{align*}
    We omit the dependence of $p^*$ on $(\bv(1),\bv(0))$ later in the proof for notational simplicity.
    By plugging this explicit form of $p^*$ into the variance $\operatorname{Var}(\hat{\tau};\bv(1),\bv(0),p)$, we can obtain the following expression:
    \begin{align*}
        T^2\cdot\operatorname{Var}(\hat{\tau};\bv(1),\bv(0),p^*)=&2\left(\sum_{t=1}^T \Delta_{t,1}(\bv(1))^2\right)^{1/2}\left(\sum_{t=1}^T \Delta_{t,0}(\bv(0))^2\right)^{1/2}\\
        &+2\sum_{t=1}^T \Delta_{t,1}(\bv(1)) \cdot \Delta_{t,0}(\bv(0))\enspace.
    \end{align*}
    \textbf{Step 2: }Prove that $\bv^*(1)=\bv_{\textrm{OLS}}(1)$, $\bv^*(0)=\bv_{\textrm{OLS}}(0)$ minimize $\operatorname{Var}(\hat{\tau};\bv(1),\bv(0),p^*)$.\\[3mm]
    For any fixed $\bv(1),\bv(0)$, we show that $\operatorname{Var}(\hat{\tau};\bv(1),\bv(0),p^*)\geq \operatorname{Var}(\hat{\tau};\bv^*(1),\bv^*(0),p^*)$.
    For $k\in\{0,1\}$, we introduce the following notations:
    \begin{align*}
        \vec{\delta}_k=&\bv(k)-\bv^*(k)\enspace,\\
        \vec{Y}_k=&(y_1(k),\ldots,y_T(k))^{\tran}\enspace,\\
        \vec{X}=&(\xv_1,\ldots,\xv_T)^{\tran}\enspace,\\
        \vec{\alpha}_k=&\vec{Y}_k-\vec{X}\bv^*(k)\enspace,\\
        \vec{z}_k=&\vec{X}\vec{\delta}_k\enspace.
    \end{align*} 
    By the explicit form of $\bv^*(1)$ and $\bv^*(0)$, for any $k_1,k_2\in\{0,1\}$ we can see that
    \begin{align*}
        \iprod{\vec{\alpha}_{k_1},\vec{z}_{k_2}}=&\vec{\delta}_{k_2}^{\tran}\vec{X}^{\tran}\left(\vec{Y}_{k_1}-\vec{X}\bv^*(k_1)\right)\\
        =&\vec{\delta}_{k_2}^{\tran}\vec{X}^{\tran}\left(\vec{I}_T-\vec{X}(\vec{X}^{\tran}\vec{X})^{-1}\vec{X}^{\tran}\right)\vec{Y}_{k_1}\\
        =&0\enspace.
    \end{align*}
    By this orthogonality, we can simplify the expressions for $\operatorname{Var}(\hat{\tau};\bv(1),\bv(0),p^*)$ and $\operatorname{Var}(\hat{\tau};\bv^*(1),\bv^*(0),p^*)$ as:
    \begin{align*}
        T^2\cdot\operatorname{Var}(\hat{\tau};\bv(1),\bv(0),p^*)=&2\left\|\vec{Y}_1-\vec{X}\bv(1)\right\|\left\|\vec{Y}_0-\vec{X}\bv(0)\right\|+2\iprod{\vec{Y}_1-\vec{X}\bv(1),\vec{Y}_0-\vec{X}\bv(0)}\\
        =&2\left(\|\vec{\alpha}_1-\vec{z}_1\|^2\right)^{1/2}\left(\|\vec{\alpha}_0-\vec{z}_0\|^2\right)^{1/2}+2\iprod{\vec{\alpha}_1-\vec{z}_1,\vec{\alpha}_0-\vec{z}_0}\\
        =&2\left(\|\vec{\alpha}_1\|^2+\|\vec{z}_1\|^2\right)^{1/2}\left(\|\vec{\alpha}_0\|^2+\|\vec{z}_0\|^2\right)^{1/2}+2\iprod{\vec{\alpha}_1,\vec{\alpha}_0}+2\iprod{\vec{z}_1,\vec{z}_0}\enspace,\\
        T^2\cdot\operatorname{Var}(\hat{\tau};\bv^*(1),\bv^*(0),p^*)=&2\left\|\vec{Y}_1-\vec{X}\bv^*(1)\right\|\left\|\vec{Y}_0-\vec{X}\bv^*(0)\right\|+2\iprod{\vec{Y}_1-\vec{X}\bv^*(1),\vec{Y}_0-\vec{X}\bv^*(0)}\\
        =&2\left(\|\vec{\alpha}_1\|^2\right)^{1/2}\left(\|\vec{\alpha}_0\|^2\right)^{1/2}+2\iprod{\vec{\alpha}_1,\vec{\alpha}_0}\enspace.
    \end{align*}
    In order to verify $\operatorname{Var}(\hat{\tau};\bv(1),\bv(0),p^*)\geq \operatorname{Var}(\hat{\tau};\bv^*(1),\bv^*(0),p^*)$, it suffices to prove that
    \begin{align*}
        \left(\|\vec{\alpha}_1\|^2+\|\vec{z}_1\|^2\right)^{1/2}\left(\|\vec{\alpha}_0\|^2+\|\vec{z}_0\|^2\right)^{1/2}+\iprod{\vec{z}_1,\vec{z}_0}\geq\left(\|\vec{\alpha}_1\|^2\right)^{1/2}\left(\|\vec{\alpha}_0\|^2\right)^{1/2}\enspace.
    \end{align*}
    Note that it holds that $\iprod{\vec{z}_1,\vec{z}_0}\geq -\|\vec{z}_1\|\|\vec{z}_0\|$ by the Cauchy-Schwarz inequality. 
    Thus, it suffices to show that
    \begin{align}\label{prop:oracle_eq1}
        &\left(\|\vec{\alpha}_1\|^2+\|\vec{z}_1\|^2\right)^{1/2}\left(\|\vec{\alpha}_0\|^2+\|\vec{z}_0\|^2\right)^{1/2}\notag\\
        \geq&\left(\|\vec{\alpha}_1\|^2\right)^{1/2}\left(\|\vec{\alpha}_0\|^2\right)^{1/2}+\|\vec{z}_1\|\|\vec{z}_0\|\enspace.
    \end{align}
    Since both sides of \eqref{prop:oracle_eq1} are nonnegative, it is equivalent to 
    \begin{align*}
        &\left(\|\vec{\alpha}_1\|^2+\|\vec{z}_1\|^2\right)\left(\|\vec{\alpha}_0\|^2+\|\vec{z}_0\|^2\right)\\
        \geq& \left(\|\vec{\alpha}_1\|^2\right)\left(\|\vec{\alpha}_0\|^2\right)+2(\|\vec{\alpha}_1\|^2\|\vec{\alpha}_0\|^2\|\vec{z}_1\|^2\|\vec{z}_0\|^2)^{1/2}+\left(\|\vec{z}_1\|^2\right)\left(\|\vec{z}_0\|^2\right)\enspace,
    \end{align*}
    which can be further simplified as
    \begin{align*}
        \|\vec{\alpha}_1\|^2\|\vec{z}_0\|^2+\|\vec{\alpha}_0\|^2\|\vec{z}_1\|^2\geq 2(\|\vec{\alpha}_1\|^2\|\vec{\alpha}_0\|^2\|\vec{z}_1\|^2\|\vec{z}_0\|^2)^{1/2}\enspace.
    \end{align*}
    This inequality is immediately verified by the AM-GM inequality, which implies that $(\bv^*(1),\bv^*(0))$ minimize $\operatorname{Var}(\hat{\tau};\bv(1),\bv(0),p^*)$. 
    By plugging $\bv(1)=\bv^*(1)$ and $\bv(0)=\bv^*(0)$ into the form of $p^*$, we can obtain the minimizer for $p$ as $p^*=(1+\mathcal{E}(0)/\mathcal{E}(1))^{-1}$.
    The corresponding oracle Neyman variance can thus be calculated.
\end{proof}

Based on the explicit form of the oracle Neyman variance in Proposition \mainref{prop:oracle}, we derive the following lemma on the decomposition of Neyman regret.

\begin{reflemma}{\mainref{lemma:regret-decomposition}}
	\regretdecomposition
\end{reflemma}

\begin{proof}
	For $k\in\{0,1\}$, denote $\Delta_{t,k}(\bv)=y_t(k)-\iprod{\xv_t,\bv}$.
    Let us recall the expression for the variance of the AIPW estimator from Proposition \mainref{prop:aipw-variance}:
    \begin{align}\label{lemma:regret-decomposition_eq1}
        &T\cdot \Var{\eate}\notag\\
        =&\e\left[\frac{1}{T}\sum_{t=1}^T\Delta_{t,1}(\bv_t(1))^2\cdot \left(\frac{1}{p_t}-1\right)\right]+\e\left[\frac{1}{T}\sum_{t=1}^T\Delta_{t,0}(\bv_t(0))^2\cdot \left(\frac{1}{1-p_t}-1\right)\right]\notag\\
        &+2\e\left[\frac{1}{T}\sum_{t=1}^T\Delta_{t,1}(\bv_t(1))\cdot\Delta_{t,0}(\bv_t(0))\right]\quad(\text{Proposition \mainref{prop:aipw-variance}})\notag\\
        =&\e\left[\frac{1}{T}\sum_{t=1}^T\left(f_t(p_t)-\Delta_{t,1}(\bv_t(1))^2-\Delta_{t,0}(\bv_t(0))^2\right)\right]+2\e\left[\frac{1}{T}\sum_{t=1}^T\Delta_{t,1}(\bv_t(1))\cdot\Delta_{t,0}(\bv_t(0))\right]\notag\\
        =&\underbrace{\e\left[\frac{1}{T}\sum_{t=1}^Tf_t(p^*)-\Delta_{t,1}(\bv_t(1))^2-\Delta_{t,0}(\bv_t(0))^2\right]+2\e\left[\frac{1}{T}\sum_{t=1}^T\Delta_{t,1}(\bv_t(1))\cdot\Delta_{t,0}(\bv_t(0))\right]}_{:=S_1}\notag\\
        &+\frac{1}{T}\E{\regretprob_T}\enspace.
    \end{align}
    By the definition of prediction regret, we can rewrite $S_1$ by adding and subtracting oracle prediction residuals:
    \begin{align}\label{lemma:regret-decomposition_eq2}
        S_1=&\e\left[\frac{1}{T}\sum_{t=1}^T\Delta_{t,1}(\bv_t(1))^2\right]\cdot \left(\frac{1}{p^*}-1\right)+\e\left[\frac{1}{T}\sum_{t=1}^T\Delta_{t,0}(\bv_t(0))^2\right]\cdot \left(\frac{1}{1-p^*}-1\right)\notag\\
        &+2\e\left[\frac{1}{T}\sum_{t=1}^T\Delta_{t,1}(\bv_t(1))\cdot \Delta_{t,0}(\bv_t(0))\right]\notag\\
        =&\frac{1}{T}\sum_{t=1}^T\Delta_{t,1}(\bv^*(1))^2\cdot \left(\frac{1}{p^*}-1\right)+\frac{1}{T}\sum_{t=1}^T\Delta_{t,0}(\bv^*(0))^2\cdot \left(\frac{1}{1-p^*}-1\right)\notag\\
        &+\frac{2}{T}\sum_{t=1}^T\Delta_{t,1}(\bv^*(1))\cdot\Delta_{t,0}(\bv^*(0))\notag\\
        &+\left(\frac{1}{p^*}-1\right)\cdot \e\left[\frac{1}{T}\sum_{t=1}^T\Delta_{t,1}(\bv_t(1))^2-\frac{1}{T}\sum_{t=1}^T\Delta_{t,1}(\bv^*(1))^2\right]\notag\\
        &+\left(\frac{1}{1-p^*}-1\right)\cdot \e\left[\frac{1}{T}\sum_{t=1}^T\Delta_{t,0}(\bv_t(0))^2-\frac{1}{T}\sum_{t=1}^T\Delta_{t,0}(\bv^*(0))^2\right]\notag\\
        &+2\e\left[\frac{1}{T}\sum_{t=1}^T\Delta_{t,1}(\bv_t(1))\cdot \Delta_{t,0}(\bv_t(0))-\frac{1}{T}\sum_{t=1}^T\Delta_{t,1}(\bv^*(1))\cdot\Delta_{t,0}(\bv^*(0))\right]\notag\\
        \intertext{Since we have $p^*=(1+\mathcal{E}(0)/\mathcal{E}(1))^{-1}$ in Proposition \mainref{prop:oracle}, it holds that $1/p^*-1=\mathcal{E}(0)/\mathcal{E}(1)$ and $1/(1-p^*)-1=\mathcal{E}(1)/\mathcal{E}(0)$. By using the definition of $\rho,\mathcal{E}(0),\mathcal{E}(1)$, we obtain}
        =&\olsres{1}^2\frac{\mathcal{E}(0)}{\mathcal{E}(1)}+\olsres{0}^2\frac{\mathcal{E}(1)}{\mathcal{E}(0)}+2\rho\mathcal{E}(1)\mathcal{E}(0) \notag\\
        &+\frac{\mathcal{E}(0)}{\mathcal{E}(1)}\cdot\e\left[\frac{1}{T}\sum_{t=1}^T\Delta_{t,1}(\bv_t(1))^2-\frac{1}{T}\sum_{t=1}^T\Delta_{t,1}(\bv^*(1))^2\right]\notag\\
        &+\frac{\mathcal{E}(1)}{\mathcal{E}(0)}\cdot \e\left[\frac{1}{T}\sum_{t=1}^T\Delta_{t,0}(\bv_t(0))^2-\frac{1}{T}\sum_{t=1}^T\Delta_{t,0}(\bv^*(0))^2\right]\notag\\
        &+2\e\left[\frac{1}{T}\sum_{t=1}^T\Delta_{t,1}(\bv_t(1))\cdot \Delta_{t,0}(\bv_t(0))-\frac{1}{T}\sum_{t=1}^T\Delta_{t,1}(\bv^*(1))\cdot\Delta_{t,0}(\bv^*(0))\right]\notag\\
        =&2(1+\rho)\mathcal{E}(1)\mathcal{E}(0)+\frac{1}{T}\e[\mathcal{R}_T^{\text{pred}}]\quad(\text{by the definition of $\mathcal{R}_T^{\text{pred}}$})\enspace,
    \end{align}
    which is equal to $T\cdot \optvar + \frac{1}{T}\e[\mathcal{R}_T^{\text{pred}}]$ by the explicit form of the oracle Neyman variance in Proposition \mainref{prop:oracle}.
    Combining \eqref{lemma:regret-decomposition_eq1} and \eqref{lemma:regret-decomposition_eq2}, the Neyman regret can finally be decomposed as:
    \begin{align*}
        \neymanregret_T=T\cdot \Var{\eate}-T\cdot \optvar=\frac{1}{T}\E{\regretprob_T}+\frac{1}{T}\E{\mathcal{R}_T^{\text{pred}}}\enspace,
    \end{align*}
    which completes the proof.
\end{proof}

\subsection{Probability Regret}\label{section:C2}
In this section, we derive an upper bound on the expected probability regret. The main result is presented in Lemma \ref{lemma:prob-regret-bound}. The bound depends on the fourth moments of the online residuals, whose uniform upper bounds will be established later in Section \ref{section:C5}.

The following lemma is a direct implication of the adaptive sequential design.

\begin{reflemma}{\mainref{lemma:unbiased-sigmoid-loss}}
	\unbiasedsigmoidloss
\end{reflemma}

\begin{proof}
    For $k\in\{0,1\}$, denote $\Delta_{t,k}(\bv)=y_t(k)-\iprod{\xv_t,\bv}$.
    Since $\bv_t(1)$, $\bv_t(0)$, $p_{t}$, and $u$ are measurable with respect to $\filt_{t-1}$, we obtain
    \begin{align*}
        &\E{\estsigloss{t}(u)\mid \filt_{t-1}}\\
        =& \e\left[\Delta_{t,1}(\bv_t(1))^2 \cdot \frac{\indicator{Z_t = 1}}{p_t} \cdot \frac{1}{\phi(u)} 
        +
        \Delta_{t,0}(\bv_t(0))^2 \cdot \frac{\indicator{Z_t = 0}}{1-p_t} \cdot \frac{1}{1-\phi(u)} \,\Big|\filt_{t-1}\right]\\
        =&\Delta_{t,1}(\bv_t(1))^2 \cdot \frac{1}{\phi(u)} 
        +
        \Delta_{t,0}(\bv_t(0))^2 \cdot \frac{1}{1-\phi(u)}\\
        =&\sigloss{t}(u)\enspace,
    \end{align*}
    which completes the proof.
\end{proof}

Using Lemma \mainref{lemma:unbiased-sigmoid-loss}, we can rewrite the expected probability regret in an equivalent form.
The resulting expression makes its connection with the proposed design explicit.

\begin{refcorollary}{\mainref{corollary:exp-prob-loss}}
	\expprobloss
\end{refcorollary}

\begin{proof}
    By Lemma \mainref{lemma:unbiased-sigmoid-loss} and the law of iterated expectations, we have
    \begin{align*}
        \E{\regretprob_T} =&\sum_{t=1}^T\e\left[f_t(p_t)-f_t(p^*)\right]\notag\\
        =&\sum_{t=1}^T\e\left[h_t(u_t)-h_t(u^*)\right]\quad(\text{sigmoid transformation})\notag\\
        =&\sum_{t=1}^T\e\left[\e\left[\widehat{h}_t(u_t)-\widehat{h}_t(u^*)\big|\filt_{t-1}\right]\right]\quad(\text{Lemma \mainref{lemma:unbiased-sigmoid-loss}})\notag\\
        =&\sum_{t=1}^T\e\left[\widehat{h}_t(u_t)-\widehat{h}_t(u^*)\right]\quad(\text{law of iterated expectations})\enspace,
    \end{align*}
    which completes the proof.
\end{proof}

By applying Corollary \ref{corollary:standard-FTRL-upper-bound}, we can establish an upper bound on the expected probability regret in the following lemma.

\begin{reflemma}{\mainref{lemma:exp-prob-regret-initial-oco-bound}}
	\expproblossub
\end{reflemma}

\begin{proof}
    By the definition of $\widehat{h}_t$, $\widetilde{H}_{t+1}$, $u_t$ and $\widetilde{u}_{t+1}$, we apply Corollary \ref{corollary:standard-FTRL-upper-bound} to obtain
    \begin{align*}
        \E{\regretprob_T}=&\e\left[\sum_{t=1}^T\widehat{h}_t(u_t)-\sum_{t=1}^T\widehat{h}_t(u^*)\right]\quad(\text{Corollary \mainref{corollary:exp-prob-loss}})\\
        \leq&\frac{1}{\eta_{T+1}} \psi(u^*)
		+ \sum_{t=1}^T \E[\Big]{ \ptprobloss{t+1}(u_t) - \ptprobloss{t+1}( \widetilde{u}_{t+1} ) }\enspace.
    \end{align*}
    Hence, the result is proved.
\end{proof}

The upper bound on the expected probability regret established in Lemma \mainref{lemma:exp-prob-regret-initial-oco-bound} consists of two components.
We first derive an upper bound for the first term in the following lemma.

\begin{lemma}\label{lemma:regularizer-upper-bound}
    Under Assumption \mainref{assumption:moments} and Condition \mainref{condition:sigmoid}, it holds that:
    \begin{align*}
        \psi(u^*)\leq \frac{8}{\bb_3^3}\left[\left(\frac{c_1}{c_0}-1\right)^3+\frac{\bb_3}{4}\left(\frac{c_1}{c_0}-1\right)^2\right]\enspace.
    \end{align*}
\end{lemma}

\begin{proof}
    Since $p^*=(1+\mathcal{E}(0)/\mathcal{E}(1))^{-1}$, by Assumption \mainref{assumption:moments} we have $p^*\in[\frac{1}{1+c_1/c_0},1-\frac{1}{1+c_1/c_0}]$.
    For $p^*\geq 1/2$ ($u^*=\phi^{-1}(p^*)\geq 0$), we have
    \begin{align*}
        \frac{1}{1-p^*}-\frac{1}{1-1/2}\leq\frac{1}{1-\left(1-\frac{1}{1+c_1/c_0}\right)}-2=\frac{c_1}{c_0}-1\enspace.
    \end{align*}
    Hence by the convexity property in Condition \mainref{condition:sigmoid}-2 and Lemma \ref{lemma:u}-1, we have
    \begin{align*}
        \frac{c_1}{c_0}-1\geq& \frac{1}{1-p^*}-\frac{1}{1-1/2}\\
        =&\frac{1}{1-\phi(u^*)}-\frac{1}{1-\phi(0)}\quad(\text{sigmoid transformation})\\
        \geq&\left(\frac{1}{1-\phi(0)}\right)^{\prime}(u^*-0)\quad(\text{$u^*\geq 0$ and convexity of $1/(1-\phi(u))$})\\
        \geq& \frac{\bb_3}{2}(u^*-0)\quad(\text{Lemma \ref{lemma:u}-1})\enspace,
    \end{align*}
    which implies that $0\leq u^*\leq \frac{2}{\bb_3}(\frac{c_1}{c_0}-1)$.
    Since $\psi(u)$ is increasing on $[0,\infty)$, we can obtain the following upper bound:
    \begin{align*}
        \psi(u^*)\leq \frac{1}{2}\left[\frac{2}{\bb_3}\left(\frac{c_1}{c_0}-1\right)\right]^2+\left[\frac{2}{\bb_3}\left(\frac{c_1}{c_0}-1\right)\right]^3\enspace.
    \end{align*}
    For $p^*\leq 1/2$, we can obtain the same upper bound by symmetry $\psi(u)=\psi(-u)$ and $\phi^{-1}(p)=-\phi^{-1}(1-p)$.
\end{proof}

Bounding the second term in Lemma \mainref{lemma:exp-prob-regret-initial-oco-bound} requires a detailed characterization of the curvature of the Bregman divergence induced by the regularizer, as discussed in the main paper.
To this end, the following lemma derives a lower bound on the Bregman divergence.

\begin{reflemma}{\mainref{lemma:breg-lb}}
	\breglb
\end{reflemma}

\begin{proof}
	We shall prove the stronger result that for all $u, v \in \Reals$:
	\[
	\breg{\psi}{v}{u} \geq \frac{1}{2}(v-u)^2 (1 + \frac{1}{2} \abs{v} + \abs{u})
	\enspace.
	\]
    Recall that $\psi(u)=\frac{1}{2} u^2+|u|^3$. Using the shorthand $h(u)=\frac{1}{2} u^2$ and $g(u)=|u|^3$, we have that $\psi(u)=h(u)+g(u)$. Using the linearity property, we can decompose the Bregman divergence as
    \begin{align*}
        \mathcal{B}_{\psi}(v|u)=\mathcal{B}_{h+g}(v|u)=\mathcal{B}_h(v|u)+\mathcal{B}_g(v|u)=\frac{1}{2}(v-u)^2+\mathcal{B}_g(v|u)\enspace.
    \end{align*}
    To complete the proof, it suffices to lower bound the second Bregman divergence as
    \begin{align}\label{lemma:breg-lb_eq1}
        \mathcal{B}_g(v|u) \geq \frac{1}{4}(v-u)^2 \cdot(|v|+2|u|)\enspace.
    \end{align}
    To this end, we derive the Bregman divergence corresponding to $g$ as
    \begin{align*}
        \mathcal{B}_g(v|u)=|v|^3-|u|^3-3|u| u(v-u)=|v|^3+2|u|^3-3 u v|u|\enspace.
    \end{align*}
    To show the inequality \eqref{lemma:breg-lb_eq1} holds, we will consider two cases.\\[3mm]
    \textbf{Case 1: $uv \geq 0$.} In this case, it is true that $uv= |u||v|$. Thus, we have that the Bregman divergence is given by
    \[
        \mathcal{B}_{g}(v|u) = |v|^3 + 2|u|^3 - 3|u|^2|v| \, \enspace.
    \]
    On the other hand, we have
    \begin{align*}
        (v-u)^2(|v|+2|u|)=&\left(v^2+u^2-2 v u\right)(|v|+2|u|)\\
        =&|v|^3+2 v^2|u|+u^2|v|+2\left|u\right|^3-2 v u(|v|+2|u|)\\
        =&|v|^3+2 v^2|u|+u^2|v|+2\left|u\right|^3-2|v||u|(|v|+2|u|)\quad(\text {because } u v\geq 0)\\
        =&|v|^3+2|u|^3-3 |u|^2|v|\enspace.
    \end{align*}
    Thus, this establishes that inequality \eqref{lemma:breg-lb_eq1} holds in this case, as we have
    \[
        \mathcal{B}_{g}(v|u) = (v - u)^2(|v| + 2|u|) \ge \frac14 (v - u)^2(|v| + 2|u|)\enspace.
    \]
    \textbf{Case 2: $uv < 0$.} This means that $vu = -|v||u|$. Before continuing, we establish a handy upper bound, which is
    \begin{align}\label{lemma:breg-lb_eq2}
        4v^2|u| \le 2\big(|v|^3 + 2|u|^3 + |v|u^2\big)\enspace.
    \end{align}
    This upper bound is verified through the AM-GM inequality, as
    \[
        4v^2|u| = 4|v|^{3/2} \cdot (|v|^{1/2}|u|) \le 2\big(|v|^3 + |v||u|^2\big) \le 2\big(|v|^3 + 2|u|^3 + |v|u^2\big)\enspace.
    \]
    With that inequality established, let us turn our attention back to the Bregman divergence.
    In this case, the Bregman divergence is
    \begin{align*}
        \mathcal{B}_{g}(v|u) &= |v|^3 + 2|u|^3 - 3uv|u| \\
        &= |v|^3 + 2|u|^3 + 3|u||v||u| \qquad \text{(because $uv < 0$)} \\
        &= |v|^3 + 2|u|^3 + 3u^2|v|\enspace.
    \end{align*}
    On the other hand, we have 
    \begin{align*}
        \frac14 (v - u)^2(|v| + 2|u|) =& \frac14 \big( v^2 + u^2 - 2vu \big)(|v| + 2|u|)\\
        =&\frac14(v^2 + u^2 + 2|v||u|)(|v| + 2|u|)\\
        =&\frac{1}{4}\left\{|v|^3+2 v^2|u|+2|u|^3+u^2|v|+2 v^2|u|+4|v| u^2\right\}\\ 
        =&\frac{1}{4}\left\{|v|^3+2|u|^3+4 v^2|u|+5|v| u^2\right\}\\
        \leq&\frac{1}{4}\left\{|v|^3+2|u|^3+2\left(|v|^3+2|u|^3+|v| u^2\right)+5|v| u^2\right\}\quad(\text{use \eqref{lemma:breg-lb_eq2}})\\ 
        =&\frac{1}{4}\left\{3|v|^3+6|u|^3+7|v| u^2\right\}\\
        \leq&|v|^3+2|u|^3+3|v| u^2\\
        =&\mathcal{B}_g(v|u)\enspace.
    \end{align*}
    Therefore, inequality \eqref{lemma:breg-lb_eq1} holds for all $u,v\in\mathbb{R}$, which completes the proof.
\end{proof}

Based on Lemma~\mainref{lemma:breg-lb} we derive the following upper bound for the second term in the regret upper bound of Lemma \mainref{lemma:exp-prob-regret-initial-oco-bound}.

\begin{lemma}\label{lemma:prob-regret-initial-upper-bound}
	The expected probability regret can be bounded by
	\[
	\sum_{t=1}^T \E[\Big]{ \ptprobloss{t+1}(u_t) - \ptprobloss{t+1}( \widetilde{u}_{t+1} )}
	\leq 
    \sum_{t=1}^T \eta_t \E[\Bigg]{ \frac{ \paren[\big]{\grad \estsigloss{t}(u_t)}^2 }{2(1+|u_t|)} }
	\enspace.
	\] 
\end{lemma}

\begin{proof}
	Using first order convexity conditions, rearranging terms, and the fact that $u_t$ minimizes $\tprobloss{t}$ so that $\grad \tprobloss{t}(u_t) = 0$, we have that
	\begin{align*}
		& \ptprobloss{t+1}(u_t) - \ptprobloss{t+1}( \widetilde{u}_{t+1} ) \\
		=& \sum_{s=1}^{t} \paren[\Big]{\estsigloss{s}( u_t ) - \estsigloss{s}( \widetilde{u}_{t+1} )}
		+ \eta_{t}^{-1} \paren[\Big]{ \psi(u_t) - \psi(\widetilde{u}_{t+1}) } \\
		\leq& \iprod[\Big]{ \sum_{s=1}^t \grad \estsigloss{s}(u_t) , u_t - \widetilde{u}_{t+1} }
		+ \eta_{t}^{-1} \paren[\Big]{ \psi(u_t) - \psi(\widetilde{u}_{t+1}) } \\
		=& \iprod[\Big]{ \grad \tprobloss{t}(u_t) + \grad \estsigloss{t}(u_t) ,  u_t - \widetilde{u}_{t+1} }
		+ \eta_{t}^{-1} \paren[\Big]{ \psi(u_t) - \psi(\widetilde{u}_{t+1}) - \iprod{ \grad \psi(u_t), u_t-\widetilde{u}_{t+1} }} \\
		=& \iprod{ \grad \estsigloss{t}(u_t) ,  u_t - \widetilde{u}_{t+1} }
		- \eta_{t}^{-1} \paren[\Big]{ \psi(\widetilde{u}_{t+1})- \psi(u_t)  - \iprod{ \grad \psi(u_t), \widetilde{u}_{t+1}-u_t }}  \\
		=& \iprod{ \grad \estsigloss{t}(u_t) ,  u_t - \widetilde{u}_{t+1} }
		- \eta_{t}^{-1} \breg{\psi}{\widetilde{u}_{t+1}}{u_t}
		\intertext{Using Lemma~\mainref{lemma:breg-lb}, the Bregman divergence can be lower bounded as}
		\leq& \iprod{ \grad \estsigloss{t}(u_t) ,  u_t - \widetilde{u}_{t+1} } - \frac{\eta_{t}^{-1}}{2} (u_t - \widetilde{u}_{t+1})^2(1 + \abs{u_t})  \\
		\leq& \frac{\eta_t}{2} \cdot \frac{ \paren{ \grad \estsigloss{t}(u_t) }^2} { 1 + \abs{u_t} }
		\enspace,
	\end{align*}
	where the final inequality follows by completing the square.
	The result follows by taking expectation.
\end{proof}

By applying the law of iterated expectations, we can further bound the quantity in Lemma \ref{lemma:prob-regret-initial-upper-bound}, formally stated in the following lemma.

\begin{lemma}\label{lemma:conditional-expectation-squared-gradient}
	Under Condition~\mainref{condition:sigmoid},
	the conditional expectation of the squared gradient term is at most
	\[
	\E[\Bigg]{ \frac{ \paren[\big]{\grad \estsigloss{t}(u_t)}^2 }{1+|u_t|} ~\Big| \filt_{t-1}  }
	\leq b_1^2\max\{b_1,2\} \paren[\Big]{
		\braces{ y_t(1)-\iprod{\xv_t,\bv_t(1)} }^4
		+\braces{ y_t(0)-\iprod{\xv_t,\bv_t(0)} }^4
	}\enspace.
	\]
	Thus, applying the law of iterated expectations leads to
	\[
	\sum_{t=1}^T \eta_t \E[\Bigg]{ \frac{ \paren[\big]{\grad \estsigloss{t}(u_t)}^2 }{2(1+|u_t|)} }
	\leq \frac{1}{2}b_1^2\max\{b_1,2\} \sum_{k \in \setb{0,1}} \E[\Bigg]{ \sum_{t=1}^T \eta_t \braces[\big]{ y_t(k) - \iprod{ \xv_t, \bv_t(k) } }^4 }
	\enspace.
	\]
\end{lemma}

\begin{proof}
    For $k\in\{0,1\}$, denote $\Delta_{t,k}(\bv)=y_t(k)-\iprod{\xv_t,\bv}$.
    The gradient of the estimated loss function can be computed as
    \begin{align*}
        \nabla \widehat{h}_t(u)=&\frac{\mathbf{1}[Z_t=1]}{p_t}\cdot \Delta_{t,1}(\bv_t(1))^2\cdot\left(\frac{1}{\phi(u)}\right)^{\prime}+\frac{\mathbf{1}[Z_t=0]}{1-p_t}\cdot \Delta_{t,0}(\bv_t(0))^2\cdot\left(\frac{1}{1-\phi(u)}\right)^{\prime}\enspace.
    \end{align*}
    Using Condition \mainref{condition:sigmoid}, and the fact that the cross term is zero (since $\mathbf{1}[Z_t=1]\mathbf{1}[Z_t=0]=0$), we can upper bound the square of the gradient at $u_t$ as
    \begin{align*}
        &(\nabla \widehat{h}_t(u_t))^2\\
        =&\frac{\mathbf{1}[Z_t=1]}{p_t^2}\cdot \Delta_{t,1}(\bv_t(1))^4\cdot\left[\left(\frac{1}{\phi(u_t)}\right)^{\prime}\right]^2+\frac{\mathbf{1}[Z_t=0]}{(1-p_t)^2}\cdot \Delta_{t,0}(\bv_t(0))^4\cdot\left[\left(\frac{1}{1-\phi(u_t)}\right)^{\prime}\right]^2\\
        =&\frac{\mathbf{1}[Z_t=1]}{p_t^2}\cdot \Delta_{t,1}(\bv_t(1))^4\cdot\left[\left(\frac{1}{\phi(u_t)}\right)^{\prime}\right]^2+\frac{\mathbf{1}[Z_t=0]}{(1-p_t)^2}\cdot \Delta_{t,0}(\bv_t(0))^4\cdot\left[\left(\frac{1}{\phi(-u_t)}\right)^{\prime}\right]^2\quad(\text{Condition \mainref{condition:sigmoid}-1})\\
        \leq&\bb_1^2\left[\frac{\mathbf{1}[Z_t=1]}{p_t}\cdot \Delta_{t,1}(\bv_t(1))^4\cdot\frac{1}{\phi(u_t)}+\frac{\mathbf{1}[Z_t=0]}{1-p_t}\cdot \Delta_{t,0}(\bv_t(0))^4\cdot\frac{1}{1-\phi(u_t)}\right]\quad(\text{Condition \mainref{condition:sigmoid}-3a, $p_t=\phi(u_t)$})\\
        \intertext{Since $\left|\frac{1}{1-\phi(u_t)}- \frac{1}{1-\phi(0)}\right|\leq \bb_1 |u_t|$ and $\left|\frac{1}{\phi(u_t)}- \frac{1}{\phi(0)}\right|=\left|\frac{1}{1-\phi(-u_t)}- \frac{1}{1-\phi(0)}\right|\leq \bb_1 |u_t|$ by Lemma \ref{lemma:u}-3, we obtain}
        \leq&\bb_1^2\left[\frac{\mathbf{1}[Z_t=1]}{p_t}\cdot \Delta_{t,1}(\bv_t(1))^4\cdot\left(\frac{1}{\phi(0)}+\bb_1|u_t|\right)+\frac{\mathbf{1}[Z_t=0]}{1-p_t}\cdot \Delta_{t,0}(\bv_t(0))^4\cdot\left(\frac{1}{1-\phi(0)}+\bb_1|u_t|\right)\right]\\
        =&\bb_1^2(2+\bb_1|u_t|)\left(\frac{\mathbf{1}[Z_t=1]}{p_t}\cdot \Delta_{t,1}(\bv_t(1))^4+\frac{\mathbf{1}[Z_t=0]}{1-p_t}\cdot \Delta_{t,0}(\bv_t(0))^4\right)\quad(\text{since $\phi(0)=1/2$})\enspace.
    \end{align*}
    Hence, we have established the following upper bound:
    \begin{align*}
        &\e\left[\frac{(\nabla \widehat{h}_t(u_t))^2}{1+|u_t|}\Big|\filt_{t-1}\right]\\
        \leq&\e\left[\frac{\bb_1^2(2+\bb_1|u_t|)}{1+|u_t|}\left(\frac{\mathbf{1}[Z_t=1]}{p_t}\cdot \Delta_{t,1}(\bv_t(1))^4+\frac{\mathbf{1}[Z_t=0]}{1-p_t}\cdot \Delta_{t,0}(\bv_t(0))^4\right)\Big|\filt_{t-1}\right]\\
        \leq&\bb_1^2\max\{\bb_1,2\}\e\left[\frac{\mathbf{1}[Z_t=1]}{p_t}\cdot \Delta_{t,1}(\bv_t(1))^4+\frac{\mathbf{1}[Z_t=0]}{1-p_t}\cdot \Delta_{t,0}(\bv_t(0))^4\Big|\filt_{t-1}\right]\\
        =&\bb_1^2\max\{\bb_1,2\}\left[\Delta_{t,1}(\bv_t(1))^4+\Delta_{t,0}(\bv_t(0))^4\right]\\
        =&\bb_1^2\max\{b_1,2\} \paren[\Big]{
		\braces{ y_t(1)-\iprod{\xv_t,\bv_t(1)} }^4
		+\braces{ y_t(0)-\iprod{\xv_t,\bv_t(0)} }^4
	}\enspace,
    \end{align*}
    which completes the proof.
\end{proof}

Combining Lemma \ref{lemma:prob-regret-initial-upper-bound} and Lemma \ref{lemma:conditional-expectation-squared-gradient}, we can prove the following lemma.

\begin{reflemma}{\mainref{lemma:probability-gradient-fourth-moment-bound}*}
	There exists a constant $C$ depending on the sigmoid Condition~\mainref{condition:sigmoid} such that the suboptimality term is bounded by the expected fourth moment of online residuals:
	\[
	\sum_{t=1}^T \E[\Big]{ \ptprobloss{t+1}(u_t) - \ptprobloss{t+1}( \widetilde{u}_{t+1} )}
	\leq 
	\frac{\bb_1^2}{2}\max\{\bb_1,2\} \sum_{k \in \setb{0,1}} \E[\Bigg]{ \sum_{t=1}^T \eta_t \paren[\big]{ y_t(k) - \iprod{ \xv_t, \bv_t(k) } }^4 }
	\enspace.
	\]
\end{reflemma}

\begin{proof}
    The result is proved by Lemma \ref{lemma:prob-regret-initial-upper-bound} and Lemma \ref{lemma:conditional-expectation-squared-gradient}.
\end{proof}

Combining Lemma \mainref{lemma:exp-prob-regret-initial-oco-bound}, Lemma \ref{lemma:regularizer-upper-bound}, and Lemma \mainref{lemma:probability-gradient-fourth-moment-bound}*, we establish the following upper bound on the expected probability regret.

\begin{lemma}\label{lemma:prob-regret-bound}
	Under Assumptions \mainref{assumption:moments}-\mainref{assumption:maximum-radius} and Condition \mainref{condition:sigmoid}, the expected probability regret is bounded as
    \begin{align*}
        \E{\regretprob_T}\leq& \frac{8T^{1/2}R_T}{\bb_3^3}\left[\left(\frac{c_1}{c_0}-1\right)^3+\frac{\bb_3}{4}\left(\frac{c_1}{c_0}-1\right)^2\right]\\
        &+ \frac{\bb_1^2}{2}\max\{\bb_1,2\}\sum_{k \in \setb{0,1}} \E[\Bigg]{ \sum_{t=1}^T \eta_t \braces[\big]{ y_t(k) - \iprod{ \xv_t, \bv_t(k) } }^4 }\enspace.
    \end{align*}
\end{lemma}

\begin{proof}
    By Lemma \mainref{lemma:exp-prob-regret-initial-oco-bound}, Lemma \ref{lemma:regularizer-upper-bound}, and Lemma \mainref{lemma:probability-gradient-fourth-moment-bound}*, we can derive the following upper bound:
    \begin{align*}
        \E{\regretprob_T}\leq&\frac{1}{\eta_{T+1}} \psi(u^*)
		+ \sum_{t=1}^T \E[\Big]{ \ptprobloss{t+1}(u_t) - \ptprobloss{t+1}( \widetilde{u}_{t+1} ) }\quad(\text{Lemma \mainref{lemma:exp-prob-regret-initial-oco-bound}})\\
        \leq&\frac{8T^{1/2}R_T}{\bb_3^3}\left[\left(\frac{c_1}{c_0}-1\right)^3+\frac{\bb_3}{4}\left(\frac{c_1}{c_0}-1\right)^2\right]\\
        &+\frac{\bb_1^2}{2}\max\{\bb_1,2\}\sum_{k \in \setb{0,1}} \E[\Bigg]{ \sum_{t=1}^T \eta_t \braces[\big]{ y_t(k) - \iprod{ \xv_t, \bv_t(k) } }^4 }\quad(\text{Lemma \ref{lemma:regularizer-upper-bound}, Lemma \mainref{lemma:probability-gradient-fourth-moment-bound}*})\enspace,
    \end{align*}
    which completes the proof.
\end{proof}

Note that this lemma differs from Proposition \mainref{prop:prob-regret-bound}.
To establish Proposition \mainref{prop:prob-regret-bound}, it remains to derive an upper bound on the fourth moment of the online residuals.
This bound will be established in Section \ref{section:C5}.

\subsection{Bounding Inverse Probability Moments}\label{section:C3}
Due to the inverse probability weighting in the AIPW estimators, it is essential to control the expectation of the inverse assignment probabilities in the analysis of the Neyman regret.
In this section, we establish an upper bound on the moments of these inverse probabilities, formally stated in Lemma \ref{lemma:p-power-moment}.
The result is obtained through a sequence of intermediate steps, relying on Lemma \mainref{lemma:effect-of-p-regularization}*, Corollary \mainref{corollary:p-moment}*, Lemma \mainref{lemma:expectation-of-estores}, and Lemma \ref{lemma:squared-residuals-random}.

Since the \ourdesign{} design solves a regularized optimization at each iteration to select the next treatment assignment probability, we begin with a detailed analysis of the relationship between the optimization problem and its corresponding minimizer.

The following lemma characterizes when the obtained probability $p$ will be on either side of $1/2$, equivalently, when the corresponding $u = \phi^{-1}(p)$ will be positive or negative.

\begin{lemma}\label{lemma:p-location}
    Under Condition~\mainref{condition:sigmoid}, consider $A, B\geq 0$ and define $p^*$ as the minimizer of the following program:
	\begin{align*}
	    p^* = \argmin_{p \in (0,1)} \frac{A}{p} + \frac{B}{1-p} + \Psi(p)
	\enspace.
    \end{align*}
	If $A\geq B$, then $p^*\geq 1/2$ and $u^*=\phi^{-1}(p^*)\geq 0$; if $A \leq B$, then $p^*\leq 1/2$ and $u^*=\phi^{-1}(p^*)\leq 0$.
\end{lemma}

\begin{proof}
    Suppose that $A \geq B$. 
    Since $u^*=\phi^{-1}(p^*)$, it is the minimizer of the following program:
    \begin{align*}
        u^*=\operatorname{argmin}_{u\in\mathbb{R}}\frac{A}{\phi(u)}+\frac{B}{1-\phi(u)}+\psi(u)\enspace.
    \end{align*}
    By the convexity property in Condition \mainref{condition:sigmoid}-2 and the fact that $A,B\geq 0$, we can verify that $\frac{A}{\phi(u)}+\frac{B}{1-\phi(u)}$ is a convex function.
    Moreover, $\psi(u)=\frac{1}{2}u^2+|u|^3$ is strictly convex.
    Hence the target function $g(u):=\frac{A}{\phi(u)}+\frac{B}{1-\phi(u)}+\psi(u)$ is strictly convex.
    Since Condition \mainref{condition:sigmoid}-1 implies that $\phi(0)=1/2$, evaluating at $u=0$, we obtain
    \begin{align*}
        g^{\prime}(u)\big|_{u=0}=\left(\frac{A}{\phi(u)}+\frac{B}{1-\phi(u)}+\psi(u)\right)^{\prime}\Big|_{u=0}=&-A\cdot\frac{\phi^{\prime}(0)}{\phi^2(0)}+B\cdot\frac{\phi^{\prime}(0)}{(1-\phi(0))^2}+\psi^{\prime}(0)\\
        =&-4(A-B)\phi^{\prime}(0)\enspace,
    \end{align*}
    which is nonpositive since $A\geq B$ and $\phi^{\prime}(0)\geq 0$ by its monotonicity.
    Since $g$ is convex and $u^*$ should satisfy the first-order equation $g^{\prime}(u^*)=0$, we can deduce that $u^*\geq 0$.
    By Condition \mainref{condition:sigmoid}-1, we also have $p^*=\phi(u^*)\geq 1/2$.
    
    Let us now suppose that $B \geq A$.
    Let $q^*=1-p^*$.
    Then $q^*$ is the minimizer of the following program:
    \begin{align*}
        q^*=\operatorname{argmin}_{q\in(0,1)}\frac{B}{q}+\frac{A}{1-q}+\Psi(1-q)\enspace.
    \end{align*}
    Since by the symmetry property in Condition \mainref{condition:sigmoid}-1, we have $\phi^{-1}(q)=-\phi^{-1}(1-q)$, this further indicates that $\Psi(q)=\Psi(1-q)$ since $\psi$ is an even function.
    Hence, by the result in case $A\geq B$, the conclusion follows.
\end{proof}

Based on Lemma \ref{lemma:p-location}, we next derive an upper bound on the inverse assignment probabilities by leveraging the structural properties of the corresponding minimization program.

\begin{reflemma}{\mainref{lemma:effect-of-p-regularization}*}
	Consider $A, B \geq 0$, $\eta>0$ and define $p^*$ as the minimizer of the following program:
	\begin{align*}
	    p^* = \argmin_{p \in (0,1)} \frac{A}{p} + \frac{B}{1-p} + \eta^{-1} \Psi(p)
	\enspace.
    \end{align*}
	Under Condition \mainref{condition:sigmoid}, the minimizer $p^*$ is bounded away from $0$ and $1$ in the following sense:
	\begin{align*}
	\frac{1}{p^*} \leq& 2 + \min\left\{\bb_1(\bb_2/6)^{1/4}\eta^{1/4}B^{1/4},\bb_1(\bb_2/\bb_3)^{1/2}(B/A)^{1/2}\right\}\enspace,\\ 
	\frac{1}{1-p^*} \leq& 2+\min\left\{\bb_1(\bb_2/6)^{1/4}\eta^{1/4}A^{1/4},\bb_1(\bb_2/\bb_3)^{1/2}(A/B)^{1/2}\right\}
	\enspace.
    \end{align*}
    For simplicity, we define $0/0$ as $0$ here.
\end{reflemma}

\begin{proof}
    We first consider the case where $A\geq B$.
    Let $u^*=\phi^{-1}(p^*)$.
    According to Lemma \ref{lemma:p-location}, $u^*\geq 0$.
    It also satisfies the following first-order equation:
    \begin{align*}
        A\left(\frac{1}{\phi(u^*)}\right)^{\prime}+B\left(\frac{1}{1-\phi(u^*)}\right)^{\prime}+\eta^{-1}(u^*+3(u^*)^2)=0\enspace.
    \end{align*}
    By Lemma \ref{lemma:u}-2 and Condition \mainref{condition:sigmoid}-1, we have
    \begin{align*}
        3\eta^{-1}(u^*)^2\leq& \eta^{-1}(u^*+3(u^*)^2)+B\left(\frac{1}{1-\phi(u^*)}\right)^{\prime}\quad(\text{$1/(1-\phi(u))$ is monotone increasing})\\
        =&-A\left(\frac{1}{\phi(u^*)}\right)^{\prime}\quad(\text{first-order equation})\\
        \leq& \frac{\bb_2}{2}\cdot \frac{A}{(1+u^*)^2}\quad(\text{Lemma \ref{lemma:u}-2})\\
        \leq& \frac{\bb_2A}{2(u^*)^2}\enspace,
    \end{align*}
    which implies that $u^*\leq (\bb_2/6)^{1/4}\eta^{1/4}A^{1/4}$. 
    On the other hand, by Lemma \ref{lemma:u}-1 and Condition \mainref{condition:sigmoid}-1, we have
    \begin{align*}
        \frac{\bb_3 B}{2}\leq& B\left(\frac{1}{1-\phi(0)}\right)^{\prime}\quad(\text{Lemma \ref{lemma:u}-1})\\
        \leq& B\left(\frac{1}{1-\phi(u^*)}\right)^{\prime}\quad(\text{$u^*\geq 0$ and $1/(1-\phi(u))$ is convex})\\
        \leq& \eta^{-1}(u^*+3(u^*)^2)+B\left(\frac{1}{1-\phi(u^*)}\right)^{\prime}\quad(\text{$u^*\geq 0$})\\
        =&-A\left(\frac{1}{\phi(u^*)}\right)^{\prime}\quad(\text{first-order equation})\\
        \leq& \frac{\bb_2A}{2(1+u^*)^2}\enspace,
    \end{align*}
    which implies that $u^*\leq (\bb_2/\bb_3)^{1/2}(A/B)^{1/2}-1$. 
    Finally, we translate the bound on $u^*$ back to a bound on $p^*$.
    By Lemma \ref{lemma:u}, we have
    \begin{align*}
        \frac{1}{1-p^*}=&\frac{1}{1-\phi(u^*)}\\
        \leq& \frac{1}{1-\phi(0)}+\bb_1u^*\quad(\text{Lemma \ref{lemma:u}-3, $u^*\geq 0$})\\
        \leq& 2+\min\left\{\bb_1(\bb_2/6)^{1/4}\eta^{1/4}A^{1/4},\bb_1(\bb_2/\bb_3)^{1/2}(A/B)^{1/2}\right\}\quad(\text{$\phi(0)=1/2$})\enspace.
    \end{align*}
    If $A<B$, a similar result follows by symmetry, i.e.,
    \begin{align*}
        \frac{1}{p^*}\leq 2+\min\left\{\bb_1(\bb_2/6)^{1/4}\eta^{1/4}B^{1/4},\bb_1(\bb_2/\bb_3)^{1/2}(B/A)^{1/2}\right\}\enspace.
    \end{align*}
    Hence the lemma is proved.
\end{proof}

Based on Lemma \mainref{lemma:effect-of-p-regularization}* and Jensen's inequality, we can bound the moments of the inverse probabilities in terms of the expected estimated squared residuals.
This result is summarized in the following corollary.

\begin{refcorollary}{\mainref{corollary:p-moment}*}
	Under Condition \mainref{condition:sigmoid}, for each iteration $t \in [T]$ and any $0 \leq k \leq 4$, the $k$th moments of the inverse probabilities are bounded as
	\begin{align*}
		\E[\Big]{ \paren[\Big]{ \frac{1}{p_t} }^k } \leq&\left(2+\bb_1(\bb_2/6)^{1/4}\eta_{t}^{1/4}\e[\widehat{A}_{t-1}(0)]^{1/4}\right)^k,
		\\
		\E[\Big]{ \paren[\Big]{ \frac{1}{1-p_t} }^k } \leq& \left(2+\bb_1(\bb_2/6)^{1/4}\eta_{t}^{1/4}\e[\widehat{A}_{t-1}(1)]^{1/4}\right)^k.
	\end{align*}
\end{refcorollary}

\begin{proof}
    For any $0 \leq k \leq 4$, note that
    \begin{align*}
        \left(2+\bb_1(\bb_2/6)^{1/4}\eta_{t}^{1/4}x^{1/4}\right)^k
    \end{align*}
    is concave (as a function of $x$) when $x\in(0,\infty)$.
    Hence, by the definition of $p_t$, Lemma \mainref{lemma:effect-of-p-regularization}* and Jensen's inequality, for any $0 \leq k \leq 4$ we have
    \begin{align*}
        \e\left[\frac{1}{p_t^k}\right]\leq&\e\left[\left(2+\bb_1(\bb_2/6)^{1/4}\eta_{t}^{1/4}\widehat{A}_{t-1}(0)^{1/4}\right)^k\right]\leq\left(2+\bb_1(\bb_2/6)^{1/4}\eta_{t}^{1/4}\e[\widehat{A}_{t-1}(0)]^{1/4}\right)^k\enspace.
    \end{align*}
    Similarly, we can prove the result for the moments of $1/(1-p_t)$.
\end{proof}

Corollary \mainref{corollary:p-moment}* implies that, to bound the moments of the inverse assignment probabilities, it suffices to control the expectation of the estimated squared residuals.
Lemma \mainref{lemma:predictor-expectation} through Lemma \ref{lemma:squared-residuals-random} are devoted to establishing the required upper bound.

The following lemma verifies the unbiasedness property of $\bv_t(k)$.

\begin{reflemma}{\mainref{lemma:predictor-expectation}}
	\predictorexpectation
\end{reflemma}

\begin{proof}
    We only prove the result for $k=1$. 
    By the law of iterated expectations and the explicit forms of $\bv_{t}(1)$ and $\bv_{t}^*(1)$, we have
    \begin{align*}
        \e\left[\bv_{t}(1)\right]=&\e\left[\left(\vec{X}_{t-1}^{\tran}\vec{X}_{t-1}+\eta_t^{-1}\vec{I}_d\right)^{-1}\sum_{s=1}^{t-1}y_s(1)\frac{\indicator{Z_s=1}}{p_s}\cdot \xv_s\right]\\
        =&\left(\vec{X}_{t-1}^{\tran}\vec{X}_{t-1}+\eta_t^{-1}\vec{I}_d\right)^{-1}\sum_{s=1}^{t-1}y_s(1)\e\left[\e\left[\frac{\indicator{Z_s=1}}{p_s}\Big|\filt_{s-1}\right]\right]\cdot \xv_s\\
        =&\left(\vec{X}_{t-1}^{\tran}\vec{X}_{t-1}+\eta_t^{-1}\vec{I}_d\right)^{-1}\sum_{s=1}^{t-1}y_s(1)\cdot \xv_s\\
        =&\bv^*_t(1)\enspace.
        \qedhere
    \end{align*}
\end{proof}

The following lemma provides the explicit form of the online ridge predictors.

\begin{reflemma}{\mainref{lemma:ridge-pred-exact-form}*}
	At each iteration $t$, the online ridge predictions admit the decompositions
	\begin{align*}
        &\iprod{ \xv_t, \bv_t(1) } = \sum_{s=1}^{t-1} \Pi_{t,s} y_s(1) \frac{\indicator{Z_s=1}}{p_s}
	    \quadand
	    \iprod{ \xv_t, \bv_t(0) } = \sum_{s=1}^{t-1} \Pi_{t,s} y_s(0) \frac{\indicator{Z_s=0}}{1-p_s}
	    \enspace,\\
        &\iprod{ \xv_t, \bv_t^*(1) } = \sum_{s=1}^{t-1} \Pi_{t,s} y_s(1) 
	    \quadand
	    \iprod{ \xv_t, \bv_t^*(0) } = \sum_{s=1}^{t-1} \Pi_{t,s} y_s(0) 
	    \enspace.
    \end{align*}
\end{reflemma}

\begin{proof}
    By the definition of $\bv_t(1)$ and $\bv_t^*(1)$, we have
    \begin{align*}
        \iprod{\xv_t,\bv_t(1)}=&\Big\langle\xv_t,\left(\vec{X}_{t-1}^{\tran}\vec{X}_{t-1}+\eta_t^{-1}\vec{I}_d\right)^{-1}\sum_{s=1}^{t-1}y_s(1)\frac{\indicator{Z_s=1}}{p_s}\cdot \xv_s\Big\rangle\\
        =&\sum_{s=1}^{t-1} \Pi_{t,s} y_s(1) \frac{\indicator{Z_s=1}}{p_s}
    \end{align*}
    and
    \begin{align*}
        \iprod{\xv_t,\bv_t^*(1)}=\Big\langle\xv_t,\left(\vec{X}_{t-1}^{\tran}\vec{X}_{t-1}+\eta_t^{-1}\vec{I}_d\right)^{-1}\sum_{s=1}^{t-1}y_s(1)\cdot \xv_s\Big\rangle=\sum_{s=1}^{t-1} \Pi_{t,s} y_s(1)\enspace.
    \end{align*}
    We can similarly derive the form of $\iprod{\xv_t,\bv_t(0)}$ and $\iprod{\xv_t,\bv_t^*(0)}$.
\end{proof}

Based on Lemma \mainref{lemma:ridge-pred-exact-form}*, the next lemma derives the explicit form of the tracking error terms that appear in the expansion of the expectation of the estimated squared residuals.

\begin{lemma}\label{lemma:tracking-term}
	The expected prediction tracking error terms satisfy
	\begin{align*}
		\E[\Big]{\iprod{ \xv_t,  \optbv_t(1) - \bv_t(1)}^2}
		=&
		\sum_{s=1}^{t-1} \Pi_{t,s}^2 y_s(1)^2 \E[\Big]{ \frac{1}{p_s} - 1 } 
		\enspace,\\
		\E[\Big]{\iprod{ \xv_t,  \optbv_t(0) - \bv_t(0)}^2}
		=&
		\sum_{s=1}^{t-1} \Pi_{t,s}^2 y_s(0)^2 \E[\Big]{ \frac{1}{1-p_s} - 1 }\enspace.
	\end{align*}
\end{lemma}

\begin{proof}
    We only prove the first equality.
    The second equality follows from the symmetry between the treated and control groups.
    By the explicit form of $\langle \xv_t,\bv_t^*(1)-\bv_t(1)\rangle$ calculated in Lemma \mainref{lemma:ridge-pred-exact-form}*, we have
    \begin{align}\label{lemma:tracking-term_eq1}
        &\mathbb{E}\left[\langle \xv_t,\bv_t^*(1)-\bv_t(1)\rangle^2\right]\notag\\
        =&\mathbb{E}\left[\left(\sum_{s=1}^{t-1}\Pi_{t,s}y_s(1)\left(\frac{\mathbf{1}[Z_s=1]}{p_s}-1\right)\right)^2\right]\notag\\
        =&\underbrace{\mathbb{E}\left[\sum_{1\leq s_1\neq s_2\leq t-1}\Pi_{t,s_1}\Pi_{t,s_2}y_{s_1}(1)y_{s_2}(1)\left(\frac{\mathbf{1}[Z_{s_1}=1]}{p_{s_1}}-1\right)\left(\frac{\mathbf{1}[Z_{s_2}=1]}{p_{s_2}}-1\right)\right]}_{:=S_1~(\text{cross terms})}\notag\\
        &+\underbrace{\mathbb{E}\left[\sum_{s=1}^{t-1}\Pi_{t,s}^2y_s(1)^2\left(\frac{\mathbf{1}[Z_s=1]}{p_s}-1\right)^2\right]}_{:=S_2~(\text{square terms})}\enspace.
    \end{align}
    By the law of iterated expectations, we can simplify $S_1$ as:
    \begin{align}\label{lemma:tracking-term_eq2}
        S_1=&2\sum_{1\leq s_1<s_2\leq t-1}\Pi_{t,s_1}\Pi_{t,s_2}y_{s_1}(1)y_{s_2}(1)\e\left[\left(\frac{\indicator{Z_{s_1}=1}}{p_{s_1}}-1\right)\left(\frac{\indicator{Z_{s_2}=1}}{p_{s_2}}-1\right)\right]\notag\\
        =&2\sum_{1\leq s_1<s_2\leq t-1}\Pi_{t,s_1}\Pi_{t,s_2}y_{s_1}(1)y_{s_2}(1)\notag\\
        &\times\e\left[\e\left[\left(\frac{\indicator{Z_{s_1}=1}}{p_{s_1}}-1\right)\left(\frac{\indicator{Z_{s_2}=1}}{p_{s_2}}-1\right)\Big|\filt_{s_2-1}\right]\right]\notag\\
        =&0\enspace.
    \end{align}
    Similarly, we can simplify $S_2$ as:
    \begin{align}\label{lemma:tracking-term_eq3}
        S_2=&\sum_{s=1}^{t-1}\Pi_{t,s}^2y_{s}(1)^2\e\left[\e\left[\left(\frac{\indicator{Z_{s}=1}}{p_{s}}-1\right)^2\Big|\filt_{s-1}\right]\right]\notag\\
        =&\sum_{s=1}^{t-1}\Pi_{t,s}^2y_s(1)^2\e\left[\frac{1}{p_s}-1\right]\enspace.
    \end{align}
    By \eqref{lemma:tracking-term_eq1}, \eqref{lemma:tracking-term_eq2} and \eqref{lemma:tracking-term_eq3}, we obtain
    \begin{align*}
        \E[\Big]{\iprod{ \xv_t,  \optbv_t(1) - \bv_t(1)}^2}
		=&
		\sum_{s=1}^{t-1} \Pi_{t,s}^2 y_s(1)^2 \E[\Big]{ \frac{1}{p_s} - 1 }\enspace.
    \end{align*}
    Hence, the result is proved.
\end{proof}

Based on Lemma \mainref{lemma:predictor-expectation} and Lemma \ref{lemma:tracking-term}, we derive the explicit expression for the expectation of the estimated squared residuals in the following lemma.

\begin{reflemma}{\mainref{lemma:expectation-of-estores}}
	\expectationestore
\end{reflemma}

\begin{proof}
    By the law of iterated expectations, the expectation of $\widehat{A}_t(1)$ can be calculated by:
    \begin{align*}
        \e[\widehat{A}_t(1)]=&\e\left[\sum_{s=1}^{t}\frac{\mathbf{1}[Z_s=1]}{p_s}\cdot\left(y_s(1)-\iprod{\xv_s,\bv_s(1)}\right)^2\right]\notag\\
        =&\sum_{s=1}^{t}\e\left[\e\left[\frac{\mathbf{1}[Z_s=1]}{p_s}\cdot\left(y_s(1)-\iprod{\xv_s,\bv_s(1)}\right)^2\Bigg| \mathcal{F}_{s-1}\right]\right]\quad(\text{law of iterated expectations})\notag\\
        =&\sum_{s=1}^{t}\e\left[\left(y_s(1)-\iprod{\xv_s,\bv_s(1)}\right)^2\right]\notag\\
        =&\sum_{s=1}^{t}\e\left[\left((y_s(1)-\iprod{\xv_s,\bv_s^*(1)})+\iprod{\xv_s,\bv_s^*(1)-\bv_s(1)}\right)^2\right]\\
        =&\sum_{s=1}^{t}(y_s(1)-\iprod{\xv_s,\bv_s^*(1)})^2+2\sum_{s=1}^{t}(y_s(1)-\iprod{\xv_s,\bv_s^*(1)})\e\left[\iprod{\xv_s,\bv_s^*(1)-\bv_s(1)}\right]\\
        &+\sum_{s=1}^{t}\e\left[\iprod{\xv_s,\bv_s^*(1)-\bv_s(1)}^2\right]\\
        \intertext{By Lemma \mainref{lemma:predictor-expectation}, the cross term equals 0. Hence, by applying Lemma \ref{lemma:tracking-term}, we have}
        =&\sum_{s=1}^{t}(y_s(1)-\iprod{\xv_s,\bv_s^*(1)})^2+\sum_{s=1}^t\sum_{r=1}^{s-1}\Pi_{s,r}^2y_r(1)^2\e\left[\frac{1}{p_r}-1\right]\\
        =&A_{t}^*(1)+\sum_{s=1}^t\sum_{r=1}^{s-1}\Pi_{s,r}^2y_r(1)^2\e\left[\frac{1}{p_r}-1\right]\quad(\text{definition of $A_t^*(1)$})\enspace.
    \end{align*}
    Similarly, we can prove that
    \begin{align*}
        \E{\widehat{A}_t(0)}=A_{t}^*(0)+\sum_{s=1}^t\sum_{r=1}^{s-1}\Pi_{s,r}^2y_r(0)^2\e\left[\frac{1}{1-p_r}-1\right]\enspace.
    \end{align*}
    This completes the proof.
\end{proof}

In order to establish upper bounds on the estimated squared residuals, we need to control the deterministic summations that appear in the explicit forms of $\E{\widehat{A}_t(1)}$ and $\E{\widehat{A}_t(0)}$.
Lemma \ref{lemma:regularity} through Corollary \ref{corollary:squared-residual-deterministic} are established to achieve this objective.

In particular, Lemma \ref{lemma:regularity} provides upper bounds for the spectral norms of the matrices involved in bounding the regression coefficients.

\begin{lemma}\label{lemma:regularity}
    Under Assumption \mainref{assumption:covariate-regularity}, for any $t\in[T]$, it holds that
    \begin{align*}
        &\left\|\left(\xM_{t-1}^{\tran}\xM_{t-1}+\eta_{t}^{-1}\mat{I}_d\right)^{-1}\right\|\leq \gamma((t-1)\vee \eta_t^{-1})^{-1},
    \end{align*}
    where $\gamma:= \max\{\gamma_0,c_2,1\}>0$ is a constant.
\end{lemma}
\begin{proof}
    For any $t\in[T]$, we have the trivial upper bound $\left\|\left(\xM_{t-1}^{\tran}\xM_{t-1}+\eta_t^{-1}\mat{I}_d\right)^{-1}\right\|\leq \eta_t$.
    Let us now consider two cases based on the iteration $t$.
    
    \textbf{Case 1:}
    Suppose that we are in an early iteration $t < \gamma_0 \cdot \eta_t^{-1}+1$.
    Rearranging this and using the definition of $\gamma$, we have that 
    $
    \eta_t^{-1} > \gamma_0^{-1} (t-1) \geq \gamma^{-1} (t-1)
    $.
    Thus, $\eta_t \leq  \gamma \paren{  \eta_t^{-1} \vee (t-1) }^{-1}$.
    Substituting this into the upper bound above completes the proof in this case.
    
    \textbf{Case 2:}
    Suppose that $t\geq \gamma_0\cdot \eta_t^{-1}+1$.
    We will show that this means that we are in a later iteration, in the sense that the covariate regularity condition of Assumption~\mainref{assumption:covariate-regularity} may be applied. 
    Because $R_t \geq d^{1/2}$, we have that the inverse step size can be lower bounded as $\eta_t^{-1}=T^{1/2}R_t \ge T^{1/2}d^{1/2}$.
    Thus, this means that $t - 1 \geq \gamma_0 \cdot T^{1/2}d^{1/2}$ so that we can apply Assumption \mainref{assumption:covariate-regularity} to obtain that
    \begin{align*}
        \norm{ \paren{\xM_{t-1}^\tran \xM_{t-1} + \eta_t^{-1} \unitM_d }^{-1}} 
    \leq&\left(\frac{t-1}{c_2}+\eta_t^{-1}\right)^{-1}\\
    \leq& \left(\min\{c_2^{-1},1\}(t-1+\eta_t^{-1})\right)^{-1}\\
    \leq&\max\{c_2,1\}((t-1)\vee \eta_t^{-1})^{-1}\enspace.
    \end{align*}
    Therefore, the proof is completed in this case.
\end{proof}

Recall that the algorithm maintains a variable $R_t$ which is initialized as $R_0 = d^{1/2}$ and updated as $R_t = \max \setb{ R_{t-1}, \norm{\xv_t} }$.
While we informally refer to this as the ``largest covariate norm seen so far'', that technically may not be true when the covariates seen so far have norm less than $d^{1/2}$.
In this sense, our proposed algorithm does not guarantee that the final variable $R_T$ is actually equal to the maximum radius, e.g. it may be larger.
However, the following lemma shows that, under Assumption \mainref{assumption:covariate-regularity}, we will have that $R_T$ is within a constant of the maximum radius $R$.
This allows us to use $R_T$ as a proxy for $R$ in all of our proofs.

\begin{lemma}\label{lemma:R}
    Under Assumption \mainref{assumption:covariate-regularity}, the final proxy radius $R_T$ is not so much bigger than the actual radius $R$: $R_T\leq (\max\{c_2,1\})^{1/2}R$.
\end{lemma}
\begin{proof}
    Let us first consider the case that the actual radius $R$ is larger than the initial proxy radius $R_0$.
    Throughout the algorithm, the proxy radius $R_t$ is increasing according to the rule $R_t = \max\setb{ R_{t-1}, \| \xv_t \| }$.
    If we have that $R_0 \leq R$, then at some iteration $t$, we will have that $R_t = R$ and thus at the end of the algorithm, we have that $R_T = R$.
    In this case, the inequality above holds.

    Next, let us consider the case where the actual maximum radius $R$ is smaller than the initial proxy radius $R_0$, i.e. $R \leq R_0 = d^{1/2}$.
    We will show that Assumption~\mainref{assumption:covariate-regularity} prevents $R$ from being too much smaller than $d^{1/2}$.
    By standard linear algebra arguments, we have that
    \begin{align*}
        \lambda_{\min}(\xM_{T}^{\tran}\xM_{T})\leq \frac{1}{d}\operatorname{tr}(\xM_{T}^{\tran}\xM_{T})=\frac{1}{d}\operatorname{tr}\left(\sum_{t=1}^T\xv_t\xv_t^{\tran}\right)=\frac{1}{d}\sum_{t=1}^T\operatorname{tr}\left(\xv_t\xv_t^{\tran}\right)=\frac{1}{d}\sum_{t=1}^T\|\xv_t\|^2\leq \frac{TR^2}{d}\enspace,
    \end{align*}
    which, together with Assumption \mainref{assumption:covariate-regularity}, implies that $R\geq c_2^{-1/2}d^{1/2}$.
    Because $R_0 = d^{1/2}$, it now follows that $R \geq c_2^{-1/2}R_0$, establishing the desired inequality.
\end{proof}

Lemma \ref{lemma:power-sum} provides upper bounds for deterministic $p$-series summations, which play a crucial role in establishing Lemma \ref{lemma:deterministic-summation-1}.

\begin{lemma}\label{lemma:power-sum}
    For any $k>1$ and fixed $t,r\in[T]$, we have $\sum_{s=t}^T((s-1)\vee \eta_r^{-1})^{-k}\leq \xi_k((t-1)\vee \eta_r^{-1})^{-(k-1)}$, where $\xi_k:=\frac{3k-2}{k-1}>0$ is a constant.
\end{lemma}
\begin{proof}
    We consider two regimes depending on $t$.
    For $t\geq \lceil\eta_r^{-1}\rceil$ and $t\geq 2$, we have
    \begin{align}\label{lemma:power-sum_eq1}
        \sum_{s=t}^T((s-1)\vee \eta_r^{-1})^{-k}\leq&\sum_{s=t-1}^{T-1}s^{-k}\notag\\
        \leq&(t-1)^{-k}+\int_{t-1}^{\infty}s^{-k}\mathrm{d}s 
        	&(\text{by integral comparison})\notag\\
        \leq&(t-1)^{-k+1}+\frac{(t-1)^{-(k-1)}}{k-1}\notag\\
        \leq&\frac{k}{k-1}\cdot (t-1)^{-(k-1)}\enspace.
    \end{align}
    For $t< \lceil\eta_r^{-1}\rceil$ or $t=1$, we have
    \begin{align}\label{lemma:power-sum_eq2}
        \sum_{s=t}^T((s-1)\vee \eta_r^{-1})^{-k}\leq&\sum_{s=\lceil\eta_r^{-1}\rceil+1}^T((s-1)\vee \eta_r^{-1})^{-k}+\sum_{s=1}^{\lceil\eta_r^{-1}\rceil}((s-1)\vee \eta_r^{-1})^{-k}\notag\\
        \leq&\sum_{s=\lceil\eta_r^{-1}\rceil+1}^T((s-1)\vee \eta_r^{-1})^{-k}+\sum_{s=1}^{\lceil\eta_r^{-1}\rceil}\eta_r^{k}\notag\\
        \leq&\frac{k}{k-1}\cdot \lceil\eta_r^{-1}\rceil^{-(k-1)}+\eta_r^{k}\lceil\eta_r^{-1}\rceil\notag\\
        \leq&\frac{k}{k-1}\eta_r^{(k-1)}+2\eta_r^{(k-1)}\quad(\text{since $\eta_r^{-1}=T^{1/2}R_r\geq 1$})\enspace.
    \end{align}
    Combining the results in \eqref{lemma:power-sum_eq1} and \eqref{lemma:power-sum_eq2}, for any fixed $t,r\in[T]$, we have
    \begin{align*}
        \sum_{s=t}^T((s-1)\vee \eta_r^{-1})^{-k}\leq \left(\frac{k}{k-1}+2\right)((t-1)\vee\eta_r^{-1})^{-(k-1)}=\xi_k((t-1)\vee\eta_r^{-1})^{-(k-1)}\enspace.
    \end{align*}
    This completes the proof.
\end{proof}

Next, we provide a lemma that contains many useful inequalities regarding the per-iteration leverage scores $\Pi_{t,s}$.
Recall that we define the \emph{per iteration leverage scores} as
\[
\Pi_{t,s} = \xv_t^\tran \paren[\Big]{ \xM_{t-1}^\tran \xM_{t-1} + \eta_t^{-1} \unitM }^{-1} \xv_s
\enspace,
\]
where $\xM_{t-1}^\tran \xM_{t-1} = \sum_{s \leq t-1} \xv_s \xv_s^\tran$.
As demonstrated by Lemma~\mainref{lemma:ridge-pred-exact-form}*, the per-iteration leverage scores arise naturally in decomposing the online linear predictors.
The first result is a bound on each of the per-iteration leverage scores in terms of the maximum observed radius and the step size.

\begin{corollary}\label{corollary:pi}
	Under Assumption \mainref{assumption:covariate-regularity}, for any $1\leq s\leq t\leq T$, it holds that $|\Pi_{t,s}|\leq\gamma R_tR_s((t-1)\vee \eta_t^{-1})^{-1}$.
\end{corollary}

\begin{proof}
	Using the definition of $\Pi_{t,s}$, the operator norm bound, Lemma \ref{lemma:regularity}, and the definition of $R_t$, we have that
	\begin{align*}
		|\Pi_{t,s}|=
		&\left|\xv_t^{\tran}\left(\xM_{t-1}^{\tran}\xM_{t-1}+\eta_t^{-1}\mat{I}_d\right)^{-1}\xv_s\right|\\
		\leq& \|\xv_t\|\|\xv_s\|\left\|\left(\xM_{t-1}^{\tran}\xM_{t-1}+\eta_{t}^{-1}\mat{I}_d\right)^{-1}\right\|\\
		\leq& \gamma R_tR_s((t-1)\vee \eta_t^{-1})^{-1}
		\enspace.
		\qedhere
	\end{align*}
\end{proof}

The bound on the per-iteration leverage score provided by Corollary~\ref{corollary:pi} can be tight for any particular $\Pi_{t,s}$.
However, tighter results can be obtained when we consider various partial sums over the leverage scores.
Bounding such partial sums will be critical for many parts of our analyses.
As such, we have gathered all of the main inequalities of this type in the following lemma, which will be referenced throughout the remainder of the appendix.

\begin{lemma}\label{lemma:deterministic-summation-1}
    Denote $c:=(1+c_3(\max\{c_2,1\})^{1/2})>0$.
    Under Assumptions \mainref{assumption:moments}-\mainref{assumption:maximum-radius}, for any $t\in[T]$ and $k\in\{0,1\}$, we have
    \begin{enumerate}
        \item[(1)] $\sum_{s=1}^{t-1}\Pi_{t,s}^2\leq \gamma R_t^2((t-1)\vee \eta_t^{-1})^{-1}$.
        \item[(2)] $\sum_{s=t+1}^{T}\Pi_{s,t}^2\leq c\gamma R_t^{2}((t-1)\vee \eta_t^{-1})^{-1}$.
        \item[(3)] $\sum_{s=1}^{t}R_s^{-\nu_1}\sum_{r=1}^{s-1}R_r^{-\nu_2}\Pi_{s,r}^2y_r(k)^2\leq c\,c_1^2\gamma \xi_2^{1/2}R_t^{3/2-\nu_1-\nu_2}T^{1/4}$ for any $\nu_1,\nu_2\geq 0$ such that $\nu_1+\nu_2\leq 1$.
    \end{enumerate}
\end{lemma}

\begin{proof}
    It suffices to prove the case $k=1$.
    The case $k=0$ is identical.
    \begin{enumerate}
        \item[(1)] By the definition of $\Pi_{t,s}$ and Corollary \ref{corollary:pi}, we have
        \begin{align*}
            \sum_{s=1}^{t-1}\Pi_{t,s}^2=&\sum_{s=1}^{t-1}\xv_t^{\tran}\left(\xM_{t-1}^{\tran}\xM_{t-1}+\eta_t^{-1}\mat{I}_d\right)^{-1}\xv_s\xv_s^{\tran}\left(\xM_{t-1}^{\tran}\xM_{t-1}+\eta_t^{-1}\mat{I}_d\right)^{-1}\xv_t\\
            =&\xv_t^{\tran}\left(\xM_{t-1}^{\tran}\xM_{t-1}+\eta_t^{-1}\mat{I}_d\right)^{-1}\xM_{t-1}^{\tran}\xM_{t-1}\left(\xM_{t-1}^{\tran}\xM_{t-1}+\eta_t^{-1}\mat{I}_d\right)^{-1}\xv_t\\
            \leq&\xv_t^{\tran}\left(\xM_{t-1}^{\tran}\xM_{t-1}+\eta_t^{-1}\mat{I}_d\right)^{-1}\left(\xM_{t-1}^{\tran}\xM_{t-1}+\eta_t^{-1}\mat{I}_d\right)\left(\xM_{t-1}^{\tran}\xM_{t-1}+\eta_t^{-1}\mat{I}_d\right)^{-1}\xv_t\\
            =&\xv_t^{\tran}\left(\xM_{t-1}^{\tran}\xM_{t-1}+\eta_t^{-1}\mat{I}_d\right)^{-1}\xv_t\\
            =&\Pi_{t,t}\\
            \leq&\gamma R_t^2((t-1)\vee \eta_t^{-1})^{-1}\quad(\text{Corollary \ref{corollary:pi}})\enspace.
        \end{align*}
        \item[(2)] Denote $\vec{A}_{s}=\xM_{s}^{\tran}\xM_{s}+\eta_{s+1}^{-1}\mat{I}_d$ and $\vec{\widebar{A}}_s=\xM_{s}^{\tran}\xM_{s}+\eta_{s}^{-1}\mat{I}_d$ for any $s\in[T]$. 
        By Lemma \ref{lemma:R} and Assumption \mainref{assumption:maximum-radius}, for any $s$, we have $\xv_s^{\tran}\vec{A}_{s-1}^{-1}\xv_s\leq R_s^2\eta_s=R_sT^{-1/2}\leq R_TT^{-1/2}\leq c_3(\max\{c_2,1\})^{1/2}$.
        Now we fix arbitrary $t\in[T]$. 
        Since $\vec{\widebar{A}}_{s}=\vec{A}_{s-1}+\xv_s\xv_s^{\tran}$, by Sherman-Morrison formula, we have
        \begin{align*}
            \vec{\widebar{A}}_{s}^{-1}=\vec{A}_{s-1}^{-1}-\frac{\vec{A}_{s-1}^{-1}\xv_s\xv_s^{\tran}\vec{A}_{s-1}^{-1}}{1+\xv_s^{\tran}\vec{A}_{s-1}^{-1}\xv_s}\enspace,
        \end{align*}
        which implies that
        \begin{align*}
            \xv_t^{\tran}\left(\vec{A}_{s-1}^{-1}-\vec{\widebar{A}}_{s}^{-1}\right)\xv_t=\frac{\xv_t^{\tran}\vec{A}_{s-1}^{-1}\xv_s\xv_s^{\tran}\vec{A}_{s-1}^{-1}\xv_t}{1+\xv_s^{\tran}\vec{A}_{s-1}^{-1}\xv_s}\geq\frac{1}{1+c_3(\max\{c_2,1\})^{1/2}}\left(\xv_t^{\tran}\vec{A}_{s-1}^{-1}\xv_s\right)^2=c^{-1}\left(\xv_t^{\tran}\vec{A}_{s-1}^{-1}\xv_s\right)^2\enspace.
        \end{align*}
        Hence by Corollary \ref{corollary:pi}, we have
        \begin{align*}
            \sum_{s=t+1}^{T}\Pi_{s,t}^2=&\sum_{s=t+1}^{T}\left(\xv_s^{\tran}\left(\xM_{s-1}^{\tran}\xM_{s-1}+\eta_s^{-1}\mat{I}_d\right)^{-1}\xv_t\right)^2\\
            =&\sum_{s=t+1}^{T}\left(\xv_t^{\tran}\vec{A}_{s-1}^{-1}\xv_s\right)^2\\
            \leq&c\sum_{s=t+1}^{T}\xv_t^{\tran}(\vec{A}_{s-1}^{-1}-\vec{\widebar{A}}_{s}^{-1})\xv_t\\
            \intertext{Since $\eta_s\geq \eta_{s+1}$ implies that $\vec{A}_s\succeq \vec{\widebar{A}}_{s}$, we can further upper bound it by}
            \leq&c\sum_{s=t+1}^{T}\xv_t^{\tran}(\vec{A}_{s-1}^{-1}-\vec{A}_{s}^{-1})\xv_t\\
            \leq&c\,\xv_t^{\tran}\vec{A}_{t}^{-1}\xv_t\\
            \leq&c\,\xv_t^{\tran}\left(\xM_{t-1}^{\tran}\xM_{t-1}+\eta_t^{-1}\mat{I}_d\right)^{-1}\xv_t\\
            \leq&c\gamma R_t^{2}((t-1)\vee \eta_t^{-1})^{-1}\quad(\text{Corollary \ref{corollary:pi}})\enspace.
        \end{align*}
        \item[(3)] By Lemma \ref{lemma:power-sum} and the result in Part 2, we have
        \begin{align*}
            &\sum_{s=1}^{t}R_s^{-\nu_1}\sum_{r=1}^{s-1}R_r^{-\nu_2}\Pi_{s,r}^2y_r(1)^2\\
            \leq&\sum_{r=1}^{t-1}R_r^{-\nu_1-\nu_2}\left(\sum_{s=r+1}^t\Pi_{s,r}^2\right)y_r(1)^2\quad(\text{swap the order of summation, $R_{r}\leq R_s$})\\
            \leq&c\gamma \sum_{r=1}^{t-1}R_r^{2-\nu_1-\nu_2}((r-1)\vee\eta_r^{-1})^{-1}y_r(1)^2\quad(\text{by result in Part 2})\\
            =&c\gamma \sum_{r=1}^{t-1}R_r^{1-\nu_1-\nu_2}((r-1)R_r^{-1}\vee T^{1/2})^{-1}y_r(1)^2\\
            \leq&c\gamma \sum_{r=1}^{t-1}R_t^{1-\nu_1-\nu_2}((r-1)R_t^{-1}\vee T^{1/2})^{-1}y_r(1)^2\quad(\text{since $R_r\leq R_t$ and $1-\nu_1-\nu_2\geq 0$})\\
            =&c\gamma \sum_{r=1}^{t-1}R_t^{2-\nu_1-\nu_2}((r-1)\vee\eta_t^{-1})^{-1}y_r(1)^2\\
            \leq&c\gamma R_t^{2-\nu_1-\nu_2}\left(\sum_{r=1}^{t-1}((r-1)\vee\eta_t^{-1})^{-2}\right)^{1/2}\left(\sum_{r=1}^{t-1}y_r(1)^4\right)^{1/2}\quad(\text{Cauchy-Schwarz})\\
            \leq&c\,c_1^2\gamma \xi_2^{1/2}R_t^{2-\nu_1-\nu_2}\eta_t^{1/2}T^{1/2}\quad(\text{Lemma \ref{lemma:power-sum}, Assumption \mainref{assumption:moments}})\\
            =&c\,c_1^2\gamma \xi_2^{1/2}R_t^{3/2-\nu_1-\nu_2}T^{1/4}\enspace.
        \end{align*}
    \end{enumerate}
    This completes the proof.
\end{proof}

The following lemma provides an upper bound on the squared norm of the OLS estimator $\bv_v^*$ for any sequence of outcomes ${v_t : t \in [T]}$ with bounded empirical second moment.

\begin{lemma}\label{lemma:ols-norm}
    Under Assumption \mainref{assumption:covariate-regularity}, for any $\{v_t:t\in[T]\}$ such that $\sum_{t=1}^Tv_t^2\leq c_1^2T$, it holds that $\norm{ \bv^*_v}^2 \leq c_1^2 c_2$.
    Here $\optbv_v:=\argmin_{\bv\in\mathbb{R}^d}\sum_{t\in[T]}(v_t-\iprod{ \xv_t, \bv })^2$.
\end{lemma}

\begin{proof}
    Denote $\vec{V}=(v_1,\ldots,v_T)^{\tran}$ and $\vec{X}=(\xv_1,\ldots,\xv_T)^{\tran}$.
    Using the explicit form of the OLS predictor $\bv^*_v$ and an operator norm bound, we have
    \[
    \norm{ \bv^*_v}^2
    =\mat{V}^{\tran}\xM \paren[\big]{ \xM^{\tran}\xM }^{-2} \xM^{\tran} \mat{V}
    \leq \norm{ \xM \paren[\big]{ \xM^{\tran}\xM }^{-2} \xM^{\tran} }  \cdot \norm{\mat{V}}^2
    \enspace.
    \]
    We use Assumption \mainref{assumption:covariate-regularity} to bound the operator norm as
    \[
    \norm{ \xM \paren[\big]{ \xM^{\tran}\xM }^{-2} \xM^{\tran} }
    = \lambda_{\max}\paren[\big]{ \paren[\big]{ \xM^{\tran}\xM }^{-1} }
    = \left(\lambda_{\min}\paren[\big]{  \xM^{\tran}\xM  }\right)^{-1}
    \leq c_2T^{-1}
    \enspace.
    \]
    On the other hand, by assumption we have $\|\mat{V}\|^2=\sum_{t=1}^Tv_t^2\leq c_1^2T$.
    Combining these bounds with the operator norm bound yields the desired result, i.e. $\norm{ \bv^*_v}^2 \leq c_1^2 c_2$.
\end{proof}

We establish the following one-step objective improvement identity in ridge regression, which is crucial in bounding the prediction regret.

\begin{lemma}\label{lemma:ridge}
    For any $\lambda>0$ and sequence $\{v_t:t\in[T]\}$, denote $G_t(\bv)=\sum_{s<t}(v_s-\iprod{\xv_s,\bv})^2+\lambda\|\bv\|^2$ for any $t\geq 1$.
    Further let $\bv_{t}$ denote the minimizer of $G_t$.
    Then for any $t\in[T]$, we have the following equality:
    \begin{align*}
        G_t(\bv_t)-G_{t+1}(\bv_{t+1})=-\frac{1}{1+\xv_t^{\tran}(\vec{X}_{t-1}^{\tran}\vec{X}_{t-1}+\lambda \mat{I}_d)^{-1}\xv_t}\cdot (v_t-\iprod{\xv_t,\bv_t})^2\enspace,
    \end{align*}
    which also implies the following equality:
    \begin{align*}
        G_{t+1}(\bv_{t})-G_{t+1}(\bv_{t+1})=\frac{\xv_t^{\tran}(\vec{X}_{t-1}^{\tran}\vec{X}_{t-1}+\lambda \mat{I}_d)^{-1}\xv_t}{1+\xv_t^{\tran}(\vec{X}_{t-1}^{\tran}\vec{X}_{t-1}+\lambda \mat{I}_d)^{-1}\xv_t}\cdot (v_t-\iprod{\xv_t,\bv_t})^2\enspace.
    \end{align*}
\end{lemma}

\begin{proof}
    For any $t\in[T]$, denote $\vec{V}_{t}=(v_1,\ldots,v_{t})^{\tran}$.
    We further introduce the following notations:
    \begin{align*}
        \mat{W}=&\xM_{t-1}^{\tran}\xM_{t-1}+\lambda\mat{I}_d\enspace,\\
        h=&\xv_{t}^{\tran}\mat{W}^{-1}\xv_t\enspace,\\
        \vec{\alpha}=&\xM_{t-1}\mat{W}^{-1}\xv_t\enspace,\\
        \mat{H}=&\xM_{t-1}\mat{W}^{-1}\xM_{t-1}^{\tran}\enspace.
    \end{align*}
    We first derive the explicit form of $G_t(\bv_t)$.
    Note that the ridge estimator $\bv_t$ can be explicitly calculated as $\bv_t=(\xM_{t-1}^{\tran}\xM_{t-1}+\lambda\mat{I}_d)^{-1}\xM_{t-1}^{\tran}\vec{V}_{t-1}=\mat{W}^{-1}\xM_{t-1}^{\tran}\vec{V}_{t-1}$.
    Hence,
    \begin{align}\label{lemma:ridge_eq1}
        G_t(\bv_t)=&\|\vec{V}_{t-1}-\xM_{t-1}\bv_t\|^2+\lambda\|\bv_t\|^2\notag\\
        =&\vec{V}_{t-1}^{\tran}\vec{V}_{t-1}-2\vec{V}_{t-1}^{\tran}\xM_{t-1}\bv_t+\underbrace{\bv_t^{\tran}\xM_{t-1}^{\tran}\xM_{t-1}\bv_t+\lambda\cdot \bv_t^{\tran}\bv_t}_{=\bv_t^{\tran}\vec{W}\bv_t}\notag\\
        =&\vec{V}_{t-1}^{\tran}\vec{V}_{t-1}-2\vec{V}_{t-1}^{\tran}\xM_{t-1}\mat{W}^{-1}\xM_{t-1}^{\tran}\vec{V}_{t-1}+\vec{V}_{t-1}^{\tran}\xM_{t-1}\vec{W}^{-1}\vec{W}\mat{W}^{-1}\xM_{t-1}^{\tran}\vec{V}_{t-1}\notag\\
        =&\vec{V}_{t-1}^{\tran}\vec{V}_{t-1}-\vec{V}_{t-1}^{\tran}\xM_{t-1}\mat{W}^{-1}\xM_{t-1}^{\tran}\vec{V}_{t-1}\notag\\
        =&\vec{V}_{t-1}^{\tran}(\mat{I}_{t-1}-\xM_{t-1}\mat{W}^{-1}\xM_{t-1}^{\tran})\vec{V}_{t-1}\enspace.
    \end{align}
    Similarly, we can calculate $G_{t+1}(\bv_{t+1})$ as
    \begin{align}\label{lemma:ridge_eq2}
        G_{t+1}(\bv_{t+1})=\vec{V}_{t}^{\tran}(\mat{I}_t-\xM_{t}(\mat{W}+\xv_t\xv_t^{\tran})^{-1}\xM_{t}^{\tran})\vec{V}_{t}\enspace.
    \end{align}
    Using the block form of vector $\vec{V}_{t}$ and matrix $\vec{X}_t$, the right-hand side of \eqref{lemma:ridge_eq2} can be rewritten as:
    \begin{align}\label{lemma:ridge_eq3}
        \begin{bmatrix}
        \vec{V}_{t-1}^{\tran}&
        v_t
    \end{bmatrix}\left(\mat{I}_t-\begin{bmatrix}
        \xM_{t-1}\\
        \xv_{t}^{\tran}
    \end{bmatrix}(\mat{W}+\xv_t\xv_t^{\tran})^{-1}\begin{bmatrix}
        \xM_{t-1}^{\tran}& \xv_{t}
    \end{bmatrix}\right)\begin{bmatrix}
         \vec{V}_{t-1}\\
        v_t
    \end{bmatrix}\enspace.
    \end{align}
    For simplicity, suppose we have the following block form for the matrix:
    \begin{align*}
    \mat{I}_t-\begin{bmatrix}
        \xM_{t-1}\\
        \xv_{t}^{\tran}
    \end{bmatrix}(\mat{W}+\xv_t\xv_t^{\tran})^{-1}\begin{bmatrix}
        \xM_{t-1}^{\tran}& \xv_{t}
    \end{bmatrix}:=\begin{bmatrix}
        \mat{R}_{1,1}&\vec{R}_{1,2}\\
        \vec{R}_{1,2}^{\tran}&R_{2,2}
    \end{bmatrix}\enspace.
    \end{align*}
    By Sherman-Morrison formula, we have
    \begin{align}\label{lemma:ridge_eq4}
        (\mat{W}+\xv_t\xv_t^{\tran})^{-1}=\mat{W}^{-1}-\frac{\mat{W}^{-1}\xv_t\xv_t^{\tran}\mat{W}^{-1}}{1+\xv_{t}^{\tran}\mat{W}^{-1}\xv_t}=\mat{W}^{-1}-\frac{\mat{W}^{-1}\xv_t\xv_t^{\tran}\mat{W}^{-1}}{1+h}\enspace.
    \end{align}
    Now we calculate the explicit form for each block:
    \begin{enumerate}
    \item[(1)] By \eqref{lemma:ridge_eq4}, $\mat{R}_{1,1}$ can be calculated as:
    \begin{align*}
        \mat{R}_{1,1}=&\vec{I}_{t-1}-\xM_{t-1}\mat{W}^{-1}\xM_{t-1}^{\tran}+\frac{\xM_{t-1}\mat{W}^{-1}\xv_t\xv_t^{\tran}\mat{W}^{-1}\xM_{t-1}^{\tran}}{1+h}=\vec{I}_{t-1}-\mat{H}+\frac{\vec{\alpha}\vec{\alpha}^{\tran}}{1+h}\enspace.
    \end{align*}
    \item[(2)] By \eqref{lemma:ridge_eq4}, $\vec{R}_{1,2}$ can be calculated as:
    \begin{align*}
        \vec{R}_{1,2}=-\xM_{t-1}\mat{W}^{-1}\xv_t+\frac{\xM_{t-1}\mat{W}^{-1}\xv_t\xv_t^{\tran}\mat{W}^{-1}\xv_t}{1+h}=-\vec{\alpha}+\frac{h}{1+h}\vec{\alpha}=-\frac{1}{1+h}\vec{\alpha}\enspace.
    \end{align*}
    \item[(3)] By \eqref{lemma:ridge_eq4}, $R_{2,2}$ can be calculated as:
    \begin{align*}
        R_{2,2}=1-\xv_t^{\tran}\mat{W}^{-1}\xv_t+\frac{\xv_t^{\tran}\mat{W}^{-1}\xv_t\xv_t^{\tran}\mat{W}^{-1}\xv_t}{1+h}=1-h+\frac{h^2}{1+h}=\frac{1}{1+h}\enspace.
    \end{align*}
    \end{enumerate}
    Hence, we can simplify \eqref{lemma:ridge_eq3} as
    \begin{align*}
        &\begin{bmatrix}
        \vec{V}_{t-1}^{\tran}&
        v_t
    \end{bmatrix}\begin{bmatrix}
        \mat{R}_{1,1}&\vec{R}_{1,2}\\
        \vec{R}_{1,2}^{\tran}&R_{2,2}
    \end{bmatrix}\begin{bmatrix}
         \vec{V}_{t-1}\\
        v_t
    \end{bmatrix}\\
    =&\begin{bmatrix}
        \vec{V}_{t-1}^{\tran}&
        v_t
    \end{bmatrix}\begin{bmatrix}
        (\mat{I}_{t-1}-\mat{H})+\frac{\vec{\alpha}\vec{\alpha}^{\tran}}{1+h}&-\frac{1}{1+h}\vec{\alpha}\\
        -\frac{1}{1+h}\vec{\alpha}^{\tran}&\frac{1}{1+h}\\
    \end{bmatrix}\begin{bmatrix}
         \vec{V}_{t-1}\\
        v_t
    \end{bmatrix}\notag\\
    =&\begin{bmatrix}
        \vec{V}_{t-1}^{\tran}&
        v_t
    \end{bmatrix}\left[\begin{bmatrix}
        \mat{I}_{t-1}-\mat{H}&\vec{0}\\
        \vec{0}&0
    \end{bmatrix}+\frac{1}{1+h}\cdot\begin{bmatrix}
         -\vec{\alpha}\\
        1
        \end{bmatrix}\begin{bmatrix}
        -\vec{\alpha}^{\tran}&1
    \end{bmatrix}\right]\begin{bmatrix}
         \vec{V}_{t-1}\\
        v_t
    \end{bmatrix}\notag\\
    =&\vec{V}_{t-1}^{\tran}(\mat{I}_{t-1}-\mat{H})\vec{V}_{t-1}+\frac{1}{1+h}(v_t-\vec{\alpha}^{\tran}\vec{V}_{t-1})^2\enspace.
    \end{align*}
    This, together with \eqref{lemma:ridge_eq1} implies that
    \begin{align*}
        G_{t+1}(\bv_{t+1})=&\vec{V}_{t-1}^{\tran}(\mat{I}_{t-1}-\xM_{t-1}\mat{W}^{-1}\xM_{t-1}^{\tran})\vec{V}_{t-1}+\frac{1}{1+h}(v_t-\xv_t^{\tran}\underbrace{\vec{W}^{-1}\xM_{t-1}^{\tran}\vec{V}_{t-1}}_{=\bv_t})^2\\
        =&G_{t}(\bv_{t})+\frac{1}{1+h}(v_t-\iprod{\xv_t,\bv_t})^2,
    \end{align*}
    which indicates that
    \begin{align*}
        G_{t+1}(\bv_{t})-G_{t+1}(\bv_{t+1})=&(v_t-\iprod{\xv_t,\bv_t})^2+G_{t}(\bv_{t})-G_{t+1}(\bv_{t+1})\\
        =&\frac{h}{1+h}(v_t-\iprod{\xv_t,\bv_t})^2\\
        =&\frac{\xv_t^{\tran}(\vec{X}_{t-1}^{\tran}\vec{X}_{t-1}+\lambda \mat{I}_d)^{-1}\xv_t}{1+\xv_t^{\tran}(\vec{X}_{t-1}^{\tran}\vec{X}_{t-1}+\lambda \mat{I}_d)^{-1}\xv_t}\cdot (v_t-\iprod{\xv_t,\bv_t})^2\enspace.
    \end{align*}
    This completes the proof.
\end{proof}

Lemma \ref{lemma:squared-residual-deterministic} provides a finite-sample upper bound on the sum of squared full-information online residuals for any arbitrary outcomes $\{v_t:t\in[T]\}$ with bounded empirical second moment.

\begin{lemma}\label{lemma:squared-residual-deterministic}
    Denote $\zeta_1:=c_1^2(1+c_3(\max\{c_2,1\})^{1/2})(1+c_2c_3(\max\{c_2,1\})^{1/2})>0$.
    Under Assumptions \mainref{assumption:covariate-regularity}-\mainref{assumption:maximum-radius}, for any $\{v_t:t\in[T]\}$ such that $\sum_{t=1}^Tv_t^2\leq c_1^2T$, it holds that $\sum_{t=1}^T \paren{ v_t - \iprod{ \xv_t, \optbv_{v,t} } }^2\leq \zeta_1 T$ for any $T$.
    Here $\optbv_{v,t}:=\argmin_{\bv\in\mathbb{R}^d}\sum_{s<t}(v_s-\iprod{ \xv_s, \bv })^2+\eta_t^{-1}\|\bv\|^2$.
\end{lemma}

\begin{proof}
    For any $t\in[T]$, denote $\Delta_{t}(\bv)=v_t-\iprod{\xv_t,\bv}$.
	Denote $\widetilde{G}_{t}(\bv):=\sum_{s<t}(v_s-\iprod{\xv_s,\bv})^2+\eta_{t-1}^{-1}\|\bv\|^2$ and $\widetilde{\bv}^*_{v,t}$ as its minimizer.
    Further denote $\bv_v^*=\arg\min_{\bv}\sum_{t\in[T]}(v_t-\iprod{\xv_t,\bv})^2$ as the OLS estimator.
    We first derive an upper bound on the difference between the online squared residuals $\sum_{t=1}^T \paren{ v_t - \iprod{ \xv_t, \optbv_{v,t} } }^2$ and the sum of squared residuals for OLS estimator, i.e., $\sum_{t=1}^{T}(v_t-\iprod{\xv_t,\bv^*_v})^2$.
    By Corollary \ref{corollary:standard-FTRL-upper-bound}, we have
	\begin{align*}
		&\sum_{t=1}^{T}\left(v_t-\iprod{\xv_t,\bv_{v,t}^*}\right)^2-\sum_{t=1}^{T}\left(v_t-\iprod{\xv_t,\bv^*_v}\right)^2\notag\\
		\leq&\eta_{T+1}^{-1}\|\bv^*_v\|^2+\sum_{t=1}^T(\widetilde{G}_{t+1}(\bv^*_{v,t})-\widetilde{G}_{t+1}(\widetilde{\bv}^*_{v,t+1}))\quad(\text{Corollary \ref{corollary:standard-FTRL-upper-bound}})\notag\\
		\intertext{By letting $\lambda=\eta_{t}^{-1}$ in Lemma \ref{lemma:ridge} and using the definition $\Pi_{t,t}=\xv_t^{\tran}(\xM_{t-1}^{\tran}\xM_{t-1}+\eta_{t}^{-1}\vec{I}_d)^{-1}\xv_t$, this leads to}
        \leq&c_1^2c_2\eta_{T}^{-1}+\sum_{t=1}^{T}\frac{\Pi_{t,t}}{1+\Pi_{t,t}}\cdot \left(v_t-\iprod{\xv_t,\bv_{v,t}^*}\right)^2\quad(\text{Lemma \ref{lemma:ols-norm}})\enspace,
	\end{align*}
    which implies that
    \begin{align*}
        &\sum_{t=1}^{T}\frac{1}{1+\Pi_{t,t}}\left(v_t-\iprod{\xv_t,\bv_{v,t}^*}\right)^2\\
        \leq& c_1^2c_2T^{1/2}R_T+\sum_{t=1}^{T}\left(v_t-\iprod{\xv_t,\bv^*_v}\right)^2\\
        \leq&c_1^2c_2c_3(\max\{c_2,1\})^{1/2}T+\sum_{t=1}^{T}\left(v_t-\iprod{\xv_t,\bv^*_v}\right)^2\quad(\text{by Assumption \mainref{assumption:maximum-radius} and Lemma \ref{lemma:R}})\\
        \leq&c_1^2c_2c_3(\max\{c_2,1\})^{1/2}T+\sum_{t=1}^{T}v_t^2\quad(\text{optimality of $\bv^*_v$})\\
        \leq&c_1^2(1+c_2c_3(\max\{c_2,1\})^{1/2})T\enspace.
    \end{align*}
    On the other hand, for any $t\in[T]$, we have
    \begin{align*}
        1+\Pi_{t,t}\leq 1+R_t^2\eta_t=1+T^{-1/2}R_t\leq 1+T^{-1/2}R_T\leq 1+c_3(\max\{c_2,1\})^{1/2}\enspace,
    \end{align*}
    where the last inequality is due to Assumption \mainref{assumption:maximum-radius} and Lemma \ref{lemma:R}.
    Hence we have proved that
    \begin{align*}
        \sum_{t=1}^{T}\left(v_t-\iprod{\xv_t,\bv_{v,t}^*}\right)^2\leq c_1^2(1+c_3(\max\{c_2,1\})^{1/2})(1+c_2c_3(\max\{c_2,1\})^{1/2})T\enspace,
    \end{align*}
    which completes the proof.
\end{proof}

The following corollary provides a uniform upper bound on the full-information squared residuals $A_T^*(k)$, which is a direct consequence of Lemma \ref{lemma:squared-residual-deterministic}.

\begin{corollary}\label{corollary:squared-residual-deterministic}
    Under Assumptions \mainref{assumption:moments}-\mainref{assumption:maximum-radius}, for $k\in\{0,1\}$, it holds that $A_T^*(k)= \sum_{t=1}^T \paren{ y_t(k) - \iprod{ \xv_t, \optbv_t(k) } }^2\leq \zeta_1 T$.
\end{corollary}

\begin{proof}
    By the Cauchy-Schwarz inequality and Assumption \mainref{assumption:moments}, 
    \begin{align*}
        \sum_{t=1}^Ty_t^2(k)\leq T^{1/2}\left(\sum_{t=1}^Ty_t^4(k)\right)^{1/2}\leq c_1^2T\enspace.
    \end{align*}
    Therefore, by Lemma \ref{lemma:squared-residual-deterministic}, the result is proved.
\end{proof}

Given the explicit expressions for the estimated squared residuals $\E{\widehat{A}_t(1)}$ and $\E{\widehat{A}_t(0)}$ derived in Lemma \mainref{lemma:expectation-of-estores}, the moment bounds for the inverse assignment probabilities in Corollary \mainref{corollary:p-moment}*, and the bounds on the associated deterministic sums established in Lemma \ref{lemma:deterministic-summation-1} and Corollary \ref{corollary:squared-residual-deterministic}, we are now prepared to derive a uniform upper bound for $\E{\widehat{A}_t(1)}$ and $\E{\widehat{A}_t(0)}$, formally stated in the following lemma.

\begin{lemma}\label{lemma:squared-residuals-random}
    Under Assumptions \mainref{assumption:moments}-\mainref{assumption:maximum-radius} and Condition \mainref{condition:sigmoid}, for any $t\in[T]$, the expectations of the estimated squared residuals can be bounded as:
    \begin{align*}
        \max\left\{\E{\widehat{A}_t(1)},\E{\widehat{A}_t(0)}\right\}\leq \varsigma\cdot T\enspace,
    \end{align*}
    where $\varsigma>0$ is defined as: 
    \begin{align*}
        \varsigma:=&2\max\big\{\zeta_1+c\,c_1^2c_3^{3/2}\gamma(\max\{c_2,1\})^{3/4}\xi_2^{1/2},\\
        &\quad\quad \quad 2^{1/3}(c\,c_1^2c_3^{5/4}\gamma(\max\{c_2,1\})^{5/8}\xi_2^{1/2}\bb_1(\bb_2/6)^{1/4})^{4/3}\big\}\enspace.
    \end{align*}
\end{lemma}

\begin{proof}
    Let $\zeta_2:=c\,c_1^2c_3^{5/4}\gamma(\max\{c_2,1\})^{5/8}\xi_2^{1/2}\bb_1(\bb_2/6)^{1/4}$ and $\zeta_3:=c\,c_1^2c_3^{3/2}\gamma(\max\{c_2,1\})^{3/4}\xi_2^{1/2}$.
    Further let $\kappa$ be the largest positive root of the following equation:
    \begin{align*}
        \kappa=(\zeta_1+\zeta_3) T+\zeta_2 T^{3/4}\kappa^{1/4}\enspace.
    \end{align*}
    \textbf{Step 1: }Prove the existence of $\kappa$ and show $\kappa\geq c_1^2T$.\\[3mm]
    Note that for any fixed $T\geq 1$, $(\zeta_1+\zeta_3) T+\zeta_2 T^{3/4}\kappa^{1/4}-\kappa$ is positive when $\kappa=0$ and tends to $-\infty$ when $\kappa$ goes to infinity.
    This guarantees the existence of such $\kappa>0$ by continuity.
    Since $\zeta_1\geq c_1^2$, we have $\kappa\geq c_1^2T$ by its definition.\\[3mm]
    \textbf{Step 2: }Prove that $\kappa\leq 2\max\{\zeta_1+\zeta_3,2^{1/3}\zeta_2^{4/3}\}T=\varsigma\cdot T$.\\[3mm]
    When $\kappa> 2\max\{\zeta_1+\zeta_3,2^{1/3}\zeta_2^{4/3}\}T$, it is straightforward to check that $\frac{1}{2}\kappa> (\zeta_1+\zeta_3) T$ and $\frac{1}{2}\kappa> \zeta_2 T^{3/4}\kappa^{1/4}$ for any $T\geq 1$.
    Thus $\kappa> (\zeta_1+\zeta_3) T+\zeta_2 T^{3/4}\kappa^{1/4}$ when $\kappa> 2\max\{\zeta_1+\zeta_3,2^{1/3}\zeta_2^{4/3}\}T$, which implies that $\kappa\leq 2\max\{\zeta_1+\zeta_3,2^{1/3}\zeta_2^{4/3}\}T=\varsigma\cdot T$ by the definition of $\kappa$.\\[3mm]
    \textbf{Step 3: }Prove that $\max\{\E{\widehat{A}_t(1)},\E{\widehat{A}_t(0)}\}\leq \kappa$ for any $t\in[T]$.\\[3mm]
    We then use induction to prove that $\max\{\E{\widehat{A}_t(1)},\E{\widehat{A}_t(0)}\}\leq \kappa$ for any $t\in[T]$.
    For $t=1$, we have $\e[\widehat{A}_1(1)]=y_1^2(1)\leq c_1^2T^{1/2}\leq c_1^2T\leq \kappa$ and $\e[\widehat{A}_1(0)]=A_1^*(0)=y_1^2(0)\leq c_1^2T^{1/2}\leq c_1^2T\leq \kappa$.
    Hence the result holds.
    Assume the result holds for $1,\ldots,t-1$.
    We then prove the result for $t$. 
    By Lemma \mainref{lemma:expectation-of-estores}, we have
    \begin{align}\label{lemma:squared-residuals-random_eq1}
        \e[\widehat{A}_t(1)]=&A_{t}^*(1)+\sum_{s=1}^t\sum_{r=1}^{s-1}\Pi_{s,r}^2y_r(1)^2\e\left[\frac{1}{p_r}-1\right]\notag\\
        \intertext{Since $A_{t}^*(1)\leq A_{T}^*(1)\leq \zeta_1T$ by Corollary \ref{corollary:squared-residual-deterministic}, using the upper bound on $\e\left[\frac{1}{p_r}-1\right]$ in Corollary \mainref{corollary:p-moment}*, we have}
        \leq&\zeta_1 T+\sum_{s=1}^t\sum_{r=1}^{s-1}\Pi_{s,r}^2y_r(1)^2\left(1+\bb_1(\bb_2/6)^{1/4}\eta_r^{1/4}\e[\widehat{A}_{r-1}(0)]^{1/4}\right)\notag\\
        \leq&\zeta_1 T+\sum_{s=1}^t\sum_{r=1}^{s-1}\Pi_{s,r}^2y_r(1)^2\left(1+\bb_1(\bb_2/6)^{1/4}\eta_r^{1/4}\kappa^{1/4}\right)\quad(\text{induction assumption})\notag\\
        =&\zeta_1 T+\underbrace{\sum_{s=1}^t\sum_{r=1}^{s-1}\Pi_{s,r}^2y_r(1)^2}_{:=S_1}+\bb_1(\bb_2/6)^{1/4}\kappa^{1/4}\underbrace{\sum_{s=1}^t\sum_{r=1}^{s-1}\Pi_{s,r}^2y_r(1)^2\eta_r^{1/4}}_{:=S_2}\enspace.
    \end{align}
    By Lemma \ref{lemma:deterministic-summation-1}  and Lemma \ref{lemma:R}, $S_1$, $S_2$ can be bounded as:
    \begin{align}\label{lemma:squared-residuals-random_eq2}
        S_1\leq& c\,c_1^2\gamma \xi_2^{1/2}R_t^{3/2}T^{1/4}\quad(\text{Lemma \ref{lemma:deterministic-summation-1}-3})\notag\\
        \leq&c\,c_1^2\gamma \xi_2^{1/2}R_T^{3/2}T^{1/4}\notag\\
        \leq&c\,c_1^2\gamma \xi_2^{1/2}(\max\{c_2,1\})^{3/4}R^{3/2}T^{1/4}\quad(\text{Lemma \ref{lemma:R}})\notag\\
        \leq&c\,c_1^2c_3^{3/2}\gamma \xi_2^{1/2}(\max\{c_2,1\})^{3/4}T\enspace,\quad(\text{Assumption \mainref{assumption:maximum-radius}})\notag\\
        S_2=&T^{-1/8}\sum_{s=1}^t\sum_{r=1}^{s-1}R_r^{-1/4}\Pi_{s,r}^2y_r(1)^2\notag\\
        \leq& c\,c_1^2\gamma \xi_2^{1/2}R_t^{5/4}T^{1/4-1/8}\quad(\text{Lemma \ref{lemma:deterministic-summation-1}-3})\notag\\
        \leq&c\,c_1^2c_3^{5/4}\gamma \xi_2^{1/2}(\max\{c_2,1\})^{5/8}T^{3/4}\quad(\text{Lemma \ref{lemma:R} and Assumption \mainref{assumption:maximum-radius}})\enspace.
    \end{align}
    By \eqref{lemma:squared-residuals-random_eq1} and \eqref{lemma:squared-residuals-random_eq2}, we have
    \begin{align*}
        \e[\widehat{A}_t(1)]\leq&\zeta_1 T+c\,c_1^2c_3^{5/4}\gamma(\max\{c_2,1\})^{5/8}\xi_2^{1/2}\cdot\bb_1(\bb_2/6)^{1/4}T^{3/4}\kappa^{1/4}+c\,c_1^2c_3^{3/2}\gamma(\max\{c_2,1\})^{3/4}\xi_2^{1/2}T\\ 
        =&\zeta_1 T+\zeta_2 T^{3/4}\kappa^{1/4}+\zeta_3 T\\
        =&\kappa\enspace.
    \end{align*}
    Similarly, we can prove that $\e[\widehat{A}_t(0)]\leq \kappa$.
    The claim is thus proved by induction and the bound on $\kappa$ derived in Step 2.
\end{proof}

Using Corollary \mainref{corollary:p-moment}* and Lemma \ref{lemma:squared-residuals-random}, we can now control the moments of the inverse assignment probabilities.
The result is formally stated in the following lemma.

\begin{lemma}\label{lemma:p-power-moment}
    Denote $\chi:=2c_3^{1/4}(\max\{c_2,1\})^{1/8}+\bb_1(\bb_2/6)^{1/4}\varsigma^{1/4}>0$.
    Under Assumptions \mainref{assumption:moments}-\mainref{assumption:maximum-radius} and Condition \mainref{condition:sigmoid}, for any $t\in[T]$ and $0\leq k\leq 4$, it holds that:
    \begin{align*}
        \max\left\{\e\left[\frac{1}{p_{t}^k}\right],\e\left[\frac{1}{(1-p_{t})^k}\right]\right\}\leq \chi^kR_t^{-k/4}T^{k/8}\enspace.
     \end{align*}
\end{lemma}

\begin{proof}
    Without loss of generality, we only bound $\e[1/p_{t}^k]$.
    By Corollary \mainref{corollary:p-moment}* and Lemma \ref{lemma:squared-residuals-random}, we have
    \begin{align*}
        \e\left[\frac{1}{p_{t}^k}\right]\leq\left(2+\bb_1(\bb_2/6)^{1/4}\eta_t^{1/4}\left(\varsigma\,T\right)^{1/4}\right)^{k}=\left(2+\bb_1(\bb_2/6)^{1/4}\varsigma^{1/4}R_t^{-1/4}T^{1/8}\right)^{k}\enspace.
    \end{align*}
    Note that by Lemma \ref{lemma:R} and Assumption \mainref{assumption:maximum-radius}, $R_t^{-1/4}T^{1/8}\geq R_T^{-1/4}T^{1/8}\geq c_3^{-1/4}(\max\{c_2,1\})^{-1/8}$. Thus
    \begin{align*}
        \left(2+\bb_1(\bb_2/6)^{1/4}\varsigma^{1/4}R_t^{-1/4}T^{1/8}\right)^{k}\leq& \left(2c_3^{1/4}(\max\{c_2,1\})^{1/8}+\bb_1(\bb_2/6)^{1/4}\varsigma^{1/4}\right)^{k}R_t^{-k/4}T^{k/8}\\
        =&\chi^kR_t^{-k/4}T^{k/8}\enspace,
    \end{align*}
    which completes the proof.
\end{proof}

\subsection{Prediction Regret}\label{section:C4}
In this section, we derive an upper bound on the expected prediction regret.
The following lemma establishes the connection between the defined prediction regret and the corresponding terms that are used in the proposed algorithm.

\begin{reflemma}{\mainref{lemma:pred-regret-unbiased}}
    \predregretunbiased
\end{reflemma}

\begin{proof}
    For $k\in\{0,1\}$, denote $\Delta_{t,k}(\bv)=y_t(k)-\iprod{\xv_t,\bv}$.
    Further denote $\widehat{\Delta}_{t,1}(\bv)=y_t(1)\cdot \frac{\indicator{Z_t=1}}{p_t}-\iprod{\xv_t,\bv}$ and $\widehat{\Delta}_{t,0}(\bv)=y_t(0)\cdot \frac{\indicator{Z_t=0}}{1-p_t}-\iprod{\xv_t,\bv}$.
    By the definition of $\widehat{\ell}_t$ and $\ell_t$, we have
    \begin{align}\label{lemma:pred-regret-unbiased_eq1}
        &\e[\widehat{\ell}_t(\bv_t(1),\bv_t(0))-\widehat{\ell}_t(\bv^*(1),\bv^*(0))|\filt_{t-1}]\notag\\
        =&\frac{\mathcal{E}(0)}{\mathcal{E}(1)}\cdot\underbrace{\e\left[\widehat{\Delta}_{t,1}(\bv_t(1))^2|\filt_{t-1}\right]}_{:=S_1}+\frac{\mathcal{E}(1)}{\mathcal{E}(0)}\cdot\underbrace{\e\left[\widehat{\Delta}_{t,0}(\bv_t(0))^2|\filt_{t-1}\right]}_{:=S_2}\notag\\
        &+2\underbrace{\e\left[\widehat{\Delta}_{t,1}(\bv_t(1))\cdot\widehat{\Delta}_{t,0}(\bv_t(0))|\filt_{t-1}\right]}_{:=S_3}-\frac{\mathcal{E}(0)}{\mathcal{E}(1)}\cdot\underbrace{\e\left[\widehat{\Delta}_{t,1}(\bv^*(1))^2|\filt_{t-1}\right]}_{:=S_4}\notag\\
        &-\frac{\mathcal{E}(1)}{\mathcal{E}(0)}\cdot\underbrace{\e\left[\widehat{\Delta}_{t,0}(\bv^*(0))^2|\filt_{t-1}\right]}_{:=S_5}-2\underbrace{\e\left[\widehat{\Delta}_{t,1}(\bv^*(1))\cdot\widehat{\Delta}_{t,0}(\bv^*(0))|\filt_{t-1}\right]}_{:=S_6}\enspace.
    \end{align}
    Using the equality $\widehat{\Delta}_{t,1}(\bv_t(1))=\Delta_{t,1}(\bv_t(1))+y_t(1)\cdot\left[\frac{\indicator{Z_t=1}}{p_t}-1\right]$, we can simplify $S_1$ as
    \begin{align}\label{lemma:pred-regret-unbiased_eq2}
        S_1=&\e\left[\widehat{\Delta}_{t,1}(\bv_t(1))^2|\filt_{t-1}\right]\notag\\
        =&\e\Bigg[\left(y_t(1)\cdot\left[\frac{\indicator{Z_t=1}}{p_t}-1\right]+\Delta_{t,1}(\bv_t(1))\right)^2\Big|\filt_{t-1}\Bigg]\notag\\
        =&\e\Bigg[y_t(1)^2\cdot\left[\frac{\indicator{Z_t=1}}{p_t}-1\right]^2\Big|\filt_{t-1}\Bigg]+\e\left[\Delta_{t,1}(\bv_t(1))^2|\filt_{t-1}\right]\notag\\
        &+2\e\Bigg[y_t(1)\cdot\left[\frac{\indicator{Z_t=1}}{p_t}-1\right]\cdot\Delta_{t,1}(\bv_t(1))\Big|\filt_{t-1}\Bigg]\notag\\
        \intertext{Since $\Delta_{t,1}(\bv_t(1))$ is measurable with respect to $\filt_{t-1}$, this simplifies to}
        =&\e\Bigg[y_t(1)^2\cdot\left[\frac{\indicator{Z_t=1}}{p_t}-1\right]^2\Big|\filt_{t-1}\Bigg]+\Delta_{t,1}(\bv_t(1))^2\notag\\
        &+2y_t(1)\cdot \Delta_{t,1}(\bv_t(1))\cdot\e\Bigg[\frac{\indicator{Z_t=1}}{p_t}-1\Big|\filt_{t-1}\Bigg]\notag\\
        =&\e\Bigg[y_t(1)^2\cdot\left[\frac{\indicator{Z_t=1}}{p_t}-1\right]^2\Big|\filt_{t-1}\Bigg]+\Delta_{t,1}(\bv_t(1))^2\enspace.
    \end{align}
    Similarly, we can prove that
    \begin{align}\label{lemma:pred-regret-unbiased_eq3}
        S_2=&\e\Bigg[y_t(0)^2\cdot\left[\frac{\indicator{Z_t=0}}{1-p_t}-1\right]^2\Big|\filt_{t-1}\Bigg]+\Delta_{t,0}(\bv_t(0))^2\enspace,\notag\\
        S_3=&\e\Bigg[y_t(1)y_t(0)\cdot\left[\frac{\indicator{Z_t=1}}{p_t}-1\right]\left[\frac{\indicator{Z_t=0}}{1-p_t}-1\right]\Big|\filt_{t-1}\Bigg]\notag\\
        &+\Delta_{t,1}(\bv_t(1))\cdot\Delta_{t,0}(\bv_t(0))\enspace,\notag\\
        S_4=&\e\Bigg[y_t(1)^2\cdot\left[\frac{\indicator{Z_t=1}}{p_t}-1\right]^2\Big|\filt_{t-1}\Bigg]+\Delta_{t,1}(\bv^*(1))^2\enspace,\notag\\
        S_5=&\e\Bigg[y_t(0)^2\cdot\left[\frac{\indicator{Z_t=0}}{1-p_t}-1\right]^2\Big|\filt_{t-1}\Bigg]+\Delta_{t,0}(\bv^*(0))^2\enspace,\notag\\
        S_6=&\e\Bigg[y_t(1)y_t(0)\cdot\left[\frac{\indicator{Z_t=1}}{p_t}-1\right]\left[\frac{\indicator{Z_t=0}}{1-p_t}-1\right]\Big|\filt_{t-1}\Bigg]\notag\\
        &+\Delta_{t,1}(\bv^*(1))\cdot \Delta_{t,0}(\bv^*(0))\enspace.
    \end{align}
    By \eqref{lemma:pred-regret-unbiased_eq1}, \eqref{lemma:pred-regret-unbiased_eq2} and \eqref{lemma:pred-regret-unbiased_eq3}, the conditional expectations cancel between $S_1,S_2,S_3$ and $S_4,S_5,S_6$.
    Hence we have
    \begin{align*}
        &\e[\widehat{\ell}_t(\bv_t(1),\bv_t(0))-\widehat{\ell}_t(\bv^*(1),\bv^*(0))|\filt_{t-1}]\\
        =&\frac{\mathcal{E}(0)}{\mathcal{E}(1)}S_1+\frac{\mathcal{E}(1)}{\mathcal{E}(0)}S_2+2S_3-\frac{\mathcal{E}(0)}{\mathcal{E}(1)}S_4-\frac{\mathcal{E}(1)}{\mathcal{E}(0)}S_5-2S_6\\
        =&\frac{\mathcal{E}(0)}{\mathcal{E}(1)}\cdot \Delta_{t,1}(\bv_t(1))^2+\frac{\mathcal{E}(1)}{\mathcal{E}(0)}\cdot \Delta_{t,0}(\bv_t(0))^2+2\Delta_{t,1}(\bv_t(1))\cdot \Delta_{t,0}(\bv_t(0))\\
        &-\frac{\mathcal{E}(0)}{\mathcal{E}(1)}\cdot \Delta_{t,1}(\bv^*(1))^2-\frac{\mathcal{E}(1)}{\mathcal{E}(0)}\cdot \Delta_{t,0}(\bv^*(0))^2-2\Delta_{t,1}(\bv^*(1))\cdot \Delta_{t,0}(\bv^*(0))\\
        =&\ell_t(\bv_t(1),\bv_t(0))-\ell_t(\bv^*(1),\bv^*(0))\enspace.
    \end{align*}
    This completes the proof.
\end{proof}

Using Lemma~\mainref{lemma:pred-regret-unbiased}, the following lemma shows that the expected prediction regret can be bounded by a regularization term and the sum of expected sub-optimality gaps.
The result follows from an application of Corollary \ref{corollary:standard-FTRL-upper-bound}. 
It remains only to verify that the two ridge estimators in the individual groups jointly minimize the target function.

\begin{reflemma}{\mainref{lemma:standard-oco-for-pred-regret}}
	\standardoco
\end{reflemma}

\begin{proof}
    For simplicity, for any $t\in[T]$, we denote $\widehat{\vec{Y}}_t(1)=\left(y_1(1)\cdot\frac{\indicator{Z_1=1}}{p_1},\ldots,y_t(1)\cdot\frac{\indicator{Z_t=1}}{p_t}\right)^{\tran}$ and $\widehat{\vec{Y}}_t(0)=\left(y_1(0)\cdot\frac{\indicator{Z_1=0}}{1-p_1},\ldots,y_t(0)\cdot\frac{\indicator{Z_t=0}}{1-p_t}\right)^{\tran}$.
    Note that the regularizer term $m(\bv(1),\bv(0))$ is nonnegative:
    \begin{align*}
        m(\bv(1),\bv(0))=&\frac{\mathcal{E}(0)}{\mathcal{E}(1)}\|\bv(1)\|^2+\frac{\mathcal{E}(1)}{\mathcal{E}(0)}\|\bv(0)\|^2+2\iprod{\bv(1),\bv(0)}\\
        =&\left\|\sqrt{\frac{\olsres{0}}{\olsres{1}}}\cdot\bv(1)+\sqrt{\frac{\olsres{1}}{\olsres{0}}}\cdot\bv(0)\right\|^2\\
        \geq& 0
    \end{align*}
    and has minimum value 0.
    By Corollary \ref{corollary:standard-FTRL-upper-bound}, it suffices to prove that $(\bv_t(1),\bv_t(0))$ minimizes $\widehat{L}_t(\bv(1),\bv(0))$ and $(\widetilde{\bv}_t(1),\widetilde{\bv}_t(0))$ minimizes $\widetilde{L}_t(\bv(1),\bv(0))$.
    
    Without loss of generality, we only prove the first result.
    For any fixed $(\bv(1),\bv(0))$, we show that $\widehat{L}_t(\bv(1),\bv(0))\geq \widehat{L}_t(\bv_t(1),\bv_t(0))$.
    By the definition of $\widehat{L}_t(\bv(1),\bv(0))$, we can rewrite it as:
    \begin{align*}
        &\widehat{L}_{t}(\bv(1),\bv(0))\\
        =&\sum_{s=1}^{t-1}\Bigg(\sqrt{\frac{\olsres{0}}{\olsres{1}}}\cdot y_s(1)\cdot \frac{\indicator{Z_s=1}}{p_s}+\sqrt{\frac{\olsres{1}}{\olsres{0}}}\cdot y_s(0)\cdot \frac{\indicator{Z_s=0}}{1-p_s}\\
        &\quad\quad-\left\langle\xv_s,\sqrt{\frac{\olsres{0}}{\olsres{1}}}\cdot\bv(1)+\sqrt{\frac{\olsres{1}}{\olsres{0}}}\cdot\bv(0)\right\rangle\Bigg)^2+\eta_t^{-1}\left\|\sqrt{\frac{\olsres{0}}{\olsres{1}}}\cdot\bv(1)+\sqrt{\frac{\olsres{1}}{\olsres{0}}}\cdot\bv(0)\right\|^2\enspace.
    \end{align*}
    We construct the artificial outcomes $v_t$ and regression coefficients $\bv_{t}$, where for any $t\in[T]$,
    \begin{align*}
        v_t:=&\sqrt{\frac{\olsres{0}}{\olsres{1}}}\cdot Y_t\cdot \frac{\indicator{Z_t=1}}{p_t}+\sqrt{\frac{\olsres{1}}{\olsres{0}}}\cdot Y_t\cdot \frac{\indicator{Z_t=0}}{1-p_t}\enspace,\\
        \bv_{t}:=&\sqrt{\frac{\olsres{0}}{\olsres{1}}}\cdot \bv_{t}(1)+\sqrt{\frac{\olsres{1}}{\olsres{0}}}\cdot \bv_{t}(0)\enspace.
    \end{align*}
    Note that the ridge estimator is linear in the response variables.
    Since $\bv_{t}(1)$ and $\bv_{t}(0)$ are the ridge estimators for
    \begin{align*}
        \left\{Y_s\cdot \frac{\indicator{Z_s=1}}{p_s}:1\leq s\leq t-1\right\}\quad\text{and}\quad\left\{Y_s\cdot \frac{\indicator{Z_s=0}}{1-p_s}:1\leq s\leq t-1\right\}\enspace,
    \end{align*}
    this implies that
    \begin{align*}
        \bv_{t}=&\arg\min_{\bv\in\mathbb{R}^d}\sum_{s=1}^{t-1}(v_s-\iprod{\xv_s,\bv})^2+\eta_{t}^{-1}\|\bv\|^2\enspace.
    \end{align*}
    Hence this implies that $(\bv_t(1),\bv_t(0))$ minimizes $\widehat{L}_t(\bv(1),\bv(0))$. 
    Similarly we can prove that $(\widetilde{\bv}_t(1),\widetilde{\bv}_t(0))$ minimizes $\widetilde{L}_t(\bv(1),\bv(0))$.
    Then the result follows from Corollary \ref{corollary:standard-FTRL-upper-bound}.
\end{proof}

The first term in the upper bound provided in Lemma \mainref{lemma:standard-oco-for-pred-regret} can be controlled through Lemma \ref{lemma:ols-norm}.
The second term can be readily calculated with the help of Lemma \ref{lemma:ridge}, which is formally stated in the following lemma.

\begin{reflemma}{\mainref{lemma:ptpredloss-diff}}
	\predlossdiff
\end{reflemma}

\begin{proof}
    By the definition of $\widetilde{L}_{t+1}(\bv(1),\bv(0))$, we can rewrite it as:
    \begin{align*}
        &\widetilde{L}_{t+1}(\bv(1),\bv(0))\\
        =&\sum_{s=1}^t\Bigg(\sqrt{\frac{\olsres{0}}{\olsres{1}}}\cdot y_s(1)\cdot \frac{\indicator{Z_s=1}}{p_s}+\sqrt{\frac{\olsres{1}}{\olsres{0}}}\cdot y_s(0)\cdot \frac{\indicator{Z_s=0}}{1-p_s}\\
        &\quad\quad-\left\langle\xv_s,\sqrt{\frac{\olsres{0}}{\olsres{1}}}\cdot\bv(1)+\sqrt{\frac{\olsres{1}}{\olsres{0}}}\cdot\bv(0)\right\rangle\Bigg)^2+\eta_t^{-1}\left\|\sqrt{\frac{\olsres{0}}{\olsres{1}}}\cdot\bv(1)+\sqrt{\frac{\olsres{1}}{\olsres{0}}}\cdot\bv(0)\right\|^2\enspace.
    \end{align*}
    We construct the artificial outcomes $v_t$ and regression coefficients $\bv_{t}$, $\widetilde{\bv}_{t+1}$, where for any $t\in[T]$,
    \begin{align*}
        v_t:=&\sqrt{\frac{\olsres{0}}{\olsres{1}}}\cdot Y_t\cdot \frac{\indicator{Z_t=1}}{p_t}+\sqrt{\frac{\olsres{1}}{\olsres{0}}}\cdot Y_t\cdot \frac{\indicator{Z_t=0}}{1-p_t}\enspace,\\
        \bv_{t}:=&\sqrt{\frac{\olsres{0}}{\olsres{1}}}\cdot \bv_{t}(1)+\sqrt{\frac{\olsres{1}}{\olsres{0}}}\cdot \bv_{t}(0)\enspace,\\
        \widetilde{\bv}_{t+1}:=&\sqrt{\frac{\olsres{0}}{\olsres{1}}}\cdot \widetilde{\bv}_{t+1}(1)+\sqrt{\frac{\olsres{1}}{\olsres{0}}}\cdot \widetilde{\bv}_{t+1}(0)\enspace.
    \end{align*}
    Note that the ridge estimator is linear in the response variables, which implies that
    \begin{align*}
        \bv_{t}=&\arg\min_{\bv\in\mathbb{R}^d}\sum_{s=1}^{t-1}(v_s-\iprod{\xv_s,\bv})^2+\eta_{t}^{-1}\|\bv\|^2\enspace,\\
        \widetilde{\bv}_{t+1}=&\arg\min_{\bv\in\mathbb{R}^d}\sum_{s=1}^{t}(v_s-\iprod{\xv_s,\bv})^2+\eta_{t}^{-1}\|\bv\|^2\enspace.
    \end{align*}
    By letting $\lambda=\eta_t^{-1}$ in Lemma \ref{lemma:ridge}, we have
    \begin{align*}
        &\widetilde{L}_{t+1}(\bv_t(1),\bv_t(0))-\widetilde{L}_{t+1}(\widetilde{\bv}_{t+1}(1),\widetilde{\bv}_{t+1}(0))\\
        =&\left(\sum_{s=1}^{t}(v_s-\iprod{\xv_s,\bv_{t}})^2+\eta_{t}^{-1}\|\bv_{t}\|^2\right)-\left(\sum_{s=1}^{t}(v_s-\iprod{\xv_s,\widetilde{\bv}_{t+1}})^2+\eta_{t}^{-1}\|\widetilde{\bv}_{t+1}\|^2\right)\\
        =&\frac{\Pi_{t,t}}{1+\Pi_{t,t}}\cdot(v_t-\iprod{\xv_t,\bv_{t}})^2\\
        =&\frac{\Pi_{t,t}}{1+\Pi_{t,t}}\cdot\paren[\Bigg]{
	       \braces[\Big]{ Y_t \cdot \frac{\indicator{Z_t=1}}{p_t} - \iprod{ \xv_t , \bv_t(1) } } \cdot \sqrt{ \frac{\olsres{0}}{\olsres{1}} }
	   +
	   \braces[\Big]{ Y_t \cdot \frac{\indicator{Z_t=0}}{1-p_t} - \iprod{ \xv_t , \bv_t(0) } } \cdot \sqrt{ \frac{\olsres{1}}{\olsres{0}} }
        }^2\\
        =&\frac{\Pi_{t,t}}{1+\Pi_{t,t}}\cdot\widehat{\ell}_t(\bv_t(1),\bv_t(0))\enspace.
    \end{align*}
    This completes the proof.
\end{proof}

Using Lemma \mainref{lemma:standard-oco-for-pred-regret} and Lemma \mainref{lemma:ptpredloss-diff}, we can finally derive the upper bound on the expected prediction regret in the following proposition.

\begin{refproposition}{\mainref{prop:pred-regret}*}
	Under Assumptions \mainref{assumption:moments}-\mainref{assumption:maximum-radius} and Condition~\mainref{condition:sigmoid}, the expected prediction regret is bounded as
    \begin{align*}
        \E{\newpredregret_T}\leq
        \left(\frac{c_1}{c_0}+\frac{c_0}{c_1}+2\right)\Big(c_1^2c_2(\max\{c_2,1\})^{1/2}+c_1^2c_3^{1/4}\gamma(\max\{c_2,1\})^{5/8}\xi_2^{1/2}\chi+(\max\{c_2,1\})^{1/2}\varsigma\Big)T^{1/2}R\enspace.
    \end{align*}
\end{refproposition}

\begin{proof}
    For $k\in\{0,1\}$, denote $\Delta_{t,k}(\bv)=y_t(k)-\iprod{\xv_t,\bv}$.
    Further denote $\widehat{\Delta}_{t,1}(\bv)=y_t(1)\cdot \frac{\indicator{Z_t=1}}{p_t}-\iprod{\xv_t,\bv}$ and $\widehat{\Delta}_{t,0}(\bv)=y_t(0)\cdot \frac{\indicator{Z_t=0}}{1-p_t}-\iprod{\xv_t,\bv}$.
    Lemma \mainref{lemma:ptpredloss-diff} shows that the successive difference appearing in Lemma \mainref{lemma:standard-oco-for-pred-regret} is expressed in terms of $\widehat{\ell}_t$.
    We first approximate $\widehat{\ell}_t$ by $\ell_t$, which is easier to control via Lemma \ref{lemma:squared-residuals-random}.
    To relate $\widehat{\ell}_t$ and $\ell_t$, we explicitly expand their difference.
    By the definition of $\widehat{\ell}_t$, we obtain the following decomposition:
    \begin{align*}
        &\widehat{\ell}_t(\bv_t(1),\bv_t(0))\\
        =&\frac{\mathcal{E}(0)}{\mathcal{E}(1)}\cdot\widehat{\Delta}_{t,1}(\bv_t(1))^2+\frac{\mathcal{E}(1)}{\mathcal{E}(0)}\cdot\widehat{\Delta}_{t,0}(\bv_t(0))^2+2\cdot\widehat{\Delta}_{t,1}(\bv_t(1))\cdot\widehat{\Delta}_{t,0}(\bv_t(0))\\
        \intertext{By rewriting $\widehat{\Delta}_{t,1}(\bv_t(1))=y_t(1)\cdot \left(\frac{\indicator{Z_t=1}}{p_t}-1\right)+\Delta_{t,1}(\bv_t(1))$ and $\widehat{\Delta}_{t,0}(\bv_t(0))=y_t(0)\cdot \left(\frac{\indicator{Z_t=0}}{1-p_t}-1\right)+\Delta_{t,0}(\bv_t(0))$, we have the following expansion:}
        =&{\ell}_t(\bv_t(1),\bv_t(0))+\frac{\mathcal{E}(0)}{\mathcal{E}(1)}\cdot y_t(1)^2\left(\frac{\indicator{Z_t=1}}{p_t}-1\right)^2+\underbrace{\frac{2\mathcal{E}(0)}{\mathcal{E}(1)}\cdot y_t(1)\cdot\Delta_{t,1}(\bv_t(1))\cdot\left(\frac{\indicator{Z_t=1}}{p_t}-1\right)}_{:=D_1}\\
        &+\frac{\mathcal{E}(1)}{\mathcal{E}(0)}\cdot y_t(0)^2\left(\frac{\indicator{Z_t=0}}{1-p_t}-1\right)^2+\underbrace{\frac{2\mathcal{E}(1)}{\mathcal{E}(0)}\cdot y_t(0)\cdot\Delta_{t,0}(\bv_t(0))\cdot\left(\frac{\indicator{Z_t=0}}{1-p_t}-1\right)}_{:=D_2}\\
        &+2y_t(1)y_t(0)\left(\frac{\indicator{Z_t=1}}{p_t}-1\right)\left(\frac{\indicator{Z_t=0}}{1-p_t}-1\right)+\underbrace{2y_t(1)\left(\frac{\indicator{Z_t=1}}{p_t}-1\right)\cdot\Delta_{t,0}(\bv_t(0))}_{:=D_3}\\
        &+\underbrace{2y_t(0)\left(\frac{\indicator{Z_t=0}}{1-p_t}-1\right)\cdot\Delta_{t,1}(\bv_t(1))}_{:=D_4}\enspace.
    \end{align*}
    Following arguments similar to those in Lemma \mainref{lemma:pred-regret-unbiased}, we can show that $D_1$, $D_2$, $D_3$ and $D_4$ have mean zero by the law of iterated expectations.
    Hence by Lemma \mainref{lemma:standard-oco-for-pred-regret} and Lemma \mainref{lemma:ptpredloss-diff}, we have
    \begin{align}\label{prop:pred-regret_eq3}
        &\E{\mathcal{R}_T^{\text{pred}}}\notag\\
        \leq&\frac{m(\bv^*(1),\bv^*(0))}{\eta_{T+1}}+\e\left[\sum_{t=1}^T(\widetilde{L}_{t+1}(\bv_t(1),\bv_t(0))-\widetilde{L}_{t+1}(\widetilde{\bv}_{t+1}(1),\widetilde{\bv}_{t+1}(0)))\right]\notag\\
        \leq&\frac{m(\bv^*(1),\bv^*(0))}{\eta_{T+1}}+\e\left[\sum_{t=1}^{T}\Pi_{t,t}\cdot\widehat{\ell}_t(\bv_t(1),\bv_t(0))\right]~~~~\text{(Lemma \mainref{lemma:ptpredloss-diff})}\notag\\
        \leq&\underbrace{\frac{m(\bv^*(1),\bv^*(0))}{\eta_{T+1}}}_{:=S_1}+\underbrace{\e\left[\sum_{t=1}^{T}\Pi_{t,t}\cdot{\ell}_t(\bv_t(1),\bv_t(0))\right]}_{:=S_2}\notag\\
        &+\frac{\mathcal{E}(0)}{\mathcal{E}(1)}\underbrace{\e\left[\sum_{t=1}^{T}\Pi_{t,t}\cdot y_t(1)^2\left(\frac{\indicator{Z_t=1}}{p_t}-1\right)^2\right]}_{:=S_3}\notag\\
        &+\frac{\mathcal{E}(1)}{\mathcal{E}(0)}\underbrace{\e\left[\sum_{t=1}^{T}\Pi_{t,t}\cdot y_t(0)^2\left(\frac{\indicator{Z_t=0}}{1-p_t}-1\right)^2\right]}_{:=S_4}\notag\\
        &+2\underbrace{\e\left[\sum_{t=1}^{T}\Pi_{t,t}\cdot y_t(1)y_t(0)\left(\frac{\indicator{Z_t=1}}{p_t}-1\right)\left(\frac{\indicator{Z_t=0}}{1-p_t}-1\right)\right]}_{:=S_5}\enspace.
    \end{align}
    We now bound each term $S_1,\dots,S_5$ separately.
    By Assumption \mainref{assumption:moments} and the Cauchy-Schwarz inequality, we can show that $\{y_t(k):t\in[T]\}$ satisfy the condition in Lemma \ref{lemma:ols-norm} for $k\in\{0,1\}$.
    Hence, we have
    \begin{align}\label{prop:pred-regret_eq4}
        \quad S_1=&\eta_{T}^{-1}\left[\frac{\mathcal{E}(0)}{\mathcal{E}(1)}\cdot\|\bv^*(1)\|^2+\frac{\mathcal{E}(1)}{\mathcal{E}(0)}\cdot\|\bv^*(0)\|^2+2\cdot\iprod{\bv^*(1),\bv^*(0)}\right]\notag\\
        \quad\leq&\eta_{T}^{-1}\left[\frac{\mathcal{E}(0)}{\mathcal{E}(1)}\cdot\|\bv^*(1)\|^2+\frac{\mathcal{E}(1)}{\mathcal{E}(0)}\cdot\|\bv^*(0)\|^2+2\|\bv^*(1)\|\|\bv^*(0)\|\right]\quad(\text{Cauchy-Schwarz})\notag\\
        \quad\leq& c_1^2c_2\eta_T^{-1}\left(\frac{\mathcal{E}(0)}{\mathcal{E}(1)}+\frac{\mathcal{E}(1)}{\mathcal{E}(0)}+2\right)\quad(\text{Lemma \ref{lemma:ols-norm}})\notag\\
        \intertext{By Assumption \mainref{assumption:moments}, we can upper bound $\frac{\olsres{0}}{\olsres{1}}+\frac{\olsres{1}}{\olsres{0}}$ by $\frac{c_1}{c_0}+\frac{c_0}{c_1}$ given the shape of function $x+\frac{1}{x}$. Then we have}
        \quad\leq&c_1^2c_2T^{1/2}R_T\left(\frac{c_1}{c_0}+\frac{c_0}{c_1}+2\right)\notag\\
        \quad\leq&c_1^2c_2(\max\{c_2,1\})^{1/2}\left(\frac{c_1}{c_0}+\frac{c_0}{c_1}+2\right)T^{1/2}R\quad(\text{Lemma \ref{lemma:R}})\enspace.
    \end{align}
    We then derive the upper bound on $S_2$.
    Using the law of iterated expectations and Lemma \ref{lemma:squared-residuals-random}, for any $k\in\{0,1\}$ we have
    \begin{align*}
        \e\left[\sum_{t=1}^{T}\Delta_{t,k}(\bv_t(k))^2\right]=\e[\widehat{A}_T(k)]\leq \varsigma\,T\enspace.
    \end{align*}
    By the Cauchy-Schwarz inequality and the definition of ${\ell}_t(\bv_t(1),\bv_t(0))$, this implies that
    \begin{align}\label{prop:pred-regret_eq1}
        \e\left[\sum_{t=1}^{T}{\ell}_t(\bv_t(1),\bv_t(0))\right]=&\frac{\olsres{0}}{\olsres{1}}\cdot\e\left[\sum_{t=1}^{T}\Delta_{t,1}(\bv_t(1))^2\right]+\frac{\olsres{1}}{\olsres{0}}\cdot\e\left[\sum_{t=1}^{T}\Delta_{t,0}(\bv_t(0))^2\right]\notag\\
        &+2\e\left[\sum_{t=1}^{T}\Delta_{t,1}(\bv_t(1))\cdot \Delta_{t,0}(\bv_t(0))\right]\notag\\
        \leq&\frac{\olsres{0}}{\olsres{1}}\cdot\e\left[\sum_{t=1}^{T}\Delta_{t,1}(\bv_t(1))^2\right]+\frac{\olsres{1}}{\olsres{0}}\cdot\e\left[\sum_{t=1}^{T}\Delta_{t,0}(\bv_t(0))^2\right]\notag\\
        &+2\left(\e\left[\sum_{t=1}^{T}\Delta_{t,1}(\bv_t(1))^2\right]\right)^{1/2}\left(\e\left[\sum_{t=1}^{T}\Delta_{t,0}(\bv_t(0))^2\right]\right)^{1/2}\notag\\
        \leq&\left(\frac{\olsres{0}}{\olsres{1}}+\frac{\olsres{1}}{\olsres{0}}+2\right)\varsigma\,T\notag\\
        \leq&\left(\frac{c_1}{c_0}+\frac{c_0}{c_1}+2\right)\varsigma\,T\quad(\text{Assumption \mainref{assumption:moments}})\enspace.
    \end{align}
    Hence, by \eqref{prop:pred-regret_eq1}, we can establish the following upper bound on $S_2$:
    \begin{align}\label{prop:pred-regret_eq5}
        S_2\leq &\e\left[\sum_{t=1}^{T}R_t^2\eta_t{\ell}_t(\bv_t(1),\bv_t(0))\right]\notag\\
        \leq&\e\left[\sum_{t=1}^{T}R_tT^{-1/2}{\ell}_t(\bv_t(1),\bv_t(0))\right]\notag\\
        \intertext{Since $R_tT^{-1/2}\leq R_TT^{-1/2}\leq (\max\{c_2,1\})^{1/2}T^{-1/2}R$ by Lemma \ref{lemma:R}, we can further upper bound it by}
        \leq&(\max\{c_2,1\})^{1/2}T^{-1/2}R\cdot\e\left[\sum_{t=1}^{T}{\ell}_t(\bv_t(1),\bv_t(0))\right]\notag\\
        \leq&(\max\{c_2,1\})^{1/2}\left(\frac{c_1}{c_0}+\frac{c_0}{c_1}+2\right)\varsigma T^{1/2}R\quad(\text{by \eqref{prop:pred-regret_eq1}}) \enspace.
    \end{align}
    By the law of iterated expectations, we can simplify $S_3$ as
    \begin{align}\label{prop:pred-regret_eq6}
        S_3=&\e\left[\sum_{t=1}^{T}\Pi_{t,t}y_t(1)^2\left(\frac{\indicator{Z_t=1}}{p_t}-1\right)^2\right]\notag\\
        =&\sum_{t=1}^{T}\Pi_{t,t}y_t(1)^2\e\left[\frac{1}{p_t}-1\right]\notag\\
        \intertext{We upper bound $\Pi_{t,t}$ using Corollary \ref{corollary:pi} and upper bound $\e\left[\frac{1}{p_t}-1\right]$ using Lemma \ref{lemma:p-power-moment}. Then we have}
        \leq&\gamma\sum_{t=1}^{T}R_t^2((t-1)\vee \eta_t^{-1})^{-1}y_t(1)^2\cdot \chi\eta_t^{1/4}T^{1/4}\notag\\
        \leq&\gamma\chi T^{1/8}\sum_{t=1}^{T}R_t^{7/4}((t-1)\vee \eta_t^{-1})^{-1}y_t(1)^2\notag\\
        \leq&\gamma\chi T^{1/8}\left(\sum_{t=1}^{T}R_t^{7/2}((t-1)\vee \eta_t^{-1})^{-2}\right)^{1/2}\left(\sum_{t=1}^{T}y_t(1)^4\right)^{1/2}~\text{(Cauchy-Schwarz)}\notag\\
        \intertext{Since $R_t^{7/2}((t-1)\vee\eta_t^{-1})^{-2}=R_t^{3/2}((t-1)R_t^{-1}\vee T^{1/2})^{-2}\leq R_T^{3/2}((t-1)R_T^{-1}\vee T^{1/2})^{-2}=R_T^{7/2}((t-1)\vee\eta_T^{-1})^{-2}$, we have}
        \leq&c_1^2\gamma\chi T^{5/8}\left(\sum_{t=1}^{T}R_T^{7/2}((t-1)\vee \eta_T^{-1})^{-2}\right)^{1/2}\quad(\text{Assumption \mainref{assumption:moments}})\notag\\
        \intertext{By Lemma \ref{lemma:power-sum}, the sum in the bracket is upper bounded by $\xi_2R_T^{7/2}\eta_T=\xi_2T^{-1/2}R_T^{5/2}$. Hence we can further upper bound by}
        \leq&c_1^2\gamma\xi_2^{1/2}\chi T^{3/8}R_T^{5/4}\notag\\
        \leq& c_1^2\gamma(\max\{c_2,1\})^{5/8}\xi_2^{1/2}\chi T^{3/8}R^{5/4}\quad(\text{Lemma \ref{lemma:R}})\notag\\
        \leq&c_1^2c_3^{1/4}\gamma(\max\{c_2,1\})^{5/8}\xi_2^{1/2}\chi T^{3/8+1/8}R\quad(\text{Assumption \mainref{assumption:maximum-radius}})\notag\\
        =&c_1^2c_3^{1/4}\gamma(\max\{c_2,1\})^{5/8}\xi_2^{1/2}\chi T^{1/2}R\enspace.
    \end{align}
    Similarly, we can prove that
    \begin{align}\label{prop:pred-regret_eq7}
        S_4\leq& c_1^2c_3^{1/4}\gamma(\max\{c_2,1\})^{5/8}\xi_2^{1/2}\chi T^{1/2}R\enspace,\notag\\
        S_5\leq&S_3^{1/2}S_4^{1/2}\leq c_1^2c_3^{1/4}\gamma(\max\{c_2,1\})^{5/8}\xi_2^{1/2}\chi T^{1/2}R\quad(\text{Cauchy-Schwarz})\enspace.
    \end{align}
    Hence by \eqref{prop:pred-regret_eq3}, \eqref{prop:pred-regret_eq4}, \eqref{prop:pred-regret_eq5}, \eqref{prop:pred-regret_eq6} and \eqref{prop:pred-regret_eq7}, we have
    \begin{align*}
        &\E{\mathcal{R}_T^{\text{pred}}}\\
        \leq&c_1^2c_2(\max\{c_2,1\})^{1/2}\left(\frac{c_1}{c_0}+\frac{c_0}{c_1}+2\right)T^{1/2}R+(\max\{c_2,1\})^{1/2}\left(\frac{c_1}{c_0}+\frac{c_0}{c_1}+2\right)\varsigma T^{1/2}R\notag\\
        &+\left(\frac{\mathcal{E}(0)}{\mathcal{E}(1)}+\frac{\mathcal{E}(1)}{\mathcal{E}(0)}+2\right)c_1^2c_3^{1/4}\gamma(\max\{c_2,1\})^{5/8}\xi_2^{1/2}\chi T^{1/2}R\notag\\
        \leq&\left(\frac{c_1}{c_0}+\frac{c_0}{c_1}+2\right)\Big(c_1^2c_2(\max\{c_2,1\})^{1/2}+c_1^2c_3^{1/4}\gamma(\max\{c_2,1\})^{5/8}\xi_2^{1/2}\chi+(\max\{c_2,1\})^{1/2}\varsigma\Big)T^{1/2}R\enspace,
    \end{align*}
    which completes the proof.
\end{proof}

\subsection{Fourth Moment of Online Residuals}\label{section:C5}
The upper bound established in Lemma \ref{lemma:prob-regret-bound} involves the fourth moments of the online residuals,
\begin{align*}
    \e\left[\sum_{t=1}^T\eta_t(y_t(k)-\iprod{\xv_t,\bv_t(k)})^4\right]\enspace.
\end{align*}
In this section, we derive an upper bound for this quantity, which is nontrivial in the design-based framework.

The proof proceeds in two steps. We first bound the fourth moments of the full-information online residuals,
\begin{align*}
    \sum_{t=1}^T\eta_t(y_t(k)-\iprod{\xv_t,\bv_t^*(k)})^4\enspace,
\end{align*}
and then control the fourth moments of the tracking error terms,
\begin{align*}
    \e\left[\sum_{t=1}^T\eta_t\iprod{\xv_t,\bv_t(k)-\bv_t^*(k)}^4\right]\enspace.
\end{align*}

The following corollary follows directly from Lemma \ref{lemma:squared-residual-deterministic} and plays a key role in bounding the spectral norm of the matrix $\breve{\mat{Q}}^{(t)}$, which will be defined prior to Lemma \ref{lemma:Q-matrix-preliminary}.

\begin{corollary}\label{corollary:squared-inner-product}
    Under Assumptions \mainref{assumption:covariate-regularity}-\mainref{assumption:maximum-radius}, for any $\{v_t:t\in[T]\}$ such that $\sum_{t=1}^Tv_t^2\leq c_1^2T$, it holds that: $\sum_{t=1}^T\iprod{\xv_t,\bv^*_{v,t}}^2\leq(\zeta_1^{1/2}+c_1)^2T$.
    Here $\optbv_{v,t}:=\argmin_{\bv\in\mathbb{R}^d}\sum_{s<t}(v_s-\iprod{ \xv_s, \bv })^2+\eta_t^{-1}\|\bv\|^2$.
\end{corollary}

\begin{proof}
    By Lemma \ref{lemma:squared-residual-deterministic} and the Cauchy-Schwarz inequality, we have
    \begin{align*}
        \zeta_1 T\geq& \sum_{t=1}^T\left(v_t-\iprod{\xv_t,\bv^*_{v,t}}\right)^2~~~~(\text{Lemma \ref{lemma:squared-residual-deterministic}})\\
        =&\sum_{t=1}^Tv_t^2-2\sum_{t=1}^Tv_t\cdot\iprod{\xv_t,\bv^*_{v,t}}+\sum_{t=1}^T\iprod{\xv_t,\bv^*_{v,t}}^2\\
        \geq&\sum_{t=1}^Tv_t^2-2\left(\sum_{t=1}^Tv_t^2\right)^{1/2}\left(\sum_{t=1}^T\iprod{\xv_t,\bv^*_{v,t}}^2\right)^{1/2}+\sum_{t=1}^T\iprod{\xv_t,\bv^*_{v,t}}^2~~~~(\text{Cauchy-Schwarz})\\
        =&\left[\left(\sum_{t=1}^Tv_t^2\right)^{1/2}-\left(\sum_{t=1}^T\iprod{\xv_t,\bv^*_{v,t}}^2\right)^{1/2}\right]^2\enspace.
    \end{align*}
    By assumption, we have $0\leq \left(\sum_{t=1}^Tv_t^2\right)^{1/2}\leq c_1T^{1/2}$. 
    Hence we have $\left(\sum_{t=1}^T\langle \xv_t,\bv_{v,t}^*\rangle^2\right)^{1/2}\leq (\zeta_1^{1/2}+c_1)T^{1/2}$, which indicates that $\sum_{t=1}^T\langle \xv_t,\bv_{v,t}^*\rangle^2\leq (\zeta_1^{1/2}+c_1)^2T$.
\end{proof}

For any $2\leq t\leq T$, let the entries of matrix $\breve{\mat{Q}}^{(t)}=(\breve{Q}_{i,j}^{(t)})\in\mathbb{R}^{(t-1)\times (t-1)}$ be $\breve{Q}_{i,j}^{(t)}=\sum_{s=i\vee j+1}^{t}\Pi_{s,i}\Pi_{s,j}$ for any $1\leq i,j\leq t-1$.
Based on Corollary \ref{corollary:squared-inner-product}, we can establish the following lemma.

\begin{lemma}\label{lemma:Q-matrix-preliminary}
    Under Assumptions \mainref{assumption:covariate-regularity}-\mainref{assumption:maximum-radius}, $\breve{\mat{Q}}^{(t)}$ is a positive semidefinite matrix and it holds that: $\|\breve{\mat{Q}}^{(t)}\|\leq (\zeta_1^{1/2}c_1^{-1}+1)^2$ for any $2\leq t\leq T$.
\end{lemma}

\begin{proof}
    For any sequence $\{v_t:t\in[T]\}$ such that $\sum_{t=1}^Tv_t^2\leq c_1^2T$, define $\optbv_{v,t}:=\argmin_{\bv\in\mathbb{R}^d}\sum_{s<t}(v_s-\iprod{ \xv_s, \bv })^2+\eta_t^{-1}\|\bv\|^2$.
    We then rewrite $\sum_{s=1}^t\iprod{\xv_s,\bv_{v,s}^*}^2$ as a quadratic form of $v_1,\ldots,v_T$.
    By expanding the squared term and swapping the summation order, we have
    \begin{align*}
        \sum_{s=1}^t\iprod{\xv_s,\bv_{v,s}^*}^2=&\sum_{s=1}^t\left(\sum_{r=1}^{s-1}\Pi_{s,r}v_r\right)^2\\
        =&\sum_{s=1}^t\sum_{1\leq r_1,r_2\leq s-1}\Pi_{s,r_1}\Pi_{s,r_2}v_{r_1}v_{r_2}\\
        =&\sum_{1\leq r_1,r_2\leq t-1}\left(\sum_{s=r_1\vee r_2+1}^t\Pi_{s,r_1}\Pi_{s,r_2}\right)v_{r_1}v_{r_2}\\
        =&\sum_{r_1=1}^{t-1}\sum_{r_2=1}^{t-1}\breve{Q}_{r_1,r_2}^{(t)}v_{r_1}v_{r_2}\enspace.
    \end{align*}
    Hence, by Corollary \ref{corollary:squared-inner-product}, we have
    \begin{align*}
        0\leq\sum_{r_1=1}^{t-1}\sum_{r_2=1}^{t-1}\breve{Q}_{r_1,r_2}^{(t)}v_{r_1}v_{r_2}=\sum_{s=1}^t\iprod{\xv_s,\bv_{v,s}^*}^2\leq\sum_{s=1}^T\iprod{\xv_s,\bv_{v,s}^*}^2\leq (\zeta_1^{1/2}+c_1)^2T.
    \end{align*}
    Since this inequality holds for any sequence $\{v_t:t\in[T]\}$ such that $\sum_{t=1}^Tv_t^2\leq c_1^2T$, we can deduce that $\breve{\mat{Q}}^{(t)}$ is positive semidefinite and $\|\breve{\mat{Q}}^{(t)}\|\leq (\zeta_1^{1/2}c_1^{-1}+1)^2$.
\end{proof}

Based on Lemma \ref{lemma:Q-matrix-preliminary}, we can extend Corollary \ref{corollary:squared-inner-product} to all indices $t=1,\ldots,T$ on the potential outcomes $\{y_t(k):t\in[T],k\in\{0,1\}\}$, which is formally stated in the following corollary.

\begin{corollary}\label{corollary:online-squared-inner-product}
    Under Assumptions \mainref{assumption:moments}-\mainref{assumption:maximum-radius}, for any $t\in[T]$ and $k\in\{0,1\}$, it holds that: $\sum_{s=1}^t\langle \xv_s,\bv_s^*(k)\rangle^2\leq (\zeta_1^{1/2}c_1^{-1}+1)^2(\sum_{s=1}^{t-1}y_s(k)^2)$.
\end{corollary}

\begin{proof}
    Without loss of generality, we only prove the result for $k=1$.
    For $t=1$, the result is proved since $\bv_1^*(1)=\vec{0}$ by construction.
    For $t\geq 2$, by Lemma \ref{lemma:Q-matrix-preliminary}, we have
    \begin{align*}
        \sum_{s=1}^t\iprod{\xv_s,\bv_s^*(1)}^2=&\sum_{s=1}^t\left(\sum_{r=1}^{s-1}\Pi_{s,r}y_r(1)\right)^2\\
        =&\sum_{s=1}^t\sum_{1\leq r_1,r_2\leq s-1}\Pi_{s,r_1}\Pi_{s,r_2}y_{r_1}(1)y_{r_2}(1)\\
        =&\sum_{1\leq r_1,r_2\leq t-1}\left(\sum_{s=r_1\vee r_2+1}^t\Pi_{s,r_1}\Pi_{s,r_2}\right)y_{r_1}(1)y_{r_2}(1)\\
        =&\sum_{r_1=1}^{t-1}\sum_{r_2=1}^{t-1}\breve{Q}_{r_1,r_2}^{(t)}y_{r_1}(1)y_{r_2}(1)\\
        \leq&\|\breve{\mat{Q}}^{(t)}\|\left(\sum_{s=1}^{t-1}y_s(1)^2\right)\\
        \leq&(\zeta_1^{1/2}c_1^{-1}+1)^2\left(\sum_{s=1}^{t-1}y_s(1)^2\right)\quad(\text{Lemma \ref{lemma:Q-matrix-preliminary}})\enspace.
    \end{align*}
    Hence the result is proved.
\end{proof}

We now state the well-known Hardy's inequality, which will be used to bound the fourth moment of the online inner products in Lemma \ref{lemma:fourth-moment-deterministic}.

\begin{proposition}[Hardy's inequality]\label{prop:hardy}
    If $a_1,\ldots,a_T$ is a sequence of nonnegative real numbers, then for every real number $p>1$ we have
    \begin{align*}
        \sum_{t=1}^T\left(\frac{1}{t}\sum_{s=1}^ta_s\right)^p\leq \left(\frac{p}{p-1}\right)^p\sum_{t=1}^Ta_t^p\enspace.
    \end{align*}
\end{proposition}

Based on Corollary \ref{corollary:online-squared-inner-product} and Hardy's inequality, we establish the following lemma which upper-bounds the fourth moment of the full-information online inner products, which serves as a key step in establishing Corollary \ref{corollary:fourth-moment-deterministic}.

\begin{lemma}\label{lemma:fourth-moment-deterministic}
    Under Assumptions \mainref{assumption:moments}-\mainref{assumption:maximum-radius}, for $k\in\{0,1\}$, it holds that:
    \begin{align*}
        \sum_{t=1}^T\eta_t\langle \xv_t,\bv_t^*(k)\rangle^4\leq 5c_1^2\gamma(\max\{c_2,1\})^{1/2}(\zeta_1^{1/2}+c_1)^2T^{1/2}R\enspace.
    \end{align*}
\end{lemma}

\begin{proof}
    Without loss of generality, we only prove the result for $k=1$.
    For simplicity, for any $t\in[T]$ we denote $S_t=\sum_{s=1}^{t}y_s(1)^2$ and $B_t=\sum_{s=1}^t\langle \xv_s,\bv_s^*(1)\rangle^2$.
    For any $2\leq t\leq T$, by the Cauchy-Schwarz inequality and Lemma \ref{lemma:deterministic-summation-1}, we have
    \begin{align}\label{lemma:fourth-moment-deterministic_eq1}
        &\eta_t\langle \xv_t,\bv_t^*(1)\rangle^2\notag\\
        =&\eta_t\left(\sum_{s=1}^{t-1}\Pi_{t,s}y_s(1)\right)^2\notag\\
        \leq&\eta_t\left(\sum_{s=1}^{t-1}\Pi_{t,s}^2\right)\left(\sum_{s=1}^{t-1}y_s(1)^2\right)\quad(\text{Cauchy-Schwarz})\notag\\
        \leq&\gamma T^{-1/2}R_t^{-1}\cdot R_t^2((t-1)\vee \eta_t^{-1})^{-1}S_{t-1}\quad(\text{Lemma \ref{lemma:deterministic-summation-1}-1})\notag\\
        \leq&\gamma T^{-1/2}R_T(t-1)^{-1}S_{t-1}\quad(\text{since $R_t\leq R_T$ and $(t-1)\vee \eta_t^{-1}\geq t-1$})\enspace.
    \end{align}
    Since $S_{t-1}\leq S_{t}$, we have the following inequality:
    \begin{align*}
        \frac{S_{t-1}}{t-1}-\frac{S_{t}}{t}\leq\frac{S_{t-1}}{t-1}-\frac{S_{t-1}}{t}=\frac{S_{t-1}}{t(t-1)}\enspace.
    \end{align*}
    Hence by Corollary \ref{corollary:online-squared-inner-product}, we can establish the following inequality:
    \begin{align}\label{lemma:fourth-moment-deterministic_eq2}
        B_t\left(\frac{S_{t-1}}{t-1}-\frac{S_{t}}{t}\right)\leq&\left(\sum_{s=1}^t\langle \xv_s,\bv_s^*(1)\rangle^2\right)\cdot \frac{S_{t-1}}{t(t-1)}\notag\\
        \leq&(\zeta_1^{1/2}c_1^{-1}+1)^2\cdot \frac{S_{t-1}(\sum_{s=1}^{t-1}y_s(1)^2)}{t(t-1)}\quad(\text{Corollary \ref{corollary:online-squared-inner-product}})\notag\\
        \leq&(\zeta_1^{1/2}c_1^{-1}+1)^2\left(\frac{S_{t-1}}{t-1}\right)^2\enspace.
    \end{align}
    By \eqref{lemma:fourth-moment-deterministic_eq1}, we have
    \begin{align}\label{lemma:fourth-moment-deterministic_eq3}
        \qquad&\sum_{t=1}^T\eta_t\langle \xv_t,\bv_t^*(1)\rangle^4\notag\\
        \leq&\gamma T^{-1/2}R_T\sum_{t=2}^T\langle \xv_t,\bv_t^*(1)\rangle^2(t-1)^{-1}S_{t-1}\quad(\text{by \eqref{lemma:fourth-moment-deterministic_eq1}})\notag\\
        =&\gamma T^{-1/2}R_T\sum_{t=2}^T\left(B_t-B_{t-1}\right)\cdot\frac{S_{t-1}}{t-1}\quad(\text{by definition of $B_t$})\notag\\
        =&\gamma T^{-1/2}R_T\Bigg[\underbrace{\sum_{t=2}^{T}B_t\left(\frac{S_{t-1}}{t-1}-\frac{S_{t}}{t}\right)}_{:=D_1}+\underbrace{\frac{S_{T}}{T}\cdot B_T}_{:=D_2}\Bigg]\quad(\text{rewrite the telescoping sum})\enspace.
    \end{align}
    By \eqref{lemma:fourth-moment-deterministic_eq2} and Hardy's inequality (Proposition \ref{prop:hardy}) using $p=2$, we have
    \begin{align}\label{lemma:fourth-moment-deterministic_eq4}
        D_1\leq& (\zeta_1^{1/2}c_1^{-1}+1)^2\sum_{t=2}^T\left(\frac{1}{t-1}\sum_{s=1}^{t-1}y_s(1)^2\right)^2\notag\\
        \leq&\left(\frac{2}{2-1}\right)^2(\zeta_1^{1/2}c_1^{-1}+1)^2\sum_{t=1}^{T}y_t(1)^4\quad(\text{Hardy's inequality})\notag\\
        \leq&4c_1^4(\zeta_1^{1/2}c_1^{-1}+1)^2T\quad(\text{Assumption \mainref{assumption:moments}})\enspace.
    \end{align}
    By Corollary \ref{corollary:squared-inner-product} and Assumption \mainref{assumption:moments}, we have
    \begin{align}\label{lemma:fourth-moment-deterministic_eq5}
        D_2\leq \frac{1}{T}\cdot (\zeta_1^{1/2}+c_1)^2T\cdot\sum_{t=1}^Ty_t(1)^2\leq c_1^2(\zeta_1^{1/2}+c_1)^2T\enspace.
    \end{align}
    By \eqref{lemma:fourth-moment-deterministic_eq3}, \eqref{lemma:fourth-moment-deterministic_eq4}, and \eqref{lemma:fourth-moment-deterministic_eq5}, we have
    \begin{align*}
        \sum_{t=1}^T\eta_t\langle \xv_t,\bv_t^*(1)\rangle^4\leq&\gamma T^{-1/2}R_T(D_1+D_2)\\
        \leq&\gamma T^{-1/2}R_T\cdot 5c_1^2(\zeta_1^{1/2}+c_1)^2T\\
        \leq&5c_1^2\gamma(\max\{c_2,1\})^{1/2}(\zeta_1^{1/2}+c_1)^2T^{1/2}R\quad(\text{Lemma \ref{lemma:R}})\enspace,
    \end{align*}
    which completes the proof.
\end{proof}

Using the proof in Lemma \ref{lemma:fourth-moment-deterministic}, we can derive the following corollary regarding the partial sum of the fourth powers of full-information online inner products.

\begin{corollary}\label{corollary:partial-fourth-moment}
    Under Assumptions \mainref{assumption:moments}-\mainref{assumption:maximum-radius}, for $k\in\{0,1\}$, any fixed $0\leq \delta\leq 2$, and any $t\in[T]$, it holds that:
    \begin{align*}
        \sum_{s=1}^tR_s^{-\delta}\langle \xv_s,\bv_s^*(k)\rangle^4\leq 5c_1^2\gamma(\zeta_1^{1/2}+c_1)^2TR_t^{2-\delta}\enspace.
    \end{align*}
\end{corollary}

\begin{proof}
    By similar arguments as in Lemma \ref{lemma:fourth-moment-deterministic}, the result is proved.
\end{proof}

Based on Lemma \ref{lemma:fourth-moment-deterministic}, we can finally establish the upper bound on the sum of fourth powers of full-information online residuals in the following corollary:

\begin{corollary}\label{corollary:fourth-moment-deterministic}
    Under Assumptions \mainref{assumption:moments}-\mainref{assumption:maximum-radius}, for $k\in\{0,1\}$, it holds that:
    \begin{align*}
        \sum_{t=1}^T\eta_t\left(y_t(k)-\iprod{\xv_t,\bv^*_t(k)}\right)^4\leq c_1^2(\max\{c_2,1\})^{1/2}(1+5^{1/4}\gamma^{1/4})^{4}(\zeta_1^{1/2}+c_1)^{2}T^{1/2}R \enspace.
    \end{align*}
\end{corollary}

\begin{proof}
    Without loss of generality, we only prove the result for $k=1$.
	By the triangle inequality in $\ell^4$ (Minkowski inequality), we have
    \begin{align}\label{corollary:fourth-moment-deterministic_eq1}
        \sum_{t=1}^T\eta_t\left(y_t(1)-\iprod{\xv_t,\bv^*_t(1)}\right)^4=&\sum_{t=1}^T|\eta_t^{1/4}y_t(1)-\eta_t^{1/4}\iprod{\xv_t,\bv^*_t(1)}|^4\notag\\
        \leq&\left[\left(\sum_{t=1}^T\eta_ty_t(1)^4\right)^{1/4}+\left(\sum_{t=1}^T\eta_t\iprod{\xv_t,\bv^*_t(1)}^4\right)^{1/4}\right]^4\enspace.
    \end{align}
    By Assumption \mainref{assumption:moments}, we have
    \begin{align}\label{corollary:fourth-moment-deterministic_eq2}
        \sum_{t=1}^T\eta_ty_t(1)^4 \leq T^{-1/2}\cdot c_1^4T=c_1^4 T^{1/2} \enspace. 
    \end{align}
    Since Lemma \ref{lemma:R} implies that $R(\max\{c_2,1\})^{1/2}\geq R_T\geq d^{1/2}\geq 1$, by Lemma \ref{lemma:fourth-moment-deterministic}, \eqref{corollary:fourth-moment-deterministic_eq1} and \eqref{corollary:fourth-moment-deterministic_eq2}, we have
    \begin{align*}
        &\sum_{t=1}^T\eta_t\left(y_t(1)-\iprod{\xv_t,\bv^*_t(1)}\right)^4\\
        \leq&\left(c_1T^{1/8}+5^{1/4}c_1^{1/2}\gamma^{1/4}(\max\{c_2,1\})^{1/8}(\zeta_1^{1/2}+c_1)^{1/2}T^{1/8}R^{1/4}\right)^4\quad(\text{Lemma \ref{lemma:fourth-moment-deterministic}})\\
        \leq&\left(c_1T^{1/8}R^{1/4}(\max\{c_2,1\})^{1/8}+5^{1/4}c_1^{1/2}\gamma^{1/4}(\max\{c_2,1\})^{1/8}(\zeta_1^{1/2}+c_1)^{1/2}T^{1/8}R^{1/4}\right)^4\\
        \leq&c_1^2(\max\{c_2,1\})^{1/2}(1+5^{1/4}\gamma^{1/4})^{4}(\zeta_1^{1/2}+c_1)^{2}T^{1/2}R\enspace.
    \end{align*}
    Therefore, the result is proved.
\end{proof}

We next derive an upper bound on the fourth moments of the tracking error terms: $\e\left[\sum_{t=1}^T\eta_t\iprod{\xv_t,\bv_t(k)-\bv_t^*(k)}^4\right]$.
As in the analysis of the estimated squared residuals, the argument requires separate control of the corresponding deterministic summations and the moments of the inverse assignment probabilities.

The following lemma provides upper bounds for moments of the inverse probability weighting terms that arise in bounding the tracking error terms.

\begin{lemma}\label{lemma:p-moments-cross-term}
    Under Assumptions \mainref{assumption:moments}-\mainref{assumption:maximum-radius} and Condition \mainref{condition:sigmoid}, it holds that:
    \begin{enumerate}
        \item[(1)] For any $s\in[T]$, $\e\left[\left(\frac{\indicator{Z_{s}=1}}{p_{s}}-1\right)^4\right]\leq 2\chi^{3}R_s^{-3/4}T^{3/8}$.
        \item[(2)] For any $1\leq s_2<s_1\leq T$, $\e\left[\left(\frac{\indicator{Z_{s_1}=1}}{p_{s_1}}-1\right)^2\left(\frac{\indicator{Z_{s_2}=1}}{p_{s_2}}-1\right)^2\right]\leq 4\chi^{9/4}R_{s_1}^{-1/4}R_{s_2}^{-5/16}T^{9/32}$.
        \item[(3)] For any $s\in[T]$, $\e\left[\left(\frac{\indicator{Z_{s}=0}}{1-p_{s}}-1\right)^4\right]\leq 2\chi^{3}R_s^{-3/4}T^{3/8}$.
        \item[(4)] For any $1\leq s_2<s_1\leq T$, $\e\left[\left(\frac{\indicator{Z_{s_1}=0}}{1-p_{s_1}}-1\right)^2\left(\frac{\indicator{Z_{s_2}=0}}{1-p_{s_2}}-1\right)^2\right]\leq 4\chi^{9/4}R_{s_1}^{-1/4}R_{s_2}^{-5/16}T^{9/32}$.
    \end{enumerate}
\end{lemma}

\begin{proof}
    Without loss of generality, we only prove (1)-(2).
    Results in (3)-(4) follow from the symmetry between the treated group and control group.
    Throughout the proof, we repeatedly use the fact that for any $s\in[T]$ and $k\geq 1$,
    \begin{align*}
        \e\left[\left(\frac{\indicator{Z_s=1}}{p_s}\right)^k\right]=\e\left[\frac{\indicator{Z_s=1}}{p_s^k}\right]=\e\left[\e\left[\frac{\indicator{Z_s=1}}{p_s^k}\Big|\filt_{s-1}\right]\right]=\e\left[\frac{1}{p_s^{k-1}}\right]\enspace.
    \end{align*}
    \textbf{Part 1:} For any $s\in [T]$, by Lemma \ref{lemma:p-power-moment}, we have
        \begin{align*}
            \e\left[\left(\frac{\indicator{Z_{s}=1}}{p_{s}}-1\right)^4\right]\leq&\e\left[\left(\frac{\indicator{Z_{s}=1}}{p_{s}}\right)^4\right]+1\quad(\text{since $\indicator{Z_{s}=1}/p_{s}$ is either 0 or larger than 1})\\
            =&\e\left[\frac{1}{p_{s}^3}\right]+1\\
            \leq&2\e\left[\frac{1}{p_{s}^3}\right]\\
            \leq&2\chi^3R_s^{-3/4}T^{3/8}\enspace.
        \end{align*}
        \textbf{Part 2:} For any $1\leq s_2< s_1 \leq T$, by the law of iterated expectations, H{\"o}lder's inequality and Lemma \ref{lemma:p-power-moment}, we have
        \begin{align*}
            &\e\left[\left(\frac{\indicator{Z_{s_1}=1}}{p_{s_1}}-1\right)^2\left(\frac{\indicator{Z_{s_2}=1}}{p_{s_2}}-1\right)^2\right]\\
            \leq&\e\left[\left(\frac{\indicator{Z_{s_1}=1}}{p_{s_1}^2}+1\right)\left(\frac{\indicator{Z_{s_2}=1}}{p_{s_2}^2}+1\right)\right]\\
            =&\e\left[\frac{\indicator{Z_{s_1}=1}}{p_{s_1}^2}\frac{\indicator{Z_{s_2}=1}}{p_{s_2}^2}\right]+\e\left[\frac{\indicator{Z_{s_1}=1}}{p_{s_1}^2}\right]+\e\left[\frac{\indicator{Z_{s_2}=1}}{p_{s_2}^2}\right]+1\\
            =&\e\left[\frac{1}{p_{s_1}}\frac{\indicator{Z_{s_2}=1}}{p_{s_2}^2}\right]+\e\left[\frac{1}{p_{s_1}}\right]+\e\left[\frac{1}{p_{s_2}}\right]+1\\
            \leq&\left(\e\left[\frac{1}{p_{s_1}^4}\right]\right)^{1/4}\left(\e\left[\frac{\indicator{Z_{s_2}=1}}{p_{s_2}^{2\cdot \frac{4}{3}}}\right]\right)^{3/4}+\e\left[\frac{1}{p_{s_1}}\right]+\e\left[\frac{1}{p_{s_2}}\right]+1\quad(\text{H{\"o}lder's inequality})\\
            \leq&\left(\e\left[\frac{1}{p_{s_1}^4}\right]\right)^{1/4}\left(\e\left[\frac{1}{p_{s_2}^{5/3}}\right]\right)^{3/4}+\chi R_{s_1}^{-1/4}T^{1/8}+\chi R_{s_2}^{-1/4}T^{1/8}+1\quad(\text{Lemma \ref{lemma:p-power-moment}})\\
            \leq&\chi R_{s_1}^{-1/4}T^{1/8}\cdot \chi^{5/4} R_{s_2}^{-5/16}T^{5/32}+\chi R_{s_1}^{-1/4}T^{1/8}+\chi R_{s_2}^{-1/4}T^{1/8}+1\quad(\text{Lemma \ref{lemma:p-power-moment}})\\
            \leq &4\chi^{9/4}R_{s_1}^{-1/4}R_{s_2}^{-5/16}T^{9/32}\enspace.
        \end{align*}
    Therefore, the lemma is proved.
\end{proof}

We next establish an upper bound on the fourth moment of the prediction tracking error terms, which quantifies the gap between the full-information fourth moments and those of the (random) online residuals.
Using Lemma \ref{lemma:deterministic-summation-1} and Lemma \ref{lemma:p-moments-cross-term}, we derive the desired upper bound in the following proposition.

\begin{refproposition}{\mainref{prop:expected-tracking}*}
    Under Assumptions \mainref{assumption:moments}-\mainref{assumption:maximum-radius} and Condition \mainref{condition:sigmoid}, the prediction tracking error for both treatments $k \in \setb{0,1}$ can be bounded as
    \begin{align*}
        \e\left[\sum_{t=1}^T\eta_t\iprod{\xv_t,\bv_t(k)-\bv_t^*(k)}^4\right]\leq \widetilde{\varUpsilon}\cdot T^{1/2}R\enspace,
    \end{align*}
    where $\widetilde{\varUpsilon}:=32cc_1^4\chi^3\gamma^2(\max\{c_2,1\})^{1/2}>0$ is a constant.
\end{refproposition}

\begin{proof}
    We only provide the proof for $k=1$ since the argument for $k=0$ is identical.
    By Lemma \mainref{lemma:ridge-pred-exact-form}*, we have the following equality:
    \begin{align}\label{prop:expected-tracking_eq1}
        \iprod{\xv_t,\bv_t(k)-\bv_t^*(k)}=\sum_{s=1}^{t-1}\underbrace{\Pi_{t,s}y_s(1)\left(\frac{\indicator{Z_s=1}}{p_s}-1\right)}_{:=M_{t,s}}\enspace.
    \end{align}
    For any fixed $t\in[T]$, $\{M_{t,s}:s=1,\ldots,t-1\}$ form a martingale difference sequence with respect to $\filt_{s}$.
    Hence, for any $s=1,\ldots,t-1$, we have
    \begin{align*}
        &\e\left[\left(\sum_{r=1}^{s}M_{t,r}\right)^4\right]\\
        =&\e\left[\left(\sum_{r=1}^{s-1}M_{t,r}+M_{t,s}\right)^4\right]\\
        =&\e\left[\left(\sum_{r=1}^{s-1}M_{t,r}\right)^4\right]+4\e\left[\left(\sum_{r=1}^{s-1}M_{t,r}\right)^3M_{t,s}\right]+6\e\left[\left(\sum_{r=1}^{s-1}M_{t,r}\right)^2M_{t,s}^2\right]\\
        &+4\e\left[\left(\sum_{r=1}^{s-1}M_{t,r}\right)M_{t,s}^3\right]+\e\left[M_{t,s}^4\right]\\
        \intertext{By the law of iterated expectations, we have $\e\left[\left(\sum_{r=1}^{s-1}M_{t,r}\right)^3M_{t,s}\right]=\e\left[\left(\sum_{r=1}^{s-1}M_{t,r}\right)^3\e\left[M_{t,s}|\filt_{s-1}\right]\right]=0$. Using the inequality $4ab^3\leq 2a^2b^2+2b^4$, we have}
        \leq&\e\left[\left(\sum_{r=1}^{s-1}M_{t,r}\right)^4\right]+8\e\left[\left(\sum_{r=1}^{s-1}M_{t,r}\right)^2M_{t,s}^2\right]+3\e\left[M_{t,s}^4\right]\\
        \leq&\e\left[\left(\sum_{r=1}^{s-1}M_{t,r}\right)^4\right]+8\e\left[\left(\sum_{r=1}^{s-1}M_{t,r}\right)^4\right]^{1/2}\e\left[M_{t,s}^4\right]^{1/2}+16\e\left[M_{t,s}^4\right]\quad(\text{Cauchy-Schwarz})\\
        =&\left(\e\left[\left(\sum_{r=1}^{s-1}M_{t,r}\right)^4\right]^{1/2}+4\e\left[M_{t,s}^4\right]^{1/2}\right)^2\enspace.
    \end{align*}
    Taking the square root on both sides leads to the following induction equation:
    \begin{align*}
        \e\left[\left(\sum_{r=1}^{s}M_{t,r}\right)^4\right]^{1/2}\leq \e\left[\left(\sum_{r=1}^{s-1}M_{t,r}\right)^4\right]^{1/2}+4\e\left[M_{t,s}^4\right]^{1/2}\enspace,
    \end{align*}
    which together with \eqref{prop:expected-tracking_eq1} implies that
    \begin{align}\label{prop:expected-tracking_eq2}
        \e\left[\iprod{\xv_t,\bv_t(k)-\bv_t^*(k)}^4\right]^{1/2}=\e\left[\left(\sum_{s=1}^{t-1}M_{t,s}\right)^4\right]^{1/2}\leq 4\sum_{s=1}^{t-1}\e\left[M_{t,s}^4\right]^{1/2}\enspace.
    \end{align}
    On the other hand, by the definition of $M_{t,s}$, it holds that
    \begin{align}\label{prop:expected-tracking_eq3}
        \e\left[M_{t,s}^4\right]=&\e\left[\left(\Pi_{t,s}y_s(1)\left(\frac{\indicator{Z_s=1}}{p_s}-1\right)\right)^4\right]\notag\\
        =&\Pi_{t,s}^4y_s(1)^4\e\left[\left(\frac{\indicator{Z_s=1}}{p_s}-1\right)^4\right]\notag\\
        \leq&2\chi^3R_s^{-3/4}T^{3/8}\Pi_{t,s}^4y_s(1)^4\quad(\text{Lemma \ref{lemma:p-moments-cross-term}--1})\enspace.
    \end{align}
    Therefore, by \eqref{prop:expected-tracking_eq2} and \eqref{prop:expected-tracking_eq3}, we have shown that
    \begin{align*}
        &\e\left[\iprod{\xv_t,\bv_t(k)-\bv_t^*(k)}^4\right]\\
        \leq& \left(4\sum_{s=1}^{t-1}\e\left[M_{t,s}^4\right]^{1/2}\right)^2\\
        \leq&32\chi^3T^{3/8}\left(\sum_{s=1}^{t-1}R_s^{-3/8}\Pi_{t,s}^2y_s(1)^2\right)^2\\
        \leq&32\chi^3T^{3/8}\left(\sum_{s=1}^{t-1}\Pi_{t,s}^2\right)\left(\sum_{s=1}^{t-1}R_s^{-3/4}\Pi_{t,s}^2y_s(1)^4\right)\quad(\text{Cauchy-Schwarz})\\
        \leq&32\chi^3\gamma T^{3/8}R_t^2((t-1)\vee \eta_t^{-1})^{-1}\left(\sum_{s=1}^{t-1}R_s^{-3/4}\Pi_{t,s}^2y_s(1)^4\right)\quad(\text{Lemma \ref{lemma:deterministic-summation-1}--1})\enspace.
    \end{align*}
    Finally, by Lemma \ref{lemma:deterministic-summation-1}, we have
    \begin{align*}
        &\e\left[\sum_{t=1}^T\eta_t\iprod{\xv_t,\bv_t(k)-\bv_t^*(k)}^4\right]\\
        \leq&32\chi^3\gamma T^{3/8}\sum_{t=1}^T\eta_tR_t^2((t-1)\vee \eta_t^{-1})^{-1}\sum_{s=1}^{t-1}R_s^{-3/4}\Pi_{t,s}^2y_s(1)^4\\
        \intertext{Since $\eta_tR_t^2((t-1)\vee \eta_t^{-1})^{-1}\leq \eta_t^2R_t^2=T^{-1}$, we have}
        \leq&32\chi^3\gamma T^{3/8-1}\sum_{t=1}^T \sum_{s=1}^{t-1}R_s^{-3/4}\Pi_{t,s}^2y_s(1)^4\\
        =&32\chi^3\gamma T^{-5/8}\sum_{s=1}^{T-1}R_s^{-3/4}y_s(1)^4\sum_{t=s+1}^T\Pi_{t,s}^2\\
        \leq&32c\chi^3\gamma^2 T^{-5/8}\sum_{s=1}^{T-1}R_s^{-3/4+2}((s-1)\vee \eta_s^{-1})^{-1}y_s(1)^4\quad(\text{Lemma \ref{lemma:deterministic-summation-1}--2})\\
        \leq&32c\chi^3\gamma^2 T^{-5/8-1/2}\sum_{s=1}^{T-1}R_s^{-3/4+2-1}y_s(1)^4\quad(\text{$((s-1)\vee \eta_s^{-1})^{-1}\leq \eta_s=R_s^{-1}T^{-1/2}$})\\
        \leq&32cc_1^4\chi^3\gamma^2 T^{-5/8-1/2+1}R_T^{1/4}\\
        \leq&32cc_1^4\chi^3\gamma^2 T^{-5/8-1/2+1}R_T\quad(\text{since $R_T\geq d^{1/2}\geq 1$})\\
        \leq&32cc_1^4\chi^3\gamma^2(\max\{c_2,1\})^{1/2}T^{1/2}R\quad(\text{Lemma \ref{lemma:R}})\enspace.
    \end{align*}
    This completes the proof.
\end{proof}

Combining the results in Corollary \ref{corollary:fourth-moment-deterministic} and Proposition \mainref{prop:expected-tracking}*, we obtain the upper bound on the fourth moment of the (random) online residuals in the following lemma.
\begin{reflemma}{\mainref{lemma:bounded-fourth-moments}*}
    Under Assumptions \mainref{assumption:moments}-\mainref{assumption:maximum-radius} and Condition \mainref{condition:sigmoid}, for each treatment $k\in\{0,1\}$, the fourth moment of the online residuals is bounded by:
    \begin{align*}
        &\e\left[\sum_{t=1}^T\eta_t\left(y_t(k)-\iprod{\xv_t,\bv_t(k)}\right)^4\right]\leq \varUpsilon \cdot T^{1/2}R,
    \end{align*}
    where $\varUpsilon>0$ is defined as:
    \begin{align*}
        \varUpsilon:=&(c_1^{1/2}(\max\{c_2,1\})^{1/8}(1+5^{1/4}\gamma^{1/4})(\zeta_1^{1/2}+c_1)^{1/2}+\widetilde{\varUpsilon}^{1/4})^4\enspace.
    \end{align*}
\end{reflemma}

\begin{proof}
    Without loss of generality, we only prove the result for $k=1$.
    By the Minkowski inequality in $L^4$, Corollary \ref{corollary:fourth-moment-deterministic} and Proposition \mainref{prop:expected-tracking}*, we have
    \begin{align*}
        &\e\left[\sum_{t=1}^T\eta_t\left(y_t(1)-\iprod{\xv_t,\bv_t(1)}\right)^4\right]\\
        \leq&\left[\left(\e\left[\sum_{t=1}^T\eta_t\iprod{\xv_t,\bv_t(1)-\bv_t^*(1)}^4\right]\right)^{1/4}+\left(\sum_{t=1}^T\eta_t\left(y_t(1)-\iprod{\xv_t,\bv_t^*(1)}\right)^4\right)^{1/4}\right]^4\quad(\text{Minkowski})\\
        \leq&(c_1^{1/2}(\max\{c_2,1\})^{1/8}(1+5^{1/4}\gamma^{1/4})(\zeta_1^{1/2}+c_1)^{1/2}+\widetilde{\varUpsilon}^{1/4})^4T^{1/2}R\\
        =&\varUpsilon\cdot T^{1/2}R\enspace.
    \end{align*}
    Hence the result is proved.
\end{proof}

\subsection{Neyman Regret Analysis (Theorem \mainref{thm:neyman-regret})}\label{section:C6}

In this section, we derive an upper bound on the expected Neyman regret.
By combining the bounds in Lemma \ref{lemma:prob-regret-bound} and Lemma \mainref{lemma:bounded-fourth-moments}*, we can establish the final upper bound on the expected probability regret in the following proposition.

\begin{refproposition}{\mainref{prop:prob-regret-bound}*}
    Under Assumptions \mainref{assumption:moments}-\mainref{assumption:maximum-radius} and Condition \mainref{condition:sigmoid}, the expected probability regret is bounded as
    \begin{align*}
        \E{\regretprob_T}\leq \left[\frac{8(\max\{c_2,1\})^{1/2}}{\bb_3^3}\left[\left(\frac{c_1}{c_0}-1\right)^3+\frac{\bb_3}{4}\left(\frac{c_1}{c_0}-1\right)^2\right]+\bb_1^2\max\{\bb_1,2\}\varUpsilon\right]T^{1/2}R\enspace.
    \end{align*}
\end{refproposition}

\begin{proof}
    By Lemma \ref{lemma:prob-regret-bound}, Lemma \mainref{lemma:bounded-fourth-moments}* and Lemma \ref{lemma:R}, we have
    \begin{align*}
        \E{\regretprob_T}\leq&\frac{8T^{1/2}R_T}{\bb_3^3}\left[\left(\frac{c_1}{c_0}-1\right)^3+\frac{\bb_3}{4}\left(\frac{c_1}{c_0}-1\right)^2\right]\\
        &+ \frac{\bb_1^2}{2}\max\{\bb_1,2\}\sum_{k \in \setb{0,1}} \E[\Bigg]{ \sum_{t=1}^T \eta_t \braces[\big]{ y_t(k) - \iprod{ \xv_t, \bv_t(k) } }^4 }\quad(\text{Lemma \ref{lemma:prob-regret-bound}})\\
        \leq&\frac{8(\max\{c_2,1\})^{1/2}T^{1/2}R}{\bb_3^3}\left[\left(\frac{c_1}{c_0}-1\right)^3+\frac{\bb_3}{4}\left(\frac{c_1}{c_0}-1\right)^2\right]\\
        &+\bb_1^2\max\{\bb_1,2\}\varUpsilon\cdot T^{1/2}R\quad(\text{Lemma \mainref{lemma:bounded-fourth-moments}* and Lemma \ref{lemma:R}})\\
        =&\left[\frac{8(\max\{c_2,1\})^{1/2}}{\bb_3^3}\left[\left(\frac{c_1}{c_0}-1\right)^3+\frac{\bb_3}{4}\left(\frac{c_1}{c_0}-1\right)^2\right]+\bb_1^2\max\{\bb_1,2\}\varUpsilon\right]T^{1/2}R\enspace.
    \end{align*}
    This completes the proof.
\end{proof}

Finally, by combining the bounds on the probability regret (Proposition \mainref{prop:prob-regret-bound}*) and the prediction regret (Proposition \mainref{prop:pred-regret}*) and the decomposition lemma (Lemma \mainref{lemma:regret-decomposition}), we obtain the Neyman regret upper bound stated in Theorem \mainref{thm:neyman-regret}.

\begin{reftheorem}{\mainref{thm:neyman-regret}*}
	Under Assumptions \mainref{assumption:moments}-\mainref{assumption:maximum-radius} and Condition \mainref{condition:sigmoid}, the Neyman regret is bounded as
	$\neymanregret_T \leq C_N\cdot T^{-1/2} R$, with constant
	\begin{align*}
        &C_N = \frac{8(\max\{c_2,1\})^{1/2}}{\bb_3^3}\left[\left(\frac{c_1}{c_0}-1\right)^3+\frac{\bb_3}{4}\left(\frac{c_1}{c_0}-1\right)^2\right]+\bb_1^2\max\{\bb_1,2\}\varUpsilon\\
        &+\left(\frac{c_1}{c_0}+\frac{c_0}{c_1}+2\right)\Big(c_1^2c_2(\max\{c_2,1\})^{1/2}+c_1^2c_3^{1/4}\gamma(\max\{c_2,1\})^{5/8}\xi_2^{1/2}\chi+(\max\{c_2,1\})^{1/2}\varsigma\Big)\enspace.
    \end{align*}
    Here, $\gamma$ is defined in Lemma \ref{lemma:regularity}; $\xi_2$ in Lemma \ref{lemma:power-sum}; $\zeta_1$ in Lemma \ref{lemma:squared-residual-deterministic}; $\varsigma$ in Lemma \ref{lemma:squared-residuals-random}; $\chi$ in Lemma \ref{lemma:p-power-moment}; and $\varUpsilon$ in Lemma \mainref{lemma:bounded-fourth-moments}*.
\end{reftheorem}

\begin{proof}
    By Lemma \mainref{lemma:regret-decomposition}, we have that the Neyman regret is decomposed into the expected probability regret and the expected prediction regret:
    \[
    \neymanregret_T = \frac{1}{T}\E{\regretprob_T}+\frac{1}{T}\e[\mathcal{R}_T^{\text{pred}}]
    \enspace.
    \]
    The proof is completed by applying Proposition \mainref{prop:prob-regret-bound}* and Proposition \mainref{prop:pred-regret}* which bound each of these expected regrets individually.
\end{proof}
	
	\section{Inference Analysis}
In this section, we aim to verify the inference procedure based on the AIPW estimator and the proposed variance bound estimator.
Section \ref{section:D1} presents the technical lemmas used in the subsequent proofs.
Section \ref{section:D2} establishes an almost sure upper bound for the inverse probabilities, which play a crucial role in verifying the central limit theorem.
Section \ref{section:D3} provides a detailed proof of the martingale central limit theorem.
Section \ref{section:D4} demonstrates the equivalence between Assumption \mainref{assumption:bounded-correlation} and the non-superefficiency condition in design-based causal inference.
Section \ref{section:D5} verifies the consistency of the proposed variance bound estimator.
Using the asymptotic normality established in Section \ref{section:D3} and the consistent variance estimator proposed in Section \ref{section:D5}, we construct Wald-type confidence intervals in Section \ref{section:D6}.

\subsection{Technical Lemmas}\label{section:D1}
In this section, we present several technical lemmas that will be used to establish the main results in the subsequent sections.

For any fixed $1\leq s< t\leq T$, let the entries of matrix $\mat{Q}^{(s,t)}=(Q_{i,j}^{(s,t)})\in\mathbb{R}^{(t-s)\times (t-s)}$ be
    \begin{align*}
        Q^{(s,t)}_{i,j}=\left(\sum_{r=i\vee j+s}^{t}\Pi_{r,i+s-1}\Pi_{r,j+s-1}\right)^2\cdot\left(\frac{R_{t}}{R_{i+s-1}}\cdot \frac{R_{t}}{R_{j+s-1}}\right)^{1/2}
    \end{align*}
for any $1\leq i,j\leq t-s$. 
The following corollary is derived from Lemma \ref{lemma:Q-matrix-preliminary}, which bounds the spectral norm of the matrix $\mat{Q}^{(s,t)}$.

\begin{corollary}\label{corollary:Q-matrix}
    Under Assumptions \mainref{assumption:covariate-regularity}-\mainref{assumption:maximum-radius}, $\mat{Q}^{(s,t)}$ is a positive semidefinite matrix, and the following holds
    \begin{align*}
        \|\mat{Q}^{(s,t)}\|\leq c\gamma(\zeta_1^{1/2}c_1^{-1}+1)^2R_t^2((s-1)\vee \eta_t^{-1})^{-1}
    \end{align*}
    for any $1\leq s< t\leq T$.
\end{corollary}

\begin{proof}
    For fixed $1\leq s< t\leq T$, let the entries in matrix $\widetilde{\mat{Q}}^{(s,t)}=(\widetilde{Q}_{i,j}^{(s,t)})\in\mathbb{R}^{(t-s)\times (t-s)}$ be $\widetilde{Q}^{(s,t)}_{i,j}=(\sum_{r=i\vee j+s}^{t}\Pi_{r,i+s-1}\Pi_{r,j+s-1})$ for any $1\leq i,j\leq t-s$ and let the entries of matrix $\widebar{\mat{Q}}^{(s,t)}=(\widebar{Q}_{i,j}^{(s,t)})\in\mathbb{R}^{(t-s)\times (t-s)}$ be
    \begin{align*}
        \widebar{Q}^{(s,t)}_{i,j}=\left(\sum_{r=i\vee j+s}^{t}\Pi_{r,i+s-1}\Pi_{r,j+s-1}\right)\cdot\left(\frac{R_{t}}{R_{i+s-1}}\cdot \frac{R_{t}}{R_{j+s-1}}\right)^{1/2}
    \end{align*}
    for any $1\leq i,j\leq t-s$.
    The two claims are then proved in the following steps.\\[3mm]
    \textbf{Step 1:} Prove that $\mat{Q}^{(s,t)}$ is positive semidefinite.\\[3mm]
    Since $\widetilde{\mat{Q}}^{(s,t)}$ is a principal submatrix of the positive semidefinite matrix $\breve{\mat{Q}}^{(t)}$ (defined before Lemma \ref{lemma:Q-matrix-preliminary}), $\widetilde{\mat{Q}}^{(s,t)}$ is also positive semidefinite. 
    By the definition of $\widebar{\mat{Q}}^{(s,t)}$, we can rewrite it as
    \begin{align*}
        \widebar{\mat{Q}}^{(s,t)}=\operatorname{diag}\left\{\left(\frac{R_{t}}{R_{s}}\right)^{1/2},\ldots,\left(\frac{R_{t}}{R_{t-1}}\right)^{1/2}\right\}\widetilde{\mat{Q}}^{(s,t)}\operatorname{diag}\left\{\left(\frac{R_{t}}{R_{s}}\right)^{1/2},\ldots,\left(\frac{R_{t}}{R_{t-1}}\right)^{1/2}\right\}\enspace,
    \end{align*}
    which implies that $\widebar{\mat{Q}}^{(s,t)}$ is also positive semidefinite. Since the Hadamard product of positive semidefinite matrices is positive semidefinite, this implies that $\mat{Q}^{(s,t)}=\widetilde{\mat{Q}}^{(s,t)}\circ \widebar{\mat{Q}}^{(s,t)}$ is positive semidefinite.\\[3mm]
    \textbf{Step 2:} Bound the spectral norm of  $\mat{Q}^{(s,t)}$.\\[3mm]
    To control the spectral norm, we aim to bound the diagonal entries of $\widebar{\mat{Q}}^{(s,t)}$.
    For any $1\leq k\leq t-s$, the $k$-th diagonal element of $\widebar{\mat{Q}}^{(s,t)}$ is bounded by:
    \begin{align}\label{corollary:Q-matrix_eq1}
        \widebar{Q}^{(s,t)}_{k,k}=&\frac{R_{t}}{R_{k+s-1}}\cdot \sum_{r=k+s}^{t}\Pi_{r,k+s-1}^2\notag\\
        \leq&c\gamma R_tR_{k+s-1}^{-1}\cdot R_{k+s-1}^2((k+s-2)\vee\eta_{k+s-1}^{-1})^{-1}\quad(\text{Lemma \ref{lemma:deterministic-summation-1}-2})\notag\\
        =&c\gamma R_t((k+s-2) R_{k+s-1}^{-1}\vee T^{1/2})^{-1}\notag\\
        \leq&c\gamma R_t((k+s-2)R_t^{-1}\vee T^{1/2})^{-1}\quad(\text{since $R_{k+s-1}\leq R_t$ for $k\leq t-s$})\notag\\
        \leq&c\gamma R_t^2((s-1)\vee\eta_t^{-1})^{-1}\enspace.
    \end{align}
    Since $\widetilde{\mat{Q}}^{(s,t)}$ is a principal submatrix of $\breve{\mat{Q}}^{(t)}$, we have $\|\widetilde{\mat{Q}}^{(s,t)}\|\leq \|\breve{\mat{Q}}^{(t)}\|\leq (\zeta_1^{1/2}c_1^{-1}+1)^2$ by Lemma \ref{lemma:Q-matrix-preliminary}. Then by \eqref{corollary:Q-matrix_eq1} and Theorem 5.3.4 in \cite{horn2012matrix}, we obtain
    \begin{align*}
        \|{\mat{Q}}^{(s,t)}\|=&\|\widetilde{\mat{Q}}^{(s,t)}\circ\widebar{\mat{Q}}^{(s,t)}\|\\
        \leq& \|\widetilde{\mat{Q}}^{(s,t)}\|\cdot\max_{1\leq k\leq t-s}\widebar{Q}^{(s,t)}_{k,k}\\
        \leq& c\gamma(\zeta_1^{1/2}c_1^{-1}+1)^2R_t^2((s-1)\vee \eta_t^{-1})^{-1}\enspace.
    \end{align*}
    Therefore, the result is proved.
\end{proof}

The following lemma provides the results for the deterministic summations that will be used in the proof of the central limit theorem.

\begin{lemma}\label{lemma:deterministic-summation-new}
    Under Assumptions \mainref{assumption:moments}-\mainref{assumption:maximum-radius-strengthened}, for any $t\in[T]$ and $k\in\{0,1\}$, we have
    \begin{enumerate}
        \item[(1)] $\sum_{r=1}^{t-1}\left(\sum_{s=r+1}^{t}\Pi_{s,r}^2\right)^2y_r(k)^4\leq c^2c_1^4\gamma^2R_t^2$.
        \item[(2)] $\sum_{s=1}^{t-1}\left(\sum_{r=s+1}^t\Pi_{r,s}\left(y_r(k)-\iprod{\xv_r,\bv_r^*(k)}\right)\right)^2R_s^{-1/4}y_s(k)^2\leq c\,c_1^2\gamma\xi_2^{1/2}\zeta_1T^{5/4}R_t^{5/4}$.
        \item[(3)] $\sum_{1\leq t_2,t_3<t_1\leq t-1}\left|\sum_{s=t_1+1}^{t}\Pi_{s,t_1}\Pi_{s,t_2}\right|\left|\sum_{s=t_1+1}^{t}\Pi_{s,t_1}\Pi_{s,t_3}\right|y_{t_1}(k)^2\left|y_{t_2}(k)\right|\left|y_{t_3}(k)\right|\left(\frac{R_t^2}{R_{t_1}^2}\cdot\frac{R_t}{R_{t_2}}\cdot\frac{R_t}{R_{t_3}}\right)^{1/4}\\
        \leq 2^{-1/2}c\,c_1^2\gamma(\zeta_1^{1/2}+c_1)^2T^{5/4}R_t^{3/2}\log^2(\eta_tT)$ for $T$ sufficiently large.
    \end{enumerate}
\end{lemma}

\begin{proof}
    Without loss of generality, we only prove the result for $k=1$.
    \begin{enumerate}
        \item[(1)] By Lemma \ref{lemma:deterministic-summation-1}, we have
        \begin{align*}
            \sum_{r=1}^{t-1}\left(\sum_{s=r+1}^{t}\Pi_{s,r}^2\right)^2y_r(1)^4\leq&c^2\gamma^2\sum_{r=1}^{t-1}R_r^4((r-1)\vee \eta_r^{-1})^{-2}y_r(1)^4\quad(\text{Lemma \ref{lemma:deterministic-summation-1}-2})\\
            \leq&c^2\gamma^2\sum_{r=1}^{t-1}R_r^4\eta_r^{2}y_r(1)^4\\
            \leq&c^2\gamma^2R_t^2T^{-1}\sum_{r=1}^{t-1}y_r(1)^4\quad(\text{since $R_r\leq R_t$})\\
            \leq&c^2c_1^4\gamma^2R_t^2\quad(\text{Assumption \mainref{assumption:moments}})\enspace.
        \end{align*}
        \item[(2)] By Lemma \ref{lemma:power-sum}, Lemma \ref{lemma:deterministic-summation-1} and Corollary \ref{corollary:squared-residual-deterministic}, we have
        \begin{align*}
             &\sum_{s=1}^{t-1}\left(\sum_{r=s+1}^t\Pi_{r,s}\left(y_r(1)-\langle \xv_r,\bv_r^*(1)\rangle\right)\right)^2R_s^{-1/4}y_s(1)^2\\
             \leq&\sum_{s=1}^{t-1}\left(\sum_{r=s+1}^t\Pi_{r,s}^2\right)\left(\sum_{r=s+1}^t\left(y_r(1)-\langle \xv_r,\bv_r^*(1)\rangle\right)^2\right)R_s^{-1/4}y_s(1)^2~~~~(\text{Cauchy-Schwarz})\\
             \leq&c\gamma\zeta_1T\sum_{s=1}^{t-1}R_s^{2}((s-1)\vee \eta_s^{-1})^{-1}R_s^{-1/4}y_s(1)^2~~~~(\text{Lemma \ref{lemma:deterministic-summation-1}-2, Corollary \ref{corollary:squared-residual-deterministic}})\\
             \leq&c\gamma\zeta_1T\left(\sum_{s=1}^{t-1}R_s^{7/2}((s-1)\vee \eta_s^{-1})^{-2}\right)^{1/2}\left(\sum_{s=1}^{t-1}y_s(1)^4\right)^{1/2}~~~~(\text{Cauchy-Schwarz})\\
             \leq&c\,c_1^2\gamma\zeta_1T^{3/2}\left(\sum_{s=1}^{t-1}R_s^{3/2}((s-1)R_s^{-1}\vee  T^{1/2})^{-2}\right)^{1/2}\quad(\text{Assumption \mainref{assumption:moments}})\\
             \leq&c\,c_1^2\gamma\zeta_1T^{3/2}\left(\sum_{s=1}^{t-1}R_t^{3/2}((s-1)R_t^{-1}\vee  T^{1/2})^{-2}\right)^{1/2}\quad(\text{since $R_s\leq R_t$})\\
             =&c\,c_1^2\gamma\zeta_1T^{3/2}\left(\sum_{s=1}^{t-1}R_t^{7/2}((s-1)\vee  \eta_t^{-1})^{-2}\right)^{1/2}\\
             \leq&c\,c_1^2\gamma\xi_2^{1/2}\zeta_1T^{3/2}\eta_t^{1/2}R_t^{7/4}\quad(\text{Lemma \ref{lemma:power-sum}})\\
             =&c\,c_1^2\gamma\xi_2^{1/2}\zeta_1T^{5/4}R_t^{5/4}\enspace.
         \end{align*}
         \item[(3)] By the definition of matrix $\mat{Q}^{(1,t)}$ before Corollary \ref{corollary:Q-matrix}, we have
         \begin{align*}
             &\sum_{1\leq t_2,t_3<t_1\leq t-1}\left|\sum_{s=t_1+1}^{t}\Pi_{s,t_1}\Pi_{s,t_2}\right|\left|\sum_{s=t_1+1}^{t}\Pi_{s,t_1}\Pi_{s,t_3}\right|y_{t_1}(1)^2\left|y_{t_2}(1)\right|\left|y_{t_3}(1)\right|\left(\frac{R_t^2}{R_{t_1}^2}\cdot\frac{R_t}{R_{t_2}}\cdot\frac{R_t}{R_{t_3}}\right)^{1/4}\\
             \leq&\frac{1}{2}\sum_{1\leq t_2,t_3<t_1\leq t-1}\Bigg[\left(\sum_{s=t_1+1}^{t}\Pi_{s,t_1}\Pi_{s,t_2}\right)^2y_{t_1}(1)^2y_{t_2}(1)^2\left(\frac{R_t^2}{R_{t_1}^2}\cdot\frac{R_t^2}{R_{t_2}^2}\right)^{1/4}\\
             &+\left(\sum_{s=t_1+1}^{t}\Pi_{s,t_1}\Pi_{s,t_3}\right)^2y_{t_1}(1)^2y_{t_3}(1)^2\left(\frac{R_t^2}{R_{t_1}^2}\cdot\frac{R_t^2}{R_{t_3}^2}\right)^{1/4}\Bigg]~~~~(\text{AM-GM inequality})\\
             =&\sum_{1\leq t_2<t_1\leq t-1}(t_1-1)\left(\sum_{s=t_1+1}^{t}\Pi_{s,t_1}\Pi_{s,t_2}\right)^2y_{t_1}(1)^2y_{t_2}(1)^2\left(\frac{R_t}{R_{t_1}}\cdot\frac{R_t}{R_{t_2}}\right)^{1/2}\quad(\text{by symmetry})\\
             \leq&\frac{1}{2}\sum_{1\leq t_1, t_2\leq t-1}(t_1\vee t_2)\left(\sum_{s=t_1\vee t_2+1}^{t}\Pi_{s,t_1}\Pi_{s,t_2}\right)^2y_{t_1}(1)^2y_{t_2}(1)^2\left(\frac{R_t}{R_{t_1}}\cdot\frac{R_t}{R_{t_2}}\right)^{1/2}\\
             =&\frac{1}{2}\sum_{1\leq t_1, t_2\leq t-1}(t_1\vee t_2)Q^{(1,t)}_{t_1,t_2}y_{t_1}(1)^2y_{t_2}(1)^2\enspace.
        \end{align*}
        We then partition the indices in this summation into $(L+1)\times(L+1)$ blocks and upper bound the summation in each block.
        Let $L=\lfloor \log^{2} (\eta_tT)\rfloor-1$. Since $\eta_t^{-1}=\bigO{T^{1/2}R_t}=o(T)$ by Assumption \mainref{assumption:maximum-radius-strengthened}, $L$ tends to infinity. For $\ell=1,\ldots,L+1$, let $C_{\ell}=\lceil\eta_t^{-1+(\ell-1)/L}T^{(\ell-1)/L}\rceil$ and let $C_{0}=0$. For $\ell=0,\ldots,L$, denote $B_\ell=\{C_\ell+1,\ldots,C_{\ell+1}\}\cap \{1,\ldots,t-1\}$. Then we have
        \begin{align*}
            &\sum_{1\leq t_2,t_3<t_1\leq t-1}\left|\sum_{s=t_1+1}^{t}\Pi_{s,t_1}\Pi_{s,t_2}\right|\left|\sum_{s=t_1+1}^{t}\Pi_{s,t_1}\Pi_{s,t_3}\right|y_{t_1}(1)^2\left|y_{t_2}(1)\right|\left|y_{t_3}(1)\right|\left(\frac{R_t^2}{R_{t_1}^2}\cdot\frac{R_t}{R_{t_2}}\cdot\frac{R_t}{R_{t_3}}\right)^{1/4}\\
            \leq&\frac{1}{2}\sum_{1\leq t_1, t_2\leq t-1}(t_1\vee t_2)Q^{(1,t)}_{t_1,t_2}y_{t_1}(1)^2y_{t_2}(1)^2\\
            =&\frac{1}{2}\sum_{\ell_1=0}^{L}\sum_{\ell_2=0}^L\underbrace{\left(\sum_{t_1\in B_{\ell_1}}\sum_{t_2\in B_{\ell_2}}(t_1\vee t_2)Q^{(1,t)}_{t_1,t_2}y_{t_1}(1)^2y_{t_2}(1)^2\right)}_{:=S_{\ell_1,\ell_2}}\enspace.
        \end{align*}
        We then bound $S_{\ell_1,\ell_2}$ for every $\ell_1,\ell_2=0,\ldots,L$.
        For $S_{0,0}$, by Corollary \ref{corollary:Q-matrix}, we have
        \begin{align*}
            S_{0,0}=&\sum_{t_1\in B_{0}}\sum_{t_2\in B_{0}}(t_1\vee t_2)Q^{(1,t)}_{t_1,t_2}y_{t_1}(1)^2y_{t_2}(1)^2\\
            \leq&C_1\sum_{t_1\in B_{0}}\sum_{t_2\in B_{0}}Q^{(1,t)}_{t_1,t_2}y_{t_1}(1)^2y_{t_2}(1)^2\\
            \leq&C_1\|\mat{Q}^{(1,t)}\|\sum_{s\in B_{0}}y_s(1)^4\\
            \leq&c\,\gamma(\zeta_1^{1/2}c_1^{-1}+1)^2R_t^2C_1\eta_t\sum_{s\in B_{0}}y_s(1)^4~~~~(\text{Corollary \ref{corollary:Q-matrix}})\enspace.
        \end{align*}
        For $\ell=1,\ldots,L$, by Corollary \ref{corollary:Q-matrix} and the Cauchy-Schwarz inequality for the bilinear form induced by a PSD matrix, we have
        \begin{align*}
            &S_{\ell,0}\\
            =&\sum_{t_1\in B_{\ell}}\sum_{t_2\in B_{0}}(t_1\vee t_2)Q^{(1,t)}_{t_1,t_2}y_{t_1}(1)^2y_{t_2}(1)^2\\
            \leq&C_{\ell+1}\sum_{t_1\in B_{\ell}}\sum_{t_2\in B_{0}}Q^{(1,t)}_{t_1,t_2}y_{t_1}(1)^2y_{t_2}(1)^2\\
            \intertext{If $C_{\ell}+1\geq t$, then the summation is zero. We only consider the case where $C_{\ell}+1< t$. Hence}
            \leq&C_{\ell+1}\left(\sum_{t_1,t_2\in B_{\ell}}Q^{(1,t)}_{t_1,t_2}y_{t_1}(1)^2y_{t_2}(1)^2\right)^{1/2}\left(\sum_{t_1,t_2\in B_{0}}Q^{(1,t)}_{t_1,t_2}y_{t_1}(1)^2y_{t_2}(1)^2\right)^{1/2}~~(\text{Cauchy-Schwarz})\\
            \leq&C_{\ell+1}\left(\left\|\mat{Q}^{(C_\ell+1,t)}\right\|\sum_{s\in B_{\ell}}y_s(1)^4\right)^{1/2}\left(\left\|\mat{Q}^{(1,t)}\right\|\sum_{s\in B_{0}}y_s(1)^4\right)^{1/2}\\
            \leq&C_{\ell+1}\left(c\,\gamma(\zeta_1^{1/2}c_1^{-1}+1)^2R_{t}^2C_{\ell}^{-1}\sum_{s\in B_{\ell}}y_s(1)^4\right)^{1/2}\left(c\,\gamma(\zeta_1^{1/2}c_1^{-1}+1)^2R_t^2\eta_t\sum_{s\in B_{0}}y_s(1)^4\right)^{1/2}~~~~(\text{Corollary \ref{corollary:Q-matrix}})\\
            =&c\,\gamma(\zeta_1^{1/2}c_1^{-1}+1)^2R_t^2\left(\frac{C_{\ell+1}}{C_\ell}\frac{C_{\ell+1}}{\eta_t^{-1}}\right)^{1/2}\left(\sum_{s\in B_{\ell}}y_s(1)^4\right)^{1/2}\left(\sum_{s\in B_0}y_s(1)^4\right)^{1/2}\enspace.
        \end{align*}
        For $\ell_1,\ell_2=1,\ldots,L$, by Corollary \ref{corollary:Q-matrix} and the Cauchy-Schwarz inequality for the bilinear form induced by a PSD matrix, we have
        \begin{align*}
            &S_{\ell_1,\ell_2}\\
            =&\sum_{t_1\in B_{\ell_1}}\sum_{t_2\in B_{\ell_2}}(t_1\vee t_2)Q_{t_1,t_2}^{(1,t)}y_{t_1}(1)^2y_{t_2}(1)^2\\
            \leq&C_{\ell_1\vee \ell_2+1}\left(\sum_{t_1,t_2\in B_{\ell_1}}Q^{(1,t)}_{t_1,t_2}y_{t_1}(1)^2y_{t_2}(1)^2\right)^{1/2}\left(\sum_{t_1,t_2\in B_{\ell_2}}Q^{(1,t)}_{t_1,t_2}y_{t_1}(1)^2y_{t_2}(1)^2\right)^{1/2}~~~~(\text{Cauchy-Schwarz})\\
            \intertext{If $C_{\ell_1}+1\geq t$ or $C_{\ell_2}+1\geq t$, then the summation is zero. We only consider the case where $\max\{C_{\ell_1},C_{\ell_2}\}+1< t$. Hence}
            \leq&C_{\ell_1\vee \ell_2+1}\left(\left\|\mat{Q}^{(C_{\ell_1}+1,t)}\right\|\sum_{s\in B_{\ell_1}}y_s(1)^4\right)^{1/2}\left(\left\|\mat{Q}^{(C_{\ell_2}+1,t)}\right\|\sum_{s\in B_{\ell_2}}y_s(1)^4\right)^{1/2}\\
            \leq&C_{\ell_1\vee \ell_2+1}\left(c\,\gamma(\zeta_1^{1/2}c_1^{-1}+1)^2R_t^2C_{\ell_1}^{-1}\sum_{s\in B_{\ell_1}}y_s(1)^4\right)^{1/2}\left(c\,\gamma(\zeta_1^{1/2}c_1^{-1}+1)^2R_t^2C_{\ell_2}^{-1}\sum_{s\in B_{\ell_2}}y_s(1)^4\right)^{1/2} 
            \intertext{where the last inequality uses Corollary~\ref{corollary:Q-matrix}.
            Continuing, we have that} 
            \leq&c\,\gamma(\zeta_1^{1/2}c_1^{-1}+1)^2R_t^2\left(\frac{C_{\ell_1\vee \ell_2+1}}{C_{\ell_1}}\frac{C_{\ell_1\vee \ell_2+1}}{C_{\ell_2}}\right)^{1/2}\left(\sum_{s\in B_{\ell_1}}y_s(1)^4\right)^{1/2}\left(\sum_{s\in B_{\ell_2}}y_s(1)^4\right)^{1/2}\enspace.
        \end{align*}
        By Assumption \mainref{assumption:maximum-radius-strengthened}, it can be shown that
        \begin{align*}
            \max_{\ell=1,\ldots,L}\frac{C_{\ell+1}}{C_{\ell}}=&\max_{\ell=1,\ldots,L}\frac{\lceil\eta_t^{-1+((\ell+1)-1)/L}T^{((\ell+1)-1)/L}\rceil}{\lceil\eta_t^{-1+(\ell-1)/L}T^{(\ell-1)/L}\rceil}\\
            \leq&\max_{\ell=1,\ldots,L}\frac{\eta_t^{-1+((\ell+1)-1)/L}T^{((\ell+1)-1)/L}+1}{\eta_t^{-1+(\ell-1)/L}T^{(\ell-1)/L}}\\
            \leq&\exp\left(\log[(T\eta_t)^{\frac{1}{L}}]\right)+\eta_t\\
            =&\exp\left(\frac{1}{L}\log(T\eta_t)\right)+\eta_t\\
            = &\exp\left(\frac{1}{\lfloor \log^{2} (T\eta_t)\rfloor-1}\log(T\eta_t)\right)+\eta_t\\
            \rightarrow&1\quad(\text{since $T\eta_t\rightarrow\infty$ and $\eta_t\rightarrow0$ by Assumption \mainref{assumption:maximum-radius-strengthened}})\enspace.
        \end{align*}
        Hence for $T$ sufficiently large, we have
        \begin{align*}
            \max_{\ell=1,\ldots,L}\frac{C_{\ell+1}}{C_\ell}\leq 2\enspace.
        \end{align*}
        By the construction of $C_1,\ldots,C_{L+1}$, we have
        \begin{align*}
            &\max\left\{C_1\eta_t,~\max_{\ell=1,\ldots,L}\left(\frac{C_{\ell+1}}{C_\ell}\frac{C_{\ell+1}}{\eta_t^{-1}}\right)^{1/2},~\max_{\ell_1,\ell_2=1,\ldots,L}\left(\frac{C_{\ell_1\vee \ell_2+1}}{C_{\ell_1}}\frac{C_{\ell_1\vee \ell_2+1}}{C_{\ell_2}}\right)^{1/2}\right\}\\
            \leq&\left(2\cdot \frac{T}{\eta_t^{-1}}\right)^{1/2}=2^{1/2}(T\eta_t)^{1/2}.
        \end{align*}
        This uniform bound implies that for any $\ell_1,\ell_2=0,\ldots,L$, it holds that
        \begin{align*}
            S_{\ell_1,\ell_2}\leq 2^{1/2}c\,\gamma(\zeta_1^{1/2}c_1^{-1}+1)^2R_t^2(T\eta_t)^{1/2}\left(\sum_{s\in B_{\ell_1}}y_s(1)^4\right)^{1/2}\left(\sum_{s\in B_{\ell_2}}y_s(1)^4\right)^{1/2}\enspace.
        \end{align*}
        Hence, for $T$ sufficiently large we have
        \begin{align*}
            &\sum_{1\leq t_2,t_3<t_1\leq t-1}\left|\sum_{s=t_1+1}^{t}\Pi_{s,t_1}\Pi_{s,t_2}\right|\left|\sum_{s=t_1+1}^{t}\Pi_{s,t_1}\Pi_{s,t_3}\right|y_{t_1}(1)^2\left|y_{t_2}(1)\right|\left|y_{t_3}(1)\right|\left(\frac{R_t^2}{R_{t_1}^2}\cdot\frac{R_t}{R_{t_2}}\cdot\frac{R_t}{R_{t_3}}\right)^{1/4}\\
            \leq&\frac{1}{2}\sum_{\ell_1=0}^{L}\sum_{\ell_2=0}^{L}S_{\ell_1,\ell_2}\\
            \leq&2^{-1/2}c\,\gamma(\zeta_1^{1/2}c_1^{-1}+1)^2R_t^2(T\eta_t)^{1/2}\left[\sum_{\ell_1=0}^{L}\sum_{\ell_2=0}^{L}\left(\sum_{s\in B_{\ell_1}}y_s(1)^4\right)^{1/2}\left(\sum_{s\in B_{\ell_2}}y_s(1)^4\right)^{1/2}\right]\\
            \leq&2^{-1/2}c\,\gamma(\zeta_1^{1/2}c_1^{-1}+1)^2R_t^2(T\eta_t)^{1/2}\left(\sum_{\ell=0}^{L}\left(\sum_{s\in B_{\ell}}y_s(1)^4\right)^{1/2}\right)^2\\
            \leq&2^{-1/2}c\,\gamma(\zeta_1^{1/2}c_1^{-1}+1)^2R_t^2(T\eta_t)^{1/2}(L+1)\left(\sum_{\ell=0}^{L}\left(\sum_{s\in B_{\ell}}y_s(1)^4\right)\right)~~~~(\text{Cauchy-Schwarz})\\
            \leq&2^{-1/2}c\,\gamma(\zeta_1^{1/2}c_1^{-1}+1)^2R_t^2(T\eta_t)^{1/2}(L+1)\sum_{s=1}^{T}y_s(1)^4\\
            \leq&2^{-1/2}c\,\gamma(\zeta_1^{1/2}c_1^{-1}+1)^2R_t^2(T\eta_t)^{1/2}(L+1)\cdot c_1^4T\quad(\text{Assumption \mainref{assumption:moments}})\\
            \leq&2^{-1/2}c\,c_1^2\gamma(\zeta_1^{1/2}+c_1)^2R_t^2\eta_t^{1/2}T^{3/2}\log^2(\eta_tT)\\
            =&2^{-1/2}c\,c_1^2\gamma(\zeta_1^{1/2}+c_1)^2T^{5/4}R_t^{3/2}\log^2(\eta_tT)\enspace.
        \end{align*}
    \end{enumerate}
    Therefore, the lemma is proved.
\end{proof}

The following lemma provides upper bounds for the moments of the inverse probability weighting terms that arise in the proof of the central limit theorem.

\begin{lemma}\label{lemma:p-moments-cross-term-new}
    Under Assumptions \mainref{assumption:moments}-\mainref{assumption:maximum-radius} and Condition \mainref{condition:sigmoid}, the following holds:
    \begin{enumerate}
        \item[(1)] For $s_2,s_3<s_1$, $s_2\neq s_3$, $\left|\operatorname{Cov}\left(\left(\frac{\indicator{Z_{s_1}=1}}{p_{s_1}}-1\right)\left(\frac{\indicator{Z_{s_2}=1}}{p_{s_2}}-1\right),\left(\frac{\indicator{Z_{s_1}=1}}{p_{s_1}}-1\right)\left(\frac{\indicator{Z_{s_3}=1}}{p_{s_3}}-1\right)\right)\right|\\
        \leq4\left(2+\bb_1(\bb_2/6)^{1/4}\eta_{s_1}^{1/4}\e[\widehat{A}_{s_1}(0)]^{1/4}\right)^{25/16}\left(\frac{R_{s_1}}{R_{s_2}}\right)^{5/64}\left(\frac{R_{s_1}}{R_{s_3}}\right)^{5/64}$.
        \item[(2)] For $s_2,s_3<s_1$, $s_2\neq s_3$, $\left|\operatorname{Cov}\left(\left(\frac{\indicator{Z_{s_1}=0}}{1-p_{s_1}}-1\right)\left(\frac{\indicator{Z_{s_2}=0}}{1-p_{s_2}}-1\right),\left(\frac{\indicator{Z_{s_1}=0}}{1-p_{s_1}}-1\right)\left(\frac{\indicator{Z_{s_3}=0}}{1-p_{s_3}}-1\right)\right)\right|\\
        \leq4\left(2+\bb_1(\bb_2/6)^{1/4}\eta_{s_1}^{1/4}\e[\widehat{A}_{s_1}(1)]^{1/4}\right)^{25/16}\left(\frac{R_{s_1}}{R_{s_2}}\right)^{5/64}\left(\frac{R_{s_1}}{R_{s_3}}\right)^{5/64}$.
    \end{enumerate}
\end{lemma}

\begin{proof}
    We only prove part 1. Part 2 follows from the symmetry between the treated group and the control group.
    We assume without loss of generality that $s_3<s_2$.
    For any $1\leq s_3<s_2<s_1\leq T$, by the law of total covariance we have
        \begin{align}\label{lemma:p-moments-cross-term-new_eq1}
            &\left|\operatorname{Cov}\left(\left(\frac{\indicator{Z_{s_1}=1}}{p_{s_1}}-1\right)\left(\frac{\indicator{Z_{s_2}=1}}{p_{s_2}}-1\right),\left(\frac{\indicator{Z_{s_1}=1}}{p_{s_1}}-1\right)\left(\frac{\indicator{Z_{s_3}=1}}{p_{s_3}}-1\right)\right)\right|\notag\\
            =&\Bigg|\e\Bigg[\operatorname{Cov}\Bigg(\left(\frac{\indicator{Z_{s_1}=1}}{p_{s_1}}-1\right)\left(\frac{\indicator{Z_{s_2}=1}}{p_{s_2}}-1\right),\left(\frac{\indicator{Z_{s_1}=1}}{p_{s_1}}-1\right)\left(\frac{\indicator{Z_{s_3}=1}}{p_{s_3}}-1\right)\Bigg|\mathcal{F}_{s_1-1}\Bigg)\Bigg]\notag\\
            +&\operatorname{Cov}\Bigg[\e\left(\left(\frac{\indicator{Z_{s_1}=1}}{p_{s_1}}-1\right)\left(\frac{\indicator{Z_{s_2}=1}}{p_{s_2}}-1\right)\Bigg|\mathcal{F}_{s_1-1}\right),\notag\\
            &\e\left(\left(\frac{\indicator{Z_{s_1}=1}}{p_{s_1}}-1\right)\left(\frac{\indicator{Z_{s_3}=1}}{p_{s_3}}-1\right)\Bigg|\mathcal{F}_{s_1-1}\right)\Bigg]\Bigg|\enspace.
        \end{align}
        By direct calculation, we obtain
        \begin{align*}
            &\operatorname{Cov}\Bigg(\left(\frac{\indicator{Z_{s_1}=1}}{p_{s_1}}-1\right)\left(\frac{\indicator{Z_{s_2}=1}}{p_{s_2}}-1\right),\left(\frac{\indicator{Z_{s_1}=1}}{p_{s_1}}-1\right)\left(\frac{\indicator{Z_{s_3}=1}}{p_{s_3}}-1\right)\Bigg|\mathcal{F}_{s_1-1}\Bigg)\\
            =&\left(\frac{\indicator{Z_{s_2}=1}}{p_{s_2}}-1\right)\left(\frac{\indicator{Z_{s_3}=1}}{p_{s_3}}-1\right)\operatorname{Var}\left(\frac{\indicator{Z_{s_1}=1}}{p_{s_1}}-1\Bigg|\mathcal{F}_{s_1-1}\right)\\
            =&\left(\frac{\indicator{Z_{s_2}=1}}{p_{s_2}}-1\right)\left(\frac{\indicator{Z_{s_3}=1}}{p_{s_3}}-1\right)\cdot \frac{1-p_{s_1}}{p_{s_1}}
        \end{align*}
        and
        \begin{align*}
            \e\left[\left(\frac{\indicator{Z_{s_1}=1}}{p_{s_1}}-1\right)\left(\frac{\indicator{Z_{s_2}=1}}{p_{s_2}}-1\right)\Bigg|\mathcal{F}_{s_1-1}\right]=\left(\frac{\indicator{Z_{s_2}=1}}{p_{s_2}}-1\right)\e\left[\frac{\indicator{Z_{s_1}=1}}{p_{s_1}}-1\Bigg|\mathcal{F}_{s_1-1}\right]=0\enspace.
        \end{align*}
        Note that $\e\left[\left(\frac{\indicator{Z_{s_2}=1}}{p_{s_2}}-1\right)\left(\frac{\indicator{Z_{s_3}=1}}{p_{s_3}}-1\right)\right]$ equals zero for $s_2\neq s_3$.
        Hence \eqref{lemma:p-moments-cross-term-new_eq1} can be further bounded as follows:
        \begin{align*}    
            &\left|\e\left[\frac{1-p_{s_1}}{p_{s_1}}\left(\frac{\indicator{Z_{s_2}=1}}{p_{s_2}}-1\right)\left(\frac{\indicator{Z_{s_3}=1}}{p_{s_3}}-1\right)\right]\right|\\
            =&\left|\e\left[\frac{1}{p_{s_1}}\left(\frac{\indicator{Z_{s_2}=1}}{p_{s_2}}-1\right)\left(\frac{\indicator{Z_{s_3}=1}}{p_{s_3}}-1\right)\right]\right|\\
            \intertext{Since $\frac{\indicator{Z_{s_2}=1}}{p_{s_2}},\frac{\indicator{Z_{s_3}=1}}{p_{s_3}}$ take either the value 0 or a value larger than 1, this can be bounded by}
            \leq&\underbrace{\e\left[\frac{1}{p_{s_1}}\right]}_{:=S_1}+\underbrace{\e\left[\frac{1}{p_{s_1}}\frac{\indicator{Z_{s_2}=1}}{p_{s_2}}\right]}_{:=S_2}+\underbrace{\e\left[\frac{1}{p_{s_1}}\frac{\indicator{Z_{s_3}=1}}{p_{s_3}}\right]}_{:=S_3}+\underbrace{\e\left[\frac{1}{p_{s_1}}\frac{\indicator{Z_{s_2}=1}}{p_{s_2}}\frac{\indicator{Z_{s_3}=1}}{p_{s_3}}\right]}_{:=S_4}\enspace.
        \end{align*}
        By Corollary \mainref{corollary:p-moment}*, we can bound $S_1$ as:
        \begin{align*}
            S_1\leq&2+\bb_1(\bb_2/6)^{1/4}\eta_{s_1}^{1/4}\e[\widehat{A}_{s_1-1}(0)]^{1/4}\leq 2+\bb_1(\bb_2/6)^{1/4}\eta_{s_1}^{1/4}\e[\widehat{A}_{s_1}(0)]^{1/4}\enspace.
        \end{align*}
        By the law of iterated expectations and H{\"o}lder's inequality, we can bound $S_2$ as:
        \begin{align*}
            S_2\leq&\left(\e\left[\frac{1}{p_{s_1}^4}\right]\right)^{1/4}\left(\e\left[\frac{\indicator{Z_{s_2}=1}}{p_{s_2}^{4/3}}\right]\right)^{3/4}\quad(\text{H{\"o}lder's inequality})\\
            =&\left(\e\left[\frac{1}{p_{s_1}^4}\right]\right)^{1/4}\left(\e\left[\frac{1}{p_{s_2}^{1/3}}\right]\right)^{3/4}\\
            \leq&\left(2+\bb_1(\bb_2/6)^{1/4}\eta_{s_1}^{1/4}\e[\widehat{A}_{s_1}(0)]^{1/4}\right)\left(2+\bb_1(\bb_2/6)^{1/4}\eta_{s_2}^{1/4}\e[\widehat{A}_{s_2}(0)]^{1/4}\right)^{1/4}\enspace.
        \end{align*}
        Similarly, $S_3$ can be bounded as:
        \begin{align*}          
            S_3\leq&\left(2+\bb_1(\bb_2/6)^{1/4}\eta_{s_1}^{1/4}\e[\widehat{A}_{s_1}(0)]^{1/4}\right)\left(2+\bb_1(\bb_2/6)^{1/4}\eta_{s_3}^{1/4}\e[\widehat{A}_{s_3}(0)]^{1/4}\right)^{1/4}\enspace.
        \end{align*}
        By H{\"o}lder's inequality, we can bound $S_4$ as:
        \begin{align*}
            S_4\leq&\left(\e\left[\frac{1}{p_{s_1}^4}\right]\right)^{1/4}\left(\e\left[\frac{\indicator{Z_{s_2}=1}}{p_{s_2}^{4/3}}\frac{\indicator{Z_{s_3}=1}}{p_{s_3}^{4/3}}\right]\right)^{3/4}\quad(\text{H{\"o}lder's inequality})\\
            \leq&\left(\e\left[\frac{1}{p_{s_1}^4}\right]\right)^{1/4}\left(\e\left[\frac{1}{p_{s_2}^{1/3}}\frac{\indicator{Z_{s_3}=1}}{p_{s_3}^{4/3}}\right]\right)^{3/4}\quad(\text{law of iterated expectations})\\
            \leq&\left(\e\left[\frac{1}{p_{s_1}^4}\right]\right)^{1/4}\left(\e\left[\frac{1}{p_{s_2}^4}\right]\right)^{\frac{3}{4}\cdot \frac{1}{12}}\left(\e\left[\frac{\indicator{Z_{s_3}=1}}{p_{s_3}^{\frac{4}{3}\cdot\frac{12}{11}}}\right]\right)^{\frac{3}{4}\cdot \frac{11}{12}}\quad(\text{H{\"o}lder's inequality})\\
            =&\left(\e\left[\frac{1}{p_{s_1}^4}\right]\right)^{1/4}\left(\e\left[\frac{1}{p_{s_2}^4}\right]\right)^{\frac{3}{4}\cdot \frac{1}{12}}\left(\e\left[\frac{1}{p_{s_3}^{\frac{5}{11}}}\right]\right)^{\frac{3}{4}\cdot \frac{11}{12}}\\
            \leq&\left(2+\bb_1(\bb_2/6)^{1/4}\eta_{s_1}^{1/4}\e[\widehat{A}_{s_1}(0)]^{1/4}\right)\left(2+\bb_1(\bb_2/6)^{1/4}\eta_{s_2}^{1/4}\e[\widehat{A}_{s_2}(0)]^{1/4}\right)^{1/4}\\
            &\times\left(2+\bb_1(\bb_2/6)^{1/4}\eta_{s_3}^{1/4}\e[\widehat{A}_{s_3}(0)]^{1/4}\right)^{5/16}\enspace.
        \end{align*}
        Since $\widehat{A}_{s_1}(0)\geq \widehat{A}_{s_2}(0)$ and $\widehat{A}_{s_1}(0)\geq \widehat{A}_{s_3}(0)$, we have
        \begin{align*}
            &\left|\operatorname{Cov}\left(\left(\frac{\indicator{Z_{s_1}=1}}{p_{s_1}}-1\right)\left(\frac{\indicator{Z_{s_2}=1}}{p_{s_2}}-1\right),\left(\frac{\indicator{Z_{s_1}=1}}{p_{s_1}}-1\right)\left(\frac{\indicator{Z_{s_3}=1}}{p_{s_3}}-1\right)\right)\right|\\
            \leq& \left(2+\bb_1(\bb_2/6)^{1/4}\eta_{s_1}^{1/4}\e[\widehat{A}_{s_1}(0)]^{1/4}\right)+\left(2+\bb_1(\bb_2/6)^{1/4}\eta_{s_1}^{1/4}\e[\widehat{A}_{s_1}(0)]^{1/4}\right)\\
            &\quad\times\Big(2+\bb_1(\bb_2/6)^{1/4}\eta_{s_2}^{1/4}\underbrace{\e[\widehat{A}_{s_1}(0)]^{1/4}}_{\text{replaces }\e[\widehat{A}_{s_2}(0)]}\Big)^{1/4}\\
            &+\left(2+\bb_1(\bb_2/6)^{1/4}\eta_{s_1}^{1/4}\e[\widehat{A}_{s_1}(0)]^{1/4}\right)\Big(2+\bb_1(\bb_2/6)^{1/4}\eta_{s_3}^{1/4}\underbrace{\e[\widehat{A}_{s_1}(0)]^{1/4}}_{\text{replaces }\e[\widehat{A}_{s_3}(0)]}\Big)^{1/4}\\
            &+\left(2+\bb_1(\bb_2/6)^{1/4}\eta_{s_1}^{1/4}\e[\widehat{A}_{s_1}(0)]^{1/4}\right)\Big(2+\bb_1(\bb_2/6)^{1/4}\eta_{s_2}^{1/4}\underbrace{\e[\widehat{A}_{s_1}(0)]^{1/4}}_{\text{replaces }\e[\widehat{A}_{s_2}(0)]}\Big)^{1/4}\\
            &\quad\times\Big(2+\bb_1(\bb_2/6)^{1/4}\eta_{s_3}^{1/4}\underbrace{\e[\widehat{A}_{s_1}(0)]^{1/4}}_{\text{replaces }\e[\widehat{A}_{s_3}(0)]}\Big)^{5/16}\\
            \leq&\left(2+\bb_1(\bb_2/6)^{1/4}\eta_{s_1}^{1/4}\e[\widehat{A}_{s_1}(0)]^{1/4}\right)\\
            &+\left(2+\bb_1(\bb_2/6)^{1/4}\eta_{s_1}^{1/4}\e[\widehat{A}_{s_1}(0)]^{1/4}\right)\Big(2+\bb_1(\bb_2/6)^{1/4}\underbrace{\eta_{s_1}^{1/4}}_{\text{replaces }\eta_{s_2}}\e[\widehat{A}_{s_1}(0)]^{1/4}\Big)^{1/4}\left(\frac{R_{s_1}}{R_{s_2}}\right)^{1/16}\\
            &+\left(2+\bb_1(\bb_2/6)^{1/4}\eta_{s_1}^{1/4}\e[\widehat{A}_{s_1}(0)]^{1/4}\right)\Big(2+\bb_1(\bb_2/6)^{1/4}\underbrace{\eta_{s_1}^{1/4}}_{\text{replaces }\eta_{s_3}}\e[\widehat{A}_{s_1}(0)]^{1/4}\Big)^{1/4}\left(\frac{R_{s_1}}{R_{s_3}}\right)^{1/16}\\
            &+\left(2+\bb_1(\bb_2/6)^{1/4}\eta_{s_1}^{1/4}\e[\widehat{A}_{s_1}(0)]^{1/4}\right)\Big(2+\bb_1(\bb_2/6)^{1/4}\underbrace{\eta_{s_1}^{1/4}}_{\text{replaces }\eta_{s_2}}\e[\widehat{A}_{s_1}(0)]^{1/4}\Big)^{1/4}\\
            &\quad\times\Big(2+\bb_1(\bb_2/6)^{1/4}\underbrace{\eta_{s_1}^{1/4}}_{\text{replaces }\eta_{s_3}}\e[\widehat{A}_{s_1}(0)]^{1/4}\Big)^{5/16}\left(\frac{R_{s_1}}{R_{s_2}}\right)^{1/16}\left(\frac{R_{s_1}}{R_{s_3}}\right)^{5/64}\\
            \leq&4\left(2+\bb_1(\bb_2/6)^{1/4}\eta_{s_1}^{1/4}\e[\widehat{A}_{s_1}(0)]^{1/4}\right)^{25/16}\left(\frac{R_{s_1}}{R_{s_2}}\right)^{5/64}\left(\frac{R_{s_1}}{R_{s_3}}\right)^{5/64}\quad(\text{since $R_{s_1}\geq R_{s_2}$, $R_{s_1}\geq R_{s_3}$})\enspace.
        \end{align*}
        Therefore, the lemma is proved.
\end{proof}

The following corollary follows directly from the proof of Proposition \mainref{prop:expected-tracking}*. 
It provides a refined asymptotic upper bound on the fourth moment of the tracking error terms and will be used in the proof of the central limit theorem.

\begin{corollary}\label{corollary:expected-tracking}
    Under Assumptions \mainref{assumption:moments}--\mainref{assumption:maximum-radius} and Condition \mainref{condition:sigmoid}, for $k\in\{0,1\}$, the following holds:
    \begin{align*}
        \e\left[\sum_{s=1}^tR_s^{-1}\iprod{\xv_s,\bv_s^*(k)-\bv_s(k)}^4\right]=\bigO{T^{3/8}R_t^{1/4}}\enspace.
    \end{align*}
\end{corollary}

\begin{proof}
    By the proof in Proposition \mainref{prop:expected-tracking}*, the result is already verified.
\end{proof}

The following lemma captures the covariance structure between inverse probability weighting terms, which is crucial in bounding the variance of the estimated squared residuals (Lemma \ref{lemma:variance-estimated-squared-residual}).
It shows that most fourth-order covariance terms vanish.

\begin{lemma}\label{lemma:covariance}
    For any $1\leq t_1\neq t_2,t_3\neq t_4\leq T$, if the largest index among $\{t_1,t_2,t_3,t_4\}$ is unique, then
    \begin{align*}
        \operatorname{Cov}\left(\left(\frac{\mathbf{1}[Z_{t_1}=1]}{p_{t_1}}-1\right)\left(\frac{\mathbf{1}[Z_{t_2}=1]}{p_{t_2}}-1\right),\left(\frac{\mathbf{1}[Z_{t_3}=1]}{p_{t_3}}-1\right)\left(\frac{\mathbf{1}[Z_{t_4}=1]}{p_{t_4}}-1\right)\right)=0\enspace.
    \end{align*}
\end{lemma}

\begin{proof}
    We assume without loss of generality that $t_1$ is the unique largest index among $\{t_1,t_2,t_3,t_4\}$.
    For simplicity, denote
    \begin{align*}
      e_1=&\left(\frac{\mathbf{1}[Z_{t_1}=1]}{p_{t_1}}-1\right)\left(\frac{\mathbf{1}[Z_{t_2}=1]}{p_{t_2}}-1\right)\enspace,\\
      e_2=&\left(\frac{\mathbf{1}[Z_{t_3}=1]}{p_{t_3}}-1\right)\left(\frac{\mathbf{1}[Z_{t_4}=1]}{p_{t_4}}-1\right)\enspace.
    \end{align*}
    By the law of total covariance, we have
    \begin{align}\label{lemma:covariance_eq1}
        \operatorname{Cov}(e_1,e_2)=\e\left[\operatorname{Cov}\left(e_1,e_2|\mathcal{F}_{t_1-1}\right)\right]+\operatorname{Cov}\left(\e\left[e_1|\mathcal{F}_{t_1-1}\right],\e\left[e_2|\mathcal{F}_{t_1-1}\right]\right)\enspace.
    \end{align}
    Note that $\max\{t_3,t_4\}<t_1$ implies that $e_2$ is measurable with respect to $\filt_{t_1-1}$.
    This indicates that the first term in \eqref{lemma:covariance_eq1} equals 0.
    On the other hand, since $t_2<t_1$, we have
    \begin{align*}
        \e\left[e_1|\mathcal{F}_{t_1-1}\right]=\left(\frac{\mathbf{1}[Z_{t_2}=1]}{p_{t_2}}-1\right)\e\left[\frac{\mathbf{1}[Z_{t_1}=1]}{p_{t_1}}-1\Bigg|\mathcal{F}_{t_1-1}\right]=0\enspace,
    \end{align*}
    which implies that the second term in \eqref{lemma:covariance_eq1} also equals 0.
    This completes the proof.
\end{proof}

\subsection{Almost Sure Bounds for Inverse Probabilities}\label{section:D2}
The following proposition derives the almost-sure bounds for the inverse probabilities.

\begin{refproposition}{\mainref{prop:as-inv-prob-bound}}
    \asinvprobbound
\end{refproposition}

\begin{proof}
    For $k\in\{0,1\}$, denote $\Delta_{t,k}(\bv)=y_t(k)-\iprod{\xv_t,\bv}$.
    Further denote $\widehat{\Delta}_{t,1}(\bv)=y_t(1)\cdot \frac{\indicator{Z_t=1}}{p_t}-\iprod{\xv_t,\bv}$ and $\widehat{\Delta}_{t,0}(\bv)=y_t(0)\cdot \frac{\indicator{Z_t=0}}{1-p_t}-\iprod{\xv_t,\bv}$.
    We prove the result for sufficiently large $T$, since the case of small $T$ can be handled by choosing the constant $K$ large enough. 
    The proof proceeds in three steps.
    First, we establish a coarse uniform bound for $1/p_t$ and $1/(1-p_t)$.
    Next, we refine this bound recursively.
    Finally, we derive time-indexed bounds for each $t\in[T]$.
    For simplicity, we denote $\omega_1:=\bb_1\bb_2^{1/4}$ and $\omega_2:=\bb_1(\bb_2/\bb_3)^{1/2}$.\\[3mm]
    \textbf{Step 1:} Obtain the coarse uniform bound.\\[3mm] Let $\mu_1>2$ denote the largest solution to the equation: $2+\omega_1c_1^{1/2}\mu_1^{3/4}T^{1/8}=\mu_1$.
    The verification of the existence of such $\mu_1$ follows arguments similar to those in the proof of Lemma \ref{lemma:squared-residuals-random}.
    We then use the induction method to prove that $\max_{t\in[T]}\left(\frac{1}{p_t}\vee\frac{1}{1-p_t}\right)\leq \mu_1$, which serves as the initial coarse uniform bound.
    The claim holds for $t=1$ since $p_1=1/2$. 
    If the result holds for $1,\ldots,t$, we then verify it at $t+1$. 
    
    We first derive an upper bound for $\widehat{A}_t(1)$. 
    By the induction assumption and the inequality $(a+b)^2\leq 2a^2+2b^2$, we have
    \begin{align}\label{prop:as-inv-prob-bound_eq}
        \widehat{A}_{t}(1)=&\sum_{s=1}^{t}\frac{\mathbf{1}[Z_s=1]}{p_s}\cdot\Delta_{s,1}(\bv_s(1))^2\notag\\
        \leq&\mu_1\sum_{s=1}^{t}\Delta_{s,1}(\bv_s(1))^2~~~~(\text{induction assumption})\notag\\
        =&\mu_1\sum_{s=1}^{t}\left[\widehat{\Delta}_{s,1}(\bv_s(1))-y_s(1)\left(\frac{\mathbf{1}[Z_s=1]}{p_s}-1\right)\right]^2\notag\\
        \leq&2\mu_1\sum_{s=1}^{t}\widehat\Delta_{s,1}(\bv_s(1))^2+2\mu_1\sum_{s=1}^{t}y_s(1)^2\left(\frac{\mathbf{1}[Z_s=1]}{p_s}-1\right)^2\enspace.
    \end{align}
    For simplicity, we introduce the following notations:
    \begin{align*}
        \widehat{G}_t^{(1)}(\bv)&=\sum_{s=1}^{t-1} \widehat{\Delta}_{s,1}(\bv)^2+\eta_t^{-1}\|\bv\|^2\enspace,\\
        \widebar{G}_t^{(1)}(\bv)&=\sum_{s=1}^{t-1} \widehat{\Delta}_{s,1}(\bv)^2+\eta_{t-1}^{-1}\|\bv\|^2\enspace.
    \end{align*}
    Denote by $\widebar{\bv}_t(1)$ the minimizer of $\widebar{G}_t^{(1)}(\bv)$.
    By Corollary \ref{corollary:standard-FTRL-upper-bound} and Lemma \ref{lemma:ridge}, we have
    \begin{align*}
        &\sum_{s=1}^t\widehat\Delta_{s,1}(\bv_s(1))^2-\sum_{s=1}^t\widehat{\Delta}_{s,1}(\bv_{t+1}(1))^2\notag\\
        \leq&\eta_{t+1}^{-1}\|\bv_{t+1}(1)\|^2+\sum_{s=1}^t(\widebar{G}_{s+1}^{(1)}(\bv_s(1))-\widebar{G}_{s+1}^{(1)}(\widebar{\bv}_{s+1}(1)))\quad(\text{Corollary \ref{corollary:standard-FTRL-upper-bound}})\notag\\
        =&\eta_{t+1}^{-1}\|\bv_{t+1}(1)\|^2+\sum_{s=1}^t\frac{\xv_s^{\tran}\left(\xM_{s-1}^{\tran}\xM_{s-1}+\eta_s^{-1}\mat{I}_d\right)^{-1}\xv_s}{1+\xv_s^{\tran}\left(\xM_{s-1}^{\tran}\xM_{s-1}+\eta_s^{-1}\mat{I}_d\right)^{-1}\xv_s}
        \cdot\widehat\Delta_{s,1}(\bv_s(1))^2\quad(\text{Lemma \ref{lemma:ridge}})\notag\\
        \leq&\eta_{t+1}^{-1}\|\bv_{t+1}(1)\|^2+\sum_{s=1}^t\xv_s^{\tran}\left(\xM_{s-1}^{\tran}\xM_{s-1}+\eta_s^{-1}\mat{I}_d\right)^{-1}\xv_s\cdot\widehat\Delta_{s,1}(\bv_s(1))^2\notag\\
        \leq&\eta_{t+1}^{-1}\|\bv_{t+1}(1)\|^2+\sum_{s=1}^tR_s^{2}\eta_s\cdot\widehat\Delta_{s,1}(\bv_s(1))^2\notag\\
        \intertext{Since $R_s^2\eta_s=T^{-1/2}R_s\leq T^{-1/2}R_T$, which is smaller than $1/2$ for sufficiently large $T$ by Lemma \ref{lemma:R} and Assumption \mainref{assumption:maximum-radius-strengthened}, we obtain}
        \leq&\eta_{t+1}^{-1}\|\bv_{t+1}(1)\|^2+\frac{1}{2}\sum_{s=1}^t\widehat\Delta_{s,1}(\bv_s(1))^2\enspace,
    \end{align*}
    which is equivalent to
    \begin{align}\label{prop:as-inv-prob-bound_eq1}
        \frac{1}{2}\sum_{s=1}^t\widehat\Delta_{s,1}(\bv_s(1))^2\leq\sum_{s=1}^{t}\widehat{\Delta}_{s,1}(\bv_{t+1}(1))^2+\eta_{t+1}^{-1}\left\|\bv_{t+1}(1)\right\|^2\enspace.
    \end{align}
    By the induction assumption and the definition of $\bv_{t+1}(1)$, we have
    \begin{align}\label{prop:as-inv-prob-bound_eq2}
        &\sum_{s=1}^{t}\widehat{\Delta}_{s,1}(\bv_{t+1}(1))^2+\eta_{t+1}^{-1}\left\|\bv_{t+1}(1)\right\|^2\notag\\
        =&\min_{\bv\in\mathbb{R}^d}\left\{\sum_{s=1}^{t}\widehat{\Delta}_{s,1}(\bv)^2+\eta_{t+1}^{-1}\left\|\bv\right\|^2\right\}\quad(\text{optimality of $\bv_{t+1}(1)$})\notag\\
        \leq&\sum_{s=1}^{t}\left(y_s(1)\cdot \frac{\mathbf{1}[Z_s=1]}{p_s}\right)^2\quad(\text{evaluation at $\bv=\vec{0}$})\notag\\
        \leq&\left(\sum_{s=1}^{t}y_s(1)^2\right)\left(\max_{s=1,\ldots,t}\frac{1}{p_s}\right)^2\notag\\
        \leq&c_1^2\mu_1^2T~~~~(\text{induction assumption and Assumption \mainref{assumption:moments}})\enspace.
    \end{align}
    Hence by \eqref{prop:as-inv-prob-bound_eq1} and \eqref{prop:as-inv-prob-bound_eq2}, we have
    \begin{align}\label{prop:as-inv-prob-bound_eq3}
        \sum_{s=1}^{t}\widehat\Delta_{s,1}(\bv_s(1))^2\leq 2\sum_{s=1}^{t}\widehat{\Delta}_{s,1}(\bv_{t+1}(1))^2+2\eta_{t+1}^{-1}\left\|\bv_{t+1}(1)\right\|^2\leq 2c_1^2\mu_1^2T\enspace.
    \end{align}
    By induction assumption and Assumption \mainref{assumption:moments}, we also have
    \begin{align}\label{prop:as-inv-prob-bound_eq4}
        \sum_{s=1}^{t}y_s(1)^2\left(\frac{\mathbf{1}[Z_s=1]}{p_s}-1\right)^2\leq \mu_1^2\sum_{s=1}^{t}y_s(1)^2\leq c_1^2\mu_1^2T\enspace.
    \end{align}
    Hence by \eqref{prop:as-inv-prob-bound_eq}, \eqref{prop:as-inv-prob-bound_eq3} and \eqref{prop:as-inv-prob-bound_eq4}, we have
    \begin{align*}
        \widehat{A}_t(1)\leq&2\mu_1\sum_{s=1}^{t}\widehat\Delta_{s,1}(\bv_s(1))^2+2\mu_1\sum_{s=1}^{t}y_s(1)^2\left(\frac{\mathbf{1}[Z_s=1]}{p_s}-1\right)^2\leq 6c_1^2\mu_1^3T\enspace.
    \end{align*}
    Then by Lemma \mainref{lemma:effect-of-p-regularization}*, we have
    \begin{align*}
        \frac{1}{1-p_{t+1}}\leq& 2+\bb_1(\bb_2/6)^{1/4}\eta_{t+1}^{1/4}\widehat{A}_t(1)^{1/4}\\
        \leq& 2+\omega_1\eta_1^{1/4}c_1^{1/2}\mu_1^{3/4}T^{1/4}\quad(\text{since $\eta_{t+1}\leq \eta_1$})\\
        \leq&2+\omega_1c_1^{1/2}\mu_1^{3/4}T^{1/8}\\
        =&\mu_1\quad(\text{definition of $\mu_1$})\enspace.
    \end{align*}
    Similarly, we can prove that $\widehat{A}_t(0)$ is bounded by $6c_1^2\mu_1^3T$ and $1/p_{t+1}$ is bounded by $\mu_1$. Then by induction, we have verified that
    \begin{align*}
        \max_{t\in[T]}\left(\frac{1}{p_t}\vee\frac{1}{1-p_t}\right)\leq \mu_1\enspace.
    \end{align*}
    \textbf{Step 2:} Refine the uniform bound recursively.\\[3mm]
    Suppose the upper bound obtained in the $r$-th round is $\mu_r$.
    Then we refine this bound to $\mu_{r+1}$ in the $(r+1)$-th round through the following procedure. 
    Repeating the arguments in Step 1 leads to $\widehat{A}_T(1)\leq 6c_1^2\mu_r^3T$ and $\widehat{A}_T(0)\leq 6c_1^2\mu_r^3T$. 
    We fix a $T$-dependent constant $\delta>0$ (to be determined later). 
    For any $t\in [T]$, if $\widehat{A}_{t-1}(1)\geq 6c_1^2\mu_r^3T\delta^{-2}$, by Lemma \mainref{lemma:effect-of-p-regularization}* we have
    \begin{align*}
        \frac{1}{p_t}\leq 2+\omega_2\left(\frac{\widehat{A}_{t-1}(0)}{\widehat{A}_{t-1}(1)}\right)^{1/2}\leq 2+\omega_2\left(\frac{6c_1^2\mu_r^3T}{6c_1^2\mu_r^3T\delta^{-2}}\right)^{1/2}=2+\omega_2\delta\enspace.
    \end{align*}
    If $\widehat{A}_{t-1}(1)\leq 6c_1^2\mu_r^3T\delta^{-2}$, then we have
    \begin{align*}
        \frac{1}{p_{t}}\leq& 2+\bb_1(\bb_2/6)^{1/4}\eta_t^{1/4}\widehat{A}_{t-1}(0)^{1/4}\\
        \intertext{The proof in Step 1 indicates that $\widehat{A}_{t-1}(0)\leq 6c_1^2T\cdot \left[\max_{s=1,\ldots,t-1}\frac{1}{1-p_s}\right]^3$. Hence}
        \leq&2+\bb_1(\bb_2/6)^{1/4}\eta_t^{1/4}\cdot \left(6c_1^2T\cdot \left[\max_{s=1,\ldots,t-1}\frac{1}{1-p_s}\right]^3\right)^{1/4}\\
        \intertext{Since $\frac{1}{1-p_s}\leq 2+\bb_1(\bb_2/6)^{1/4}\eta_s^{1/4}\widehat{A}_{s-1}(1)^{1/4}$ by Lemma \mainref{lemma:effect-of-p-regularization}* and $\eta_t\leq \eta_1$, we obtain}
        \leq&2+\bb_1\bb_2^{1/4}c_1^{1/2}\eta_1^{1/4}T^{1/4}\max_{s=1,\ldots,t-1}\left(2+\bb_1(\bb_2/6)^{1/4}\eta_s^{1/4}\widehat{A}_{s-1}(1)^{1/4}\right)^{3/4}\\
        \leq&2+\omega_1c_1^{1/2}\eta_1^{1/4}T^{1/4}\Big(2+\bb_1(\bb_2/6)^{1/4}\underbrace{\eta_1^{1/4}}_{\eta_{s}\leq \eta_1}\underbrace{\widehat{A}_{t-1}(1)^{1/4}}_{\widehat{A}_{s-1}(1)\leq \widehat{A}_{t-1}(1)}\Big)^{3/4}\\
        \leq&2+\omega_1c_1^{1/2}\eta_1^{1/4}T^{1/4}\left[2+\bb_1(\bb_2/6)^{1/4}\eta_1^{1/4}\left(6c_1^2\mu_r^3T\delta^{-2}\right)^{1/4}\right]^{3/4}\\
        \leq&2+\omega_1c_1^{1/2}\eta_1^{1/4}T^{1/4}\left(2+\omega_1\eta_1^{1/4}c_1^{1/2}\mu_r^{3/4}T^{1/4}\delta^{-1/2}\right)^{3/4}\\
        \leq&2+\omega_1c_1^{1/2}T^{1/8}\left(2+\omega_1c_1^{1/2}\mu_r^{3/4}T^{1/8}\delta^{-1/2}\right)^{3/4}\enspace.
    \end{align*}
    Now we choose $\delta>0$ to be the supremum of the solution set of the equation:
    \begin{align}\label{prop:as-inv-prob-bound_eq5}
        2+\omega_2\delta=2+\omega_1c_1^{1/2}T^{1/8}\left(2+\omega_1c_1^{1/2}\mu_r^{3/4}T^{1/8}\delta^{-1/2}\right)^{3/4}\enspace.
    \end{align}
    We first prove the existence of such $\delta$.
    The left-hand side of equation \eqref{prop:as-inv-prob-bound_eq5} attains value $2$ at $\delta=0$ and tends to $+\infty$ when $\delta\rightarrow\infty$. 
    The right-hand side of equation \eqref{prop:as-inv-prob-bound_eq5} tends to infinity as $\delta\downarrow 0$ and tends to $2+2^{3/4}\omega_1c_1^{1/2}T^{1/8}$ when $\delta\rightarrow\infty$ for any fixed $T$. 
    Hence the solution to \eqref{prop:as-inv-prob-bound_eq5} exists by the intermediate value theorem for continuous functions.
    By continuity, the supremum $\delta$ is itself a solution.
    Then for any $t\in[T]$, we have $1/p_t\leq 2+\omega_2\delta:= \mu_{r+1}$. 
    Similarly, we can prove that for any $t\in[T]$, we have $1/(1-p_t)\leq 2+\omega_2\delta=\mu_{r+1}$.
    
    Hence, we have obtained a new uniform bound $\mu_{r+1}$ for $1/p_t$ and $1/(1-p_t)$. 
    By such recursion, we can obtain a sequence of upper bounds: $\{\mu_r:r=1,2,\ldots\}$.
    Our goal is to prove that this sequence is decreasing and converges to some limit, which serves as the refined uniform bound.
    
    Suppose $\mu>2$ is the supremum of the solution set of the equation:
    \begin{align}\label{prop:as-inv-prob-bound_eq6}
        \mu=2+\omega_1c_1^{1/2}T^{1/8}\left(2+\omega_1c_1^{1/2}\mu^{3/4}T^{1/8}(\omega_2^{-1}(\mu-2))^{-1/2}\right)^{3/4}\enspace.
    \end{align}
    We first prove the existence of such solutions. 
    The left-hand side of \eqref{prop:as-inv-prob-bound_eq6} attains value $2$ at $\mu=2$, while the right-hand side tends to infinity when $\mu\downarrow 2$. 
    Meanwhile, for any fixed $T$, the right-hand side is of order $\bigO{\mu^{3/16}}$ while the left-hand side is of order $\bigO{\mu}$. 
    Hence by the intermediate value theorem for continuous functions, there exist solutions to \eqref{prop:as-inv-prob-bound_eq6} on $(2,\infty)$ for any given $T$. 
    By continuity, the supremum $\mu$ is itself a solution.
    Then we have
    \begin{align*}
        \mu=&2+\omega_1c_1^{1/2}T^{1/8}\left(2+\omega_1c_1^{1/2}\mu^{3/4}T^{1/8}(\omega_2^{-1}(\mu-2))^{-1/2}\right)^{3/4}\\
        \asymp&1+T^{1/8}\left[1+(\mu^{1/4}T^{1/8})^{3/4}\right]\\
        \asymp&T^{7/32}\mu^{3/16}\enspace.
    \end{align*}
    We can verify that $\mu_1=\Theta(T^{1/2})$ and $\mu=\bigO{T^{7/26}}$, which indicates that $\mu=o(\mu_1)$. 
    Hence for sufficiently large $T$, we may assume without loss of generality that $\mu_1>\mu$. 
    Now we prove that the sequence $\{\mu_r:r=1,2,\ldots\}$ is monotone decreasing and satisfies
    \begin{align*}
        \mu_r>2+\omega_1c_1^{1/2}T^{1/8}\left(2+\omega_1c_1^{1/2}\mu_r^{3/4}T^{1/8}(\omega_2^{-1}(\mu_r-2))^{-1/2}\right)^{3/4}
    \end{align*}
    through induction.
    By the previous discussions on the asymptotic orders of both sides of \eqref{prop:as-inv-prob-bound_eq6} and the definition of $\mu$, for any $\delta >\mu$, by continuity we have
    \begin{align*}
        \delta>2+\omega_1c_1^{1/2}T^{1/8}\left(2+\omega_1c_1^{1/2}\delta^{3/4}T^{1/8}(\omega_2^{-1}(\delta-2))^{-1/2}\right)^{3/4}\enspace.
    \end{align*}
    Since $\mu_1>\mu$, we have already verified the result when $r=1$. Suppose the result holds for $1,\ldots,r$. 
    By the definition of $\mu_{r+1}$, the following holds:
    \begin{align}\label{prop:as-inv-prob-bound_eq7}
        \mu_{r+1}=2+\omega_1c_1^{1/2}T^{1/8}\left(2+\omega_1c_1^{1/2}\mu_r^{3/4}T^{1/8}(\omega_2^{-1}(\mu_{r+1}-2))^{-1/2}\right)^{3/4}\enspace.\quad
    \end{align}
    By induction assumption, we have
    \begin{align*}
        \mu_r>2+\omega_1c_1^{1/2}T^{1/8}\left(2+\omega_1c_1^{1/2}\mu_r^{3/4}T^{1/8}(\omega_2^{-1}(\mu_r-2))^{-1/2}\right)^{3/4}\enspace.
    \end{align*}
    Hence, when evaluated at $\delta=\omega_2^{-1}(\mu_{r+1}-2)$, the left-hand side of \eqref{prop:as-inv-prob-bound_eq5} equals the right-hand side due to \eqref{prop:as-inv-prob-bound_eq7}, whereas when evaluated at $\delta=\omega_2^{-1}(\mu_{r}-2)$, the left-hand side exceeds the right-hand side. 
    Since the left-hand side of \eqref{prop:as-inv-prob-bound_eq5} is increasing in $\delta$ while the right-hand side is decreasing, it follows that $\mu_{r+1}<\mu_r$.
    Hence, by \eqref{prop:as-inv-prob-bound_eq7} we also have
    \begin{align*}
        \mu_{r+1}=&2+\omega_1c_1^{1/2}T^{1/8}\left(2+\omega_1c_1^{1/2}\mu_r^{3/4}T^{1/8}(\omega_2^{-1}(\mu_{r+1}-2))^{-1/2}\right)^{3/4}\\
        >&2+\omega_1c_1^{1/2}T^{1/8}\left(2+\omega_1c_1^{1/2}\mu_{r+1}^{3/4}T^{1/8}(\omega_2^{-1}(\mu_{r+1}-2))^{-1/2}\right)^{3/4}\enspace.
    \end{align*}
    
    Therefore, we have proved that $\{\mu_r:r=1,2,\ldots\}$ is decreasing by induction. 
    Hence the limit of the sequence exists by the monotone bounded sequence theorem.
    This limit is also a solution to \eqref{prop:as-inv-prob-bound_eq6} by the continuity property. 
    Since $\mu>2$ is the largest solution to equation \eqref{prop:as-inv-prob-bound_eq6} by definition, it follows that $\lim_{r\rightarrow \infty}\mu_r\leq \mu$. 
    Hence we can prove that
    \begin{align*}
        \max_{t\in[T]}\left(\frac{1}{p_t}\vee\frac{1}{1-p_t}\right)\leq \mu=\bigO{T^{7/26}}\enspace.
    \end{align*}
    \textbf{Step 3:} Derive time-indexed bounds for each $t\in[T]$.\\[3mm]
    If $\widehat{A}_{t-1}(1)\geq 6c_1^2\mu^3T\delta^{-2}$, by Lemma \mainref{lemma:effect-of-p-regularization}* we have
    \begin{align*}
        \frac{1}{p_t}\leq 2+\omega_2\left(\frac{\widehat{A}_{t-1}(0)}{\widehat{A}_{t-1}(1)}\right)^{1/2}\leq 2+\omega_2\left(\frac{6c_1^2\mu^3T}{6c_1^2\mu^3T\delta^{-2}}\right)^{1/2}=2+\omega_2\delta\enspace.
    \end{align*}
    If $\widehat{A}_{t-1}(1)\leq 6c_1^2\mu^3T\delta^{-2}$, by arguments similar to those in Step 2, we have
    \begin{align*}
        \frac{1}{p_{t}}\leq& 2+\bb_1(\bb_2/6)^{1/4}\eta_t^{1/4}\widehat{A}_{t-1}(0)^{1/4}\\
        \leq&2+\bb_1(\bb_2/6)^{1/4}\eta_t^{1/4}\cdot \left(6c_1^2T\cdot \left[\max_{s=1,\ldots,t-1}\frac{1}{1-p_s}\right]^3\right)^{1/4}\\
        \leq&2+\omega_1c_1^{1/2}\eta_t^{1/4}T^{1/4}\max_{s=1,\ldots,t-1}\left(2+\bb_1(\bb_2/6)^{1/4}\eta_s^{1/4}\widehat{A}_{s-1}(1)^{1/4}\right)^{3/4}\\
        \leq&2+\omega_1c_1^{1/2}\eta_t^{1/4}T^{1/4}\left(2+\bb_1(\bb_2/6)^{1/4}\eta_1^{1/4}\widehat{A}_{t-1}(1)^{1/4}\right)^{3/4}\\
        \leq&2+\omega_1c_1^{1/2}\eta_t^{1/4}T^{1/4}\left[2+\bb_1(\bb_2/6)^{1/4}\eta_1^{1/4}\left(6c_1^2\mu^3T\delta^{-2}\right)^{1/4}\right]^{3/4}\\
        \leq&2+\omega_1c_1^{1/2}\eta_t^{1/4}T^{1/4}\left(2+\omega_1c_1^{1/2}\mu^{3/4}T^{1/8}\delta^{-1/2}\right)^{3/4}\enspace.
    \end{align*}
    Let $\delta_t>0$ be the supremum of the solution set of the equation:
    \begin{align}\label{prop:as-inv-prob-bound_eq8}
        2+\omega_2\delta=2+\omega_1c_1^{1/2}\eta_t^{1/4}T^{1/4}\left(2+\omega_1c_1^{1/2}\mu^{3/4}T^{1/8}\delta^{-1/2}\right)^{3/4}\enspace.
    \end{align}
    We can prove that $\delta_t$ is the largest solution to \eqref{prop:as-inv-prob-bound_eq8} by similar arguments.
    Consider $\widetilde{\delta}=(2\omega_1^{7/4}\omega_2^{-1}c_1^{7/8}\eta_t^{1/4}T^{11/32}\mu^{9/16})^{8/11}$.
    We can verify that both $\mu^{3/4}T^{1/8}\widetilde{\delta}^{-1/2}$ and $\eta_t^{1/4}T^{1/4}$ tend to infinity when $T$ goes to infinity.
    Hence, for sufficiently large $T$, we have
    \begin{align*}
        &(2+\omega_2\widetilde{\delta})-\left(2+\omega_1c_1^{1/2}\eta_t^{1/4}T^{1/4}\left(2+\omega_1c_1^{1/2}\mu^{3/4}T^{1/8}\widetilde{\delta}^{-1/2}\right)^{3/4}\right)\\
        \geq&\omega_2\widetilde{\delta}-2\omega_1c_1^{1/2}\eta_t^{1/4}T^{1/4}\left(\omega_1c_1^{1/2}\mu^{3/4}T^{1/8}\widetilde{\delta}^{-1/2}\right)^{3/4}\quad(\text{for sufficiently large $T$})\\
        =&\omega_2\widetilde{\delta}-2\omega_1^{7/4}c_1^{7/8}\eta_t^{1/4}T^{11/32}\mu^{9/16}\widetilde{\delta}^{-3/8}\\
        =&0\enspace.
    \end{align*}
    Since the left-hand side of \eqref{prop:as-inv-prob-bound_eq8} is monotone increasing in $\delta$ while the right-hand side is monotone decreasing in $\delta$, this implies that $\delta_t\leq \widetilde{\delta}$.
    Hence we have
    \begin{align*}
        \max\left\{\frac{1}{p_t},\frac{1}{1-p_t}\right\}\leq 2+\omega_2\delta_t\leq& 2+\omega_2\widetilde{\delta}\\
        =&2+\omega_2(2\omega_1^{7/4}\omega_2^{-1}c_1^{7/8}\eta_t^{1/4}T^{11/32}\mu^{9/16})^{8/11}\\
        =&2+\omega_2(2\omega_1^{7/4}\omega_2^{-1}c_1^{7/8}T^{7/32}R_t^{-1/4}\mu^{9/16})^{8/11}\\
        =&2+\omega_2(2\omega_1^{7/4}\omega_2^{-1}c_1^{7/8})^{8/11}T^{7/44}\mu^{9/22}R_t^{-2/11}\enspace.
    \end{align*}
    Since $\mu=\bigO{T^{7/26}}$, we can show that the right-hand side is of order $\bigO{T^{7/26}R_t^{-2/11}}$, which implies that there exists a constant $K>0$ (independent of $t$ and $T$) such that
    \begin{align*}
	   \Pr[\Big]{ \max \braces[\Big]{ \frac{1}{p_t}, \frac{1}{1-p_t} } \leq K \cdot T^{7/26} R_t^{-2/11} \ \text{ for all } t \in [T]} = 1 \enspace.
	\end{align*}
    This completes the proof.
\end{proof}

\subsection{Central Limit Theorem}\label{section:D3}
In this section, we prove a central limit theorem for the AIPW estimator.
We begin by stating the martingale central limit theorem, which is the main tool for our proof.

\begin{reflemma}{\mainref{lemma:martingale-clt}}[\cite{helland1982central}]
	\martingaleclt
\end{reflemma}

Throughout this section, we omit the subscript $T$ and write $X_{t,T}$ and $\filt_{t,T}$ as $X_t$ and $\filt_t$, respectively.
For $k\in\{0,1\}$, denote $\Delta_{t,k}(\bv)=y_t(k)-\iprod{\xv_t,\bv}$.
In Section 5, we have already defined the martingale difference sequence:

\begin{align*}
X_{t} = \frac{ \eate_t - \ate_t }{T \sqrt{ \Var{\eate} } }=\frac{1}{T \sqrt{ \Var{\eate} } }\left[\Delta_{t,1}(\bv_t(1)) \cdot \left(\frac{ \indicator{Z_t=1} }{p_t}-1\right)- \Delta_{t,0}(\bv_t(0)) \cdot \left(\frac{ \indicator{Z_t=0} }{1-p_t}-1\right)\right]
\enspace.
\end{align*}

In the next two subsections, we verify the two conditions in Lemma \mainref{lemma:martingale-clt}, i.e., the conditional variance condition and the conditional Lyapunov condition. 
We first state, without proof, the non-superefficiency condition in the following corollary. Its proof is deferred to Section \ref{section:D4}.

\begin{refcorollary}{\mainref{corollary:non-superefficiency}}[Non-Superefficiency]
	\nonsuperefficiency
\end{refcorollary}

\subsubsection{Conditional Variance Convergence}
In this section, we verify the conditional variance convergence for the martingale central limit theorem.
We begin by deriving a simplified form for the conditional variance condition.
Following the derivation in Proposition \mainref{prop:aipw-variance}, we obtain
\begin{align*}
    &T\cdot\operatorname{Var}(\hat{\tau})=\frac{1}{T}\sum_{t=1}^T\e\Bigg[\Delta_{t,1}(\bv_t(1))^2\cdot \left(\frac{1}{p_t}-1\right)+\Delta_{t,0}(\bv_t(0))^2\cdot \left(\frac{1}{1-p_t}-1\right)\\
    &\hspace{1in}+2\Delta_{t,1}(\bv_t(1))\cdot\Delta_{t,0}(\bv_t(0))\Bigg]\enspace,\\
    &\sum_{t=1}^T\operatorname{Var}(\hat{\tau}_t|\filt_{t-1})=\sum_{t=1}^T\Bigg[\Delta_{t,1}(\bv_t(1))^2\cdot \left(\frac{1}{p_t}-1\right)+\Delta_{t,0}(\bv_t(0))^2\cdot \left(\frac{1}{1-p_t}-1\right)\\
    &\hspace{1.1in}+2\Delta_{t,1}(\bv_t(1))\cdot\Delta_{t,0}(\bv_t(0))\Bigg]\enspace.
\end{align*}
Hence, the term $V_T^2$ can be simplified as:
\begin{align*}
    &V_T^2\\
    =&\frac{1}{T^2\operatorname{Var}(\hat{\tau})}\sum_{t=1}^T\left[\Delta_{t,1}(\bv_t(1))^2\cdot \left(\frac{1}{p_t}-1\right)+\Delta_{t,0}(\bv_t(0))^2\cdot \left(\frac{1}{1-p_t}-1\right)+2\Delta_{t,1}(\bv_t(1))\cdot\Delta_{t,0}(\bv_t(0))\right]\\
    =&1+\frac{1}{T^2\operatorname{Var}(\hat{\tau})}\Bigg(\sum_{t=1}^T\Delta_{t,1}(\bv_t(1))^2\cdot \left(\frac{1}{p_t}-1\right)-\e\left[\sum_{t=1}^T\Delta_{t,1}(\bv_t(1))^2\cdot \left(\frac{1}{p_t}-1\right)\right]\\
    &\hspace{0.9in}+\sum_{t=1}^T\Delta_{t,0}(\bv_t(0))^2\cdot \left(\frac{1}{1-p_t}-1\right)-\e\left[\sum_{t=1}^T\Delta_{t,0}(\bv_t(0))^2\cdot \left(\frac{1}{1-p_t}-1\right)\right]\\
    &\hspace{0.9in}+2\sum_{t=1}^T\Delta_{t,1}(\bv_t(1))\cdot\Delta_{t,0}(\bv_t(0))-2\e\left[\sum_{t=1}^T\Delta_{t,1}(\bv_t(1))\cdot\Delta_{t,0}(\bv_t(0))\right]\Bigg)\enspace.
\end{align*}

Under the non-superefficiency condition (Corollary \mainref{corollary:non-superefficiency}), we have $\operatorname{Var}(\hat{\tau})=\Omega(T^{-1})$.
In order to prove that $V_T^2 \xrightarrow{p} 1$, it suffices to establish the following three convergence results:
\begin{align}\label{simplified-stable-1}
    \qquad&\frac{1}{T}\left[\sum_{t=1}^T\Delta_{t,1}(\bv_t(1))^2\cdot \left(\frac{1}{p_t}-1\right)-\e\left[\sum_{t=1}^T\Delta_{t,1}(\bv_t(1))^2\cdot \left(\frac{1}{p_t}-1\right)\right]\right]\xrightarrow{p}0\enspace,\notag\\
    &\frac{1}{T}\left[\sum_{t=1}^T\Delta_{t,0}(\bv_t(0))^2\cdot \left(\frac{1}{1-p_t}-1\right)-\e\left[\sum_{t=1}^T\Delta_{t,0}(\bv_t(0))^2\cdot \left(\frac{1}{1-p_t}-1\right)\right]\right]\xrightarrow{p}0\enspace,\notag\\
    &\frac{1}{T}\left[\sum_{t=1}^T\Delta_{t,1}(\bv_t(1))\cdot \Delta_{t,0}(\bv_t(0))-\e\left[\sum_{t=1}^T\Delta_{t,1}(\bv_t(1))\cdot\Delta_{t,0}(\bv_t(0))\right]\right]\xrightarrow{p}0\enspace.
\end{align}
These are the treated-group squared term, the control-group squared term, and the cross term.

We begin with the first convergence result in \eqref{simplified-stable-1}. Using the equality $\Delta_{t,1}(\bv_t(1))-\Delta_{t,1}(\bv_t^*(1))=\iprod{\xv_t,\bv_t^*(1)-\bv_t(1)}$, we obtain the following decomposition:
\begin{align*}
    &\frac{1}{T}\sum_{t=1}^T\Delta_{t,1}(\bv_t(1))^2\cdot\left(\frac{1}{p_t}-1\right)\\
    =&\frac{1}{T}\sum_{t=1}^T\Delta_{t,1}(\bv_t^*(1))^2\cdot\left(\frac{1}{p_t}-1\right)+\frac{1}{T}\sum_{t=1}^T\left[\Delta_{t,1}(\bv_t(1))^2-\Delta_{t,1}(\bv_t^*(1))^2\right]\left(\frac{1}{p_t}-1\right)\\
    \intertext{Using equality $\frac{1}{p_t}-1=\left(\frac{1}{p_t}-\frac{1}{\widebar{p}_t}\right)+\left(\frac{1}{\widebar{p}_t}-1\right)$ and $a^2-b^2=(a+b)(a-b)$, we obtain}
    =&\frac{1}{T}\sum_{t=1}^T\Delta_{t,1}(\bv_t^*(1))^2\cdot\left(\frac{1}{p_t}-\frac{1}{\widebar{p}_t}\right)+\frac{1}{T}\sum_{t=1}^T\Delta_{t,1}(\bv_t^*(1))^2\cdot\left(\frac{1}{\widebar{p}_t}-1\right)\\
    &+\frac{1}{T}\sum_{t=1}^T\underbrace{(2\Delta_{t,1}(\bv_t^*(1))+\iprod{\xv_t,\bv_t^*(1)-\bv_t(1)})}_{=\Delta_{t,1}(\bv_t^*(1))+\Delta_{t,1}(\bv_t(1))}\cdot \iprod{\xv_t,\bv_t^*(1)-\bv_t(1)}\cdot\left(\frac{1}{p_t}-1\right)\\
    =&\frac{1}{T}\sum_{t=1}^T\Delta_{t,1}(\bv_t^*(1))^2\cdot\left(\frac{1}{p_t}-\frac{1}{\widebar{p}_t}\right)+\frac{1}{T}\sum_{t=1}^T\Delta_{t,1}(\bv_t^*(1))^2\cdot\left(\frac{1}{\widebar{p}_t}-1\right)\\
    &+\frac{2}{T}\sum_{t=1}^T\Delta_{t,1}(\bv_t^*(1))\cdot\iprod{\xv_t,\bv_t^*(1)-\bv_t(1)}\cdot\left(\frac{1}{p_t}-\frac{1}{\widebar{p}_t}\right)\\
    &+\frac{2}{T}\sum_{t=1}^T\Delta_{t,1}(\bv_t^*(1))\cdot\iprod{\xv_t,\bv_t^*(1)-\bv_t(1)}\cdot\left(\frac{1}{\widebar{p}_t}-1\right)\\
    &+\frac{1}{T}\sum_{t=1}^T\iprod{\xv_t,\bv_t^*(1)-\bv_t(1)}^2\left(\frac{1}{p_t}-\frac{1}{\widebar{p}_t}\right)+\frac{1}{T}\sum_{t=1}^T\iprod{\xv_t,\bv_t^*(1)-\bv_t(1)}^2\left(\frac{1}{\widebar{p}_t}-1\right)\enspace.
\end{align*}
Hence, it suffices to show that each centered term converges to zero in probability, i.e.,

(1) Full-information online residuals multiplied by the difference between the random inverse assignment probability and its fixed comparator.
\begin{align*}
    \frac{1}{T}\sum_{t=1}^T\Delta_{t,1}(\bv_t^*(1))^2\cdot\left(\frac{1}{p_t}-\frac{1}{\widebar{p}_t}\right)-\e\left[\frac{1}{T}\sum_{t=1}^T\Delta_{t,1}(\bv_t^*(1))^2\cdot\left(\frac{1}{p_t}-\frac{1}{\widebar{p}_t}\right)\right]\xrightarrow{p}0\enspace.
\end{align*}

(2) Full-information online residuals multiplied by the fixed comparator term.
\begin{align*}
    \frac{1}{T}\sum_{t=1}^T\Delta_{t,1}(\bv_t^*(1))^2\cdot\left(\frac{1}{\widebar{p}_t}-1\right)-\e\left[\frac{1}{T}\sum_{t=1}^T\Delta_{t,1}(\bv_t^*(1))^2\cdot\left(\frac{1}{\widebar{p}_t}-1\right)\right]\xrightarrow{p}0\enspace.
\end{align*}

(3) Full-information online residuals combined with tracking error terms, multiplied by the difference between the random inverse assignment probability and its fixed comparator.
\begin{align*}
    &\frac{1}{T}\sum_{t=1}^T\Delta_{t,1}(\bv_t^*(1))\cdot\iprod{\xv_t,\bv_t(1)-\bv^*_t(1)}\cdot\left(\frac{1}{p_t}-\frac{1}{\widebar{p}_t}\right)\notag\\
    &-\e\left[\frac{1}{T}\sum_{t=1}^T\Delta_{t,1}(\bv_t^*(1))\cdot\iprod{\xv_t,\bv_t(1)-\bv^*_t(1)}\cdot\left(\frac{1}{p_t}-\frac{1}{\widebar{p}_t}\right)\right]\xrightarrow{p}0\enspace.
\end{align*}

(4) Full-information online residuals combined with tracking error terms, multiplied by the fixed comparator term.
\begin{align*}
    &\frac{1}{T}\sum_{t=1}^T\Delta_{t,1}(\bv_t^*(1))\cdot\iprod{\xv_t,\bv_t(1)-\bv^*_t(1)}\cdot\left(\frac{1}{\widebar{p}_t}-1\right)\notag\\
    &-\e\left[\frac{1}{T}\sum_{t=1}^T\Delta_{t,1}(\bv_t^*(1))\cdot\iprod{\xv_t,\bv_t(1)-\bv^*_t(1)}\cdot\left(\frac{1}{\widebar{p}_t}-1\right)\right]\xrightarrow{p}0\enspace.
\end{align*}

(5) Tracking error terms multiplied by the difference between the random inverse assignment probability and its fixed comparator.
\begin{align*}
    \frac{1}{T}\sum_{t=1}^T\iprod{\xv_t,\bv_t(1)-\bv^*_t(1)}^2\left(\frac{1}{p_t}-\frac{1}{\widebar{p}_t}\right)-\e\left[\frac{1}{T}\sum_{t=1}^T\iprod{\xv_t,\bv_t(1)-\bv^*_t(1)}^2\left(\frac{1}{p_t}-\frac{1}{\widebar{p}_t}\right)\right]\xrightarrow{p}0\enspace.
\end{align*}

(6) Tracking error terms multiplied by the fixed comparator term.
\begin{align*}
    \frac{1}{T}\sum_{t=1}^T\iprod{\xv_t,\bv_t(1)-\bv^*_t(1)}^2\left(\frac{1}{\widebar{p}_t}-1\right)-\e\left[\frac{1}{T}\sum_{t=1}^T\iprod{\xv_t,\bv_t(1)-\bv^*_t(1)}^2\left(\frac{1}{\widebar{p}_t}-1\right)\right]\xrightarrow{p}0\enspace.
\end{align*}

The convergence result in (2) holds trivially, since all terms within it are deterministic.
To establish the remaining convergence results, we simplify the arguments using the following lemma:
\begin{lemma}\label{lemma:markov}
    Let $\{X_n\}$ be a sequence of random variables. If $\e[|X_n|] \to 0$, then $X_n - \e[X_n] \xrightarrow{p} 0$.
\end{lemma}

\begin{proof}
    For any fixed $\epsilon > 0$, Markov's inequality leads to
    \begin{align*}
        \operatorname{Pr}(|X_n| > \epsilon) \le \frac{\e[|X_n|]}{\epsilon} \to 0\enspace,
    \end{align*}
    which implies $|X_n| \xrightarrow{p} 0$, and hence $X_n \xrightarrow{p} 0$.
    Moreover, by Jensen's inequality,
    \begin{align*}
        |\e[X_n]| \le \e[|X_n|] \to 0\enspace.
    \end{align*}
    Therefore, $X_n - \e[X_n] \xrightarrow{p} 0$.
\end{proof}

Given Lemma \ref{lemma:markov}, it suffices to verify the following convergence results:
\begin{align}\label{simplified-stable-3}
    &\e\left[\frac{1}{T}\sum_{t=1}^T\Delta_{t,1}(\bv_t^*(1))^2\cdot\left|\frac{1}{p_t}-\frac{1}{\widebar{p}_t}\right|\right]\rightarrow 0\enspace,\notag\\
    &\e\left[\frac{1}{T}\sum_{t=1}^T|\Delta_{t,1}(\bv_t^*(1))|\cdot|\iprod{\xv_t,\bv_t(1)-\bv^*_t(1)}|\cdot\left|\frac{1}{p_t}-\frac{1}{\widebar{p}_t}\right|\right]\rightarrow 0\enspace,\notag\\
    &\e\left[\frac{1}{T}\sum_{t=1}^T|\Delta_{t,1}(\bv_t^*(1))|\cdot|\iprod{\xv_t,\bv_t(1)-\bv^*_t(1)}|\cdot\left|\frac{1}{\widebar{p}_t}-1\right|\right]\rightarrow 0\enspace,\notag\\
    &\e\left[\frac{1}{T}\sum_{t=1}^T\iprod{\xv_t,\bv_t(1)-\bv^*_t(1)}^2\cdot\left|\frac{1}{p_t}-\frac{1}{\widebar{p}_t}\right|\right]\rightarrow 0\enspace,\notag\\
    &\e\left[\frac{1}{T}\sum_{t=1}^T\iprod{\xv_t,\bv_t(1)-\bv^*_t(1)}^2\cdot\left|\frac{1}{\widebar{p}_t}-1\right|\right]\rightarrow 0\enspace.
\end{align}

Before proving the first convergence result in \eqref{simplified-stable-3}, we establish several supporting lemmas that control the difference between the inverse probabilities associated with $p_t$ and $\widebar{p}_t$.
The following two lemmas bound the difference between the inverse probabilities in terms of the differences between their corresponding first-order equations.
Lemma \ref{lemma:difference-inverse-probability-main} covers the two cases where both probabilities lie on the same side of $1/2$ (both larger than $1/2$ or both smaller than $1/2$).

\begin{lemma}\label{lemma:difference-inverse-probability-main}
    Let $A,B>0$ and $\widetilde{A},\widetilde{B}\geq 0$. 
    Suppose $p,\widetilde{p}$ satisfy
    \begin{align*}
        -\frac{A}{p^2}+\frac{B}{(1-p)^2}+\Psi^{\prime}(p)&=0\enspace,\\
        -\frac{\widetilde{A}}{\widetilde{p}^2}+\frac{\widetilde{B}}{(1-\widetilde{p})^2}+\Psi^{\prime}(\widetilde{p})&=0\enspace.
    \end{align*}
    \begin{enumerate}
        \item[(1)] If $A\geq B,\widetilde{A}\geq \widetilde{B}$, under Condition \mainref{condition:sigmoid} we can obtain the following upper bounds:
        \begin{align*}
            \left|\frac{1}{1-p}-\frac{1}{1-\widetilde{p}}\right|\leq&\frac{\bb_1\bb_2}{\bb_3}\left(1+\frac{2}{\bb_3}\right)\frac{|A-\widetilde{A}|}{A}\frac{p}{1-p}+\frac{2\bb_1^2\bb_2}{\bb_3^2}\left(1+\frac{2}{\bb_3}\right)\frac{|B-\widetilde{B}|}{B}\frac{\widetilde{p}}{1-\widetilde{p}}\enspace,\\
            \left|\frac{1}{p}-\frac{1}{\widetilde{p}}\right|\leq&\frac{\bb_2^2}{2\bb_3}\frac{|A-\widetilde{A}|}{A}+\frac{\bb_1\bb_2}{\bb_3}\left(1+\frac{2}{\bb_3}\right)\frac{|B-\widetilde{B}|}{A}\frac{p}{1-p}\enspace.
        \end{align*}
        \item[(2)] If $A\leq B,\widetilde{A}\leq \widetilde{B}$, under Condition \mainref{condition:sigmoid} we can obtain the following upper bounds:
        \begin{align*}
            \left|\frac{1}{p}-\frac{1}{\widetilde{p}}\right|\leq&\frac{\bb_1\bb_2}{\bb_3}\left(1+\frac{2}{\bb_3}\right)\frac{|B-\widetilde{B}|}{B}\frac{1-p}{p}+\frac{2\bb_1^2\bb_2}{\bb_3^2}\left(1+\frac{2}{\bb_3}\right)\frac{|A-\widetilde{A}|}{A}\frac{1-\widetilde{p}}{\widetilde{p}}\enspace,\\
            \left|\frac{1}{1-p}-\frac{1}{1-\widetilde{p}}\right|\leq&\frac{\bb_2^2}{2\bb_3}\frac{|B-\widetilde{B}|}{B}+\frac{\bb_1\bb_2}{\bb_3}\left(1+\frac{2}{\bb_3}\right)\frac{|A-\widetilde{A}|}{B}\frac{1-p}{p}\enspace.
        \end{align*}
    \end{enumerate}
\end{lemma}

\begin{proof}
    \textbf{Part 1:} We can verify that $p$ and $\widetilde{p}$ are the minimizers of the following optimization programs, respectively.
        \begin{align*}
            p=&\operatorname{argmin}_{p\in(0,1)}\frac{A}{p}+\frac{B}{1-p}+\Psi(p)\enspace,\\
            \widetilde{p}=&\operatorname{argmin}_{p\in(0,1)}\frac{\widetilde{A}}{p}+\frac{\widetilde{B}}{1-p}+\Psi(p)\enspace.
        \end{align*}
        For $u=\phi^{-1}(p)$ and $\widetilde{u}=\phi^{-1}(\widetilde{p})$, one can verify that $u$ and $\widetilde{u}$ are the minimizers of the corresponding transformed optimization problems and therefore satisfy the following first-order conditions:
        \begin{align}\label{lemma:difference-inverse-probability-main_eq1}
            A\left(\frac{1}{\phi(u)}\right)^{\prime}+B\left(\frac{1}{1-\phi(u)}\right)^{\prime}+u+3u|u|=&0\enspace,\notag\\
            \widetilde{A}\left(\frac{1}{\phi(\widetilde{u})}\right)^{\prime}+\widetilde{B}\left(\frac{1}{1-\phi(\widetilde{u})}\right)^{\prime}+\widetilde{u}+3\widetilde{u}|\widetilde{u}|=&0\enspace.
        \end{align}
        Since $A\geq B,\widetilde{A}\geq \widetilde{B}$, by Lemma \ref{lemma:p-location} we have $u,\widetilde{u}\geq 0$. 
        By subtracting the second equation in \eqref{lemma:difference-inverse-probability-main_eq1} from the first equation, we have
        \begin{align}\label{lemma:difference-inverse-probability-main_eq2}
            &A\left[\left(\frac{1}{\phi(u)}\right)^{\prime}-\left(\frac{1}{\phi(\widetilde{u})}\right)^{\prime}\right]+B\left[\left(\frac{1}{1-\phi(u)}\right)^{\prime}-\left(\frac{1}{1-\phi(\widetilde{u})}\right)^{\prime}\right]+(u-\widetilde{u})(1+3u+3\widetilde{u})\notag\\
            &+(A-\widetilde{A})\left(\frac{1}{\phi(\widetilde{u})}\right)^{\prime}+(B-\widetilde{B})\left(\frac{1}{1-\phi(\widetilde{u})}\right)^{\prime}=0\enspace.
        \end{align}
        By the convexity property in Condition \mainref{condition:sigmoid}-2, $A\left[\left(\frac{1}{\phi(u)}\right)^{\prime}-\left(\frac{1}{\phi(\widetilde{u})}\right)^{\prime}\right]$, $B\left[\left(\frac{1}{1-\phi(u)}\right)^{\prime}-\left(\frac{1}{1-\phi(\widetilde{u})}\right)^{\prime}\right]$ and $(u-\widetilde{u})(1+3u+3\widetilde{u})$ always have the same sign as $u-\widetilde{u}$ for $u,\widetilde{u}\geq 0$. 
        Hence by Lemma \ref{lemma:u}-5 and \eqref{lemma:difference-inverse-probability-main_eq2}, we have
        \begin{align}\label{lemma:difference-inverse-probability-main_eq3}
            &\left|(A-\widetilde{A})\left(\frac{1}{\phi(\widetilde{u})}\right)^{\prime}+(B-\widetilde{B})\left(\frac{1}{1-\phi(\widetilde{u})}\right)^{\prime}\right|\notag\\
            =&\left|A\left[\left(\frac{1}{\phi(u)}\right)^{\prime}-\left(\frac{1}{\phi(\widetilde{u})}\right)^{\prime}\right]+B\left[\left(\frac{1}{1-\phi(u)}\right)^{\prime}-\left(\frac{1}{1-\phi(\widetilde{u})}\right)^{\prime}\right]+(u-\widetilde{u})(1+3u+3\widetilde{u})\right|\notag\\
            \geq&A\left|\left(\frac{1}{\phi(u)}\right)^{\prime}-\left(\frac{1}{\phi(\widetilde{u})}\right)^{\prime}\right|\notag\\
            \geq&\frac{\bb_3}{2}\cdot\frac{A}{(1+\widetilde{u})(1+u)(1+\widetilde{u}\land u)}\cdot |u-\widetilde{u}|\enspace.
        \end{align}
        By the monotonicity assumption in Condition \mainref{condition:sigmoid}-1, we have $\left(1/(1-\phi(u))\right)^{\prime}> 0$ and $\left(1/\phi(u)\right)^{\prime}<0$.
        Hence by \eqref{lemma:difference-inverse-probability-main_eq1}, we have
        \begin{align*}
            A=\frac{B\left(\frac{1}{1-\phi(u)}\right)^{\prime}+u+3u^2}{-\left(\frac{1}{\phi(u)}\right)^{\prime}}=\frac{B\left(\frac{1}{1-\phi(u)}\right)^{\prime}+u+3u^2}{\left|\left(\frac{1}{\phi(u)}\right)^{\prime}\right|}\geq B\frac{\left(\frac{1}{1-\phi(u)}\right)^{\prime}}{\left|\left(\frac{1}{\phi(u)}\right)^{\prime}\right|}\enspace,
        \end{align*}
        which implies that
        \begin{align}\label{lemma:difference-inverse-probability-main_eq4}
            \frac{1}{A}\leq \frac{\left|\left(\frac{1}{\phi(u)}\right)^{\prime}\right|}{B\left(\frac{1}{1-\phi(u)}\right)^{\prime}}\leq&\frac{1}{B\left(\frac{1}{1-\phi(0)}\right)^{\prime}}\cdot \frac{\bb_2}{2(1+u)^2}\quad(\text{Lemma \ref{lemma:u}-2, the monotonicity of $\left(1/(1-\phi(u))\right)^{\prime}$})\notag\\
            \leq&\frac{2}{\bb_3B}\cdot \frac{\bb_2}{2(1+u)^2}\quad(\text{Lemma \ref{lemma:u}-1})\notag\\
            =&\frac{\bb_2}{\bb_3B(1+u)^2}\enspace.
        \end{align}
        Hence, by \eqref{lemma:difference-inverse-probability-main_eq3} and \eqref{lemma:difference-inverse-probability-main_eq4}, we have
        \begin{align}\label{lemma:difference-inverse-probability-main_eq5}
            |u-\widetilde{u}|\leq& \frac{2}{\bb_3}\cdot\frac{(1+\widetilde{u})(1+u)(1+\widetilde{u}\land u)}{A}\cdot\left|(A-\widetilde{A})\left(\frac{1}{\phi(\widetilde{u})}\right)^{\prime}+(B-\widetilde{B})\left(\frac{1}{1-\phi(\widetilde{u})}\right)^{\prime}\right|\notag\\
            \leq&\frac{2}{\bb_3}\frac{|A-\widetilde{A}|}{A}(1+\widetilde{u})^2(1+u)\left|\left(\frac{1}{\phi(\widetilde{u})}\right)^{\prime}\right|+\frac{2}{\bb_3}\frac{|B-\widetilde{B}|}{A}(1+\widetilde{u})(1+u)^2\left|\left(\frac{1}{1-\phi(\widetilde{u})}\right)^{\prime}\right|\notag\\
            \intertext{Since $\left|\left(\frac{1}{\phi(\widetilde{u})}\right)^{\prime}\right|\leq \frac{\bb_2}{2(1+\widetilde{u})^2}$ by Lemma \ref{lemma:u}-2, $\left|\left(\frac{1}{1-\phi(\widetilde{u})}\right)^{\prime}\right|=\left|\left(\frac{1}{\phi(-\widetilde{u})}\right)^{\prime}\right|\leq \bb_1$ by Condition \mainref{condition:sigmoid}-3, and $1/A\leq \frac{\bb_2}{\bb_3B(1+u)^2}$ by \eqref{lemma:difference-inverse-probability-main_eq4}, we can further bound this by}
            \leq&\frac{\bb_2}{\bb_3}\frac{|A-\widetilde{A}|}{A}(1+u)+\frac{2\bb_1\bb_2}{\bb_3^2}\frac{|B-\widetilde{B}|}{B}(1+\widetilde{u})\enspace.
        \end{align}
        Hence, by Lemma \ref{lemma:u}-3, we have
        \begin{align}\label{lemma:difference-inverse-probability-main_eq6}
            \left|\frac{1}{1-p}-\frac{1}{1-\widetilde{p}}\right|\leq\bb_1|u-\widetilde{u}|\leq&\frac{\bb_1\bb_2}{\bb_3}\frac{|A-\widetilde{A}|}{A}(1+u)+\frac{2\bb_1^2\bb_2}{\bb_3^2}\frac{|B-\widetilde{B}|}{B}(1+\widetilde{u})\enspace.
        \end{align}
        Moreover, by Lemma \ref{lemma:u}-1 and Lemma \ref{lemma:u}-3, and the convexity property in Condition \mainref{condition:sigmoid}-2, we obtain
        \begin{align}\label{lemma:difference-inverse-probability-main_eq7}
            \frac{1}{1-p}-\frac{1}{1-1/2}=&\frac{1}{1-\phi(u)}-\frac{1}{1-\phi(0)}\geq \left(\frac{1}{1-\phi(0)}\right)^{\prime}(u-0)\geq \frac{\bb_3u}{2}\enspace,\notag\\
            \frac{1}{1-\widetilde{p}}-\frac{1}{1-1/2}\geq &\frac{\bb_3\widetilde{u}}{2}\enspace.
        \end{align}
        We can finally obtain from \eqref{lemma:difference-inverse-probability-main_eq6} and \eqref{lemma:difference-inverse-probability-main_eq7} that
        \begin{align*}
            \left|\frac{1}{1-p}-\frac{1}{1-\widetilde{p}}\right|\leq&\frac{\bb_1\bb_2}{\bb_3}\frac{|A-\widetilde{A}|}{A}(1+u)+\frac{2\bb_1^2\bb_2}{\bb_3^2}\frac{|B-\widetilde{B}|}{B}(1+\widetilde{u})\quad(\text{by \eqref{lemma:difference-inverse-probability-main_eq6}})\\
            \leq&\frac{\bb_1\bb_2}{\bb_3}\frac{|A-\widetilde{A}|}{A}\left[1+\frac{2}{\bb_3}\left(\frac{1}{1-p}-2\right)\right]+\frac{2\bb_1^2\bb_2}{\bb_3^2}\frac{|B-\widetilde{B}|}{B}\left[1+\frac{2}{\bb_3}\left(\frac{1}{1-\widetilde{p}}-2\right)\right]\quad(\text{by \eqref{lemma:difference-inverse-probability-main_eq7}})\\
            =&\frac{\bb_1\bb_2}{\bb_3}\frac{|A-\widetilde{A}|}{A}\left[1+\frac{2}{\bb_3}\left(\frac{p}{1-p}-1\right)\right]+\frac{2\bb_1^2\bb_2}{\bb_3^2}\frac{|B-\widetilde{B}|}{B}\left[1+\frac{2}{\bb_3}\left(\frac{\widetilde{p}}{1-\widetilde{p}}-1\right)\right]\\
            \intertext{Since $p\geq 1/2$ by Lemma \ref{lemma:p-location}, we have $\frac{p}{1-p}\geq 1$ and $1+\frac{2}{\bb_3}\left(\frac{p}{1-p}-1\right)\leq 1+\frac{2}{\bb_3}\cdot\frac{p}{1-p}\leq \left(1+\frac{2}{\bb_3}\right)\frac{p}{1-p}$. Hence,}
            \leq&\frac{\bb_1\bb_2}{\bb_3}\left(1+\frac{2}{\bb_3}\right)\frac{|A-\widetilde{A}|}{A}\frac{p}{1-p}+\frac{2\bb_1^2\bb_2}{\bb_3^2}\left(1+\frac{2}{\bb_3}\right)\frac{|B-\widetilde{B}|}{B}\frac{\widetilde{p}}{1-\widetilde{p}}\enspace.
        \end{align*}
        By the first two steps in \eqref{lemma:difference-inverse-probability-main_eq5}, we can also derive
        \begin{align*}
            &\left|\frac{1}{p}-\frac{1}{\widetilde{p}}\right|\\
            \leq&\frac{\bb_2}{2}\cdot\frac{|u-\widetilde{u}|}{(1+u)(1+\widetilde{u})}\quad(\text{Lemma \ref{lemma:u}-4})\\
            \leq&\frac{\bb_2}{2(1+u)(1+\widetilde{u})}\left(\frac{2}{\bb_3}\frac{|A-\widetilde{A}|}{A}(1+\widetilde{u})^2(1+u)\left|\left(\frac{1}{\phi(\widetilde{u})}\right)^{\prime}\right|+\frac{2}{\bb_3}\frac{|B-\widetilde{B}|}{A}(1+\widetilde{u})(1+u)^2\left|\left(\frac{1}{1-\phi(\widetilde{u})}\right)^{\prime}\right|\right)\\
            \leq&\frac{\bb_2}{2(1+u)(1+\widetilde{u})}\left(\frac{\bb_2}{\bb_3}\frac{|A-\widetilde{A}|}{A}(1+u)+\frac{2}{\bb_3}\frac{|B-\widetilde{B}|}{A}(1+\widetilde{u})(1+u)^2\left|\left(\frac{1}{1-\phi(\widetilde{u})}\right)^{\prime}\right|\right)\quad(\text{Lemma \ref{lemma:u}-2})\\
            \intertext{Since $\left|\left(\frac{1}{1-\phi(\widetilde{u})}\right)^{\prime}\right|=\left|\left(\frac{1}{\phi(-\widetilde{u})}\right)^{\prime}\right|\leq \bb_1$ by Condition \mainref{condition:sigmoid}-3, this can be further bounded as}
            \leq&\frac{\bb_2^2}{2\bb_3}\frac{|A-\widetilde{A}|}{A}+\frac{\bb_1\bb_2}{\bb_3}\frac{|B-\widetilde{B}|}{A}(1+u)\\
            \leq&\frac{\bb_2^2}{2\bb_3}\frac{|A-\widetilde{A}|}{A}+\frac{\bb_1\bb_2}{\bb_3}\left(1+\frac{2}{\bb_3}\right)\frac{|B-\widetilde{B}|}{A}\frac{p}{1-p}\quad(\text{by \eqref{lemma:difference-inverse-probability-main_eq7}})\enspace,
        \end{align*}
        which completes the proof in part 1.\\[3mm]
    \textbf{Part 2:}  Let $q=1-p$ and $\widetilde{q}=1-\widetilde{p}$. 
    Then $q$ and $\widetilde{q}$ satisfy:
    \begin{align*}
        -\frac{A}{(1-q)^2}+\frac{B}{q^2}+\Psi^{\prime}(1-q)&=0\enspace,\\
        -\frac{\widetilde{A}}{(1-\widetilde{q})^2}+\frac{\widetilde{B}}{\widetilde{q}^2}+\Psi^{\prime}(1-\widetilde{q})&=0\enspace.
    \end{align*}
    Since by Condition \mainref{condition:sigmoid}-1, we have $\phi^{-1}(p)=-\phi^{-1}(1-p)$, which further indicates that $\Psi(p)=\Psi(1-p)$ since $\psi$ is an even function. 
    This implies that $\Psi^{\prime}(1-q)=-\Psi^{\prime}(q)$ and $\Psi^{\prime}(1-\widetilde{q})=-\Psi^{\prime}(\widetilde{q})$, i.e.,
    \begin{align*}
        -\frac{B}{q^2}+\frac{A}{(1-q)^2}+\Psi^{\prime}(q)=&0\enspace,\\
        -\frac{\widetilde{B}}{\widetilde{q}^2}+\frac{\widetilde{A}}{(1-\widetilde{q})^2}+\Psi^{\prime}(\widetilde{q})=&0\enspace.
    \end{align*}
    Hence by the claim in part 1, the result is verified.\qedhere
\end{proof}

When $p_t$ and $\widebar{p}_t$ do not lie on the same side of 1/2, it is not possible to directly bound the difference between their inverse probabilities.
However, when the two first-order equations are close, both probabilities should be close to $1/2$.
Therefore, we use $1/2$ as an intermediate comparator and bound their difference in the following lemma.

\begin{lemma}\label{lemma:difference-inverse-probability-half}
    Let $A, B$ be positive numbers. 
    Suppose $p$ satisfies:
    \begin{align*}
        -\frac{A}{p^2}+\frac{B}{(1-p)^2}+\Psi^{\prime}(p)&=0\enspace.
    \end{align*}
    \begin{enumerate}
        \item[(1)] If $B\leq A$, then under Condition \mainref{condition:sigmoid}, $p$ satisfies $0\leq\frac{1}{1-p}-2\leq \frac{\bb_1\bb_2}{\bb_3}\cdot\frac{A-B}{B}$.
        \item[(2)] If $A\leq B$, then under Condition \mainref{condition:sigmoid}, $p$ satisfies $0\leq \frac{1}{p}-2\leq \frac{\bb_1\bb_2}{\bb_3}\cdot\frac{B-A}{A}$.
    \end{enumerate}
\end{lemma}
\begin{proof}
    \textbf{Part 1:} When $A\geq B$, by Lemma \ref{lemma:p-location} we have $p\geq1/2$ and $u=\phi^{-1}(p)\geq 0$. 
        Let $\widetilde{u}=\phi^{-1}(1/2)=0$.
        By the symmetry property in Condition \mainref{condition:sigmoid}, we have
        \begin{align*}
            \left(\frac{1}{\phi(\widetilde{u})}\right)^{\prime}=-\left(\frac{1}{1-\phi(-\widetilde{u})}\right)^{\prime}=-\left(\frac{1}{1-\phi(\widetilde{u})}\right)^{\prime}\quad(\text{since $\widetilde{u}=0$})\enspace.
        \end{align*}
        Hence, $u$ and $\widetilde{u}$ satisfy the following equations:
        \begin{align}\label{lemma:difference-inverse-probability-half_eq1}
            A\left(\frac{1}{\phi(u)}\right)^{\prime}+B\left(\frac{1}{1-\phi(u)}\right)^{\prime}+u+3u^2=&0\enspace,\notag\\
            B\left(\frac{1}{\phi(\widetilde{u})}\right)^{\prime}+B\left(\frac{1}{1-\phi(\widetilde{u})}\right)^{\prime}+\widetilde{u}+3\widetilde{u}^2=&0\enspace.
        \end{align}
        Subtracting the second equation in \eqref{lemma:difference-inverse-probability-half_eq1} from the first equation, we can obtain
        \begin{align}\label{lemma:difference-inverse-probability-half_eq2}
            &B\left[\left(\frac{1}{\phi(u)}\right)^{\prime}-\left(\frac{1}{\phi(\widetilde{u})}\right)^{\prime}\right]+B\left[\left(\frac{1}{1-\phi(u)}\right)^{\prime}-\left(\frac{1}{1-\phi(\widetilde{u})}\right)^{\prime}\right]+(u-\widetilde{u})(1+3u+3\widetilde{u})\notag\\
            &+(A-B)\left(\frac{1}{\phi(u)}\right)^{\prime}=0\enspace.
        \end{align}
        The convexity property in Condition \mainref{condition:sigmoid}-2 guarantees that $B\left[\left(\frac{1}{\phi(u)}\right)^{\prime}-\left(\frac{1}{\phi(\widetilde{u})}\right)^{\prime}\right]$, $B\left[\left(\frac{1}{1-\phi(u)}\right)^{\prime}-\left(\frac{1}{1-\phi(\widetilde{u})}\right)^{\prime}\right]$ and $(u-\widetilde{u})(1+3u+3\widetilde{u})$ have the same sign as $u-\widetilde{u}$, which is nonnegative.
        Then by Lemma \ref{lemma:u}-5 and \eqref{lemma:difference-inverse-probability-half_eq2}, we have the following inequality:
        \begin{align}\label{lemma:difference-inverse-probability-half_eq3}
            &(A-B)\left|\left(\frac{1}{\phi(u)}\right)^{\prime}\right|\notag\\
            =&\left|B\left[\left(\frac{1}{\phi(u)}\right)^{\prime}-\left(\frac{1}{\phi(\widetilde{u})}\right)^{\prime}\right]+B\left[\left(\frac{1}{1-\phi(u)}\right)^{\prime}-\left(\frac{1}{1-\phi(\widetilde{u})}\right)^{\prime}\right]+(u-\widetilde{u})(1+3u+3\widetilde{u})\right|\notag\\
            \geq&B\left|\left(\frac{1}{\phi(u)}\right)^{\prime}-\left(\frac{1}{\phi(\widetilde{u})}\right)^{\prime}\right|\notag\\
            \geq&\frac{\bb_3}{2}\cdot\frac{B}{(1+\widetilde{u})(1+u)(1+\widetilde{u}\land u)}\cdot (u-\widetilde{u})\quad(\text{Lemma \ref{lemma:u}-5})\enspace.
        \end{align}
        Hence by Lemma \ref{lemma:u} and \eqref{lemma:difference-inverse-probability-half_eq3}, we obtain
        \begin{align*}
            \frac{1}{1-p}-\frac{1}{1-1/2}\leq& \bb_1(u-\widetilde{u})\quad(\text{Lemma \ref{lemma:u}-3})\\
            \leq&\bb_1\cdot\frac{2}{\bb_3}\cdot\frac{(1+\widetilde{u})(1+u)(1+\widetilde{u}\land u)}{B}\cdot(A-B)\left|\left(\frac{1}{\phi(u)}\right)^{\prime}\right|\quad(\text{by \eqref{lemma:difference-inverse-probability-half_eq2} and \eqref{lemma:difference-inverse-probability-half_eq3}})\\
            \leq&\frac{\bb_1\bb_2}{\bb_3}\cdot\frac{A-B}{B}(1+u)^2(1+\widetilde{u})(1+u)^{-2}\quad(\text{Lemma \ref{lemma:u}-2})\\
            =&\frac{\bb_1\bb_2}{\bb_3}\cdot\frac{A-B}{B}\quad(\text{since $\widetilde{u}=0$})\enspace.
        \end{align*}
        \textbf{Part 2:}  Let $q=1-p$. 
        Following the same arguments as in part 2 of Lemma \ref{lemma:difference-inverse-probability-main}, $q$ satisfies
        \begin{align*}
            -\frac{B}{q^2}+\frac{A}{(1-q)^2}+\Psi^{\prime}(q)&=0\enspace.
        \end{align*}
        By the claim in part 1, we can prove that
        \begin{align*}
            \frac{1}{p}-2=\frac{1}{1-q}-2\leq \frac{\bb_1\bb_2}{\bb_3}\cdot \frac{B-A}{A}\enspace.
        \end{align*}
        This completes the proof. \qedhere
\end{proof}

Lemma \ref{lemma:difference-inverse-probability-main} and Lemma \ref{lemma:difference-inverse-probability-half} imply that controlling the differences between inverse probabilities reduces to bounding the differences between the corresponding first-order equations.
By the definitions of $p_t$ and $\widebar{p}_t$, this further reduces to bounding the difference between the estimated squared residuals and their expectations, which can in turn be achieved by controlling the variance of the estimated squared residuals.
The following lemma establishes such variance bounds.
For simplicity, we denote $\widehat{a}_{t,k}=\eta_t\widehat{A}_{t}(k)$ and $a_{t,k}=\eta_t\e[\widehat{A}_{t}(k)]$ for $k\in\{0,1\}$ and any $t\in[T]$.

\begin{lemma}\label{lemma:variance-estimated-squared-residual}
    Suppose $T$ is sufficiently large. Under Assumptions \mainref{assumption:moments}-\mainref{assumption:maximum-radius-strengthened}, and Condition \mainref{condition:sigmoid}, for any $t\in[T]$, the variance of $\widehat{a}_{t,1}$ and $\widehat{a}_{t,0}$ can be bounded by:
    \begin{align*}
        \Var{\widehat{a}_{t,1}}\leq&C(\eta_tT)^{1/2}\log^2(\eta_t T)(2+\bb_1(\bb_2/6)^{1/4}a_{t,0}^{1/4})^{25/16}\enspace,\\
        \Var{\widehat{a}_{t,0}}\leq&C(\eta_tT)^{1/2}\log^2(\eta_t T)(2+\bb_1(\bb_2/6)^{1/4}a_{t,1}^{1/4})^{25/16}\enspace,
    \end{align*}
    where $C:=c\,c_1^2\gamma\left(8\xi_2^{1/2}\zeta_1+46(\zeta_1^{1/2}+c_1)^2\right)>0$ is a constant.
\end{lemma}

\begin{proof}
    Without loss of generality, we only verify the first inequality. 
    Throughout the proof, we omit irrelevant constants in small order terms for clarity.
    For $k\in\{0,1\}$, denote $\Delta_{t,k}(\bv)=y_t(k)-\iprod{\xv_t,\bv}$.
    Further denote $\widehat{\Delta}_{t,1}(\bv)=y_t(1)\cdot \frac{\indicator{Z_t=1}}{p_t}-\iprod{\xv_t,\bv}$ and $\widehat{\Delta}_{t,0}(\bv)=y_t(0)\cdot \frac{\indicator{Z_t=0}}{1-p_t}-\iprod{\xv_t,\bv}$.
    For any $t\in[T]$, denote $e_t:=\left(\frac{\mathbf{1}[Z_{t}=1]}{p_{t}}-1\right)\cdot\Delta_{t,1}(\bv_{t}(1))^2$.
    Then the variance of $\widehat{A}_t(1)$ can be bounded by:
    \begin{align}\label{lemma:variance-estimated-squared-residual_eq1}
        \Var{\widehat{A}_t(1)}=&\operatorname{Var}\Big[\sum_{s=1}^{t}\frac{\mathbf{1}[Z_s=1]}{p_s}\cdot\Delta_{s,1}(\bv_{s}(1))^2\Big]\notag\\
        =&\operatorname{Var}\left[\sum_{s=1}^{t}e_s+\sum_{s=1}^{t}\Delta_{s,1}(\bv_{s}(1))^2\right]\notag\\
        \intertext{Using the Cauchy-Schwarz inequality $|\operatorname{Cov}(X,Y)|\leq \operatorname{Var}(X)^{1/2}\operatorname{Var}(Y)^{1/2}$, we can bound this as}
        \leq&\underbrace{\operatorname{Var}\left[\sum_{s=1}^{t}e_s\right]}_{:=S_1}+\underbrace{\operatorname{Var}\left[\sum_{s=1}^{t}\Delta_{s,1}(\bv_{s}(1))^2\right]}_{:=S_2}\notag\\
        &+2\underbrace{\left(\operatorname{Var}\left[\sum_{s=1}^{t}e_s\right]\right)^{1/2}\left(\operatorname{Var}\left[\sum_{s=1}^{t}\Delta_{s,1}(\bv_{s}(1))^2\right]\right)^{1/2}}_{=S_1^{1/2}S_2^{1/2}}\enspace.
    \end{align}
    We then derive upper bounds on $S_1$ and $S_2$ separately.\\[3mm]
    \textbf{Step 1: }Derive an upper bound on $S_1$.\\[3mm]
    For term $S_1$, by expanding the variance term, we obtain
    \begin{align}\label{lemma:variance-estimated-squared-residual_eq2}
        S_1=\operatorname{Var}\left[\sum_{s=1}^{t}e_s\right]=\sum_{s=1}^{t}\operatorname{Var}\left[e_s\right]+2\sum_{1\leq s_1<s_2\leq t}\operatorname{Cov}(e_{s_1},e_{s_2})\enspace.
    \end{align}
    Using the law of total variance, the first term in \eqref{lemma:variance-estimated-squared-residual_eq2} can be simplified as
    \begin{align}\label{lemma:variance-estimated-squared-residual_eq3}
        \sum_{s=1}^{t}\operatorname{Var}\left[e_s\right]=&\sum_{s=1}^{t}\operatorname{Var}\left[\e\left[e_s|\mathcal{F}_{s-1}\right]\right]+\sum_{s=1}^{t}\e\left[\operatorname{Var}\left[e_s|\mathcal{F}_{s-1}\right]\right]\notag\\
        \intertext{Since $\e[e_s|\filt_{s-1}]=\Delta_{s,1}(\bv_{s}(1))^2\cdot\e[\frac{\mathbf{1}[Z_s=1]}{p_s}-1|\filt_{s-1}]=0$, this can be further simplified as}
        =&\sum_{s=1}^{t}\e\left[\operatorname{Var}\left[\Delta_{s,1}(\bv_{s}(1))^2\cdot\left(\frac{\mathbf{1}[Z_s=1]}{p_s}-1\right)\Bigg|\mathcal{F}_{s-1}\right]\right]\notag\\
        =&\sum_{s=1}^{t}\e\left[\Delta_{s,1}(\bv_{s}(1))^4\cdot\operatorname{Var}\left[\left(\frac{\mathbf{1}[Z_s=1]}{p_s}-1\right)\Bigg|\mathcal{F}_{s-1}\right]\right]\notag\\
        =&\sum_{s=1}^{t}\e\left[\frac{1-p_s}{p_s}\cdot\Delta_{s,1}(\bv_{s}(1))^4\right]\enspace.
    \end{align}
    For $s_1<s_2$, we can show that $e_{s_1}$ is $\filt_{s_2-1}$ measurable and $\e[e_{s_2}|\filt_{s_2-1}]=0$.
    Hence by using the law of total covariance, the second term in \eqref{lemma:variance-estimated-squared-residual_eq2} can be simplified as
    \begin{align*}
        2\sum_{1\leq s_1<s_2\leq t}\operatorname{Cov}(e_{s_1},e_{s_2})=2\sum_{1\leq s_1<s_2\leq t}(\operatorname{Cov}(\e[e_{s_1}|\mathcal{F}_{s_2-1}],\e[e_{s_2}|\mathcal{F}_{s_2-1}])+\e[\operatorname{Cov}(e_{s_1},e_{s_2}|\mathcal{F}_{s_2-1})])=0\enspace.
    \end{align*}
    Hence by using inequality $(a+b)^4\leq 8a^4+8b^4$, \eqref{lemma:variance-estimated-squared-residual_eq2} and \eqref{lemma:variance-estimated-squared-residual_eq3}, $S_1$ can be simplified and bounded as
    \begin{align}\label{lemma:variance-estimated-squared-residual_eq4}
        S_1=&\sum_{s=1}^{t}\e\left[\frac{1-p_s}{p_s}\cdot\Delta_{s,1}(\bv_{s}(1))^4\right]\notag\\
        \leq&\sum_{s=1}^{t}\e\left[\frac{1}{p_s}\cdot\Delta_{s,1}(\bv_{s}(1))^4\right]\notag\\
        =&\sum_{s=1}^{t}\e\left[\frac{1}{p_s}\cdot\left\{\Delta_{s,1}(\bv_{s}^*(1))+\iprod{\xv_s,\bv_s^*(1)-\bv_s(1)}\right\}^4\right]\notag\\
        \leq&8\underbrace{\sum_{s=1}^{t}\Delta_{s,1}(\bv_{s}^*(1))^4\cdot\e\left[\frac{1}{p_s}\right]}_{:=S_{1,1}}+8\underbrace{\sum_{s=1}^{t}\e\left[\frac{1}{p_s}\cdot\iprod{\xv_s,\bv_s(1)-\bv^*_s(1)}^4\right]}_{:=S_{1,2}}\enspace.
    \end{align}
    We then bound $S_{1,1}$ and $S_{1,2}$, corresponding to the full-information online residuals and tracking error terms, respectively.
    By Corollary \ref{corollary:partial-fourth-moment} and Lemma \ref{lemma:p-power-moment}, we separately control the deterministic residual term and the random inverse probability term in $S_{1,1}$, leading to:
    \begin{align}\label{lemma:variance-estimated-squared-residual_eq5}
        S_{1,1}=&\sum_{s=1}^{t}(y_s(1)-\iprod{\xv_s,\bv_s^*(1)})^4\cdot\underbrace{\e\left[\frac{1}{p_s}\right]}_{\text{Lemma \ref{lemma:p-power-moment}}}\notag\\
        \lesssim&T^{1/8}\underbrace{\sum_{s=1}^{t}R_s^{-1/4}(y_s(1)-\iprod{\xv_s,\bv_s^*(1)})^4}_{\text{Corollary \ref{corollary:partial-fourth-moment}}}\notag\\
        \lesssim&T^{1/8}\cdot TR_t^{2-1/4}\notag\\
        =&T^{5/4}R_t^{3/2}\cdot (TR_t^{-2})^{-1/8}\quad(\text{rewrite in the target scale})\notag\\
        =&o(T^{5/4}R_t^{3/2}\log^2(\eta_tT))\quad(\text{Assumption \mainref{assumption:maximum-radius-strengthened} and Lemma \ref{lemma:R}})\enspace.
    \end{align}
    By Proposition \mainref{prop:as-inv-prob-bound} and Corollary \ref{corollary:expected-tracking}, $S_{1,2}$ can be bounded as:
    \begin{align}\label{lemma:variance-estimated-squared-residual_eq6}
        S_{1,2}=&\sum_{s=1}^{t}\e\left[\frac{1}{p_s}\iprod{\xv_s,\bv_s(1)-\bv^*_s(1)}^4\right]\notag\\
        \lesssim&T^{7/26}\sum_{s=1}^{t}R_s^{9/11}\cdot R_s^{-1}\e\left[\iprod{\xv_s,\bv_s(1)-\bv^*_s(1)}^4\right]~~~~(\text{Proposition \mainref{prop:as-inv-prob-bound}})\notag\\
        \lesssim&T^{7/26}R_t^{9/11}\cdot T^{3/8}R_t^{1/4}~~~~(\text{using the monotonicity $R_s \leq R_t$ and Corollary \ref{corollary:expected-tracking}})\notag\\
        =&T^{5/4}R_t^{3/2}\cdot T^{-63/104}R_t^{-19/44}\quad(\text{rewrite in the target scale})\notag\\
        =&o\,(T^{5/4}R_t^{3/2}\log^2(\eta_tT))\quad(\text{Assumption \mainref{assumption:maximum-radius-strengthened} and Lemma \ref{lemma:R}})\enspace.
    \end{align}
    Hence by \eqref{lemma:variance-estimated-squared-residual_eq4}, \eqref{lemma:variance-estimated-squared-residual_eq5} and \eqref{lemma:variance-estimated-squared-residual_eq6}, $S_1$ can be bounded as
    \begin{align}\label{lemma:variance-estimated-squared-residual_eq7}
        S_1\leq& 8S_{1,1}+8S_{1,2}=o\,(T^{5/4}R_t^{3/2}\log^2(\eta_tT))\enspace.
    \end{align}
    \textbf{Step 2: }Derive an upper bound on $S_2$.\\[3mm]
    For $S_2$, we have
    \begin{align}\label{lemma:variance-estimated-squared-residual_eq8}
        S_2=&\operatorname{Var}\left[\sum_{s=1}^{t}\Delta_{s,1}(\bv_{s}(1))^2\right]\notag\\
        =&\operatorname{Var}\left[\sum_{s=1}^{t}\left(\Delta_{s,1}(\bv_{s}(1))^2-\Delta_{s,1}(\bv_{s}^*(1))^2\right)\right]\quad(\text{subtract nonrandom terms})\notag\\
        \leq&\e\left[\left(\sum_{s=1}^{t}\left(\Delta_{s,1}(\bv_{s}(1))^2-\Delta_{s,1}(\bv_{s}^*(1))^2\right)\right)^2\right]\quad(\text{using $\operatorname{Var}(X)\leq \e [X^2]$})\notag\\
        \intertext{Using equality $a^2 - b^2 = (a+b)(a-b)$ and equality $\Delta_{s,1}(\bv_s(1))-\Delta_{s,1}(\bv_s^*(1))=\iprod{\xv_s,\bv_s^*(1)-\bv_s(1)}$, this can be simplified as}
        =&\e\left[\left(\sum_{s=1}^{t}\left(\Delta_{s,1}(\bv_s(1))+\Delta_{s,1}(\bv_s^*(1))\right)\cdot\iprod{\xv_s,\bv_s^*(1)-\bv_s(1)}\right)^2\right]\notag\\
        =&\e\left[\left(\sum_{s=1}^{t}\left(2\Delta_{s,1}(\bv_s^*(1))+\iprod{\xv_s,\bv_s^*(1)-\bv_s(1)}\right)\cdot\iprod{\xv_s,\bv_s^*(1)-\bv_s(1)}\right)^2\right]\notag\\
        =&\e\Bigg[\Bigg(2\sum_{s=1}^{t}\Delta_{s,1}(\bv_s^*(1))\cdot\iprod{\xv_s,\bv_s^*(1)-\bv_s(1)}+\sum_{s=1}^{t}\iprod{\xv_s,\bv_s^*(1)-\bv_s(1)}^2\Bigg)^2\Bigg]\notag\\
        \intertext{Using the inequality $\e[(X+Y)^2]\leq 2\e [X^2]+2\e [Y^2]$, we have}
        \leq&8\underbrace{\e\left[\left(\sum_{s=1}^{t}\Delta_{s,1}(\bv_s^*(1))\cdot\iprod{\xv_s,\bv_s(1)-\bv^*_s(1)}\right)^2\right]}_{:=S_{2,1}}\notag\\
        &+2\underbrace{\e\left[\left(\sum_{s=1}^{t}\iprod{\xv_s,\bv_s(1)-\bv^*_s(1)}^2\right)^2\right]}_{:=S_{2,2}}\enspace.
    \end{align}
    We then derive upper bounds on $S_{2,1}$ and $S_{2,2}$ separately.\\[3mm]
    \textbf{Step 2.1: }Derive an upper bound on $S_{2,1}$.\\[3mm]
    Since $\iprod{\xv_s,\bv_s(1)-\bv^*_s(1)}$ has mean zero by Lemma \mainref{lemma:predictor-expectation}, by plugging in the explicit form of $\iprod{\xv_s,\bv_s(1)-\bv^*_s(1)}$, we have
    \begin{align}\label{lemma:variance-estimated-squared-residual_eq9}
        \quad S_{2,1}=&\e\left[\left(\sum_{s=1}^{t}\Delta_{s,1}(\bv_s^*(1))\cdot\iprod{\xv_s,\bv_s(1)-\bv^*_s(1)}\right)^2\right]\notag\\
        =&\operatorname{Var}\left[\sum_{s=1}^{t}\Delta_{s,1}(\bv_s^*(1))\cdot\iprod{\xv_s,\bv_s(1)-\bv^*_s(1)}\right]\notag\\
        =&\operatorname{Var}\left[\sum_{s=1}^{t}\Delta_{s,1}(\bv_s^*(1))\cdot\sum_{r=1}^{s-1}\Pi_{s,r}y_r(1)\left(\frac{\mathbf{1}[Z_{r}=1]}{p_{r}}-1\right)\right]\notag\\
        \intertext{Swapping the summation order leads to}
        =&\operatorname{Var}\left[\sum_{r=1}^{t-1}\left(\sum_{s=r+1}^t\Pi_{s,r}\cdot\Delta_{s,1}(\bv_s^*(1))\right)y_r(1)\left(\frac{\mathbf{1}[Z_{r}=1]}{p_{r}}-1\right)\right]\notag\\
        \intertext{Since the terms $\left(\sum_{s=r+1}^t\Pi_{s,r}\cdot\Delta_{s,1}(\bv_s^*(1))\right)y_r(1)\left(\frac{\mathbf{1}[Z_{r}=1]}{p_{r}}-1\right)$ form an $L^2$ martingale difference sequence with respect to $\filt_r$, the variance of the sum equals the sum of variances.}
        =&\sum_{r=1}^{t-1}\operatorname{Var}\left[\left(\sum_{s=r+1}^t\Pi_{s,r}\cdot\Delta_{s,1}(\bv_s^*(1))\right)y_r(1)\left(\frac{\mathbf{1}[Z_{r}=1]}{p_{r}}-1\right)\right]\notag\\
        =&\sum_{r=1}^{t-1}\left(\sum_{s=r+1}^t\Pi_{s,r}\cdot\Delta_{s,1}(\bv_s^*(1))\right)^2y_r(1)^2\operatorname{Var}\left(\frac{\mathbf{1}[Z_{r}=1]}{p_{r}}-1\right)\notag\\
        \intertext{Applying the law of total variance gives $\operatorname{Var}\left(\frac{\mathbf{1}[Z_{r}=1]}{p_{r}}-1\right)=\operatorname{Var}\left(\e\left[\frac{\mathbf{1}[Z_{r}=1]}{p_{r}}-1|\filt_{r-1}\right]\right)+\e\left[\operatorname{Var}\left(\frac{\mathbf{1}[Z_{r}=1]}{p_{r}}-1|\filt_{r-1}\right)\right]=\e\left[\frac{1}{p_r}-1\right]$. Then we can further simplify this expression as}
        =&\sum_{r=1}^{t-1}\left(\sum_{s=r+1}^t\Pi_{s,r}\cdot\Delta_{s,1}(\bv_s^*(1))\right)^2y_r(1)^2\underbrace{\e\left[\frac{1}{p_r}-1\right]}_{\text{Corollary \mainref{corollary:p-moment}*}}\notag\\
        \leq&\sum_{r=1}^{t-1}\left(\sum_{s=r+1}^t\Pi_{s,r}\cdot\Delta_{s,1}(\bv_s^*(1))\right)^2y_r(1)^2\left(2+\bb_1(\bb_2/6)^{1/4}\eta_r^{1/4}\e[\widehat{A}_{r-1}(0)]^{1/4}\right)\notag\\
        \intertext{Since $r \leq t$ and $\widehat A_s(0)$ is nondecreasing in $s$, we can bound $\e[\widehat{A}_{r-1}(0)]\leq \e[\widehat{A}_{t}(0)]$. Using the fact that $R_r\leq R_t$ and $\eta_r/\eta_t=R_t/R_r\geq 1$, we can further bound this by}
        \leq&\sum_{r=1}^{t-1}\left(\sum_{s=r+1}^t\Pi_{s,r}\cdot\Delta_{s,1}(\bv_s^*(1))\right)^2y_r(1)^2\left(\frac{R_t}{R_r}\right)^{1/4}\left(2+\bb_1(\bb_2/6)^{1/4}\eta_t^{1/4}\e[\widehat{A}_t(0)]^{1/4}\right)\notag\\
        \leq&R_t^{1/4}\underbrace{\sum_{r=1}^{t-1}\left(\sum_{s=r+1}^t\Pi_{s,r}\cdot\Delta_{s,1}(\bv_s^*(1))\right)^2R_r^{-1/4}y_r(1)^2}_{\text{bound by Lemma \ref{lemma:deterministic-summation-new}-2}}\left(2+\bb_1(\bb_2/6)^{1/4}a_{t,0}^{1/4}\right)\notag\\
        \leq&c\,c_1^2\gamma\xi_2^{1/2}\zeta_1T^{5/4}R_t^{5/4+1/4}\left(2+\bb_1(\bb_2/6)^{1/4}a_{t,0}^{1/4}\right)\notag\\
        =&c\,c_1^2\gamma\xi_2^{1/2}\zeta_1\cdot T^{5/4}R_t^{3/2}\left(2+\bb_1(\bb_2/6)^{1/4}a_{t,0}^{1/4}\right)\enspace.
    \end{align}
    \textbf{Step 2.2: }Derive an upper bound on $S_{2,2}$.\\[3mm]
    The overall strategy for bounding $S_{2,2}$ is to expand the square, separate diagonal and cross terms, and control each part separately.
    \begin{align}\label{lemma:variance-estimated-squared-residual_eq10}
        &S_{2,2}\notag\\
        =&\e\left[\left(\sum_{s=1}^{t}\iprod{\xv_s,\bv_s(1)-\bv^*_s(1)}^2\right)^2\right]\notag\\
        =&\e\left[\left(\sum_{s=1}^t\left(\sum_{r=1}^{s-1}\Pi_{s,r}y_r(1)\left[\frac{\mathbf{1}[Z_r=1]}{p_r}-1\right]\right)^2\right)^2\right]\quad(\text{Lemma \mainref{lemma:ridge-pred-exact-form}*})\notag\\
        =&\e\left[\left(\sum_{s=1}^{t}\sum_{1\leq t_1,t_2\leq s-1}\Pi_{s,t_1}\Pi_{s,t_2}y_{t_1}(1)y_{t_2}(1)\left(\frac{\mathbf{1}[Z_{t_1}=1]}{p_{t_1}}-1\right)\left(\frac{\mathbf{1}[Z_{t_2}=1]}{p_{t_2}}-1\right)\right)^2\right]\quad(\text{expand the squares})\notag\\
        \intertext{Separating diagonal ($t_1=t_2$) and cross ($t_1\ne t_2$) terms gives}
        =&\e\Bigg[\Bigg(\underbrace{\sum_{s=1}^{t}\sum_{1\leq r\leq s-1}\Pi_{s,r}^2y_r(1)^2\left(\frac{\mathbf{1}[Z_{r}=1]}{p_{r}}-1\right)^2}_{\text{diagonal terms}}\notag\\
        &+\underbrace{\sum_{s=1}^{t}\sum_{1\leq t_1\neq t_2\leq s-1}\Pi_{s,t_1}\Pi_{s,t_2}y_{t_1}(1)y_{t_2}(1)\left(\frac{\mathbf{1}[Z_{t_1}=1]}{p_{t_1}}-1\right)\left(\frac{\mathbf{1}[Z_{t_2}=1]}{p_{t_2}}-1\right)}_{\text{cross terms}}\Bigg)^2\Bigg]\notag\\
        \intertext{Using inequality $\e[(X+Y)^2]\leq 2\e [X^2]+2\e [Y^2]$, this can further be bounded as}
        \leq&2\e\left[\left(\sum_{s=1}^{t}\sum_{r=1}^{s-1}\Pi_{s,r}^2y_r(1)^2\left(\frac{\mathbf{1}[Z_{r}=1]}{p_{r}}-1\right)^2\right)^2\right]\notag\\
        &+2\e\left[\left(\sum_{s=1}^{t}\sum_{1\leq t_1\neq t_2\leq s-1}\Pi_{s,t_1}\Pi_{s,t_2}y_{t_1}(1)y_{t_2}(1)\left(\frac{\mathbf{1}[Z_{t_1}=1]}{p_{t_1}}-1\right)\left(\frac{\mathbf{1}[Z_{t_2}=1]}{p_{t_2}}-1\right)\right)^2\right]\notag\\
        \intertext{Swapping the summation orders in both terms leads to}
        =&2\e\left[\left(\sum_{r=1}^{t-1}\left[\sum_{s=r+1}^{t}\Pi_{s,r}^2\right]y_r(1)^2\left(\frac{\mathbf{1}[Z_{r}=1]}{p_{r}}-1\right)^2\right)^2\right]\notag\\
        &+2\operatorname{Var}\left[\sum_{1\leq t_1\neq t_2\leq t-1}\left(\sum_{s=t_1\vee t_2+1}^{t}\Pi_{s,t_1}\Pi_{s,t_2}\right)y_{t_1}(1)y_{t_2}(1)\left(\frac{\mathbf{1}[Z_{t_1}=1]}{p_{t_1}}-1\right)\left(\frac{\mathbf{1}[Z_{t_2}=1]}{p_{t_2}}-1\right)\right]\notag\\
        \intertext{By expanding the squares in the first term and grouping diagonal/cross terms, this expression can be further decomposed as}
        =&2\underbrace{\sum_{r=1}^{t-1}\left(\sum_{s=r+1}^{t}\Pi_{s,r}^2\right)^2y_r(1)^4\e\left[\left(\frac{\mathbf{1}[Z_{r}=1]}{p_{r}}-1\right)^4\right]}_{:=B_1~(\text{diagonal terms})}\notag\\
        &+2\underbrace{\sum_{1\leq t_1\neq t_2\leq t-1}\left(\sum_{s=t_1+1}^{t}\Pi_{s,t_1}^2\right)\left(\sum_{s=t_2+1}^{t}\Pi_{s,t_2}^2\right)y_{t_1}(1)^2y_{t_2}(1)^2\e\left[\left(\frac{\mathbf{1}[Z_{t_1}=1]}{p_{t_1}}-1\right)^2\left(\frac{\mathbf{1}[Z_{t_2}=1]}{p_{t_2}}-1\right)^2\right]}_{:=B_2~(\text{cross terms})}\notag\\
        &+2\underbrace{\operatorname{Var}\left[\sum_{1\leq t_1\neq t_2\leq t-1}\left(\sum_{s=t_1\vee t_2+1}^{t}\Pi_{s,t_1}\Pi_{s,t_2}\right)y_{t_1}(1)y_{t_2}(1)\left(\frac{\mathbf{1}[Z_{t_1}=1]}{p_{t_1}}-1\right)\left(\frac{\mathbf{1}[Z_{t_2}=1]}{p_{t_2}}-1\right)\right]}_{:=B_3}\enspace.
    \end{align}
    We then bound $B_1$, $B_2$ and $B_3$ separately.
    By Lemma \ref{lemma:deterministic-summation-new} and Lemma \ref{lemma:p-moments-cross-term}, we can bound the random and nonrandom parts in $B_1$ separately, leading to:
    \begin{align}\label{lemma:variance-estimated-squared-residual_eq11}
        B_1=&\sum_{r=1}^{t-1}\left(\sum_{s=r+1}^{t}\Pi_{s,r}^2\right)^2y_r(1)^4\underbrace{\e\left[\left(\frac{\mathbf{1}[Z_{r}=1]}{p_{r}}-1\right)^4\right]}_{\text{Lemma \ref{lemma:p-moments-cross-term}-1}}\notag\\
        \lesssim&\underbrace{\sum_{r=1}^{t-1}\left(\sum_{s=r+1}^{t}\Pi_{s,r}^2\right)^2y_r(1)^4}_{\text{Lemma \ref{lemma:deterministic-summation-new}-1}}\cdot T^{3/8}\notag\\
        \lesssim&R_t^2\cdot T^{3/8}\notag\\
        =&T^{5/4}R_t^{3/2}\cdot (TR_t^{-2})^{-1/4}\cdot T^{-5/8}\quad(\text{rewrite in the target scale})\notag\\
        =&o(T^{5/4}R_t^{3/2}\log^2(\eta_tT))\quad(\text{Assumption \mainref{assumption:maximum-radius-strengthened} and Lemma \ref{lemma:R}})\enspace.
    \end{align}
    We next bound $B_2$.
    By Lemma \ref{lemma:deterministic-summation-1} and Lemma \ref{lemma:p-moments-cross-term}, we can bound $B_2$ by:
    \begin{align}\label{lemma:variance-estimated-squared-residual_eq12}
        \qquad B_2=&\sum_{1\leq t_1\neq t_2\leq t-1}\left(\sum_{s=t_1+1}^{t}\Pi_{s,t_1}^2\right)\left(\sum_{s=t_2+1}^{t}\Pi_{s,t_2}^2\right)y_{t_1}(1)^2y_{t_2}(1)^2\notag\\
        &\times\e\left[\left(\frac{\mathbf{1}[Z_{t_1}=1]}{p_{t_1}}-1\right)^2\left(\frac{\mathbf{1}[Z_{t_2}=1]}{p_{t_2}}-1\right)^2\right]\notag\\
        \lesssim&\sum_{1\leq t_1< t_2\leq t-1}\left(\sum_{s=t_1+1}^{t}\Pi_{s,t_1}^2\right)\left(\sum_{s=t_2+1}^{t}\Pi_{s,t_2}^2\right)y_{t_1}(1)^2y_{t_2}(1)^2\notag\\
        &\times\underbrace{\e\left[\left(\frac{\mathbf{1}[Z_{t_1}=1]}{p_{t_1}}-1\right)^2\left(\frac{\mathbf{1}[Z_{t_2}=1]}{p_{t_2}}-1\right)^2\right]}_{\text{Lemma \ref{lemma:p-moments-cross-term}-2}}\notag\\
        \lesssim&\left(\sum_{1\leq t_1< t_2\leq t-1}R_{t_1}^{-5/16}R_{t_2}^{-1/4}\left(\sum_{s=t_1+1}^{t}\Pi_{s,t_1}^2\right)\left(\sum_{s=t_2+1}^{t}\Pi_{s,t_2}^2\right)y_{t_1}(1)^2y_{t_2}(1)^2\right)\cdot T^{9/32}\notag\\
        \leq&\Bigg(\underbrace{\sum_{1\leq t_1, t_2\leq t-1}}_{\text{relax }t_1\neq t_2}R_{t_1}^{-5/16}R_{t_2}^{-1/4}\left(\sum_{s=t_1+1}^{t}\Pi_{s,t_1}^2\right)\left(\sum_{s=t_2+1}^{t}\Pi_{s,t_2}^2\right)y_{t_1}(1)^2y_{t_2}(1)^2\Bigg)\cdot T^{9/32}\notag\\
        =&\left(\sum_{s=1}^{t-1}R_s^{-1/4}\left(\sum_{r=s+1}^t\Pi_{r,s}^2\right)y_s(1)^2\right)\left(\sum_{s=1}^{t-1}R_s^{-5/16}\left(\sum_{r=s+1}^t\Pi_{r,s}^2\right)y_s(1)^2\right)\cdot T^{9/32}\notag\\
        =&\underbrace{\left(\sum_{r=1}^{t}\sum_{s=1}^{r-1}R_s^{-1/4}\Pi_{r,s}^2y_s(1)^2\right)\left(\sum_{r=1}^{t}\sum_{s=1}^{r-1}R_s^{-5/16}\Pi_{r,s}^2y_s(1)^2\right)}_{\text{Lemma \ref{lemma:deterministic-summation-1}-3}}\cdot T^{9/32}\quad(\text{swap the summation order})\notag\\
        \lesssim&R_t^{3/2-1/4}T^{1/4}\cdot R_t^{3/2-5/16}T^{1/4}\cdot T^{9/32}\notag\\
        =&T^{5/4}R_t^{3/2}\cdot (TR_t^{-2})^{-15/32}\quad(\text{rewrite in the target scale})\notag\\
        =&o(T^{5/4}R_t^{3/2}\log^2(\eta_tT))\quad(\text{Assumption \mainref{assumption:maximum-radius-strengthened} and Lemma \ref{lemma:R}})\enspace.
    \end{align}
    It remains to control $B_3$, which involves cross terms.
    We therefore expand the variance expression and identify the index configurations that lead to nonzero contributions.
    Lemma \ref{lemma:covariance} indicates that the only nonzero terms in the expansion of $B_3$ have four indices where exactly two indices attain the maximal index among the four.
    Up to symmetry, there are two cases that satisfy this requirement: (i) $t_1=t_3$ and $t_2=t_4$. (ii) $t_1=t_3$ with $t_1>t_2$, $t_3>t_4$, and $t_2\neq t_4$.
    Hence we can expand $B_3$ as
    \begin{align}\label{lemma:variance-estimated-squared-residual_eq13}
        B_3=&\sum_{1\leq t_1\neq t_2,t_3\neq t_4\leq t-1}\operatorname{Cov}\Bigg(\left(\sum_{s=t_1\vee t_2+1}^{t}\Pi_{s,t_1}\Pi_{s,t_2}\right)y_{t_1}(1)y_{t_2}(1)\left(\frac{\mathbf{1}[Z_{t_1}=1]}{p_{t_1}}-1\right)\left(\frac{\mathbf{1}[Z_{t_2}=1]}{p_{t_2}}-1\right),\notag\\
        &\quad\quad\quad\quad\quad\quad\quad\quad\left(\sum_{s=t_3\vee t_4+1}^{t}\Pi_{s,t_3}\Pi_{s,t_4}\right)y_{t_3}(1)y_{t_4}(1)\left(\frac{\mathbf{1}[Z_{t_3}=1]}{p_{t_3}}-1\right)\left(\frac{\mathbf{1}[Z_{t_4}=1]}{p_{t_4}}-1\right)\Bigg)\notag\\
        =&2\underbrace{\sum_{1\leq t_1\neq t_2\leq t-1}\left(\sum_{s=t_1\vee t_2+1}^{t}\Pi_{s,t_1}\Pi_{s,t_2}\right)^2y_{t_1}(1)^2y_{t_2}(1)^2\operatorname{Var}\left[\left(\frac{\mathbf{1}[Z_{t_1}=1]}{p_{t_1}}-1\right)\left(\frac{\mathbf{1}[Z_{t_2}=1]}{p_{t_2}}-1\right)\right]}_{:=B_{3,1}}\notag\\
        &+4\sum_{1\leq t_2\neq t_3<t_1\leq t-1}\left(\sum_{s=t_1+1}^{t}\Pi_{s,t_1}\Pi_{s,t_2}\right)\left(\sum_{s=t_1+1}^{t}\Pi_{s,t_1}\Pi_{s,t_3}\right)y_{t_1}(1)^2y_{t_2}(1)y_{t_3}(1)\notag\\
        &\qquad\underbrace{\times\operatorname{Cov}\left(\left(\frac{\mathbf{1}[Z_{t_1}=1]}{p_{t_1}}-1\right)\left(\frac{\mathbf{1}[Z_{t_2}=1]}{p_{t_2}}-1\right),\left(\frac{\mathbf{1}[Z_{t_1}=1]}{p_{t_1}}-1\right)\left(\frac{\mathbf{1}[Z_{t_3}=1]}{p_{t_3}}-1\right)\right)}_{:=B_{3,2}}\enspace,
    \end{align}
    where $B_{3,1}$ and $B_{3,2}$ correspond to case (i) and (ii), respectively.
    We then bound $B_{3,1}$ and $B_{3,2}$ separately.
    By Lemma \ref{lemma:deterministic-summation-1} and Lemma \ref{lemma:p-moments-cross-term}, $B_{3,1}$ can be bounded by
    \begin{align}\label{lemma:variance-estimated-squared-residual_eq14}
        \qquad B_{3,1}\leq& \sum_{1\leq t_1\neq t_2\leq t-1}\left(\sum_{s=t_1\vee t_2+1}^{t}\Pi_{s,t_1}^2\right)\left(\sum_{s=t_1\vee t_2+1}^{t}\Pi_{s,t_2}^2\right)y_{t_1}(1)^2y_{t_2}(1)^2\notag\\
        &\times\operatorname{Var}\left[\left(\frac{\mathbf{1}[Z_{t_1}=1]}{p_{t_1}}-1\right)\left(\frac{\mathbf{1}[Z_{t_2}=1]}{p_{t_2}}-1\right)\right]\notag\\
        \lesssim& \sum_{1\leq t_1< t_2\leq t-1}\left(\sum_{s=t_1\vee t_2+1}^{t}\Pi_{s,t_1}^2\right)\left(\sum_{s=t_1\vee t_2+1}^{t}\Pi_{s,t_2}^2\right)y_{t_1}(1)^2y_{t_2}(1)^2\notag\\
        &\times\underbrace{\e\left[\left(\frac{\mathbf{1}[Z_{t_1}=1]}{p_{t_1}}-1\right)^2\left(\frac{\mathbf{1}[Z_{t_2}=1]}{p_{t_2}}-1\right)^2\right]}_{\text{Lemma \ref{lemma:p-moments-cross-term}-2}}\notag\\
        \lesssim&\left(\sum_{1\leq t_1< t_2\leq t-1}R_{t_1}^{-5/16}R_{t_2}^{-1/4}\left(\sum_{s=t_1+1}^{t}\Pi_{s,t_1}^2\right)\left(\sum_{s=t_2+1}^{t}\Pi_{s,t_2}^2\right)y_{t_1}(1)^2y_{t_2}(1)^2\right)\cdot T^{9/32}\notag\\
        \leq&\Bigg(\underbrace{\sum_{1\leq t_1, t_2\leq t-1}}_{\text{relax }t_1\neq t_2}R_{t_1}^{-5/16}R_{t_2}^{-1/4}\left(\sum_{s=t_1+1}^{t}\Pi_{s,t_1}^2\right)\left(\sum_{s=t_2+1}^{t}\Pi_{s,t_2}^2\right)y_{t_1}(1)^2y_{t_2}(1)^2\Bigg)\cdot T^{9/32}\notag\\
        =&\left(\sum_{s=1}^{t-1}R_s^{-1/4}\left(\sum_{r=s+1}^t\Pi_{r,s}^2\right)y_s(1)^2\right)\left(\sum_{s=1}^{t-1}R_s^{-5/16}\left(\sum_{r=s+1}^t\Pi_{r,s}^2\right)y_s(1)^2\right)\cdot T^{9/32}\notag\\
        =&\underbrace{\left(\sum_{r=1}^{t}\sum_{s=1}^{r-1}R_s^{-1/4}\Pi_{r,s}^2y_s(1)^2\right)\left(\sum_{r=1}^{t}\sum_{s=1}^{r-1}R_s^{-5/16}\Pi_{r,s}^2y_s(1)^2\right)}_{\text{Lemma \ref{lemma:deterministic-summation-1}-3}}\cdot T^{9/32}\quad(\text{swap the summation order})\notag\\
        \lesssim&R_t^{3/2-1/4}T^{1/4}\cdot R_t^{3/2-5/16}T^{1/4}\cdot T^{9/32}\notag\\
        =&T^{5/4}R_t^{3/2}\cdot (TR_t^{-2})^{-15/32}\quad(\text{rewrite in the target scale})\notag\\
        =&o(T^{5/4}R_t^{3/2}\log^2(\eta_tT))\quad(\text{Assumption \mainref{assumption:maximum-radius-strengthened} and Lemma \ref{lemma:R}})\enspace.
    \end{align}
    Applying Lemma \ref{lemma:p-moments-cross-term-new}-1 to bound the covariance term, we can bound $B_{3,2}$ as:
    \begin{align}\label{lemma:variance-estimated-squared-residual_eq15}
        |B_{3,2}|\leq&\sum_{1\leq t_2\neq t_3<t_1\leq t-1}\left|\left(\sum_{s=t_1+1}^{t}\Pi_{s,t_1}\Pi_{s,t_2}\right)\left(\sum_{s=t_1+1}^{t}\Pi_{s,t_1}\Pi_{s,t_3}\right)y_{t_1}(1)^2\,y_{t_2}(1)\,y_{t_3}(1)\right|\notag\\
        &\times\underbrace{\left|\operatorname{Cov}\left(\left(\frac{\mathbf{1}[Z_{t_1}=1]}{p_{t_1}}-1\right)\left(\frac{\mathbf{1}[Z_{t_2}=1]}{p_{t_2}}-1\right),\left(\frac{\mathbf{1}[Z_{t_1}=1]}{p_{t_1}}-1\right)\left(\frac{\mathbf{1}[Z_{t_3}=1]}{p_{t_3}}-1\right)\right)\right|}_{\text{Lemma \ref{lemma:p-moments-cross-term-new}-1}}\notag\\
        \leq&4\sum_{1\leq t_2\neq t_3<t_1\leq t-1}\left|\left(\sum_{s=t_1+1}^{t}\Pi_{s,t_1}\Pi_{s,t_2}\right)\left(\sum_{s=t_1+1}^{t}\Pi_{s,t_1}\Pi_{s,t_3}\right)y_{t_1}(1)^2\,y_{t_2}(1)\,y_{t_3}(1)\right|\notag\\
        &\times\left(2+\bb_1(\bb_2/6)^{1/4}\eta_{t_1}^{1/4}\e[\widehat{A}_{t_1}(0)]^{1/4}\right)^{25/16}\left(\frac{R_{t_1}}{R_{t_2}}\right)^{5/64}\left(\frac{R_{t_1}}{R_{t_3}}\right)^{5/64}\notag\\
        \intertext{Since $t_1\leq t$, we have $R_{t_1}\leq R_t$. By the fact that $\eta_{t_1}/\eta_t=R_t/R_{t_1}\geq 1$ and $\widehat{A}_{t_1}(0)\leq \widehat{A}_{t}(0)$, we can further bound it as}
        \leq&4\sum_{1\leq t_2\neq t_3<t_1\leq t-1}\left|\left(\sum_{s=t_1+1}^{t}\Pi_{s,t_1}\Pi_{s,t_2}\right)\left(\sum_{s=t_1+1}^{t}\Pi_{s,t_1}\Pi_{s,t_3}\right)y_{t_1}(1)^2\,y_{t_2}(1)\,y_{t_3}(1)\right|\notag\\
        &\times\left(2+\bb_1(\bb_2/6)^{1/4}\eta_{t}^{1/4}\e[\widehat{A}_{t}(0)]^{1/4}\right)^{25/16}\left(\frac{R_t}{R_{t_1}}\right)^{25/64}\left(\frac{R_{t_1}}{R_{t_2}}\right)^{5/64}\left(\frac{R_{t_1}}{R_{t_3}}\right)^{5/64}\notag\\
        =&4\sum_{1\leq t_2\neq t_3<t_1\leq t-1}\left|\left(\sum_{s=t_1+1}^{t}\Pi_{s,t_1}\Pi_{s,t_2}\right)\left(\sum_{s=t_1+1}^{t}\Pi_{s,t_1}\Pi_{s,t_3}\right)y_{t_1}(1)^2\,y_{t_2}(1)\,y_{t_3}(1)\right|\notag\\
        &\times\left(2+\bb_1(\bb_2/6)^{1/4}a_{t,0}^{1/4}\right)^{25/16}\left(\frac{R_t}{R_{t_1}}\right)^{15/64}\left(\frac{R_t}{R_{t_2}}\right)^{5/64}\left(\frac{R_t}{R_{t_3}}\right)^{5/64}\notag\\
        \intertext{Here we drop the restriction $t_2 \neq t_3$, which leads to the increase in the summation:}
        \leq&4\sum_{1\leq t_2,t_3<t_1\leq t-1}\underbrace{\left|\left(\sum_{s=t_1+1}^{t}\Pi_{s,t_1}\Pi_{s,t_2}\right)\left(\sum_{s=t_1+1}^{t}\Pi_{s,t_1}\Pi_{s,t_3}\right)y_{t_1}(1)^2\,y_{t_2}(1)\,y_{t_3}(1)\left(\frac{R_t^2}{R_{t_1}^2}\cdot\frac{R_t}{R_{t_2}}\cdot\frac{R_t}{R_{t_3}}\right)^{1/4}\right|}_{\text{Lemma \ref{lemma:deterministic-summation-new}-3}}\notag\\
        &\times\left(2+\bb_1(\bb_2/6)^{1/4}a_{t,0}^{1/4}\right)^{25/16}~~~~(\text{since $R_{t_1},R_{t_2},R_{t_3}\leq R_t$})\notag\\
        \leq& 2^{3/2}c\,c_1^2\gamma(\zeta_1^{1/2}+c_1)^2\cdot R_t^{3/2}T^{5/4}\log^2(\eta_tT)\left(2+\bb_1(\bb_2/6)^{1/4}a_{t,0}^{1/4}\right)^{25/16}\enspace.
    \end{align}
    Combining the upper bounds in \eqref{lemma:variance-estimated-squared-residual_eq8}, \eqref{lemma:variance-estimated-squared-residual_eq9}, \eqref{lemma:variance-estimated-squared-residual_eq10}, \eqref{lemma:variance-estimated-squared-residual_eq11}, \eqref{lemma:variance-estimated-squared-residual_eq12}, \eqref{lemma:variance-estimated-squared-residual_eq13}, \eqref{lemma:variance-estimated-squared-residual_eq14} and \eqref{lemma:variance-estimated-squared-residual_eq15}, $S_2$ is bounded as:
    \begin{align}\label{lemma:variance-estimated-squared-residual_eq16}
        \quad\qquad&S_2\notag\\
        \leq& 8S_{2,1}+2S_{2,2}\notag\\
        \leq&8S_{2,1}+4B_{1}+4B_{2}+8B_{3,1}+16|B_{3,2}|\notag\\
        \leq&8c\,c_1^2\gamma\xi_2^{1/2}\zeta_1T^{5/4}R_t^{3/2}\left(2+\bb_1(\bb_2/6)^{1/4}a_{t,0}^{1/4}\right)\notag\\
        &+2^{11/2}c\,c_1^2\gamma(\zeta_1^{1/2}+c_1)^2R_t^{3/2}T^{5/4}\log^2(\eta_tT)\cdot \left(2+\bb_1(\bb_2/6)^{1/4}a_{t,0}^{1/4}\right)^{25/16}\notag\\
        &+o(T^{5/4}R_t^{3/2}\log^2(\eta_tT))\notag\\
        \leq&c\,c_1^2\gamma\left(8\xi_2^{1/2}\zeta_1+2^{11/2}(\zeta_1^{1/2}+c_1)^2\right)R_t^{3/2}T^{5/4}\log^2(\eta_tT)\cdot \left(2+\bb_1(\bb_2/6)^{1/4}a_{t,0}^{1/4}\right)^{25/16}\notag\\
        &+o(T^{5/4}R_t^{3/2}\log^2(\eta_tT))\enspace.
    \end{align}
    \textbf{Step 3: }Combine the results in Step 1 and Step 2 to obtain the final bound.\\[3mm]
    For sufficiently large $T$, by \eqref{lemma:variance-estimated-squared-residual_eq1}, \eqref{lemma:variance-estimated-squared-residual_eq7} and \eqref{lemma:variance-estimated-squared-residual_eq16}, we have
    \begin{align*}
        \Var{\widehat{A}_t(1)}\leq&(S_1^{1/2}+S_2^{1/2})^2\\
        \leq& c\,c_1^2\gamma\left(8\xi_2^{1/2}\zeta_1+2^{11/2}(\zeta_1^{1/2}+c_1)^2\right)R_t^{3/2}T^{5/4}\log^2(\eta_tT)\cdot \left(2+\bb_1(\bb_2/6)^{1/4}a_{t,0}^{1/4}\right)^{25/16}\\
        &+o\left(T^{5/4}R_t^{3/2}\log^2(\eta_tT)\cdot \left(2+\bb_1(\bb_2/6)^{1/4}a_{t,0}^{1/4}\right)^{25/16}\right)\enspace.
    \end{align*}
    Using the identity $R_t^{3/2}T^{5/4}=\eta_t^{-2}(\eta_tT)^{1/2}$ and $2^{11/2}<46$, this indicates that
    \begin{align*}
        \Var{\widehat{a}_{t,1}}=\Var{\eta_t\widehat{A}_t(1)}\leq&C(\eta_tT)^{1/2}\log^2(\eta_tT)\left(2+\bb_1(\bb_2/6)^{1/4}a_{t,0}^{1/4}\right)^{25/16}\enspace
    \end{align*}
    for sufficiently large $T$.
\end{proof}

Lemma \ref{lemma:variance-estimated-squared-residual} establishes that the variances of $\widehat{a}_{t,1}$ (or $\widehat{A}_t(1)$) and $\widehat{a}_{t,0}$ (or $\widehat{A}_t(0)$) are bounded through a mutually normalizing manner; specifically, the variance of $\widehat{a}_{t,1}$ is bounded in terms of $a_{t,0}$, while that of $\widehat{a}_{t,0}$ is bounded in terms of $a_{t,1}$.
Such dependence can be undesirable, as the estimated squared residuals for each group may fail to concentrate around their respective expectations when the treatment group and control group are imbalanced. 
Nevertheless, as implied by Lemma \ref{lemma:squared-residuals-random}, once the expected squared residuals exceed a threshold, the maximal ratio between the two groups is controlled, thereby ensuring concentration despite the mutual normalization.
The result is stated formally in the following lemma.

\begin{lemma}\label{lemma:expected-difference-inverse-probability}
    Suppose $T$ is sufficiently large, and let $0 < \delta < \frac{1}{4}$ be an arbitrary constant. 
    Let $\widetilde{C} > 0$ be the constant defined by
    \begin{align*}
        \widetilde{C}:=&\max\Bigg\{\frac{2^{-7/32}C^{1/2}\bb_1^{25/32}\bb_2^{281/128}}{\bb_3},\frac{2^{57/32}C^{1/2}\bb_1^{89/32}\bb_2^{185/128}}{\bb_3}\left(1+\frac{2}{\bb_3}\right),\\
        &\quad\quad\quad 2^{283/64}C^{3/4}\bb_1^{139/64}\bb_2^{139/256},~\frac{2^{89/32}C^{1/2}\bb_1^{121/32}\bb_2^{217/128}\varsigma^{89/128}}{\bb_3^{5/2}}\left(1+\frac{2}{\bb_3}\right),\\
        &\quad\quad\quad 2^{347/64}C^{3/4}\bb_1^{139/64}\bb_2^{139/256}\varsigma^{139/256},~\frac{8\bb_1\bb_2}{\bb_3},~33C^{3/4}\bb_1^{139/64}\bb_2^{139/256}\Bigg\}\enspace.
    \end{align*}
    Then under Assumptions \mainref{assumption:moments}-\mainref{assumption:maximum-radius-strengthened} and Condition \mainref{condition:sigmoid}, the following bounds hold:
    \begin{enumerate}
        \item[(1)] For any $t\in[T]$ such that $a_{t,1}\geq 7\bb_1^{-4}\bb_2^{-1}$, the following holds:
        \begin{align*}
            &\e\left[\left|\frac{1}{p_{t+1}}-\frac{1}{\widebar{p}_{t+1}}\right|\right]\\
            \leq&\widetilde{C}\max\Bigg\{(\eta_tT)^{1/4}\log(\eta_tT)a_{t,1}^{-103/128}+(\eta_tT)^{1/4}\log(\eta_tT)a_{t,1}^{-71/128}+(\eta_tT)^{3/8}\log^{3/2}(\eta_tT)\delta^{-3/2}a_{t,1}^{-245/256},\notag\\
            &(\eta_tT)^{1/4}\log(\eta_tT)a_{t,1}^{-71/128}+(\eta_tT)^{121/128}\log(\eta_tT)a_{t,1}^{-3/2}+(\eta_tT)^{235/256}\log^{3/2}(\eta_tT)\delta^{-3/2}a_{t,1}^{-3/2},\notag\\
        &\delta+(\eta_tT)^{3/8}\log^{3/2}(\eta_tT)\delta^{-3/2}a_{t,1}^{-245/256}\Bigg\}\enspace.
        \end{align*}
        \item[(2)] For any $t\in[T]$ such that $a_{t,0}\geq 7\bb_1^{-4}\bb_2^{-1}$, the following holds:
        \begin{align*}
            &\e\left[\left|\frac{1}{1-p_{t+1}}-\frac{1}{1-\widebar{p}_{t+1}}\right|\right]\\
            \leq&\widetilde{C}\max\Bigg\{(\eta_tT)^{1/4}\log(\eta_tT)a_{t,0}^{-103/128}+(\eta_tT)^{1/4}\log(\eta_tT)a_{t,0}^{-71/128}+(\eta_tT)^{3/8}\log^{3/2}(\eta_tT)\delta^{-3/2}a_{t,0}^{-245/256},\notag\\
            &(\eta_tT)^{1/4}\log(\eta_tT)a_{t,0}^{-71/128}+(\eta_tT)^{121/128}\log(\eta_tT)a_{t,0}^{-3/2}+(\eta_tT)^{235/256}\log^{3/2}(\eta_tT)\delta^{-3/2}a_{t,0}^{-3/2},\notag\\
            &\delta+(\eta_tT)^{3/8}\log^{3/2}(\eta_tT)\delta^{-3/2}a_{t,0}^{-245/256}\Bigg\}\enspace.
        \end{align*}
    \end{enumerate}
\end{lemma}

\begin{proof}
    Without loss of generality, we prove only the first part. 
    For any $t\in[T]$ satisfying $a_{t,1}\geq 7\bb_1^{-4}\bb_2^{-1}$, we have
    \begin{align}\label{lemma:expected-difference-inverse-probability_eq1}
        2+\bb_1(\bb_2/6)^{1/4}a_{t,1}^{1/4}\leq 2\bb_1\bb_2^{1/4}a_{t,1}^{1/4}\enspace.
    \end{align}
    This inequality is used to simplify expressions throughout the proof.
    
    The high-level intuition for the proof is as follows.
    When the ratio between $a_{t,1}$ and $a_{t,0}$ is bounded away from $1$, we can lower-bound the probability that $(\widehat{a}_{t,1}, \widehat{a}_{t,0})$ has the same ordering as $(a_{t,1}, a_{t,0})$ using Chebyshev's inequality and the variance bounds in Lemma \ref{lemma:variance-estimated-squared-residual}.
    This ordering implies that $p_{t+1}$ and $\widebar{p}_{t+1}$ lie on the same side of $1/2$ by Lemma \ref{lemma:p-location}.
    Then the difference between the corresponding inverse probabilities can be bounded using Lemmas \ref{lemma:difference-inverse-probability-main} and \ref{lemma:variance-estimated-squared-residual}.
    On the other hand, when the ratio between $a_{t,1}$ and $a_{t,0}$ is close to $1$, we can directly bound the difference between $\widebar{p}_{t+1}$ and $1/2$, and similarly bound the difference between $p_{t+1}$ and $1/2$ using Lemma \ref{lemma:difference-inverse-probability-half}.
    We consider the following three cases:\\[3mm]
    \textbf{Case 1:} $a_{t,0}\leq (1-\delta)a_{t,1}$.\\[3mm]
        By splitting the expectation according to the events, we have the following decomposition:
         \begin{align}\label{lemma:expected-difference-inverse-probability_eq2}
            &\e\left[\left|\frac{1}{p_{t+1}}-\frac{1}{\widebar{p}_{t+1}}\right|\right]\notag\\
            =&\underbrace{\e\left[\left|\frac{1}{p_{t+1}}-\frac{1}{\widebar{p}_{t+1}}\right|\mathbf{1}\left[\widehat{a}_{t,1}\geq \widehat{a}_{t,0}\right]\right]}_{:=D_1}+\underbrace{\e\left[\left|\frac{1}{p_{t+1}}-\frac{1}{\widebar{p}_{t+1}}\right|\mathbf{1}\left[\widehat{a}_{t,1}<\widehat{a}_{t,0}\right]\right]}_{:=D_2}\enspace. 
        \end{align}
        The two terms correspond to whether the ordering of $\widehat a_{t,1}$ and $\widehat a_{t,0}$ matches that of $a_{t,1}$ and $a_{t,0}$.\\[3mm]
        \textbf{Step 1.1:} Derive an upper bound on $D_1$.\\[3mm]
        On the event $\widehat{a}_{t,1}\geq \widehat{a}_{t,0}$, the ordering between $\widehat{a}_{t,1},\widehat{a}_{t,0}$ is the same as $a_{t,1},a_{t,0}$.
        Hence Lemma \ref{lemma:difference-inverse-probability-main} can be applied to bound the difference between $1/p_{t+1}$ and $1/\widebar{p}_{t+1}$.
        
        On the event $\widehat{a}_{t,1}\geq \widehat{a}_{t,0}$, we have $p_{t+1},\widebar{p}_{t+1}\geq 1/2$ by Lemma \ref{lemma:p-location}. By definition, $p_{t+1}$ and $\widebar{p}_{t+1}$ satisfy the following first-order equations:
        \begin{align*}
            -\frac{\e[\widehat{A}_{t}(1)]}{\widebar{p}_{t+1}^2}+\frac{\e[\widehat{A}_{t}(0)]}{(1-\widebar{p}_{t+1})^2}+\eta_{t+1}^{-1}\Psi^{\prime}(\widebar{p}_{t+1})&=0\enspace,\\
            -\frac{\widehat{A}_{t}(1)}{p_{t+1}^2}+\frac{\widehat{A}_{t}(0)}{(1-p_{t+1})^2}+\eta_{t+1}^{-1}\Psi^{\prime}(p_{t+1})&=0\enspace.
        \end{align*}
        Using the notations $\widehat{a}_{t,k}=\eta_t\widehat{A}_{t}(k)$ and ${a}_{t,k}=\eta_t\e[\widehat{A}_{t}(k)]$ for $k\in\{0,1\}$, we can rewrite the equations as:
        \begin{align*}
            -\frac{a_{t,1}\cdot \frac{\eta_{t+1}}{\eta_t}}{\widebar{p}_{t+1}^2}+\frac{a_{t,0}\cdot \frac{\eta_{t+1}}{\eta_t}}{(1-\widebar{p}_{t+1})^2}+\Psi^{\prime}(\widebar{p}_{t+1})&=0\enspace,\\
            -\frac{\widehat{a}_{t,1}\cdot \frac{\eta_{t+1}}{\eta_t}}{p_{t+1}^2}+\frac{\widehat{a}_{t,0}\cdot \frac{\eta_{t+1}}{\eta_t}}{(1-p_{t+1})^2}+\Psi^{\prime}(p_{t+1})&=0\enspace.
         \end{align*}
         If $a_{t,0}=0$, i.e., $\e[\widehat{a}_{t,0}]=0$, then $\widehat{a}_{t,0}=0$ almost surely by its nonnegativity.
         Together with the assumption that $a_{t,1}>0$, this means that we can always apply Lemma \ref{lemma:difference-inverse-probability-main} whether $a_{t,0}=0$ or not.
         By Lemma \mainref{lemma:effect-of-p-regularization}* and Lemma \ref{lemma:difference-inverse-probability-main}, we have
         \begin{align}\label{lemma:expected-difference-inverse-probability_eq3}
             \left|\frac{1}{p_{t+1}}-\frac{1}{\widebar{p}_{t+1}}\right|\leq&\frac{\bb_2^2}{2\bb_3}\frac{|a_{t,1}-\widehat{a}_{t,1}|\cdot \frac{\eta_{t+1}}{\eta_t}}{a_{t,1}\cdot \frac{\eta_{t+1}}{\eta_t}}+\frac{\bb_1\bb_2}{\bb_3}\left(1+\frac{2}{\bb_3}\right)\frac{\widebar{p}_{t+1}}{1-\widebar{p}_{t+1}}\cdot\frac{|a_{t,0}-\widehat{a}_{t,0}|\cdot \frac{\eta_{t+1}}{\eta_t}}{a_{t,1}\cdot \frac{\eta_{t+1}}{\eta_t}}\quad(\text{Lemma \ref{lemma:difference-inverse-probability-main}})\notag\\
             \leq&\frac{\bb_2^2}{2\bb_3}\frac{|a_{t,1}-\widehat{a}_{t,1}|\cdot \frac{\eta_{t+1}}{\eta_t}}{a_{t,1}\cdot \frac{\eta_{t+1}}{\eta_t}}\notag\\
             &+\frac{\bb_1\bb_2}{\bb_3}\left(1+\frac{2}{\bb_3}\right)\frac{\left(1+\bb_1(\bb_2/6)^{1/4}a_{t,1}^{1/4}\cdot \left(\frac{\eta_{t+1}}{\eta_t}\right)^{1/4}\right)|a_{t,0}-\widehat{a}_{t,0}|}{a_{t,1}}\quad(\text{Lemma \mainref{lemma:effect-of-p-regularization}*})\notag\\
             \intertext{Since $\eta_{t+1}\leq \eta_t$, by \eqref{lemma:expected-difference-inverse-probability_eq1} we have $1+\bb_1(\bb_2/6)^{1/4}a_{t,1}^{1/4}\cdot \left(\frac{\eta_{t+1}}{\eta_t}\right)^{1/4}\leq 2+\bb_1(\bb_2/6)^{1/4}a_{t,1}^{1/4}\leq 2\bb_1\bb_2^{1/4}a_{t,1}^{1/4}$. Hence}
             \leq&\frac{\bb_2^2}{2\bb_3}\frac{|a_{t,1}-\widehat{a}_{t,1}|}{a_{t,1}}+\frac{2\bb_1^2\bb_2^{5/4}}{\bb_3}\left(1+\frac{2}{\bb_3}\right)\frac{|a_{t,0}-\widehat{a}_{t,0}|}{a_{t,1}^{3/4}}\enspace.
         \end{align}
        Hence by \eqref{lemma:expected-difference-inverse-probability_eq3} and the Cauchy-Schwarz inequality, we have
        \begin{align}\label{lemma:expected-difference-inverse-probability_eq4}
            D_1\leq& \e\left[\frac{\bb_2^2}{2\bb_3}\frac{|a_{t,1}-\widehat{a}_{t,1}|}{a_{t,1}}+\frac{2\bb_1^2\bb_2^{5/4}}{\bb_3}\left(1+\frac{2}{\bb_3}\right)\frac{|a_{t,0}-\widehat{a}_{t,0}|}{a_{t,1}^{3/4}}\right]\notag\\
            \leq&\frac{\bb_2^2}{2\bb_3}a_{t,1}^{-1}\operatorname{Var}^{1/2}(\widehat{a}_{t,1})+\frac{2\bb_1^2\bb_2^{5/4}}{\bb_3}\left(1+\frac{2}{\bb_3}\right)a_{t,1}^{-3/4}\operatorname{Var}^{1/2}(\widehat{a}_{t,0})\quad(\text{Cauchy-Schwarz})\notag\\
            \intertext{Applying Lemma \ref{lemma:variance-estimated-squared-residual} to $\operatorname{Var}(\widehat{a}_{t,1})$ and $\operatorname{Var}(\widehat{a}_{t,0})$ yields}
            \leq&\frac{C^{1/2}\bb_2^2}{2\bb_3}a_{t,1}^{-1}(\eta_tT)^{1/4}\log(\eta_tT) \left(2+\bb_1(\bb_2/6)^{1/4}a_{t,0}^{1/4}\right)^{25/32}\notag\\
            &+\frac{2C^{1/2}\bb_1^2\bb_2^{5/4}}{\bb_3}\left(1+\frac{2}{\bb_3}\right)a_{t,1}^{-3/4}(\eta_tT)^{1/4}\log(\eta_tT)\left(2+\bb_1(\bb_2/6)^{1/4}a_{t,1}^{1/4}\right)^{25/32}\notag\\
            \intertext{Since $a_{t,0}\leq a_{t,1}$, \eqref{lemma:expected-difference-inverse-probability_eq1} implies that $2+\bb_1(\bb_2/6)^{1/4}a_{t,0}^{1/4}\leq 2+\bb_1(\bb_2/6)^{1/4}a_{t,1}^{1/4}\leq 2\bb_1\bb_2^{1/4}a_{t,1}^{1/4}$. Hence, this can be further bounded as}
            \leq&\frac{2^{-7/32}C^{1/2}\bb_1^{25/32}\bb_2^{281/128}}{\bb_3}a_{t,1}^{-1}(\eta_tT)^{1/4}\log(\eta_tT)a_{t,1}^{25/128}\notag\\
            &+\frac{2^{57/32}C^{1/2}\bb_1^{89/32}\bb_2^{185/128}}{\bb_3}\left(1+\frac{2}{\bb_3}\right)a_{t,1}^{-3/4}(\eta_tT)^{1/4}\log(\eta_tT)a_{t,1}^{25/128}\notag\\
            \leq&\frac{2^{-7/32}C^{1/2}\bb_1^{25/32}\bb_2^{281/128}}{\bb_3}(\eta_tT)^{1/4}\log(\eta_tT)a_{t,1}^{-103/128}\notag\\
            &+\frac{2^{57/32}C^{1/2}\bb_1^{89/32}\bb_2^{185/128}}{\bb_3}\left(1+\frac{2}{\bb_3}\right)(\eta_tT)^{1/4}\log(\eta_tT)a_{t,1}^{-71/128}\enspace.
        \end{align}
        \textbf{Step 1.2:} Derive an upper bound on $D_2$.\\[3mm]
        On the event $\widehat{a}_{t,1}<\widehat{a}_{t,0}$, since $a_{t,0}\leq (1-\delta)a_{t,1}$, then either $|\widehat{a}_{t,1}-a_{t,1}|\geq \frac{\delta}{2} a_{t,1}$ or $|\widehat{a}_{t,0}-a_{t,0}|\geq \frac{\delta}{2} a_{t,1}$ holds.
        Hence, by the union bound, we obtain
        \begin{align}\label{lemma:expected-difference-inverse-probability_eq5}
            \operatorname{Pr}\left(\widehat{a}_{t,1}<\widehat{a}_{t,0}\right)\leq& \operatorname{Pr}\left(\frac{|\widehat{a}_{t,1}-a_{t,1}|}{a_{t,1}}+\frac{|\widehat{a}_{t,0}-a_{t,0}|}{a_{t,1}}\geq \delta\right)\notag\\
            \leq&\operatorname{Pr}\left(\frac{|\widehat{a}_{t,1}-a_{t,1}|}{a_{t,1}}\geq \frac{\delta}{2}\right)+\operatorname{Pr}\left(\frac{|\widehat{a}_{t,0}-a_{t,0}|}{a_{t,1}}\geq \frac{\delta}{2}\right)\notag\\
            \intertext{Applying Chebyshev's inequality together with Lemma \ref{lemma:variance-estimated-squared-residual} yields}
            \leq&\frac{4\operatorname{Var}\left(\widehat{a}_{t,1}\right)}{\delta^2a_{t,1}^2}+\frac{4\operatorname{Var}\left(\widehat{a}_{t,0}\right)}{\delta^2a_{t,1}^2}\quad(\text{Chebyshev's inequality})\notag\\
            \leq&\frac{4C(\eta_tT)^{1/2}\log^2(\eta_tT)\left(2+\bb_1(\bb_2/6)^{1/4}a_{t,0}^{1/4}\right)^{25/16}}{\delta^2a_{t,1}^2}\notag\\
            &+\frac{4C(\eta_tT)^{1/2}\log^2(\eta_tT)\left(2+\bb_1(\bb_2/6)^{1/4}a_{t,1}^{1/4}\right)^{25/16}}{\delta^2a_{t,1}^2}\quad(\text{Lemma \ref{lemma:variance-estimated-squared-residual}})\notag\\
            \leq&\frac{8C(\eta_tT)^{1/2}\log^2(\eta_tT)\left(2+\bb_1(\bb_2/6)^{1/4}a_{t,1}^{1/4}\right)^{25/16}}{\delta^2a_{t,1}^2}\quad(\text{since $a_{t,0}\leq a_{t,1}$})\notag\\
            \leq&\frac{8\cdot 2^{25/16}C\bb_1^{25/16}\bb_2^{25/64}(\eta_tT)^{1/2}\log^2(\eta_tT)a_{t,1}^{25/64}}{\delta^2a_{t,1}^2}\quad(\text{using \eqref{lemma:expected-difference-inverse-probability_eq1}})\notag\\
            \leq&\frac{2^{73/16}C\bb_1^{25/16}\bb_2^{25/64}(\eta_tT)^{1/2}\log^2(\eta_tT)}{\delta^2a_{t,1}^{103/64}}\enspace.
        \end{align}
        Note that we have $p_{t+1}\leq 1/2$ and $\widebar{p}_{t+1}\geq 1/2$ by Lemma \ref{lemma:p-location}. Hence by \eqref{lemma:expected-difference-inverse-probability_eq5}, Corollary \mainref{corollary:p-moment}* and H{\"o}lder's inequality, we have
        \begin{align}\label{lemma:expected-difference-inverse-probability_eq6}
            \qquad D_2=&\e\left[\left|\frac{1}{p_{t+1}}-\frac{1}{\widebar{p}_{t+1}}\right|\mathbf{1}\left[\widehat{a}_{t,1}<\widehat{a}_{t,0}\right]\right]\notag\\
            \leq&\e\left[\frac{1}{p_{t+1}}\mathbf{1}\left[\widehat{a}_{t,1}<\widehat{a}_{t,0}\right]\right]\quad(\text{since $1/p_{t+1}\geq 1/\widebar{p}_{t+1}\geq 0$})\notag\\
            \leq&\underbrace{\left(\e\left[\frac{1}{p_{t+1}^4}\right]\right)^{1/4}}_{\text{Corollary \mainref{corollary:p-moment}*}}\underbrace{\left(\operatorname{Pr}\left(\widehat{a}_{t,1}<\widehat{a}_{t,0}\right)\right)^{3/4}}_{\text{bound by }\eqref{lemma:expected-difference-inverse-probability_eq5}}\quad(\text{H{\"o}lder's inequality})\notag\\
            \leq&\left(2+\bb_1(\bb_2/6)^{1/4}a_{t,0}^{1/4}\right)\left(\frac{2^{73/16}C\bb_1^{25/16}\bb_2^{25/64}(\eta_tT)^{1/2}\log^2(\eta_tT)}{\delta^2a_{t,1}^{103/64}}\right)^{3/4}\notag\\
            \leq&2\bb_1\bb_2^{1/4}a_{t,1}^{1/4}\left(\frac{2^{73/16}C\bb_1^{25/16}\bb_2^{25/64}(\eta_tT)^{1/2}\log^2(\eta_tT)}{\delta^2a_{t,1}^{103/64}}\right)^{3/4}\quad(\text{using $a_{t,0}\leq a_{t,1}$ and \eqref{lemma:expected-difference-inverse-probability_eq1}})\notag\\
            \leq&2^{283/64}C^{3/4}\bb_1^{139/64}\bb_2^{139/256}(\eta_tT)^{3/8}\log^{3/2}(\eta_tT)\delta^{-3/2}a_{t,1}^{-245/256}\enspace.
        \end{align}
        \textbf{Step 1.3:} Combine the results in Step 1.1 and Step 1.2.\\[3mm]
        By combining the upper bounds in \eqref{lemma:expected-difference-inverse-probability_eq2}, \eqref{lemma:expected-difference-inverse-probability_eq4} and \eqref{lemma:expected-difference-inverse-probability_eq6}, we have
        \begin{align}\label{lemma:expected-difference-inverse-probability_eq7}
            \e\left[\left|\frac{1}{p_{t+1}}-\frac{1}{\widebar{p}_{t+1}}\right|\right]\leq&\frac{2^{-7/32}C^{1/2}\bb_1^{25/32}\bb_2^{281/128}}{\bb_3}(\eta_tT)^{1/4}\log(\eta_tT)a_{t,1}^{-103/128}\notag\\
            &+\frac{2^{57/32}C^{1/2}\bb_1^{89/32}\bb_2^{185/128}}{\bb_3}\left(1+\frac{2}{\bb_3}\right)(\eta_tT)^{1/4}\log(\eta_tT)a_{t,1}^{-71/128}\notag\\
            &+2^{283/64}C^{3/4}\bb_1^{139/64}\bb_2^{139/256}(\eta_tT)^{3/8}\log^{3/2}(\eta_tT)\delta^{-3/2}a_{t,1}^{-245/256}\enspace.
        \end{align}
        \textbf{Case 2:}  $a_{t,1}\leq (1-\delta)a_{t,0}$.\\[3mm]
        By splitting the expectation according to the events, we have the following decomposition:
         \begin{align}\label{lemma:expected-difference-inverse-probability_eq8}
            \e\left[\left|\frac{1}{p_{t+1}}-\frac{1}{\widebar{p}_{t+1}}\right|\right]=&\underbrace{\e\left[\left|\frac{1}{p_{t+1}}-\frac{1}{\widebar{p}_{t+1}}\right|\mathbf{1}\left[|\widehat{a}_{t,0}-a_{t,0}|< \delta a_{t,0}/2\right]\cdot\mathbf{1}\left[|\widehat{a}_{t,1}-a_{t,1}|< \delta a_{t,1}/2\right]\right]}_{:=D_1}\notag\\
            &+\underbrace{\e\left[\left|\frac{1}{p_{t+1}}-\frac{1}{\widebar{p}_{t+1}}\right|\mathbf{1}\left[\{|\widehat{a}_{t,0}-a_{t,0}|\geq \delta a_{t,0}/2\}\cup \{|\widehat{a}_{t,1}-a_{t,1}|\geq \delta a_{t,1}/2\}\right]\right]}_{:=D_2}\enspace. 
        \end{align}
        \textbf{Step 2.1:} Derive an upper bound on $D_1$.\\[3mm]
        On the events $\{|\widehat{a}_{t,0}-a_{t,0}|< \frac{\delta}{2}a_{t,0}\}$ and $\{|\widehat{a}_{t,1}-a_{t,1}|< \frac{\delta}{2}a_{t,1}\}$, we have
        \begin{align*}
            \widehat{a}_{t,0}-\widehat{a}_{t,1}\geq& (a_{t,0}-a_{t,1})-|\widehat{a}_{t,0}-a_{t,0}|-|\widehat{a}_{t,1}-a_{t,1}|\\
            \geq&\left(1-(1-\delta)\right)a_{t,0}-\frac{\delta}{2}a_{t,0}-\frac{\delta}{2}a_{t,1}\\
            \geq&\delta a_{t,0}-\frac{\delta}{2}a_{t,0}-\frac{\delta}{2}a_{t,0}\quad(\text{since $a_{t,0}\geq a_{t,1}$})\\
            =&0\enspace.
        \end{align*}
        Hence, the ordering between $\widehat{a}_{t,1},\widehat{a}_{t,0}$ is the same as $a_{t,1},a_{t,0}$.
        In this case, we have $p_{t+1},\widebar{p}_{t+1}\leq 1/2$ by Lemma \ref{lemma:p-location}.
        By definition, $p_{t+1}$ and $\widebar{p}_{t+1}$ satisfy the following first-order equations:
         \begin{align*}
            -\frac{a_{t,1}\cdot \frac{\eta_{t+1}}{\eta_t}}{\widebar{p}_{t+1}^2}+\frac{a_{t,0}\cdot \frac{\eta_{t+1}}{\eta_t}}{(1-\widebar{p}_{t+1})^2}+\Psi^{\prime}(\widebar{p}_{t+1})&=0\enspace,\\
            -\frac{\widehat{a}_{t,1}\cdot \frac{\eta_{t+1}}{\eta_t}}{p_{t+1}^2}+\frac{\widehat{a}_{t,0}\cdot \frac{\eta_{t+1}}{\eta_t}}{(1-p_{t+1})^2}+\Psi^{\prime}(p_{t+1})&=0\enspace.
         \end{align*}
         Because $\delta<1/4$ and $a_{t,0}>a_{t,1}>0$, by Lemma \ref{lemma:difference-inverse-probability-main} and the proof in Lemma \mainref{lemma:effect-of-p-regularization}*, we have
         \begin{align}\label{lemma:expected-difference-inverse-probability_eq9}
             \left|\frac{1}{p_{t+1}}-\frac{1}{\widebar{p}_{t+1}}\right|\leq&\frac{\bb_1\bb_2}{\bb_3}\left(1+\frac{2}{\bb_3}\right)\frac{1-\widebar{p}_{t+1}}{\widebar{p}_{t+1}}\frac{|\widehat{a}_{t,0}-a_{t,0}|\cdot \frac{\eta_{t+1}}{\eta_t}}{a_{t,0}\cdot \frac{\eta_{t+1}}{\eta_t}}\notag\\
             &+\frac{2\bb_1^2\bb_2}{\bb_3^2}\left(1+\frac{2}{\bb_3}\right)\frac{1-p_{t+1}}{p_{t+1}}\frac{|\widehat{a}_{t,1}-a_{t,1}|\cdot \frac{\eta_{t+1}}{\eta_t}}{a_{t,1}\cdot \frac{\eta_{t+1}}{\eta_t}}\quad(\text{Lemma \ref{lemma:difference-inverse-probability-main}})\notag\\
             \intertext{By Lemma \mainref{lemma:effect-of-p-regularization}*, we have $\frac{1-\widebar{p}_{t+1}}{\widebar{p}_{t+1}}=\frac{1}{\widebar{p}_{t+1}}-1\leq 1+\bb_1(\bb_2/6)^{1/4}\cdot \left(\frac{\eta_{t+1}}{\eta_t}\cdot a_{t,0}\right)^{1/4}$ and $\frac{1-p_{t+1}}{p_{t+1}}=\frac{1}{p_{t+1}}-1\leq \bb_1(\bb_2/\bb_3)^{1/2}\left(\frac{\widehat{a}_{t,0}\cdot \frac{\eta_{t+1}}{\eta_t}}{\widehat{a}_{t,1}\cdot \frac{\eta_{t+1}}{\eta_t}}\right)^{1/2}=\bb_1(\bb_2/\bb_3)^{1/2}\left(\frac{\widehat{a}_{t,0}}{\widehat{a}_{t,1}}\right)^{1/2}$. Hence}
             \leq&\frac{\bb_1\bb_2}{\bb_3}\left(1+\frac{2}{\bb_3}\right)\frac{\left(1+\bb_1(\bb_2/6)^{1/4}a_{t,0}^{1/4}\cdot \left(\frac{\eta_{t+1}}{\eta_t}\right)^{1/4}\right)|\widehat{a}_{t,0}-a_{t,0}|}{a_{t,0}}\notag\\
             &+\frac{2\bb_1^2\bb_2}{\bb_3^2}\left(1+\frac{2}{\bb_3}\right)\cdot \bb_1(\bb_2/\bb_3)^{1/2}\left(\frac{\widehat{a}_{t,0}}{\widehat{a}_{t,1}}\right)^{1/2}\frac{|\widehat{a}_{t,1}-a_{t,1}|}{a_{t,1}}\notag\\
             \intertext{Since $\eta_{t+1}\leq \eta_t$ and $7\bb_1^{-4}\bb_2^{-1}\leq a_{t,1}\leq a_{t,0}$, we have $1+\bb_1(\bb_2/6)^{1/4}a_{t,0}^{1/4}\cdot \left(\frac{\eta_{t+1}}{\eta_t}\right)^{1/4}\leq 2+\bb_1(\bb_2/6)^{1/4}a_{t,0}^{1/4}\leq 2\bb_1\bb_2^{1/4}a_{t,0}^{1/4}$. Moreover, $|\widehat{a}_{t,0}-a_{t,0}|\leq \frac{\delta}{2}a_{t,0}$ implies that $\widehat{a}_{t,0}\leq (1+\delta/2)a_{t,0}$. Similarly, we have $\widehat{a}_{t,1}\geq (1-\delta/2)a_{t,1}$. Hence}
             \leq&\frac{2\bb_1^2\bb_2^{5/4}}{\bb_3}\left(1+\frac{2}{\bb_3}\right)\frac{a_{t,0}^{1/4}|\widehat{a}_{t,0}-a_{t,0}|}{a_{t,0}}\notag\\
             &+\frac{2\bb_1^3\bb_2^{3/2}}{\bb_3^{5/2}}\left(1+\frac{2}{\bb_3}\right)\left(\frac{1+\delta/2}{1-\delta/2}\right)^{1/2}\left(\frac{a_{t,0}}{a_{t,1}}\right)^{1/2}\frac{|\widehat{a}_{t,1}-a_{t,1}|}{a_{t,1}}\notag\\
             \intertext{We have $a_{t,0}=\eta_t\e[\widehat{A}_{t}(0)]\leq \varsigma\cdot \eta_tT$ by Lemma \ref{lemma:squared-residuals-random} and $\left(\frac{1+\delta/2}{1-\delta/2}\right)^{1/2}\leq 2$ since $\delta<1/4$. Hence}
             \leq&\frac{2\bb_1^2\bb_2^{5/4}}{\bb_3}\left(1+\frac{2}{\bb_3}\right)\frac{|\widehat{a}_{t,0}-a_{t,0}|}{a_{t,0}^{3/4}}+\frac{4\bb_1^3\bb_2^{3/2}\varsigma^{1/2}}{\bb_3^{5/2}}\left(1+\frac{2}{\bb_3}\right)\frac{(\eta_tT)^{1/2}|\widehat{a}_{t,1}-a_{t,1}|}{a_{t,1}^{3/2}}\notag\\
             \leq&\frac{2\bb_1^2\bb_2^{5/4}}{\bb_3}\left(1+\frac{2}{\bb_3}\right)\frac{|\widehat{a}_{t,0}-a_{t,0}|}{a_{t,1}^{3/4}}+\frac{4\bb_1^3\bb_2^{3/2}\varsigma^{1/2}}{\bb_3^{5/2}}\left(1+\frac{2}{\bb_3}\right)\frac{(\eta_tT)^{1/2}|\widehat{a}_{t,1}-a_{t,1}|}{a_{t,1}^{3/2}}\enspace.
        \end{align}
        The last step is because $a_{t,1}\leq a_{t,0}$.
        By Lemma \ref{lemma:variance-estimated-squared-residual}, \eqref{lemma:expected-difference-inverse-probability_eq9} and the Cauchy-Schwarz inequality, we have
        \begin{align}\label{lemma:expected-difference-inverse-probability_eq10}
            D_1\leq& \e\left[\frac{2\bb_1^2\bb_2^{5/4}}{\bb_3}\left(1+\frac{2}{\bb_3}\right)\frac{|\widehat{a}_{t,0}-a_{t,0}|}{a_{t,1}^{3/4}}+\frac{4\bb_1^3\bb_2^{3/2}\varsigma^{1/2}}{\bb_3^{5/2}}\left(1+\frac{2}{\bb_3}\right)\frac{(\eta_tT)^{1/2}|\widehat{a}_{t,1}-a_{t,1}|}{a_{t,1}^{3/2}}\right]\notag\\
            \leq&\frac{2\bb_1^2\bb_2^{5/4}}{\bb_3}\left(1+\frac{2}{\bb_3}\right)a_{t,1}^{-3/4}\underbrace{\operatorname{Var}(\widehat{a}_{t,0})^{1/2}}_{\text{Lemma \ref{lemma:variance-estimated-squared-residual}}}+\frac{4\bb_1^3\bb_2^{3/2}\varsigma^{1/2}}{\bb_3^{5/2}}\left(1+\frac{2}{\bb_3}\right)(\eta_tT)^{1/2}a_{t,1}^{-3/2}\underbrace{\operatorname{Var}(\widehat{a}_{t,1})^{1/2}}_{\text{Lemma \ref{lemma:variance-estimated-squared-residual}}}\notag\\
            \leq&\frac{2\bb_1^2\bb_2^{5/4}}{\bb_3}\left(1+\frac{2}{\bb_3}\right)a_{t,1}^{-3/4}C^{1/2}(\eta_tT)^{1/4}\log(\eta_tT)\Big(\underbrace{2+\bb_1(\bb_2/6)^{1/4}a_{t,1}^{1/4}}_{\leq 2\bb_1\bb_2^{1/4}a_{t,1}^{1/4}}\Big)^{25/32}\notag\\
            &+\frac{4\bb_1^3\bb_2^{3/2}\varsigma^{1/2}}{\bb_3^{5/2}}\left(1+\frac{2}{\bb_3}\right)(\eta_tT)^{1/2} a_{t,1}^{-3/2}C^{1/2}(\eta_tT)^{1/4}\log(\eta_tT) \Big(\underbrace{2+\bb_1(\bb_2/6)^{1/4}a_{t,0}^{1/4}}_{\leq 2\bb_1\bb_2^{1/4}a_{t,0}^{1/4}}\Big)^{25/32}\notag\\
            \leq&\frac{2^{57/32}C^{1/2}\bb_1^{89/32}\bb_2^{185/128}}{\bb_3}\left(1+\frac{2}{\bb_3}\right)a_{t,1}^{-3/4}(\eta_tT)^{1/4}\log(\eta_tT)a_{t,1}^{25/128}\notag\\
            &+\frac{2^{89/32}C^{1/2}\bb_1^{121/32}\bb_2^{217/128}\varsigma^{1/2}}{\bb_3^{5/2}}\left(1+\frac{2}{\bb_3}\right) a_{t,1}^{-3/2}(\eta_tT)^{3/4}\log(\eta_tT)a_{t,0}^{25/128}\notag\\
            \intertext{Since $a_{t,0}=\eta_t\e[\widehat{A}_t(0)]\leq \varsigma\cdot\eta_tT$ by Lemma \ref{lemma:squared-residuals-random}, this can further be bounded as}
            \leq&\frac{2^{57/32}C^{1/2}\bb_1^{89/32}\bb_2^{185/128}}{\bb_3}\left(1+\frac{2}{\bb_3}\right)(\eta_tT)^{1/4}\log(\eta_tT)a_{t,1}^{-71/128}\notag\\
            &+\frac{2^{89/32}C^{1/2}\bb_1^{121/32}\bb_2^{217/128}\varsigma^{89/128}}{\bb_3^{5/2}}\left(1+\frac{2}{\bb_3}\right)(\eta_tT)^{121/128}\log(\eta_tT)a_{t,1}^{-3/2}\enspace.
        \end{align}
        \textbf{Step 2.2:} Derive an upper bound on $D_2$.\\[3mm]
        We first bound the probabilities of event $\{|\widehat{a}_{t,0}-a_{t,0}|\geq \frac{\delta}{2}a_{t,0}\}$ and event $\{|\widehat{a}_{t,1}-a_{t,1}|\geq \frac{\delta}{2}a_{t,1}\}$.
        By Chebyshev's inequality and Lemma \ref{lemma:variance-estimated-squared-residual}, we have
        \begin{align}\label{lemma:expected-difference-inverse-probability_eq12}
            \operatorname{Pr}\left(\frac{|\widehat{a}_{t,0}-a_{t,0}|}{a_{t,0}}\geq \frac{\delta}{2}\right)\leq&\frac{4\operatorname{Var}\left(\widehat{a}_{t,0}\right)}{\delta^2a_{t,0}^2}\quad(\text{Chebyshev's inequality})\notag\\
            \leq&\frac{4C(\eta_tT)^{1/2}\log^2(\eta_tT) \left(2+\bb_1(\bb_2/6)^{1/4}a_{t,1}^{1/4}\right)^{25/16}}{\delta^2a_{t,0}^2}\quad(\text{Lemma \ref{lemma:variance-estimated-squared-residual}})\notag\\
            \intertext{Since $a_{t,1}\leq a_{t,0}$, we have $2+\bb_1(\bb_2/6)^{1/4}a_{t,1}^{1/4}\leq 2\bb_1\bb_2^{1/4}a_{t,1}^{1/4}\leq 2\bb_1\bb_2^{1/4}a_{t,0}^{1/4}$. Hence}
            \leq&\frac{4\cdot 2^{25/16}C\bb_1^{25/16}\bb_2^{25/64}(\eta_tT)^{1/2}\log^2(\eta_tT)a_{t,0}^{25/64}}{\delta^2a_{t,0}^2}\notag\\
            \leq&\frac{2^{57/16}C\bb_1^{25/16}\bb_2^{25/64}(\eta_tT)^{1/2}\log^2(\eta_tT)}{\delta^2a_{t,0}^{103/64}}\notag\\
            \leq&\frac{2^{57/16}C\bb_1^{25/16}\bb_2^{25/64}(\eta_tT)^{1/2}\log^2(\eta_tT)}{\delta^2a_{t,1}^{103/64}}\quad(\text{since $a_{t,0}\geq a_{t,1}$})\enspace.
        \end{align}
        Since $7\bb_1^{-4}\bb_2^{-1}\leq a_{t,1}\leq a_{t,0}$, we can show that $2+\bb_1(\bb_2/6)^{1/4}a_{t,0}^{1/4}\leq 2\bb_1\bb_2^{1/4}a_{t,0}^{1/4}$.
        Moreover, we have $a_{t,0}\leq \varsigma\cdot \eta_tT$, which leads to
        \begin{align}\label{lemma:expected-difference-inverse-probability_eq13}
            \operatorname{Pr}\left(\frac{|\widehat{a}_{t,1}-a_{t,1}|}{a_{t,1}}\geq \frac{\delta}{2}\right)\leq&\frac{4\operatorname{Var}\left(\widehat{a}_{t,1}\right)}{\delta^2a_{t,1}^2}\quad(\text{Chebyshev's inequality})\notag\\
            \leq&\frac{4C(\eta_tT)^{1/2}\log^2(\eta_tT) \left(2+\bb_1(\bb_2/6)^{1/4}a_{t,0}^{1/4}\right)^{25/16}}{\delta^2a_{t,1}^2}\quad(\text{Lemma \ref{lemma:variance-estimated-squared-residual}})\notag\\
            \leq&\frac{2^{57/16}C\bb_1^{25/16}\bb_2^{25/64}(\eta_tT)^{1/2}\log^2(\eta_tT)(\varsigma\cdot\eta_tT)^{25/64}}{\delta^2a_{t,1}^2}\enspace.
        \end{align}
        By \eqref{lemma:expected-difference-inverse-probability_eq12} and \eqref{lemma:expected-difference-inverse-probability_eq13}, the union bound gives
        \begin{align}\label{lemma:expected-difference-inverse-probability_eq13_new}
            &\operatorname{Pr}\left(\left\{\frac{|\widehat{a}_{t,0}-a_{t,0}|}{a_{t,0}}\geq \frac{\delta}{2}\right\}\cup\left\{\frac{|\widehat{a}_{t,1}-a_{t,1}|}{a_{t,1}}\geq \frac{\delta}{2}\right\}\right)\notag\\
            \leq&\operatorname{Pr}\left(\frac{|\widehat{a}_{t,0}-a_{t,0}|}{a_{t,0}}\geq \frac{\delta}{2}\right)+\operatorname{Pr}\left(\frac{|\widehat{a}_{t,1}-a_{t,1}|}{a_{t,1}}\geq \frac{\delta}{2}\right)\notag\\
            \leq&2^{57/16}C\bb_1^{25/16}\bb_2^{25/64}(\eta_tT)^{1/2}\log^2(\eta_tT)\delta^{-2}\left(a_{t,1}^{-2}(\varsigma\cdot\eta_tT)^{25/64}+a_{t,1}^{-103/64}\right)\notag\\
            \intertext{Since $a_{t,1}\leq \varsigma\cdot\eta_tT$ by Lemma \ref{lemma:squared-residuals-random}, we have $a_{t,1}^{-103/64}\leq a_{t,1}^{-2}(\varsigma\cdot\eta_tT)^{25/64}$. Hence we obtain}
            \leq&2^{57/16}C\bb_1^{25/16}\bb_2^{25/64}(\eta_tT)^{1/2}\log^2(\eta_tT)\delta^{-2}\cdot 2a_{t,1}^{-2}(\varsigma\cdot\eta_tT)^{25/64}\notag\\
            \leq&2^{73/16}C\bb_1^{25/16}\bb_2^{25/64}\varsigma^{25/64}(\eta_tT)^{57/64}\log^2(\eta_tT)\delta^{-2}a_{t,1}^{-2}\enspace.
        \end{align}
        Therefore, by H{\"o}lder's inequality, \eqref{lemma:expected-difference-inverse-probability_eq13_new}, Lemma \mainref{lemma:effect-of-p-regularization}*, and Corollary \mainref{corollary:p-moment}*, we obtain
        \begin{align}\label{lemma:expected-difference-inverse-probability_eq14}
            D_2=&\e\left[\left|\frac{1}{p_{t+1}}-\frac{1}{\widebar{p}_{t+1}}\right|\mathbf{1}\left[\{|\widehat{a}_{t,0}-a_{t,0}|\geq \delta a_{t,0}/2\}\cup\{|\widehat{a}_{t,1}-a_{t,1}|\geq \delta a_{t,1}/2\}\right]\right]\notag\\
            \leq&\e\left[\left(\frac{1}{p_{t+1}}+\frac{1}{\widebar{p}_{t+1}}\right)\mathbf{1}\left[\{|\widehat{a}_{t,0}-a_{t,0}|\geq \delta a_{t,0}/2\}\cup\{|\widehat{a}_{t,1}-a_{t,1}|\geq \delta a_{t,1}/2\}\right]\right]\notag\\
            \leq&\Bigg(\underbrace{\e\left[\frac{1}{p_{t+1}^4}\right]^{1/4}+\e\left[\frac{1}{\widebar{p}_{t+1}^4}\right]^{1/4}}_{\text{Corollary \mainref{corollary:p-moment}*, Lemma \mainref{lemma:effect-of-p-regularization}*}}\Bigg)\underbrace{\left(\operatorname{Pr}\left[\left\{\frac{|\widehat{a}_{t,0}-a_{t,0}|}{a_{t,0}}\geq\frac{\delta}{2}\right\}\cup\left\{\frac{|\widehat{a}_{t,1}-a_{t,1}|}{a_{t,1}}\geq\frac{\delta}{2}\right\}\right]\right)^{3/4}}_{\text{bound by \eqref{lemma:expected-difference-inverse-probability_eq13_new}}}\notag\\
            \leq&2\underbrace{\Big(2+\bb_1(\bb_2/6)^{1/4}a_{t,0}^{1/4}\Big)}_{\leq 2\bb_1\bb_2^{1/4}a_{t,0}^{1/4}\leq 2\bb_1\bb_2^{1/4}(\varsigma\cdot \eta_tT)^{1/4}}\left(2^{73/16}C\bb_1^{25/16}\bb_2^{25/64}\varsigma^{25/64}(\eta_tT)^{57/64}\log^2(\eta_tT)\delta^{-2}a_{t,1}^{-2}\right)^{3/4}\notag\\
            \leq&4\bb_1\bb_2^{1/4}(\varsigma\cdot \eta_tT)^{1/4}\cdot 2^{219/64}C^{3/4}\bb_1^{75/64}\bb_2^{75/256}\varsigma^{75/256}(\eta_tT)^{171/256}\log^{3/2}(\eta_tT)\delta^{-3/2}a_{t,1}^{-3/2}\notag\\
            =&2^{347/64}C^{3/4}\bb_1^{139/64}\bb_2^{139/256}\varsigma^{139/256}(\eta_tT)^{235/256}\log^{3/2}(\eta_tT)\delta^{-3/2}a_{t,1}^{-3/2}\enspace.
        \end{align}
        \textbf{Step 2.3:} Combine the results in Step 2.1 and Step 2.2.\\[3mm]
        By combining the upper bounds in \eqref{lemma:expected-difference-inverse-probability_eq8}, \eqref{lemma:expected-difference-inverse-probability_eq10} and \eqref{lemma:expected-difference-inverse-probability_eq14}, we have
        \begin{align}\label{lemma:expected-difference-inverse-probability_eq15}
            \e\left[\left|\frac{1}{p_{t+1}}-\frac{1}{\widebar{p}_{t+1}}\right|\right]\leq &\frac{2^{57/32}C^{1/2}\bb_1^{89/32}\bb_2^{185/128}}{\bb_3}\left(1+\frac{2}{\bb_3}\right)(\eta_tT)^{1/4}\log(\eta_tT)a_{t,1}^{-71/128}\notag\\
            &+\frac{2^{89/32}C^{1/2}\bb_1^{121/32}\bb_2^{217/128}\varsigma^{89/128}}{\bb_3^{5/2}}\left(1+\frac{2}{\bb_3}\right)(\eta_tT)^{121/128}\log(\eta_tT)a_{t,1}^{-3/2}\notag\\
            &+2^{347/64}C^{3/4}\bb_1^{139/64}\bb_2^{139/256}\varsigma^{139/256}(\eta_tT)^{235/256}\log^{3/2}(\eta_tT)\delta^{-3/2}a_{t,1}^{-3/2}\enspace.
        \end{align}
        \textbf{Case 3:} $a_{t,0}\geq (1-\delta)a_{t,1}$ and $a_{t,1}\geq (1-\delta)a_{t,0}$.\\[3mm]
        By splitting the expectation according to the events, we have the following upper bound:
        \begin{align}\label{lemma:expected-difference-inverse-probability_eq16}
            \e\left[\left|\frac{1}{p_{t+1}}-\frac{1}{\widebar{p}_{t+1}}\right|\right]\leq&\underbrace{\e\left[\left|\frac{1}{p_{t+1}}-\frac{1}{\widebar{p}_{t+1}}\right|\mathbf{1}\left[|\widehat{a}_{t,1}-a_{t,1}|\leq \delta a_{t,0}\right]\cdot\mathbf{1}\left[|\widehat{a}_{t,0}-a_{t,0}|\leq \delta a_{t,1}\right]\right]}_{:=D_1}\notag\\
            &+\underbrace{\e\left[\left|\frac{1}{p_{t+1}}-\frac{1}{\widebar{p}_{t+1}}\right|\mathbf{1}\left[|\widehat{a}_{t,0}-a_{t,0}|\geq \delta a_{t,1}\right]\right]}_{:=D_2}\notag\\
            &+\underbrace{\e\left[\left|\frac{1}{p_{t+1}}-\frac{1}{\widebar{p}_{t+1}}\right|\mathbf{1}\left[|\widehat{a}_{t,1}-a_{t,1}|\geq \delta a_{t,0}\right]\right]}_{:=D_3}\enspace.
        \end{align}
        \textbf{Step 3.1:} Derive an upper bound on $D_1$.\\[3mm]
        Since $a_{t,1}>0$ by assumption, we also have $a_{t,0}\geq (1-\delta)a_{t,1}>0$.
        On the event $|\widehat{a}_{t,1}-a_{t,1}|\leq \delta a_{t,0}$ and $|\widehat{a}_{t,0}-a_{t,0}|\leq \delta a_{t,1}$, we have
        \begin{align}\label{lemma:expected-difference-inverse-probability_eqq}
            \frac{\widehat{a}_{t,0}}{\widehat{a}_{t,1}}\leq \frac{a_{t,0}+\delta a_{t,1}}{a_{t,1}-\delta a_{t,0}}\leq\frac{\left(1+\frac{\delta}{1-\delta}\right)a_{t,0}}{(1-2\delta)a_{t,0}}=\frac{1}{(1-\delta)(1-2\delta)}\enspace,\notag\\
            \frac{\widehat{a}_{t,1}}{\widehat{a}_{t,0}}\leq \frac{a_{t,1}+\delta a_{t,0}}{a_{t,0}-\delta a_{t,1}}\leq\frac{\left(1+\frac{\delta}{1-\delta}\right)a_{t,1}}{(1-2\delta)a_{t,1}}=\frac{1}{(1-\delta)(1-2\delta)}\enspace.
        \end{align}
        We then discuss four cases regarding the ordering among $a_{t,1},a_{t,0}$ and $\widehat{a}_{t,1},\widehat{a}_{t,0}$.
        On the events $\{a_{t,1}\leq a_{t,0}\}$ and $\{\widehat{a}_{t,1}\leq \widehat{a}_{t,0}\}$, it holds that $p_{t+1},\widebar{p}_{t+1}\leq 1/2$ by Lemma \ref{lemma:p-location}.
        Hence, by Lemma \ref{lemma:difference-inverse-probability-half} and \eqref{lemma:expected-difference-inverse-probability_eqq}, we have
        \begin{align*}
            \left|\frac{1}{p_{t+1}}-\frac{1}{\widebar{p}_{t+1}}\right|\leq& \left|\frac{1}{p_{t+1}}-2\right|+\left|\frac{1}{\widebar{p}_{t+1}}-2\right|\\
            \leq&\frac{\bb_1\bb_2}{\bb_3}\cdot \frac{|a_{t,1}-a_{t,0}|}{a_{t,1}}+\frac{\bb_1\bb_2}{\bb_3}\cdot \frac{|\widehat{a}_{t,1}-\widehat{a}_{t,0}|}{\widehat{a}_{t,1}}\quad(\text{Lemma \ref{lemma:difference-inverse-probability-half}})\\
            \leq&\frac{\bb_1\bb_2}{\bb_3}\left|\frac{a_{t,0}}{a_{t,1}}-1\right|+\frac{\bb_1\bb_2}{\bb_3}\left|\frac{\widehat{a}_{t,0}}{\widehat{a}_{t,1}}-1\right|\\
            \leq&\frac{\bb_1\bb_2}{\bb_3}\left(\frac{1}{1-\delta}-1\right)+\frac{\bb_1\bb_2}{\bb_3}\left(\frac{1}{(1-\delta)(1-2\delta)}-1\right)\quad(\text{by \eqref{lemma:expected-difference-inverse-probability_eqq}})\\
            =&\frac{4\bb_1\bb_2}{\bb_3}\frac{\delta}{1-2\delta}\\
            \leq&\frac{8\bb_1\bb_2}{\bb_3}\delta\quad(\text{since $\delta<1/4$})\enspace.
        \end{align*}
        Similarly, on the events $\{a_{t,1}\geq a_{t,0}\}$ and $\{\widehat{a}_{t,1}\geq \widehat{a}_{t,0}\}$, it holds that $p_{t+1},\widebar{p}_{t+1}\geq 1/2$ by Lemma \ref{lemma:p-location}. 
        Hence, we have
        \begin{align*}
            \left|\frac{1}{p_{t+1}}-\frac{1}{\widebar{p}_{t+1}}\right|\leq& \left|\frac{1}{p_{t+1}}-2\right|+\left|\frac{1}{\widebar{p}_{t+1}}-2\right|\\
            \leq&\left|\frac{\frac{1}{2}-p_{t+1}}{\frac{1}{2}p_{t+1}}\right|+\left|\frac{\frac{1}{2}-\widebar{p}_{t+1}}{\frac{1}{2}\widebar{p}_{t+1}}\right|\\
            \leq&\left|\frac{\frac{1}{2}-p_{t+1}}{\left(1-\frac{1}{2}\right)\left(1-p_{t+1}\right)}\right|+\left|\frac{\frac{1}{2}-\widebar{p}_{t+1}}{\left(1-\frac{1}{2}\right)\left(1-\widebar{p}_{t+1}\right)}\right|\\
            \leq&\left|\frac{1}{1-p_{t+1}}-2\right|+\left|\frac{1}{1-\widebar{p}_{t+1}}-2\right|\\
            \leq&\frac{8\bb_1\bb_2}{\bb_3}\delta\enspace.
        \end{align*}
        The last inequality is verified through Lemma \ref{lemma:difference-inverse-probability-half} and similar arguments as in the first case.
        Similarly, we can prove that under the remaining two cases, it holds that
        \begin{align*}
            \left|\frac{1}{p_{t+1}}-\frac{1}{\widebar{p}_{t+1}}\right|\leq \frac{8\bb_1\bb_2}{\bb_3}\delta\enspace.
        \end{align*}
        Hence, by combining the upper bounds in these four cases, $D_1$ can be bounded as:
        \begin{align}\label{lemma:expected-difference-inverse-probability_eq17}
            D_1\leq \frac{8\bb_1\bb_2}{\bb_3}\delta\enspace.
        \end{align}
        \textbf{Step 3.2:} Derive upper bounds on $D_2$ and $D_3$.\\[3mm]
        By Chebyshev's inequality and Lemma \ref{lemma:variance-estimated-squared-residual}, we have
        \begin{align}\label{lemma:expected-difference-inverse-probability_eq18}
            \operatorname{Pr}\left(\frac{|\widehat{a}_{t,0}-a_{t,0}|}{a_{t,1}}\geq \delta\right)\leq&\frac{\operatorname{Var}\left(\widehat{a}_{t,0}\right)}{\delta^2a_{t,1}^2}\quad(\text{Chebyshev's inequality})\notag\\
            \leq&\frac{C(\eta_tT)^{1/2}\log^2(\eta_tT)\left(2+\bb_1(\bb_2/6)^{1/4}a_{t,1}^{1/4}\right)^{25/16}}{\delta^2a_{t,1}^2}\quad(\text{Lemma \ref{lemma:variance-estimated-squared-residual}})\notag\\
            \leq&\frac{2^{25/16}C\bb_1^{25/16}\bb_2^{25/64}(\eta_tT)^{1/2}\log^2(\eta_tT)a_{t,1}^{25/64}}{\delta^2a_{t,1}^2}\quad(\text{using \eqref{lemma:expected-difference-inverse-probability_eq1}})\notag\\
            =&\frac{2^{25/16}C\bb_1^{25/16}\bb_2^{25/64}(\eta_tT)^{1/2}\log^2(\eta_tT)}{\delta^2a_{t,1}^{103/64}}\enspace.
        \end{align}
        Since $\delta<1/4$, we have $a_{t,0}\geq (1-\delta)a_{t,1}\geq \frac{3}{4}\cdot 7\bb_1^{-4}\bb_2^{-1}$, which leads to $2+\bb_1(\bb_2/6)^{1/4}a_{t,0}^{1/4}\leq 2\bb_1\bb_2^{1/4}a_{t,0}^{1/4}$.
        Hence, by similar arguments as in \eqref{lemma:expected-difference-inverse-probability_eq18}, we have
        \begin{align}\label{lemma:expected-difference-inverse-probability_eq19}
            \operatorname{Pr}\left(\frac{|\widehat{a}_{t,1}-a_{t,1}|}{a_{t,0}}\geq \delta\right)\leq&\frac{2^{25/16}C\bb_1^{25/16}\bb_2^{25/64}(\eta_tT)^{1/2}\log^2(\eta_tT)}{\delta^2a_{t,0}^{103/64}}\notag\\
            \leq&\frac{2^{25/16}C\bb_1^{25/16}\bb_2^{25/64}(\eta_tT)^{1/2}\log^2(\eta_tT)}{\delta^2(1-\delta)^{103/64}a_{t,1}^{103/64}}\quad(\text{since $a_{t,0}\geq (1-\delta)a_{t,1}$})\notag\\
            \leq&\frac{2^{41/16}C\bb_1^{25/16}\bb_2^{25/64}(\eta_tT)^{1/2}\log^2(\eta_tT)}{\delta^2a_{t,1}^{103/64}}\quad(\text{since $\delta<1/4$})\enspace.
        \end{align}
        By \eqref{lemma:expected-difference-inverse-probability_eq18}, Corollary \mainref{corollary:p-moment}*, Lemma \mainref{lemma:effect-of-p-regularization}*, and H{\"o}lder's inequality,
        \begin{align}\label{lemma:expected-difference-inverse-probability_eq20}
            \qquad D_2\leq&\e\left[\left(\frac{1}{p_{t+1}}+\frac{1}{\widebar{p}_{t+1}}\right)\mathbf{1}\left[|\widehat{a}_{t,0}-a_{t,0}|\geq \delta a_{t,1}\right]\right]\notag\\
            \leq&\underbrace{\left(\e\left[\frac{1}{p_{t+1}^4}\right]^{1/4}+\e\left[\frac{1}{\widebar{p}_{t+1}^4}\right]^{1/4}\right)}_{\text{Corollary \mainref{corollary:p-moment}*, Lemma \mainref{lemma:effect-of-p-regularization}*}}\underbrace{\left(\operatorname{Pr}\left(|\widehat{a}_{t,0}-a_{t,0}|\geq \delta a_{t,1}\right)\right)^{3/4}}_{\text{bound by \eqref{lemma:expected-difference-inverse-probability_eq18}}}\quad(\text{H{\"o}lder's inequality})\notag\\
            \leq&\left(\frac{2^{25/16}C\bb_1^{25/16}\bb_2^{25/64}(\eta_tT)^{1/2}\log^2(\eta_tT)}{\delta^2a_{t,1}^{103/64}}\right)^{3/4}\cdot 2\underbrace{\Big(2+\bb_1(\bb_2/6)^{1/4}a_{t,0}^{1/4}\Big)}_{\leq 2\bb_1\bb_2^{1/4}a_{t,0}^{1/4}}\notag\\
            \leq&2^{203/64}\bb_1\bb_2^{1/4}a_{t,0}^{1/4}\cdot \frac{C^{3/4}\bb_1^{75/64}\bb_2^{75/256}(\eta_tT)^{3/8}\log^{3/2}(\eta_tT)}{\delta^{3/2}a_{t,1}^{309/256}}\enspace.
        \end{align}
        By \eqref{lemma:expected-difference-inverse-probability_eq19}, Lemma \mainref{lemma:effect-of-p-regularization}*, Corollary \mainref{corollary:p-moment}*, and H{\"o}lder's inequality,
        \begin{align}\label{lemma:expected-difference-inverse-probability_eq21}
            \qquad D_3\leq&\e\left[\left(\frac{1}{p_{t+1}}+\frac{1}{\widebar{p}_{t+1}}\right)\mathbf{1}\left[|\widehat{a}_{t,1}-a_{t,1}|\geq \delta a_{t,0}\right]\right]\notag\\
            \leq&\underbrace{\left(\e\left[\frac{1}{p_{t+1}^4}\right]^{1/4}+\e\left[\frac{1}{\widebar{p}_{t+1}^4}\right]^{1/4}\right)}_{\text{Corollary \mainref{corollary:p-moment}*, Lemma \mainref{lemma:effect-of-p-regularization}*}}\underbrace{\left(\operatorname{Pr}\left(|\widehat{a}_{t,1}-a_{t,1}|\geq \delta a_{t,0}\right)\right)^{3/4}}_{\text{bound by \eqref{lemma:expected-difference-inverse-probability_eq19}}}\quad(\text{H{\"o}lder's inequality})\notag\\
            \leq&\left(\frac{2^{41/16}C\bb_1^{25/16}\bb_2^{25/64}(\eta_tT)^{1/2}\log^2(\eta_tT)}{\delta^2a_{t,1}^{103/64}}\right)^{3/4}\cdot 2\underbrace{\Big(2+\bb_1(\bb_2/6)^{1/4}a_{t,0}^{1/4}\Big)}_{\leq 2\bb_1\bb_2^{1/4}a_{t,0}^{1/4}}\notag\\
            \leq&2^{251/64}\bb_1\bb_2^{1/4}a_{t,0}^{1/4}\cdot \frac{C^{3/4}\bb_1^{75/64}\bb_2^{75/256}(\eta_tT)^{3/8}\log^{3/2}(\eta_tT)}{\delta^{3/2}a_{t,1}^{309/256}}\enspace.
        \end{align}
        \textbf{Step 3.3:} Combine the results in Step 3.1 and Step 3.2.\\[3mm]
        Since $\delta<1/4$, by the upper bounds in \eqref{lemma:expected-difference-inverse-probability_eq16}, \eqref{lemma:expected-difference-inverse-probability_eq17}, \eqref{lemma:expected-difference-inverse-probability_eq20} and \eqref{lemma:expected-difference-inverse-probability_eq21}, we obtain
        \begin{align}\label{lemma:expected-difference-inverse-probability_eq22}
            \qquad \qquad&\e\left[\left|\frac{1}{p_{t+1}}-\frac{1}{\widebar{p}_{t+1}}\right|\right]\notag\\
            \leq&\frac{8\bb_1\bb_2}{\bb_3}\delta+(2^{251/64}+2^{203/64})\cdot\bb_1\bb_2^{1/4}a_{t,0}^{1/4}\cdot \frac{C^{3/4}\bb_1^{75/64}\bb_2^{75/256}(\eta_tT)^{3/8}\log^{3/2}(\eta_tT)}{\delta^{3/2}a_{t,1}^{309/256}}\notag\\
            \leq&\frac{8\bb_1\bb_2}{\bb_3}\delta+2^{315/64}a_{t,0}^{1/4}C^{3/4}\bb_1^{139/64}\bb_2^{139/256}(\eta_tT)^{3/8}\log^{3/2}(\eta_tT)\delta^{-3/2}a_{t,1}^{-309/256}\notag\\
            \leq&\frac{8\bb_1\bb_2}{\bb_3}\delta+\underbrace{2^{315/64}(1-\delta)^{-1/4}}_{<33\text{ since }\delta<1/4}C^{3/4}\bb_1^{139/64}\bb_2^{139/256}(\eta_tT)^{3/8}\log^{3/2}(\eta_tT)\delta^{-3/2}a_{t,1}^{1/4-309/256}\notag\\
            \leq&\frac{8\bb_1\bb_2}{\bb_3}\delta+33C^{3/4}\bb_1^{139/64}\bb_2^{139/256}(\eta_tT)^{3/8}\log^{3/2}(\eta_tT)\delta^{-3/2}a_{t,1}^{-245/256}\enspace.
        \end{align}
        \textbf{Final Step:} Combine the results in Case 1, Case 2 and Case 3.\\[3mm]
        By the upper bounds in \eqref{lemma:expected-difference-inverse-probability_eq7}, \eqref{lemma:expected-difference-inverse-probability_eq15} and \eqref{lemma:expected-difference-inverse-probability_eq22}, we have
    \begin{align*}
        &\e\left[\left|\frac{1}{p_{t+1}}-\frac{1}{\widebar{p}_{t+1}}\right|\right]\\
        \leq&\widetilde{C}\max\Bigg\{(\eta_tT)^{1/4}\log(\eta_tT)a_{t,1}^{-103/128}+(\eta_tT)^{1/4}\log(\eta_tT)a_{t,1}^{-71/128}+(\eta_tT)^{3/8}\log^{3/2}(\eta_tT)\delta^{-3/2}a_{t,1}^{-245/256},\notag\\
        &(\eta_tT)^{1/4}\log(\eta_tT)a_{t,1}^{-71/128}+(\eta_tT)^{121/128}\log(\eta_tT)a_{t,1}^{-3/2}+(\eta_tT)^{235/256}\log^{3/2}(\eta_tT)\delta^{-3/2}a_{t,1}^{-3/2},\notag\\
        &\delta+(\eta_tT)^{3/8}\log^{3/2}(\eta_tT)\delta^{-3/2}a_{t,1}^{-245/256}\Bigg\}\enspace,
    \end{align*}
    where
    \begin{align*}
        \widetilde{C}=&\max\Bigg\{\frac{2^{-7/32}C^{1/2}\bb_1^{25/32}\bb_2^{281/128}}{\bb_3},\frac{2^{57/32}C^{1/2}\bb_1^{89/32}\bb_2^{185/128}}{\bb_3}\left(1+\frac{2}{\bb_3}\right),\\
        &\quad\quad\quad 2^{283/64}C^{3/4}\bb_1^{139/64}\bb_2^{139/256},~\frac{2^{89/32}C^{1/2}\bb_1^{121/32}\bb_2^{217/128}\varsigma^{89/128}}{\bb_3^{5/2}}\left(1+\frac{2}{\bb_3}\right),\\
        &\quad\quad\quad 2^{347/64}C^{3/4}\bb_1^{139/64}\bb_2^{139/256}\varsigma^{139/256},~\frac{8\bb_1\bb_2}{\bb_3},~33C^{3/4}\bb_1^{139/64}\bb_2^{139/256}\Bigg\}\enspace.
    \end{align*}
    This completes the proof.
\end{proof}

By imposing a threshold on $a_{t,1}$ or $a_{t,0}$ in Lemma~\ref{lemma:expected-difference-inverse-probability} and choosing $\delta>0$ appropriately, we can obtain the following corollary.

\begin{corollary}\label{corollary:expected-difference-inverse-probability}
	Suppose $T$ is sufficiently large. 
    Under Assumptions \mainref{assumption:moments}-\mainref{assumption:maximum-radius-strengthened}, and Condition \mainref{condition:sigmoid}, there exists a constant $K_2>0$ such that:
	\begin{enumerate}
		\item[(1)] For any $t\in[T]$ such that $a_{t,1}\geq  (\eta_tT)^{7/10}$, the following holds:
		\begin{align*}
			&\e\left[\left|\frac{1}{p_{t+1}}-\frac{1}{\widebar{p}_{t+1}}\right|\right]\leq K_2(\eta_tT)^{-1/20}\enspace.
		\end{align*}
		\item[(2)] For any $t\in[T]$ such that $a_{t,0}\geq (\eta_tT)^{7/10}$, the following holds:
		\begin{align*}
			&\e\left[\left|\frac{1}{1-p_{t+1}}-\frac{1}{1-\widebar{p}_{t+1}}\right|\right]\leq K_2(\eta_tT)^{-1/20}\enspace.
		\end{align*}
	\end{enumerate}
\end{corollary}

\begin{proof}
    We only prove the first part. 
    When $T$ is sufficiently large, $(\eta_tT)^{7/10}\geq (\eta_TT)^{7/10}\geq 7 \bb_1^{-4}\bb_2^{-1}$ by Assumption \mainref{assumption:maximum-radius-strengthened} and Lemma \ref{lemma:R}. For $0<\delta<1/4$, by Lemma \ref{lemma:expected-difference-inverse-probability}, we have
    \begin{align*}
        &\e\left[\left|\frac{1}{p_{t+1}}-\frac{1}{\widebar{p}_{t+1}}\right|\right]\\
        \leq&\widetilde{C}\max\Bigg\{\underbrace{(\eta_tT)^{1/4}\log(\eta_tT)a_{t,1}^{-103/128}+(\eta_tT)^{1/4}\log(\eta_tT)a_{t,1}^{-71/128}+(\eta_tT)^{3/8}\log^{3/2}(\eta_tT)\delta^{-3/2}a_{t,1}^{-245/256}}_{:=S_1},\notag\\
        &\underbrace{(\eta_tT)^{1/4}\log(\eta_tT)a_{t,1}^{-71/128}+(\eta_tT)^{121/128}\log(\eta_tT)a_{t,1}^{-3/2}+(\eta_tT)^{235/256}\log^{3/2}(\eta_tT)\delta^{-3/2}a_{t,1}^{-3/2}}_{:=S_2},\notag\\
        &\underbrace{\delta+(\eta_tT)^{3/8}\log^{3/2}(\eta_tT)\delta^{-3/2}a_{t,1}^{-245/256}}_{:=S_3}\Bigg\}\enspace.
    \end{align*}
    We then choose $\delta=(\eta_tT)^{-1/20}$ to balance the terms:
    \begin{align*}
        S_1\lesssim&\log(\eta_tT)(\eta_tT)^{\frac{1}{4}-\frac{103}{128}\cdot \frac{7}{10}}+\log(\eta_tT)(\eta_tT)^{\frac{1}{4}-\frac{71}{128}\cdot \frac{7}{10}}+\log^{3/2}(\eta_tT)(\eta_tT)^{\frac{3}{2}\cdot \frac{1}{20}+\frac{3}{8}-\frac{245}{256}\cdot \frac{7}{10}}\\
        \lesssim&(\eta_tT)^{-\frac{1}{20}}\enspace,\\
        S_2\lesssim&\log(\eta_tT)(\eta_tT)^{\frac{1}{4}-\frac{71}{128}\cdot \frac{7}{10}}+\log(\eta_tT)(\eta_tT)^{\frac{121}{128}-\frac{3}{2}\cdot \frac{7}{10}}+\log^{3/2}(\eta_tT)(\eta_tT)^{\frac{3}{2}\cdot \frac{1}{20}-\frac{3}{2}\cdot \frac{7}{10}+\frac{235}{256}}\\
        \lesssim&(\eta_tT)^{-\frac{1}{20}}\enspace,\\
        S_3\lesssim&(\eta_tT)^{-\frac{1}{20}}+\log^{3/2}(\eta_tT)(\eta_tT)^{\frac{3}{2}\cdot \frac{1}{20}+\frac{3}{8}-\frac{245}{256}\cdot \frac{7}{10}}\\
        \lesssim&(\eta_tT)^{-\frac{1}{20}}\enspace.
    \end{align*}
    Hence there exists a constant $K_2>0$ (independent of $t$ and $T$) such that
    \begin{align*}
        \e\left[\left|\frac{1}{p_{t+1}}-\frac{1}{\widebar{p}_{t+1}}\right|\right]\leq \widetilde{C}\max\{S_1,S_2,S_3\}\leq K_2(\eta_tT)^{-1/20}\enspace,
    \end{align*}
    which completes the proof.
\end{proof}

Based on Corollary \ref{corollary:expected-difference-inverse-probability}, the following lemma verifies the first convergence result in \eqref{simplified-stable-3}.

\begin{lemma}\label{lemma:stable-1}
    Denote $\Delta_{t,1}(\bv)=y_t(1)-\iprod{\xv_t,\bv}$.
    Under Assumptions \mainref{assumption:moments}-\mainref{assumption:maximum-radius-strengthened}, and Condition \mainref{condition:sigmoid}, the following holds:
    \begin{align*}
        \e\left[\frac{1}{T}\sum_{t=1}^T\Delta_{t,1}(\bv_t^*(1))^2\cdot\left|\frac{1}{p_t}-\frac{1}{\widebar{p}_t}\right|\right]\rightarrow 0\enspace.
    \end{align*}
\end{lemma}

\begin{proof}
    The high-level intuition of the proof is as follows.
    We split the sum according to whether $a_{t-1,1}$ is large or small.
    On the large set, we apply Corollary~\ref{corollary:expected-difference-inverse-probability} to bound $\e[| 1/p_t - 1/\widebar p_t |]$.
    On the small set, we control the deterministic part of the sum using a telescoping argument.

    We let $\mathcal{G}=\left\{2\leq t\leq T:a_{t-1,1}\geq (\eta_{t-1}T)^{7/10}\right\}$ and obtain the following decomposition:
    \begin{align}\label{lemma:stable-1_eq1}
        &\e\left(\frac{1}{T}\sum_{t=1}^{T}\Delta_{t,1}(\bv_t^*(1))^2\cdot \left|\frac{1}{p_t}-\frac{1}{\widebar{p}_t}\right|\right)\notag\\
        =&\underbrace{\frac{1}{T}\sum_{t\in\mathcal{G}^c}\Delta_{t,1}(\bv_t^*(1))^2\cdot \e\left[\left|\frac{1}{p_t}-\frac{1}{\widebar{p}_t}\right|\right]}_{:=S_1}+\underbrace{\frac{1}{T}\sum_{t\in{\mathcal{G}}}\Delta_{t,1}(\bv_t^*(1))^2\cdot \e\left[\left|\frac{1}{p_t}-\frac{1}{\widebar{p}_t}\right|\right]}_{:=S_2}\enspace.
    \end{align}
    \textbf{Step 1: }Derive an upper bound on $S_1$.\\[3mm]
    Before bounding $S_1$, we first provide an upper bound on $\iprod{\xv_t,\bv_t^*(1)}^2$.
    By Lemma \mainref{lemma:ridge-pred-exact-form}*, we have
    \begin{align}\label{lemma:stable-1_eq2}
        \iprod{\xv_t,\bv_t^*(1)}^2=&\left(\sum_{s=1}^{t-1}\Pi_{t,s}y_s(1)\right)^2\quad(\text{Lemma \mainref{lemma:ridge-pred-exact-form}*})\notag\\
        \leq&\left(\sum_{s=1}^{t-1}\Pi_{t,s}^2\right)\left(\sum_{s=1}^{t-1}y_s(1)^2\right)\quad(\text{Cauchy-Schwarz})\notag\\
        \intertext{The first term is bounded by $\gamma R_t^2((t-1)\vee\eta_t^{-1})^{-1}$ using Lemma \ref{lemma:deterministic-summation-1}-1. Applying the Cauchy-Schwarz inequality to the second term leads to}
        \lesssim&R_t^2((t-1)\vee\eta_t^{-1})^{-1}(t-1)^{1/2}\left(\sum_{s=1}^{t-1}y_s(1)^4\right)^{1/2}\notag\\
        \lesssim&R_t^2((t-1)\vee\eta_t^{-1})^{-1/2}T^{1/2}\quad(\text{Assumption \mainref{assumption:moments}})\notag\\
        \leq&R_t^2\eta_t^{1/2}T^{1/2}\notag\\
        =&R_t^{3/2}T^{1/4}\enspace.
    \end{align}
    In order to apply the telescoping technique, we first specify the elements in $\mathcal{G}^c$.
    Define $\eta_0=T^{-1/2}$, $\widehat{A}_{t_0}(1)=0$, and suppose $\mathcal{G}^c=[T]\,\setminus\,\mathcal{G}=\{t_1,\ldots,t_N\}$, where $1=t_1<\ldots<t_N$. 
    For any $t\in[T]$, by the definition of $\widehat{A}_{t}(1)$, we have
    \begin{align*}
        &\E{\widehat{A}_{t}(1)}-\E{\widehat{A}_{t-1}(1)}\\
        =&\e\left[\frac{\mathbf{1}[Z_{t}=1]}{p_{t}}\cdot(y_{t}(1)-\iprod{\xv_t,\bv_t(1)})^2\right]\\
        =&\e\left[(y_{t}(1)-\iprod{\xv_t,\bv_t(1)})^2\right]\quad(\text{by the law of iterated expectations})\\
        \intertext{Following the proof in Lemma \mainref{lemma:expectation-of-estores}, we can simplify this as}
        =&(y_{t}(1)-\iprod{\xv_t,\bv_t^*(1)})^2+\sum_{s=1}^{t-1}\Pi_{t,s}^2y_s(1)^2\e\left[\frac{1}{p_s}-1\right]\\
        \intertext{We bound the first term using inequality $(a+b)^2\leq 2a^2+2b^2$. Since $\Pi_{t,s}\lesssim R_tR_s((t-1)\vee \eta_t^{-1})^{-1}$ by Corollary \ref{corollary:pi} and $\e\left[\frac{1}{p_s}-1\right]\lesssim (\eta_sT)^{1/4}$ by Lemma \ref{lemma:p-power-moment}, we have}
        \lesssim&\underbrace{y_{t}^2(1)}_{\text{Assumption }  \mainref{assumption:moments}}+\underbrace{\iprod{\xv_t,\bv_t^*(1)}^2}_{\text{bound by }\eqref{lemma:stable-1_eq2}}+R_{t}\underbrace{((t-1)\vee\eta_{t}^{-1})^{-1}}_{\leq \eta_t}\sum_{s=1}^{t-1}|\Pi_{t,s}|R_sy_s(1)^2\cdot \eta_s^{1/4}T^{1/4}\\
        \lesssim&T^{1/2}+R_{t}^{3/2}T^{1/4}+R_{t}\eta_{t}T^{1/8}\sum_{s=1}^{t-1}|\Pi_{t,s}|R_s\cdot R_s^{-1/4}y_s(1)^2~~~~(\text{by Assumption \mainref{assumption:moments} and \eqref{lemma:stable-1_eq2}})\\
        \intertext{Since $R_t^{3/2}T^{1/4}=T^{1/2}R_t\cdot(TR_t^{-2})^{-1/4}\lesssim T^{1/2}R_t$ by Assumption \mainref{assumption:maximum-radius-strengthened} and Lemma \ref{lemma:R}, $R_s\leq R_t$ for $s\leq t$, we have}
        \lesssim&T^{1/2}R_t+T^{1/8-1/2}R_{t}^{1-1+1-1/4}\left(\sum_{s=1}^{t-1}\Pi_{t,s}^2\right)^{1/2}\left(\sum_{s=1}^{t-1}y_s(1)^4\right)^{1/2}~~~~(\text{Cauchy-Schwarz})\\
        \intertext{Since $\sum_{s=1}^{t-1}\Pi_{t,s}^2\lesssim R_t^2((t-1)\vee\eta_t^{-1})^{-1}\leq R_t^2\eta_t=R_tT^{-1/2}$ by Lemma \ref{lemma:deterministic-summation-1}-1 and $\sum_{s=1}^{t-1}y_s(1)^4\lesssim T$ by Assumption \mainref{assumption:moments}, we have}
        \lesssim&T^{1/2}R_t+T^{-3/8}R_{t}^{3/4}\cdot (R_tT^{-1/2})^{1/2}\cdot T^{1/2}\\
        \lesssim&T^{1/2}R_t+T^{1/2}R_t\cdot (TR_t^{-2})^{-1/8}\cdot T^{-1/2}\\
        \lesssim&T^{1/2}R_t\quad(\text{Lemma \ref{lemma:R} and Assumption \mainref{assumption:maximum-radius-strengthened}})\enspace.
    \end{align*}
    Hence there exists a constant $K_3>0$, such that
    \begin{align}\label{lemma:stable-1_eq3}
        \E{\widehat{A}_{t}(1)}-\E{\widehat{A}_{t-1}(1)}\leq K_3T^{1/2}R_t
    \end{align}
    for any $t\in[T]$. 
    Moreover, the previous proof also establishes that $(y_t(1)-\langle \xv_t,\bv_t^*(1)\rangle)^2\leq \E{\widehat{A}_t(1)}-\E{\widehat{A}_{t-1}(1)}$.
    By Lemma \ref{lemma:p-power-moment}, Corollary \mainref{corollary:p-moment}*, and Lemma \ref{lemma:squared-residuals-random}, we have $\e[1/p_t]\lesssim (\eta_tT)^{1/4}$ and $1/\widebar{p}_t\lesssim (\eta_tT)^{1/4}$.
    Hence, we have
    \begin{align}\label{lemma:stable-1_eq4}
        \qquad S_1\lesssim& T^{-1}\sum_{t\in\mathcal{G}^c}\left(\E{\widehat{A}_t(1)}-\E{\widehat{A}_{t-1}(1)}\right)(\eta_tT)^{1/4}\notag\\
        =&T^{-1}\sum_{k=1}^N\underbrace{\left(\E{\widehat{A}_{t_k}(1)}-\E{\widehat{A}_{t_{k}-1}(1)}\right)}_{\leq \E{\widehat{A}_{t_k}(1)}-\E{\widehat{A}_{t_{k-1}}(1)}\text{ by monotonicity}}(\eta_{t_k}T)^{1/4}\notag\\
        \leq&T^{-1}\sum_{k=1}^N\left(\E{\widehat{A}_{t_k}(1)}-\E{\widehat{A}_{t_{k-1}}(1)}\right)(\eta_{t_k}T)^{1/4}\notag\\
        \intertext{Rewriting the telescoping sum leads to}
        =&T^{-1}\sum_{k=1}^{N-1}\E{\widehat{A}_{t_k}(1)}\left((\eta_{t_k}T)^{1/4}-(\eta_{t_{k+1}}T)^{1/4}\right)\notag\\
        &+T^{-1}\E{\widehat{A}_{t_N}(1)}\cdot(\eta_{t_N}T)^{1/4}-T^{-1}\E{\widehat{A}_{t_0}(1)}\cdot(\eta_{t_1}T)^{1/4}\notag\\
        \intertext{Note that $\E{\widehat{A}_{t_0}(1)}=0$ by definition. By adding and subtracting $\E{\widehat{A}_{t_1-1}(1)},\ldots,\E{\widehat{A}_{t_N-1}(1)}$, we have}
        \leq&T^{-1}\sum_{k=1}^{N-1}\left(\eta_{t_k-1}^{-1}\cdot \eta_{t_k-1}\E{\widehat{A}_{t_k-1}(1)}+\left[\E{\widehat{A}_{t_k}(1)}-\E{\widehat{A}_{t_{k}-1}(1)}\right]\right)\notag\\
        &\times\left((\eta_{t_k}T)^{1/4}-(\eta_{t_{k+1}}T)^{1/4}\right)\notag\\
        &+T^{-1}\left(\eta_{t_N-1}^{-1}\cdot \eta_{t_N-1}\E{\widehat{A}_{t_N-1}(1)}+\left[\E{\widehat{A}_{t_N}(1)}-\E{\widehat{A}_{t_{N}-1}(1)}\right]\right)(\eta_{t_N}T)^{1/4}\notag\\
        \leq&T^{-1}\sum_{k=1}^{N-1}\left(\eta_{t_k-1}^{-1}\cdot \eta_{t_k-1}\E{\widehat{A}_{t_k-1}(1)}+K_3T^{1/2}R_{t_k}\right)\left((\eta_{t_k}T)^{1/4}-(\eta_{t_{k+1}}T)^{1/4}\right)\notag\\
        &+T^{-1}\left(\eta_{t_N-1}^{-1}\cdot \eta_{t_N-1}\E{\widehat{A}_{t_N-1}(1)}+K_3T^{1/2}R_{t_N}\right)(\eta_{t_N}T)^{1/4}\quad(\text{by \eqref{lemma:stable-1_eq3}})\notag\\
        \intertext{Rearranging the terms in the summation leads to}
        =&K_3T^{-1+1/2+1/2}\left[\sum_{k=1}^{N-1}(\eta_{t_k}T)^{-1}\left((\eta_{t_k}T)^{1/4}-(\eta_{t_{k+1}}T)^{1/4}\right)+(\eta_{t_N}T)^{-1}(\eta_{t_N}T)^{1/4}\right]\notag\\
        &+T^{-1}\Bigg[\sum_{k=1}^{N-1}\eta_{t_k-1}^{-1}\cdot \eta_{t_k-1}\E{\widehat{A}_{t_k-1}(1)}\left((\eta_{t_k}T)^{1/4}-(\eta_{t_{k+1}}T)^{1/4}\right)\notag\\
        &\qquad\qquad +\eta_{t_N-1}^{-1}\cdot \eta_{t_N-1}\E{\widehat{A}_{t_N-1}(1)}\cdot (\eta_{t_{N}}T)^{1/4}\Bigg]\notag\\
        \intertext{By the definition of $\mathcal{G}^c$, for any $k=1,\ldots,N$, it holds that $\eta_{t_k-1}^{-1}\cdot\eta_{t_k-1}\E{\widehat{A}_{t_k-1}(1)}=\eta_{t_k-1}^{-1}a_{t_k-1,1}\leq \eta_{t_k-1}^{-1}(\eta_{t_k-1}T)^{7/10}\leq \eta_{t_k}^{-1}(\eta_{t_k}T)^{7/10}$ since $\eta_{t_k-1}\geq \eta_{t_k}$. Hence we have}
        \leq&K_3\left[\sum_{k=1}^{N-1}\left((\eta_{t_k}T)^{1/4}\right)^{-4}\left((\eta_{t_k}T)^{1/4}-(\eta_{t_{k+1}}T)^{1/4}\right)+(\eta_{t_N}T)^{-3/4}\right]\notag\\
        &+T^{-1}\left[\sum_{k=1}^{N-1}\eta_{t_k}^{-1}(\eta_{t_k}T)^{7/10}\left((\eta_{t_k}T)^{1/4}-(\eta_{t_{k+1}}T)^{1/4}\right)+\eta_{t_N}^{-1}(\eta_{t_N}T)^{7/10}(\eta_{t_N}T)^{1/4}\right]\notag\\
        \lesssim&\sum_{k=1}^{N-1}\left((\eta_{t_k}T)^{1/4}\right)^{-4}\left((\eta_{t_k}T)^{1/4}-(\eta_{t_{k+1}}T)^{1/4}\right)+\sum_{k=1}^{N-1}\left((\eta_{t_k}T)^{1/4}\right)^{-6/5}\left((\eta_{t_k}T)^{1/4}-(\eta_{t_{k+1}}T)^{1/4}\right)\notag\\
        &+(\eta_{t_N}T)^{-3/4}+(\eta_{t_N}T)^{-1+7/10+1/4}\notag\\
        \intertext{Since $\eta_{t_1}T\geq\ldots \geq \eta_{t_N}T$, by an integral comparison argument we have}
        \leq&\int_{(\eta_{t_{N}}T)^{1/4}}^{(\eta_{t_{1}}T)^{1/4}}x^{-4}\mathrm{d}x+\int_{(\eta_{t_{N}}T)^{1/4}}^{(\eta_{t_{1}}T)^{1/4}}x^{-6/5}\mathrm{d}x+(\eta_{T}T)^{-3/4}+(\eta_{T}T)^{-1/20}\notag\\
        \lesssim&o(1)+(\eta_{t_N}T)^{-\frac{1}{4}\cdot 3}+(\eta_{t_N}T)^{-\frac{1}{4}\cdot \frac{1}{5}}~~~~(\text{Assumption \mainref{assumption:maximum-radius-strengthened} and Lemma \ref{lemma:R}})\notag\\
        \leq&o(1)+(\eta_TT)^{-3/4}+(\eta_TT)^{-1/20}\notag\\
        \rightarrow&0~~~~(\text{Assumption \mainref{assumption:maximum-radius-strengthened} and Lemma \ref{lemma:R}})\enspace.
    \end{align}
    \textbf{Step 2: }Derive an upper bound on $S_2$.\\[3mm]
    By the definition of the set $\mathcal{G}$ and Corollary \ref{corollary:expected-difference-inverse-probability}-1, we have
    \begin{align*}
        \max_{t\in\mathcal{G}}\e\left[\left|\frac{1}{p_t}-\frac{1}{\widebar{p}_t}\right|\right]\lesssim&\max_{t\in\mathcal{G}}(\eta_{t-1}T)^{-1/20}\leq (\eta_TT)^{-1/20}\enspace.
    \end{align*}
    Hence by Corollary \ref{corollary:squared-residual-deterministic}, we can bound $S_2$ as:
    \begin{align}\label{lemma:stable-1_eq5}
        S_2=&\frac{1}{T}\sum_{t\in\mathcal{G}}\Delta_{t,1}(\bv_t^*(1))^2\cdot \e\left[\left|\frac{1}{p_t}-\frac{1}{\widebar{p}_t}\right|\right]\notag\\
        \lesssim&T^{-1}\underbrace{\sum_{t=1}^T\Delta_{t,1}(\bv_t^*(1))^2}_{\text{Corollary \ref{corollary:squared-residual-deterministic}}}\cdot (\eta_TT)^{-1/20}\notag\\
        \lesssim&T^{-1}\cdot T\cdot (\eta_TT)^{-1/20}\notag\\
        \leq&(\eta_TT)^{-1/20}\notag\\
        \rightarrow&0~~~~(\text{Assumption \mainref{assumption:maximum-radius-strengthened} and Lemma \ref{lemma:R}})\enspace.
    \end{align}
    \textbf{Step 3: }Combine the results in Step 1 and Step 2.\\[3mm]
    By \eqref{lemma:stable-1_eq1}, \eqref{lemma:stable-1_eq4}, \eqref{lemma:stable-1_eq5}, we have
    \begin{align}\label{thm5.1_eq10}
        \e\left[\frac{1}{T}\sum_{t=1}^{T}\Delta_{t,1}(\bv_t^*(1))^2\cdot \left|\frac{1}{p_t}-\frac{1}{\widebar{p}_t}\right|\right]=S_1+S_2\rightarrow 0\enspace.
    \end{align}
    Hence the result is proved.
\end{proof}

The remaining convergence results in \eqref{simplified-stable-3} follow from bounding the tracking-error terms and the inverse-probability moments.

\begin{lemma}\label{lemma:stable-2}
    Under Assumptions \mainref{assumption:moments}-\mainref{assumption:maximum-radius-strengthened}, and Condition \mainref{condition:sigmoid}, the following holds:
    \begin{align*}
        \e\left[\frac{1}{T}\sum_{t=1}^T\iprod{\xv_t,\bv_t(1)-\bv^*_t(1)}^2\cdot\left|\frac{1}{p_t}-\frac{1}{\widebar{p}_t}\right|\right]\rightarrow 0\enspace.
    \end{align*}
\end{lemma}

\begin{proof}
    By Corollary \ref{corollary:expected-tracking} and Lemma \ref{lemma:p-power-moment}, we have 
    \begin{align*}
        &\e\Bigg[\frac{1}{T}\sum_{t=1}^T\iprod{\xv_t,\bv_t(1)-\bv^*_t(1)}^2\cdot\left|\frac{1}{p_t}-\frac{1}{\widebar{p}_t}\right|\Bigg]\\
        =&T^{-1}\e\Bigg[\sum_{t=1}^TR_t^{-1/2}\iprod{\xv_t,\bv_t(1)-\bv^*_t(1)}^2\cdot R_t^{1/2}\left|\frac{1}{p_t}-\frac{1}{\widebar{p}_t}\right|\Bigg]\\
        \leq&T^{-1}\e\left[\sum_{t=1}^TR_t^{-1}\iprod{\xv_t,\bv_t(1)-\bv^*_t(1)}^4\right]^{1/2}\e\left[\sum_{t=1}^TR_t\left(\frac{1}{p_t}-\frac{1}{\widebar{p}_t}\right)^2\right]^{1/2}\quad(\text{Cauchy-Schwarz})\\
        \intertext{Since $\e[1/p_t^2]$ and $\e[1/\widebar{p}_t^2]$ can be bounded using Lemma \ref{lemma:p-power-moment}, by the inequality $(a+b)^2\leq 2a^2+2b^2$, we have}
        \lesssim&T^{-1}\cdot (T^{3/8}R_T^{1/4})^{1/2}\left[\sum_{t=1}^TR_t(\eta_tT)^{1/2}\right]^{1/2}\\
        =&T^{-1+3/16}R_T^{1/8}\left[\sum_{t=1}^TR_t^{1-1/2}T^{1/4}\right]^{1/2}\\
        \leq&T^{-1+3/16+1/2+1/8}R_T^{1/8+1/4}\quad(\text{since $R_t\leq R_T$})\\
        =&T^{-3/16}R_T^{3/8}\\
        =&o(1)~~~~(\text{Assumption \mainref{assumption:maximum-radius-strengthened} and Lemma \ref{lemma:R}})\enspace,
    \end{align*}
    which completes the proof.
\end{proof}

\begin{lemma}\label{lemma:stable-3}
    Under Assumptions \mainref{assumption:moments}-\mainref{assumption:maximum-radius-strengthened}, and Condition \mainref{condition:sigmoid}, the following holds:
    \begin{align*}
        \e\left[\frac{1}{T}\sum_{t=1}^T\iprod{\xv_t,\bv_t(1)-\bv^*_t(1)}^2\cdot\left|\frac{1}{\widebar{p}_t}-1\right|\right]\rightarrow 0\enspace.
    \end{align*}
\end{lemma}

\begin{proof}
    Note that the arguments in Lemma \ref{lemma:stable-2} still apply to this lemma, which completes the proof.
\end{proof}

\begin{lemma}\label{lemma:stable-4}
    Denote $\Delta_{t,1}(\bv)=y_t(1)-\iprod{\xv_t,\bv}$.
    Under Assumptions \mainref{assumption:moments}-\mainref{assumption:maximum-radius-strengthened}, and Condition \mainref{condition:sigmoid}, the following holds:
    \begin{align*}
        \e\left[\frac{1}{T}\sum_{t=1}^T|\Delta_{t,1}(\bv_t^*(1))|\cdot|\iprod{\xv_t,\bv_t(1)-\bv^*_t(1)}|\cdot\left|\frac{1}{p_t}-\frac{1}{\widebar{p}_t}\right|\right]\rightarrow 0\enspace.
    \end{align*}
\end{lemma}

\begin{proof}
    Applying the Cauchy-Schwarz inequality to the product term and using Lemma \ref{lemma:stable-1}, Lemma \ref{lemma:stable-2} lead to
    \begin{align*}
        &\e\left[\frac{1}{T}\sum_{t=1}^T|\Delta_{t,1}(\bv_t^*(1))|\cdot|\iprod{\xv_t,\bv_t(1)-\bv^*_t(1)}|\cdot\left|\frac{1}{p_t}-\frac{1}{\widebar{p}_t}\right|\right]\\
        \leq&\underbrace{\left(\e\left[\frac{1}{T}\sum_{t=1}^T\Delta_{t,1}(\bv_t^*(1))^2\cdot\left|\frac{1}{p_t}-\frac{1}{\widebar{p}_t}\right|\right]\right)^{1/2}}_{\text{Lemma } \ref{lemma:stable-1}}\cdot\underbrace{\left(\e\left[\frac{1}{T}\sum_{t=1}^T\iprod{\xv_t,\bv_t(1)-\bv^*_t(1)}^2\cdot\left|\frac{1}{p_t}-\frac{1}{\widebar{p}_t}\right|\right]\right)^{1/2}}_{\text{Lemma } \ref{lemma:stable-2}}\\
        \rightarrow&0\enspace,
    \end{align*}
    which completes the proof.
\end{proof}

\begin{lemma}\label{lemma:stable-5}
    Denote $\Delta_{t,1}(\bv)=y_t(1)-\iprod{\xv_t,\bv}$.
    Under Assumptions \mainref{assumption:moments}-\mainref{assumption:maximum-radius-strengthened}, and Condition \mainref{condition:sigmoid}, the following holds:
    \begin{align*}
        \e\left[\frac{1}{T}\sum_{t=1}^T|\Delta_{t,1}(\bv_t^*(1))|\cdot|\iprod{\xv_t,\bv_t(1)-\bv^*_t(1)}|\cdot\left|\frac{1}{\widebar{p}_t}-1\right|\right]\rightarrow 0\enspace.
    \end{align*}
\end{lemma}

\begin{proof}
    By Lemma \mainref{lemma:effect-of-p-regularization}*, Lemma \ref{lemma:squared-residuals-random}, and the Cauchy-Schwarz inequality, we have
    \begin{align}\label{thm5.1_eq19}
        &\e\Bigg[\frac{1}{T}\sum_{t=1}^T|\Delta_{t,1}(\bv_t^*(1))|\cdot|\iprod{\xv_t,\bv_t(1)-\bv^*_t(1)}|\cdot\left|\frac{1}{\widebar{p}_t}-1\right|\Bigg]\notag\\
        \lesssim&\e\left[\frac{1}{T}\sum_{t=1}^T|\Delta_{t,1}(\bv_t^*(1))|\cdot|\iprod{\xv_t,\bv_t(1)-\bv^*_t(1)}|\cdot\eta_t^{1/4}T^{1/4}\right]\quad(\text{Lemma \mainref{lemma:effect-of-p-regularization}* and Lemma \ref{lemma:squared-residuals-random}})\notag\\
        \lesssim&T^{1/8}\e\left[\frac{1}{T}\sum_{t=1}^TR_t^{-1/4}|\Delta_{t,1}(\bv_t^*(1))|\cdot|\iprod{\xv_t,\bv_t(1)-\bv^*_t(1)}|\right]\notag\\
        \lesssim&T^{-1+1/8}\left(\sum_{t=1}^T\Delta_{t,1}(\bv_t^*(1))^2\right)^{1/2}\left(\e\left[\sum_{t=1}^TR_t^{-1/2}\iprod{\xv_t,\bv_t(1)-\bv^*_t(1)}^2\right]\right)^{1/2}~~~~(\text{Cauchy-Schwarz})\notag\\
        \intertext{$\sum_{t=1}^T\Delta_{t,1}(\bv_t^*(1))^2$ is of order $\bigO{T}$ using Corollary \ref{corollary:squared-residual-deterministic}. $\e\left[\sum_{t=1}^TR_t^{-1/2}\iprod{\xv_t,\bv_t(1)-\bv^*_t(1)}^2\right]$ can be simplified using Lemma \ref{lemma:tracking-term}. Hence we have}
        \lesssim&T^{-1+1/8+1/2}\Bigg(\sum_{t=1}^TR_t^{-1/2}\sum_{s=1}^{t-1}\Pi_{t,s}^2y_s(1)^2\underbrace{\e\left[\frac{1}{p_s}-1\right]}_{\text{Lemma \ref{lemma:p-power-moment}}}\Bigg)^{1/2}\notag\\
        \lesssim&T^{-3/8}\left(\sum_{t=1}^TR_t^{-1/2}\sum_{s=1}^{t-1}\Pi_{t,s}^2y_s(1)^2\cdot(\eta_sT)^{1/4}\right)^{1/2}\notag\\
        \lesssim&T^{-3/8+1/16}\underbrace{\left(\sum_{t=1}^TR_t^{-1/2}\sum_{s=1}^{t-1}R_s^{-1/4}\Pi_{t,s}^2y_s(1)^2\right)^{1/2}}_{\text{Lemma \ref{lemma:deterministic-summation-1}-3}}\notag\\
        \lesssim&T^{-5/16}\cdot \left(R_T^{3/2-1/2-1/4}T^{1/4}\right)^{1/2}\notag\\
        \lesssim&(TR_T^{-2})^{-3/16}\notag\\
        \rightarrow&0\quad(\text{Assumption \mainref{assumption:maximum-radius-strengthened} and Lemma \ref{lemma:R}})\enspace,
    \end{align}
    which completes the proof.
\end{proof}

Using Lemma \ref{lemma:stable-1}-\ref{lemma:stable-5}, we establish the first convergence result in \eqref{simplified-stable-1} by the following lemma.

\begin{lemma}\label{lemma:stable-variance-1}
    For $k\in\{0,1\}$, denote $\Delta_{t,k}(\bv)=y_t(k)-\iprod{\xv_t,\bv}$.
    Under Assumptions \mainref{assumption:moments}-\mainref{assumption:maximum-radius-strengthened}, and Condition \mainref{condition:sigmoid}, the following holds:
    \begin{align*}
        &\frac{1}{T}\left[\sum_{t=1}^T\Delta_{t,1}(\bv_t(1))^2\cdot \left(\frac{1}{p_t}-1\right)-\e\left[\sum_{t=1}^T\Delta_{t,1}(\bv_t(1))^2\cdot \left(\frac{1}{p_t}-1\right)\right]\right]\xrightarrow{p}0\enspace,\\
        &\frac{1}{T}\left[\sum_{t=1}^T\Delta_{t,0}(\bv_t(0))^2\cdot \left(\frac{1}{1-p_t}-1\right)-\e\left[\sum_{t=1}^T\Delta_{t,0}(\bv_t(0))^2\cdot \left(\frac{1}{1-p_t}-1\right)\right]\right]\xrightarrow{p}0\enspace.
    \end{align*}
\end{lemma}

\begin{proof}
    By Lemma \ref{lemma:stable-1} and Lemma \ref{lemma:markov}, we have established that
    \begin{align*}
        \frac{1}{T}\sum_{t=1}^T\Delta_{t,1}(\bv_t^*(1))^2\cdot\left(\frac{1}{p_t}-\frac{1}{\widebar{p}_t}\right)-\e\left[\frac{1}{T}\sum_{t=1}^T\Delta_{t,1}(\bv_t^*(1))^2\cdot\left(\frac{1}{p_t}-\frac{1}{\widebar{p}_t}\right)\right]\xrightarrow{p}0\enspace.
    \end{align*}
    By Lemma \ref{lemma:markov}, Lemma \ref{lemma:stable-2}, Lemma \ref{lemma:stable-3}, Lemma \ref{lemma:stable-4}, and Lemma \ref{lemma:stable-5}, we can similarly prove the remaining four convergence results, which together imply the following:
    \begin{align*}
        \frac{1}{T}\left[\sum_{t=1}^T\Delta_{t,1}(\bv_t(1))^2\cdot \left(\frac{1}{p_t}-1\right)-\e\left[\sum_{t=1}^T\Delta_{t,1}(\bv_t(1))^2\cdot \left(\frac{1}{p_t}-1\right)\right]\right]\xrightarrow{p}0\enspace.
    \end{align*}
    The second convergence result is similarly verified due to the symmetry.
\end{proof}

We can similarly verify the final convergence result in \eqref{simplified-stable-1} in the following lemma.

\begin{lemma}\label{lemma:stable-variance-2}
    For $k\in\{0,1\}$, denote $\Delta_{t,k}(\bv)=y_t(k)-\iprod{\xv_t,\bv}$.
    Under Assumptions \mainref{assumption:moments}-\mainref{assumption:maximum-radius-strengthened}, and Condition \mainref{condition:sigmoid}, the following holds:
    \begin{align*}
        \frac{1}{T}\left[\sum_{t=1}^T\Delta_{t,1}(\bv_t(1))\cdot\Delta_{t,0}(\bv_t(0))-\e\left[\sum_{t=1}^T\Delta_{t,1}(\bv_t(1))\cdot\Delta_{t,0}(\bv_t(0))\right]\right]\xrightarrow{p}0\enspace.
    \end{align*}
\end{lemma}

\begin{proof}
    We prove the result using a method similar to that in Lemma \ref{lemma:stable-variance-1}.
    First, we decompose the random term:
    \begin{align*}
        &\frac{1}{T}\sum_{t=1}^T\Delta_{t,1}(\bv_t(1))\cdot \Delta_{t,0}(\bv_t(0))\\
        =&\frac{1}{T}\sum_{t=1}^T\Delta_{t,1}(\bv_t(1))\cdot \Delta_{t,0}(\bv_t^*(0))+\frac{1}{T}\sum_{t=1}^T\Delta_{t,1}(\bv_t(1))\cdot\iprod{\xv_t,\bv_t^*(0)-\bv_t(0)}\\
        =&\frac{1}{T}\sum_{t=1}^T\Delta_{t,1}(\bv_t^*(1))\cdot \Delta_{t,0}(\bv_t^*(0))+\frac{1}{T}\sum_{t=1}^T\Delta_{t,1}(\bv_t(1))\cdot\iprod{\xv_t,\bv_t^*(0)-\bv_t(0)}\\
        &+\frac{1}{T}\sum_{t=1}^T\Delta_{t,0}(\bv_t^*(0))\cdot\iprod{\xv_t,\bv_t^*(1)-\bv_t(1)}\enspace.
    \end{align*}
    Note that the first term is nonrandom. 
    By Lemma \ref{lemma:markov}, it suffices to show that
    \begin{align*}
        &\e\left[\frac{1}{T}\sum_{t=1}^T\left|\Delta_{t,1}(\bv_t(1))\right|\cdot|\iprod{\xv_t,\bv_t^*(0)-\bv_t(0)}|\right]\rightarrow 0\enspace,\\
        &\e\left[\frac{1}{T}\sum_{t=1}^T\left|\Delta_{t,0}(\bv_t^*(0))\right|\cdot|\iprod{\xv_t,\bv_t^*(1)-\bv_t(1)}|\right]\rightarrow 0\enspace.
    \end{align*}
    By Lemma \ref{lemma:tracking-term}, Lemma \ref{lemma:deterministic-summation-1} and Lemma \ref{lemma:p-power-moment}, we have
    \begin{align}\label{lemma:stable-variance-2_eq1}
        \e\left[\sum_{t=1}^T\langle \xv_t,\bv_t^*(1)-\bv_t(1)\rangle^2\right]=&\sum_{t=1}^T\sum_{s=1}^{t-1}\Pi_{t,s}^2y_s(1)^2\e\left[\frac{1}{p_s}-1\right]\quad(\text{Lemma \ref{lemma:tracking-term}})\notag\\
        \lesssim&T^{1/8}\sum_{t=1}^T\sum_{s=1}^{t-1}R_s^{-1/4}\Pi_{t,s}^2y_s(1)^2\quad(\text{Lemma \ref{lemma:p-power-moment}})\notag\\
        \lesssim&T^{1/8}\cdot R_T^{3/2-1/4}T^{1/4}\quad(\text{Lemma \ref{lemma:deterministic-summation-1}-3})\notag\\
        =&T^{3/8}R_T^{5/4}\enspace.
    \end{align}
    A similar result can be derived for $\e\left[\sum_{t=1}^T\langle \xv_t,\bv_t^*(0)-\bv_t(0)\rangle^2\right]$.
    Hence by \eqref{lemma:stable-variance-2_eq1} and Lemma \ref{lemma:squared-residuals-random}, we have
    \begin{align*}
        &\e\left[\frac{1}{T}\sum_{t=1}^T\left|\Delta_{t,1}(\bv_t(1))\right|\cdot|\iprod{\xv_t,\bv_t^*(0)-\bv_t(0)}|\right]\notag\\
        \leq&\frac{1}{T}\underbrace{\left(\e\left[\sum_{t=1}^T\Delta_{t,1}(\bv_t(1))^2\right]\right)^{1/2}}_{\text{Lemma \ref{lemma:squared-residuals-random}}}\cdot\underbrace{\left(\e\left[\sum_{t=1}^T\iprod{\xv_t,\bv_t^*(0)-\bv_t(0)}^2\right]\right)^{1/2}}_{\text{bound by \eqref{lemma:stable-variance-2_eq1}}}~~~~(\text{Cauchy-Schwarz})\notag\\
        \lesssim&T^{-1}\cdot T^{1/2}\cdot (T^{3/8}R_T^{5/4})^{1/2}\notag\\
        =&(TR_T^{-2})^{-5/16}\notag\\
        \rightarrow&0\quad(\text{Assumption \mainref{assumption:maximum-radius-strengthened} and Lemma \ref{lemma:R}})\enspace.
    \end{align*}
    By \eqref{lemma:stable-variance-2_eq1} and Corollary \ref{corollary:squared-residual-deterministic}, we have
    \begin{align*}
        &\e\left[\frac{1}{T}\sum_{t=1}^T\left|\Delta_{t,0}(\bv_t^*(0))\right|\cdot|\iprod{\xv_t,\bv_t^*(1)-\bv_t(1)}|\right]\notag\\
        \leq&\frac{1}{T}\underbrace{\left(\e\left[\sum_{t=1}^T\Delta_{t,0}(\bv_t^*(0))^2\right]\right)^{1/2}}_{\text{Corollary \ref{corollary:squared-residual-deterministic}}}\cdot\underbrace{\left(\e\left[\sum_{t=1}^T\iprod{\xv_t,\bv_t^*(1)-\bv_t(1)}^2\right]\right)^{1/2}}_{\text{bound by \eqref{lemma:stable-variance-2_eq1}}}~~~~(\text{Cauchy-Schwarz})\notag\\
        \lesssim&T^{-1}\cdot T^{1/2}\cdot (T^{3/8}R_T^{5/4})^{1/2}\notag\\
        =&(TR_T^{-2})^{-5/16}\notag\\
        \rightarrow&0\quad(\text{Assumption \mainref{assumption:maximum-radius-strengthened} and Lemma \ref{lemma:R}})\enspace.
    \end{align*}
    Combining these two convergence results with Lemma \ref{lemma:markov} completes the proof.
\end{proof}

By Lemma \ref{lemma:stable-variance-1} and Lemma \ref{lemma:stable-variance-2}, we can verify the conditional variance convergence in the following lemma.

\begin{lemma}\label{lemma:stable-variance}
    Under Assumptions \mainref{assumption:moments}-\mainref{assumption:bounded-correlation} and Condition \mainref{condition:sigmoid}, it holds that $V_{T}^2\xrightarrow{p}1$.
\end{lemma}

\begin{proof}
    By Lemma \ref{lemma:stable-variance-1}, Lemma \ref{lemma:stable-variance-2} and the simplified form of $V_T^2$ derived at the beginning of the section, the result is proved.
\end{proof}

\subsubsection{Conditional Lyapunov Condition}
In this section, we verify the conditional Lyapunov condition for the martingale central limit theorem.
We begin by deriving a simplified form of the conditional Lyapunov condition.
For $k\in\{0,1\}$, denote $\Delta_{t,k}(\bv)=y_t(k)-\iprod{\xv_t,\bv}$.
Setting $\delta = 2$ and using the inequality $(a+b)^4\leq 8a^4+8b^4$, we obtain:
\begin{align*}
    \sum_{t=1}^T\E{X_{t}^4|\filt_{t-1}}\lesssim&\frac{1}{T^4(\Var{\eate})^2}\sum_{t=1}^T\e\left[\Delta_{t,1}(\bv_t(1))^4\cdot \left(\frac{ \indicator{Z_t=1} }{p_t}-1\right)^4\Big|\filt_{t-1}\right]\\
    &+\frac{1}{T^4(\Var{\eate})^2}\sum_{t=1}^T\e\left[\Delta_{t,0}(\bv_t(0))^4\cdot \left(\frac{ \indicator{Z_t=0} }{1-p_t}-1\right)^4\Big|\filt_{t-1}\right]\\
    \lesssim&\frac{1}{T^4(\Var{\eate})^2}\sum_{t=1}^T\left(\Delta_{t,1}(\bv_t(1))^4\cdot \frac{1}{p_t^3}+\Delta_{t,0}(\bv_t(0))^4\cdot \frac{1}{(1-p_t)^3}\right)\enspace.
\end{align*}

The non-superefficiency condition (Corollary \mainref{corollary:non-superefficiency}) guarantees that $\Var{\eate}=\Omega(T^{-1})$.
In order to prove that $\sum_{t=1}^T\E{X_{t}^4|\filt_{t-1}}$ converges to 0 in probability, it suffices to verify that
\begin{align*}
    &\frac{1}{T^2}\sum_{t=1}^T\e\left[\Delta_{t,1}(\bv_t(1))^4\cdot \frac{1}{p_t^3}\right]\rightarrow 0\quad \text{and}\quad\frac{1}{T^2}\sum_{t=1}^T\e\left[\Delta_{t,0}(\bv_t(0))^4\cdot \frac{1}{(1-p_t)^3}\right]\rightarrow 0\enspace.
\end{align*}
By symmetry between the treated ($k=1$) and control groups ($k=0$), we only prove the first convergence result.
Using the inequality $(a+b)^4\leq 8a^4+8b^4$ and the equality $\Delta_{t,1}(\bv_t(1))=\Delta_{t,1}(\bv_t^*(1))+\iprod{ \xv_t, \bv_t^*(1)-\bv_t(1) }$, it suffices to show that
\begin{align}\label{simplified-lyapunov}
    \qquad \frac{1}{T^2}\sum_{t=1}^T\e\left[\Delta_{t,1}(\bv_t^*(1))^4\cdot \frac{1}{p_t^3}\right]\rightarrow 0~ \text{and}~\frac{1}{T^2}\sum_{t=1}^T\e\left[\iprod{ \xv_t, \bv_t(1)-\bv_t^*(1) } ^4\cdot \frac{1}{p_t^3}\right]\rightarrow 0\enspace.
\end{align}
The two convergence results in \eqref{simplified-lyapunov} are established in Lemma \ref{lemma:Lyapunov-1} and Lemma \ref{lemma:Lyapunov-2} separately.

\begin{lemma}\label{lemma:Lyapunov-1}
    Denote $\Delta_{t,1}(\bv)=y_t(1)-\iprod{\xv_t,\bv}$.
    Under Assumptions \mainref{assumption:moments}-\mainref{assumption:maximum-radius-strengthened}, and Condition \mainref{condition:sigmoid}, the following holds:
    \begin{align*}
        \frac{1}{T^2}\sum_{t=1}^T\e\left[\Delta_{t,1}(\bv_t^*(1))^4\cdot \frac{1}{p_t^3}\right]\rightarrow 0\enspace.
    \end{align*}
\end{lemma}

\begin{proof}
    By Lemma \ref{lemma:p-power-moment} and Corollary \ref{corollary:fourth-moment-deterministic}, we obtain
    \begin{align*}
        &\frac{1}{T^2}\sum_{t=1}^T\e\left[\Delta_{t,1}(\bv_t^*(1))^4\cdot \frac{1}{p_t^3}\right]\\
        \lesssim&T^{-2}\sum_{t=1}^T\Delta_{t,1}(\bv_t^*(1))^4\cdot(\eta_tT)^{3/4}\quad(\text{Lemma \ref{lemma:p-power-moment}})\\
        \leq&T^{-2+3/4}\sum_{t=1}^T\eta_t^{-1/4}\cdot\eta_t\cdot\Delta_{t,1}(\bv_t^*(1))^4\\
        \leq&T^{-5/4}\eta_T^{-1/4}\cdot R_TT^{1/2}\quad(\text{using Corollary \ref{corollary:fourth-moment-deterministic} and monotonicity $\eta_T\leq \eta_t$})\\
        =&(TR_T^{-2})^{-5/8}\\
        \rightarrow&0\quad(\text{Assumption \mainref{assumption:maximum-radius-strengthened} and Lemma \ref{lemma:R}})\enspace,
    \end{align*}
    which completes the proof.
\end{proof}

\begin{lemma}\label{lemma:Lyapunov-2}
    Under Assumptions \mainref{assumption:moments}-\mainref{assumption:maximum-radius-strengthened}, and Condition \mainref{condition:sigmoid}, the following holds:
    \begin{align*}
        \frac{1}{T^2}\sum_{t=1}^T\e\left[\iprod{ \xv_t, \bv_t(1)-\bv_t^*(1) } ^4\cdot \frac{1}{p_t^3}\right]\rightarrow 0\enspace.
    \end{align*}
\end{lemma}

\begin{proof}
    By Proposition \mainref{prop:as-inv-prob-bound} and Corollary \ref{corollary:expected-tracking}, we have
    \begin{align*}
        &\frac{1}{T^2}\sum_{t=1}^T\e\left[\frac{1}{p_t^3}\iprod{\xv_t,\bv_t(1)-\bv_t^*(1)}^4\right]\\
        \lesssim&T^{-2+21/26}\sum_{t=1}^T\e\left[R_t^{-6/11}\iprod{\xv_t,\bv_t(1)-\bv_t^*(1)}^4\right]\quad(\text{Proposition \mainref{prop:as-inv-prob-bound}})\\
        \lesssim&T^{-31/26}R_T^{5/11}\sum_{t=1}^T\e\left[R_t^{-1}\iprod{\xv_t,\bv_t(1)-\bv_t^*(1)}^4\right]\quad(\text{since $R_t\leq R_T$})\\
        \lesssim&T^{-31/26}R_T^{5/11}\cdot T^{3/8}R_T^{1/4}\quad(\text{Corollary \ref{corollary:expected-tracking}})\\
        =&(TR_T^{-2})^{-31/88}\cdot T^{-133/286}\\
        \rightarrow&0\quad(\text{Assumption \mainref{assumption:maximum-radius-strengthened} and Lemma \ref{lemma:R}})\enspace,
    \end{align*}
    which completes the proof.
\end{proof}

Based on Lemma \ref{lemma:Lyapunov-1} and Lemma \ref{lemma:Lyapunov-2}, we now verify the conditional Lyapunov condition.

\begin{lemma}\label{lemma:Lyapunov}
    Under Assumptions \mainref{assumption:moments}-\mainref{assumption:bounded-correlation} and Condition \mainref{condition:sigmoid}, it holds that $\sum_{t=1}^T\E{X_{t}^4|\filt_{t-1}}\xrightarrow{p} 0$.
\end{lemma}

\begin{proof}
    For $k\in\{0,1\}$, denote $\Delta_{t,k}(\bv)=y_t(k)-\iprod{\xv_t,\bv}$.
    By Corollary \mainref{corollary:non-superefficiency} and the inequality $(a+b)^4\leq 8a^4+8b^4$, we have
    \begin{align*}
        &\e\left[\sum_{t=1}^T\E{X_{t}^4|\filt_{t-1}}\right]\\
        \lesssim&\frac{1}{T^4(\Var{\eate})^2}\e\left[\sum_{t=1}^T\e\left[\Delta_{t,1}(\bv_t(1))^4\cdot \frac{1}{p_t^3}+\Delta_{t,0}(\bv_t(0))^4\cdot \frac{1}{(1-p_t)^3}\Big|\filt_{t-1}\right]\right]\\
        \lesssim&\frac{1}{T^2}\cdot\sum_{t=1}^T\e\left[\Delta_{t,1}(\bv_t(1))^4\cdot \frac{1}{p_t^3}+\Delta_{t,0}(\bv_t(0))^4\cdot \frac{1}{(1-p_t)^3}\right]\quad(\text{Corollary \mainref{corollary:non-superefficiency}})\\
        \lesssim&\frac{1}{T^2}\sum_{t=1}^T\e\left[\Delta_{t,1}(\bv_t^*(1))^4\cdot \frac{1}{p_t^3} \right]+\frac{1}{T^2}\sum_{t=1}^T\e\left[\iprod{ \xv_t, \bv_t(1)-\bv_t^*(1) } ^4\cdot \frac{1}{p_t^3} \right]\\
        &+\frac{1}{T^2}\sum_{t=1}^T\e\left[\Delta_{t,0}(\bv_t^*(0))^4\cdot \frac{1}{(1-p_t)^3} \right]+\frac{1}{T^2}\sum_{t=1}^T\e\left[\iprod{ \xv_t, \bv_t(0)-\bv_t^*(0) } ^4\cdot \frac{1}{(1-p_t)^3} \right]\enspace,
    \end{align*}
    which converges to 0 by Lemma \ref{lemma:Lyapunov-1}, Lemma \ref{lemma:Lyapunov-2} and symmetric results for the control group ($k=0$). 
    Markov's inequality therefore implies $\sum_{t=1}^T\e[X_{t}^4|\filt_{t-1}] \xrightarrow{p} 0$, which completes the proof.
\end{proof}

We have finally verified the conditional variance condition
and the conditional Lyapunov condition in Lemma \ref{lemma:stable-variance} and Lemma \ref{lemma:Lyapunov}.
The martingale central limit theorem (Lemma \mainref{lemma:martingale-clt}) therefore applies, leading to the following theorem.

\begin{reftheorem}{\mainref{thm:clt}}
    \ourclt
\end{reftheorem}

\begin{proof}
    The result can be readily shown by Lemma \mainref{lemma:martingale-clt}, Lemma \ref{lemma:stable-variance}, and Lemma \ref{lemma:Lyapunov}.
\end{proof}

\subsection{Non-Superefficiency}\label{section:D4}
In this section, we verify the non-superefficiency condition in Corollary \mainref{corollary:non-superefficiency}.
For $k\in\{0,1\}$, denote $\Delta_{t,k}(\bv)=y_t(k)-\iprod{\xv_t,\bv}$.
We further introduce the following notation for simplicity:
\begin{align*}
    \widetilde{A}_T(1)=&\sum_{t=1}^T\frac{\mathbf{1}[Z_t=1]}{p_t}\cdot \Delta_{t,1}(\bv_t^*(1))^2\enspace,\\
    \widetilde{A}_T(0)=&\sum_{t=1}^T\frac{\mathbf{1}[Z_t=0]}{1-p_t}\cdot \Delta_{t,0}(\bv_t^*(0))^2\enspace,\\
    \widetilde{p}=&\operatorname{argmin}_{p\in(0,1)}\sum_{t=1}^T\Bigg(\frac{\mathbf{1}[Z_t=1]}{p_t}\cdot\Delta_{t,1}(\bv_t^*(1))^2\cdot \frac{1}{p}+\frac{\mathbf{1}[Z_t=0]}{1-p_t}\cdot\Delta_{t,0}(\bv_t^*(0))^2\cdot \frac{1}{1-p}\Bigg)+\eta_T^{-1}\Psi(p)\\
    =&\operatorname{argmin}_{p\in(0,1)}\frac{\widetilde{A}_T(1)}{p}+\frac{\widetilde{A}_T(0)}{1-p}+\eta_T^{-1}\Psi(p)\enspace,\\
    \breve{p}=&\operatorname{argmin}_{p\in(0,1)}\sum_{t=1}^T\Bigg(\Delta_{t,1}(\bv_t^*(1))^2\cdot \frac{1}{p}+\Delta_{t,0}(\bv_t^*(0))^2\cdot \frac{1}{1-p}\Bigg)+\eta_T^{-1}\Psi(p)\\
    =&\operatorname{argmin}_{p\in(0,1)}\frac{{A}_T^*(1)}{p}+\frac{{A}_T^*(0)}{1-p}+\eta_T^{-1}\Psi(p)\enspace.
\end{align*}

In Proposition \mainref{prop:pred-regret}*, we showed that the expected prediction regret admits a sublinear upper bound.
However, this result does not exclude the possibility that the expected prediction regret could be substantially negative (for example, $\E{\newpredregret_T} = -T$).
The following lemma further establishes that the expected prediction regret is sublinear in magnitude.

\begin{lemma}\label{lemma:pred-regret-exact-order}
    Under Assumptions \mainref{assumption:moments}-\mainref{assumption:maximum-radius-strengthened}, and Condition \mainref{condition:sigmoid}, it holds that $|\E{\newpredregret_T}|=o(T)$.
\end{lemma}

\begin{proof}
    For $k\in\{0,1\}$, denote $\Delta_{t,k}(\bv)=y_t(k)-\iprod{\xv_t,\bv}$.
    By Lemma \mainref{lemma:standard-oco-for-pred-regret} and Assumption \mainref{assumption:maximum-radius-strengthened}, it has already been shown that the positive part
    $\max\{\E{\mathcal{R}_T^{\text{pred}}},0\}$ is $o(T)$.
    Therefore, it suffices to prove that the negative part $\max\{-\E{\mathcal{R}_T^{\text{pred}}},0\}$ is also $o(T)$.
    For notational simplicity, we use the concatenation $\bv_t=(\bv_t(1),\bv_t(0))$ throughout the proof.
    We first decompose the difference between $\ell_t(\bv_t)$ and the deterministic quantity $\ell_t(\bv_t^*)$.
    By using the equality $\Delta_{t,k}(\bv_t(k))=\Delta_{t,k}(\bv_t^*(k))+\iprod{ \xv_t, \bv_t^*(k)-\bv_t(k) }$ and expanding the squares, we isolate linear and quadratic errors, leading to: 
    \begin{align}\label{lemma:pred-regret-exact-order_eq1}
        &\ell_t(\bv_t)\notag\\
        =&\frac{\mathcal{E}(0)}{\mathcal{E}(1)}\cdot\Delta_{t,1}(\bv_t(1))^2+\frac{\mathcal{E}(1)}{\mathcal{E}(0)}\cdot\Delta_{t,0}(\bv_t(0))^2+2\Delta_{t,1}(\bv_t(1))\cdot\Delta_{t,0}(\bv_t(0))\notag\\
        =&\frac{\mathcal{E}(0)}{\mathcal{E}(1)}\cdot[\Delta_{t,1}(\bv_t^*(1))+\iprod{\xv_t,\bv_t^*(1)-\bv_t(1)}]^2+\frac{\mathcal{E}(1)}{\mathcal{E}(0)}\cdot[\Delta_{t,0}(\bv_t^*(0))+\iprod{\xv_t,\bv_t^*(0)-\bv_t(0)}]^2\notag\\
        &+2[\Delta_{t,1}(\bv_t^*(1))+\iprod{\xv_t,\bv_t^*(1)-\bv_t(1)}]\cdot [\Delta_{t,0}(\bv_t^*(0))+\iprod{\xv_t,\bv_t^*(0)-\bv_t(0)}]\notag\\
        =&\frac{\mathcal{E}(0)}{\mathcal{E}(1)}\cdot\Delta_{t,1}(\bv_t^*(1))^2+\frac{2\mathcal{E}(0)}{\mathcal{E}(1)}\cdot\Delta_{t,1}(\bv_t^*(1))\cdot\iprod{\xv_t,\bv_t^*(1)-\bv_t(1)}+\frac{\mathcal{E}(0)}{\mathcal{E}(1)}\cdot\iprod{\xv_t,\bv_t^*(1)-\bv_t(1)}^2\notag\\
        &+\frac{\mathcal{E}(1)}{\mathcal{E}(0)}\cdot\Delta_{t,0}(\bv_t^*(0))^2+\frac{2\mathcal{E}(1)}{\mathcal{E}(0)}\cdot\Delta_{t,0}(\bv_t^*(0))\cdot\iprod{\xv_t,\bv_t^*(0)-\bv_t(0)}+\frac{\mathcal{E}(1)}{\mathcal{E}(0)}\cdot\iprod{\xv_t,\bv_t^*(0)-\bv_t(0)}^2\notag\\
        &+2\Delta_{t,1}(\bv_t^*(1))\cdot\Delta_{t,0}(\bv_t^*(0))+2\Delta_{t,1}(\bv_t^*(1))\cdot\iprod{\xv_t,\bv_t^*(0)-\bv_t(0)}\notag\\
        &+2\Delta_{t,0}(\bv_t^*(0))\cdot\iprod{\xv_t,\bv_t^*(1)-\bv_t(1)}+2\iprod{\xv_t,\bv_t^*(1)-\bv_t(1)}\cdot\iprod{\xv_t,\bv_t^*(0)-\bv_t(0)}\notag\\
        =&\Bigg[\frac{2\mathcal{E}(0)}{\mathcal{E}(1)}\cdot \Delta_{t,1}(\bv_t^*(1))\cdot\iprod{\xv_t,\bv_t^*(1)-\bv_t(1)}+\frac{2\mathcal{E}(1)}{\mathcal{E}(0)}\cdot\Delta_{t,0}(\bv_t^*(0))\cdot\iprod{\xv_t,\bv_t^*(0)-\bv_t(0)}\notag\\
        &\underbrace{\quad\quad\quad\quad\quad+2\Delta_{t,1}(\bv_t^*(1))\cdot\iprod{\xv_t,\bv_t^*(0)-\bv_t(0)}+2\Delta_{t,0}(\bv_t^*(0))\cdot\iprod{\xv_t,\bv_t^*(1)-\bv_t(1)}\Bigg]}_{:=B_{t,1}}\notag\\
        &+\underbrace{\Bigg[\frac{\mathcal{E}(0)}{\mathcal{E}(1)}\cdot\iprod{\xv_t,\bv_t^*(1)-\bv_t(1)}^2+\frac{\mathcal{E}(1)}{\mathcal{E}(0)}\cdot\iprod{\xv_t,\bv_t^*(0)-\bv_t(0)}^2+2\iprod{\xv_t,\bv_t^*(1)-\bv_t(1)}\cdot \iprod{\xv_t,\bv_t^*(0)-\bv_t(0)}\Bigg]}_{:=B_{t,2}}\notag\\
        &+\ell_t(\bv_t^*(1),\bv_t^*(0))\enspace,
    \end{align}
    where $B_{t,1}$ collects the cross terms, $B_{t,2}$ collects the square terms, and $\ell_t(\bv_t^*(1),\bv_t^*(0))$ is the full-information residual.
    Lemma \mainref{lemma:predictor-expectation} guarantees that $B_{t,1}$ has mean zero.
    By \eqref{lemma:stable-variance-2_eq1} and the AM-GM inequality, we obtain the bound
    \begin{align}\label{lemma:pred-regret-exact-order_eq2}
        \sum_{t=1}^{T}\e[|B_{t,2}|]\leq&\frac{2\mathcal{E}(0)}{\mathcal{E}(1)}\cdot\e\left[\sum_{t=1}^{T}\iprod{\xv_t,\bv_t^*(1)-\bv_t(1)}^2\right]+\frac{2\mathcal{E}(1)}{\mathcal{E}(0)}\cdot\e\left[\sum_{t=1}^{T}\iprod{\xv_t,\bv_t^*(0)-\bv_t(0)}^2\right]\notag\\
        \lesssim&T^{3/8}R_T^{5/4}\quad(\text{by Assumption \mainref{assumption:moments} and \eqref{lemma:stable-variance-2_eq1}})\notag\\
        =&o(T)\quad(\text{Assumption \mainref{assumption:maximum-radius-strengthened} and Lemma \ref{lemma:R}})\enspace.
    \end{align}
    Denote
    \begin{align*}
        \widebar{L}_t(\bv)=&\sum_{s=1}^{t-1}\ell_s(\bv)+\eta_{t-1}^{-1}m(\bv)\enspace,\\
        L^*_t(\bv)=&\sum_{s=1}^{t-1}\ell_s(\bv)+\eta_{t}^{-1}m(\bv)
    \end{align*}
    and let $\widetilde{\bv}^*_t$ denote the minimizer of $\widebar{L}_t$.
    By \eqref{lemma:pred-regret-exact-order_eq1}, \eqref{lemma:pred-regret-exact-order_eq2} and Lemma \ref{lemma:standard-FTRL-decomposition}, we decompose the regret and lower bound it as follows:
    \begin{align}\label{lemma:pred-regret-exact-order_eq3}
        \E{\mathcal{R}_T^{\text{pred}}}=&\e\left[\sum_{t=1}^T\ell_t(\bv_t)-\sum_{t=1}^T\ell_t(\bv^*)\right]\notag\\
        =&\sum_{t=1}^T\ell_t(\bv_t^*)-\sum_{t=1}^T\ell_t(\bv^*)+\sum_{t=1}^T\e[B_{t,1}]+\sum_{t=1}^T\e[B_{t,2}]\quad(\text{by \eqref{lemma:pred-regret-exact-order_eq1}})\notag\\
        =&o(T)+\sum_{t=1}^T\ell_t(\bv_t^*)-\sum_{t=1}^T\ell_t(\bv^*)\quad(\text{by \eqref{lemma:pred-regret-exact-order_eq2}})\notag\\
        \intertext{Applying Lemma \ref{lemma:standard-FTRL-decomposition} leads to}
        =&o(T)+\eta_{T}^{-1}m(\bv^*)+{L}_{T+1}^*(\bv_{T+1}^*)-{L}_{T+1}^*(\bv^*)+\sum_{t=1}^T\left(\widebar{L}_{t+1}(\bv_t^*)-{L}_{t+1}^*(\bv_{t+1}^*)\right)\notag\\
        \intertext{By adding and subtracting $\widebar{L}_{t+1}(\bv_{t+1}^*)$, we can further decompose this as}
        =&o(T)+\eta_{T}^{-1}m(\bv^*)+{L}_{T+1}^*(\bv_{T+1}^*)-{L}_{T+1}^*(\bv^*)+\sum_{t=1}^T\underbrace{\left(\widebar{L}_{t+1}(\bv_{t+1}^*)-{L}_{t+1}^*(\bv_{t+1}^*)\right)}_{=\left(\eta_{t}^{-1}-\eta_{t+1}^{-1}\right)\cdot m(\bv_{t+1}^*)}\notag\\
        &+\sum_{t=1}^T(\widebar{L}_{t+1}(\bv_t^*(1),\bv_t^*(0))-\widebar{L}_{t+1}(\bv_{t+1}^*(1),\bv_{t+1}^*(0)))\notag\\
        \intertext{Since $\widetilde{\bv}_{t+1}^*$ is the minimizer of $\widebar{L}_{t+1}$, we have}
        \geq&o(T)+\underbrace{\eta_{T}^{-1}m(\bv^*)+{L}_{T+1}^*(\bv_{T+1}^*)-{L}_{T+1}^*(\bv^*)}_{:=S_1}+\underbrace{\sum_{t=1}^T\left(\eta_{t}^{-1}-\eta_{t+1}^{-1}\right)m(\bv_{t+1}^*)}_{:=S_2}\notag\\
        &+\underbrace{\sum_{t=1}^T(\widebar{L}_{t+1}(\widetilde{\bv}^*_{t+1})-\widebar{L}_{t+1}(\bv_{t+1}^*))}_{:=S_3}\enspace.
    \end{align}
    \textbf{Step 1:} Derive a lower bound on $S_1$.\\[3mm]
    By the definition of $L_{T+1}(\bv_{T+1}^*)$, we have
    \begin{align}\label{lemma:pred-regret-exact-order_eq4}
        L_{T+1}^*(\bv_{T+1}^*)=&\frac{\olsres{0}}{\olsres{1}}\cdot\left[\sum_{t=1}^T\Delta_{t,1}(\bv_{T+1}^*(1))^2+\eta_T^{-1}\|\bv_{T+1}^*(1)\|^2\right]\notag\\
        &+\frac{\olsres{1}}{\olsres{0}}\cdot\left[\sum_{t=1}^T\Delta_{t,0}(\bv_{T+1}^*(0))^2+\eta_T^{-1}\|\bv_{T+1}^*(0)\|^2\right]\notag\\
        +&2\left[\sum_{t=1}^T\Delta_{t,1}(\bv_{T+1}^*(1))\cdot\Delta_{t,0}(\bv_{T+1}^*(0))+\eta_T^{-1}\iprod{\bv_{T+1}^*(1),\bv_{T+1}^*(0)}\right]\enspace.
    \end{align}
    For simplicity, we denote $\vec{Y}_t(1)=(y_1(1),\ldots,y_t(1))^{\tran}$, $\vec{Y}_t(0)=(y_1(0),\ldots,y_t(0))^{\tran}$ for $t\in[T]$ and $\vec{X}=(\xv_1,\ldots,\xv_T)^{\tran}$. Using the ridge regression identity derived in Lemma \ref{lemma:ridge}, the loss under the ridge estimator can be calculated as
    \begin{align}\label{lemma:pred-regret-exact-order_eq5}
        &\sum_{t=1}^T\Delta_{t,1}(\bv_{T+1}^*(1))^2+\eta_T^{-1}\|\bv_{T+1}^*(1)\|^2\notag\\
        =&\vec{Y}_T^{\tran}(1)\left(\vec{I}_T-\vec{X}(\vec{X}^{\tran}\vec{X}+\eta_T^{-1}\vec{I}_d)^{-1}\vec{X}^{\tran}\right)\vec{Y}_T(1)\enspace.
    \end{align}
    Similarly, we can obtain
    \begin{align}\label{lemma:pred-regret-exact-order_eq6}
        &\sum_{t=1}^T\Delta_{t,0}(\bv_{T+1}^*(0))^2+\eta_T^{-1}\|\bv_{T+1}^*(0)\|^2\notag\\
        =&\vec{Y}_T^{\tran}(0)\left(\vec{I}_T-\vec{X}(\vec{X}^{\tran}\vec{X}+\eta_T^{-1}\vec{I}_d)^{-1}\vec{X}^{\tran}\right)\vec{Y}_T(0)\enspace,\notag\\
        &\sum_{t=1}^T\Delta_{t,1}(\bv_{T+1}^*(1))\cdot\Delta_{t,0}(\bv_{T+1}^*(0))+\eta_T^{-1}\iprod{\bv_{T+1}^*(1),\bv_{T+1}^*(0)}\notag\\
        =&\vec{Y}_T^{\tran}(1)\left(\vec{I}_T-\vec{X}(\vec{X}^{\tran}\vec{X}+\eta_T^{-1}\vec{I}_d)^{-1}\vec{X}^{\tran}\right)\vec{Y}_T(0)\enspace.
    \end{align}
    Substituting \eqref{lemma:pred-regret-exact-order_eq5} and \eqref{lemma:pred-regret-exact-order_eq6} into \eqref{lemma:pred-regret-exact-order_eq4}, we have
    \begin{align}\label{lemma:pred-regret-exact-order_eq7}
        L_{T+1}^*(\bv_{T+1}^*)=&\frac{\olsres{0}}{\olsres{1}}\cdot\vec{Y}_T^{\tran}(1)\left(\vec{I}_T-\vec{X}(\vec{X}^{\tran}\vec{X}+\eta_T^{-1}\vec{I}_d)^{-1}\vec{X}^{\tran}\right)\vec{Y}_T(1)\notag\\
        &+\frac{\olsres{1}}{\olsres{0}}\cdot\vec{Y}_T^{\tran}(0)\left(\vec{I}_T-\vec{X}(\vec{X}^{\tran}\vec{X}+\eta_T^{-1}\vec{I}_d)^{-1}\vec{X}^{\tran}\right)\vec{Y}_T(0)\notag\\
        &+2\cdot\vec{Y}_T^{\tran}(1)\left(\vec{I}_T-\vec{X}(\vec{X}^{\tran}\vec{X}+\eta_T^{-1}\vec{I}_d)^{-1}\vec{X}^{\tran}\right)\vec{Y}_T(0)\enspace.
    \end{align}
    On the other hand, since $L_{T+1}^*(\bv^*)-\eta_{T+1}^{-1}m(\bv^*)$ corresponds to the squared loss under the OLS estimator, we have
    \begin{align}\label{lemma:pred-regret-exact-order_eq8}
        L_{T+1}^*(\bv^*)-\eta_{T+1}^{-1}m(\bv^*)=&\frac{\olsres{0}}{\olsres{1}}\cdot\vec{Y}_T^{\tran}(1)\left(\vec{I}_T-\vec{X}(\vec{X}^{\tran}\vec{X})^{-1}\vec{X}^{\tran}\right)\vec{Y}_T(1)\notag\\
        &+\frac{\olsres{1}}{\olsres{0}}\cdot\vec{Y}_T^{\tran}(0)\left(\vec{I}_T-\vec{X}(\vec{X}^{\tran}\vec{X})^{-1}\vec{X}^{\tran}\right)\vec{Y}_T(0)\notag\\
        &+2\cdot\vec{Y}_T^{\tran}(1)\left(\vec{I}_T-\vec{X}(\vec{X}^{\tran}\vec{X})^{-1}\vec{X}^{\tran}\right)\vec{Y}_T(0)\enspace.
    \end{align}
    Denote
    \begin{align*}
        \vec{U}=\sqrt{\frac{\mathcal{E}(0)}{\mathcal{E}(1)}}\mat{Y}_T(1)+\sqrt{\frac{\mathcal{E}(1)}{\mathcal{E}(0)}}\mat{Y}_T(0)\enspace.
    \end{align*}
    Combining the results in \eqref{lemma:pred-regret-exact-order_eq7} and \eqref{lemma:pred-regret-exact-order_eq8}, we obtain
    \begin{align}\label{lemma:pred-regret-exact-order_eq9}
        S_1=&\eta_{T}^{-1}m(\bv^*)+{L}_{T+1}^*(\bv_{T+1}^*)-{L}_{T+1}^*(\bv^*)\notag\\
        =&\vec{U}^{\tran}\vec{X}(\vec{X}^{\tran}\vec{X})^{-1}\vec{X}^{\tran}\vec{U}-\vec{U}^{\tran}\vec{X}(\vec{X}^{\tran}\vec{X}+\eta_T^{-1}\vec{I}_d)^{-1}\vec{X}^{\tran}\vec{U}\notag\\
        \geq& 0\enspace,
    \end{align}
    where the last step is due to $(\vec{X}^{\tran}\vec{X})^{-1}\succeq (\vec{X}^{\tran}\vec{X}+\eta_T^{-1}\vec{I}_d)^{-1}$.\\[3mm]
    \textbf{Step 2:} Derive a lower bound on $S_2$.\\[3mm]
    To control $S_2$, we first derive an upper bound on the norm of the full-information predictor.
    For any $t\in[T]$, denote $\vec{X}_t=(\xv_1,\ldots,\xv_t)^{\tran}$.
    By Lemma \ref{lemma:regularity}, we have
    \begin{align*}
        \|\bv^*_{t+1}(1)\|^2=&\mat{Y}_{t}^{\tran}(1)\xM_{t}\left(\xM_{t}^{\tran}\xM_t+\eta_{t+1}^{-1}\mat{I}_d\right)^{-2}\xM_{t}^{\tran}\mat{Y}_t(1)\\
        \leq&\left(\sum_{s=1}^ty_s(1)^2\right)\left\|\xM_{t}\left(\xM_{t}^{\tran}\xM_t+\eta_{t+1}^{-1}\mat{I}_d\right)^{-2}\xM_{t}^{\tran}\right\|\\
        \leq&\underbrace{\left(\sum_{s=1}^ty_s(1)^2\right)}_{\text{Cauchy-Schwarz}}\underbrace{\left\|\left(\xM_{t}^{\tran}\xM_t+\eta_{t+1}^{-1}\mat{I}_d\right)^{-1}\right\|}_{\text{Lemma \ref{lemma:regularity}}}\\
        \leq&\left(\sum_{s=1}^ty_s(1)^4\right)^{1/2}t^{1/2}\gamma(t\vee \eta_{t+1}^{-1})^{-1}\\
        \lesssim&T^{1/2}(t\vee \eta_{t+1}^{-1})^{-1/2}~~~~\text{(Assumption \mainref{assumption:moments})}\\
        \leq&T^{1/2}(\eta_{t+1}^{-1})^{-1/2}\quad(\text{since $t\vee \eta_{t+1}^{-1}\geq \eta_{t+1}^{-1}$})\enspace.
    \end{align*}
    Similar results can be obtained for $\bv_{t+1}^*(0)$.
    Hence, we upper bound $-S_2$, which yields a lower bound for $S_2$:
    \begin{align}\label{lemma:pred-regret-exact-order_eq10}
        \quad-S_2=&\sum_{t=1}^T\left(\eta_{t+1}^{-1}-\eta_{t}^{-1}\right)\left(\frac{\mathcal{E}(0)}{\mathcal{E}(1)}\|\bv^*_{t+1}(1)\|^2+\frac{\mathcal{E}(1)}{\mathcal{E}(0)}\|\bv^*_{t+1}(0)\|^2+2\iprod{\bv^*_{t+1}(1),\bv^*_{t+1}(0)}\right)\notag\\
        \lesssim&\sum_{t=1}^T\left(\eta_{t+1}^{-1}-\eta_{t}^{-1}\right)\left(\frac{\mathcal{E}(0)}{\mathcal{E}(1)}\|\bv^*_{t+1}(1)\|^2+\frac{\mathcal{E}(1)}{\mathcal{E}(0)}\|\bv^*_{t+1}(0)\|^2\right)\quad(\text{AM-GM inequality})\notag\\
        \lesssim&\sum_{t=1}^T\left(\eta_{t+1}^{-1}-\eta_{t}^{-1}\right)\left(\|\bv^*_{t+1}(1)\|^2+\|\bv^*_{t+1}(0)\|^2\right)\quad(\text{Assumption \mainref{assumption:moments}})\notag\\
        \lesssim&T^{1/2}\sum_{t=1}^T\left(\eta_{t+1}^{-1}-\eta_{t}^{-1}\right)(\eta_{t+1}^{-1})^{-1/2}\notag\\
        \intertext{Since $\eta_1\geq\ldots\geq\eta_{T+1}$, by an integral comparison argument we have}
        \leq&T^{1/2}\int_{\eta_1^{-1}}^{\eta_{T+1}^{-1}}x^{-1/2}\mathrm{d}x\notag\\
        \lesssim&T^{1/2}\eta_{T+1}^{-1/2}\notag\\
        =&T\cdot (T\eta_{T+1})^{-1/2}\notag\\
        =&o(T)\quad(\text{Assumption \mainref{assumption:maximum-radius-strengthened} and Lemma \ref{lemma:R}})\enspace.
    \end{align}
    \textbf{Step 3:} Derive a lower bound on $S_3$.\\[3mm]
    We first fix $t\in[T]$.
    For notational simplicity, we denote $a=\eta_{t}^{-1}$ and $b=\eta_{t+1}^{-1}$.
    By the explicit form of $\bv^*_{t+1}(1)$ and $\widetilde{\bv}^*_{t+1}(1)$, we have
    \begin{align*}
        &\left\|\mat{Y}_t(1)-\xM_t\bv^*_{t+1}(1)\right\|^2+\eta_t^{-1}\|\bv^*_{t+1}(1)\|^2-\left(\left\|\mat{Y}_t(1)-\xM_t\widetilde{\bv}^*_{t+1}(1)\right\|^2+\eta_t^{-1}\|\widetilde{\bv}^*_{t+1}(1)\|^2\right)\\
        =&\left(\mat{Y}_t^{\tran}(1)\mat{Y}_t(1)-2\mat{Y}_t^{\tran}(1)\xM_t\bv^*_{t+1}(1)+(\bv^*_{t+1}(1))^{\tran}\left(\xM_t^{\tran}\xM_t+a\mat{I}_d\right)\bv^*_{t+1}(1)\right)\\
        &-\left(\mat{Y}_t^{\tran}(1)\mat{Y}_t(1)-2\mat{Y}_t^{\tran}(1)\xM_t\widetilde{\bv}^*_{t+1}(1)+(\widetilde{\bv}^*_{t+1}(1))^{\tran}\left(\xM_t^{\tran}\xM_t+a\mat{I}_d\right)\widetilde{\bv}^*_{t+1}(1)\right)\\
        =&-2\underbrace{\mat{Y}_t^{\tran}(1)\xM_t\left(\xM_t^{\tran}\xM_t+b\mat{I}_d\right)^{-1}\xM_t^{\tran}\mat{Y}_t(1)}_{\mat{Y}_t^{\tran}(1)\xM_t\bv^*_{t+1}(1)}+2\underbrace{\mat{Y}_t^{\tran}(1)\xM_t\left(\xM_t^{\tran}\xM_t+a\mat{I}_d\right)^{-1}\xM_t^{\tran}\mat{Y}_t(1)}_{\mat{Y}_t^{\tran}(1)\xM_t\widetilde{\bv}^*_{t+1}(1)}\\
        &+\underbrace{\mat{Y}_t^{\tran}(1)\xM_t\left(\xM_t^{\tran}\xM_t+b\mat{I}_d\right)^{-1}\left(\xM_t^{\tran}\xM_t+a\mat{I}_d\right)\left(\xM_t^{\tran}\xM_t+b\mat{I}_d\right)^{-1}\xM_t^{\tran}\mat{Y}_t(1)}_{(\bv^*_{t+1}(1))^{\tran}\left(\xM_t^{\tran}\xM_t+a\mat{I}_d\right)\bv^*_{t+1}(1)}\\
        &-\underbrace{\mat{Y}_t^{\tran}(1)\xM_t\left(\xM_t^{\tran}\xM_t+a\mat{I}_d\right)^{-1}\left(\xM_t^{\tran}\xM_t+a\mat{I}_d\right)\left(\xM_t^{\tran}\xM_t+a\mat{I}_d\right)^{-1}\xM_t^{\tran}\mat{Y}_t(1)}_{(\widetilde{\bv}^*_{t+1}(1))^{\tran}\left(\xM_t^{\tran}\xM_t+a\mat{I}_d\right)\widetilde{\bv}^*_{t+1}(1)}\\
        =&\mat{Y}_t^{\tran}(1)\Big(-2\xM_t\left(\xM_t^{\tran}\xM_t+b\mat{I}_d\right)^{-1}\xM_t^{\tran}+2\xM_t\left(\xM_t^{\tran}\xM_t+a\mat{I}_d\right)^{-1}\xM_t^{\tran}\\
        &+\xM_t\left(\xM_t^{\tran}\xM_t+b\mat{I}_d\right)^{-1}\left(\xM_t^{\tran}\xM_t+a\mat{I}_d\right)\left(\xM_t^{\tran}\xM_t+b\mat{I}_d\right)^{-1}\xM_t^{\tran}-\xM_t\left(\xM_t^{\tran}\xM_t+a\mat{I}_d\right)^{-1}\xM_t^{\tran}\Big)\mat{Y}_t(1)\enspace.
    \end{align*}
    By the same method, we also have
    \begin{align*}
        &\left\|\mat{Y}_t(0)-\xM_t\bv^*_{t+1}(0)\right\|^2+\eta_t^{-1}\|\bv^*_{t+1}(0)\|^2-\left(\left\|\mat{Y}_t(0)-\xM_t\widetilde{\bv}^*_{t+1}(0)\right\|^2+\eta_t^{-1}\|\widetilde{\bv}^*_{t+1}(0)\|^2\right)\\
        =&\mat{Y}_t^{\tran}(0)\Big(-2\xM_t\left(\xM_t^{\tran}\xM_t+b\mat{I}_d\right)^{-1}\xM_t^{\tran}+2\xM_t\left(\xM_t^{\tran}\xM_t+a\mat{I}_d\right)^{-1}\xM_t^{\tran}\\
        &+\xM_t\left(\xM_t^{\tran}\xM_t+b\mat{I}_d\right)^{-1}\left(\xM_t^{\tran}\xM_t+a\mat{I}_d\right)\left(\xM_t^{\tran}\xM_t+b\mat{I}_d\right)^{-1}\xM_t^{\tran}-\xM_t\left(\xM_t^{\tran}\xM_t+a\mat{I}_d\right)^{-1}\xM_t^{\tran}\Big)\mat{Y}_t(0)\enspace.
    \end{align*}
    and
    \begin{align*}
        &\iprod{\mat{Y}_t(1)-\xM_t\bv^*_{t+1}(1),\mat{Y}_t(0)-\xM_t\bv^*_{t+1}(0)}+\eta_t^{-1}\iprod{\bv^*_{t+1}(1),\bv^*_{t+1}(0)}\\
        &-\left(\iprod{\mat{Y}_t(1)-\xM_t\widetilde{\bv}^*_{t+1}(1),\mat{Y}_t(0)-\xM_t\widetilde{\bv}^*_{t+1}(0)}+\eta_t^{-1}\iprod{\widetilde{\bv}^*_{t+1}(1),\widetilde{\bv}^*_{t+1}(0)}\right)\\
        =&\mat{Y}_t^{\tran}(1)\Big(-2\xM_t\left(\xM_t^{\tran}\xM_t+b\mat{I}_d\right)^{-1}\xM_t^{\tran}+2\xM_t\left(\xM_t^{\tran}\xM_t+a\mat{I}_d\right)^{-1}\xM_t^{\tran}\\
        &+\xM_t\left(\xM_t^{\tran}\xM_t+b\mat{I}_d\right)^{-1}\left(\xM_t^{\tran}\xM_t+a\mat{I}_d\right)\left(\xM_t^{\tran}\xM_t+b\mat{I}_d\right)^{-1}\xM_t^{\tran}-\xM_t\left(\xM_t^{\tran}\xM_t+a\mat{I}_d\right)^{-1}\xM_t^{\tran}\Big)\mat{Y}_t(0)\enspace.
    \end{align*}
    Hence
    \begin{align*}
        &\widebar{L}_{t+1}(\bv_{t+1}^*(1),\bv_{t+1}^*(0))-\widebar{L}_{t+1}(\widetilde{\bv}^*_{t+1}(1),\widetilde{\bv}^*_{t+1}(0))\\
        =&\left[\left(\frac{\mathcal{E}(0)}{\mathcal{E}(1)}\right)^{1/2}\mat{Y}_t(1)+\left(\frac{\mathcal{E}(1)}{\mathcal{E}(0)}\right)^{1/2}\mat{Y}_t(0)\right]^{\tran}\Big(-2\xM_t\left(\xM_t^{\tran}\xM_t+b\mat{I}_d\right)^{-1}\xM_t^{\tran}+2\xM_t\left(\xM_t^{\tran}\xM_t+a\mat{I}_d\right)^{-1}\xM_t^{\tran}\\
        &+\xM_t\left(\xM_t^{\tran}\xM_t+b\mat{I}_d\right)^{-1}\left(\xM_t^{\tran}\xM_t+a\mat{I}_d\right)\left(\xM_t^{\tran}\xM_t+b\mat{I}_d\right)^{-1}\xM_t^{\tran}-\xM_t\left(\xM_t^{\tran}\xM_t+a\mat{I}_d\right)^{-1}\xM_t^{\tran}\Big)\\
        &\times\left[\left(\frac{\mathcal{E}(0)}{\mathcal{E}(1)}\right)^{1/2}\mat{Y}_t(1)+\left(\frac{\mathcal{E}(1)}{\mathcal{E}(0)}\right)^{1/2}\mat{Y}_t(0)\right]\\
        \intertext{By Assumption \mainref{assumption:moments} and the inequality $\|\vec{\alpha}+\vec{\beta}\|^2\leq 2\|\vec{\alpha}\|^2+2\|\vec{\beta}\|^2$, we have $\left\|\left(\frac{\mathcal{E}(0)}{\mathcal{E}(1)}\right)^{1/2}\mat{Y}_t(1)+\left(\frac{\mathcal{E}(1)}{\mathcal{E}(0)}\right)^{1/2}\mat{Y}_t(0)\right\|^2\lesssim \left(\frac{\mathcal{E}(0)}{\mathcal{E}(1)}\right)\|\mat{Y}_t(1)\|^2+\left(\frac{\mathcal{E}(1)}{\mathcal{E}(0)}\right)\|\mat{Y}_t(0)\|^2\lesssim \|\mat{Y}_t(1)\|^2+\|\mat{Y}_t(0)\|^2$. Hence we have}
        \lesssim&\left(\sum_{s=1}^ty_s(1)^2+\sum_{s=1}^ty_s(0)^2\right)\Big\|\xM_t\Big(-2\left(\xM_t^{\tran}\xM_t+b\mat{I}_d\right)^{-1}+2\left(\xM_t^{\tran}\xM_t+a\mat{I}_d\right)^{-1}-\left(\xM_t^{\tran}\xM_t+a\mat{I}_d\right)^{-1}\\
        &+\left(\xM_t^{\tran}\xM_t+b\mat{I}_d\right)^{-1}\left(\xM_t^{\tran}\xM_t+a\mat{I}_d\right)\left(\xM_t^{\tran}\xM_t+b\mat{I}_d\right)^{-1}\Big)\xM_t^{\tran}\Big\|\\
        =&\left(\sum_{s=1}^ty_s(1)^2+\sum_{s=1}^ty_s(0)^2\right)\Big\|\xM_t^{\tran}\xM_t\Big(-2\left(\xM_t^{\tran}\xM_t+b\mat{I}_d\right)^{-1}+2\left(\xM_t^{\tran}\xM_t+a\mat{I}_d\right)^{-1}-\left(\xM_t^{\tran}\xM_t+a\mat{I}_d\right)^{-1}\\
        &\qquad\qquad\qquad\qquad\qquad\quad\underbrace{+\left(\xM_t^{\tran}\xM_t+b\mat{I}_d\right)^{-1}\left(\xM_t^{\tran}\xM_t+a\mat{I}_d\right)\left(\xM_t^{\tran}\xM_t+b\mat{I}_d\right)^{-1}\Big)\Big\|\quad\quad\qquad\qquad}_{:=\|V\|}\enspace.
    \end{align*}
    To bound the spectral norm of $\mat V$, we compute its eigenvalues explicitly.
    Since $\xM_t^{\tran}\xM_t$, $\xM_t^{\tran}\xM_t+a\mat{I}_d$ and $\xM_t^{\tran}\xM_t+b\mat{I}_d$ share the same eigenvectors and hence can be simultaneously diagonalized, if $\lambda$ is an eigenvalue of $\xM_t^{\tran}\xM_t$, then the corresponding eigenvalue of $\mat{V}$ is given by:
    \begin{align*}
        &\lambda\left(-\frac{2}{\lambda+b}+\frac{2}{\lambda+a}-\frac{1}{\lambda+a}+\frac{\lambda+a}{(\lambda+b)^2}\right)\\
        =&\lambda\cdot \frac{-2(\lambda+b)(\lambda+a)+(\lambda+b)^2+(\lambda+a)^2}{(\lambda+b)^2(\lambda+a)}\\
        =&\lambda\cdot \frac{(\lambda+b-\lambda-a)^2}{(\lambda+b)^2(\lambda+a)}\\
        \leq&(b-a)^2\frac{1}{(\lambda+b)^2}
        \quad(\text{since $a
        \geq 0$})\enspace.
    \end{align*}
    By the same argument as in Lemma \ref{lemma:regularity}, we can show that $(\lambda+b)^{-1}=(\lambda+\eta_{t+1}^{-1})^{-1}\leq \gamma(t\vee \eta_{t+1}^{-1})^{-1}$, then we have $\|\mat{V}\|\leq \gamma^2(\eta_{t+1}^{-1}-\eta_{t}^{-1})^2(t\vee \eta_{t+1}^{-1})^{-2}$. 
    Hence, we can upper bound $-S_3$ as follows, which leads to a lower bound for $S_3$:
    \begin{align}\label{lemma:pred-regret-exact-order_eq11}
        -S_3=&\sum_{t=1}^T\left(\widebar{L}_{t+1}(\bv_{t+1}^*)-\widebar{L}_{t+1}(\widetilde{\bv}^*_{t+1})\right)\notag\\
        \lesssim&\sum_{t=1}^T\underbrace{\left(\sum_{s=1}^ty_s(1)^2+\sum_{s=1}^ty_s(0)^2\right)}_{\text{Cauchy-Schwarz}}\|\mat{V}\|\notag\\
        \lesssim&\sum_{t=1}^T\left[\left(\sum_{s=1}^ty_s(1)^4\right)^{1/2}+\left(\sum_{s=1}^ty_s(0)^4\right)^{1/2}\right]t^{1/2}(\eta_{t+1}^{-1}-\eta_{t}^{-1})^2(t\vee \eta_{t+1}^{-1})^{-2}\notag\\
        \lesssim&T^{1/2}\sum_{t=1}^T(\eta_{t+1}^{-1}-\eta_{t}^{-1})^2(t\vee \eta_{t+1}^{-1})^{-3/2}\quad(\text{Assumption \mainref{assumption:moments}})\notag\\
        \intertext{Since $(t\vee \eta_{t+1}^{-1})^{-1}(\eta_{t+1}^{-1}-\eta_{t}^{-1})\leq (t\vee \eta_{t+1}^{-1})^{-1}\eta_{t+1}^{-1}\leq 1$, we have}
        \leq&T^{1/2}\sum_{t=1}^T(\eta_{t+1}^{-1}-\eta_{t}^{-1})(\eta_{t+1}^{-1})^{-1/2}\notag\\
        \intertext{Since $\eta_1\geq\ldots\geq\eta_{T+1}$, by an integral comparison argument we have}
        \leq&T^{1/2}\int_{\eta_1^{-1}}^{\eta_{T+1}^{-1}}x^{-1/2}\mathrm{d}x\notag\\
        \lesssim&T^{1/2}\eta_{T+1}^{-1/2}\notag\\
        =&T\cdot (T\eta_{T+1})^{-1/2}\notag\\
        =&o(T)\quad(\text{Assumption \mainref{assumption:maximum-radius-strengthened} and Lemma \ref{lemma:R}})\enspace.
    \end{align}
    Hence by \eqref{lemma:pred-regret-exact-order_eq3}, \eqref{lemma:pred-regret-exact-order_eq9}, \eqref{lemma:pred-regret-exact-order_eq10} and \eqref{lemma:pred-regret-exact-order_eq11}, we have proved that $\max\{-\e[\mathcal{R}_T^{\text{pred}}],0\}$ is of order $o(T)$, which completes the proof of Lemma~\ref{lemma:pred-regret-exact-order}.
\end{proof}

The following corollary follows from the proofs of Proposition \mainref{prop:pred-regret}* and Lemma~\ref{lemma:pred-regret-exact-order}.
It shows that the sum of squared full-information online residuals is asymptotically equivalent to the squared loss under the OLS estimator.

\begin{corollary}\label{corollary:squared-residual-deterministic-exact-order}
    Under Assumptions \mainref{assumption:moments}-\mainref{assumption:maximum-radius-strengthened}, for $k\in\{0,1\}$, it holds that $A_T^*(k)-T\cdot\olsres{k}^2=o(T)$ and $A_T^*(k)=\Theta(T)$.
\end{corollary}

\begin{proof}
    By definition, we have
    \begin{align*}
        A_T^*(k)-T\cdot\olsres{k}^2=\sum_{t=1}^T(y_t(k)-\iprod{\xv_t,\bv_t^*(k)})^2-\sum_{t=1}^T(y_t(k)-\iprod{\xv_t,\bv^*(k)})^2\enspace.
    \end{align*}
    Following the same argument as in Proposition \mainref{prop:pred-regret}*, we can show that the positive part of $A_T^*(k)-T\cdot\olsres{k}^2$ is $o(T)$.
    Following the same argument as in Lemma \ref{lemma:pred-regret-exact-order}, we can show that the negative part of $A_T^*(k)-T\cdot\olsres{k}^2$ is also $o(T)$.
    Hence $A_T^*(k)-T\cdot\olsres{k}^2=o(T)$.
    Since Assumption \mainref{assumption:moments} implies that $T\cdot\olsres{k}^2=\Theta(T)$, we can deduce that $A_T^*(k)=\Theta(T)$.
\end{proof}

The following lemma provides an upper bound on the variability of the inverse probability at the terminal time $T$.

\begin{lemma}\label{lemma:difference-inverse-probability-last}
    Under Assumptions \mainref{assumption:moments}-\mainref{assumption:maximum-radius-strengthened}, and Condition \mainref{condition:sigmoid}, it holds that
    \begin{align*}
        \widetilde{A}_T(1)\left|\frac{1}{\widetilde{p}}-\frac{1}{\breve{p}}\right|&\lesssim |\widetilde{A}_T(1)-A_T^*(1)|+|\widetilde{A}_T(0)-A_T^*(0)|\enspace,\\
        \widetilde{A}_T(0)\left|\frac{1}{1-\widetilde{p}}-\frac{1}{1-\breve{p}}\right|&\lesssim |\widetilde{A}_T(1)-A_T^*(1)|+|\widetilde{A}_T(0)-A_T^*(0)|\enspace.
    \end{align*}
\end{lemma}

\begin{proof}
    Without loss of generality, we only prove the first inequality.
    Denote $\delta=|\widetilde{A}_T(1)-A_T^*(1)|+|\widetilde{A}_T(0)-A_T^*(0)|$.
    Corollary \ref{corollary:squared-residual-deterministic-exact-order} indicates that there exist constants $m,M>0$ such that $mT\leq A_T^*(1)\leq MT$ and $mT\leq A_T^*(0)\leq MT$.
    Hence, we have $1/\breve{p}=\bigO{1}$ by Lemma \mainref{lemma:effect-of-p-regularization}*.
    We then discuss the following two cases.\\[3mm]
    \textbf{Case 1:} $\delta\leq mT/2$.\\[3mm]
    In this case, it holds that $\widetilde{A}_T(1)=\Theta(T)$ and $\widetilde{A}_T(0)=\Theta(T)$.
    Then we have $1/\widetilde{p}=\bigO{1}$ by Lemma \mainref{lemma:effect-of-p-regularization}*.
    If $\widetilde{A}_T(1)\geq \widetilde{A}_T(0)$, $A_T^*(1)\geq A_T^*(0)$ or $\widetilde{A}_T(1)\leq \widetilde{A}_T(0)$, $A_T^*(1)\leq A_T^*(0)$, then Lemma \ref{lemma:difference-inverse-probability-main} indicates that
    \begin{align*}
        \left|\frac{1}{\widetilde{p}}-\frac{1}{\breve{p}}\right|\lesssim \frac{|\widetilde{A}_T(1)-A_T^*(1)|}{T}+\frac{|\widetilde{A}_T(0)-A_T^*(0)|}{T}\lesssim \frac{|\widetilde{A}_T(1)-A_T^*(1)|}{\widetilde{A}_T(1)}+\frac{|\widetilde{A}_T(0)-A_T^*(0)|}{\widetilde{A}_T(1)}\enspace.
    \end{align*}
    If $\widetilde{A}_T(1)\geq \widetilde{A}_T(0)$, $A_T^*(1)\leq A_T^*(0)$ or $\widetilde{A}_T(1)\leq \widetilde{A}_T(0)$, $A_T^*(1)\geq A_T^*(0)$, since $\widetilde{p}$ and $\breve{p}$ are bounded away from 0 and 1, by Lemma \ref{lemma:difference-inverse-probability-half} we have
    \begin{align*}
        \left|\frac{1}{\widetilde{p}}-2\right|&\lesssim \frac{|\widetilde{A}_T(1)-\widetilde{A}_T(0)|}{T}\enspace,\\
        \left|\frac{1}{\breve{p}}-2\right|&\lesssim \frac{|A^*_T(1)-A^*_T(0)|}{T}\enspace.
    \end{align*}
    Since $\{\widetilde{A}_T(1),\widetilde{A}_T(0)\}$ and $\{A^*_T(1),A^*_T(0)\}$ have different orderings, it holds that
    \begin{align*}
        \left|\frac{1}{\widetilde{p}}-\frac{1}{\breve{p}}\right|\lesssim&\left|\frac{1}{\widetilde{p}}-2\right|+\left|\frac{1}{\breve{p}}-2\right|\\
        \lesssim&\frac{|\widetilde{A}_T(1)-\widetilde{A}_T(0)|+|A^*_T(1)-A^*_T(0)|}{T}\\
        =&\frac{|(\widetilde{A}_T(1)-\widetilde{A}_T(0))-(A^*_T(1)-A^*_T(0))|}{T}\\
        \leq&\frac{|\widetilde{A}_T(1)-A^*_T(1)|+|\widetilde{A}_T(0)-A^*_T(0)|}{T}\\
        \lesssim&\frac{|\widetilde{A}_T(1)-A_T^*(1)|}{\widetilde{A}_T(1)}+\frac{|\widetilde{A}_T(0)-A_T^*(0)|}{\widetilde{A}_T(1)}\enspace.
    \end{align*}
    \textbf{Case 2:} $\delta\geq mT/2$.\\[3mm]
    The optimality of $\widetilde{p}$ indicates that
    \begin{align*}
        \frac{\widetilde{A}_T(1)}{\widetilde{p}}\leq\frac{\widetilde{A}_T(1)}{\widetilde{p}}+\frac{\widetilde{A}_T(0)}{1-\widetilde{p}}+\eta_T^{-1}\Psi(\widetilde{p})\leq \frac{\widetilde{A}_T(1)}{1/2}+\frac{\widetilde{A}_T(0)}{1-1/2}+\eta_T^{-1}\Psi(1/2)=2(\widetilde{A}_T(1)+\widetilde{A}_T(0))\enspace.
    \end{align*}
    Hence, we have that
    \begin{align*}
        \widetilde{A}_T(1)\left|\frac{1}{\widetilde{p}}-\frac{1}{\breve{p}}\right|\leq& \frac{\widetilde{A}_T(1)}{\widetilde{p}}+\frac{\widetilde{A}_T(1)}{\breve{p}}\\
        \lesssim& (\widetilde{A}_T(1)+\widetilde{A}_T(0))+\widetilde{A}_T(1)\quad(\text{since $1/\breve{p}=\bigO{1}$})\\
        \leq&(A_T^*(1)+A_T^*(0))+A_T^*(1)+2\delta\\
        \lesssim&\delta\enspace,
    \end{align*}
    where the last step is because $A_T^*(1)=\bigO{T}, A_T^*(0)=\bigO{T}$ and $\delta\geq mT/2$.
    Therefore, we have proved that
    \begin{align*}
        \widetilde{A}_T(1)\left|\frac{1}{\widetilde{p}}-\frac{1}{\breve{p}}\right| \lesssim |\widetilde{A}_T(1)-A_T^*(1)|+|\widetilde{A}_T(0)-A_T^*(0)|\enspace.
    \end{align*}
    This completes the proof.
\end{proof}

The following lemma characterizes the difference in the inverse probabilities when different step sizes $\eta$ are used.
It is proved by a similar argument to that of Lemma \ref{lemma:difference-inverse-probability-main}.

\begin{lemma}\label{lemma:difference-inverse-probability-step}
    Let $A, B>0$ and $\eta\geq \widetilde{\eta}>0$.
    Suppose $p,\widetilde{p}$ satisfy
    \begin{align*}
        -\frac{A}{p^2}+\frac{B}{(1-p)^2}+\eta^{-1}\Psi^{\prime}(p)&=0\enspace,\\
        -\frac{A}{\widetilde{p}^2}+\frac{B}{(1-\widetilde{p})^2}+\widetilde{\eta}^{-1}\Psi^{\prime}(\widetilde{p})&=0\enspace.
    \end{align*}
    \begin{enumerate}
        \item[(1)] If $A\geq B$, then under Condition \mainref{condition:sigmoid} we have the following upper bound:
        \begin{align*}
            \left|\frac{1}{1-p}-\frac{1}{1-\widetilde{p}}\right|\leq&\frac{6\bb_1}{\bb_3}\left(1+\frac{2}{\bb_3}\right)^5\frac{|\eta^{-1}-\widetilde{\eta}^{-1}|}{A}\left(\frac{p}{1-p}\right)\left(\frac{\widetilde{p}}{1-\widetilde{p}}\right)^4\enspace.
        \end{align*}
        \item[(2)] If $A\leq B$, then under Condition \mainref{condition:sigmoid} we have the following upper bound:
        \begin{align*}
            \left|\frac{1}{p}-\frac{1}{\widetilde{p}}\right|\leq&\frac{6\bb_1}{\bb_3}\left(1+\frac{2}{\bb_3}\right)^5\frac{|\eta^{-1}-\widetilde{\eta}^{-1}|}{B}\left(\frac{1-p}{p}\right)\left(\frac{1-\widetilde{p}}{\widetilde{p}}\right)^4\enspace.
        \end{align*}
    \end{enumerate}
\end{lemma}

\begin{proof}
    \textbf{Part 1:} Since $A\geq B$, we have $p,\widetilde{p}\geq 1/2$ by Lemma \ref{lemma:p-location}.
    For $u=\phi^{-1}(p)\geq 0$ and $\widetilde{u}=\phi^{-1}(\widetilde{p})\geq 0$, they should satisfy
    \begin{align}\label{lemma:difference-inverse-probability-step_eq1}
        A\left(\frac{1}{\phi(u)}\right)^{\prime}+B\left(\frac{1}{1-\phi(u)}\right)^{\prime}+\eta^{-1}(u+3u^2)=&0\enspace,\notag\\
        A\left(\frac{1}{\phi(\widetilde{u})}\right)^{\prime}+B\left(\frac{1}{1-\phi(\widetilde{u})}\right)^{\prime}+\widetilde{\eta}^{-1}(\widetilde{u}+3\widetilde{u}^2)=&0\enspace.
    \end{align}
    By subtracting the second equation in \eqref{lemma:difference-inverse-probability-step_eq1} from the first equation, we obtain
    \begin{align}\label{lemma:difference-inverse-probability-step_eq2}
        &A\left[\left(\frac{1}{\phi(u)}\right)^{\prime}-\left(\frac{1}{\phi(\widetilde{u})}\right)^{\prime}\right]+B\left[\left(\frac{1}{1-\phi(u)}\right)^{\prime}-\left(\frac{1}{1-\phi(\widetilde{u})}\right)^{\prime}\right]\notag\\
        &+\eta^{-1}(u-\widetilde{u})(1+3u+3\widetilde{u})+(\eta^{-1}-\widetilde{\eta}^{-1})(\widetilde{u}+3\widetilde{u}^2)=0\enspace.
    \end{align}
    The convexity assumption in Condition \mainref{condition:sigmoid}-2 implies that $A\left[\left(\frac{1}{\phi(u)}\right)^{\prime}-\left(\frac{1}{\phi(\widetilde{u})}\right)^{\prime}\right]$, $B\left[\left(\frac{1}{1-\phi(u)}\right)^{\prime}-\left(\frac{1}{1-\phi(\widetilde{u})}\right)^{\prime}\right]$ and $\eta^{-1}(u-\widetilde{u})(1+3u+3\widetilde{u})$ should have the same sign as $u-\widetilde{u}$. 
    Hence by \eqref{lemma:difference-inverse-probability-step_eq2} and Lemma \ref{lemma:u}, we have
    \begin{align}\label{lemma:difference-inverse-probability-step_eq3}
        &|(\eta^{-1}-\widetilde{\eta}^{-1})(\widetilde{u}+3\widetilde{u}^2)|\notag\\
        =&\left|A\left[\left(\frac{1}{\phi(u)}\right)^{\prime}-\left(\frac{1}{\phi(\widetilde{u})}\right)^{\prime}\right]+B\left[\left(\frac{1}{1-\phi(u)}\right)^{\prime}-\left(\frac{1}{1-\phi(\widetilde{u})}\right)^{\prime}\right]+\eta^{-1}(u-\widetilde{u})(1+3u+3\widetilde{u})\right|\notag\\
        \geq&A\left|\left(\frac{1}{\phi(u)}\right)^{\prime}-\left(\frac{1}{\phi(\widetilde{u})}\right)^{\prime}\right|\notag\\
        \geq&\frac{\bb_3}{2}\frac{A}{(1+\widetilde{u})(1+u)(1+\widetilde{u}\land u)}\cdot |u-\widetilde{u}|\quad(\text{Lemma \ref{lemma:u}-5})\enspace.
    \end{align}
    Then by \eqref{lemma:difference-inverse-probability-step_eq3}, we have
    \begin{align*}
        |u-\widetilde{u}|\leq& \frac{2}{\bb_3}\cdot\frac{(1+\widetilde{u})(1+u)(1+\widetilde{u}\land u)}{A}\cdot|\eta^{-1}-\widetilde{\eta}^{-1}|(\widetilde{u}+3\widetilde{u}^2)\quad(\text{by \eqref{lemma:difference-inverse-probability-step_eq3}})\\
        \leq&\frac{2}{\bb_3}\frac{|\eta^{-1}-\widetilde{\eta}^{-1}|}{A}(1+u)(1+\widetilde{u})^2(\widetilde{u}+3\widetilde{u}^2)
        \quad(\text{since $1+\widetilde{u}\land u\leq 1+\widetilde{u}$})\\
        \leq&\frac{6}{\bb_3}\frac{|\eta^{-1}-\widetilde{\eta}^{-1}|}{A}(1+u)(1+\widetilde{u})^4\quad(\text{since $\widetilde{u}+3\widetilde{u}^2\leq 3(1+\widetilde{u})^2$ for $\widetilde{u}\geq 0$})\\
        \intertext{Since we have $\frac{1}{1-p}-\frac{1}{1-1/2}\geq \frac{\bb_3}{2}|u-0|$ and $\frac{1}{1-\widetilde{p}}-\frac{1}{1-1/2}\geq \frac{\bb_3}{2}|\widetilde{u}-0|$ by Lemma \ref{lemma:u}-3, we have}
        \leq&\frac{6}{\bb_3}\frac{|\eta^{-1}-\widetilde{\eta}^{-1}|}{A}\left[\frac{2}{\bb_3}\left(\frac{1}{1-p}-\frac{1}{1-1/2}\right)+1\right]\left[\frac{2}{\bb_3}\left(\frac{1}{1-\widetilde{p}}-\frac{1}{1-1/2}\right)+1\right]^4\\
        \leq&\frac{6}{\bb_3}\left(1+\frac{2}{\bb_3}\right)^5\frac{|\eta^{-1}-\widetilde{\eta}^{-1}|}{A}\left(\frac{p}{1-p}\right)\left(\frac{\widetilde{p}}{1-\widetilde{p}}\right)^4\enspace.
    \end{align*}
    Hence we can prove by Lemma \ref{lemma:u}-3 that
    \begin{align*}
        \left|\frac{1}{1-p}-\frac{1}{1-\widetilde{p}}\right|\leq \bb_1|u-\widetilde{u}|\leq \frac{6\bb_1}{\bb_3}\left(1+\frac{2}{\bb_3}\right)^5\frac{|\eta^{-1}-\widetilde{\eta}^{-1}|}{A}\left(\frac{p}{1-p}\right)\left(\frac{\widetilde{p}}{1-\widetilde{p}}\right)^4\enspace.
    \end{align*}
    \textbf{Part 2:}  Let $q=1-p$ and $\widetilde{q}=1-\widetilde{p}$. 
        Following the same arguments as in part 2 of Lemma \ref{lemma:difference-inverse-probability-main}, $q$ and $\widetilde{q}$ satisfy
        \begin{align*}
            -\frac{B}{q^2}+\frac{A}{(1-q)^2}+\eta^{-1}\Psi^{\prime}(q)&=0\enspace,\\
            -\frac{B}{\widetilde{q}^2}+\frac{A}{(1-\widetilde{q})^2}+\widetilde{\eta}^{-1}\Psi^{\prime}(\widetilde{q})&=0\enspace.
        \end{align*}
        By the claim in part 1, we can prove that
        \begin{align*}
            \left|\frac{1}{p}-\frac{1}{\widetilde{p}}\right|=\left|\frac{1}{1-q}-\frac{1}{1-\widetilde{q}}\right|\leq& \frac{6\bb_1}{\bb_3}\left(1+\frac{2}{\bb_3}\right)^5\frac{|\eta^{-1}-\widetilde{\eta}^{-1}|}{B}\left(\frac{q}{1-q}\right)\left(\frac{\widetilde{q}}{1-\widetilde{q}}\right)^4\\
            =&\frac{6\bb_1}{\bb_3}\left(1+\frac{2}{\bb_3}\right)^5\frac{|\eta^{-1}-\widetilde{\eta}^{-1}|}{B}\left(\frac{1-p}{p}\right)\left(\frac{1-\widetilde{p}}{\widetilde{p}}\right)^4\enspace.
        \end{align*}
        This completes the proof. \qedhere
\end{proof}

Using Lemma \ref{lemma:pred-regret-exact-order}, Lemma \ref{lemma:difference-inverse-probability-last} and Lemma \ref{lemma:difference-inverse-probability-step}, we can now derive the exact form for the asymptotic variance of $\hat{\tau}$ in the following theorem.

\begin{reftheorem}{\mainref{thm:exact-asymptotic-variance}}
	\asymptoticvariance
\end{reftheorem}

\begin{proof}
    For $k\in\{0,1\}$, denote $\Delta_{t,k}(\bv)=y_t(k)-\iprod{\xv_t,\bv}$.
    By Lemma \mainref{lemma:regret-decomposition}, we have the following decomposition for the variance:
    \begin{align}\label{thm:exact-asymptotic-variance_eq1}
        T\cdot \operatorname{Var}(\hat{\tau})=&T\cdot \optvar+\frac{1}{T}\E{\regretprob_T}+\frac{1}{T}\e[\mathcal{R}_T^{\text{pred}}]\enspace.
    \end{align}
    The proof proceeds by showing that all regret terms in the variance decomposition vanish asymptotically.
    In Lemma \ref{lemma:pred-regret-exact-order}, we have shown that $\frac{1}{T}\e[\mathcal{R}_T^{\text{pred}}]=o(1)$.
    Then it suffices to prove that $\frac{1}{T}\E{\regretprob_T}=o(1)$.
    
    For simplicity, we introduce the following notation:
    \begin{align*}
        \widehat{f}_t(p)=&\frac{\mathbf{1}[Z_t=1]}{p_t}\cdot\Delta_{t,1}(\bv_t(1))^2\cdot \frac{1}{p}+\frac{\mathbf{1}[Z_t=0]}{1-p_t}\cdot\Delta_{t,0}(\bv_t(0))^2\cdot \frac{1}{1-p}\enspace,\\
        \widehat{F}_{t}(p)=&\sum_{s=1}^{t-1}\widehat{f}_s(p)+\frac{1}{\eta_t}\Psi(p)\enspace,\\
        \widetilde{F}_t(p)=&\sum_{s=1}^{t-1}\widehat{f}_s(p)+\frac{1}{\eta_{t-1}}\Psi(p)\enspace,\\
        \widehat{h}_t(u)=&\frac{\mathbf{1}[Z_t=1]}{p_t}\cdot\Delta_{t,1}(\bv_t(1))^2\cdot\frac{1}{\phi(u)}+\frac{\mathbf{1}[Z_t=0]}{1-p_t}\cdot\Delta_{t,0}(\bv_t(0))^2\cdot \frac{1}{1-\phi(u)}\enspace,\\
        \widehat{H}_t(u)=&\sum_{s=1}^{t-1}\widehat{h}_s(u)+\frac{1}{\eta_t}\psi(u)\enspace,\\
        \widetilde{H}_t(u)=&\sum_{s=1}^{t-1}\widehat{h}_s(u)+\frac{1}{\eta_{t-1}}\psi(u)\enspace.
    \end{align*}
    Let $\widetilde{p}_{t}$, $u_t$, $\widetilde{u}_{t}$ be the minimizers of the functions $\widetilde{F}_{t}$, $\widehat{H}_t$ and $\widetilde{H}_t$, respectively. 
    Using the analogous version of Lemma \mainref{lemma:unbiased-sigmoid-loss} on the $p$-space, we can decompose the expected probability regret as:
    \begin{align}\label{thm:exact-asymptotic-variance_eq2}
        \frac{1}{T}\E{\regretprob_T}=&\frac{1}{T}\e\left[\sum_{t=1}^Tf_t(p_t)-\sum_{t=1}^Tf_t(p^*)\right]\notag\\
        =&\frac{1}{T}\e\left[\sum_{t=1}^T\widehat{f}_t(p_t)-\sum_{t=1}^T\widehat{f}_t(p^*)\right]\quad(\text{Lemma \mainref{lemma:unbiased-sigmoid-loss} and law of iterated expectations})\notag\\
        \intertext{Applying Lemma \ref{lemma:standard-FTRL-decomposition} leads to}
        =&\frac{1}{T}\e\left[\frac{1}{\eta_{T}}\Psi(p^{*})+\sum_{t=1}^T\left(\widetilde{F}_{t+1}(p_t)-\widehat{F}_{t+1}(p_{t+1})\right)+\widehat{F}_{T+1}(p_{T+1})-\widehat{F}_{T+1}(p^*)\right]\notag\\
        \intertext{By adding and subtracting $\widetilde{F}_{t+1}(\widetilde{p}_{t+1})$ and $\widetilde{F}_{t+1}(p_{t+1})$, we can further decompose it as}
        =&\frac{1}{T\eta_{T}}\Psi(p^{*})+\frac{1}{T}\sum_{t=1}^T\e\left[\widetilde{F}_{t+1}(p_t)-\widetilde{F}_{t+1}(\widetilde{p}_{t+1})\right]+\frac{1}{T}\sum_{t=1}^T\e\left[\widetilde{F}_{t+1}(\widetilde{p}_{t+1})-\widetilde{F}_{t+1}(p_{t+1})\right]\notag\\
        &+\frac{1}{T}\sum_{t=1}^T\e\left[\widetilde{F}_{t+1}(p_{t+1})-\widehat{F}_{t+1}({p}_{t+1})\right]+\frac{1}{T}\e\left[\widehat{F}_{T+1}(p_{T+1})-\widehat{F}_{T+1}(p^*)\right]\notag\\
        \intertext{Proposition \mainref{prop:prob-regret-bound} has already shown that $0\leq \frac{1}{T\eta_{T}}\Psi(p^{*})+\frac{1}{T}\sum_{t=1}^T\e\left[\widetilde{F}_{t+1}(p_t)-\widetilde{F}_{t+1}(\widetilde{p}_{t+1})\right]=o(1)$. Hence we have}
        =&o(1)+\frac{1}{T}\sum_{t=1}^T\e\left[\widetilde{F}_{t+1}(\widetilde{p}_{t+1})-\widetilde{F}_{t+1}(p_{t+1})\right]+\frac{1}{T}\sum_{t=1}^T\e\underbrace{\Big[\widetilde{F}_{t+1}(p_{t+1})-\widehat{F}_{t+1}({p}_{t+1})\Big]}_{=-\left(\eta_{t+1}^{-1}-\eta_t^{-1}\right)\Psi(p_{t+1})}\notag\\
        &+\frac{1}{T}\e\left[\widehat{F}_{T+1}(p_{T+1})-\widehat{F}_{T+1}(p^*)\right]\notag\\
        \intertext{By adding and subtracting $\widehat{F}_{T+1}(\widetilde{p})$ and $\widehat{F}_{T+1}(\breve{p})$, we have}
        =&o(1)-\underbrace{\frac{1}{T}\sum_{t=1}^T\e\left[\widetilde{F}_{t+1}(p_{t+1})-\widetilde{F}_{t+1}(\widetilde{p}_{t+1})\right]}_{:=S_1}-\underbrace{\frac{1}{T}\sum_{t=1}^T\left(\eta_{t+1}^{-1}-\eta_t^{-1}\right)\e\left[\Psi(p_{t+1})\right]}_{:=S_2}\notag\\
        &+\underbrace{\frac{1}{T}\e\left[\widehat{F}_{T+1}(p_{T+1})-\widehat{F}_{T+1}(\widetilde{p})\right]}_{:=S_3}+\underbrace{\frac{1}{T}\e\left[\widehat{F}_{T+1}(\widetilde{p})-\widehat{F}_{T+1}(\breve{p})\right]}_{:=S_4}\notag\\
        &+\underbrace{\frac{1}{T}\e\left[\widehat{F}_{T+1}(\breve{p})-\widehat{F}_{T+1}(p^*)\right]}_{:=S_5}\enspace.
    \end{align}
    We then prove that each of $S_1$, $S_2$, $S_3$, $S_4$ and $S_5$ is of order $o(1)$.\\[3mm]
    \textbf{Step 1:} Prove $S_1=o(1)$.\\[3mm]
    By the definition of $\widetilde{F}$ and $\widetilde{H}$, we have
    \begin{align*}
        0\leq&\frac{1}{T}\e\left[\sum_{t=1}^T\left(\widetilde{F}_{t+1}(p_{t+1})-\widetilde{F}_{t+1}(\widetilde{p}_{t+1})\right)\right]\quad(\text{optimality of $\widetilde{p}_{t+1}$})\\
        =&\frac{1}{T}\sum_{t=1}^T\e\left[\left(\widetilde{H}_{t+1}(u_{t+1})-\widetilde{H}_{t+1}(\widetilde{u}_{t+1})\right)\right]\quad(\text{sigmoid transformation})\\
        =&\frac{1}{T}\e\left[\sum_{t=1}^T\left[\langle \nabla \widetilde{H}_{t+1}(\widetilde{u}_{t+1}),u_{t+1}-\widetilde{u}_{t+1}\rangle+\mathcal{B}_{\,\widetilde{\!H}_{t+1}}(u_{t+1}|\widetilde{u}_{t+1})\right]\right]\quad(\text{definition of Bregman divergence})\\
        =&\frac{1}{T}\e\left[\sum_{t=1}^T\mathcal{B}_{\,\widetilde{\!H}_{t+1}}(u_{t+1}|\widetilde{u}_{t+1})\right]\quad(\text{optimality of $\widetilde{u}_{t+1}$})\enspace.
    \end{align*}
    For any $t\in[T]$, if $\widehat{A}_t(1)\geq \widehat{A}_t(0)$, then $u_{t+1},\widetilde{u}_{t+1}\geq 0$ by Lemma \ref{lemma:p-location}. 
    Since $\eta_{t+1}\leq \eta_{t}$, we can also show that $\widetilde{u}_{t+1}\geq u_{t+1}$ and $\widetilde{p}_{t+1}\geq p_{t+1}$ by similar arguments as in Lemma \ref{lemma:p-location} (using the monotonicity property of the gradient).
    We then control the Bregman divergence terms appearing in $\widetilde{H}_{t+1}$.
    For $f_1(u):=1/\phi(u)$ and $f_0(u):=1/(1-\phi(u))$ (which are both convex functions by Condition \mainref{condition:sigmoid}-2), by Lemma \ref{lemma:u}--6, we have
    \begin{align}\label{thm:exact-asymptotic-variance_eq3}
        \mathcal{B}_{f_1}(u_{t+1}|\widetilde{u}_{t+1})=&f_1(u_{t+1})-f_1(\widetilde{u}_{t+1})-f_1^{\prime}(\widetilde{u}_{t+1})(u_{t+1}-\widetilde{u}_{t+1})\notag\\
        =&\frac{1}{\phi(u_{t+1})}-\frac{1}{\phi(\widetilde{u}_{t+1})}-\left(\frac{1}{\phi(\widetilde{u}_{t+1})}\right)^{\prime}(u_{t+1}-\widetilde{u}_{t+1})\notag\\
        \leq&\bb_2\cdot \frac{\widetilde{u}_{t+1}-u_{t+1}}{(1+u_{t+1})(1+\widetilde{u}_{t+1})}\enspace,\notag\\
        \mathcal{B}_{f_0}(u_{t+1}|\widetilde{u}_{t+1})=&f_0(u_{t+1})-f_0(\widetilde{u}_{t+1})-f_0^{\prime}(\widetilde{u}_{t+1})(u_{t+1}-\widetilde{u}_{t+1})\notag\\
        =&\frac{1}{1-\phi(u_{t+1})}-\frac{1}{1-\phi(\widetilde{u}_{t+1})}-\left(\frac{1}{1-\phi(\widetilde{u}_{t+1})}\right)^{\prime}(u_{t+1}-\widetilde{u}_{t+1})\notag\\
        \leq&\bb_2\cdot \frac{\widetilde{u}_{t+1}-u_{t+1}}{(1+u_{t+1})(1+\widetilde{u}_{t+1})}\enspace,\notag\\
        \mathcal{B}_{\psi}(u_{t+1}|\widetilde{u}_{t+1})=&\psi(u_{t+1})-\psi(\widetilde{u}_{t+1})-\psi^{\prime}(\widetilde{u}_{t+1})(u_{t+1}-\widetilde{u}_{t+1})\notag\\
        =&\frac{1}{2}u_{t+1}^2+u_{t+1}^3-\frac{1}{2}\widetilde{u}_{t+1}^2-\widetilde{u}_{t+1}^3-(\widetilde{u}_{t+1}+3\widetilde{u}_{t+1}^2)(u_{t+1}-\widetilde{u}_{t+1})\notag\\
        =&\frac{1}{2}(u_{t+1}-\widetilde{u}_{t+1})^2(1+2u_{t+1}+4\widetilde{u}_{t+1})\enspace.
    \end{align} 
    On the other hand, by the definition of $p_{t+1}$ and $\widetilde{p}_{t+1}$, they should satisfy the following first-order equations:
    \begin{align*}
        -\frac{\widehat{A}_{t}(1)}{\widetilde{p}_{t+1}^2}+\frac{\widehat{A}_{t}(0)}{(1-\widetilde{p}_{t+1})^2}+\eta_t^{-1}\Psi^{\prime}(\widetilde{p}_{t+1})&=0\enspace,\\
        -\frac{\widehat{A}_t(1)}{p_{t+1}^2}+\frac{\widehat{A}_t(0)}{(1-p_{t+1})^2}+\eta_{t+1}^{-1}\Psi^{\prime}(p_{t+1})&=0\enspace.
    \end{align*}
    By Lemma \ref{lemma:u}-3, we have $\widetilde{u}_{t+1}-u_{t+1}\lesssim \frac{1}{1-\widetilde{p}_{t+1}}-\frac{1}{1-p_{t+1}}$, $u_{t+1}-0\geq \frac{1}{\bb_1}\left(\frac{1}{1-p_{t+1}}-\frac{1}{1-1/2}\right)$ and $u_{t+1}-0\leq \frac{2}{\bb_3}\left(\frac{1}{1-p_{t+1}}-\frac{1}{1-1/2}\right)$, which implies that
    \begin{align*}
        1+\widetilde{u}_{t+1}\leq& 1+\frac{2}{\bb_3}\left(\frac{1}{1-\widetilde{p}_{t+1}}-2\right)\lesssim 1+\left(\frac{1}{1-\widetilde{p}_{t+1}}-2\right)=\frac{\widetilde{p}_{t+1}}{1-\widetilde{p}_{t+1}}\enspace,\\
        1+\widetilde{u}_{t+1}\geq& 1+\frac{1}{\bb_1}\left(\frac{1}{1-\widetilde{p}_{t+1}}-2\right)\gtrsim 1+\left(\frac{1}{1-\widetilde{p}_{t+1}}-2\right)=\frac{\widetilde{p}_{t+1}}{1-\widetilde{p}_{t+1}}\enspace,\\
        1+u_{t+1}\gtrsim&\frac{p_{t+1}}{1-p_{t+1}}\enspace.
    \end{align*}
    Note that $\widetilde{H}_{t+1}=\widehat{A}_t(1)\cdot f_1+\widehat{A}_t(0)\cdot f_0+\eta_t^{-1}\cdot \psi$.
    By the linearity of Bregman divergence, \eqref{thm:exact-asymptotic-variance_eq3}, Lemma \ref{lemma:u}, Lemma \ref{lemma:difference-inverse-probability-main} and Lemma \ref{lemma:difference-inverse-probability-step}, we have
    \begin{align*}
        &\mathcal{B}_{\,\widetilde{\!H}_{t+1}}(u_{t+1}|\widetilde{u}_{t+1})\\
        =&\widehat{A}_t(1)\cdot \mathcal{B}_{f_1}(u_{t+1}|\widetilde{u}_{t+1})+\widehat{A}_t(0)\cdot \mathcal{B}_{f_0}(u_{t+1}|\widetilde{u}_{t+1})+\eta_t^{-1}\mathcal{B}_{\psi}(u_{t+1}|\widetilde{u}_{t+1})\\
        \lesssim &\widehat{A}_t(1)\cdot \frac{\widetilde{u}_{t+1}-u_{t+1}}{(1+u_{t+1})(1+\widetilde{u}_{t+1})}+\underbrace{\widehat{A}_t(0)}_{\leq \widehat{A}_t(1)}\cdot \frac{\widetilde{u}_{t+1}-u_{t+1}}{(1+u_{t+1})(1+\widetilde{u}_{t+1})}+\eta_t^{-1}(\widetilde{u}_{t+1}-u_{t+1})^2(1+\widetilde{u}_{t+1}+\underbrace{u_{t+1}}_{\leq \widetilde{u}_{t+1}})\\
        \lesssim&\widehat{A}_t(1)\cdot \frac{\widetilde{u}_{t+1}-u_{t+1}}{(1+u_{t+1})(1+\widetilde{u}_{t+1})}+\eta_t^{-1}(\widetilde{u}_{t+1}-u_{t+1})^2(1+\widetilde{u}_{t+1})\\
        \intertext{Using the bounds $\widetilde{u}_{t+1}-u_{t+1}\lesssim \frac{1}{1-\widetilde{p}_{t+1}}-\frac{1}{1-p_{t+1}}$, $1+u_{t+1}\gtrsim p_{t+1}/(1-p_{t+1})$, $1+\widetilde{u}_{t+1}\gtrsim \widetilde{p}_{t+1}/(1-\widetilde{p}_{t+1})$ and $1+\widetilde{u}_{t+1}\lesssim \widetilde{p}_{t+1}/(1-\widetilde{p}_{t+1})$, we obtain}
        \lesssim&\widehat{A}_t(1)\cdot \left(\frac{\widetilde{p}_{t+1}}{1-\widetilde{p}_{t+1}}\right)^{-1}\left(\frac{p_{t+1}}{1-p_{t+1}}\right)^{-1}\underbrace{\left(\frac{1}{1-\widetilde{p}_{t+1}}-\frac{1}{1-p_{t+1}}\right)}_{\text{Lemma \ref{lemma:difference-inverse-probability-step}}}+\eta_t^{-1}\underbrace{\left(\frac{1}{1-\widetilde{p}_{t+1}}-\frac{1}{1-p_{t+1}}\right)^2}_{\text{Lemma \ref{lemma:difference-inverse-probability-main}}}\left(\frac{\widetilde{p}_{t+1}}{1-\widetilde{p}_{t+1}}\right)\\
        \lesssim&\widehat{A}_t(1)\cdot \left(\frac{\widetilde{p}_{t+1}}{1-\widetilde{p}_{t+1}}\right)^{-1}\left(\frac{p_{t+1}}{1-p_{t+1}}\right)^{-1}\cdot \frac{\eta_{t+1}^{-1}-\eta_{t}^{-1}}{\widehat{A}_t(1)}\left(\frac{\widetilde{p}_{t+1}}{1-\widetilde{p}_{t+1}}\right)\left(\frac{p_{t+1}}{1-p_{t+1}}\right)^{4}\\
        &+\eta_t^{-1}\left(\frac{\widetilde{p}_{t+1}}{1-\widetilde{p}_{t+1}}\frac{\eta_{t}\widehat{A}_t(1)-\eta_{t+1}\widehat{A}_t(1)}{\eta_{t}\widehat{A}_t(1)}+\frac{p_{t+1}}{1-p_{t+1}}\frac{\eta_{t}\widehat{A}_t(0)-\eta_{t+1}\widehat{A}_t(0)}{\eta_{t}\widehat{A}_t(0)}\right)^2\left(\frac{\widetilde{p}_{t+1}}{1-\widetilde{p}_{t+1}}\right)\\
        \intertext{We have $\frac{\eta_{t}\widehat{A}_t(1)-\eta_{t+1}\widehat{A}_t(1)}{\eta_{t}\widehat{A}_t(1)}=\frac{\eta_{t}\widehat{A}_t(0)-\eta_{t+1}\widehat{A}_t(0)}{\eta_{t}\widehat{A}_t(0)}=\frac{\eta_{t}-\eta_{t+1}}{\eta_t}=\frac{\eta_{t+1}^{-1}-\eta_{t}^{-1}}{\eta_{t+1}^{-1}}$. Since $\widetilde{p}_{t+1}\geq p_{t+1}$, we have $\frac{p_{t+1}}{1-p_{t+1}}\leq\frac{\widetilde{p}_{t+1}}{1-\widetilde{p}_{t+1}}$. Hence we can obtain}
        \lesssim&(\eta_{t+1}^{-1}-\eta_{t}^{-1})\left(\frac{p_{t+1}}{1-p_{t+1}}\right)^{3}+\eta_t^{-1}\left(\frac{\eta_{t+1}^{-1}-\eta_{t}^{-1}}{\eta_{t+1}^{-1}}\right)^2\left(\frac{\widetilde{p}_{t+1}}{1-\widetilde{p}_{t+1}}\right)^3\enspace.
    \end{align*}
    If $\widehat{A}_t(1)\leq \widehat{A}_t(0)$, by symmetry we can similarly prove that
    \begin{align*}
        \mathcal{B}_{\,\widetilde{\!H}_{t+1}}(u_{t+1}|\widetilde{u}_{t+1})\lesssim(\eta_{t+1}^{-1}-\eta_{t}^{-1})\left(\frac{1-p_{t+1}}{p_{t+1}}\right)^{3}+\eta_t^{-1}\left(\frac{\eta_{t+1}^{-1}-\eta_{t}^{-1}}{\eta_{t+1}^{-1}}\right)^2\left(\frac{1-\widetilde{p}_{t+1}}{\widetilde{p}_{t+1}}\right)^3\enspace.
    \end{align*}
    By applying Lemma \ref{lemma:p-power-moment}, Corollary \mainref{corollary:p-moment}*, and Lemma \ref{lemma:squared-residuals-random}, we have
    \begin{align*}
        &\frac{1}{T}\e\left[\sum_{t=1}^T\mathcal{B}_{\,\widetilde{\!H}_{t+1}}(u_{t+1}|\widetilde{u}_{t+1})\right]\\
        \lesssim&\frac{1}{T}\sum_{t=1}^T(\eta_{t+1}^{-1}-\eta_{t}^{-1})\underbrace{\e\left[\left(\frac{1}{1-p_{t+1}}\right)^{3}+\left(\frac{1}{p_{t+1}}\right)^{3}\right]}_{\text{Lemma \ref{lemma:p-power-moment}}}+\frac{1}{T}\sum_{t=1}^T\eta_t^{-1}\left(\frac{\eta_{t+1}^{-1}-\eta_{t}^{-1}}{\eta_{t+1}^{-1}}\right)^2\underbrace{\e\left[\left(\frac{1}{1-\widetilde{p}_{t+1}}\right)^3+\left(\frac{1}{\widetilde{p}_{t+1}}\right)^3\right]}_{\text{Corollary \mainref{corollary:p-moment}* and Lemma \ref{lemma:squared-residuals-random}}}\\
        \lesssim&\frac{1}{T}\sum_{t=1}^T(\eta_{t+1}^{-1}-\eta_{t}^{-1})\eta_{t+1}^{3/4}T^{3/4}+\frac{1}{T}\sum_{t=1}^T\eta_t^{-1}\underbrace{\left(\frac{\eta_{t+1}^{-1}-\eta_{t}^{-1}}{\eta_{t+1}^{-1}}\right)^2}_{\leq \frac{\eta_{t+1}^{-1}-\eta_{t}^{-1}}{\eta_{t+1}^{-1}}}\eta_{t}^{3/4}T^{3/4}\\
        \leq&T^{-1/4}\left[\sum_{t=1}^T(\eta_{t+1}^{-1}-\eta_{t}^{-1})(\eta_{t+1}^{-1})^{-3/4}+\sum_{t=1}^T\eta_t^{-1/4}(\eta_{t+1}^{-1}-\eta_{t}^{-1})\eta_{t+1}\right]\\
        \leq&T^{-1/4}\left[\sum_{t=1}^T(\eta_{t+1}^{-1}-\eta_{t}^{-1})(\eta_{t+1}^{-1})^{-3/4}+\sum_{t=1}^T(\eta_{t+1}^{-1}-\eta_{t}^{-1})(\eta_{t+1}^{-1})^{-3/4}\right]~~~~(\text{since $\eta_t\geq \eta_{t+1}$})\\
        \leq&2T^{-1/4}\int_{\eta_{1}^{-1}}^{\eta_{T+1}^{-1}}x^{-3/4}\mathrm{d}x\quad(\text{by integral comparison})\\
        \lesssim&T^{-1/4}\eta_{T+1}^{-1/4}\\
        \rightarrow&0~~~~(\text{Assumption \mainref{assumption:maximum-radius-strengthened} and Lemma \ref{lemma:R}})\enspace.
    \end{align*}
    Therefore,
    \begin{align}\label{thm:exact-asymptotic-variance_eq4}
        |S_{1}|=\frac{1}{T}\e\left[\sum_{t=1}^T\mathcal{B}_{\,\widetilde{\!H}_{t+1}}(u_{t+1}|\widetilde{u}_{t+1})\right]=o(1)\enspace.
    \end{align}
    \textbf{Step 2:} Prove $S_2=o(1)$.\\[3mm]
    If $p_{t+1}\geq 1/2$, by Lemma \ref{lemma:u}-3, we have
    \begin{align*}
        0\leq \phi^{-1}(p_{t+1})=u_{t+1}-0\leq \frac{2}{\bb_3}\left(\frac{1}{1-p_{t+1}}-\frac{1}{1-1/2}\right)\leq \frac{2}{\bb_3}\cdot \frac{1}{1-p_{t+1}}\lesssim \frac{1}{1-p_{t+1}}\enspace.
    \end{align*}
    Similarly, we can prove that if $p_{t+1}\leq 1/2$, we have $0\leq -\phi^{-1}(p_{t+1})\lesssim 1/p_{t+1}$. 
    Hence by Lemma \ref{lemma:p-power-moment}, we have
    \begin{align}\label{thm:exact-asymptotic-variance_eq5}
        |S_{2}|=&\frac{1}{T}\sum_{t=1}^T\left(\eta_{t+1}^{-1}-\eta_t^{-1}\right)\e\left[\Psi(p_{t+1})\right]\notag\\
        =&\frac{1}{T}\sum_{t=1}^T\left(\eta_{t+1}^{-1}-\eta_t^{-1}\right)\e\left[\frac{1}{2}(\phi^{-1}(p_{t+1}))^2+|\phi^{-1}(p_{t+1})|^3\right]\notag\\
        \lesssim&\frac{1}{T}\sum_{t=1}^T\left(\eta_{t+1}^{-1}-\eta_t^{-1}\right)\underbrace{\left(\e\left[\frac{1}{p_{t+1}^2}+\frac{1}{(1-p_{t+1})^2}\right]+\e\left[\frac{1}{p_{t+1}^3}+\frac{1}{(1-p_{t+1})^3}\right]\right)}_{\text{Lemma \ref{lemma:p-power-moment}}}\notag\\
        \lesssim&\frac{1}{T}\sum_{t=1}^T\left(\eta_{t+1}^{-1}-\eta_t^{-1}\right)\eta_{t+1}^{3/4}T^{3/4}\notag\\
        \leq&T^{-1/4}\sum_{t=1}^T(\eta_{t+1}^{-1}-\eta_{t}^{-1})(\eta_{t+1}^{-1})^{-3/4}\notag\\
        \leq&T^{-1/4}\int_{\eta_{1}^{-1}}^{\eta_{T+1}^{-1}}x^{-3/4}\mathrm{d}x\quad(\text{integral comparison})\notag\\
        \lesssim&T^{-1/4}\eta_{T+1}^{-1/4}\notag\\
        \rightarrow&0~~~~(\text{Assumption \mainref{assumption:maximum-radius-strengthened} and Lemma \ref{lemma:R}})\enspace.
    \end{align}
    \textbf{Step 3:} Prove $S_3=o(1)$.\\[3mm]
    Let $\widetilde{u}=\phi^{-1}(\widetilde{p})$. Then by Lemma \mainref{lemma:breg-lb}, we have
    \begin{align}\label{thm:exact-asymptotic-variance_eq6}
        \qquad 0\leq&\frac{1}{T}\e\left[\widehat{F}_{T+1}(\widetilde{p})-\widehat{F}_{T+1}(p_{T+1})\right]\quad(\text{optimality of $p_{T+1}$})\notag\\
        =&\frac{1}{T}\e\left[\widehat{H}_{T+1}(\widetilde{u})-\widehat{H}_{T+1}(u_{T+1})\right]\quad(\text{sigmoid transformation})\notag\\
        =&\frac{1}{T}\e\left[\langle \nabla \widehat{H}_{T+1}(\widetilde{u}),\widetilde{u}-u_{T+1}\rangle-\mathcal{B}_{\,\widehat{\!H}_{T+1}}(u_{T+1}|\widetilde{u})\right]\notag\\
        \leq&\frac{1}{T}\e\left[\langle \nabla \widehat{H}_{T+1}(\widetilde{u}),\widetilde{u}-u_{T+1}\rangle-\eta_T^{-1}\mathcal{B}_{\psi}(u_{T+1}|\widetilde{u})\right]\quad(\text{since $\widehat{H}_{T+1}-\eta_T^{-1}\psi$ is convex})\notag\\
        \leq&\frac{1}{T}\e\left[\langle \nabla \widehat{H}_{T+1}(\widetilde{u}),\widetilde{u}-u_{T+1}\rangle-\frac{\eta_T^{-1}}{2}(u_{T+1}-\widetilde{u})^2(1+|\widetilde{u}|)\right]\quad(\text{Lemma \mainref{lemma:breg-lb}})\notag\\
        \leq&\frac{1}{T}\e\left[\frac{\eta_T}{2(1+|\widetilde{u}|)}(\nabla \widehat{H}_{T+1}(\widetilde{u}))^2\right]\quad(\text{complete the square})\notag\\
        \lesssim&T^{-1}\eta_{T}\e\left[(\nabla \widehat{H}_{T+1}(\widetilde{u}))^2\right]\enspace.
    \end{align}
    The goal is therefore to upper bound $\e[\nabla \widehat{H}_{T+1}(\widetilde{u})^2]$.
    By definition, $\widetilde{u}=\phi^{-1}(\widetilde{p})$ should satisfy the following first-order equation:
    \begin{align}\label{thm:exact-asymptotic-variance_eq7}
        \qquad~~&\sum_{t=1}^T\Bigg[\frac{\mathbf{1}[Z_t=1]}{p_t}\cdot\Delta_{t,1}(\bv_t^*(1))^2\cdot\left(\frac{1}{\phi(\widetilde{u})}\right)^{\prime}+\frac{\mathbf{1}[Z_t=0]}{1-p_t}\cdot\Delta_{t,0}(\bv_t^*(0))^2\cdot \left(\frac{1}{1-\phi(\widetilde{u})}\right)^{\prime}\Bigg]\notag\\
        &+\eta_{T}^{-1}\psi^{\prime}(\widetilde{u})=0\enspace.
    \end{align}
    For any $t\in[T]$ and $k\in\{0,1\}$, define $\widetilde{\Delta}_{t,k}=\Delta_{t,k}(\bv_t(k))^2-\Delta_{t,k}(\bv_t^*(k))^2$.
    By the definition of $\widehat{H}_{T+1}$, we have
    \begin{align}\label{thm:exact-asymptotic-variance_eq8}
        &\e\left[(\nabla \widehat{H}_{T+1}(\widetilde{u}))^2\right]\notag\\
        =&\e\left[\left(\sum_{t=1}^T\nabla \widehat{h}_{t}(\widetilde{u})+\eta_{T}^{-1}\psi^{\prime}(\widetilde{u})\right)^2\right]\notag\\
        =&\e\Bigg[\Bigg(\sum_{t=1}^T\frac{\mathbf{1}[Z_t=1]}{p_t}\cdot\Delta_{t,1}(\bv_t(1))^2\cdot\left(\frac{1}{\phi(\widetilde{u})}\right)^{\prime}\notag\\
        &+\sum_{t=1}^T\frac{\mathbf{1}[Z_t=0]}{1-p_t}\cdot\Delta_{t,0}(\bv_t(0))^2\cdot\left(\frac{1}{1-\phi(\widetilde{u})}\right)^{\prime}+\eta_{T}^{-1}\psi^{\prime}(\widetilde{u})\Bigg)^2\Bigg]\notag\\
        \intertext{By subtracting the first-order equation \eqref{thm:exact-asymptotic-variance_eq7} from the terms in the bracket, we can simplify this as}
        =&\e\Bigg[\Bigg(\sum_{t=1}^T\frac{\mathbf{1}[Z_t=1]}{p_t}\cdot\widetilde{\Delta}_{t,1}\cdot\left(\frac{1}{\phi(\widetilde{u})}\right)^{\prime}+\sum_{t=1}^T\frac{\mathbf{1}[Z_t=0]}{1-p_t}\cdot\widetilde{\Delta}_{t,0}\cdot\left(\frac{1}{1-\phi(\widetilde{u})}\right)^{\prime}\Bigg)^2\Bigg]\notag\\
        \intertext{Using the inequality $(a+b)^2\leq 2a^2+2b^2$ and pulling out the common factor $\left[\left(\frac{1}{\phi(\widetilde{u})}\right)^{\prime}\right]^2$ and $\left[\left(\frac{1}{1-\phi(\widetilde{u})}\right)^{\prime}\right]^2$, we obtain}
        \lesssim&\e\left[\left[\left(\frac{1}{\phi(\widetilde{u})}\right)^{\prime}\right]^2\left(\sum_{t=1}^T\frac{\mathbf{1}[Z_t=1]}{p_t}\cdot\widetilde{\Delta}_{t,1}\right)^2\right]\notag\\
        &+\e\left[\left[\left(\frac{1}{1-\phi(\widetilde{u})}\right)^{\prime}\right]^2\left(\sum_{t=1}^T\frac{\mathbf{1}[Z_t=0]}{1-p_t}\cdot\widetilde{\Delta}_{t,0}\right)^2\right]\notag\\
        \intertext{Since $\left|\left(\frac{1}{\phi(\widetilde{u})}\right)^{\prime}\right|\leq \bb_1$ and $\left|\left(\frac{1}{1-\phi(\widetilde{u})}\right)^{\prime}\right|=\left|\left(\frac{1}{\phi(-\widetilde{u})}\right)^{\prime}\right|\leq \bb_1$ by Condition \mainref{condition:sigmoid}, we obtain}
        \lesssim&\e\left[\left(\sum_{t=1}^T\frac{\mathbf{1}[Z_t=1]}{p_t}\cdot\widetilde{\Delta}_{t,1}\right)^2\right]+\e\left[\left(\sum_{t=1}^T\frac{\mathbf{1}[Z_t=0]}{1-p_t}\cdot\widetilde{\Delta}_{t,0}\right)^2\right]\notag\\
        \intertext{Using the decomposition $\frac{\mathbf{1}[Z_t=1]}{p_t}=\left(\frac{\mathbf{1}[Z_t=1]}{p_t}-1\right)+1$ and the inequality $(a+b)^2\leq 2a^2+2b^2$, we obtain}
        \lesssim&\underbrace{\e\left[\left(\sum_{t=1}^T\left(\frac{\mathbf{1}[Z_t=1]}{p_t}-1\right)\cdot \widetilde{\Delta}_{t,1}\right)^2\right]}_{:=B_1}+\underbrace{\e\left[\left(\sum_{t=1}^T\left(\frac{\mathbf{1}[Z_t=0]}{1-p_t}-1\right)\cdot \widetilde{\Delta}_{t,0}\right)^2\right]}_{:=B_2}\notag\\
        &+\underbrace{\e\left[\left(\sum_{t=1}^T\widetilde{\Delta}_{t,1}\right)^2\right]}_{:=B_3}+\underbrace{\e\left[\left(\sum_{t=1}^T\widetilde{\Delta}_{t,0}\right)^2\right]}_{:=B_4}\enspace.
    \end{align}
    Since the terms $\left(\frac{\mathbf{1}[Z_t=1]}{p_t}-1\right)\cdot\widetilde{\Delta}_{t,1}=\left(\frac{\mathbf{1}[Z_t=1]}{p_t}-1\right)\cdot\left[\Delta_{t,1}(\bv_t(1))^2-\Delta_{t,1}(\bv_t^*(1))^2\right]$ form a martingale difference sequence and are orthogonal in $L^2$, we obtain
    \begin{align}\label{thm:exact-asymptotic-variance_eq9}
        \quad B_1=&\e\left[\left(\sum_{t=1}^T\left(\frac{\mathbf{1}[Z_t=1]}{p_t}-1\right)\cdot\widetilde{\Delta}_{t,1}\right)^2\right]\notag\\
        =&\operatorname{Var}\left(\sum_{t=1}^T\left(\frac{\mathbf{1}[Z_t=1]}{p_t}-1\right)\cdot \widetilde{\Delta}_{t,1}\right)\notag\\
        =&\sum_{t=1}^T\operatorname{Var}\left[\left(\frac{\mathbf{1}[Z_t=1]}{p_t}-1\right)\cdot\widetilde{\Delta}_{t,1}\right]\quad(\text{property of $L^2$ martingale difference sequence})\notag\\
        \intertext{Using the inequality $\operatorname{Var}(X+Y)\leq 2\operatorname{Var}(X)+2\operatorname{Var}(Y)$ and $\widetilde{\Delta}_{t,1}=\Delta_{t,1}(\bv_t(1))^2-\Delta_{t,1}(\bv_t^*(1))^2$, we have}
        \lesssim&\underbrace{\sum_{t=1}^T\operatorname{Var}\left(\left(\frac{\mathbf{1}[Z_t=1]}{p_t}-1\right)\cdot\Delta_{t,1}(\bv_t(1))^2\right)}_{:=B_{1,1}}\notag\\
        &+\underbrace{\sum_{t=1}^T\operatorname{Var}\left(\left(\frac{\mathbf{1}[Z_t=1]}{p_t}-1\right)\cdot \Delta_{t,1}(\bv_t^*(1))^2\right)}_{:=B_{1,2}}\enspace.
    \end{align}
    In the proof of Lemma \ref{lemma:variance-estimated-squared-residual}, we have shown that
    \begin{align}\label{thm:exact-asymptotic-variance_eq10}
        B_{1,1}=& o(T^{5/4}R_T^{3/2})=o(T^{3/2}R_T)\enspace,
    \end{align}
    where the last steps follow from Assumption \mainref{assumption:maximum-radius-strengthened} and Lemma \ref{lemma:R}.
    On the other hand, by direct calculation, we have
    \begin{align}\label{thm:exact-asymptotic-variance_eq100}
        B_{1,2}=&\sum_{t=1}^T\operatorname{Var}\left(\left(\frac{\mathbf{1}[Z_t=1]}{p_t}-1\right)\cdot \Delta_{t,1}(\bv_{t}^*(1))^2\right)\notag\\
        =&\sum_{t=1}^T\Delta_{t,1}(\bv_{t}^*(1))^4\cdot\operatorname{Var}\left(\frac{\mathbf{1}[Z_t=1]}{p_t}-1\right)\notag\\
        =&\sum_{t=1}^T\Delta_{t,1}(\bv_{t}^*(1))^4\cdot\e\left[\frac{1-p_t}{p_t}\right]\notag\\
        \lesssim&\sum_{t=1}^T\Delta_{t,1}(\bv_{t}^*(1))^4\cdot(\eta_tT)^{1/4}\quad(\text{Lemma \ref{lemma:p-power-moment}})\notag\\
        \leq&\eta_T^{-3/4}T^{1/4}\cdot\sum_{t=1}^T\eta_t\Delta_{t,1}(\bv_{t}^*(1))^4\quad(\text{since $\eta_t$ is nonincreasing})\notag\\
        \lesssim&\eta_T^{-3/4}T^{1/4}\cdot T^{1/2}R_T\quad(\text{Corollary \ref{corollary:fourth-moment-deterministic}})\notag\\
        \leq&T^{9/8}R_T^{7/4}\notag\\
        =&o(T^{3/2}R_T)\quad(\text{\mainref{assumption:maximum-radius-strengthened} and Lemma \ref{lemma:R}})\enspace.
    \end{align}
    Hence by \eqref{thm:exact-asymptotic-variance_eq9}, \eqref{thm:exact-asymptotic-variance_eq10}, and \eqref{thm:exact-asymptotic-variance_eq100}, we obtain
    \begin{align*}
        B_1\lesssim B_{1,1}+B_{1,2}=o(T^{3/2}R_T)\enspace.
    \end{align*}
    Similarly, we can prove that $B_2=o(T^{3/2}R_T)$.
    In the proof of Lemma \ref{lemma:variance-estimated-squared-residual}, we have already verified that
    \begin{align*}
        B_{3}=&\e\left[\left(\sum_{t=1}^T\left[\Delta_{t,1}(\bv_t(1))^2-\Delta_{t,1}(\bv_t^*(1))^2\right]\right)^2\right]\\
        \lesssim& \eta_T^{-2}(\eta_TT)^{1/2}\log^{2}(\eta_TT)(\eta_TT)^{25/64}\\
        =& T^{3/2}R_T\cdot (\eta_TT)^{-1+1/2+25/64}\log^{2}(\eta_TT)\\
        =&o(T^{3/2}R_T)\enspace,
    \end{align*}
    where the last inequality follows from Assumption \mainref{assumption:maximum-radius-strengthened} and Lemma \ref{lemma:R}.
    Similarly, we can prove that $B_4=o(T^{3/2}R_T)$.
    Hence by \eqref{thm:exact-asymptotic-variance_eq8}, we have $\e\left[\nabla \widehat{H}_{T+1}(\widetilde{u})^2\right]\lesssim B_1+B_2+B_3+B_4=o(T^{3/2}R_T)$.
    Then by \eqref{thm:exact-asymptotic-variance_eq6}, we have
    \begin{align}\label{thm:exact-asymptotic-variance_eq11}
        |S_{3}|\lesssim  T^{-1}\eta_{T}\cdot \e\left[\nabla \widehat{H}_{T+1}(\widetilde{u})^2\right]=T^{-1}T^{-1/2}R_T^{-1}\cdot o(T^{3/2}R_T)=o(1)\enspace.
    \end{align}
    \textbf{Step 4:} Prove $S_4=o(1)$.\\[3mm]
    We have the following decomposition for $S_4$:
    \begin{align}\label{thm:exact-asymptotic-variance_eq12}
        S_{4}=&\frac{1}{T}\sum_{t=1}^T\e\left[\frac{\mathbf{1}[Z_t=1]}{p_t}\cdot\Delta_{t,1}(\bv_t(1))^2\cdot \frac{1}{\widetilde{p}}+\frac{\mathbf{1}[Z_t=0]}{1-p_t}\cdot\Delta_{t,0}(\bv_t(0))^2\cdot \frac{1}{1-\widetilde{p}}\right]\notag\\
        &-\frac{1}{T}\sum_{t=1}^T\e\left[\frac{\mathbf{1}[Z_t=1]}{p_t}\cdot\Delta_{t,1}(\bv_t(1))^2\cdot \frac{1}{\breve{p}}+\frac{\mathbf{1}[Z_t=0]}{1-p_t}\cdot\Delta_{t,0}(\bv_t(0))^2\cdot \frac{1}{1-\breve{p}}\right]\notag\\
        &+(\eta_{T}T)^{-1}\e\left[\Psi(\widetilde{p})-\Psi(\breve{p})\right]\notag\\
        =&\underbrace{\frac{1}{T}\sum_{t=1}^T\e\left[\frac{\mathbf{1}[Z_t=1]}{p_t}\cdot\Delta_{t,1}(\bv_t(1))^2\cdot \left(\frac{1}{\widetilde{p}}-\frac{1}{\breve{p}}\right)\right]}_{:=S_{4,1}}\notag\\
        &+\underbrace{\frac{1}{T}\sum_{t=1}^T\e\left[\frac{\mathbf{1}[Z_t=0]}{1-p_t}\cdot\Delta_{t,0}(\bv_t(0))^2\cdot \left(\frac{1}{1-\widetilde{p}}-\frac{1}{1-\breve{p}}\right)\right]}_{:=S_{4,2}}\notag\\
        &+\underbrace{(\eta_{T}T)^{-1}\e\left[\Psi(\widetilde{p})-\Psi(\breve{p})\right]}_{:=S_{4,3}}\enspace.
    \end{align}
    We then bound $S_{4,1}$, $S_{4,2}$ and $S_{4,3}$ separately.
    Using the inequality $(a+b)^2\leq 2a^2+2b^2$, we can bound and decompose $|S_{4,1}|$ as:
    \begin{align}\label{thm:exact-asymptotic-variance_eq13}
        |S_{4,1}|=&\left|\frac{1}{T}\sum_{t=1}^T\e\left[\frac{\mathbf{1}[Z_t=1]}{p_t}\cdot\Delta_{t,1}(\bv_t(1))^2\cdot \left(\frac{1}{\widetilde{p}}-\frac{1}{\breve{p}}\right)\right]\right|\notag\\
        \leq&\frac{1}{T}\sum_{t=1}^T\e\left[\frac{\mathbf{1}[Z_t=1]}{p_t}\cdot\Delta_{t,1}(\bv_t(1))^2\cdot \left|\frac{1}{\widetilde{p}}-\frac{1}{\breve{p}}\right|\right]\notag\\
        =&\frac{1}{T}\sum_{t=1}^T\e\left[\frac{\mathbf{1}[Z_t=1]}{p_t}\cdot\left(\Delta_{t,1}(\bv_t^*(1))+\iprod{\xv_t,\bv_t^*(1)-\bv_t(1)}\right)^2\cdot \left|\frac{1}{\widetilde{p}}-\frac{1}{\breve{p}}\right|\right]\notag\\
        \lesssim&\underbrace{\frac{1}{T}\sum_{t=1}^T\e\left[\frac{\mathbf{1}[Z_t=1]}{p_t}\cdot\Delta_{t,1}(\bv_t^*(1))^2\cdot \left|\frac{1}{\widetilde{p}}-\frac{1}{\breve{p}}\right|\right]}_{:=D_1}\notag\\
        &+\underbrace{\frac{1}{T}\sum_{t=1}^T\e\left[\frac{\mathbf{1}[Z_t=1]}{p_t}\cdot\iprod{\xv_t,\bv_t(1)-\bv^*_t(1)}^2\cdot \left|\frac{1}{\widetilde{p}}-\frac{1}{\breve{p}}\right|\right]}_{:=D_2}\enspace.
    \end{align}
    Before bounding $D_1$, we first bound the variance of $\widetilde{A}_T(1)$ and $\widetilde{A}_T(0)$.
    By \eqref{thm:exact-asymptotic-variance_eq100}, we have
    \begin{align}\label{lemma:difference-inverse-probability-last_eq1}
        \Var{\widetilde{A}_T(1)}=&\operatorname{Var}\left(\sum_{t=1}^T\frac{\mathbf{1}[Z_t=1]}{p_t}\cdot \Delta_{t,1}(\bv_{t}^*(1))^2\right)\notag\\
        =&\operatorname{Var}\left(\sum_{t=1}^T\left(\frac{\mathbf{1}[Z_t=1]}{p_t}-1\right)\cdot \Delta_{t,1}(\bv_{t}^*(1))^2\right)\notag\\
        \intertext{Since the terms $\left(\frac{\mathbf{1}[Z_t=1]}{p_t}-1\right)\cdot (y_t(1)-\iprod{\xv_t,\bv_t^*(1)})^2$ form a martingale difference sequence and are orthogonal in $L^2$, the variance of the sum equals the sum of variances.}
        =&\sum_{t=1}^T\operatorname{Var}\left(\left(\frac{\mathbf{1}[Z_t=1]}{p_t}-1\right)\cdot \Delta_{t,1}(\bv_{t}^*(1))^2\right)\notag\\
        \lesssim&T^{9/8}R_T^{7/4}\quad(\text{by \eqref{thm:exact-asymptotic-variance_eq100}})\enspace.
    \end{align}
    Similarly, we can show that $\Var{\widetilde{A}_T(0)}\lesssim T^{9/8}R_T^{7/4}$.
    By Lemma \ref{lemma:difference-inverse-probability-last} and the definition of $\widetilde{A}_T(1)$, we have
    \begin{align}\label{thm:exact-asymptotic-variance_eq14}
        D_1=&\frac{1}{T}\e\left[\widetilde{A}_T(1)\left|\frac{1}{\widetilde{p}}-\frac{1}{\breve{p}}\right|\right]\notag\\
        \lesssim&\frac{1}{T}\e\left[|\widetilde{A}_T(1)-A^*_T(1)|\right]+\frac{1}{T}\e\left[|\widetilde{A}_T(0)-A^*_T(0)|\right]\quad(\text{Lemma \ref{lemma:difference-inverse-probability-last}})\notag\\
        \lesssim&\frac{1}{T}\left(\e\left[(\widetilde{A}_T(1)-A^*_T(1))^2\right]^{1/2}+\e\left[(\widetilde{A}_T(0)-A^*_T(0))^2\right]^{1/2}\right)\quad(\text{Jensen's inequality})\notag\\
        =&\frac{1}{T}\left(\Var{\widetilde{A}_T(1)}^{1/2}+\Var{\widetilde{A}_T(0)}^{1/2}\right)\quad(\text{since $\e[\widetilde{A}_T(k)]=A_T^*(k)$})\notag\\
        \lesssim&T^{-1}\cdot T^{9/16}R_T^{7/8}\quad(\text{by \eqref{lemma:difference-inverse-probability-last_eq1}})\notag\\
        =&o(T^{-1+9/16+7/16})\quad(\text{Assumption \mainref{assumption:maximum-radius-strengthened} and Lemma \ref{lemma:R}})\notag\\
        =&o(1)\enspace.
    \end{align}
    By H{\"o}lder's inequality, we have
    \begin{align}\label{thm:exact-asymptotic-variance_eq15}
        D_{2}\leq&\frac{1}{T}\e\left[\sum_{t=1}^T\frac{\mathbf{1}[Z_t=1]}{p_t}\cdot\iprod{\xv_t,\bv_t(1)-\bv^*_t(1)}^2\cdot \left(\frac{1}{\widetilde{p}}+\frac{1}{\breve{p}}\right)\right]\notag\\
        \leq&\frac{1}{T}\e\left[\left(\sum_{t=1}^T\frac{\mathbf{1}[Z_t=1]}{p_t}\cdot\iprod{\xv_t,\bv_t(1)-\bv^*_t(1)}^2\right)^{4/3}\right]^{3/4}\left(\e\left[\frac{1}{\widetilde{p}^4}\right]^{1/4}+\e\left[\frac{1}{\breve{p}^4}\right]^{1/4}\right)\quad(\text{H{\"o}lder's inequality})\notag\\
        \lesssim&T^{-1}\e\left[\left(\sum_{t=1}^T\frac{\mathbf{1}[Z_t=1]}{p_t}\cdot\iprod{\xv_t,\bv_t(1)-\bv^*_t(1)}^2\right)\right]^{1/2}\notag\\
        &\times\e\left[\left(\sum_{t=1}^T\frac{\mathbf{1}[Z_t=1]}{p_t}\cdot\iprod{\xv_t,\bv_t(1)-\bv^*_t(1)}^2\right)^2\right]^{1/4}\times \eta_T^{1/4}\left(\e[\widetilde{A}_T(0)]+A_T^*(0)\right)^{1/4}\quad(\text{Lemma \mainref{lemma:effect-of-p-regularization}*})\notag\\
        \intertext{Since $\e[\widetilde{A}_T(0)]=A_T^*(0)$ and $A_T^*(0)=\Theta(T)$, this leads to}
        \lesssim&T^{-1+1/8}R_T^{-1/4}\underbrace{\e\left[\left(\sum_{t=1}^T\frac{\mathbf{1}[Z_t=1]}{p_t}\cdot\iprod{\xv_t,\bv_t(1)-\bv^*_t(1)}^2\right)\right]^{1/2}}_{:=D_{2,1}}\notag\\
        &\times\underbrace{\e\left[\left(\sum_{t=1}^T\frac{\mathbf{1}[Z_t=1]}{p_t}\cdot\iprod{\xv_t,\bv_t(1)-\bv^*_t(1)}^2\right)^2\right]^{1/4}}_{:=D_{2,2}}\enspace.
    \end{align}
    By the law of iterated expectations and \eqref{lemma:stable-variance-2_eq1}, we have
    \begin{align*}
        D_{2,1}=\left(\sum_{t=1}^T\e\left[\iprod{\xv_t,\bv_t(1)-\bv_t^*(1)}^2\right]\right)^{1/2}\lesssim T^{3/16}R_T^{5/8}\enspace.
    \end{align*}
    Since $\iprod{\xv_t,\bv_t(1)-\bv_t^*(1)}=\Delta_{t,1}(\bv_t^*(1))-\Delta_{t,1}(\bv_t(1))$, using the inequality $\max\{(a-b)^2,(a+b)^2\}\leq 2a^2+2b^2$, we have
    \begin{align*}
        &\e\left[\left(\sum_{t=1}^T\frac{\mathbf{1}[Z_t=1]}{p_t}\cdot\iprod{\xv_t,\bv_t(1)-\bv^*_t(1)}^2\right)^2\right]\\
        \lesssim& \e\left[\left(\sum_{t=1}^T\frac{\mathbf{1}[Z_t=1]}{p_t}\cdot \left[\Delta_{t,1}(\bv_t^*(1))^2+\Delta_{t,1}(\bv_t(1))^2\right]\right)^2\right]\\
        =&\e\left[\left(\widehat{A}_T(1)+\widetilde{A}_T(1)\right)^2\right]\\
        \lesssim&\e[\widehat{A}_T(1)^2]+\e[\widetilde{A}_T(1)^2]\\
        \leq&\operatorname{Var}(\widehat{A}_T(1))+\e[\widehat{A}_T(1)]^2+\operatorname{Var}(\widetilde{A}_T(1))+\e[\widetilde{A}_T(1)]^2
        \intertext{Lemma \ref{lemma:variance-estimated-squared-residual} implies that $\operatorname{Var}(\widehat{A}_T(1))\lesssim T^2(\eta_TT)^{-2+1/2}\log^2(\eta_TT)(\eta_TT)^{1/4\cdot 25/16}$, which is of order $o(T^2)$ by Assumption \mainref{assumption:maximum-radius-strengthened} and Lemma \ref{lemma:R}. Lemma \ref{lemma:squared-residuals-random} implies that $\e[\widehat{A}_T(1)]^2\lesssim T^2$. \eqref{lemma:difference-inverse-probability-last_eq1} implies that $\operatorname{Var}(\widetilde{A}_T(1))\lesssim T^{9/8}R_T^{7/4}=o(T^2)$ by Assumption \mainref{assumption:maximum-radius-strengthened} and Lemma \ref{lemma:R}. Corollary \ref{corollary:squared-residual-deterministic-exact-order} implies that $\e[\widetilde{A}_T(1)]^2=A_T^*(1)^2\lesssim T^2$. Hence this can be further bounded by}
        \lesssim& T^2\enspace,
    \end{align*}
    which implies that $D_{2,2}\lesssim T^{1/2}$.
    Hence, by \eqref{thm:exact-asymptotic-variance_eq15} and Assumption \mainref{assumption:maximum-radius-strengthened}, we have
    \begin{align}\label{thm:exact-asymptotic-variance_eqqq}
        D_2\lesssim T^{-1+1/8}R_T^{-1/4}\cdot T^{3/16}R_T^{5/8}\cdot T^{1/2}=T^{-3/16}R_T^{3/8}=o(1)\enspace.
    \end{align}
    By \eqref{thm:exact-asymptotic-variance_eq13}, \eqref{thm:exact-asymptotic-variance_eq14} and \eqref{thm:exact-asymptotic-variance_eqqq}, we have $S_{4,1}=o(1)$. Similarly, we can prove that $S_{4,2}=o(1)$. For $S_{4,3}$, by similar arguments as in bounding $S_{2}$, Corollary \ref{corollary:squared-residual-deterministic} and Corollary \mainref{corollary:p-moment}*, we have
    \begin{align*}
        |S_{4,3}|\leq&(\eta_TT)^{-1}\e\left[\Psi(\widetilde{p})\right]+(\eta_TT)^{-1}\e\left[\Psi(\breve{p})\right]\\
        \lesssim&(\eta_TT)^{-1}\underbrace{\e\left[\frac{1}{\widetilde{p}^3}+\frac{1}{(1-\widetilde{p})^3}\right]}_{\text{Corollary \mainref{corollary:p-moment}*}}+(\eta_TT)^{-1}\underbrace{\e\left[\frac{1}{\breve{p}^3}+\frac{1}{(1-\breve{p})^3}\right]}_{\text{Corollary \mainref{corollary:p-moment}*}}\\
        \lesssim&(\eta_TT)^{-1}\left[\eta_T^{3/4}\e[\widetilde{A}_T(0)]^{3/4}+\eta_T^{3/4}\e[\widetilde{A}_T(1)]^{3/4}+\eta_T^{3/4}A_T^{*}(0)^{3/4}+\eta_T^{3/4}A_T^{*}(1)^{3/4}\right]~~~~(\text{Corollary \mainref{corollary:p-moment}*})\\
        \lesssim&(\eta_TT)^{-1}\left[\eta_T^{3/4}(A_T^{*}(0))^{3/4}+\eta_T^{3/4}(A_T^{*}(1))^{3/4}\right]\\
        \lesssim&(\eta_TT)^{-1}\cdot (\eta_TT)^{3/4}~~~~(\text{Corollary \ref{corollary:squared-residual-deterministic}})\\
        =&(\eta_TT)^{-1/4}\\
        =&o(1)\quad(\text{Assumption \mainref{assumption:maximum-radius-strengthened} and Lemma \ref{lemma:R}})\enspace.
    \end{align*}
    Hence by \eqref{thm:exact-asymptotic-variance_eq12}, we have
    \begin{align}\label{thm:exact-asymptotic-variance_eq16}
        S_{4}=S_{4,1}+S_{4,2}+S_{4,3}=o(1)\enspace.
    \end{align}
    \textbf{Step 5:} Prove $S_5=o(1)$.\\[3mm]
    We have the following decomposition for $S_{5}$:
    \begin{align}\label{thm:exact-asymptotic-variance_eq17}
        \quad S_{5}=&\frac{1}{T}\sum_{t=1}^T\e\left[\frac{\mathbf{1}[Z_t=1]}{p_t}\cdot\Delta_{t,1}(\bv_t(1))^2\cdot \frac{1}{\breve{p}}+\frac{\mathbf{1}[Z_t=0]}{1-p_t}\cdot\Delta_{t,0}(\bv_t(0))^2\cdot \frac{1}{1-\breve{p}}\right]\notag\\
        &-\frac{1}{T}\sum_{t=1}^T\e\left[\frac{\mathbf{1}[Z_t=1]}{p_t}\cdot\Delta_{t,1}(\bv_t(1))^2\cdot \frac{1}{p^*}+\frac{\mathbf{1}[Z_t=0]}{1-p_t}\cdot\Delta_{t,0}(\bv_t(0))^2\cdot \frac{1}{1-p^*}\right]\notag\\
        &+(\eta_{T}T)^{-1}\left[\Psi(\breve{p})-
        \Psi(p^*)\right]\notag\\
        \intertext{Since $\breve{p}$ and $p^*$ are nonrandom, by applying the law of iterated expectations, we obtain}
        =&\frac{1}{T}\sum_{t=1}^T\e\left[\Delta_{t,1}(\bv_t(1))^2\cdot \frac{1}{\breve{p}}+\Delta_{t,0}(\bv_t(0))^2\cdot \frac{1}{1-\breve{p}}\right]\notag\\
        &-\frac{1}{T}\sum_{t=1}^T\e\left[\Delta_{t,1}(\bv_t(1))^2\cdot \frac{1}{p^*}+\Delta_{t,0}(\bv_t(0))^2\cdot \frac{1}{1-p^*}\right]+(\eta_{T}T)^{-1}\left[\Psi(\breve{p})-
        \Psi(p^*)\right]\notag\\
        =&\underbrace{(\eta_{T}T)^{-1}\left[\Psi(\breve{p})-
        \Psi(p^*)\right]}_{:=S_{5,1}}+\underbrace{\frac{1}{T}\sum_{t=1}^T\e\left[\Delta_{t,1}(\bv_t(1))^2\right]\cdot \left(\frac{1}{\breve{p}}-\frac{1}{p^*}\right)}_{:=S_{5,2}}\notag\\
        &+\underbrace{\frac{1}{T}\sum_{t=1}^T\e\left[\Delta_{t,0}(\bv_t(0))^2\right]\cdot \left(\frac{1}{1-\breve{p}}-\frac{1}{1-p^*}\right)}_{:=S_{5,3}}\enspace.
    \end{align}
    By similar arguments as in bounding $S_{4,3}$, we can prove that $S_{5,1}=o(1)$. 
    For $S_{5,2}$, we consider the following four cases:\\[3mm]
    \textbf{Case 1:} $\mathcal{E}(1)\geq \mathcal{E}(0)$ and $A_T^*(1)\geq A_T^*(0)$.\\[3mm]
        In this case, by definition, $\breve{p}$ and $p^*$ should satisfy:
        \begin{align*}
            -\frac{A_{T}^*(1)}{\breve{p}^2}+\frac{A_{T}^*(0)}{(1-\breve{p})^2}+\eta_{T}^{-1}\Psi^{\prime}(\breve{p})&=0\enspace,\\
            -\frac{T\cdot\olsres{1}^2}{(p^*)^2}+\frac{T\cdot\olsres{0}^2+\Delta}{(1-p^*)^2}+\eta_{T}^{-1}\Psi^{\prime}(p^*)&=0\enspace,
        \end{align*}
        where $\Delta:=-\eta_{T}^{-1}(1-p^*)^2\Psi^{\prime}(p^*)$. 
        Since $p^*=(1+\mathcal{E}(0)/\mathcal{E}(1))^{-1}$ is bounded away from 0 and 1 by Assumption \mainref{assumption:moments}, we have $|\Delta|=\bigO{\eta_{T}^{-1}}=o(T)$ by Assumption \mainref{assumption:maximum-radius-strengthened} and Lemma \ref{lemma:R}. 
        Since $A_{T}^*(1)=\Theta(T)$ and $A_{T}^*(0)=\Theta(T)$ by Corollary \ref{corollary:squared-residual-deterministic-exact-order}, we have $1/\breve{p}=\bigO{1}$ and $1/(1-\breve{p})=\bigO{1}$ by Lemma \mainref{lemma:effect-of-p-regularization}*. Then by Lemma \ref{lemma:difference-inverse-probability-main} and Corollary \ref{corollary:squared-residual-deterministic-exact-order}, we have
        \begin{align*}
            \left|\frac{1}{\breve{p}}-\frac{1}{p^*}\right|\lesssim&\frac{|A_{T}^*(1)-T\cdot\olsres{1}^2|}{T\cdot\olsres{1}^2}+\frac{|A_{T}^*(0)-T\cdot\olsres{0}^2-\Delta|}{T\cdot\olsres{1}^2}\quad(\text{Lemma \ref{lemma:difference-inverse-probability-main}})\notag\\
            \lesssim&T^{-1}\cdot o(T)+T^{-1}\cdot o(T)\quad(\text{Assumption \mainref{assumption:moments}, Corollary \ref{corollary:squared-residual-deterministic-exact-order}})\notag\\
            =&o(1)\enspace.
        \end{align*}
        \textbf{Case 2:} $\mathcal{E}(1)\leq\mathcal{E}(0)$ and $A_T^*(1)\leq A_T^*(0)$.\\[3mm]
        We can prove that $\left|\frac{1}{\breve{p}}-\frac{1}{p^*}\right|=o(1)$ by the same method as in Case 1.\\[3mm]
        \textbf{Case 3:} $\mathcal{E}(1)\geq\mathcal{E}(0)$ and $A_T^*(1)\leq A_T^*(0)$.\\[3mm]
        In this case, $\breve{p}$ and $p^*$ satisfy:
        \begin{align*}
            -\frac{A_{T}^*(1)}{\breve{p}^2}+\frac{A_{T}^*(0)}{(1-\breve{p})^2}+\eta_{T}^{-1}\Psi^{\prime}(\breve{p})&=0\enspace,\\
            -\frac{T\cdot\olsres{1}^2}{(p^*)^2}+\frac{T\cdot\olsres{0}^2+\Delta}{(1-p^*)^2}+\eta_{T}^{-1}\Psi^{\prime}(p^*)&=0\enspace,
        \end{align*}
        where $\Delta:=-\eta_{T}^{-1}(1-p^*)^2\Psi^{\prime}(p^*)=o(T)$. Then by Corollary \ref{corollary:squared-residual-deterministic-exact-order}, we have
        \begin{align*}
            0\leq A_{T}^*(0)-A_T^*(1)\leq& \underbrace{T\cdot\olsres{0}^2-T\cdot\olsres{1}^2}_{\leq 0}+|A_T^*(1)-T\cdot\olsres{1}^2|+|A_{T}^*(0)-T\cdot\olsres{0}^2|\\
            \leq&|A_T^*(1)-T\cdot\olsres{1}^2|+|A_{T}^*(0)-T\cdot\olsres{0}^2|\\
            =&o(T)\quad(\text{Corollary \ref{corollary:squared-residual-deterministic-exact-order}})\enspace.
        \end{align*}
        Then by Lemma \ref{lemma:difference-inverse-probability-half} and Corollary \ref{corollary:squared-residual-deterministic-exact-order}, we have
        \begin{align*}
            \left|\frac{1}{\breve{p}}-2\right|\leq \frac{\bb_1\bb_2}{\bb_3}\cdot\frac{A_{T}^*(0)-A_T^*(1)}{A_T^*(1)}=o(1)\enspace.
        \end{align*}
        Similarly we can prove that $T\cdot\olsres{1}^2-T\cdot\olsres{0}^2=o(T)$. Then by Assumption \mainref{assumption:moments} and Lemma \ref{lemma:difference-inverse-probability-half}, we have
        \begin{align*}
            \left|\frac{1}{p^*}-2\right|=\left|\frac{p^*-(1/2)}{p^*\cdot (1/2)}\right|\leq\left|\frac{p^*-(1/2)}{(1-p^*)\cdot (1/2)}\right|=\left|\frac{1}{1-p^*}-2\right|\lesssim \frac{T\cdot\olsres{1}^2-T\cdot\olsres{0}^2-\Delta}{T\cdot\olsres{0}^2+\Delta}=o(1)\enspace.
        \end{align*}
        Hence we have
        \begin{align*}
            \left|\frac{1}{\breve{p}}-\frac{1}{p^*}\right|\leq \left|\frac{1}{\breve{p}}-2\right|+\left|\frac{1}{p^*}-2\right|=o(1)\enspace.
        \end{align*}
        \textbf{Case 4:} $\mathcal{E}(1)\leq\mathcal{E}(0)$ and $A_T^*(1)\geq A_T^*(0)$.\\[3mm]
        We can prove that $\left|\frac{1}{\breve{p}}-\frac{1}{p^*}\right|=o(1)$ by the same method as in Case 3.
    
    By the analysis in these four cases, we have proved that $\left|\frac{1}{\breve{p}}-\frac{1}{p^*}\right|=o(1)$. 
    Then by Lemma \ref{lemma:squared-residuals-random}, we have
    \begin{align*}
        |S_{5,2}|\leq \frac{1}{T}\sum_{t=1}^T\e\left[\Delta_{t,1}(\bv_t(1))^2\right]\cdot \left|\frac{1}{\breve{p}}-\frac{1}{p^*}\right|\lesssim T^{-1}\cdot T\cdot o(1)=o(1)\enspace.
    \end{align*}
    Similarly, we can show that $S_{5,3}=o(1)$.
    Hence by \eqref{thm:exact-asymptotic-variance_eq17}, we have
    \begin{align}\label{thm:exact-asymptotic-variance_eq18}
        S_5=S_{5,1}+S_{5,2}+S_{5,3}=o(1)\enspace.
    \end{align}
    \textbf{Final Step:} Combine the results in Steps 1-5.\\[3mm]
    Combining the results in \eqref{thm:exact-asymptotic-variance_eq2}, \eqref{thm:exact-asymptotic-variance_eq4}, \eqref{thm:exact-asymptotic-variance_eq5}, \eqref{thm:exact-asymptotic-variance_eq11}, \eqref{thm:exact-asymptotic-variance_eq16}, and \eqref{thm:exact-asymptotic-variance_eq18}, we obtain
    \begin{align*}
        \frac{1}{T}\E{\regretprob_T}=o(1)-S_{1}-S_{2}+S_{3}+S_{4}+S_{5}=o(1)\enspace.
    \end{align*}
    Since $T\cdot \optvar=2(1+\rho)\mathcal{E}(1)\mathcal{E}(0)$ by Proposition \mainref{prop:oracle}, by \eqref{thm:exact-asymptotic-variance_eq1} we have established the asymptotic variance of the adaptive AIPW estimator:
    \begin{align*}
        T\cdot \operatorname{Var}(\hat{\tau})=T\cdot \optvar+\frac{1}{T}\E{\regretprob_T}+\frac{1}{T}\e[\mathcal{R}_T^{\text{pred}}]=2(1+\rho)\mathcal{E}(1)\mathcal{E}(0)+o(1)\enspace.
    \end{align*}
    Therefore, the result is proved.
\end{proof}

The explicit form of the asymptotic variance in Theorem \mainref{thm:exact-asymptotic-variance} establishes the equivalence between Assumption \mainref{assumption:bounded-correlation} and the non-superefficiency condition, which is stated in the following corollary.

\begin{refcorollary}{\mainref{corollary:non-superefficiency}}[Non-Superefficiency]
	\nonsuperefficiency
\end{refcorollary}

\begin{proof}
    It holds that $-1\leq \rho\leq 1$ by the Cauchy-Schwarz inequality. Since $c_0\leq \mathcal{E}(0),\mathcal{E}(1)\leq c_1$ by Assumption \mainref{assumption:moments}, the result is proved by Assumption \mainref{assumption:bounded-correlation} and Theorem \mainref{thm:exact-asymptotic-variance}.
\end{proof}

\subsection{Variance Estimation}\label{section:D5}
Recall that $\sqolsres{k} = \olsres{k}^2$, and the estimated squared residuals are defined as:
\begin{align*}
\estsqolsres{1}
&= \frac{1}{T} \sum_{t=1}^T Q_{t,t} Y_t^2 \frac{ \indicator{Z_t =1} }{p_t} 
+ \frac{1}{T} \sum_{t=1}^T \sum_{s \neq t} Q_{t,s} Y_t Y_s \frac{ \indicator{Z_t = 1, Z_s = 1} }{p_s p_t} 
\quadand \\
\estsqolsres{0}
&=  \frac{1}{T} \sum_{t=1}^T Q_{t,t} Y_t^2 \frac{ \indicator{Z_t = 0} }{1-p_t} 
+ \frac{1}{T} \sum_{t=1}^T \sum_{s \neq t} Q_{t,s} Y_t Y_s \frac{ \indicator{Z_t = 0, Z_s = 0} }{(1-p_s) (1-p_t)}
\enspace.
\end{align*}
In this section, we establish the consistency of the proposed variance estimator and characterize its rate of convergence.
The following lemma provides upper bounds on the deterministic terms involved in controlling the variance of the estimator.

\begin{lemma}\label{lemma:deterministic-summation-VB}
    Under Assumptions \mainref{assumption:moments}-\mainref{assumption:covariate-regularity}, for $k\in\{0,1\}$, the following holds:
    \begin{enumerate}
        \item[(1)] $\sum_{t=1}^TQ_{t,t}^2y_t(k)^4\leq c_1^4T$.
        \item[(2)] $\sum_{1\leq t_1\neq t_2\leq T}Q_{t_1,t_2}^2y_{t_1}(k)^2y_{t_2}(k)^2\left(\frac{R_{T}}{R_{t_1}}\right)^{1/2}\left(\frac{R_{T}}{R_{t_2}}\right)^{1/2}\leq c_1^4c_2R_T^2$.
    \end{enumerate}
\end{lemma}

\begin{proof}
    We only prove the result for $k=1$.
    Denote $\vec{X}=(\xv_1,\ldots,\xv_T)^{\tran}$.
    \begin{enumerate} 
        \item[(1)] We have $\mat{Q}=\mat{I}_T-\xM\left(\xM^{\tran}\xM\right)^{-1}\xM^{\tran}$ by definition. Since $\mat{Q}\preceq \mat{I}_T$, we have $0\leq Q_{t,t}\leq 1$ for any $t\in[T]$. Hence by Assumption \mainref{assumption:moments}, we have
        \begin{align*}
            \sum_{t=1}^TQ_{t,t}^2y_t(1)^4\leq c_1^4T\enspace.
        \end{align*}
        \item[(2)] Denote $\mat{H}=\xM\left(\xM^{\tran}\xM\right)^{-1}\xM^{\tran}$ and
        \begin{align*}
            \widetilde{\!\mat{H}}=\operatorname{diag}\left\{\left(\frac{R_{T}}{R_{1}}\right)^{1/2},\ldots,\left(\frac{R_{T}}{R_{T}}\right)^{1/2}\right\}\mat{H}\operatorname{diag}\left\{\left(\frac{R_{T}}{R_{1}}\right)^{1/2},\ldots,\left(\frac{R_{T}}{R_{T}}\right)^{1/2}\right\}\enspace.
        \end{align*}
        It follows that $\widetilde{\!\mat{H}}$ is positive semidefinite since $\mat{H}$ is positive semidefinite.
        Denote $\mat{Y}=(y_1(1)^2,\ldots,y_T(1)^2)^{\tran}$.
        For any $1\leq t_1\neq t_2\leq T$, by definition, it holds that $H_{t_1,t_2}=-Q_{t_1,t_2}$. Hence, we have
        \begin{align*}
            &\sum_{1\leq t_1\neq t_2\leq T}Q_{t_1,t_2}^2y_{t_1}(1)^2y_{t_2}(1)^2\left(\frac{R_{T}}{R_{t_1}}\right)^{1/2}\left(\frac{R_{T}}{R_{t_2}}\right)^{1/2}\\
            =&\sum_{1\leq t_1\neq t_2\leq T}H_{t_1,t_2}^2y_{t_1}(1)^2y_{t_2}(1)^2\left(\frac{R_{T}}{R_{t_1}}\right)^{1/2}\left(\frac{R_{T}}{R_{t_2}}\right)^{1/2}\\
            \leq&\sum_{1\leq t_1,t_2\leq T}H_{t_1,t_2}^2y_{t_1}(1)^2y_{t_2}(1)^2\left(\frac{R_{T}}{R_{t_1}}\right)^{1/2}\left(\frac{R_{T}}{R_{t_2}}\right)^{1/2}\\
            =&\mat{Y}^{\tran}(\mat{H}\circ \widetilde{\!\mat{H}})\mat{Y}\enspace,
        \end{align*}
        where $\circ$ denotes the Hadamard product.
        By Assumption \mainref{assumption:covariate-regularity}, the $t$-th diagonal element of $\widetilde{\!\mat{H}}$ satisfies:
        \begin{align*}
            \widetilde{H}_{t,t}=\left(\frac{R_{T}}{R_{t}}\right)H_{t,t}\leq \left(\frac{R_{T}}{R_{t}}\right)\|\xv_t\|^2\left\|\left(\xM^{\tran}\xM\right)^{-1}\right\|\leq c_2T^{-1}\left(\frac{R_{T}}{R_{t}}\right)R_t^2\leq c_2R_T^2T^{-1}\enspace.
        \end{align*}
        Then by Theorem 5.3.4 in \cite{horn2012matrix}, we have
        \begin{align*}
            \|\mat{H}\circ\widetilde{\!\mat{H}}\|\leq& \|\mat{H}\|\cdot\max_{t=1,\ldots,T}\widetilde{H}_{t,t}\\
            \leq&\left\|\xM\left(\xM^{\tran}\xM\right)^{-1}\xM^{\tran}\right\|\cdot c_2R_T^2T^{-1}\\
            \leq&\left\|\xM^{\tran}\xM\left(\xM^{\tran}\xM\right)^{-1}\right\|\cdot c_2R_T^2T^{-1}\\
            \leq&c_2R_T^2T^{-1}\enspace.
        \end{align*}
        Hence, we have
        \begin{align*}
            &\sum_{1\leq t_1\neq t_2\leq T}Q_{t_1,t_2}^2y_{t_1}(1)^2y_{t_2}(1)^2\left(\frac{R_{T}}{R_{t_1}}\right)^{1/2}\left(\frac{R_{T}}{R_{t_2}}\right)^{1/2}\\
            \leq& \mat{Y}^{\tran}(\mat{H}\circ \widetilde{\!\mat{H}})\mat{Y}\\
            \leq& \|\mat{H}\circ \widetilde{\!\mat{H}}\|\cdot\sum_{t=1}^Ty_t(1)^4\\
            \leq& c_1^4c_2R_T^2\quad(\text{Assumption \mainref{assumption:moments}})\enspace.
        \end{align*}
    \end{enumerate}
    Hence the lemma is proved.
\end{proof}

We now establish the following result, which can be viewed as an extension of Lemma \ref{lemma:squared-residuals-random}.
The result follows directly from Lemma \ref{lemma:variance-estimated-squared-residual} and Lemma \ref{lemma:squared-residuals-random}, which will be used to derive higher-moment bounds for the inverse assignment probabilities.

\begin{corollary}\label{corollary:estimated-squared-residuals-second-moment}
    Under Assumptions \mainref{assumption:moments}-\mainref{assumption:maximum-radius-strengthened}, and Condition \mainref{condition:sigmoid}, there exists a constant $\widetilde{\zeta}>0$ such that, for all $t\in[T]$,
    \begin{align*}
        \max\left\{\e[\widehat{A}_t(1)^2],\e[\widehat{A}_t(0)^2]\right\}\leq \widetilde{\zeta}\cdot T^2\enspace.
    \end{align*}
\end{corollary}

\begin{proof}
    We only prove the result for $\e[\widehat{A}_t(1)^2]$.
    By Lemma \ref{lemma:squared-residuals-random} and Lemma \ref{lemma:variance-estimated-squared-residual}, we have
    \begin{align*}
        \Var{\widehat{A}_t(1)}\lesssim&\eta_t^{-2}(\eta_tT)^{1/2}\log^2(\eta_tT)\left(2+\bb_1(\bb_2/6)^{1/4}\eta_t^{1/4}\e[\widehat{A}_t(0)]^{1/4}\right)^{25/16}\quad(\text{Lemma \ref{lemma:variance-estimated-squared-residual}})\\
        \lesssim&T^2\cdot (\eta_tT)^{-2}(\eta_tT)^{1/2}\log^2(\eta_tT)(\eta_tT)^{\frac{1}{4}\cdot \frac{25}{16}}\quad(\text{Lemma \ref{lemma:squared-residuals-random}})\\
        =&T^2\cdot(\eta_tT)^{-71/64}\log^2(\eta_tT)\\
        =&o(T^2)\quad(\text{Assumption \mainref{assumption:maximum-radius-strengthened} and Lemma \ref{lemma:R}})\enspace.
    \end{align*}
    Hence by Lemma \ref{lemma:squared-residuals-random} we have
    \begin{align*}
        \e[\widehat{A}_t(1)^2]=\e[\widehat{A}_t(1)]^2+\Var{\widehat{A}_t(1)}\leq \varsigma^2T^2+o(T^2)=(\varsigma^2+o(1))T^2\enspace,
    \end{align*}
    which completes the proof.
\end{proof}

Based on Corollary \ref{corollary:estimated-squared-residuals-second-moment}, we establish the following lemma.
The lemma characterizes the covariance structure of the inverse probability weighting terms, which is crucial for bounding the variance of the estimated squared OLS residuals in Theorem \mainref{thm:sq-ols-res-estimates}.
In particular, it shows that most fourth-order covariance terms vanish and provides upper bounds for the remaining nonzero covariance terms.

\begin{lemma}\label{lemma:refined-p-bound}
    Under Assumptions \mainref{assumption:moments}-\mainref{assumption:maximum-radius-strengthened}, and Condition \mainref{condition:sigmoid}, for any $s_1,s_2,s_3,s_4\in[T]$ such that  $s_1>s_2$ and $s_3>s_4$, we have the following results:
    \begin{itemize}
        \item[(1)] If $s_1,s_2,s_3,s_4$ are distinct, we have
        \begin{align*}
            \operatorname{Cov}\left(\frac{\mathbf{1}[Z_{s_1}=1]}{p_{s_1}}\frac{\mathbf{1}[Z_{s_2}=1]}{p_{s_2}},\frac{\mathbf{1}[Z_{s_3}=1]}{p_{s_3}}\frac{\mathbf{1}[Z_{s_4}=1]}{p_{s_4}}\right)=0\enspace.
        \end{align*}
        \item[(2)] If $s_1=s_3$ and $s_2=s_4$, we have
        \begin{align*}
            \left|\operatorname{Cov}\left(\frac{\mathbf{1}[Z_{s_1}=1]}{p_{s_1}}\frac{\mathbf{1}[Z_{s_2}=1]}{p_{s_2}},\frac{\mathbf{1}[Z_{s_3}=1]}{p_{s_3}}\frac{\mathbf{1}[Z_{s_4}=1]}{p_{s_4}}\right)\right|\leq C\cdot T^{9/32}\enspace.
        \end{align*}
        \item[(3)] If there are exactly three distinct indices among $\{s_1,s_2,s_3,s_4\}$, we have
        \begin{align*}
            &\left|\operatorname{Cov}\left(\frac{\mathbf{1}[Z_{s_1}=1]}{p_{s_1}}\frac{\mathbf{1}[Z_{s_2}=1]}{p_{s_2}},\frac{\mathbf{1}[Z_{s_3}=1]}{p_{s_3}}\frac{\mathbf{1}[Z_{s_4}=1]}{p_{s_4}}\right)\right|\\
            \leq&C\cdot\left(\eta_TT\right)^{81/256}\cdot\left(\frac{R_T}{R_{s_1}}\right)^{1/4}\left(\frac{R_T}{R_{s_2}}\right)^{1/4}\left(\frac{R_T}{R_{s_3}}\right)^{1/4}\left(\frac{R_T}{R_{s_4}}\right)^{1/4}\enspace.
        \end{align*}
    \end{itemize}
    Here $C>0$ is a universal constant.
\end{lemma}

\begin{proof}
    By Corollary \ref{corollary:estimated-squared-residuals-second-moment} and similar arguments as in Corollary \mainref{corollary:p-moment}*, we can prove that for any $0\leq k\leq 8$ and any $t\in[T]$, the following holds:
    \begin{align*}
        \max\left\{\e\left[\frac{1}{p_t^k}\right],\e\left[\frac{1}{(1-p_t)^k}\right]\right\}\lesssim (\eta_tT)^{k/4}\enspace.
    \end{align*}
    By this result, we can sharpen the use of H{\"o}lder's inequality when bounding the moments of inverse probabilities.
    \begin{itemize}
        \item[(1)] If $s_1,s_2,s_3,s_4$ are distinct.
        By the law of iterated expectations, we have
        \begin{align*}
        &\operatorname{Cov}\left(\frac{\mathbf{1}[Z_{s_1}=1]}{p_{s_1}}\frac{\mathbf{1}[Z_{s_2}=1]}{p_{s_2}},\frac{\mathbf{1}[Z_{s_3}=1]}{p_{s_3}}\frac{\mathbf{1}[Z_{s_4}=1]}{p_{s_4}}\right)\\
        =&\e\left[\frac{\mathbf{1}[Z_{s_1}=1]}{p_{s_1}}\frac{\mathbf{1}[Z_{s_2}=1]}{p_{s_2}}\frac{\mathbf{1}[Z_{s_3}=1]}{p_{s_3}}\frac{\mathbf{1}[Z_{s_4}=1]}{p_{s_4}}\right]\\
        &-\e\left[\frac{\mathbf{1}[Z_{s_1}=1]}{p_{s_1}}\frac{\mathbf{1}[Z_{s_2}=1]}{p_{s_2}}\right]\e\left[\frac{\mathbf{1}[Z_{s_3}=1]}{p_{s_3}}\frac{\mathbf{1}[Z_{s_4}=1]}{p_{s_4}}\right]\\
        =&1-1\\
        =&0\enspace.
        \end{align*}
        \item[(2)] By the proof of part (2) in Lemma \ref{lemma:p-moments-cross-term}, we have
        \begin{align*}
            &\left|\operatorname{Cov}\left(\frac{\mathbf{1}[Z_{s_1}=1]}{p_{s_1}}\frac{\mathbf{1}[Z_{s_2}=1]}{p_{s_2}},\frac{\mathbf{1}[Z_{s_3}=1]}{p_{s_3}}\frac{\mathbf{1}[Z_{s_4}=1]}{p_{s_4}}\right)\right|\\
            =&\Bigg|\e\left[\frac{\mathbf{1}[Z_{s_1}=1]}{p_{s_1}^2}\frac{\mathbf{1}[Z_{s_2}=1]}{p_{s_2}^2}\right]-\underbrace{\left(\e\left[\frac{\mathbf{1}[Z_{s_1}=1]}{p_{s_1}}\frac{\mathbf{1}[Z_{s_2}=1]}{p_{s_2}}\right]\right)^2}_{=1}\Bigg|\\
            \leq&\e\left[\frac{\mathbf{1}[Z_{s_1}=1]}{p_{s_1}^2}\frac{\mathbf{1}[Z_{s_2}=1]}{p_{s_2}^2}\right]+1\\
            \lesssim&T^{9/32}+1\\
            \lesssim&T^{9/32}\enspace.
        \end{align*}
        \item[(3)] If there are exactly three distinct indices among $s_1,s_2,s_3,s_4$, we consider the following cases:\\[3mm]
        \textbf{Case 1: }$s_1>s_2=s_3>s_4$ or $s_3>s_4=s_1>s_2$.\\[3mm]
        By symmetry, we only consider the case where $s_1>s_2=s_3>s_4$.
        Using the law of iterated expectations, we have
        \begin{align*}
            &\left|\operatorname{Cov}\left(\frac{\mathbf{1}[Z_{s_1}=1]}{p_{s_1}}\frac{\mathbf{1}[Z_{s_2}=1]}{p_{s_2}},\frac{\mathbf{1}[Z_{s_3}=1]}{p_{s_3}}\frac{\mathbf{1}[Z_{s_4}=1]}{p_{s_4}}\right)\right|\\
            =&\left|\e\left[\frac{\mathbf{1}[Z_{s_1}=1]}{p_{s_1}}\frac{\mathbf{1}[Z_{s_2}=1]}{p_{s_2}^2}\frac{\mathbf{1}[Z_{s_4}=1]}{p_{s_4}}\right]-1\right|\\
            =&\left|\e\left[\frac{\mathbf{1}[Z_{s_2}=1]}{p_{s_2}^2}\frac{\mathbf{1}[Z_{s_4}=1]}{p_{s_4}}\right]-1\right|\\
            \leq&\e\left[\frac{1}{p_{s_2}}\frac{\mathbf{1}[Z_{s_4}=1]}{p_{s_4}}\right]+1\\
            \leq&\left(\e\left[\frac{1}{p_{s_2}^8}\right]\right)^{1/8}\left(\e\left[\frac{\mathbf{1}[Z_{s_4}=1]}{p_{s_4}^{8/7}}\right]\right)^{7/8}
            +1\quad(\text{H{\"o}lder's inequality})\\
            =&\left(\e\left[\frac{1}{p_{s_2}^8}\right]\right)^{1/8}\left(\e\left[\frac{1}{p_{s_4}^{1/7}}\right]\right)^{7/8}
            +1\\
            \lesssim&(\eta_{s_2}T)^{\frac{1}{4}\cdot 8\cdot \frac{1}{8}}(\eta_{s_4}T)^{\frac{1}{4}\cdot\frac{1}{7}\cdot\frac{7}{8}}\\
            =&(\eta_TT)^{9/32}\cdot\left(\frac{R_T}{R_{s_2}}\right)^{1/4}\left(\frac{R_T}{R_{s_4}}\right)^{1/32}\\
            \leq&(\eta_TT)^{9/32}\cdot\left(\frac{R_T}{R_{s_1}}\right)^{1/4}\left(\frac{R_T}{R_{s_2}}\right)^{1/4}\left(\frac{R_T}{R_{s_3}}\right)^{1/4}\left(\frac{R_T}{R_{s_4}}\right)^{1/4}\quad(\text{since $R_T\geq R_{s_1},R_{s_2},R_{s_3},R_{s_4}$})\enspace.
        \end{align*}
        \textbf{Case 2: }$s_1=s_3$ and $s_2\neq s_4$.\\[3mm]
        By symmetry, we only consider the case where $s_2>s_4$.
        Using the law of iterated expectations, we have
        \begin{align*}
        &\left|\operatorname{Cov}\left(\frac{\mathbf{1}[Z_{s_1}=1]}{p_{s_1}}\frac{\mathbf{1}[Z_{s_2}=1]}{p_{s_2}},\frac{\mathbf{1}[Z_{s_3}=1]}{p_{s_3}}\frac{\mathbf{1}[Z_{s_4}=1]}{p_{s_4}}\right)\right|\\
        =&\left|\e\left[\frac{\mathbf{1}[Z_{s_1}=1]}{p_{s_1}^2}\frac{\mathbf{1}[Z_{s_2}=1]}{p_{s_2}}\frac{\mathbf{1}[Z_{s_4}=1]}{p_{s_4}}\right]-1\right|\\
        \leq&\e\left[\frac{1}{p_{s_1}}\frac{\indicator{Z_{s_2}=1}}{p_{s_2}}\frac{\indicator{Z_{s_4}=1}}{p_{s_4}}\right]+1\\
        \leq&\left(\e\left[\frac{1}{p_{s_1}^8}\right]\right)^{1/8}\left(\e\left[\frac{\indicator{Z_{s_2}=1}}{p_{s_2}^{8/7}}\frac{\indicator{Z_{s_4}=1}}{p_{s_4}^{8/7}}\right]\right)^{7/8}+1\quad(\text{H{\"o}lder's inequality})\\
        =&\left(\e\left[\frac{1}{p_{s_1}^8}\right]\right)^{1/8}\left(\e\left[\frac{1}{p_{s_2}^{1/7}}\frac{\indicator{Z_{s_4}=1}}{p_{s_4}^{8/7}}\right]\right)^{7/8}+1\\
        \leq&\left(\e\left[\frac{1}{p_{s_1}^8}\right]\right)^{1/8}\left(\e\left[\frac{1}{p_{s_2}^8}\right]\right)^{\frac{7}{8}\cdot \frac{1}{56}}\left(\e\left[\frac{\indicator{Z_{s_4}=1}}{p_{s_4}^{\frac{8}{7}\cdot\frac{56}{55}}}\right]\right)^{\frac{7}{8}\cdot \frac{55}{56}}+1\quad(\text{H{\"o}lder's inequality})\\
        =&\left(\e\left[\frac{1}{p_{s_1}^8}\right]\right)^{1/8}\left(\e\left[\frac{1}{p_{s_2}^8}\right]\right)^{\frac{7}{8}\cdot \frac{1}{56}}\left(\e\left[\frac{1}{p_{s_4}^{\frac{9}{55}}}\right]\right)^{\frac{7}{8}\cdot \frac{55}{56}}+1\\
        \lesssim&(\eta_{s_1}T)^{\frac{1}{4}(8\cdot \frac{1}{8})}(\eta_{s_2}T)^{\frac{1}{4}(8\cdot \frac{7}{8}\cdot \frac{1}{56})}(\eta_{s_4}T)^{\frac{1}{4}(\frac{9}{55}\cdot \frac{7}{8}\cdot \frac{55}{56})}\\
        =&\left(\eta_TT\right)^{81/256}\cdot \left(\frac{R_{T}}{R_{s_1}}\right)^{1/4}\left(\frac{R_{T}}{R_{s_2}}\right)^{1/32}\left(\frac{R_{T}}{R_{s_4}}\right)^{9/256}\\
        \leq&\left(\eta_TT\right)^{81/256}\cdot\left(\frac{R_T}{R_{s_1}}\right)^{1/4}\left(\frac{R_T}{R_{s_2}}\right)^{1/4}\left(\frac{R_T}{R_{s_3}}\right)^{1/4}\left(\frac{R_T}{R_{s_4}}\right)^{1/4}\quad(\text{since $R_T\geq R_{s_1},R_{s_2},R_{s_3},R_{s_4}$})\enspace.
    \end{align*}
    \textbf{Case 3: }$s_2=s_4$ and $s_1\neq s_3$.\\[3mm]
    By symmetry, we only consider the case where $s_1>s_3$.
    Using the law of iterated expectations, we have
    \begin{align*}
        &\left|\operatorname{Cov}\left(\frac{\mathbf{1}[Z_{s_1}=1]}{p_{s_1}}\frac{\mathbf{1}[Z_{s_2}=1]}{p_{s_2}},\frac{\mathbf{1}[Z_{s_3}=1]}{p_{s_3}}\frac{\mathbf{1}[Z_{s_4}=1]}{p_{s_4}}\right)\right|\\
        =&\left|\e\left[\frac{\mathbf{1}[Z_{s_1}=1]}{p_{s_1}}\frac{\mathbf{1}[Z_{s_3}=1]}{p_{s_3}}\frac{\mathbf{1}[Z_{s_4}=1]}{p_{s_4}^2}\right]-1\right|\\
        =&\left|\e\left[\frac{\mathbf{1}[Z_{s_3}=1]}{p_{s_3}}\frac{\mathbf{1}[Z_{s_4}=1]}{p_{s_4}^2}\right]-1\right|\\
        =&\left|\e\left[\frac{\mathbf{1}[Z_{s_4}=1]}{p_{s_4}^2}\right]-1\right|\\
        \leq&\e\left[\frac{1}{p_{s_4}}\right]+1\\
        \lesssim&(\eta_{s_4}T)^{1/4}\\
        =&(\eta_TT)^{1/4}\cdot\left(\frac{R_T}{R_{s_4}}\right)^{1/4}\\
        \leq&\left(\eta_TT\right)^{81/256}\cdot\left(\frac{R_T}{R_{s_1}}\right)^{1/4}\left(\frac{R_T}{R_{s_2}}\right)^{1/4}\left(\frac{R_T}{R_{s_3}}\right)^{1/4}\left(\frac{R_T}{R_{s_4}}\right)^{1/4}\quad(\text{since $R_T\geq R_{s_1},R_{s_2},R_{s_3},R_{s_4}$})\enspace.
    \end{align*}
    The discussions in these three cases complete the proof. \qedhere
    \end{itemize}
\end{proof}

Based on Lemma \ref{lemma:deterministic-summation-VB} and Lemma \ref{lemma:refined-p-bound}, we now bound the variance of the estimated squared OLS residuals in the following theorem.

\begin{reftheorem}{\mainref{thm:sq-ols-res-estimates}}
	\sqolsresestimates
\end{reftheorem}

\begin{proof}
    We only prove the result for $k=1$. 
    We first verify the unbiasedness of $\widehat{A}(1)$.
    By the law of iterated expectations, we have
    \begin{align*}
        &\E{\widehat{A}(1)}\\
        =&\frac{1}{T}\e\left[\sum_{t=1}^TQ_{t,t}y_t(1)^2\cdot \frac{\mathbf{1}[Z_t=1]}{p_t}+\sum_{1\leq t_1\neq t_2\leq T}Q_{t_1,t_2}y_{t_1}(1)y_{t_2}(1)\cdot \frac{\mathbf{1}[Z_{t_1}=1]}{p_{t_1}}\frac{\mathbf{1}[Z_{t_2}=1]}{p_{t_2}}\right]\\
        =&\frac{1}{T}\sum_{t=1}^T\e\left[Q_{t,t}y_t(1)^2\cdot \frac{\mathbf{1}[Z_t=1]}{p_t}\right]+\frac{2}{T}\sum_{1\leq t_1<t_2\leq T}\e\left[\e\left[Q_{t_1,t_2}y_{t_1}(1)y_{t_2}(1)\cdot \frac{\mathbf{1}[Z_{t_1}=1]}{p_{t_1}}\frac{\mathbf{1}[Z_{t_2}=1]}{p_{t_2}}\Big|\filt_{t_2-1}\right]\right]\\
        =&\frac{1}{T}\sum_{t=1}^TQ_{t,t}y_t(1)^2+\frac{2}{T}\sum_{1\leq t_1<t_2\leq T}\e\left[Q_{t_1,t_2}y_{t_1}(1)y_{t_2}(1)\cdot \frac{\mathbf{1}[Z_{t_1}=1]}{p_{t_1}}\right]\\
        =&\frac{1}{T}\sum_{t=1}^TQ_{t,t}y_t(1)^2+\frac{2}{T}\sum_{1\leq t_1<t_2\leq T}Q_{t_1,t_2}y_{t_1}(1)y_{t_2}(1)\\
        =&\frac{1}{T}\vec{y}(1)^\tran \mat{Q} \vec{y}(1)\\
        =&A(1)\enspace.
    \end{align*}
    Now we turn to bounding the variance of $\widehat{A}(1)$. 
    Using the inequality $\operatorname{Var}(X+Y)\leq 2\operatorname{Var}(X)+2\operatorname{Var}(Y)$, we have:
    \begin{align*}
        &\Var{\widehat{A}(1)}\\
        =&\frac{1}{T^2}\operatorname{Var}\left(\sum_{t=1}^TQ_{t,t}y_t(1)^2\cdot \frac{\mathbf{1}[Z_t=1]}{p_t}+\sum_{1\leq t_1\neq t_2\leq T}Q_{t_1,t_2}y_{t_1}(1)y_{t_2}(1)\cdot \frac{\mathbf{1}[Z_{t_1}=1]}{p_{t_1}}\frac{\mathbf{1}[Z_{t_2}=1]}{p_{t_2}}\right)\\
        \lesssim&\underbrace{\frac{1}{T^2}\operatorname{Var}\left(\sum_{t=1}^TQ_{t,t}y_t(1)^2\cdot \frac{\mathbf{1}[Z_t=1]}{p_t}\right)}_{:=S_1}+\underbrace{\frac{1}{T^2}\operatorname{Var}\left(\sum_{1\leq t_1\neq t_2\leq T}Q_{t_1,t_2}y_{t_1}(1)y_{t_2}(1)\cdot \frac{\mathbf{1}[Z_{t_1}=1]}{p_{t_1}}\frac{\mathbf{1}[Z_{t_2}=1]}{p_{t_2}}\right)}_{:=S_2}\enspace.
    \end{align*}
    We first derive an upper bound on $S_1$.
    By similar arguments as in the proof of Lemma \ref{lemma:p-moments-cross-term}, we can show that
    \begin{align*}
        &\max_{1\leq t\leq T}\operatorname{Var}\left(\frac{\mathbf{1}[Z_{t}=1]}{p_{t}}-1\right)=\bigO{T^{1/8}}\enspace.
    \end{align*}
    Then by Lemma \ref{lemma:deterministic-summation-VB}, we obtain
    \begin{align*}
        S_1\lesssim& T^{-2}\cdot T\cdot T^{1/8}=T^{-7/8}\enspace.
    \end{align*}
    We then derive an upper bound on $S_2$ by expanding the variance. 
    Since Lemma \ref{lemma:refined-p-bound} implies that the covariance equals 0 when $t_1,t_2,t_3,t_4$ are all distinct, we have
    \begin{align*}
        &S_2\\
        =&\frac{1}{T^2}\sum_{1\leq t_1\neq t_2\leq T}\sum_{1\leq t_3\neq t_4\leq T}\operatorname{Cov}\Bigg(Q_{t_1,t_2}y_{t_1}(1)y_{t_2}(1)\cdot \frac{\mathbf{1}[Z_{t_1}=1]}{p_{t_1}}\frac{\mathbf{1}[Z_{t_2}=1]}{p_{t_2}},\\
        &\hspace{1.8in}Q_{t_3,t_4}y_{t_3}(1)y_{t_4}(1)\cdot \frac{\mathbf{1}[Z_{t_3}=1]}{p_{t_3}}\frac{\mathbf{1}[Z_{t_4}=1]}{p_{t_4}}\Bigg)\\
        \lesssim&\underbrace{\frac{1}{T^2}\sum_{1\leq t_1\neq t_2\leq T}Q_{t_1,t_2}^2y_{t_1}(1)^2y_{t_2}(1)^2\cdot\operatorname{Var}\left(\frac{\mathbf{1}[Z_{t_1}=1]}{p_{t_1}}\frac{\mathbf{1}[Z_{t_2}=1]}{p_{t_2}}\right)}_{:=S_{2,1}~(\text{exactly two distinct indices among }t_1,t_2,t_3,t_4)}\\
        +&\underbrace{\frac{1}{T^2}\sum_{1\leq t_1\neq t_2\neq t_3\leq T}|Q_{t_1,t_2}||Q_{t_1,t_3}|y_{t_1}(1)^2|y_{t_2}(1)||y_{t_3}(1)|\left|\operatorname{Cov}\left(\frac{\mathbf{1}[Z_{t_1}=1]}{p_{t_1}}\frac{\mathbf{1}[Z_{t_2}=1]}{p_{t_2}},\frac{\mathbf{1}[Z_{t_1}=1]}{p_{t_1}}\frac{\mathbf{1}[Z_{t_3}=1]}{p_{t_3}}\right)\right|}_{:=S_{2,2}~(\text{exactly three distinct indices among }t_1,t_2,t_3,t_4)}\enspace.
    \end{align*}
    We then derive upper bounds on $S_{2,1}$ and $S_{2,2}$ separately.
    By Lemma \ref{lemma:deterministic-summation-VB} and Lemma \ref{lemma:refined-p-bound}, we obtain
    \begin{align}
        S_{2,1}\lesssim& \frac{1}{T^2}\cdot R_T^2\cdot T^{9/32}=T^{-55/32}R_T^2\enspace.
    \end{align}
    By Lemma \ref{lemma:deterministic-summation-VB} and Lemma \ref{lemma:refined-p-bound} we obtain
    \begin{align*}
        |S_{2,2}|\leq &\frac{1}{T^2}\sum_{1\leq t_1\neq t_2\neq t_3\leq T}|Q_{t_1,t_2}||Q_{t_1,t_3}|y_{t_1}(1)^2|y_{t_2}(1)||y_{t_3}(1)|\\
        &\times\underbrace{\left|\operatorname{Cov}\left(\frac{\mathbf{1}[Z_{t_1}=1]}{p_{t_1}}\frac{\mathbf{1}[Z_{t_2}=1]}{p_{t_2}},\frac{\mathbf{1}[Z_{t_1}=1]}{p_{t_1}}\frac{\mathbf{1}[Z_{t_3}=1]}{p_{t_3}}\right)\right|}_{\text{bound by Lemma \ref{lemma:refined-p-bound}}}\\
        \lesssim&\frac{1}{T^2}\sum_{1\leq t_1\neq t_2\neq t_3\leq T}|Q_{t_1,t_2}||Q_{t_1,t_3}|y_{t_1}(1)^2|y_{t_2}(1)||y_{t_3}(1)|\cdot(\eta_TT)^{81/256}\left(\frac{R_{T}}{R_{t_1}}\right)^{1/2}\left(\frac{R_{T}}{R_{t_2}}\right)^{1/4}\left(\frac{R_{T}}{R_{t_3}}\right)^{1/4}\\
        \lesssim&T^{-2}(\eta_TT)^{81/256}\cdot T\cdot\underbrace{\sum_{1\leq t_1\neq t_2\leq T}Q_{t_1,t_2}^2y_{t_1}(1)^2y_{t_2}(1)^2\left(\frac{R_{T}}{R_{t_1}}\right)^{1/2}\left(\frac{R_{T}}{R_{t_2}}\right)^{1/2}}_{\text{bound by Lemma \ref{lemma:deterministic-summation-VB}}}\quad(\text{AM-GM inequality})\\
        \lesssim&T^{-1}(\eta_TT)^{81/256}R_T^2\\
        =&T^{-1}R_T^2(TR_T^{-2})^{81/512}\\
        \lesssim&T^{-1}R_T^2(TR_T^{-2})^{1/6}\quad(\text{Assumption \mainref{assumption:maximum-radius-strengthened} and Lemma \ref{lemma:R}})\enspace.
    \end{align*}
    Hence we have
    \begin{align*}
        \Var{\widehat{A}(1)}\lesssim& S_1+S_{2,1}+S_{2,2}\\
        \lesssim&T^{-7/8}+T^{-55/32}R_T^2+T^{-1}R_T^2(TR_T^{-2})^{1/6}\\
        =&T^{-1}R_T^2(TR_T^{-2})^{1/6}\cdot \left(T^{-1/24}R_T^{-5/3}+T^{-23/32}(TR_T^{-2})^{-1/6}+1\right)\\
        \lesssim&T^{-1}R_T^2(TR_T^{-2})^{1/6}\quad(\text{Lemma \ref{lemma:R} and Assumption \mainref{assumption:maximum-radius-strengthened}})\\
        =&T^{-5/6}R_T^{5/3}\\
        \lesssim&(T^{-5/12}R^{5/6})^2\quad(\text{Lemma \ref{lemma:R}})\enspace.
    \end{align*}
    This completes the proof.
\end{proof}

The following corollary follows from Theorem \mainref{thm:sq-ols-res-estimates} and establishes the consistency of the variance bound estimator.

\begin{refcorollary}{\mainref{corollary:vb-is-consistent}}
	\vbconsistent
\end{refcorollary}

\begin{proof}
    Recall that $T\cdot \vb=4\mathcal{E}(1)\mathcal{E}(0)$ and $T\cdot \evb=4\widehat{\mathcal{E}}(1)\widehat{\mathcal{E}}(0)$. By the definition of $\widehat{\mathcal{E}}(k)$ and Theorem \mainref{thm:sq-ols-res-estimates}, we have
    \begin{align*}
        |\widehat{\mathcal{E}}(k)^2-\olsres{k}^2|\leq |\widehat{A}(k)-\olsres{k}^2|=|\widehat{A}(k)-A(k)|=\mathcal{O}_p(T^{-5/12}R^{5/6})\enspace.
    \end{align*} 
    Since $\mathcal{E}(k)=\Theta(1)$ by Assumption \mainref{assumption:moments}, it follows that $|\widehat{\mathcal{E}}(k)-\mathcal{E}(k)|=\mathcal{O}_p(T^{-5/12}R^{5/6})$ using a Taylor expansion of the square root. Then by Assumption \mainref{assumption:moments}, we obtain
    \begin{align*}
        |\widehat{\mathcal{E}}(1)\widehat{\mathcal{E}}(0)-\mathcal{E}(1)\mathcal{E}(0)|\leq&\widehat{\mathcal{E}}(1)\cdot|\widehat{\mathcal{E}}(0)-\mathcal{E}(0)|+\mathcal{E}(0)\cdot|\widehat{\mathcal{E}}(1)-\mathcal{E}(1)|\\
        \leq&(\underbrace{\mathcal{E}(1)}_{\bigO{1}}+\underbrace{|\widehat{\mathcal{E}}(1)-\mathcal{E}(1)|}_{o_p(1)})\cdot|\widehat{\mathcal{E}}(0)-\mathcal{E}(0)|+\underbrace{\mathcal{E}(0)}_{\bigO{1}}\cdot|\widehat{\mathcal{E}}(1)-\mathcal{E}(1)|\\
        =&\mathcal{O}_p(T^{-5/12}R^{5/6})\enspace,
    \end{align*}
    which implies that $T \cdot \evb-T \cdot \vb=\mathcal{O}_p(T^{-5/12}R^{5/6})$.
\end{proof}

\subsection{Wald-type Confidence Intervals}\label{section:D6}

Under the central limit theorem (Theorem \mainref{thm:clt}) and the consistency of the variance estimator (Corollary \mainref{corollary:vb-is-consistent}), we obtain the following result on the asymptotic coverage of the Wald-type confidence intervals.

\begin{refcorollary}{\mainref{corollary:wald-coverage}}
	\wald
\end{refcorollary}

\begin{proof}
    The result can be readily shown by Theorem \mainref{thm:clt} and Corollary \mainref{corollary:vb-is-consistent}.
\end{proof}

	\section{Additional Simulation Results} \label{sec:supp-simulations}

\begin{figure}[t]
	\centering  
	\begin{subfigure}{0.32\textwidth}
		\centering
		\includegraphics[width=\linewidth]{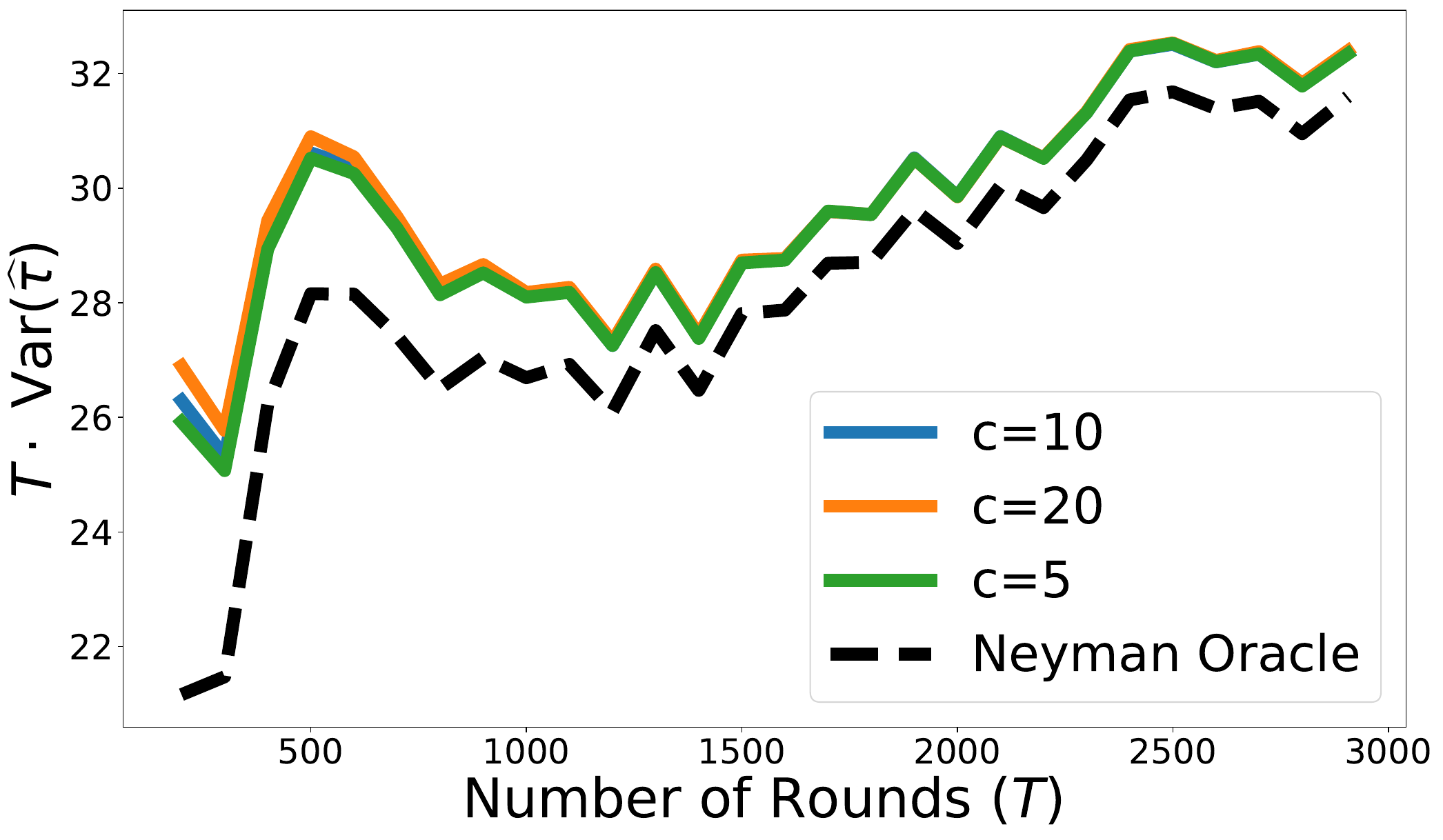}
		\caption{Intercept only}\label{fig:intercept_stepsize}
	\end{subfigure}
	\hfill
	\begin{subfigure}{0.32\textwidth}
		\centering
		\includegraphics[width=\linewidth]{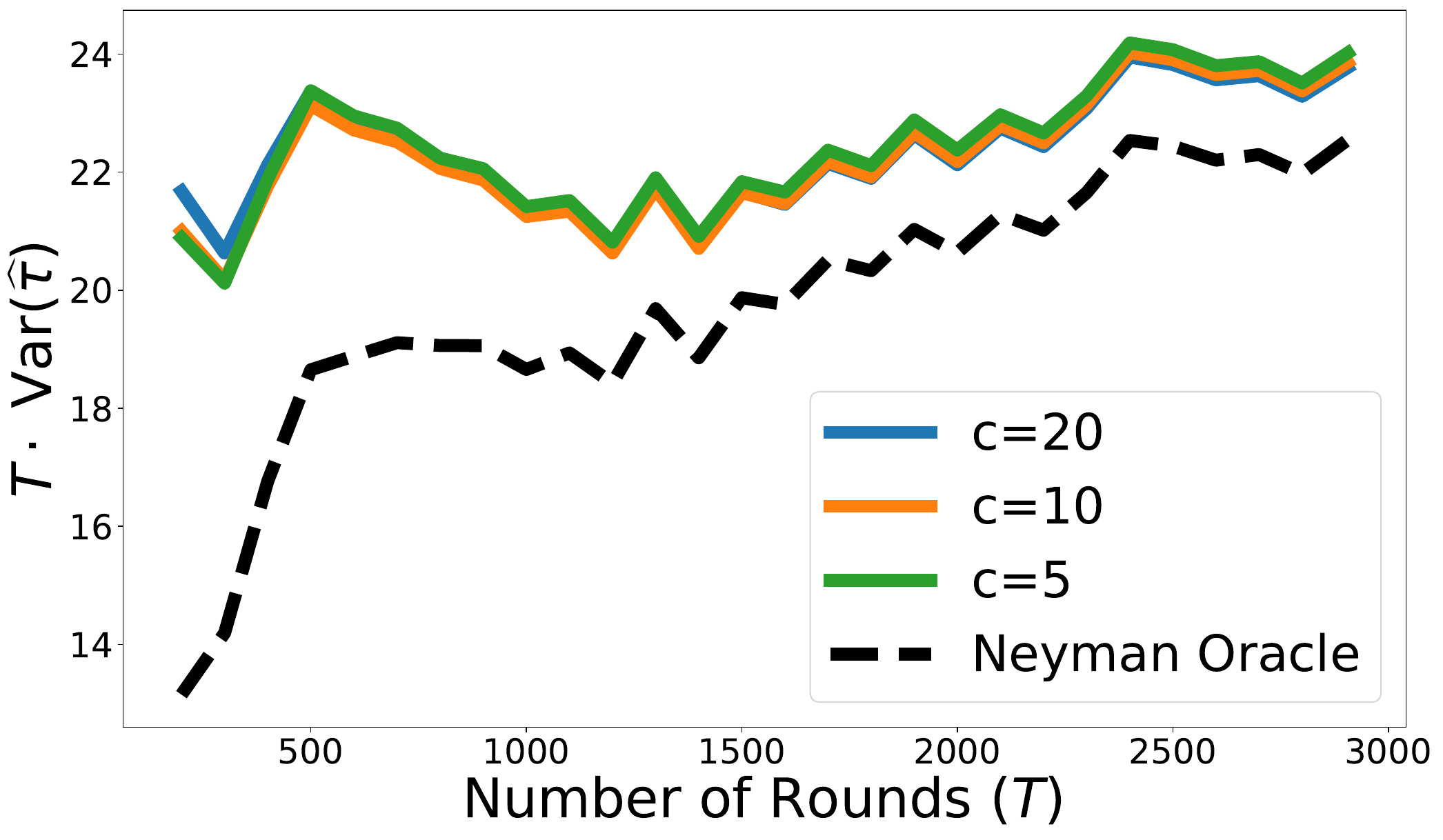}
		\caption{Three covariates}\label{fig:three_covariates_stepsize}
	\end{subfigure}
	\hfill
	\begin{subfigure}{0.32\textwidth}
		\centering
		\includegraphics[width=\linewidth]{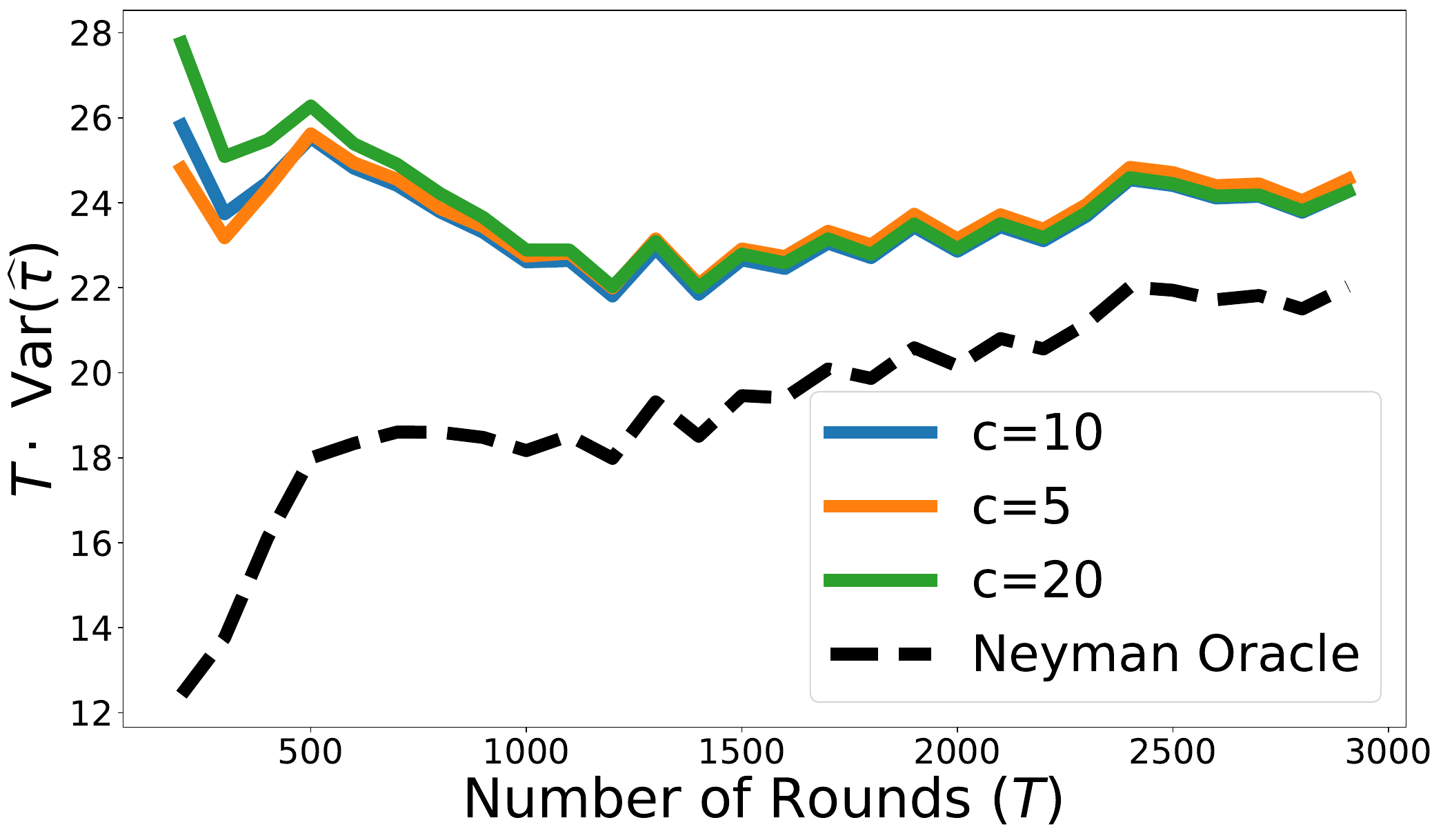}
		\caption{All covariates}\label{fig:all_covariates_stepsize}
	\end{subfigure}
	\caption{
        We analyze the sensitivity of \ourdesign{} to the choice of constant $c$ in the step size $\eta_t = c \cdot T^{-1/2}R_t^{-1}$.
		We evaluate this sensitivity across the three choices of covariates.
		Figure~\ref{fig:intercept_stepsize} shows the intercept alone, Figure~\ref{fig:three_covariates_stepsize} shows three most informative covariates, and Figure~\ref{fig:all_covariates_stepsize} shows all the covariates.}
	\label{fig:stepsize-comparison}
\end{figure}

In this section, we present additional simulation results.
In particular, we investigate the sensitivity of \ourdesign{} to the choice of the constant $c$ in the step size $\eta_t = c \cdot T^{-1/2}R_t^{-1}$.
We use the same simulation set-up as described in Section~\mainref{sec:simulations} of the main body.
In particular, we consider three choices of covariates: (1) the intercept alone, (2) the three most informative covariates (``expectloslessthan3000'', ``gender'', ``profits1''), and (3)  all 12 covariates.

Figure~\ref{fig:stepsize-comparison} shows the effect of using $c \in \setb{5, 10, 20}$ in the step size $\eta_t$ in these three settings.
Across all settings, we see while there is a very small difference in performance for small sample sizes, the choice of step size plays virtually no role in the actual variance once the sample size is at least $T \geq 500$.

\end{document}